\documentclass[a4paper,11pt]{article}

\usepackage[T1]{fontenc}
\usepackage{lmodern}

\usepackage[latin1]{inputenc}
\usepackage{amsmath}
\usepackage{amsthm}
\usepackage{amssymb}
\usepackage{mathrsfs}
\usepackage{graphicx}
\usepackage[all]{xy}
\usepackage{hyperref}

\usepackage{makeidx}

\usepackage{abstract}

\newenvironment{changemargin}[2]{\begin{list}{}{%
\setlength{\topsep}{0pt}%
\setlength{\leftmargin}{0pt}%
\setlength{\rightmargin}{0pt}%
\setlength{\listparindent}{\parindent}%
\setlength{\itemindent}{\parindent}%
\setlength{\parsep}{0pt plus 1pt}%
\addtolength{\leftmargin}{#1}%
\addtolength{\rightmargin}{#2}%
}\item }{\end{list}}

\newtheorem{thm}{Theorem}[section]
\newtheorem{cor}[thm]{Corollary}
\newtheorem{stat}[thm]{Statement}
\newtheorem{claim}[thm]{Claim}
\newtheorem{fact}[thm]{Fact}
\newtheorem{construction}[thm]{Construction}
\newtheorem{lemma}[thm]{Lemma}
\newtheorem{prop}[thm]{Proposition}

\theoremstyle{definition}
\newtheorem{definition}[thm]{Definition}
\newtheorem{ex}[thm]{Example}
\newtheorem{remark}[thm]{Remark}
\newtheorem{question}[thm]{Question}

\newtheorem{problem}[thm]{Problem}

\newcommand{\link}{\operatorname{link}}
\newcommand{\st}{\operatorname{star}}

\newcommand{\supp}{\operatorname{supp}}

\newcommand{\X}{X(\Gamma,\mathcal{G})}

\makeindex

\title{Cubical-like geometry of quasi-median graphs and applications to geometric group theory}
\date{\today}
\author{Anthony Genevois}

\begin{document}

\maketitle

\begin{abstract}
The class of quasi-median graphs is a generalisation of median graphs, or equivalently of CAT(0) cube complexes. The purpose of this thesis is to introduce these graphs in geometric group theory. In the first part of our work, we extend the definition of hyperplanes from CAT(0) cube complexes, and we show that the geometry of a quasi-median graph essentially reduces to the combinatorics of its hyperplanes. In the second part, we exploit the specific structure of the hyperplanes to state combination results. The main idea is that if a group acts in a suitable way on a quasi-median graph so that clique-stabilisers satisfy some non-positively curved property $\mathcal{P}$, then the whole group must satisfy $\mathcal{P}$ as well. The properties we are interested in are mainly (relative) hyperbolicity, (equivariant) $\ell^p$-compressions, CAT(0)-ness and cubicality. In the third part, we apply our general criteria to several classes of groups, including graph products, Guba and Sapir's diagram products, some wreath products, and some graphs of groups. Graph products are our most natural examples, where the link between the group and its quasi-median graph is particularly strong and explicit; in particular, we are able to determine precisely when a graph product is relatively hyperbolic.
\end{abstract}

\begin{small}

\medskip \noindent
\begin{center}
\textbf{R\'esum\'e}
\end{center}

\end{small}

\begin{footnotesize}

\medskip
\begin{changemargin}{1cm}{1cm}
La classe des graphes quasi-m\'edians est une g\'en\'eralisation des graphes m\'edians, ou de mani\`ere \'equivalente, des complexes cubiques CAT(0). L'objectif de cette th\`ese est d'introduire ces graphes dans le monde de la th\'eorie g\'eom\'etrique des groupes. Dans un premier temps, nous \'etendons la notion d'hyperplan d\'efinie dans les complexes cubiques CAT(0), et nous montrons que la g\'eom\'etrie d'un graphe quasi-m\'edian se r\'eduit essentiellement \`a la combinatoire de ses hyperplans. Dans la deuxi\`eme partie de notre texte, qui est le c\oe{}ur de la th\`ese, nous exploitons la structure particuli\`ere des hyperplans pour d\'emontrer des r\'esultats de combinaison. L'id\'ee principale est que si un groupe agit d'une bonne mani\`ere sur un graphe quasi-m\'edian de sorte que les stabilisateurs de cliques satisfont une certaine propri\'et\'e $\mathcal{P}$ de courbure n\'egative ou nulle, alors le groupe tout entier doit satisfaire $\mathcal{P}$ \'egalement. Les propri\'et\'es que nous consid\'erons incluent : l'hyperbolicit\'e (\'eventuellement relative), les compressions $\ell^p$ (\'equivariantes), la g\'eom\'etrie CAT(0) et la g\'eom\'etrie cubique. Finalement, la troisi\`eme et derni\`ere partie de la th\`ese est consacr\'ee \`a l'application des crit\`eres g\'en\'eraux d\'emontr\'es pr\'ec\'edemment \`a certaines classes de groupes particuli\`eres, incluant les produits graph\'es, les groupes de diagrammes introduits par Guba et Sapir, certains produits en couronne (permutationnels), et certains graphes de groupes. Les produits graph\'es constituent notre application la plus naturelle, o\`u le lien entre le groupe et son graphe quasi-m\'edian associ\'e est particuli\`erement fort et explicite; en particulier, nous sommes capables de d\'eterminer pr\'ecis\'ement quand un produit graph\'e est relativement hyperbolique. 
\end{changemargin}

\end{footnotesize}

\newpage

\tableofcontents

\newpage

\section{Introduction}

\hspace{0.5cm} CAT(0) cube complexes were introduced by Gromov in his seminal paper \cite{Gromov1987} as a convenient source of examples of CAT(0) and CAT(-1) groups. But their strength really appeared with the recognition of the central role played by the combinatorics of their hyperplanes, initiated by Sageev in his thesis \cite{MR1347406}. Since then, several still open conjectures for CAT(0) spaces were verified for CAT(0) cube complexes, including the (bi)automaticity of cubulated groups \cite{NibloReeves}, the Tits Alternative for groups acting freely on finite-dimensional CAT(0) cube complexes \cite{alternative}, and the Rank Rigidity Conjecture \cite{MR2827012}. Recently, CAT(0) cube complexes were also crucial in the proof of the famous virtual Haken conjecture \cite{MR3104553}. 

Independently, Chepo\"{i} \cite{mediangraphs} and Roller \cite{Roller} realised that the class of CAT(0) cube complexes can be naturally identified with the class of the so-called \emph{median graphs}. These graphs were known by graph theoretists for a long time since they were introduced by Nebesk\'y in 1971 \cite{NebeskyMedian}. Since then, several classes of graphs were introduced as generalisations of median graphs (see for instance \cite{weaklymoduloar} and references therein), including the main subjet of this article, \emph{quasi-median graphs}, which were introduced by Mulder in 1980 \cite{Mulder} and more extensively studied by Bandelt, Mulder and Wilkeit in 1994 \cite{quasimedian}. The ambition of this article is to introduce quasi-median graphs into the world of geometric group theory as natural and powerful objects.

The first part of our work consists in extending the definition of hyperplanes from CAT(0) cube complexes to quasi-median graphs, and then showing that the geometry of quasi-median graphs essentially reduces to the combinatorics of their hyperplanes. The results we prove are quite similar to those which hold for CAT(0) cube complexes. The main difference is that cutting along a hyperplane in a quasi-median graph may produce arbitrarily many connected components (but always at least two), whereas cutting along a hyperplane in a CAT(0) cube complex produces exactly two connected components. In fact, the analogy between quasi-median graphs and CAT(0) cube complexes is so strong that a precise dictionnary is possible, allowing to translate any statement which holds for CAT(0) cube complexes into a statement which is expected to hold for quasi-median graphs; Table \ref{tab:Dictionnary} is an attempt for such a dictionnary. For instance, we prove:

\begin{prop}
Let $X$ be a quasi-median graph.
\begin{itemize}
	\item A hyperplane separates $X$ into at least two connected components, called \emph{sectors}; they are gated subgraphs.
	\item The carrier of a hyperplane decomposes as a Cartesian product of a clique and a quasi-median subgraph.
	\item For every finite collection of pairwise transverse hyperplanes, there exists a prism whose dual hyperplanes are precisely those hyperplanes.
	\item A path in $X$ is a geodesic if and only if it intersects each hyperplane at most once; as a consequence, the distance between two vertices is equal to the number of hyperplanes separating them.
	\item If a vertex $x$ does not belong to a gated subgraph $C$, then there exists a hyperplane separating $x$ from $C$.
\end{itemize}
\end{prop}

Therefore, if one keeps CAT(0) cube complexes in mind, we do not prove any surprising result, but this is precisely this similarity which turns out to be surprising. In fact, several results we prove are already known by graph theorists, possibly in a different language. The main objective of the first part of our work is to introduce a common formalism as the foundation of the rest of our work, but also of our current knowledge about quasi-median graphs.

\begin{table}
	\centering
		\begin{tabular}{|c|c|c|}
		\hline
			\textbf{CAT(0) cube complexes / Median graphs} & \textbf{Quasi-median graphs} \\ 
		\hline
			Hyperplanes & Hyperplanes [\ref{def:hyperplanes}] \\ \hline
			Halfspaces & Sectors [\ref{def:sectorclique}, \ref{def:sectorhyp}] \\ \hline
			Combinatorially convex subcomplexes & Gated subgraphs [\ref{def:gated}] \\ \hline
			Finite subcomplexes & Cubically finite subgraphs [\ref{def:cubicallyfinite}, \ref{cor:cubicallyfinite}] \\ \hline
			Cubes & Prisms [\ref{def:prism}] \\ \hline
			Median vertex & Quasi-median triangle [\ref{def:qmtriangle}] \\ \hline
			Pocsets & Popsets [\ref{def:popset}] \\ \hline
			Spaces with walls & Spaces with partitions [\ref{def:spacepartitions}] \\ \hline
		\end{tabular}
	\caption{Dictionnary between quasi-median graphs and CAT(0) cube complexes}
	\label{tab:Dictionnary}
\end{table}

Once quasi-median graphs are studied, our purpose is to find information on groups acting on them. However, by considering pairs $\{ S,S^c\}$ where $S$ is a sector, quasi-median graphs are naturally spaces with walls, so that a CAT(0) cube complex is associated to any quasi-median graph. As a consequence, admitting a ``nice'' action (eg. a geometric action) on a quasi-median graph produces a similar action on a CAT(0) cube complex (see Proposition \ref{prop:quasimedianimplycubical} for a precise statement). So why do we care about quasi-median graphs? The first reason is that their geometries may be easier to handle than those of the corresponding CAT(0) cube complexes, allowing us to exploit the action further in order to deduce interesting properties of the group. Graph products are the typical examples where such a situation happens; more details will be given below. The second reason is that ``suitable'' actions on quasi-median graphs, which are far from being proper, lead to combination theorems. More precisely, given a group $G$ acting on some quasi-median graph $X$, our strategy is the following: 
\begin{description}
	\item[Step 1.] Fix a convenient set $\mathcal{C}$ of representatives of cliques of $X$.
	\item[Step 2.] For every $C \in \mathcal{C}$, use the action $\mathrm{stab}(C) \curvearrowright C$ to transfer structures from $\mathrm{stab}(C)$ to $C$ (eg. a metric). The most convenient situation is when this action is transitive and free, so that an orbit map provides a bijection. Typically, $\mathcal{C}$ decomposes into two parts $\mathcal{C}= \mathcal{C}_1 \sqcup \mathcal{C}_2$: in the first one, the actions we are interested in are indeed transitive and free; and in the second part, we have no control on the actions, but the cliques are finite. Therefore, we are able to define structures on the cliques of $\mathcal{C}_1$ from their stabilisers, and usually we put trivial structures (eg. discrete metrics) on the cliques of $\mathcal{C}_2$; because they are finite, this does not cause any trouble.
	\item[Step 3.] Next, if the action of $G$ on $X$ is ``well-behaved'', it is possible to extend our collection of structures defined on the cliques of $\mathcal{C}$ to a system of structures, defined on each clique of $X$, which is compatible with the graph structure and which is $G$-invariant. This is not always possible, so we need to define carefully what a ``suitable action'' means, but when the extension is possible, it is unique.
	\item[Step 4.]  From such a system, we define a global structure on $X$ which extends the ``local structures'', and study how the global structure inherits its properties from the local ones.
	\item[Step 5.] Use the action of $G$ on $X$ endowed with the global structure to find information about $G$ (eg. existence of a good action on some metric space). 
\end{description}
We introduce \emph{$\mathcal{C}$-transitive} actions on quasi-median graphs in order to make the extension mentionned in Step 2 possible; and similarly, \emph{topical} actions for the extension mentionned in Step 3. 

\begin{remark}\label{rem:transitivestrong}
In fact, $\mathcal{C}$-transitive actions are stronger than what it is described in Step 2, because we require some control on the finite cliques of $\mathcal{C}_2$. The first reason for this choice is that this generality is sufficient to cover all our applications, so we chose to give a presentation as simple as possible. As a consequence, the results proved in Section \ref{section:topicalactionsI} hold in a framework which is slightly more general. The second reason is that with our stronger hypotheses, we get a canonical way to identify an arbitrary clique of $X$ to a clique of our set of representatives $\mathcal{C}$, and this will be fundamental in a construction described below, we call \emph{inflating the hyperplanes}. 
\end{remark}

The second part of our work develops this strategy. Typically, the results we obtain have the following form: given a property of groups $\mathcal{P}$ and a group $G$ acting topically-transitively on a quasi-median graph, such that the action satisfies some convenient finiteness conditions (depending on $\mathcal{P}$), if clique-stabilisers satisfy $\mathcal{P}$ then so does $G$. The structures we are interested in are mainly metrics and collections of walls, in order to study the geometry of a group and its cubical properties; but we also consider measured wallspaces, spaces with labelled partitions and Lipschitz maps to $L^p$-spaces in order to study a-T-menability, a-$\mathcal{B}$-menability and (equivariant) $\ell^p$-compressions respectively. Thus, we proved criteria for
\begin{itemize}
	\item acting metrically properly on a CAT(0) cube complex (Proposition \ref{prop:CAT0metricallyproper});
	\item acting geometrically on a CAT(0) cube complex (Proposition \ref{prop:cubulatinggeometrically});
	\item being a-T-menable (Proposition \ref{aTmenablegroups});
	\item being a-$\mathcal{B}$-menable (Proposition \ref{aBmenablegroups});
	\item bounding below $\ell^p$-compressions (Proposition \ref{prop:equicompression});
	\item being relatively hyperbolic (Theorem \ref{thm:qmrelativelyhyp}).
\end{itemize}
However, our strategy does not work for many interesting properties, because in Step 2 we define structures on cliques from their stabilisers. But many properties cannot be read directly from the group, for instance being CAT(0). The trick is the following. Let $G$ be a group acting on a quasi-median graph $X$. If a clique-stabiliser $\mathrm{stab}(C)$ acts on another space $Y_C$ (with a point of trivial stabiliser and such that the action $\mathrm{stab}(C) \curvearrowright C$ is transitive and free), then we can identify $C$ with a subset of $Y_C$. Next, we add to $C$ the missing points of $Y_C$, and doing this process in a suitable way for every clique of $X$ produces a new quasi-median graph $Y$ on which $G$ acts. Now, the cliques of $Y$ are naturally identified with some spaces, so that they are naturally endowed with the corresponding structures. Finally, reproducing our strategy from Step 3 produces other combination results. We find criteria for:
\begin{itemize}
	\item acting properly discontinuously on a CAT(0) cube complex (Proposition \ref{prop:properlydiscontinuous});
	\item acting geometrically and virtually specially on some CAT(0) cube complex (Proposition \ref{prop:cocompactspecial});
	\item being CAT(0) (Theorem \ref{thm:producingCAT0groups}).
\end{itemize}
The third and last part of this work is dedicated to the application of the previous criteria to specific classes of groups, which we now describe.

\paragraph{Graph products.} Given a simplicial graph $\Gamma$ and a collection of groups $\mathcal{G}=\{ G_u \mid u \in V(\Gamma) \}$ indexed by the vertices of $\Gamma$, Green \cite{GreenGP} defined the \emph{graph product} $\Gamma \mathcal{G}$ as the quotient 
$$\left( \underset{u \in V(\Gamma)}{\ast} G_u \right) / \langle \langle [g,h]=1, g \in G_u, h \in G_v \ \text{if} \ (u,v) \in E(\Gamma) \rangle \rangle.$$
For instance, the geometry and the combinatorics of graph products were studied in \cite{GPDehnFunction, AntolinDreesen, GPAsymptotic, Bourdon1997, GPcommensurability, HermillerMeier, HsuWise, MeierGP, MeierGPtopology, GPsurfacesub, BNSinvGP, TitsAltGP}. In our work, we are interested in the Cayley graph $\X$ defined from the generating set $\bigsqcup \mathcal{G}$, which turns out to be a quasi-median graph. Moreover, thanks to the normal form proved by Green in her thesis, we understand precisely the geodesics of $\X$, making $\X$ a very convenient geometric model of $\Gamma \mathcal{G}$. It is worth noticing that, for graph products of finitely many finite groups, $\X$ is a Cayley graph obtained from a finite generating set, so that it turns out to be a complete geometric realisation of the graph product $\Gamma \mathcal{G}$. Although we deduce from this observation that $\Gamma \mathcal{G}$ must act geometrically on some CAT(0) cube complex, the precise description of $\X$ allows us to prove easily that $\Gamma \mathcal{G}$ is virtually cocompact special, that it embeds quasi-isometrically into a product of $\chi(\Gamma)$ trees (see Theorem \ref{thm:embedxtrees}), and also to determine exactly when it is hyperbolic, so that we are able to reprove \cite[Theorem A]{MeierGP} which characterizes the hyperbolic graph products (of arbitrary groups). In fact, applying our general criterion on relative hyperbolicity allows us to determine when a graph product is hyperbolic relatively to its factors. Thanks to our description of $\X$, the argument can be improved further, by following our argument in \cite[Theorem 5.24]{coningoff}, in order to characterize precisely when a graph product is relatively hyperbolic. More precisely, we associate to any finite simplicial graph $\Gamma$ and to any collection of groups $\mathcal{G}$ labelled by $V(\Gamma)$, a collection $\mathfrak{I}(\Gamma,\mathcal{G})$ of subgraphs of $\Gamma$. We emphasize the fact that, the construction of $\mathfrak{I}(\Gamma,\mathcal{G})$ is explicit, algorithmic, and does not depend heavily on the vertex-groups: we only need to know which groups of $\mathcal{G}$ are finite, other information about the vertex-groups being unnecessary. Then

\begin{thm}
Let $\Gamma$ be a finite graph and $\mathcal{G}$ a collection of finitely generated groups indexed by $V(\Gamma)$. The graph product $\Gamma \mathcal{G}$ is relatively hyperbolic if and only if $\mathfrak{I}(\Gamma, \mathcal{G}) \neq \{ \Gamma \}$. If so, $\Gamma \mathcal{G}$ is hyperbolic relatively to $\{ \Lambda \mathcal{G} \mid \Lambda \in \mathfrak{I}(\Gamma ,\mathcal{G}) \}$.
\end{thm}

\noindent
Next, by applying our general criteria, we immediately find that the properties of
\begin{itemize}
	\item being CAT(0) (Theorem \ref{thm:GPCAT0});
	\item acting geometrically on a CAT(0) cube complex (Theorem \ref{thm:GPcc});
	\item acting geometrically and virtually specially on a CAT(0) cube complex (Theorem \ref{thm:GPcc});
\end{itemize}
are stable under graph products along finite graphs; that the properties of
\begin{itemize}
	\item acting metrically properly on a CAT(0) cube complex (Proposition \ref{prop:GPmetricallyproper});
	\item being a-T-menable (Corollary \ref{cor:aTaB});
	\item being a-$\mathcal{B}$-menable (Corollary \ref{cor:aTaB});
\end{itemize}
are stable under graph products along countable graphs; and finally that acting properly discontinuously on a CAT(0) cube complex is stable under arbitrary graph products (Theorem \ref{thm:GPccproper}). By applying our criterion about $\ell^p$-compressions, we are also able to reprove \cite[Corollary 4.4]{AntolinDreesen}:

\begin{thm}
Let $\Gamma$ be a finite simplicial graph and $\mathcal{G}$ a collection of finitely generated groups indexed by $V(\Gamma)$. For every $p \geq 1$,
$$\alpha^*_p(\Gamma \mathcal{G}) \geq \min \left( \frac{1}{p}, \min\limits_{G \in \mathcal{G}} \alpha_p^*(G) \right).$$
\end{thm}

Some of the results we prove are already known, but the point is that we are able to include all of them in a unique point of view, the study of quasi-median graphs, and that we are sometimes able to simplify or even improve them. In our opinion, quasi-median geometry is the most natural way to study graph products.

\paragraph{Wreath products.} The \emph{wreath product} of two groups $G$ and $H$ is defined as the semidirect product
$$G \wr H = \left( \bigoplus\limits_{h \in H} G \right) \rtimes H,$$
where $H$ acts on the direct sum by permuting the coordinates. Wreath products are important in geometric group theory because they lead to interesting counterexamples \cite{WreathLip, CornulierFPWP, PerturbationofWP, QIvsLIP}, but their geometries remain widely unknown appart from a few particular cases \cite{RigiditySolLamp}. Motivated by the cubulation of $\mathbb{Z} \wr \mathbb{Z}$ as a diagram group (see \cite[Example 10]{MR1725439} and \cite{MR1978047}), we associate to any action of $H$ on a CAT(0) cube complex (containing a vertex of trivial stabiliser) a quasi-median graph $\mathfrak{W}$, called the \emph{graph of wreaths}, on which $G \wr H$ acts topically-transitively. We refer to the introduction of Section \ref{section:appli3} for a description of the idea of the construction. First of all, by applying our criterion producing proper actions on cube complexes, we are able to reprove \cite[Theorem 5.C.3]{CornulierCommensurated}:

\begin{thm}
Acting properly discontinuously on a CAT(0) cube complex is stable under wreath products. 
\end{thm}

It is worth noticing that, contrary to Cornulier's proof, we get an explicit construction of the CAT(0) cube complex on which the wreath product acts (see Remark \ref{remark:explicitW}). But the new results on wreath products we prove in this paper concern equivariant $\ell^p$-compressions. The paper \cite{MR2271228} written by Arzhantseva, Guba and Sapir, where they show (in particular) that the Hilbert space compression of $\mathbb{Z} \wr \mathbb{Z}$ belongs to the interval $[1/2,3/4]$, motivated a lot of works on finding estimates on the $\ell^p$-compressions of wreath products; see for instance \cite{NaorPeres, NaorYuval1, LiWreathProducts, TesseraWP, WPZwrZ}. In this paper, we adapt ideas from \cite{propertyA}, which were formulated for CAT(0) cube complexes, to quasi-median graphs, and in particular to the graph of wreaths $\mathfrak{W}$ for the study of wreath products. Our main result is the following:

\begin{thm}
Let $G,H$ be two finitely generated groups. Suppose that $H$ acts on a CAT(0) cube complex $X$ with an orbit map which has compression $\alpha$. Then
$$\alpha_p^*(G \wr H) \geq \alpha \cdot TS(X) \cdot \min \left( \frac{1}{p} , \alpha_p^*(G) \right).$$
\end{thm}

This result is not an immediate application of our general criterion, because wreath products may not be quasi-isometrically embedded into the geometric models we construct. Understanding the distorsion of this embedding leads to the introduction of the constant $TS(X)$, depending on the geometry of the CAT(0) cube complex $X$. By noticing that $TS(X)=1$ for any uniformly locally finite hyperbolic CAT(0) cube complex $X$ (see Proposition \ref{prop:TS1}), we deduce our main application:

\begin{cor}
Let $H$ be a hyperbolic group acting geometrically on some CAT(0) cube complex. For every finitely generated group $G$ and every $p \geq 1$, 
$$\alpha_p^*(G \wr H) \geq \min \left( \frac{1}{p}, \alpha_p^*(G) \right),$$
with equality if $H$ is non elementary and $p \in [1,2]$.
\end{cor}

\paragraph{Diagram products.} In \cite{MR1725439}, Guba and Sapir defined a \emph{diagram product} $D(\mathcal{P},\mathcal{G},w)$ as the fundamental group of a 2-complex of groups associated to a semigroup presentation $\mathcal{P}= \langle \Sigma \mid \mathcal{R} \rangle$, a collection of groups $\mathcal{G}$ indexed by the alphabet $\Sigma$, and a base word $w \in \Sigma^+$. By noticing that diagram products can be described as diagrams whose edges are labelled by a new alphabet $\Sigma(\mathcal{G})= \{ (s,g) \mid s \in \Sigma, g \in G_s \}$, we generalize the construction due to Farley of CAT(0) cube complexes associated to diagram groups \cite{MR1978047} in order to produce a quasi-median graph on which $D(\mathcal{P},\mathcal{G},w)$ acts topically-transitively. This allows us to apply our different combination results.

So let $\mathcal{P}= \langle \Sigma \mid \mathcal{R} \rangle$ be a semigroup presentation, $\mathcal{G}$ a collection of groups indexed by $\Sigma$ and $w \in \Sigma^+$ a base word. 
\begin{itemize}
	\item If the groups of $\mathcal{G}$ act properly on CAT(0) cube complexes, then so does the diagram product $D(\mathcal{P}, \mathcal{G},w)$.
\end{itemize}
Moreover, assuming that $\mathcal{P}$ is a finite presentation,
\begin{itemize}
	\item if the groups of $\mathcal{G}$ act metrically properly on CAT(0) cube complexes, then so does the diagram product $D(\mathcal{P}; \mathcal{G},w)$;
	\item if the groups of $\mathcal{G}$ are a-T-menable, then so is the diagram product $D(\mathcal{P}, \mathcal{G},w)$;
	\item if the groups of $\mathcal{G}$ are a-$\mathcal{B}$-menable, then so is the diagram product $D(\mathcal{P}, \mathcal{G},w)$.
\end{itemize} 
And finally, assuming that $\mathcal{P}$ is a finite presentation and that the class $[w]_{\mathcal{P}}$ of $w$ modulo $\mathcal{P}$ is finite,
\begin{itemize}
	\item if the groups of $\mathcal{G}$ act geometrically on CAT(0) cube complexes, then so does the diagram product $D(\mathcal{P}, \mathcal{G},w)$;
	\item if the groups of $\mathcal{G}$ are CAT(0), then so is the diagram product $D(\mathcal{P}, \mathcal{G},w)$.
\end{itemize} 

We are also able to estimate the $\ell^p$-compressions of diagram products when the class of the base word is finite; we will refer to this class of diagram products as \emph{cocompact diagram products}, because this assumption amounts to say that the the quasi-median graph we define contains only finitely many orbits of cliques.

\begin{thm}
Let $\mathcal{P}= \langle \Sigma \mid \mathcal{R} \rangle$ be a semigroup presentation, $\mathcal{G}$ a collection of groups indexed by $\Sigma$ and $w \in \Sigma^+$ a base word. Suppose that $\mathcal{P}$ is a finite presentation, that $[w]_{\mathcal{P}}$ is finite, and that the groups of $\mathcal{G}$ are finitely generated. For every $p \geq 1$, the diagram product $D(\mathcal{P}, \mathcal{G},w)$ satifies
$$\alpha_p^*(D(\mathcal{P}, \mathcal{G},w)) \geq \min \left( \frac{1}{p}, \min\limits_{G \in \mathcal{G}} \alpha_p^*(G) \right).$$
\end{thm}

We also deduce from our general combination results a characterisation of hyperbolic cocompact diagram products. We suspect that our statement holds without the cocompactness assumption; see Question \ref{question:twistedGPhyp} and the related discussion.

\begin{thm}
Let $\mathcal{P}= \langle \Sigma \mid \mathcal{R} \rangle$ be a semigroup presentation, $\mathcal{G}$ a collection of finitely generated groups indexed by $\Sigma$, and $w \in \Sigma^+$ a base word whose class $[w]_{\mathcal{P}}$ is finite. The diagram product $D(\mathcal{P}, \mathcal{G},w)$ is hyperbolic if and only if, for every non empty words $u,v \in \Sigma^+$ such that $w$ is equal to $uv$ modulo $\mathcal{P}$, at least one of the two diagram products $D(\mathcal{P}, \mathcal{G}, u)$ and $D(\mathcal{P}, \mathcal{G}, v)$ is finite. 
\end{thm}

Finally, in order to better understand the structure of diagram groups, we introduce and study another type of actions on quasi-median graphs, called \emph{rotative actions}. Our main result shows that a group admitting such an action decomposes as a semidirect product with a normal factor which is a graph product (see Theorem \ref{thm:rotativeactions} for a precise statement). This structure result is of independent interest, and applies to our other classes of groups acting on quasi-median graphs, but in these cases we do not get new information. Once applied to diagram products, we prove

\begin{thm}
Let $\mathcal{P}= \langle \Sigma \mid \mathcal{R} \rangle$ be a semigroup presentation, $\mathcal{G}$ a collection of finitely generated groups indexed by $\Sigma$, and $w \in \Sigma^+$ a base word. The diagram product decomposes as
$$D(\mathcal{P},\mathcal{G},w) = \Gamma(\mathcal{P}, \mathcal{G},w) \rtimes D(\mathcal{P},w),$$
where $D(\mathcal{P},w)$ is the underlying diagram group and $\Gamma(\mathcal{P},\mathcal{G},w)$ a graph product whose vertex-groups are elements of $\mathcal{G}$.
\end{thm}

We refer to Section \ref{section:DPsemidirect} for a precise description of the graph product $\Gamma(\mathcal{P},\mathcal{G},w)$. As a consequence, when the underlying diagram group is trivial and when the groups of $\mathcal{G}$ are either trivial or infinite cyclic, our theorem produces new examples of right-angled Artin groups which are also diagram groups; concrete examples are given at the end of Section \ref{section:DPsemidirect}. The general question of determining which right-angled Artin groups are diagram groups remains open.

\paragraph{Right-angled graph of groups.} Our last class of groups acting on quasi-median graphs generalises the graph products. A \emph{right-angled graph of groups} is a graph of groups $\mathfrak{G}$ whose vertices are graph products, whose edges are ``subgraph products'', and whose monomorphisms are canonical embeddings (see Section \ref{section:appli4} for a precise definition). Concrete examples are given in Section \ref{section:RAGGex}. We use the normal form of fundamental groupoids of graph of groups proved in \cite{Higgins} to deduce that a natural (connected component of a) Cayley graph of the fundamental groupoid of $\mathfrak{G}$ turns out to be a quasi-median graph. Thus, we make the fundamental group of $\mathfrak{G}$ act on a quasi-median graph. In particular, if all the factors of this product are finite, we find that the fundamental group acts geometrically on a CAT(0) cube complex, and we are able to determine precisely when it is hyperbolic (see Proposition \ref{prop:RAGGfinitehyp} for a precise statement). When some factors are infinite, we would like to show that the action is topical-transitive and then apply our criteria. However, it turns out that our action may not be topical-transitive. This is also a motivation for introducing this class of groups: having examples to test possible future extensions of our work. The obstruction for the action to be topical-transitive is clearly identified: by concatenating monomorphisms of our graph of groups, every factor $G$ arises with a collection of automorphisms $\Phi(G)$, and our action is topical-transitive precisely when $\Phi(G)$ is reduced to $\{ \operatorname{Id} \}$ for every factor $G$. Loosely speaking, the dynamics of the action of the fundamental group of $\mathfrak{G}$ on its associated quasi-median graph is related to the dynamics of collections of automorphisms on factors, and the topical-transitive situation corresponds to the simplest case. We suspect that our results extend to the situation where the $\Phi(G)$'s are finite collections of finite-order automorphisms. 

Let $\mathfrak{G}$ be a right-angled graph of groups such that $\Phi(G)= \{ \operatorname{Id} \}$ for every factor $G$. By applying our criteria, we prove that
\begin{itemize}
	\item if the factors act properly on CAT(0) cube complexes, then so does the fundamental group of $\mathfrak{G}$;
\end{itemize}
assuming that the underlying abstract graph is locally finite and that the simplicial graphs defining the graph products are finite, we are able to deduce that:
\begin{itemize}
	\item if the factors act metrically properly on CAT(0) cube complexes, then so does the fundamental group of $\mathfrak{G}$;
	\item if the factors are a-T-menable, then so is the fundamental group of $\mathfrak{G}$;
	\item if the factors are a-$\mathcal{B}$-menable, then so is the fundamental group of $\mathfrak{G}$;
\end{itemize} 
and finally, if the underlying abstract graph and the simplicial graphs defining the graph products are all finite, then
\begin{itemize}
	\item If the factors act geometrically on CAT(0) cube complexes, then so does the fundamental group of $\mathfrak{G}$;
	\item if the factors are CAT(0), then so is the fundamental group of $\mathfrak{G}$.
\end{itemize} 
(All these statements are proved in Section \ref{section:RAGGtopic}.) Moreover, because the orbit map from the fundamental group of $\mathfrak{G}$ to its associated quasi-median graph turns out to be a quasi-isometric embedding, it is also possible to estimate the equivariant $\ell^p$-compressions of the group. More precisely,

\begin{thm}
Let $\mathfrak{G}$ be a right-angled graph of groups such that $\Phi(G)= \{ \operatorname{Id} \}$ for every factor $G$. Suppose that the underlying abstract graph and the simplicial graphs defining the graph products are all finite, and that the factors are finitely generated. For every $p \geq 1$, the fundamental group $\mathfrak{F}_{\omega}$ of $\mathfrak{G}$ satifies
$$\alpha_p^*(\mathfrak{F}_{\omega}) \geq \min \left( \frac{1}{p}, \min\limits_{\text{$G$ factor}} \alpha_p^*(G) \right).$$
\end{thm}

Finally, we exploit the quasi-median geometry of the fundamental group of $\mathfrak{G}$ in order to determine precisely when it is hyperbolic. The characterisation we obtain is the following (we refer to Section \ref{section:RAGGgeometry} for a definition of the \emph{link} of a factor).

\begin{thm}
Let $\mathfrak{G}$ be a right-angled graph of groups such that $\Phi(G)= \{ \operatorname{Id} \}$ for every factor $G$. Suppose that the factors are finitely generated, and that the underlying abstract graph and the simplicial graphs defining the graph products are all finite. The fundamental group of $\mathfrak{G}$ is hyperbolic if and only if the following conditions are satisfied:
\begin{itemize}
	\item any simplicial graph defining one of our graph products is square-free;
	\item  there do not exist a loop $p$ in the abstract graph, based at some vertex $v \in V$, and two non adjacent vertex-groups $G,H$ in the graph product $G_v$, such that $\varphi_p^G(G)=G$ and $\varphi_p^H(H)=H$;
	\item the link of every factor is finite;
	\item the factors are hyperbolic.
\end{itemize}
\end{thm}

\paragraph{Organization of this paper.} Section \ref{section:QMmain} contains the first part of our work, dedicated to the study of quasi-median graphs from hyperplanes. Section \ref{section:metrizingQM} and Section \ref{section:cubulatingQM} describe how to perform the fourth step of the strategy described above for metrics and wallspaces respectively, while the second and third steps are described in Section \ref{section:topicalactionsI}, in which we also state and prove our first general criteria. Section \ref{section:inflatingintro} describes how to inflate the hyperplanes of a quasi-median graph, and we use this construction in Section \ref{section:topicalactionsII} to find our second family of general criteria. The next four sections are dedicated to applications to graph products, wreath products, diagram groups, and right-angled graph of groups. We conclude this paper by some open problems and questions in Section \ref{section:open}.

\paragraph{Acknowledgement.} I am grateful to J\'er\'emie Chalopin and Victor Chepo\"{i}, for a useful discussion on generalisations of median graphs, which lead me to state and prove the arguments I wrote in an early draft in the context of quasi-median graphs; of course to my advisor, Peter Ha\"{i}ssinsky, for all our discussions during the last three years and for all his comments on my manuscript; and to Vincent Guirardel and Michah Sageev for having accepted to refer my thesis, with a special thought for Vincent and all his useful comments. Finally, I would like to thank Goulnara Arzhantseva, Indira Chatterji, Yves Cornulier and Victor Chepo\"{i} for having accepted to participate to my jury.

\section{Cubical-like geometry of quasi-median graphs}\label{section:QMmain}

\noindent
In this section, we introduce quasi-median graphs and study their geometries. The main idea we want to stress out is that one may define hyperplanes in quasi-median graphs so that arguments holding for CAT(0) cube complexes can be naturally translated to quasi-median graphs.

\subsection{Quasi-median graphs}

\noindent
\emph{Quasi-median graphs} are the main objects studied in this work. Among all the equivalent definitions of these graphs, it seems that the more convenient definition for our purposes is to see quasi-median graphs as particular \emph{weakly modular graphs}. Weakly modular graphs were introduced in \cite{ChepoiTriangles} and \cite{BandeltChepoi}; see also \cite{weaklymoduloar}. Quasi-median graphs appeared independently under different definitions in several places in the litterature; see \cite{quasimedian} and references therein for more information.

\begin{definition}
A graph is \emph{weakly modular}\index{Weakly modular graphs} if it satisfies the following two conditions:
\begin{description}
	\item[(triangle condition)] for any vertex $u$ and any two adjacent vertices $v,w$ at distance $k$ from $u$, there exists a common neighbor $x$ of $v,w$ at distance $k-1$ from $u$;
	\item[(quadrangle condition)] for any vertices $u,z$ at distance $k$ appart and any two neighbors $v,w$ of $z$ at distance $k-1$ from $u$, there exists a common neighbor $x$ of $v,w$ at distance $k-2$ from $u$.
\end{description}
A graph is \emph{quasi-median}\index{Quasi-median graphs} if it weakly modular and does not contain $K_4^-$ and $K_{3,2}$ as induced subgraphs (see Figure \ref{figure7}).
\end{definition}
\begin{figure}
\begin{center}
\includegraphics[scale=0.6]{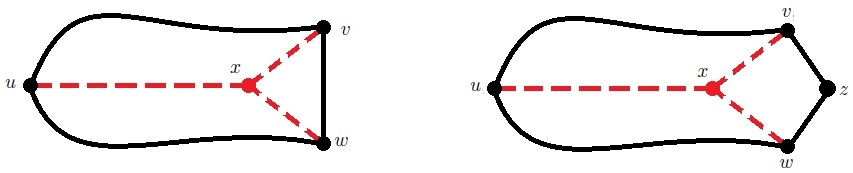}
\end{center}
\caption{The triangle and quadrangle conditions.}
\label{figure8}
\end{figure}
\begin{figure}
\begin{center}
\includegraphics[scale=0.6]{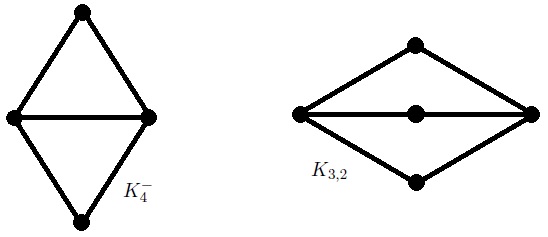}
\end{center}
\caption{The graphs $K_4^-$ and $K_{3,2}$.}
\label{figure7}
\end{figure}
\noindent
In the context of CAT(0) cube complexes, combinatorially convex subcomplexes play an important role. For quasi-median graphs, this role will be played by \emph{gated subgraphs}. 

\begin{definition}\label{def:gated}
Let $X$ be a graph and $Y \subset X$ a subgraph. Fixing a vertex $x \in X$, we say that a vertex $p \in Y$ is a \emph{gate} for $x$ if, for every $y \in Y$, there exists a geodesic between $x$ and $y$ passing through $p$. If any vertex of $X$ admits a gate in $Y$, we say that $Y$ is \emph{gated}\index{Gated subgraphs}.
\end{definition}

\noindent
It is worth noticing that:

\begin{fact}
Let $X$ be a graph and $Y \subset X$ a subgraph.
\begin{itemize}
	\item A gate in $Y$ of some vertex $x \in X$, when it exists, is the unique vertex of $Y$ minimising the distance to $x$.
	\item If $Y$ is gated, then it is convex.
\end{itemize}
\end{fact}

\begin{proof}
Fix a vertex $x \in X$, and suppose that it admits a gate $y \in Y$. For every $z \in Y$, we know that
$$d(x,z) =d(x,y)+d(y,z) \geq d(x,y).$$
Therefore, $y$ minimises the distance to $x$ in $Y$. It also follows from this inequality that, if $z \in Y$ is another vertex minimising the distance to $x$, necessarily $d(y,z)=0$, ie., $y=z$. This proves the first assertion of our lemma. Now, suppose that $Y$ is gated, and fix three vertices $x,y,z \in X$ such that $x,y$ are vertices of $Y$ and $z$ belongs to a geodesic between $x,y$. If $p \in Y$ denotes the gate of $z$, then
$$d(x,y)=d(x,z)+d(z,y)=d(x,p)+2d(p,z)+d(p,y).$$
On the other hand, we know from the triangle inequality that $d(x,y) \leq d(x,p)+d(p,y)$. Consequently, $d(p,z)$ must be zero, so that $z=p \in Y$. 
\end{proof}

\noindent
In general, it is difficult to determine whether or not a subgraph is gated just by applying the definition. Therefore, our first goal is to find a criterion which is easy to verify. We begin with a few definitions.

\begin{definition}
Let $X$ be a graph and $Y \subset X$ a subgraph. If any \emph{triangle} of $X$, ie., any three pairwise adjacent vertices of $X$, which shares an edge with $Y$ is contained into $Y$, we say that $Y$ \emph{contains its triangles}. 
\end{definition}

\begin{definition}
Let $X$ be a graph and $Y \subset X$ a subgraph. We say that $Y$ is \emph{locally convex} if any square in $X$ with two adjacent edges contained in $Y$ is necessarily included into $Y$. 
\end{definition}

\noindent
The following proposition is our main criterion to prove that a subgraph is gated. The statement also follows from different results proved in \cite{ChepoiTriangles}. We give a direct proof for completeness.

\begin{prop}\label{prop:gated}
Let $X$ be a weakly modular graph and $Y \subset Y$ a connected subgraph. Then $Y$ is gated if and only if it is locally convex and contains its triangles.
\end{prop}

\noindent
We begin by proving a weaker statement.

\begin{lemma}\label{gated:step1}
Let $X$ be a weakly modular graph and $Y \subset X$ a connected subgraph. If $Y$ is locally convex and contains its triangles, then it is geodesic (ie., between any two vertices of $Y$ there exists a geodesic in $X$ lying in $Y$).
\end{lemma}

\begin{proof}
Because $Y$ is connected, it is sufficient to prove that, for every vertices $x,y,z \in Y$ with $y,z$ adjacent in $Y$ such that there exists a geodesic between $x$ and $y$ lying in $Y$, necessarily there exists a geodesic between $x$ and $z$ lying in $Y$. We will argue by induction on $d(x,y)$. Of course, for $d(x,y)=0$, there is nothing to prove.

\medskip \noindent
Let $\gamma \subset Y$ be a geodesic between $x$ and $y$. For convenience, set $\ell=d(x,y)$. Notice that $|d(x,z)-d(x,y)| \leq 1$, so three cases may happen.
\begin{enumerate}
	\item Suppose that $d(x,z)=\ell+1$. Then, the concatenation of $\gamma$ with the edge $(y,z)$ between $y$ and $z$ defines a geodesic between $x$ and $z$ lying in $Y$.
	\item Suppose that $d(x,z)=\ell$. By the triangle condition, there exists a vertex $a \in X$ adjacent to both $y$ and $z$ such that $d(x,a)=\ell-1$. Notice that, because $Y$ contains its triangles, the triangle defined by $y,z,a$ must be included into $Y$. Since $d(x,m)<d(x,y)$, we can apply our induction hypothesis to a find a geodesic between $x$ and $a$ which lies in $Y$; if we concatenate it with the edge $(a,z)$, we find a geodesic between $x$ and $z$ lying in $Y$.
	\item Suppose that $d(x,z)=\ell-1$. Let $a$ be the vertex of $\gamma$ defined by $d(x,a)=\ell-1$. Notice first that, if $a=z$, then taking a subsegment of $\gamma$ suffices to produce a geodesic between $x$ and $z$ lying in $Y$. So we will suppose that $a \neq z$. Then the quadrangle condition implies that there exists a vertex $m \in X$ adjacent to both $a$ and $z$ such that $d(x,m)= \ell-2$. Notice that, since $Y$ is locally convex, the square defined by $a,m,y,z$ must be included into $Y$. By applying our induction hypothesis, there exists a geodesic between $x$ and $m$ lying in $Y$, so that, by concatenating it with the edge $(m,z)$, we find a geodesic between $x$ and $z$ which lies in $Y$.
\end{enumerate}
We conclude that $Y$ is indeed geodesic.
\end{proof}

\begin{proof}[Proof of Proposition \ref{prop:gated}.]
It is clear that a gated subgraph must be locally convex and must contain its triangles. So we assume that our subgraph $Y$ is locally convex and that it contains its triangles, and we want to prove that it turns out to be a gated subgraph.

\medskip \noindent
Suppose by contradiction that $Y$ is not gated. Let $x \in X$ be a vertex which has no gate in $Y$. Thus, if $y \in Y$ denotes a vertex minimizing in $Y$ the distance to $x$, there exists a vertex $z \in Y$ such that no geodesic between $x$ and $z$ passes through $y$. Without loss of generality, we may require our counterexample to minimize the quantity $d(y,z)$, so that, for every $w \in Y$ satisfying $d(y,w)<d(y,z)$, there exists a geodesic between $x$ and $w$ passing through $y$. In particular, if we fix a geodesic $(u_1, \ldots, u_r,u_{r+1})$ from $u_1=y$ to $u_{r+1}=z$, then $d(x,u_i)=d(x,y)+d(y,u_i)$ for every $1 \leq i \leq r$. According to the previous lemma, we can choose our geodesic between $y$ and $z$ in $Y$. Because $|d(x,z)-d(x,u_r)| \leq d(z,u_r) =1$, three cases may happen. 

\medskip \noindent
Suppose that $d(x,z)-d(x,u_r)=1$. Then the concatenation of a geodesic between $x$ and $y$ with $(u_1, \ldots, u_{r+1})$ would produce a geodesic between $x$ and $z$ passing through $y$, contradicting our choice of $z$.

\medskip \noindent
Suppose that $d(x,z)-d(x,u_r)=-1$. Because $d(x,u_{r+1})=d(x,u_r)-1=d(x,u_{r-1})$, the quadrangle condition implies that there exists a vertex $v_r$ adjacent to both $u_{r-1}$ and $u_{r+1}$ such that $d(x,v_r)=d(x,u_r)-2$. Moreover, since $Y$ is locally convex, the square defined by $u_r,u_{r+1},v_r,u_{r-1}$ must be included into $Y$, so that $v_r \in Y$. Similarly, because $d(x,v_r)= d(x,u_r)-2 = d(x,u_{r-2})$, the quadrangle condition implies that there exists a vertex $v_{r-1}$ adjacent to both $u_{r-2}$ and $v_r$ such that $d(x,v_{r-1})=d(x,u_{r-1})-2$; moreover, since $Y$ is locally convex, the square defined by $u_{r-1},v_r,v_{r-1},u_{r-2}$ must be included into $Y$, so that $v_{r-1} \in Y$. Then, it is possible to iterate the argument with $u_{r-3}$ and $v_{r-1}$, and so on. Thus, we have constructed a sequence of vertices $v_2, \ldots, v_r \in Y$ satisfying $d(x,v_i)=d(x,u_i)-2$. It is worth noticing that, because $y$ minimizes the distance to $x$ in $Y$, necessarily $r \geq 2$, so at least $v_2$ exists. But then $$d(x,v_2)=d(x,u_2)-2=d(x,u_1)-1=d(x,y)-1,$$ contradicting the choice of $y$.

\medskip \noindent
Finally, suppose that $d(x,z)-d(x,u_r)=0$. So $d(x,u_r)=d(x,u_{r+1})$, and we deduce from the triangle condition that there exists a vertex $m \in X$ adjacent to both $u_{r}$ and $u_{r+1}$ such that $d(x,m)=d(x,u_r)-1$. Notice that, since $Y$ contains its triangles, necessarily $m \in Y$. Finally, we get a contradiction by replacing $u_{r+1}$ with $m$ in the previous argument.
\end{proof}

\noindent
We conclude this section by proving some general properties of quasi-median graphs which will be useful in the sequel. We begin by noticing that gated subgraphs satisfy the Helly's property.

\begin{prop}\label{prop:Helly}
Let $X$ be a graph and $Y_1, \ldots, Y_n \subset X$ a collection of gated subgraphs. If $Y_p \cap Y_q \neq \emptyset$ for every $1 \leq p,q \leq n$, then $\bigcap\limits_{i=1}^n Y_i \neq \emptyset$. 
\end{prop}

\begin{proof}
We will suppose that $n=3$. The general case follows by induction. Fix some vertices $x \in Y_1 \cap Y_2$, $y \in Y_2 \cap Y_3$, $z \in Y_3 \cap Y_1$, and let $m \in Y_3$ be the gate of $x$ in $Y_3$. Because $m$ belongs to a geodesic between $x$ and $y$, we deduce from the convexity of $Y_2$ that $m \in Y_2$; similarly, because $m$ belongs to a geodesic between $x$ and $z$, we deduce from the convexity of $Y_1$ that $m \in Y_3$. Thus, $m$ belongs to $Y_1 \cap Y_2 \cap Y_3$, so that $Y_1 \cap Y_2 \cap Y_3 \neq \emptyset$. 
\end{proof}

\noindent
The following lemma follows easily from the definition of quasi-median graphs. We leave the proof as an exercice.

\begin{lemma}\label{lem:convexquasimedian}
Let $X$ be a quasi-median graph and $Y \subset X$ a subgraph. If $Y$ is a convex subgraph, then it is a quasi-median graph on its own right.
\end{lemma}

\noindent
For our next lemma, recall that a \emph{clique}\index{Cliques} is a maximal complete subgraph.

\begin{lemma}\label{lem:cliquecontainstriangles}
In a quasi-median graph, a clique contains its triangles.
\end{lemma}

\begin{proof}
Let $C$ be a clique, and $x,y,z$ three pairwise adjacent vertices with $x,y \in C$. If $C$ contains no other vertices than $x,y$ and possibly $z$, then in fact $z$ must belong to $C$ since $x,y,z$ define a complete subgraph containing $C$. Otherwise, suppose that there exists some vertex $v \in C$ different from $x,y,z$. Then the subgraph generated by $v,x,y,z$ contains a $K_4^-$; since this subgraph cannot be induced by the definition of quasi-median graphs, we deduce that $v$ and $z$ must be adjacent. Thus, we conclude that $z$ is ajdacent to any vertex of $C$, which implies that $z$ belongs to $C$.
\end{proof}

\noindent
In the next section, we will use this lemma to deduce that cliques are gated subgraphs (see Lemma \ref{lem:cliquegated}). 

\begin{lemma}\label{lem:cliqueinter}
In a quasi-median graph, the intersection between two different cliques is either empty or a single vertex.
\end{lemma}

\begin{proof}
Let $C_1,C_2$ be two cliques intersecting along an edge $e$. Fix two vertices $x_1 \in C_1$, $x_2 \in C_2$ which are not endpoints of $e$. The subgraph generated by $e$, $x_1$ and $x_2$ contains a $K_4^-$; since such a subgraph cannot be induced by definition of a quasi-median graph, we deduce that $x_1$ and $x_2$ must be adjacent. Therefore, the vertices of $C_1$ and $C_2$ are pairwise adjacent, hence $C_1=C_2$.
\end{proof}

\noindent
A direct consequence of the previous lemma is:

\begin{cor}\label{cor:clique}
In a quasi-median graph, there exists a unique clique containing a given edge.
\end{cor}

\noindent
Essentially, our last lemma states that two cliques containing parallel edges must be parallel themselves.

\begin{lemma}\label{lem:parallelcliques}
Let $C_1,C_2 \subset X$ be two distinct cliques. Suppose that there exist two edges $e_1 \subset C_1$, $e_2 \subset C_2$ which are opposite sides of some square in $X$. Then the subgraph generated by $C_1 \cup C_2$ is isomorphic to $C_1 \times [0,1]$, where $C_1= C_1 \times \{0 \}$ and $C_2= C_1 \times \{ 1 \}$. 
\end{lemma}

\begin{proof}
Let $e_i=(x_i,y_i)$ for $i=1,2$, where $x_1$ and $x_2$ (resp. $y_1$ and $y_2$) are adjacent in the square defined by $e_1$ and $e_2$. 

\medskip \noindent
First, notice that $C_1$ and $C_2$ are disjoint. Indeed, suppose that there exists a vertex $v \in C_1 \cap C_2$. Then the vertices $v$, $x_1$ and $x_2$ are pairwise adjacent, hence $x_2 \in C_1$. Similarly, we show that $y_2 \in C_1$, so that the edge $e_2$ belongs both to $C_1$ and $C_2$. We deduce from Corollary \ref{cor:clique} that $C_1=C_2$. Therefore, because we supposed $C_1$ and $C_2$ distinct, they are necessarily disjoint.

\medskip \noindent
Then, we notice that any vertex of $C_1$ is adjacent to some vertex of $C_2$. Let $v \in C_1$ be a vertex, which we will suppose different from $x_1$ and $y_1$ since we already know that $x_1$ and $y_1$ are adjacent to some vertices of $C_2$. We claim that $d(v,x_2)=d(v,y_2)=2$. Indeed, if $v$ is adjacent to $x_2$ or $y_2$, say to $x_2$, then the vertices $v,x_1,x_2$ define a triangle which must be included into $C_1$ since a clique contains its triangles; however, we know that $C_1$ and $C_2$ are disjoint, so this is impossible. Thus, the triangle condition implies that there exists a vertex $m \in X$ adjacent to both $x_2$ and $y_2$ such that $d(v,m)=1$. Because $C_2$ contains its triangles, the triangle defined by the vertices $m,x_2,y_2$ must be included into $C_2$, and in particular $m \in C_2$. We have proved that $v$ is adjacent to some vertex of $C_2$.

\medskip \noindent
On the other hand, because $C_2$ contains its triangles according to Lemma \ref{lem:cliquecontainstriangles} and that $C_1$ and $C_2$ are disjoint, a vertex of $C_1$ may be adjacent to at most one vertex of $C_2$. Therefore, any vertex of $C_1$ is adjacent to exactly one vertex of $C_2$. Similarly, we prove by symmetry that any vertex of $C_2$ is adjacent to exacly one vertex of $C_1$. Our claim follows. 
\end{proof}

\subsection{Hyperplanes and sectors}

\noindent
In this section, we define a notion of hyperplanes in quasi-median graphs by following the definition introduced by Sageev in \cite{MR1347406} for CAT(0) cube complexes. Our goal is to prove that hyperplanes in CAT(0) cube complexes and hyperplanes in quasi-median graphs have essentially the same behaviour. A posteriori, we know that our definition coincides with the transitive closure of the $\theta$-equivalence introduced by Djokovic in \cite{Djokovic}; in particular, a few of the results proved in this section were also proved in \cite{quasimedian} using the $\theta$-equivalence. 

\begin{definition}\label{def:hyperplanes}
A \emph{hyperplane}\index{Hyperplanes} is an equivalence class of edges, where two edges $e$ and $e'$ are said \emph{equivalent} whenever there exists a sequence of edges $e_0=e,e_1, \ldots, e_{n-1},e_n=e'$ such that, for every $1 \leq i \leq n-1$, either $e_i$ and $e_{i+1}$ are opposite sides of some square or they are two sides of some triangle.

\noindent
Alternatively, if we say that two cliques are \emph{parallel} whenever they respectively contain two opposite sides of some square, then a hyperplane is the collection of edges of some class of cliques with respect to the transitive closure of being parallel.

\noindent
One says that an edge or a clique is \emph{dual} to a given hyperplane if it belongs to the associated class of edges. Of course, because two distinct equivalence classes are necessarily disjoint, an edge or a clique is dual to a unique hyperplane. 
\end{definition}

\noindent
From Lemma \ref{lem:parallelcliques}, we know that two cliques $C_1,C_2$ are parallel if and only if there exists an induced subgraph $C \times [0,1]$ where $C\times \{0\} = C_1$ and $C \times \{ 1 \}= C_2$. 

\begin{figure}
\begin{center}
\includegraphics[scale=0.6]{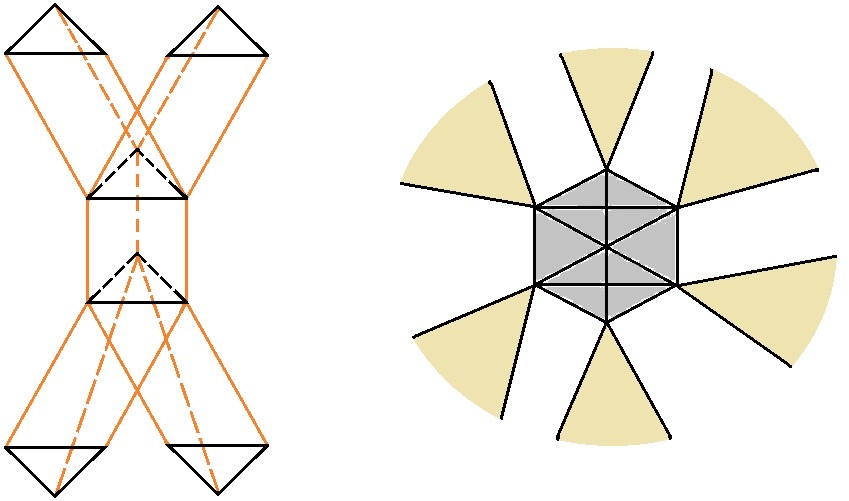}
\end{center}
\caption{The carrier and the sector-decomposition of a hyperplane.}
\label{figure2}
\end{figure}

\medskip \noindent
We sum up the results of this section in the next proposition (see also Figure \ref{figure2}); the needed definitions are given progressively below. The third point is contained in Corollary \ref{cor:sector} and Corollary \ref{cor:sectorgated}; the first point is contained in Lemma \ref{lem:carriergated} and Lemma \ref{lem:carrierproduct}; and finally the second point is contained in Lemma \ref{lem:ccpartialJ}. We also use these results to characterize geodesics, see Proposition \ref{prop:geodesichyp}.

\begin{prop}\label{prop:hypsumup}
Let $X$ be a quasi-median graph and $J$ a hyperplane. 
\begin{itemize}
	\item The carrier $N(J)$ of $J$ is a gated subgraph isomorphic to $F(J) \times C$, where $F(J)$ is the \emph{main fiber} of $J$ and $C$ is a clique. 
	\item Every connected component of $\partial J$, called a \emph{fiber}\index{Fibers of hyperplanes}, is a gated subgraph isomorphic to $F(J)$; in particular, $F(J)$ is a quasi-median graph on its own right.
	\item The hyperplane $J$ separates $X$ into at least two connected components, called \emph{sectors}. They are gated subgraphs.
\end{itemize}
\end{prop}

\noindent
First, we want to define what is a sector.

\begin{lemma}\label{lem:cliquegated}
In a quasi-median graph, any clique is gated.
\end{lemma}

\begin{proof}
Let $X$ be a quasi-median graph and $C \subset X$ a clique. According to Proposition \ref{prop:gated} and Lemma \ref{lem:cliquecontainstriangles}, it is sufficient to prove that $C$ is locally convex in order to deduce that $C$ is gated.

\medskip \noindent
Let $\Pi$ be a square in $X$ with two adjacent edges in $C$. Because $C$ is complete, this implies that $\Pi$ admits a diagonal in $X$, defining a subgraph isomorphic to $K_4^-$; then, since such a subgraph cannot be induced, the second diagonal of $\Pi$ must exist in $X$, so the subgraph generated by $\Pi$ is complete. Because we already know that $C$ contains its triangles, we deduce that $\Pi \subset C$. Therefore, $C$ is locally convex.
\end{proof}

\begin{definition}\label{def:sectorclique}
Let $X$ be a quasi-median graph, $C \subset X$ a clique and $v \in C$ a vertex. The subgraph $[C,v]$ generated by the set of all the vertices of $X$ whose gate in $C$ is $v$ is a \emph{sector}. We refer to the collection $\{ [C,v] \mid v \in C \}$ as the \emph{sector decomposition} defined by $C$. 
\end{definition}

\begin{lemma}\label{lem:sectordecomposition}
Let $C_1,C_2 \subset X$ be two cliques. If $C_1$ and $C_2$ are dual to the same hyperplane, their sector decompositions are the same.
\end{lemma}

\begin{proof}
If $C_1=C_2$ there is nothing to prove. Otherwise, there exists a sequence of cliques between $C_1$ and $C_2$ such that any two successive cliques are parallel. We will assume that $C_1$ and $C_2$ are parallel, the general case following by induction.

\medskip \noindent
We want to prove the following statement, which is sufficient to conclude. For every vertex $x \in C_1$, if $p$ denotes the vertex of $C_2$ opposite to $x$, then $[C_1,x]=[C_2,p]$. In fact, we only need to prove the inclusion $[C_1,x] \subset [C_2,p]$, the reverse inclusion following by symmetry. More precisely, for any $a \in [C_1,x]$ and $q \in C_2 \backslash \{ p \}$, we claim that $d(a,p)<d(a,q)$, meaning that $p$ is the gate of $a$ in $C_2$. 

\medskip \noindent
Let $y \in C_1$ denote the vertex opposite to $q$. Notice that neither $x$ and $q$ nor $y$ and $p$ are adjacent, since a clique must contain its triangles and that $C_1$ and $C_2$ are disjoint. For convenience, let $k=d(a,x)$. Because $|d(a,p)-d(a,x)| \leq d(p,x)=1$, we know that $d(a,p) \in \{k-1,k,k+1 \}$. Notice moreover that $$d(a,x)+1=d(a,x)+d(x,y)= d(a,y) \leq d(a,q)+d(q,y)=d(a,q)+1$$ implies $d(a,q) \geq d(a,x)=k$. Therefore, if we suppose that $d(a,q) \leq d(a,p)$, only three cases may happen.

\medskip \noindent
Suppose that $d(a,p)=k$ and $d(a,q)=k$. Because $d(a,x)=k=d(a,q)$ and $d(a,y)=k+1$, the quadrangle condition implies that there exists a vertex $m \in X$ adjacent to both $x$ and $q$ such that $d(a,m)=k-1$. In particular, $d(a,m)=k-1$ implies that $m \notin \{ p,q \}$. Notice that the vertices $x,y,p,q,m$ define a subgraph isomorphic to $K_{2,3}$, so it cannot be induced. However, we know that $x$ and $q$ are not adjacent, and $$k+1 = d(a,y) \leq d(a,m)+d(m,y) = k-1+d(m,y)$$ implies $d(m,y)\geq 2$, ie., $m$ and $y$ are not adjacent. Therefore, either $y$ and $p$ or $m$ and $p$ must be adjacent. If $y$ and $p$ are adjacent, then the vertices $p,q,x,y$ define an induced subgraph isomorphic to $K_4^-$; similarly, if $m$ and $p$ are adjacent, then the vertices $p,q,x,m$ define an induced subgraph isomorphic to $K_4^-$. Thus, we get a contradiction.

\medskip \noindent
Suppose that $d(a,p)=k+1$ and $d(a,q)=k$. Because $d(a,x)=k=d(a,q)$ and $d(a,y)=k+1$, the quadrangle condition implies that there exists a vertex $m \in X$ adjacent to both $x$ and $q$ such that $d(a,m)=k-1$. Then exactly the same argument as above produces a contradiction.

\medskip \noindent
Suppose that $d(a,p)=k+1$ and $d(a,q)=k+1$. By applying the triangle condition, we find a vertex $r$ adjacent to both $p$ and $q$ such that $d(a,r)=k$. Since the clique $C_2$ contains its triangles, necessarily $r \in C_2$. Let $s$ denote the vertex of $C_1$ opposite to $r$. A fortiori, $d(a,s)=d(a,x)+d(x,s)=k+1$. Thus, by replacing $y$ and $q$ with $s$ and $r$ respectively, we produce the siutation that occurs in our first case. In particular, we deduce a contradiction as well.

\medskip \noindent
As a consequence, we conclude that necessarily $d(a,q) > d(a,p)$.
\end{proof}

\begin{definition}\label{def:sectorhyp}
Let $X$ be a quasi-median graph and $J$ a hyperplane. If $C$ is a clique dual to $J$, the \emph{sectors delimited by $J$}\index{Sectors} are the sectors delimited by $C$. 
\end{definition}

\noindent
According to Lemma \ref{lem:sectordecomposition}, this definition does not depend on the choice of the clique $C$. 

\begin{lemma}\label{lem:samesector}
Let $X$ be a quasi-median graph, $J$ a hyperplane and $e \subset X$ an edge. The endpoints of $e$ belong to the same sector delimited by $J$ if and only if $e \notin J$. 
\end{lemma}

\begin{proof}
Let $C$ be a clique dual to $J$, $x,y \in X$ the endpoints of $e$, and $a,b \in C$ their gates in $C$ respectively. Suppose that $x$ and $y$ does not belong to the same sector, so that $a$ and $b$ are two distinct adjacent vertices. Notice that $$1+d(x,a)=d(x,b) \leq d(x,y)+d(y,b)=1+d(y,b),$$ hence $d(x,a) \leq d(y,b)$. Similarly, we show that $d(y,b) \leq d(x,a)$, hence $d(x,a)=d(y,b)$; let $\ell$ denote this common value. Let $x_0,x_1, \ldots, x_m$ be a geodesic between $x_0=x$ and $x_m=a$, and set $y_0=y$. Notice that $x_1 \neq y$ because $x_1$ and $y$ have different gates in $C$, $d(b,y)= \ell$, $d(b,x)=\ell+1$ and $$d(b,x_1)=d(b,a)+d(a,x_1)=1+ \ell-1= \ell,$$ so the quadrangle condition implies that there exists a vertex $y_1 \in X$ adjacent to both $x_1$ and $y$ such that $d(b,y_1)=d(b,x)-2=d(b,y)-1$. Morever, the gate of $y_1$ in $C$ is again $b$. Indeed, $$1+d(y_1,b)= d(y,b)=d(y,C) \leq 1+d(y_1,C),$$ so $b$ minimizes in $C$ the distance to $y_1$. Then, we can define $y_2$ similarly from $x_2$ and $y_1$, and so on. Thus, we have constructed a sequence $y_0,y_1, \ldots, y_m \in X$ such that the gate of $y_i$ in $C$ is $b$, and $x_i,y_i,y_{i+1},x_{i+1}$ define a square, and $d(b,y_i)=d(b,y)-i$. In particular, $(x_m,y_m)=(a,b)$ so we conclude that the edge $e=(x,y)$ is dual to $J$.

\medskip \noindent
Conversely, suppose that $e \in J$. Let $C'$ denote the unique clique containing the edge $(x,y)$. Clearly, $x$ and $y$ belong to different sectors delimited by $C'$, so, because $J$ is dual to $C'$, they have to belong to different sectors delimited by $J$. 
\end{proof}

\begin{cor}\label{cor:sector}
Let $X$ be a quasi-median graph and $J$ a hyperplane. Let $X \backslash J$ denote the graph obtained from $X$ by removing the interiors of the edges of $J$. The connected components of $X \backslash J$ are precisely the sectors delimited by $J$. 
\end{cor}

\begin{proof}
Let $C$ be a clique dual to $J$. According to Lemma \ref{lem:samesector}, we know that any edge of $X$ either belongs to a sector delimited by $C$ or is dual to $J$, so we only have to show that two sectors delimited by $C$ are separated by $J$ and that a sector is connected. 

\medskip \noindent
Let $u,v \in C$ be two distinct vertices. If $\gamma$ is a path between two vertices of $[C,u]$ and $[C,v]$, then there necessarily exists an edge $e \subset \gamma$ whose endpoints have different gates in $C$. It follows from Lemma \ref{lem:samesector} that $e$ is dual to $J$, hence $\gamma \nsubseteq X \backslash J$. Therefore, two sectors delimited by $C$ are separated by $J$.

\medskip \noindent
The connectedness of a sector follows from the observation that any vertex can be joined by a geodesic to its gate in $C$, and any vertex of this geodesic has necessarily the same gate in $C$. Such a path does not contain any edge of $J$ according to Lemma \ref{lem:samesector}.
\end{proof}

\begin{cor}\label{cor:sectorgated}
In a quasi-median graph, a sector is gated.
\end{cor}

\begin{proof}
Let $C$ be a clique and $v \in C$ a vertex, and let $J$ denote the hyperplane dual to $C$. According to Proposition \ref{prop:gated}, it is sufficient to prove that $[C,v]$ is locally convex and contains its triangles in order to deduce that it is gated. 

\medskip \noindent
Let $(a,b,c)$ be a triangle with $(a,b) \subset [C,v]$. Because $a$ and $b$ belong to the same sector delimited by $J$, it follows from Lemma \ref{lem:samesector} that $(a,b) \notin J$. Necessarily, $(a,c) \notin J$. Lemma \ref{lem:samesector} implies that $a$ and $c$ belong to the same sector delimited by $J$, hence $c \in [C,v]$, and finally $(a,b,c) \subset [C,v]$.

\medskip \noindent
Let $(a,b,c,d)$ be a square with $(a,b),(b,c) \subset [C,v]$. Because $a$ and $b$ belong to the same sector delimited by $J$, it follows from Lemma \ref{lem:samesector} that $(a,b) \notin J$. A fortiori, $(c,d) \notin J$. Therefore, Lemma \ref{lem:samesector} implies that $c$ and $d$ belong to the same sector delimited by $J$, hence $d \in [C,v]$, and finally $(a,b,c,d) \subset [C,v]$.
\end{proof}

\begin{definition}
The \emph{carrier}\index{Carriers of hyperplanes} of a hyperplane $J$ of $X$, denoted by $N(J)$, is the subgraph generated by the endpoints of the edges which belong to $J$. The \emph{boundary} of $J$, denoted by $\partial J$, is the graph obtained from $N(J)$ by removing the interiors of the edges which belong to $J$. 
\end{definition}

\begin{lemma}\label{lem:carriergated}
In a quasi-median graph, the carrier of a hyperplane is gated.
\end{lemma}

\noindent
We begin by proving two preliminary lemmas. Using the vocabulary of CAT(0) cube complexes, the first lemma below states that hyperplanes do not self-intersect nor self-osculate.

\begin{lemma}\label{lem:cliqueshyp}
Let $C_1,C_2 \subset X$ be two distinct cliques dual to the same hyperplane. Either $C_1=C_2$ or $C_1 \cap C_2= \emptyset$.
\end{lemma}

\begin{proof}
Suppose that there exists some vertex $u \in C_1 \cap C_2$. Let $v \in C_1 \backslash \{ u \}$ be a vertex; or course, $u$ and $v$ are adjacent. If $d(v,C_2)=1$, then $u$ must be the gate of $v$ in $C_2$, so that $u$ and $v$ belong to the same sector delimited by $C_2$, contradicting Lemma \ref{lem:samesector} because the edge $(u,v)$ is dual to the same hyperplane as the clique $C_2$. Therefore, $d(v,C_2)=0$. This implies that the edge $(u,v)$ is included into $C_1 \cap C_2$, hence $C_1=C_2$ as a consequence of Lemma \ref{lem:cliqueinter}.
\end{proof}

\begin{lemma}\label{lem:parallelcliqueshyp}
Let $C_1,C_2 \subset X$ be two distinct cliques dual to the same hyperplane. If there exists an edge between $C_1$ and $C_2$, then $C_1$ and $C_2$ are parallel.
\end{lemma}

\begin{proof}
Let $u \in C_1$ and $v \in C_2$ be two adjacent vertices and let $w \in C_1 \backslash \{ u\}$. Notice that, as a consequence of Lemma \ref{lem:cliqueshyp}, $C_1$ and $C_2$ must be disjoint. First, we claim that $d(v,w)=2$. Indeed, because $C_1 \cap C_2 = \emptyset$, necessarily $d(v,w) \geq 1$; moreover, because a clique contains its triangles, $d(v,w) =1$ would imply $v \in C_1 \cap C_2$. Thus, we have $d(v,w)=2$, and in particular $d(w,C_2) \leq 2$. On the other hand, once again because $C_1$ and $C_2$ are disjoint, necessarily $d(w,C_2) \geq 1$; and $d(w,C_2)=2$ would imply that $v$ is the gate of $u$ and $w$ in $C_2$, so that $u$ and $w$ would belong to the same sector delimited by $C_2$, contradiction Lemma \ref{lem:samesector}. Therefore, $d(w,C_2)=1$, ie., there exists some vertex $x \in C_2$ adjacent to $w$. This proves that $C_1$ and $C_2$ are parallel according to Lemma \ref{lem:parallelcliques}.
\end{proof}

\begin{proof}[Proof of Lemma \ref{lem:carriergated}.]
According to Proposition \ref{prop:gated}, it is sufficient to prove that $N(J)$ is locally convex and contains its triangles in order to deduce that it is gated. 

\medskip \noindent
Let $(a,b,c)$ be a triangle with $(a,b) \subset N(J)$. If $(a,b) \in J$ then $(a,c) \in J$, which implies $c \in N(J)$. So we will suppose that $(a,b) \notin J$. Thus, there exist two edges $(a,e),(b,f) \in J$. Notice that, because $(a,b) \notin J$, necessarily $a$ and $f$, and similarly $b$ and $e$ are not adjacent. A fortiori, $(a,e)$ and $(b,f)$ belong to distinct cliques of $J$, and Lemma \ref{lem:parallelcliqueshyp} implies that $e$ and $f$ are adjacent. Once again because $(a,b) \notin J$, we know that $c$ and $e$, and similarly $c$ and $f$, are not adjacent. Therefore, $d(c,e)=d(c,f)=2$. The triangle condition implies that there exists a vertex $m \in X$ adjacent to $e,f,c$. Clearly, $(c,m) \in J$ hence $c \in N(J)$, and finally $(a,b,c) \subset N(J)$.

\medskip \noindent
Let $(a,b,c,d)$ be a square with $(a,b),(b,c) \subset N(J)$. If $(a,b) \in J$ then $(d,c) \in J$, which implies $d \in N(J)$; similarly, if $(b,c) \in J$ then $(a,d) \in J$, which implies $d \in N(J)$. So we suppose that $(a,b),(b,c) \notin J$. In particular, we know that there exist three edges $(a,e),(b,f),(c,g) \in J$. Notice that, because $(a,b) \notin J$, necessarily $a$ is not adjacent to $f$ and $b$ is not adjacent to $e$. Therefore, $(a,e)$ and $(b,f)$ belong to distinct cliques of $J$ so that Lemma \ref{lem:parallelcliqueshyp} implies that $e$ and $f$ are adjacent. Similarly, we show that $f$ and $g$ are adjacent. Now, notice that $(b,c) \notin J$ implies that $d$ and $e$ are not adjacent, hence $d(e,d)=2$; similarly, we show that $d(g,d)=2$. Let $J_1$ denote the hyperplane dual to $(a,b)$ and $J_2$ the hyperplane dual to $(b,c)$. Because $(a,b) \notin J$ and $(b,c) \notin J$, we know that $J_1 \neq J$ and $J_2 \neq J$. If $J_1=J_2$, then this hyperplane is dual to the two edges $(c,d)$ and $(c,b)$, and we deduce from Lemma \ref{lem:cliqueshyp} that these edges must belong to the same clique, so that $b$ and $d$ must be adjacent. On the other hand, we already know that $N(J)$ contains its triangles, so $(a,b) \subset N(J)$ implies $d \in N(J)$. So suppose that $J_1 \neq J_2$. It is worth noticing that $d$ and $f$ are separated by at least three (distinct) hyperplanes, namely $J,J_1,J_2$. Indeed, it follows from Lemma \ref{lem:samesector} that $d$ and $b$ belong to the same sector delimited by $J$ and that $b$ and $d$ belong to different sectors delimited by $J$, so that $J$ must separate $f$ and $d$; similarly, $J_1$ and $J_2$ separate $d$ and $f$. As a consequence, necessarily $d(f,d) \geq 3$, and in fact $d(f,d)=3$ since there exists a path of length three between $d$ and $f$ by construction. Now we are able to deduce from the quadrangle condition that there exists a vertex $m \in X$ adjacent to $d,e,g$. Clearly, the edge $(d,m)$ belongs to $J$, hence $d \in N(J)$. Consequently, $(a,b,c,d) \subset N(J)$. 
\end{proof}

\begin{definition}
The \emph{main fiber}\index{Main fiber of a hyperplane} of a hyperplane $J$, denoted by $F(J)$, is the graph whose vertices are the cliques of $J$ and whose edges link two parallel cliques.
\end{definition}

\noindent
Since any vertex of $N(J)$ belongs to a unique clique of $J$, as a consequence of Lemma \ref{lem:cliqueshyp}, there exists a canonical projection $c : N(J) \to F(J)$. Furthermore, if we fix some clique $C$ dual to $J$ and denote by $p : X \to C$ the application which associates to a vertex of $X$ its gate in $C$, we can define a natural map $$ \Psi : \left\{ \begin{array}{ccc} N(J) & \longrightarrow & F(J) \times C \\ v & \longmapsto & (c(v), p(v)) \end{array} \right.$$

\begin{lemma}\label{lem:carrierproduct}
$\Psi$ defines an isomorphism between the graphs $N(J)$ and $F(J) \times C$. 
\end{lemma}

\begin{proof}
First, we prove that $\Psi$ sends adjacent vertices to adjacent vertices. Let $e=(x,y) \subset N(J)$ be an edge. If $e \in J$ then $x$ and $y$ belong to the same clique of $J$, ie., $c(x)=c(y)$; moreover, Lemma \ref{lem:samesector} implies that $p(x) \neq p(y)$, so that $p(x)$ and $p(y)$ are necessarily adjacent. A fortiori, $\Psi(x)$ and $\Psi(y)$ are adjacent. If $e \notin J$, Lemma \ref{lem:samesector} implies that $p(x)=p(y)$. Moreover, $x$ and $y$ necessarily belong to distinct cliques of $J$, which are linked by an edge since $x$ and $y$ are adjacent; so it follows from Lemma \ref{lem:parallelcliqueshyp} that these cliques are parallel, ie., $c(x)$ and $c(y)$ are adjacent in $F(J)$. A fortiori, $\Psi(x)$ and $\Psi(y)$ are adjacent.

\medskip \noindent
We claim that $\Psi$ is injective. Let $u,v \in N(J)$ be two distinct vertices. If $c(u) \neq c(v)$, necessarily $\Psi(u) \neq \Psi(v)$, so we suppose that $c(u)=c(v)$, ie., $u$ and $v$ belong to the same clique of $J$. In particular, $u$ and $v$ are adjacent, and the edge $(u,v)$ belongs to $J$. It follows from Lemma \ref{lem:samesector} that $p(u) \neq p(v)$. A fortiori, $\Psi(u) \neq \Psi(v)$.

\medskip \noindent
We claim that $\Psi$ is surjective. Let $(C',x) \in F(J) \times C$. According to Lemma \ref{lem:sectordecomposition}, there exists a vertex $u \in C'$ such that $[C',u]=[C,x]$. In particular, we know that $u \in [C,x]$, hence $p(u)=x$. Finally, $\Psi(u)=(C',x)$.

\medskip \noindent
In order to conclude that $\Psi$ is an isomorphism, it remains to verify that, if $u,v \in N(J)$ are vertices, then $\Psi(u)$ adjacent to $\Psi(v)$ implies that $u$ is adjacent to $v$. So suppose that $\Psi(u)$ and $\Psi(v)$ are adjacent. Two cases may happen. First, suppose that $c(u)=c(v)$ and that $p(u)$ is adjacent to $p(v)$. Then $u$ and $v$ belong to the same clique of $J$, and $p(u) \neq p(v)$ implies $u \neq v$, so that $u$ is adjacent to $v$. Now, suppose that $p(u)=p(v)$ and that $c(u)$ is adjacent to $c(v)$. Thus, $u$ and $v$ belong to the same sector delimited by $J$, and to two parallel cliques of $J$, say $u \in C_1$ and $v \in C_2$. Let $x$ denote the vertex of $C_1$ opposite to $v$. According to Lemma \ref{lem:cliqueshyp}, the edge $(x,v)$ does not belong to $J$, so that $x$ and $v$ belong to the same sector delimited by $J$ according to Lemma \ref{lem:samesector}. A fortiori, $u$ and $x$ necessarily belong to the same sector delimited by $J$. Because $u$ and $x$ belong to the same clique of $J$, we deduce from Lemma \ref{lem:samesector} that $u=x$. As a consequence, $u$ and $v$ are adjacent.
\end{proof}

\begin{lemma}\label{lem:ccpartialJ}
Let $F$ be a connected component of $\partial J$. Then $F$ is a gated subgraph isomorphic to the main fiber $F(J)$. More precisely, there exists a vertex $v \in C$ such that $F= \Psi^{-1}(F(J) \times \{ v \})$.
\end{lemma}

\begin{proof}
Noticing that $\Psi(J)= \bigcup\limits_{K \in F(J)} \{ K \} \times C$, we deduce that $\Psi(\partial J)= \bigcup\limits_{v \in C} F(J) \times \{ v \}$. As a consequence, there exists a vertex $v \in C$ such that $F= \Psi^{-1}(F(J) \times \{ v \})$. In particular, $F$ is clearly isomorphic to $F(J)$. Moreover, it is worth noticing that the vertices of $F$ are precisely the vertices of $N(J)$ whose image by $p$ is $v$, ie., $F= N(J) \cap [C,v]$. Since the intersection of two gated subgraphs is gated, we deduce from Corollary \ref{cor:sectorgated} and Lemma \ref{lem:carriergated} that $F$ is a gated subgraph.
\end{proof}

\noindent
All the results above imply Proposition \ref{prop:hypsumup}, as already mentionned. Now, we want to characterize geodesics from the hyperplanes they intersect. The analogue result for CAT(0) cube complexes is fundamental.

\begin{prop}\label{prop:geodesichyp}
Let $X$ be a quasi-median graph. A path in $X$ is a geodesic if and only if it intersects any hyperplane at most once. In particular, the distance between two vertices of $X$ is equal to the number of hyperplanes separating them.
\end{prop}

\begin{proof}
Let $p$ be a path which intersects a hyperplane $J$ twice. Let $u$ be the first vertex of $p$ which belongs to $N(J)$, and similarly $v$ the last vertex of $p$ which belongs to $N(J)$. Let $q$ denote the subpath of $p$ between $u$ and $v$. Since $N(J)$ is convex, according to Lemma \ref{lem:carriergated}, if $q \nsubseteq N(J)$ then $p$ cannot be a geodesic. So suppose that $q \subset N(J)$. Let $\Psi : N(J) \to F(J) \times C$ denote the isomorphism given by Lemma \ref{lem:carrierproduct}. Clearly, $\Psi(q)$ is not a geodesic in $F(J) \times C$ since $q$ contains two edges dual to $J$, so $q= \Psi(\Psi^{-1}(q))$ cannot be geodesic in $X$. A fortiori, $p$ cannot be a geodesic as well. Thus, we have proved that a geodesic must intersect any hyperplane at most once.

\medskip \noindent
Let $x,y \in X$ be two vertices and $\gamma$ a geodesic between them. If $J$ is a hyperplane which does not separate $x$ and $y$, it follows from the convexity of the sectors given by Corollary \ref{cor:sectorgated} that $\gamma$ must be disjoint from $J$. Therefore, any hyperplane intersecting $\gamma$ necessarily separates $x$ and $y$. On the other hand, we know that $\gamma$ intersects at most once each of these hyperplanes, so its length is bounded above by the number of hyperplanes separating $x$ and $y$. Conversely, $\gamma$ necessarily intersects any hyperplane separating $x$ and $y$, so its length is bounded below by the number of hyperplanes separating $x$ and $y$. Therefore, the distance between $x$ and $y$ is precisely the number of hyperplanes separating them.

\medskip \noindent
Now, let $p$ be a path between two vertices $x,y \in X$ such that $p$ intersects any hyperplane at most once. If $J$ is a hyperplane which does not separate $x$ and $y$, necessarily $p$ must be disjoint from $J$ since otherwise it would have to intersect it at least twice because $x$ and $y$ lie in different sectors delimited by $J$. Therefore, $p$ intersects only the hyperplanes separating $x$ and $y$. Because $p$ intersects any hyperplane at most once, we deduce that its length is bounded above by the number of hyperplanes separating $x$ and $y$, which correponds exactly to the distance between $x$ and $y$. Therefore, $p$ has to be a geodesic.
\end{proof}

\begin{remark}
Most of the results contained in this section were proved for CAT(0) cube complexes by Sageev in his thesis \cite{MR1347406}. See also \cite{MR2413337}.
\end{remark}

\subsection{Projections onto gated subgraphs}

\noindent
In this section, we generalize the results proved in \cite[Section 2]{article3} on combinatorial projections in CAT(0) cube complexes. 

\begin{definition}
Let $X$ be a quasi-median graph and $Y \subset X$ a gated subgraph. The application $\mathrm{proj}_Y : X \to Y$ which sends a vertex of $X$ to its gate in $Y$ will be referred to as the \emph{projection onto $Y$}\index{Projections onto gated subgraphs}. 
\end{definition}

\noindent
Our first main result describes the hyperplanes separating the projections of two vertices.

\begin{prop}\label{prop:2projseparate}
Let $Y \subset X$ be a gated subgraph and let $p : X \to Y$ denote the projection onto $Y$. For every vertices $x,y \in X$, the hyperplanes separating $p(x)$ and $p(y)$ are precisely the hyperplanes which separate $x$ and $y$ and intersect $Y$.
\end{prop}

\noindent
Our proposition will be a consequence of the following lemma:

\begin{lemma}\label{lem:projseparate}
Let $Y \subset X$ be a gated subgraph and $x \in X$ a vertex. Any hyperplane separating $x$ and its projection onto $Y$ separates $x$ and $Y$.
\end{lemma}

\begin{proof}
Suppose by contradiction that there exists a hyperplane separating $x$ and its projection $p$ onto $Y$ which intersects $Y$. In particular, $J$ must separate $p$ from some vertex $y \in Y$. As a consequence, if we fix a geodesic $\gamma$ between $x$ and $y$ passing through $p$, necessarily $\gamma$ has to intersect $\gamma$ at least twice. This contradicts Proposition \ref{prop:geodesichyp}.
\end{proof}

\begin{proof}[Proof of Proposition \ref{prop:2projseparate}.]
Let $J$ be a hyperplane separating $p(x)$ and $p(y)$. In particular, $J$ intersects $Y$. Because the previous lemma implies that any hyperplane separating $x$ and $p(x)$ or $y$ and $p(y)$ must be disjoint from $Y$, we deduce that $J$ separates neither $x$ and $p(x)$ nor $y$ and $p(y)$. Thus, $J$ must separate $x$ and $y$. Conversely, suppose that $J$ is a hyperplane which separates $x$ and $y$ and intersects $Y$. The same argument shows that $J$ cannot separate $x$ and $p(x)$ or $y$ and $p(y)$, so that $J$ must separate $p(x)$ and $p(y)$.
\end{proof}

\begin{cor}\label{cor:projectionlip}
Let $Y \subset X$ be a gated subgraph and let $p : X \to Y$ denote the projection onto $Y$. Then $p$ is $1$-Lipschitz, ie., $d(p(x),p(y)) \leq d(x,y)$ for every vertices $x,y \in X$.
\end{cor}

\begin{proof}
It follows from Proposition \ref{prop:geodesichyp} and Proposition \ref{prop:2projseparate} that
$$\begin{array}{lcl} d(p(x),p(y)) & = & \# \{ \text{hyperplanes separating $p(x)$ and $p(y)$} \} \\ & \leq & \# \{ \text{hyperplanes separating $x$ and $y$} \} = d(x,y), \end{array}$$
which concludes the proof.
\end{proof}

\noindent
Although Lemma \ref{lem:projseparate} has been introduced to prove Proposition \ref{prop:2projseparate}, it will turn out to be quite useful on its own right. For instance, it leads to the following result:

\begin{lemma}\label{lem:2min}
Let $Y_1,Y_2 \subset X$ be two gated subgraphs. If $x \in Y_1$ and $y \in Y_2$ are two vertices minimizing the distance between $Y_1$ and $Y_2$, then the hyperplanes separating $x$ and $y$ are precisely the hyperplanes separating $Y_1$ and $Y_2$.
\end{lemma}

\begin{proof}
Let $p : X \to Y_1$ and $q : X \to Y_2$ denote respectively the projections onto $Y_1$ and $Y_2$. In particular, we have $y=p(x)$ and $x=q(y)$. By applying Lemma \ref{lem:projseparate} twice, we deduce that any hyperplane separating $x$ and $y$ must be disjoint from $Y_1$ and $Y_2$; a fortiori, it has to separate $Y_1$ and $Y_2$. Conversely, it is clear that any hyperplane separating $Y_1$ and $Y_2$ must separate $x$ and $y$.
\end{proof}

\noindent
Our second main result is the following:

\begin{prop}
Let $Y_1,Y_2 \subset X$ be two gated subgraphs, and let $p : X \to Y_2$ denote the projection onto $Y_2$. Then $p(Y_1)$ is a geodesic subgraph of $Y_2$, and the hyperplanes intersecting $p(Y_1)$ are precisely the hyperplanes intersecting both $Y_1$ and $Y_2$. 
\end{prop}

\begin{proof}
According to Corollary \ref{cor:projectionlip}, $p$ is 1-Lipschitz, so $p$ sends an edge to either an edge or a vertex; in particular, $p$ sends a path to a path. Fix two vertices $x,y \in X$ and a geodesic $\gamma$ between them; let $\bar{\gamma} \subset p(Y_1)$ be the path which is the image of $\gamma$ by $p$. According to Proposition \ref{prop:geodesichyp}, it is sufficient to prove that no hyperplane intersects twice $\bar{\gamma}$ in order to justify that $\bar{\gamma}$ is a geodesic. If $e_1$ and $e_2$ are two edges of $\bar{\gamma}$, choose two lifts $f_1,f_2 \subset \gamma$. It follows from Proposition \ref{prop:2projseparate} that the hyperplane dual to $e_1$ (resp. $e_2$) is also dual to $f_1$ (resp. $f_2$). Therefore, a hyperplane dual to both $e_1$ and $e_2$ would be dual to both $f_1$ and $f_2$, and so would intersect $\gamma$ twice, which contradicts the fact that $\gamma$ is a geodesic. Thus, we have proved that $p(Y_1)$ is a geodesic subgraph.

\medskip \noindent
If a hyperplane intersects $p(Y_1) \subset Y_2$, it cleary intersects $Y_2$ and it must intersect $Y_1$ according to Proposition \ref{prop:2projseparate}. Conversely, if a hyperplane $J$ intersects both $Y_1$ and $Y_2$, it is dual to an edge $e \subset Y_1$, and it follows from Proposition \ref{prop:2projseparate} that $J$ separates the projections onto $Y_2$ of the endpoints of $e$; a fortiori, $J$ intersects $p(Y_1)$.
\end{proof}

\noindent
We conclude this section by proving the following result which will be useful later.

\begin{lemma}\label{lem:projinter1}
Let $X$ be a quasi-median graph and $Y_1,Y_2 \subset X$ two intersecting gated subgraphs. Then $\mathrm{proj}_{Y_2} \circ \mathrm{proj}_{Y_1} = \mathrm{proj}_{Y_1 \cap Y_2}$.
\end{lemma}

\noindent
We begin by proving the following observation.

\begin{lemma}\label{lem:projinter2}
Let $X$ be a quasi-median graph, $Y_1,Y_2 \subset X$ two intersecting gated subgraphs, and $y \in Y_1$ a vertex. If $z$ denotes the projection of $y$ onto $Y_2$, then $z \in Y_1 \cap Y_2$. 
\end{lemma}

\begin{proof}
If we fix some vertex $w \in Y_1 \cap Y_2$, then $z$ must belong to a geodesic between $y$ and $w$, hence $z \in I(y,w)$. On the other hand, because $Y_1$ is convex, necessarily $I(y,w) \subset Y_1$, hence $z \in Y_1$. We already know that $z \in Y_2$, so $z \in Y_1 \cap Y_2$.
\end{proof}

\begin{proof}[Proof of Lemma \ref{lem:projinter1}.]
Let $x \in X$ be a vertex. Let $y$ denote its projection onto $Y_1$ and $z$ the projection of $y$ onto $Y_2$. We want to prove that $z$ is the projection of $x$ onto $Y_1 \cap Y_2$. Notice that $z$ belongs to $Y_1 \cap Y_2$ according to the previous lemma.

\medskip \noindent
Now, let $w \in Y_1 \cap Y_2$ be any vertex. Noticing that
$$\begin{array}{lcl} d(x,w) & = & d(x,y)+d(y,w) \\ & = & d(x,y)+d(y,z)+d(z,w) \\ & \geq & d(x,y)+d(y,z)=d(x,z), \end{array}$$
we deduce that $z$ is the vertex of $Y_1 \cap Y_2$ minimizing the distance to $x$, ie., $z$ is the projection of $x$ onto $Y_1 \cap Y_2$.
\end{proof}

\noindent
As an immediate corollary, we deduce the following statement.

\begin{cor}\label{cor:projnested}
Let $X$ be a quasi-median graph and $Y_1 \subset Y_2$ two gated subgraphs. Then $\mathrm{proj}_{Y_1} \circ \mathrm{proj}_{Y_2} = \mathrm{proj}_{Y_1}$. 
\end{cor}

\subsection{Quasi-cubulating popsets}\label{section:spaceswithpartitions}

\noindent
Based on works of Dunwoody and Sageev \cite{MR1347406}, Roller describes in \cite{Roller} how to construct CAT(0) cube complexes from some special class of posets, called \emph{pocsets}; see also \cite{SageevCAT(0)}. We generalize the construction to obtain quasi-median graphs. It is worth noticing that, in this paper, \emph{quasi-cubulating popsets} will be used only to prove Proposition \ref{prop:quasicubulatingqm}, from which we will prove Corollary \ref{cor:finitesubgraph}, Lemma \ref{lem:prismhyp} and Lemma \ref{lem:productingflatrectangle} below; and later to show how to cubulate a \emph{space with walls} containing duplicates (which can also be done using pocsets). Therefore, the technicalities below can be skipped in a first lecture. 

\noindent
Before giving the definition of a popset, let us mention an important example to keep in mind. Let $X$ be a quasi-median graph, and let $\mathcal{S}(X)$ denote the set of the sectors of $X$. For every hyperplane $J$ of $X$, let $\mathcal{P}(J)$ denote the set of sectors delimited by $J$. Notice that $\mathfrak{P} = \{ \mathcal{P}(J) \mid J \ \text{hyperplane} \}$ defines a partition of $\mathcal{S}(X)$. The data $(\mathcal{S}(X), \subset, \mathfrak{P})$ is an example of what we call a \emph{popset}. 

\begin{definition}\label{def:popset}
A \emph{popset}\index{Popsets} $(X, < , \mathfrak{P})$ is the data of a poset $(X,<)$ with a partition $\mathfrak{P}$ of $X$ such that:
\begin{itemize}
	\item for every $\mathcal{P} \in \mathfrak{P}$, $\# \mathcal{P} \geq 2$;
	\item for every $\mathcal{P} \in \mathfrak{P}$, no two elements of $\mathcal{P}$ are $<$-comparable;
	\item for every $\mathcal{P}_1, \mathcal{P}_2 \in \mathfrak{P}$, if there exist $A_1 \in \mathcal{P}_1$ and $A_2 \in \mathcal{P}_2$ such that $A_1 < A_2$, then there exists some $B_1 \in \mathcal{P}_1$ such that $A < A_2$ for every $A \in \mathcal{P}_1 \backslash \{ B_1 \}$ and $B < B_1$ for every $B \in \mathcal{P}_2 \backslash \{ A_2 \}$. If so, we say that $\mathcal{P}_1$ and $\mathcal{P}_2$ are \emph{nested}. (Notice that being nested is a symmetric relation.)
\end{itemize}
An element of $\mathfrak{P}$ is referred to as a \emph{wall}, and an element of a wall as a \emph{sector}; notice that a sector is an element of $X$. Two walls which are not nested are \emph{transverse}. If $\mathcal{P}_1, \mathcal{P}_2 \in \mathfrak{P}$ are two nested walls and $A_1 \in \mathcal{P}_1$ is a sector, we say that $A_1$ is a \emph{sector delimited by $\mathcal{P}_1$ which contains $\mathcal{P}_2$} if there exists some $A_2 \in \mathcal{P}_2$ such that $D<A_1$ for every $D \in \mathcal{P}_2 \backslash \{ A_2 \}$. 
\end{definition}

\noindent
First, let us notice that there exists a unique sector of $\mathcal{P}_1$ which contains $\mathcal{P}_2$ when $\mathcal{P}_1$ and $\mathcal{P}_2$ are nested.

\begin{lemma}\label{lem:noncomparable}
Let $\mathcal{P}_1, \mathcal{P}_2 \in \mathfrak{P}$ be two nested walls and $A_1$ (resp. $A_2$) a sector delimited by $\mathcal{P}_1$ (resp. $\mathcal{P}_2$) which contains $\mathcal{P}_2$ (resp. $\mathcal{P}_1$). Then $A_1$ and $A_2$ are not $<$-comparable.
\end{lemma}

\begin{proof}
Suppose by contradiction that $A_1<A_2$. Then $D<A_2$ for every $D \in \mathcal{P}_1$, and there exists a sector $B_1 \in \mathcal{P}_1$ such that $A< B_1$ for every $A \in \mathcal{P}_2 \backslash \{ A_2 \}$. Therefore, if $D \in \mathcal{P}_2 \backslash \{ A_2 \}$, then $D<B_1<A_2$, which is impossible since two sectors of $\mathcal{P}_2$ are not $<$-comparable.
\end{proof}

\begin{cor}
Let $\mathcal{P}_1, \mathcal{P}_2 \in \mathcal{P}$ be two nested walls. There exists a unique sector delimited by $\mathcal{P}_1$ which contains $\mathcal{P}_2$.
\end{cor}

\begin{proof}
Let $A_1$ (resp. $A_2$) be a sector delimited by $\mathcal{P}_1$ (resp. $\mathcal{P}_2$) which contains $\mathcal{P}_2$ (resp. $\mathcal{P}_1$). If $A \in \mathcal{P}_1 \backslash \{ A_1 \}$, then $A<A_2$. We deduce from the previous lemma that $A$ cannot be a sector delimited by $\mathcal{P}_1$ which contains $\mathcal{P}_2$.
\end{proof}

\noindent
An \emph{orientation}\index{Orientations} $\sigma$, as defined below, may be thought of as ``choice function'': given a wall $\mathcal{P} \in \mathfrak{P}$, $\sigma$ chooses the sector $\sigma(\mathcal{P})$ delimited by $\mathcal{P}$. Our definition mimic the definition of \emph{ultrafilters} introduced in \cite{Roller}.

\begin{definition}
Let $(X, < , \mathfrak{P})$ be a popset. An \emph{orientation} is a map $\sigma : \mathfrak{P} \to X$ such that
\begin{itemize}
	\item $\sigma(\mathcal{P}) \in \mathcal{P}$ for every $\mathcal{P}\in \mathfrak{P}$;
	\item if $A_1 \in \mathcal{P}_1$ and $A_2 \in \mathcal{P}_2$ satisfy $A_1 < A_2$, then $A_1= \sigma(\mathcal{P}_1)$ implies $A_2= \sigma(\mathcal{P}_2)$.
\end{itemize}
\end{definition}

\begin{definition}
Given a popset $(X,<, \mathfrak{P})$, we define its \emph{quasi-cubulation} $C(X,<, \mathfrak{P})$ as the graph whose vertices are the orientations of $(X,<, \mathfrak{P})$ and whose edges link two orientations which differ only on a single wall. 
\end{definition}

\noindent
Our purpose is first to show that this construction yields quasi-median graphs, and next to understand the combinatorics of its hyperplanes (see Theorem \ref{thm:quasicubulation}). The first step, which is Proposition \ref{prop:qcqm} below, needs several preliminary lemmas.

\begin{prop}\label{prop:qcqm}
Let $(X,<, \mathfrak{P})$ be a popset. Any connected component of its quasi-cubulation $C(X,<, \mathfrak{P})$ is a quasi-median graph.
\end{prop}

\noindent
For convenience, we introduce some notation. If $\sigma$ is an orientation and $D$ a sector of some wall $\mathcal{P} \in \mathfrak{P}$, we define the map $$[\sigma,D] : \left\{ \begin{array}{ccc} \mathfrak{P} & \longrightarrow & X \\ \mathcal{Q} & \longmapsto & \left\{ \begin{array}{cl} \sigma(\mathcal{Q}) & \text{if} \ \mathcal{Q} \neq \mathcal{P} \\ D & \text{otherwise} \end{array} \right. \end{array} \right..$$ In particular, two orientations $\sigma_1$ and $\sigma_2$ are adjacent in $C(X, \mathfrak{P})$ if and only if $\sigma_2=[\sigma_1,D]$ for some sector $D \neq \sigma_1(\mathcal{P})$ (if so, notice that we also have $\sigma_1=[\sigma_2,\sigma_1(\mathcal{P})]$ where $\mathcal{P}$ denotes the wall delimiting $D$). For convenience, we will write $[\sigma,D_1,D_2]$ instead of $[[\sigma,D_1],D_2]$, $[\sigma,D_1,D_2,D_3]$ instead of $[[\sigma,D_1,D_2],D_3]$, and so on.

\medskip \noindent
Thought of as choice functions, a crucial property satisfied by any orientation is that, if it is possible to modify the choice on some wall to get an orientation, then choosing any sector will produce an orientation. This statement is made precise by Corollary \ref{cor:nochoice}.

\begin{lemma}\label{lem:orientationmin}
Let $\sigma$ be an orientation and $A \neq \sigma(\mathcal{P})$ a sector delimited by some wall $\mathcal{P}$. Then $[\sigma,A]$ defines an orientation if and only if $\sigma(\mathcal{P})$ is minimal in $\sigma(\mathfrak{P})$. 
\end{lemma}

\begin{proof}
Suppose that $\sigma(\mathcal{P})$ is minimal in $\sigma(\mathfrak{P})$. To prove that $[\sigma,A]$ defines an orientation, it is sufficient to show that, if $A_1,A_2$ are two sectors, respectively delimited by the walls $\mathcal{P}_1, \mathcal{P}_2$, satisfying $A_1<A_2$ and $A_1= [\sigma,A](\mathcal{P}_1)$, then $A_2= [\sigma,A](\mathcal{P}_2)$. Notice that $\mathcal{P}_1$ and $\mathcal{P}_2$ must be nested, and that $A_2$ is the sector delimited by $\mathcal{P}_2$ which contains $\mathcal{P}_1$; let $B_1$ denote the sector delimited by $\mathcal{P}_1$ which contains $\mathcal{P}_2$. Three cases may happen.

\medskip \noindent
Suppose that $\mathcal{P}=\mathcal{P}_1$. If $\sigma(\mathcal{P}_1) \neq B_1$, then $\sigma(\mathcal{P}_1)<A_2$ which implies that $A_2= \sigma(\mathcal{P}_2)$ because $\sigma$ is an orientation. If $\sigma(\mathcal{P}_1)=B_1$, notice that, for every $D \in \mathcal{P}_2 \backslash \{ A_2 \}$, we have $D< B_1 = \sigma(\mathcal{P}_1)= \sigma(\mathcal{P})$, so that, since $\sigma(\mathcal{P})$ is minimal in $\sigma(\mathfrak{P})$ by assumption, we deduce that $D \neq \sigma(\mathcal{P}_2)$. A fortiori, $\sigma(\mathcal{P}_2)= A_2$. We conclude that $A_2= \sigma(\mathcal{P}_2)= [\sigma,A](\mathcal{P}_2)$ since $\mathcal{P}_2 \neq \mathcal{P}_1=\mathcal{P}$.

\medskip \noindent
Suppose that $\mathcal{P}= \mathcal{P}_2$. We have $A_1= [\sigma,A](\mathcal{P}_1)= \sigma(\mathcal{P}_1)$ since $\mathcal{P}_1 \neq \mathcal{P}_2 = \mathcal{P}$. Therefore, because $\sigma$ is an orientation, $A_1< A_2$ implies $A_2= \sigma(\mathcal{P}_2)$. We deduce that $\sigma(\mathcal{P}_1)=A_1 < A_2= \sigma(\mathcal{P}_2) = \sigma(\mathcal{P})$, which is impossible since $\sigma(\mathcal{P})$ is minimal in $\sigma(\mathfrak{P})$ by assumption.

\medskip \noindent
Suppose that $\mathcal{P} \neq \mathcal{P}_1, \mathcal{P}_2$. Then $A_1=[\sigma,A](\mathcal{P}_1)= \sigma(\mathcal{P}_1)$, so, because $\sigma$ is an orientation, $A_1 < A_2$ implies $A_2= \sigma(\mathcal{P}_2)= [\sigma,A](\mathcal{P}_2)$.

\medskip \noindent
Thus, we have proved that $[\sigma,A]$ defines an orientation. Conversely, if $\sigma(\mathcal{P})$ is not minimal in $\sigma(\mathfrak{P})$, then there exists some $\mathcal{Q} \in \mathfrak{P} \backslash \{ \mathcal{P} \}$ satisfying $\sigma(\mathcal{Q})< \sigma(\mathcal{P})$. Noticing that $[\sigma,A](\mathcal{Q})= \sigma(\mathcal{Q})< \sigma(\mathcal{P})$ but $[\sigma,A](\mathcal{P})=A \neq \sigma(\mathcal{P})$, we deduce that $[\sigma,A]$ does not define an orientation.
\end{proof}

\noindent
As an immediate consequence, we deduce the following statement.

\begin{cor}\label{cor:nochoice}
Let $\sigma$ be an orientation and $D \neq \sigma(\mathcal{P})$ a sector delimited by some partition $\mathcal{P}$. If $[\sigma,D]$ defines an orientation, then so does $[\sigma,D']$ for every $D' \in \mathcal{P}$.
\end{cor}

\noindent
The first step in understanding the geometry of $C(X,<, \mathfrak{P})$ is to understand its geodesics. This is the purpose of Lemma \ref{lem:orientationgeod} below, but before we need a technical preliminary result. 

\begin{lemma}\label{lem:orientationcommute}
Let $\sigma$ be an orientation and $D_1, \ldots, D_n$ a collection of sectors delimited by $\mathcal{P}_1, \ldots, \mathcal{P}_n$ respectively, such that $[\sigma,D_1, \ldots, D_k]$ defines an orientation for every $1 \leq k \leq n$. If $\mathcal{P}_n$ is transverse to $\mathcal{P}_1, \ldots, \mathcal{P}_{n-1}$, then $[\sigma,D_1, \ldots, D_n]=[\sigma,D_n, D_1, \ldots, D_{n-1}]$ and $[\sigma,D_n, D_1, \ldots, D_k]$ defines an orientation for every $1 \leq k \leq n-1$.
\end{lemma}

\begin{proof}
Notice that, for every $1 \leq k \leq n-1$ and $\mathcal{P} \in \mathfrak{P}$, we have
$$[\sigma,D_n,D_1, \ldots, D_k](\mathcal{P})= \left\{ \begin{array}{cl} [\sigma,D_1, \ldots, D_k](\mathcal{P}) & \text{if} \ \mathcal{P} \neq \mathcal{P}_n \\ D_n & \text{if} \ \mathcal{P}= \mathcal{P}_n \end{array} \right. = [\sigma,D_1, \ldots, D_k,D_n](\mathcal{P}),$$
hence $[\sigma,D_n,D_1, \ldots, D_k]= [\sigma,D_1, \ldots, D_k,D_n]$. In particular, if we set $\mu_k=[\sigma,D_1, \ldots, D_k]$, it is sufficient to prove that $[\mu_k,D_n]$ defines an orientation to deduce that $[\sigma, D_n,D_1, \ldots, D_k]$ defines an orientation as well, since $\mu_k$ is an orientation by assumption and $[\mu_k,D_n]= [\sigma,D_n,D_1, \ldots, D_k]$ by our previous observation. Precisely, we want to prove that, if $A_1 \in \mathcal{Q}_1$ and $A_2 \in \mathcal{Q}_2$ satisfy $A_1< A_2$ and $A_1= [\mu_k,D_n](\mathcal{Q}_1)$, then $A_2=[\mu_k,D_n](\mathcal{Q}_2)$. Notice that our assumption implies that $\mathcal{Q}_1$ and $\mathcal{Q}_2$ are nested. Two cases may happen.

\medskip \noindent
Suppose that $\mathcal{Q}_1 \neq \mathcal{P}_n$. As a consequence, $A_1=[\mu_k,D_n](\mathcal{Q}_1)=\mu_k(\mathcal{Q}_1)$, so that $A_1<A_2$ implies $A_2 = [\sigma,D_1, \ldots, D_k](\mathcal{Q}_2)$. If $\mathcal{Q}_2 \neq \mathcal{P}_n$, then $A_2=[\sigma, D_1, \ldots, D_k, D_n](\mathcal{Q}_2)$, and we are done. So suppose that $\mathcal{Q}_2= \mathcal{P}_n$. In particular, since $\mathcal{Q}_2= \mathcal{P}_n$ is not transverse to $\mathcal{Q}_1$, we deduce that $\mathcal{Q}_1 \neq \mathcal{P}_i$ for every $1 \leq i \leq n$, hence $A_1= [\sigma,D_1, \ldots, D_n](\mathcal{Q}_1)$. Since $[\sigma,D_1, \ldots, D_n]$ is an orientation, $A_1 < A_2$ implies that $A_2=[\sigma,D_1, \ldots, D_n](\mathcal{Q}_2)=D_n = [\sigma,D_1, \ldots, D_k,D_n](\mathcal{Q}_2)$.

\medskip \noindent
Suppose that $\mathcal{Q}_1 = \mathcal{P}_n$. Then $A_1= [\sigma, D_1, \ldots, D_k,D_n](\mathcal{Q}_1)=D_n$. Noticing that $\mathcal{Q}_1 = \mathcal{P}_n$ is not transverse to $\mathcal{Q}_2$, necessarily $\mathcal{Q}_2 \neq \mathcal{P}_i$ for every $1 \leq i \leq n$. On the other hand, $[\sigma,D_1, \ldots, D_n](\mathcal{Q}_1)=D_n=A_1< A_2$ implies $A_2= [\sigma,D_1,\ldots, D_n](\mathcal{Q}_2)$. Therefore, $A_2= [\sigma,D_1, \ldots, D_n](\mathcal{Q}_2) = \sigma(\mathcal{Q}_2)=[\sigma,D_1, \ldots, D_k,D_n](\mathcal{Q}_2)$. This concludes the proof.
\end{proof}

\begin{lemma}\label{lem:orientationgeod}
Let $\sigma$ be an orientation and $D_1, \ldots, D_n$ a collection of sectors, delimited by the walls $\mathcal{P}_1, \ldots, \mathcal{P}_n$ respectively, such that $$\sigma, \ [\sigma,D_1], \ [\sigma,D_1,D_2], \ldots, [\sigma,D_1, \ldots, D_n]$$ defines a path in $C(X, \mathfrak{P})$. This path is a geodesic if and only if, for every $i \neq j$, $\mathcal{P}_i \neq \mathcal{P}_j$ and $D_i \neq \sigma(\mathcal{P}_i)$.
\end{lemma}

\begin{proof}
Suppose that there exist some $1 \leq i<k \leq n$ such that $\mathcal{P}_i = \mathcal{P}_k$. We choose $i$ and $k$ such that $\mathcal{P}_r \neq \mathcal{P}_s$ for every $1 \leq r < s \leq k-1$. 

\medskip \noindent
Suppose by contradiction that there exists some $i< j< k$ such that $\mathcal{P}_j$ and $\mathcal{P}_i$ are nested; let $A$ be the sector of $\mathcal{P}_i=\mathcal{P}_k$ which contains $\mathcal{P}_j$, and $B$ the sector of $\mathcal{P}_j$ which contains $\mathcal{P}_i=\mathcal{P}_k$. If $D_i = D_k$, then the fact that $\mathcal{P}_r \neq \mathcal{P}_s$ for every $1 \leq r < s \leq k-1$ implies that $[\sigma,D_1, \ldots, D_k](\mathcal{P})=[\sigma,D_1, \ldots, D_{k-1}](\mathcal{P})$ for every wall $\mathcal{P} \in \mathfrak{P}$, hence $[\sigma,D_1, \ldots, D_k]=[\sigma,D_1, \ldots, D_{k-1}]$. As a consequence, we can shorten our path as
$$\sigma, \ [\sigma,D_1], \ldots, [\sigma,D_1, \ldots, D_{k-1}], \ [\sigma,D_1, \ldots, D_{k-1},D_{k+1}], \ldots, [\sigma,D_1, \ldots, D_{k-1},D_{k+1}, \ldots, D_n].$$
Similarly, if $D_j= \sigma(\mathcal{P}_j)$, then $[\sigma,D_1, \ldots, D_j]=[\sigma,D_1, \ldots, D_{j-1}]$ so that we are able to shorten out path. The same argument holds if $D_i = \sigma(\mathcal{P}_i)$. From now on, we suppose that $D_i \neq D_k$, $D_i \neq \sigma(\mathcal{P}_i)$ and $D_j \neq \sigma(\mathcal{P}_j)$. Now, we distinguish two cases.

\medskip \noindent
Suppose that $D_i \neq A$. Then $D_i< B$ with $D_i=[\sigma,D_1, \ldots, D_i](\mathcal{P}_i)$, hence $B=[\sigma,D_1, \ldots, D_i](\mathcal{P}_j)=\sigma(\mathcal{P}_j)$. On the other hand, $D_j \neq \sigma(\mathcal{P}_j) =B$ implies $D_j< A$, so that we deduce from $D_j=[\sigma,D_1, \ldots, D_j](\mathcal{P}_j)$ that $A=[\sigma,D_1, \ldots, D_j](\mathcal{P}_i)=D_i$, a contradiction.

\medskip \noindent
Suppose that $D_i=A$. In particular, $D_k \neq A$ so $D_k < B$. On the other hand, we know that $D_k= [\sigma, D_1, \ldots, D_k](\mathcal{P}_k)$ hence $B=[\sigma,D_1, \ldots, D_k](\mathcal{P}_j)=D_j$. But $\sigma(\mathcal{P}_j) \neq D_j=B$ implies $\sigma(\mathcal{P}_j)< A$, so that we deduce that $D_i=A= \sigma(\mathcal{P}_i)$, a contradiction.

\medskip \noindent
Thus, we have proved that $\mathcal{P}_k$ is transverse to $\mathcal{P}_j$ for every $i< j< k$. As a consequence of Lemma \ref{lem:orientationcommute}, we deduce that
$$[\sigma,D_1, \ldots, D_k]=[\sigma,D_1, \ldots, D_i,D_k,D_{i+1}, \ldots, D_{k-1}]= [\sigma,D_1, \ldots, D_{i-1},D_k,D_{i+1}, \ldots, D_{k-1}]$$
and that $[\sigma,D_1, \ldots, D_{i-1},D_k, D_{i+1}, \ldots, D_{i+s}]$ defines an orientation for every $0 \leq s \leq k-i$. Thus, it is possible to shorten our path by replacing the subsegment
$$\sigma, \ [\sigma,D_1], \ [\sigma,D_1, D_2], \ldots, [\sigma,D_1, \ldots, D_k]$$
with
$$\sigma, [\sigma,D_1], \ldots, [\sigma,D_1, \ldots, D_{i-1}], [\sigma,D_1, \ldots, D_{i-1},D_k,], [\sigma,D_1, \ldots, D_{i-1},D_k,D_{i+1}], \ldots,$$ 
$$[\sigma,D_1, \ldots, D_{i-1},D_k,D_{i+1}, \ldots, D_{k-1}]$$
Thus, we have prove that, if there exist $1 \leq i,j \leq n$ such that $\mathcal{P}_i= \mathcal{P}_j$, then our path is not a geodesic. Now, suppose that $\mathcal{P}_i \neq \mathcal{P}_j$ for every $i \neq j$ but $\sigma(\mathcal{P}_i)=D_i$ for some $1 \leq i \leq n$. Then
$$[\sigma,D_1, \ldots, D_{i-1}](\mathcal{P}_i)= \sigma(\mathcal{P}_i)=D_i=[\sigma, D_1, \ldots, D_i](\mathcal{P}_i),$$
hence $[\sigma,D_1, \ldots, D_{i-1}]=[\sigma,D_1, \ldots, D_i]$. Therefore, we can shorten our path as
$$\sigma, \ [\sigma,D_1], \ldots, [\sigma,D_1, \ldots, D_{i-1}], \ [\sigma,D_1, \ldots, D_{i-1},D_{i+1}], \ldots, [\sigma,D_1, \ldots, D_{i-1}, D_{i+1}, \ldots, D_n].$$
In particular, our path was not a geodesic.

\medskip \noindent
Conversely, it is clear that the distance between two orientations is at least equal to the number of walls on which they differ. On the other hand, if we suppose that, for every $i \neq j$, $\mathcal{P}_i \neq \mathcal{P}_j$ and $\mathcal{P}_i \neq \sigma(\mathcal{P}_i)$, then $\sigma$ and $[\sigma, D_1, \ldots, D_n]$ differ precisely on $\mathcal{P}_1, \ldots, \mathcal{P}_n$, so that our path must be a geodesic.
\end{proof}

\noindent
Notice that we have proved, in the last paragraph of the previous proof, that:

\begin{cor}\label{cor:orientationdist}
The distance between two vertices of $C(X, \mathfrak{P})$ is equal to the number of walls on which they differ.
\end{cor}

\noindent
The following lemma will be fundamental in the proof of Proposition \ref{prop:qcqm}.

\begin{lemma}\label{lem:orientationpermutation}
Let $\sigma\in C(X, \mathfrak{P})$ be an orientation and $D_1, \ldots, D_n$ some sectors such that
$$\sigma, \ [\sigma,D_1], \ [\sigma,D_1,D_2], \ldots, [\sigma,D_1, \ldots, D_n]$$ 
defines a geodesic in $C(X, \mathfrak{P})$. Fix some $1 \leq i \leq n$. If there is no $1 \leq j \leq n$ different from $i$ such that $D_j < D_i$, then there exists a permutation $\varphi$ of $\{ 1, \ldots, n \}$ such that 
$$\sigma, \ [\sigma, D_{\varphi(1)}], \ [\sigma, D_{\varphi(1)}, D_{\varphi(2)}], \ldots, [\sigma, D_{\varphi(1)}, \ldots, D_{\varphi(n)}]$$
defines a geodesic as well, with $\varphi(n)=i$. 
\end{lemma}

\begin{proof}
Let $\mathcal{P}_1, \ldots, \mathcal{P}_n$ denote the walls delimiting the sectors $D_1, \ldots, D_n$ respectively. We argue by induction over $n$. If $n=1$, there is nothing to prove. From now on, suppose that $n \geq 2$. 

\medskip \noindent
First, we claim that there exists some $k \neq i$ such that $\sigma(\mathcal{P}_k)$ is minimal in $\sigma(\mathfrak{P})$. Suppose by contradiction that this is not the case. Because we already know that $\sigma(\mathcal{P}_1)$ is minimal in $\sigma(\mathfrak{P})$, since $[\sigma,D_1]$ defines an orientation, necessarily $i=1$. Then, notice that, for every wall $\mathcal{Q} \in \mathfrak{P}$ and any $1\leq j \leq n$, $\sigma(\mathcal{Q}) < \sigma(\mathcal{P}_j)$ implies $\mathcal{Q} \in \{ \mathcal{P}_1, \ldots, \mathcal{P}_n \}$. Indeed, if $\mathcal{Q} \notin \{ \mathcal{P}_1, \ldots, \mathcal{P}_n \}$, we deduce from $[\sigma,D_1, \ldots, D_n](\mathcal{Q})= \sigma(\mathcal{Q})< \sigma(\mathcal{P}_j)$ that $\sigma(\mathcal{P}_j)=[\sigma,D_1, \ldots, D_n](\mathcal{P}_j)=D_j$, which is impossible according to Lemma \ref{lem:orientationgeod}. As a consequence, since for every $2 \leq j \leq n$ the sector $\sigma(\mathcal{P}_j)$ is not minimal in $\sigma(\mathfrak{P})$, we have $\sigma(\mathcal{P}_1)< \sigma(\mathcal{P}_j)$.

\medskip \noindent
In particular, $\mathcal{P}_1$ and $\mathcal{P}_j$ must be nested; let $A_j$ denote the sector delimited by $\mathcal{P}_1$ which contains $\mathcal{P}_j$, and notice that $\sigma(\mathcal{P}_j)$ is the sector delimited by $\mathcal{P}_j$ which contains $\mathcal{P}_1$. If $A_j=D_1$, then, because we know that $\sigma(\mathcal{P}_j) \neq D_j$ thanks to Lemma \ref{lem:orientationgeod}, we deduce that $D_j< A_j=D_1$, which is impossible. Otherwise, if $A_j \neq D_1$, we have $D_1 < \sigma(\mathcal{P}_j)$, so that $[\sigma,D_1, \ldots, D_n](\mathcal{P}_1)=D_1< \sigma(\mathcal{P}_j)$ implies that $\sigma(\mathcal{P}_j) = [\sigma, D_1, \ldots, D_n](\mathcal{P}_j)=D_j$, which is impossible according to Lemma \ref{lem:orientationgeod}.

\medskip \noindent
Thus, we have proved that there exists some $k \neq i$ such that $\sigma(\mathcal{P}_k)$ is minimal in $\sigma(\mathfrak{P})$. As a consequence, $\sigma'=[\sigma,D_k]$ is an orientation. To conclude the proof by applying our induction hypothesis, it is sufficient to show that the path
$$\sigma', [\sigma',D_1], \ldots, [\sigma',D_1, \ldots, D_{k-1}], [\sigma',D_1, \ldots, D_{k-1},D_{k+1}], \ldots, [\sigma',D_1, \ldots, D_{k-1},D_{k+1}, \ldots, D_n]$$
defines a geodesic in $C(X,<, \mathfrak{P})$. We first need to verify that it defines a path in $C(X,<, \mathfrak{P})$. Notice that, for every $0 \leq r \leq n-k$ and every $\mathcal{P} \in \mathfrak{P}$, 
$$\begin{array}{lcl} [\sigma',D_1, \ldots, D_{k-1},D_{k+1},\ldots, D_{k+r}](\mathcal{P}) & = & \left\{ \begin{array}{cl} \sigma'(\mathcal{P})=\sigma(\mathcal{P}) & \text{if} \ \mathcal{P} \notin \{ \mathcal{P}_1, \ldots, \mathcal{P}_{k+r} \} \\ D_j & \text{if} \ \mathcal{P}= \mathcal{P}_j, \ j \in \{ 1, \ldots, k+r \} \backslash \{ k \} \\ \sigma'(\mathcal{P}_k)=D_k & \text{if} \ \mathcal{P}= \mathcal{P}_k \end{array} \right. \\ \\ & = & [\sigma,D_1, \ldots, D_{k+r}](\mathcal{P}), \end{array}$$
so $[\sigma',D_1, \ldots, D_{k-1}, D_{k+1},\ldots, D_{k+r}]$ is an orientation. Similarly, if $i<k$, we have $[\sigma',D_1, \ldots, D_i]=[\sigma,D_1, \ldots, D_i,D_k]$. Because we know that $[\sigma,D_1, \ldots, D_i]$ is an orientation, it is sufficient to show that $[\sigma,D_1, \ldots, D_i](\mathcal{P}_k)= \sigma(\mathcal{P}_k)$ is minimal in $[\sigma,D_1, \ldots, D_i](\mathfrak{P})$ in order to deduce that $[\sigma',D_1, \ldots, D_i]$ defines an orientation. So let $\mathcal{Q} \in \mathfrak{P}$ be a wall satisfying $[\sigma, D_1, \ldots, D_i](\mathcal{Q}) < \sigma(\mathcal{P}_k)$, and suppose by contradiction that $\mathcal{Q} \neq \mathcal{P}_k$. Because $\sigma(\mathcal{P}_k)$ is minimal in $\sigma(\mathfrak{P})$, necessarily $\mathcal{Q} \in \{ \mathcal{P}_1, \ldots, \mathcal{P}_i \}$, so that $D_j < \sigma(\mathcal{P}_k)$ for some $1 \leq j \leq i$. Now, we deduce from $[\sigma,D_1, \ldots, D_k](\mathcal{P}_j)=D_j<\sigma( \mathcal{P}_k)$ that $\sigma(\mathcal{P}_k) = [\sigma,D_1, \ldots, D_k](\mathcal{P}_k) =D_k$, contradicting Lemma \ref{lem:orientationgeod}.

\medskip \noindent
Finally, we can apply Lemma \ref{lem:orientationgeod} to conclude that our path is a geodesic. Indeed, we know that $\mathcal{P}_i \neq \mathcal{P}_j$ for every $i \neq j$, and, for every $i \in \{ 1, \ldots, n \} \backslash \{k \}$, we have $\sigma'(\mathcal{P}_i) = \sigma(\mathcal{P}_i) \neq D_i$. This concludes the proof.
\end{proof}

\noindent
It is clear from the definition of $C(X,<,\mathcal{P})$ that any of its edges is naturally labelled by a wall of $\mathfrak{P}$. Our last step before proving Proposition \ref{prop:qcqm} is to understand how behave these labels in the triangles and squares of $C(X,<, \mathfrak{P})$.

\begin{lemma}\label{lem:qctriangle}
The edges of a triangle in $C(X,<,\mathfrak{P})$ are labelled by the same wall.
\end{lemma}

\begin{proof}
Let $\alpha, \beta, \gamma \in C(X,<, \mathfrak{P})$ be three pairwise adjacent vertices. Let $\mathcal{A,B,C}$ denote the walls labelling the edges $(\beta,\gamma), (\alpha,\gamma),(\alpha,\beta)$ respectively. If $\mathcal{A} \neq \mathcal{B}$, then $\alpha$ and $\beta$ differ on two walls, but this is impossible according to Corollary \ref{cor:orientationdist} since $d(\alpha, \beta)=1$. Therefore, $\mathcal{A}= \mathcal{B}$. Similarly, we show that $\mathcal{B}= \mathcal{C}$, concluding the proof.
\end{proof}

\begin{lemma}\label{lem:qcsquare}
Two opposite edges of some square in $C(X,<,\mathfrak{P})$ are labelled by the same wall of $\mathfrak{P}$. Moreover, the two walls labelling the edges of some induced square are transverse.
\end{lemma}

\begin{proof}
Let $(\sigma,\mu,\nu, \xi)$ be a square in $C(X,<, \mathfrak{P})$. If $\mu$ and $\xi$, or $\sigma$ and $\nu$, are adjacent, then we deduce from Lemma \ref{lem:qctriangle} that all the walls labelling the edges of our square are identical. From now on, we will suppose that our square is induced. As a consequence, $\sigma, \mu, \nu$ and $\sigma, \xi,\nu$ define two geodesics between $\sigma$ and $\nu$. According to Lemma \ref{lem:orientationgeod}, there exist $A_1 \in \mathcal{P}_1$, $A_2 \in \mathcal{P}_2$, $B_1 \in \mathcal{Q}_1$ and $B_2 \in \mathcal{Q}_2$ such that $\mu= [\sigma,B_1]$, $\nu = [\sigma,B_1,B_2]$, $\xi = [\sigma,A_1]$ and $\nu=[\sigma,A_1,A_2]$, with $\mathcal{P}_1 \neq \mathcal{P}_2$, $\mathcal{Q}_1 \neq \mathcal{Q}_2$, $A_1 \neq \sigma(\mathcal{P}_1)$, $A_2 \neq \sigma(\mathcal{P}_2)$, $B_1 \neq \sigma(\mathcal{Q}_1)$ and $B_2 \neq \sigma(\mathcal{Q}_2)$. In particular, we deduce that $\sigma$ and $\nu = [\sigma,A_1,A_2] = [\sigma,B_1,B_2]$ differ on the walls $\{\mathcal{P}_1, \mathcal{P}_2, \mathcal{Q}_1, \mathcal{Q}_2 \}$. On the other hand, we know that they differ on only two walls according to Corollary \ref{cor:orientationdist} since $d(\sigma,\nu)=2$, so we deduce that $\{ \mathcal{P}_1, \mathcal{P}_2 \} = \{ \mathcal{Q}_1, \mathcal{Q}_2 \}$. Next, we know similarly that $\mu$ and $\xi$ must differ on exactly two walls, so we deduce from
$$[\sigma,A_1](\mathcal{P})= \left\{ \begin{array}{cl} \sigma(\mathcal{P}) & \text{if} \ \mathcal{P} \neq \mathcal{P}_1 \\ A_1 & \text{if} \ \mathcal{P}= \mathcal{P}_1 \end{array} \right. \ \text{and} \ [\sigma,B_1](\mathcal{P}) = \left\{ \begin{array}{cl} \sigma(\mathcal{P}) & \text{if} \ \mathcal{P} \neq \mathcal{Q}_1 \\ B_1 & \text{if} \ \mathcal{P}= \mathcal{Q}_1 \end{array} \right.$$
that $\mathcal{P}_1$ and $\mathcal{Q}_1$ must be different. Therefore, $\mathcal{Q}_1= \mathcal{P}_2$ and $ \mathcal{Q}_2= \mathcal{P}_1$. As a consequence, notice that 
$$B_1=[\sigma,B_1,B_2](\mathcal{P}_2)=\nu(\mathcal{P}_2)=[\sigma,A_1,A_2](\mathcal{P}_2)=A_2,$$
and similarly $B_2=A_1$. Thus, $[\sigma,A_1,A_2]=[\sigma,A_2,A_1]$. 

\medskip \noindent
Now, suppose by contradiction that $\mathcal{P}_1$ and $\mathcal{P}_2$ are nested, and let $C_1$ (resp. $C_2$) denote the sector delimited by $\mathcal{P}_1$ (resp. $\mathcal{P}_2$) which contains $\mathcal{P}_2$ (resp. $\mathcal{P}_1$). Two cases may happen.

\medskip \noindent
Suppose that $A_1=C_1$ and $A_2=C_2$. Then $\sigma(\mathcal{P}_1) \neq A_1$ implies $\sigma(\mathcal{P}_1)<C_2=A_2$, hence $\sigma(\mathcal{P}_2)=A_2$, a contradiction.

\medskip \noindent
Suppose that either $A_1 \neq C_1$ or $A_2 \neq C_2$. Because the two possibilities are symmetric, say that $A_2 \neq C_2$. We deduce from $[\sigma,A_1,A_2](\mathcal{P}_2)=A_2< C_1$ that $C_1= [\sigma,A_1,A_2](\mathcal{P}_1)= A_1$. As a consequence, $\sigma(\mathcal{P}_1) \neq A_1 = C_1$ implies $\sigma(\mathcal{P}_1)< C_2$, hence $C_2= \sigma(\mathcal{P}_2)$. But $\sigma(\mathcal{P}_1)< \sigma(\mathcal{P}_2)$ is impossible since $\sigma(\mathcal{P}_2)$ must be minimal in $\sigma(\mathfrak{P})$.
\end{proof}

\begin{proof}[Proof of Proposition \ref{prop:qcqm}.]
First, we want to prove that $C(X,<,\mathfrak{P})$ satisfies the triangle condition.
Let $\mu, \nu \in C(X,<, \mathfrak{P})$ be two adjacent orientations and $\sigma \in C(X,<, \mathfrak{P})$ a third orientation satisfying $d(\sigma,\mu)=d(\sigma, \nu)=k$. Because $d(\mu,\nu)=1$, the orientations $\mu$ and $\nu$ differ on a single wall $\mathcal{P}_0$, say $\nu=[\mu,D]$ where $D$ is a sector delimited by $\mathcal{P}_0$. We write $\mu= [\sigma,D_1, \ldots, D_k]$, where $$\sigma, \ [\sigma,D_1], \ [\sigma,D_1,D_2], \ldots, [\sigma,D_1, \ldots, D_k]$$ defines a geodesic between $\sigma$ and $\mu$. Notice that there exists some $1 \leq i \leq k$ such that the underlying wall of $D_i$ is $\mathcal{P}_0$, since otherwise the path $$\sigma, \ [\sigma,D_1], \ldots, [\sigma,D_1, \ldots, D_k], \ [\sigma,D_1, \ldots, D_k,D]$$ would define a geodesic of length $k+1$ between $\sigma$ and $\nu$, according to Lemma \ref{lem:orientationgeod}. Notice that $$ \nu = [\mu,D] = [\sigma,D_1, \ldots, D_k,D].$$ In particular, because $D_i$ and $D$ have the same underlying wall (and $\mathcal{P}_0$ is not the underlying wall of any other $D_j$), $D_i=[\sigma,D_1, \ldots, D_k](\mathcal{P}_0)$ must be minimal in $[\sigma,D_1, \ldots, D_k](\mathfrak{P})$ according to Lemma \ref{lem:orientationmin}, so that $D_j \subset D_i$ for no $j \neq i$. It follows from Lemma \ref{lem:orientationpermutation} that there exists a permutation $\varphi$ of $\{ 1, \ldots, k \}$ satisfying $\varphi(k)=i$ such that 
$$\sigma, \ [\sigma,D_{\varphi(1)}], \ [\sigma,D_{\varphi(1)},D_{\varphi(2)}], \ldots, [\sigma,D_{\varphi(1)}, \ldots, D_{\varphi(k)}]$$ 
defines a geodesic between $\sigma$ and $\mu$. Let $\xi = [\sigma,D_{\varphi(1)}, \ldots, D_{\varphi(k-1)}]$. Then $d(\sigma, \xi)=d(\sigma,\mu)-1=k-1$, and $[\xi,D_i]=\mu$, and $[\xi,D]=\nu$. Therefore, $\xi$ is the orientation we are looking for.

\medskip \noindent
Now, we want to prove that $C(X,<,\mathfrak{P})$ satisfies the quadrangle condition.
Let $\alpha,\beta \in C(X,<, \mathfrak{P})$ be two orientations both adjacent to a third one $\gamma \in C(X,<, \mathfrak{P})$, and let $\sigma \in C(X,<, \mathfrak{P})$ be a last orientation satisfying $d(\sigma,\alpha)= d(\sigma,\beta)=k$, $d(\sigma,\gamma)=k+1$. Because we already know that the triangle condition holds, we will suppose that $\alpha$ and $\beta$ are not adjacent, ie., $d(\alpha, \beta)=2$. Fix some geodesic 
$$\sigma, \ [\sigma,D_1], \ [\sigma,D_1,D_2], \ldots, [\sigma,D_1, \ldots, D_k]=\alpha$$
between $\sigma$ and $\alpha$. Let $S_1,S_2$ be some sectors such that $\gamma=[\alpha,S_1]$ and $\beta=[\gamma,S_2]$, and let $\mathcal{P}_1$ and $\mathcal{P}_2$ denote the walls delimiting $S_1$ and $S_2$ respectively. Because $d(\sigma,\gamma)=k+1$, the concatenation of a geodesic between $\sigma$ and $\alpha$ with the edge between $\alpha$ and $\gamma$ must be a geodesic, so we deduce from Lemma \ref{lem:orientationgeod} that, for every $1 \leq i \leq k$, the underlying wall of $D_i$ is different from $\mathcal{P}_1$. Now, concatenating this geodesic with the edge $(\gamma,\beta)$ produces a path which cannot be a geodesic, so that Lemma \ref{lem:orientationgeod} implies that $\mathcal{P}_2$ is the underlying wall of $S_1$ or of some $D_i$. Notice that, because $\alpha$ and $\beta$ are not adjacent, $\mathcal{P}_1 \neq \mathcal{P}_2$, so there exists some $1 \leq j \leq k$ such that the underlying wall of $D_j$ is $\mathcal{P}_2$. Since we know that
$$\beta= [\gamma,S_2]=[\alpha,S_1,S_2]=[\sigma,D_1, \ldots, D_k,S_1,S_2],$$
and that $\beta$ and $\gamma$ are two orientations, necessarily $D_j=[\sigma,D_1, \ldots, D_k,S_1](\mathcal{P}_2)$ is minimal in $[\sigma,D_1, \ldots, D_k,S_1] (\mathfrak{P})$. As a consequence, if $D_j$ is not minimal in $[\sigma,D_1, \ldots, D_k](\mathfrak{P})$, necessarily
$$\sigma(\mathcal{P}_1)=[\sigma,D_1, \ldots, D_k](\mathcal{P}_1) < D_j,$$ hence
$D_j= \sigma(\mathcal{Q}_j)$, where $\mathcal{Q}_j$ denotes the underlying wall of $D_j$, which contradicts Lemma \ref{lem:orientationgeod}. Therefore, $D_j$ is minimal in $[\sigma,D_1, \ldots, D_k](\mathfrak{P})$, so that $D_j > D_i$ does not hold for any $1 \leq i \leq k$. According to Lemma \ref{lem:orientationpermutation}, there exists a permutation $\varphi$ of $\{1, \ldots, k\}$ satisfying $\varphi(k)=j$ such that 
$$\sigma, \ [\sigma, D_{\varphi(1)}], \ [\sigma, D_{\varphi(1)}, D_{\varphi(2)}], \ldots, [\sigma, D_{\varphi(1)}, \ldots, D_{\varphi(k)}]$$
defines a geodesic between $\sigma$ and $\beta$. Set $\xi=[\sigma,D_{\varphi(1)}, \ldots,D_{\varphi(k-1)}]$. Notice that $d(\sigma,\xi)=k-1$ and that $\xi$ is adjacent to $\alpha$ since $[\xi,D_j]= \alpha$. To conclude that $\xi$ is the orientation we are looking for, it is sufficient to show that $\xi$ and $\beta$ are adjacent.

\medskip \noindent
We know that $d(\sigma,\gamma)=k+1$, ie., $\sigma$ and $\gamma$ differ on $k+1$ walls, so, because $\beta=[\gamma,S_2]$, $\sigma(\mathcal{P}_2) \neq S_2$ would imply that $\sigma$ and $\beta$ differ on $k+1$ walls, contradicting the fact that $d(\sigma,\beta)= k$. Therefore, $\sigma(\mathcal{P}_2)=S_2$. As a consequence, $[\xi,S_1]$ is equal to $\beta = [\alpha,S_1,S_2]= [\xi,D_j,S_1,S_2]$. Indeed, these two orientations may only differ on $\mathcal{P}_1$ and $\mathcal{P}_2$, but $$[\xi,S_1](\mathcal{P}_1) = S_1 = [\xi,D_j,S_1,S_2](\mathcal{P}_1),$$ and $$[\xi,S_1](\mathcal{P}_2)=\xi(\mathcal{P}_2)=\sigma(\mathcal{P}_2)=S_2 = [\xi,D_j,S_1,S_2](\mathcal{P}_2).$$ This concludes the proof of the quadrangle condition.

\medskip \noindent
Finally, we need to verify that $C(X,<,\mathfrak{P})$ does not contain induced subgraphs isomorphic to $K_4^-$ or $K_{2,3}$. If $Y \subset X$ is a subgraph isomorphic to $K_4^-$, it follows from Lemma \ref{lem:qctriangle} that all its edges are labelled by the same wall, which implies that all its vertices must be pairwise adjacent; in particular, $Y$ is not an induced subgraph. Next, if $Y$ is a subgraph isomorphic to $K_{2,3}$, we deduce from Lemma \ref{lem:qcsquare} that there exist two non-adjacent vertices of $Y$ which share a common neighbor along two edges labelled by the same wall; necessarily, these two vertices have to be adjacent in $X$, so that $Y$ is not an induced subgraph. This concludes the proof that any connected component of the quasi-cubulation $C(X,<, \mathfrak{P})$ is a quasi-median graph.
\end{proof}

\noindent
It is worth noticing that $C(X,<,\mathfrak{P})$ is in general not connected. Even worse, it may happen that no natural choice of a connected component is possible. The example to keep in mind is the following. Let $S$ be an infinite set. We consider the partition $\mathfrak{P}$ of the power set $2^S$ containing the walls
$$\{ \{ A \subset S \mid x \in A \}, \ \{ A \subset S \mid x \notin A \} \}$$
for every $x \in S$. The quasi-cubulation $C(2^S,\subset,\mathfrak{P})$ is naturally isometric to the graph whose vertices are the sequences which belong to $\{0,1\}^S$ and whose edges link two sequences which differ on a single coordinate, and its connected components are all isomorphic. Nevertheless, we show below that some specific connected components, if they exist, are more strongly related to the initial popset than others.  

\begin{definition}
Let $(X,<, \mathfrak{P})$ be a popset. An orientation $\sigma$ is \emph{well-founded}\index{Well-founded orientations} if, for every wall $\mathcal{P} \in \mathfrak{P}$, the set $\{ \mathcal{Q} \in \mathfrak{P} \mid \sigma(\mathcal{Q}) < \sigma(\mathcal{P}) \}$ is finite.
\end{definition}

\noindent
Unfortunately, a popset does not always admit a well-founded orientation; see \cite[Example 9.7]{Roller}. Nevertheless, following respectively \cite[Theorem 9.6]{Roller} and \cite[Proposition 9.4]{Roller}, it can be proved that a well-founded orientation exists whenever our popset $(X,<, \mathfrak{P})$ is countable and \emph{discrete} (ie., for every $A,B \in X$, the set $\{ C \mid A< C< B\}$ is finite) or whenever it is discrete and \emph{$\omega$-dimensional} (ie., every wall $\mathcal{P} \in \mathfrak{P}$ is finite and there does not exist an infinite collection of pairwise transverse walls). Moreover, in the latter case, the well-founded orientations correspond precisely to the orientations satisfying the \emph{descending chain condition}, as introduced in \cite{SageevCAT(0)} for pocsets. 

\begin{thm}\label{thm:quasicubulation}
Let $(X,<, \mathfrak{P})$ be a popset admitting a well-founded orientation, and let $Y$ denote a connected component of the quasi-cubulation $C(X,<,\mathfrak{P})$ which contains such an orientation. Then $Y$ is a quasi-median graph, and there is a natural bijection between the walls of $\mathfrak{P}$ and the hyperplanes of $Y$ which respects transversality and tangency.
\end{thm}

\noindent
The following definition gives a precise meaning of \emph{tangent walls}\index{Tangent walls (in a popset)}, extending the notion of tangent hyperplanes in CAT(0) cube complexes.

\begin{definition}
Let $(X,<, \mathfrak{P})$ be a popset. Two walls $\mathcal{P}_1, \mathcal{P}_2 \in \mathfrak{P}$ are \emph{tangent} if they are nested  and if there do not exist $A_1 \in \mathcal{P}_1$, $A_2 \in \mathcal{P}_2$, $\mathcal{P} \in \mathfrak{P}$ and $A \in \mathcal{P}$ such that $A_1 < A < A_2$ or $A_2 < A < A_1$.
\end{definition}

\begin{proof}[Proof of Theorem \ref{thm:quasicubulation}.]
We already know that $Y$ is a quasi-median graph thanks to Proposition \ref{prop:qcqm}.

\medskip \noindent
Notice that Lemma \ref{lem:qctriangle} and Lemma \ref{lem:qcsquare} imply that the edges of a given hyperplane of $C(X,<, \mathfrak{P})$ are labelled by the same wall of $\mathfrak{P}$, so that the hyperplanes of $C(X,<, \mathfrak{P})$ are naturally labelled by the walls of $\mathfrak{P}$. 

\medskip \noindent
We claim that two distinct hyperplanes of $C(X,<, \mathfrak{P})$ are labelled by distinct walls of $\mathfrak{P}$. More precisely, we will prove that two edges labelled by the same wall of $\mathfrak{P}$ are dual to the same hyperplane. 

\medskip \noindent
Let $e,f'$ be two edges such that $e=(\mu, [\mu,A])$ and $f'=(\nu',[\nu',B])$ for some $A,B \in \mathcal{P}$ and $\mathcal{P} \in \mathfrak{P}$. Notice that, setting $\nu=[\nu',\mu(\mathcal{P})]$, the edge $f=( \nu, [\nu,A])$ is dual to the same hyperplane as $f'$, since the edges
$$f'=(\nu',[\nu',B]), \ (\nu', [\nu', \mu(\mathcal{P})]), \ ([\nu', \mu(\mathcal{P})], [\nu', \mu(\mathcal{P}),A])=f$$
successively belong to the same triangle. Thus, it is sufficient to prove that $e$ and $f$ are dual to the same hyperplane. Let
$$\mu, \ [\mu,D_1], \ [\mu,D_1,D_2], \ldots, [\mu,D_1, \ldots, D_n]=\nu$$
be a geodesic between $\mu$ and $\nu$. For convenience, let $\mathcal{P}_i$ denote the underlying wall of $D_i$ for $1 \leq i \leq n$. According to Lemma \ref{lem:orientationgeod}, $\mathcal{P}_i \neq \mathcal{P}_j$ for every $i \neq j$; moreover, because $\mu(\mathcal{P})= \nu (\mathcal{P})$, necessarily $\mathcal{P} \neq \mathcal{P}_i$ for every $i$. 

\medskip \noindent
We argue by induction on $n$. If $n=0$, then $e=f$ and the hyperplanes dual to $e$ and $f$ are obviously the same. From now on, suppose that $n \geq 1$ and set $\sigma= [\mu,D_1, \ldots, D_{n-1}]$. Let $\mathcal{Q} \in \mathfrak{P}$ satisfy $\sigma(\mathcal{Q}) < \sigma(\mathcal{P})$. If $\mathcal{Q} \neq \mathcal{P}_n$, then
$$\nu(\mathcal{Q})= [\sigma,D_n](\mathcal{Q})= \sigma(\mathcal{Q}) < \sigma(\mathcal{P})= [\sigma,D_n](\mathcal{P})= \nu(\mathcal{P}).$$
Because we already know that $\nu(\mathcal{P})$ is minimal in $\nu(\mathfrak{P})$, we deduce that $\mathcal{Q}= \mathcal{P}$. On the other hand, 
$$\sigma(\mathcal{P}_n)=[\mu,D_1,\ldots, D_{n-1}](\mathcal{P}_n)= \mu(\mathcal{P}_n) \nless \mu(\mathcal{P})= [\mu,D_1,\ldots, D_{n-1}](\mathcal{P})= \sigma(\mathcal{P})$$
since we already know that $\mu(\mathcal{P})$ is minimal in $\mu(\mathfrak{P})$. Therefore, $\sigma(\mathcal{P})$ must be minimal in $\sigma(\mathfrak{P})$, so that $[\sigma,A]$ defines an orientation. Now, notice that, for every $\mathcal{Q} \notin \{ \mathcal{P}, \mathcal{P}_n \}$, we have $[\sigma,A]( \mathcal{Q}) = \sigma(\mathcal{Q})$ and
$$[\nu,A](\mathcal{Q})= \nu (\mathcal{Q})= [\sigma,D_n](\mathcal{Q})= \sigma(\mathcal{Q});$$
and $[\sigma,A](\mathcal{P})=A=[\nu,A](\mathcal{P})$; and finally, using Lemma \ref{lem:orientationgeod},
$$[\sigma,A](\mathcal{P}_n) = \mu(\mathcal{P}_n) \neq D_n= \nu(\mathcal{P}_n) = [\nu,A](\mathcal{P}_n).$$
Therefore, $[\sigma,A]$ and $[\nu,A]$ differ on a single wall of $\mathfrak{P}$, so that they must be adjacent. Thus, we have proved that the edges $(\sigma, [\sigma,A])$ and $f=(\nu,[\nu,A])$ are opposite sides of some square; on the other hand, our induction hypothesis implies that $e$ and $(\sigma,[\sigma,A])$ are dual to the same hyperplane. A fortiori, $e$ and $f$, and so $e$ and $f'$, are dual to the same hyperplane, concluding the proof of our claim.

\medskip \noindent
We claim that any wall $\mathcal{P} \in \mathfrak{P}$ labels some edge of $Y$. Let $\sigma$ denote a well-founded orientation which belongs to $Y$. Because $\sigma$ is well-founded, the set
$$I= \{ \mathcal{Q} \mid \sigma(\mathcal{Q}) < \sigma(\mathcal{P}) \} \backslash \{ \mathcal{P} \}$$
is finite, say $I= \{ \mathcal{Q}_1, \ldots, \mathcal{Q}_n \}$. Notice that, for every $1 \leq i \leq n$, $\sigma(\mathcal{Q}_i)< \sigma(\mathcal{P})$ implies that $\mathcal{Q}_i$ and $\mathcal{P}$ are nested; let $A_i$ denote the sector delimited by $\mathcal{Q}_i$ which contains $\mathcal{P}$. Notice that, for every $1 \leq i \leq n$, $\sigma(\mathcal{P})$ is the sector delimited by $\mathcal{P}$ which contains $\mathcal{Q}_i$.

\medskip \noindent
Setting $\mu=[\sigma,A_1, \ldots, A_n]$, we want to prove that $\mu$ is an orientation. So let $B_1 \in \mathcal{R}_1$ and $B_2 \in \mathcal{R}_2$ be two sectors respectively delimited by two walls $\mathcal{R}_1, \mathcal{R}_2$, satisfying $B_1< B_2$ and $B_1= \mu(\mathcal{R}_1)$. Our goal is to show that $B_2= \mu(\mathcal{R}_2)$. We distinguish two cases.

\medskip \noindent
Suppose that $\mathcal{R}_1 \notin \{ \mathcal{Q}_1, \ldots, \mathcal{Q}_n\}$. Then we deduce from $\sigma(\mathcal{R}_1) = \mu(\mathcal{R}_1)=B_1 < B_2$ that $B_2= \sigma(\mathcal{R}_2)$. If $\mathcal{R}_2 \notin \{ \mathcal{Q}_1, \ldots, \mathcal{Q}_n \}$, we conclude that $B_2= \sigma(\mathcal{R}_2)= \mu(\mathcal{R}_2)$. Otherwise, say $\mathcal{R}_2= \mathcal{Q}_i$, we deduce from $\sigma(\mathcal{R}_1)=B_1<B_2= \sigma(\mathcal{R}_2)=\sigma(\mathcal{Q}_i)< \sigma(\mathcal{P})$, hence $\mathcal{R}_1 \in \{ \mathcal{Q}_1, \ldots, \mathcal{Q}_n \}$, a contradiction.

\medskip \noindent
Suppose that $\mathcal{R}_1= \mathcal{Q}_i$ for some $1 \leq i \leq n$. Notice that $A_i= \mu(\mathcal{Q}_i)=\mu(\mathcal{R}_1)= B_1 <B_2$. Therefore, for every $D \in \mathcal{P} \backslash \{ \sigma(\mathcal{P}) \}$, we have $D<A_i<B_2$, so that $B_2$ is the sector of $\mathcal{R}_2$ which contains $\mathcal{P}$. If $\mathcal{R}_2= \mathcal{Q}_j$ for some $1 \leq j \leq n$, this means that $B_2=A_j$, hence $B_2=A_j=\mu(\mathcal{Q}_j)=\mu(\mathcal{R}_2)$. Now, suppose that $\mathcal{R}_2 \notin \{ \mathcal{Q}_1, \ldots, \mathcal{Q}_n \}$, ie., $\sigma(\mathcal{R}_2) \nless \sigma(\mathcal{P})$. Therefore, $\sigma(\mathcal{R}_2)$ must be the sector of $\mathcal{R}_2$ which contains $\mathcal{P}$, hence $B_2= \sigma(\mathcal{R}_2)= \mu(\mathcal{R}_2)$. 

\medskip \noindent
Thus, we have proved that $\mu$ is an orientation. Now, we notice that $\mu(\mathcal{P})= \sigma(\mathcal{P})$ is minimal in $\mu(\mathfrak{P})$. Indeed, if $\mathcal{Q} \in \mathfrak{P}$ is a wall satisfying $\mu(\mathcal{Q})< \sigma(\mathcal{P})$, then either $\mathcal{Q} \notin \{ \mathcal{Q}_1, \ldots, \mathcal{Q}_n \}$, so that $\sigma(\mathcal{Q})=\mu(\mathcal{Q})< \sigma(\mathcal{P})$ which implies $\mathcal{Q}= \mathcal{P}$; or $\mathcal{Q}= \mathcal{Q}_i$ for some $1 \leq i \leq n$, so that $A_i= \mu(\mathcal{Q}_i)= \mu(\mathcal{Q})< \sigma(\mathcal{P})$, which is impossible since $A_i$ and $\sigma(\mathcal{P})$ are not $<$-comparable according to Lemma \ref{lem:noncomparable}. Therefore, if we fix some $D \in \mathcal{P} \backslash \{ \sigma(\mathcal{P}) \}$, then the two orientations $\mu$ and $[\mu,D]$ define two vertices of $Y$ linked by an edge, which is clearly labelled by $\mathcal{P}$. This concludes the proof of our claim.

\medskip \noindent
We claim that two hyperplanes $J_1,J_2$ of $C(X,<,\mathfrak{P})$ are transverse if and only if the walls $\mathcal{P}_1, \mathcal{P}_2$ which label them are transverse as well.

\medskip \noindent
Suppose that $J_1$ and $J_2$ are transverse. Then there exists some square whose dual hyperplanes are $J_1$ and $J_2$. We deduce from Lemma \ref{lem:qcsquare} that $\mathcal{P}_1$ and $\mathcal{P}_2$ must be transverse. Conversely, suppose that $\mathcal{P}_1$ and $\mathcal{P}_2$ are transverse. Let $\sigma$ be a well-founded orientation which belongs to $Y$. Using exactly the same argument as above, we show that the set 
$$I= \{ \mathcal{Q} \mid \sigma(\mathcal{Q})< \sigma(\mathcal{P}_1) \} \backslash \{ \mathcal{P}_1 \}$$
is finite, say $I= \{ \mathcal{A}_1, \ldots, \mathcal{A}_n \}$, that $\mu=[\sigma,A_1, \ldots, A_n]$ is an orientation if $A_i$ denotes the sector of $\mathcal{A}_i$ containing $\mathcal{P}_1$ for every $1 \leq i \leq n$, and finally that $\sigma(\mathcal{P}_1)$ is minimal in $\mu(\mathfrak{P})$. Applying this argument once again, we show that the set
$$J= \{ \mathcal{Q} \mid \mu(\mathcal{Q}) < \mu(\mathcal{P}_2) \} \backslash \{ \mathcal{P}_2 \}$$
is finite, say $J= \{ \mathcal{B}_1, \ldots, \mathcal{B}_m \}$, that $\nu= [\mu,B_1, \ldots,B_m ]$ is an orientation if $B_i$ denotes the sector of $\mathcal{B}_i$ containing $\mathcal{P}_2$ for every $1 \leq i \leq m$, and finally that $\mu(\mathcal{P}_2)$ is minimal in $\nu(\mathfrak{P})$. Notice that, because $\mathcal{P}_1$ and $\mathcal{P}_2$ are transverse, they do not belong to $I \cup J$. As a consequence, $\nu(\mathcal{P}_2)= \mu(\mathcal{P}_2)$ is minimal in $\nu(\mathfrak{P})$. Now, we want to prove that $\nu(\mathcal{P}_1)=\sigma(\mathcal{P}_1)$ is minimal in $\nu(\mathfrak{P})$. 

\medskip \noindent
First, we notice that $I \cap J= \emptyset$. Indeed, suppose by contradiction that there exists a wall $\mathcal{Q} \in \mathfrak{P} \backslash \{ \mathcal{P}_1, \mathcal{P}_2 \}$ satisfying $\sigma(\mathcal{Q})< \sigma(\mathcal{P}_1)$ and $\mu(\mathcal{Q}) < \mu(\mathcal{P}_2)$. Let $D \in \mathcal{P}_1 \backslash \{ \sigma(\mathcal{P}_1) \}$. Because $\mu(\mathcal{Q})$ is the sector delimited by $\mathcal{Q}$ which contains $\mathcal{P}_1$, we deduce that $D < \mu(\mathcal{Q})<\mu(\mathcal{P}_2)$, so that $\mathcal{P}_1$ and $\mathcal{P}_2$ must be nested, a contradiction.

\medskip \noindent
Let $\mathcal{Q} \in \mathfrak{P}$ be a wall satisfying $\nu(\mathcal{Q}) < \sigma(\mathcal{P}_1)$. Our goal is to prove that $\mathcal{Q}= \mathcal{P}_1$. If $\mathcal{Q} \notin J$, then we deduce from $\mu(\mathcal{Q})= \nu(\mathcal{Q}) < \sigma(\mathcal{P}_1)$ that $\mathcal{Q}= \mathcal{P}_1$ because we already know that $\sigma(\mathcal{P}_1)$ is minimal in $\mu(\mathfrak{P})$. Now, suppose that $Q \in J$. Let $D \in \mathcal{P}_1 \backslash \{ \nu ( \mathcal{P}_1) \}$. From $\nu(\mathcal{Q})< \sigma(\mathcal{P}_1)$, we deduce that $\nu(\mathcal{P}_1)=\sigma(\mathcal{P}_1)$ is the sector delimited by $\mathcal{P}_1$ which contains $\mathcal{Q}$, so that $D$ must be included into the sector delimited by $\mathcal{Q}$ which contains $\mathcal{P}_1$. On the other hand, we observed that $I \cap J = \emptyset$, so that $\mathcal{Q} \notin I$, hence $\mu(\mathcal{Q}) = \sigma(\mathcal{Q}) \nless \sigma(\mathcal{P}_1)=\nu(\mathcal{P}_1)$. A fortiori, $\mu(\mathcal{Q})$ is the sector delimited by $\mathcal{Q}$ which contains $\mathcal{P}_1$, hence $D< \mu(\mathcal{Q})$. But $\mathcal{Q} \in J$ implies $\mu(\mathcal{Q}) < \mu(\mathcal{P}_2)$, so that $D< \mu(\mathcal{Q}) < \mu(\mathcal{P}_2)$. We conclude that $\mathcal{P}_1$ and $\mathcal{P}_2$ must be nested, a contradiction. Thus, we have proved that $\sigma(\mathcal{P}_1)$ is minimal in $\nu(\mathfrak{P})$.

\medskip \noindent
We conclude that, if we fix two sectors $D_1 \in \mathcal{P}_1 \backslash \{ \nu(\mathcal{P}_1) \}$ and $D_2 \in \mathcal{P}_2 \backslash \{ \nu(\mathcal{P}_2) \}$, then $[\nu,D_1]$ and $[\nu,D_2]$ are two orientations. Moreover, since $\mathcal{P}_1$ and $\mathcal{P}_2$ are transverse, we deduce that $\nu(\mathcal{P}_1)$ is minimal in $[\nu,D_2](\mathcal{P})$, so that $[\nu,D_1,D_2]$ is an orientation as well. Finally, the four orientations $\nu$, $[\nu,D_1]$, $[\nu,D_2]$ and $[\nu,D_1,D_2]=[\nu,D_2,D_1]$ define a square in $Y$ whose dual hyperplanes are $J_1$ and $J_2$. A fortiori, $J_1$ and $J_2$ must be transverse.

\medskip \noindent
Finally, we want to prove that two hyperplanes of $C(X,<,\mathfrak{P})$ are tangent if and only if the walls labelling them are tangent as well. 

\medskip \noindent
First of all, let us notice that, if $J$ is a hyperplane of $C(X,<,\mathfrak{P})$ labelled by some wall $\mathcal{P} \in \mathfrak{P}$, then the sectors delimited by $J$ are precisely the $\{ \sigma \mid \sigma(\mathcal{P})=A \}$, where $A \in \mathcal{P}$. Indeed, fix two orientations $\mu,\nu$ and some geodesic
$$\mu, \ [\mu, A_1], \ [\mu, A_1,A_2], \ldots, [\mu,A_1, \ldots, A_n] =\nu$$
between them. If $\mathcal{P}_1, \ldots, \mathcal{P}_n$ are the walls underlying $A_1, \ldots, A_n$ respectively, we deduce from Lemma \ref{lem:orientationgeod} that the walls on which $\mu$ and $\nu$ differ are precisely $\mathcal{P}_1, \ldots, \mathcal{P}_n$. On the other hand, $\mu$ and $\nu$ belong to the same sector delimited by $J$ if and only if $J$ intersects this geodesic, which is equivalent to $\mathcal{P}= \mathcal{P}_i$ for some $1 \leq i \leq n$. Therefore, $\mu$ and $\nu$ belong to the same sector delimited by $J$ if and only if $\mu(\mathcal{P})= \nu(\mathcal{P})$. This proves our claim. 

\medskip \noindent
As a consequence, we are able to prove that, if $\mathcal{P}_1$ and $\mathcal{P}_2$ are two walls which are not tangent, then the associated hyperplanes $J_1$ and $J_2$ respectively are not tangent as well. Indeed, if $\mathcal{P}_1$ and $\mathcal{P}_2$ are not tangent, either $\mathcal{P}_1$ and $\mathcal{P}_2$ are transverse or there exist $A_1 \in \mathcal{P}_1$, $A_2 \in \mathcal{P}_2$, $\mathcal{P} \in \mathfrak{P}$ and $A \in \mathcal{P}$ such that $A_1<A<A_2$ (up to switching $\mathcal{P}_1$ and $\mathcal{P}_2$). In the former case, we know that $J_1$ and $J_2$ must be transverse, so that they cannot be tangent. In the latter case, we deduce that $S_1 \subset S \subset S_2$, where $S_1 = \{ \sigma \mid \sigma(\mathcal{P}_1)=A_1 \}$ is a sector delimited by $J_1$, $S_2 = \{ \sigma \mid \sigma(\mathcal{P}_2)=A_2 \}$ a sector delimited by $J_2$ and $S = \{ \sigma \mid \sigma(\mathcal{P})=A \}$ a sector delimited by $J$, which proves that $J_1$ and $J_2$ are not tangent. Indeed, if $\sigma \in S_1$, then $\sigma(\mathcal{P}_1)=A_1< A$ implies that $\sigma(\mathcal{P}) = A$, ie., $\sigma \in S$; and similarly, if $\sigma \in S$, then $\sigma(\mathcal{P})=A< A_2$ implies that $\sigma(\mathcal{P}_2) = A_2$, ie., $\sigma \in S_2$.

\medskip \noindent
Conversely, we claim that, if $S_1 \subsetneq S_2$ are two sectors of $C(X,<, \mathfrak{P})$ such that $S_1 = \{ \sigma \mid \sigma(\mathcal{P}_1)=A_1\}$ and $S_2 = \{ \sigma \mid \sigma(\mathcal{P}_2)=A_2\}$ for some $\mathcal{P}_1, \mathcal{P}_2 \in \mathfrak{P}$, $A_1 \in \mathcal{P}_1$ and $A_2 \in \mathcal{P}_2$, then $A_1 < A_2$. For this purpose, notice that $\mathcal{P}_1$ and $\mathcal{P}_2$ must be nested and let $B_1$ (resp. $B_2$) denote the sector of $\mathcal{P}_1$ (resp. $\mathcal{P}_2$) containing $\mathcal{P}_2$ (resp. $\mathcal{P}_1$). Fix three orientations $\sigma \in S_1$, $\mu \in S_2 \backslash S_1$ and $\nu \notin S_2$. We distinguish three cases.
\begin{itemize}
	\item Suppose that $A_1 \neq B_1$. Then $\sigma(\mathcal{P}_1) = A_1 < B_2$, which implies that $B_2= \sigma(\mathcal{P}_2)=A_2$. Therefore, $A_1< B_2=A_2$.
	\item Suppose that $A_1= B_1$ and $A_2=B_2$. Then $\nu(\mathcal{P}_2) \neq A_2= B_2$ implies that $\nu(\mathcal{P}_2)< B_1=A_1$, hence $A_1 = \nu(\mathcal{P}_1)$, which is impossible since $\nu \notin S_1$.
	\item Suppose that $A_1 = B_1$ but $A_2 \neq B_2$. Then $\mu(\mathcal{P}_2)=A_2 \neq B_2$ implies that $\mu(\mathcal{P}_2) < B_1 = A_1$, hence $A_1= \mu( \mathcal{P}_1)$, which is impossible since $\mu \notin S_1$.
\end{itemize}
This concludes the proof of our claim.

\medskip \noindent
As a consequence, if $\mathcal{P}_1$ and $\mathcal{P}_2$ are two walls labelling two hyperplanes $J_1$ and $J_2$ which are not tangent, then $\mathcal{P}_1$ and $\mathcal{P}_2$ cannot be tangent. Indeed, if $J_1$ and $J_2$ are not tangent, either they are transverse or there exists a third hyperplane $J$ delimiting some sector $S$ satisfying $S_1 \subset S \subset S_2$ for some sectors $S_1,S_2$ delimited by $J_1,J_2$ respectively. In the former case, we know that $\mathcal{P}_1$ and $\mathcal{P}_2$ must be transverse as well, so that they cannot be transverse. In the latter case, if $\mathcal{P}$ denote the wall labelling $J$, we deduce from our previous claim that there exist $A_1 \in \mathcal{P}_1$, $A_2 \in \mathcal{P}_2$ and $A \in \mathcal{P}$ satisfying $A_1 < A < A_2$. A fortiori, $\mathcal{P}_1$ and $\mathcal{P}_2$ are not tangent. 
\end{proof}

\begin{remark}
Suppose that a group $G$ acts on $X$ by preserving the order $<$ and the partition $\mathfrak{P}$. Then $G$ naturally acts on the quasi-cubulation $C(X,<, \mathfrak{P})$ by isometries. Unfortunately, if $Y$ denotes a connected component of $C(X,<, \mathfrak{P})$ which contains a well-founded orientation, then $Y$ is not necessarily $G$-invariant. Nevertheless, if this is the case, then the bijection between the walls of $\mathfrak{P}$ and the hyperplanes of $Y$ produced by the previous theorem is $G$-equivariant.
\end{remark}

\noindent
Besides cubulating pocsets, a related method to construct CAT(0) cube complexes is cubulating \emph{spaces with walls}; see \cite{NicaWallspaces, ChatterjiNibloWallspaces}. We also generalize this construction. Before defining what we call a \emph{space of partitions}, we need to introduce some terminology. Fix some set $X$ and some collection of partitions $\mathfrak{P}$. Elements of $\mathfrak{P}$ will be referred to as \emph{partitions} and elements of a partition will be referred to as \emph{sectors}. Two elements $\mathcal{P}_1,\mathcal{P}_2$ are
\begin{itemize}
	\item \emph{indistinguishable} if they represent the same partition of $X$, otherwise they are \emph{distinguishable};
	\item \emph{nested} if there exist two sectors $A_1\in \mathcal{P}_1$ and $A_2 \in \mathcal{P}_2$ such that $D \subset A_1$ for every $D \in \mathcal{P}_2 \backslash \{ A_2 \}$ and $D \subset A_2$ for every $D \in \mathcal{P}_1 \backslash \{ A_1 \}$. 
\end{itemize}
Finally, we say that a partition $\mathcal{P} \in \mathfrak{P}$ separates two points of $X$ if they belong to different sectors delimited by $\mathcal{P}$.

\begin{definition}\label{def:spacepartitions}
A \emph{space with partitions}\index{Spaces with partitions} $(X,\mathfrak{P})$ is the data of a set $X$ and a collection of partitions $\mathfrak{P}$ satisfying the following conditions:
\begin{itemize}
	\item every $\mathcal{P} \in \mathfrak{P}$ satisfies $\# \mathcal{P} \geq 2$ and $\emptyset \notin \mathcal{P}$;
	\item for every distinguishable partitions $\mathcal{P}_1, \mathcal{P}_2 \in \mathfrak{P}$, if there exist two sectors $A_1 \in \mathcal{P}_1$, $A_2 \in \mathcal{P}_2$ such that $A_1 \subset A_2$, then $\mathcal{P}_1$ and $\mathcal{P}_2$ are nested.
	\item two points of $X$ are separated by finitely many partitions of $\mathfrak{P}$.
\end{itemize}
\end{definition}

\noindent
The third assumption allows us to define the pseudo-distance $d_{\mathfrak{P}}$ on $X$, counting the number of partitions separating two points of $X$. A quasi-median graph $X$ is naturally a space with partitions, without indistinguishable partitions: for the set of partitions $\mathfrak{P}$, consider the sector-decompositions induced by the hyperplanes of $X$. In this case, the associated pseudo-distance coincides with the initial distance defined on $X$.

\begin{figure}
\begin{center}
\includegraphics[scale=0.6]{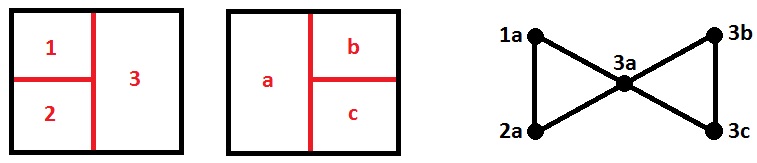}
\end{center}
\caption{Two partitions on the square, and the associated quasi-cubulation.}
\label{figure5}
\end{figure}

\medskip \noindent
We associate to any space with partitions $(X, \mathfrak{P})$ a popset in the following way. Let $\mathfrak{D}$ denote the set of all the sectors delimited by the partitions of $\mathfrak{P}$; by convention, a sector delimited by two distinct partitions is counted twice. In particular, $\mathfrak{P}$ defines a partition of $\mathfrak{D}$. The order $<$ we define on $\mathfrak{D}$ is: for every sectors $D_1,D_2 \in \mathfrak{D}$, we set $D_1<D_2$ if $D_1 \subset D_2$ and if the underlying partitions of $D_1$ and $D_2$ are distinguishable. Then it essentially follows from the definition of spaces with partitions that $(\mathfrak{D},<, \mathfrak{P})$ is a popset. An orientation of $(\mathfrak{D},<, \mathfrak{P})$ can be thought of as a map which chooses, for each partition of $\mathfrak{P}$, a sector it delimits in such a way that, whenever $\mathcal{P}_1$ and $\mathcal{P}_2$ are two distinguishable partitions delimiting respectively two sectors $A_1,A_2$ satisfying $A_1 <  A_2$, then $A_1= \sigma(\mathcal{P}_1)$ implies $A_2= \sigma(\mathcal{P}_2)$. A simple example of a space with partitions and the quasi-cubulation of its associated popset is given in Figure \ref{figure5}.

\medskip \noindent
If $x \in X$, we can define the orientation $\sigma_x$ associating to any partition of $\mathfrak{P}$ the sector which contains $x$. Such an orientation is called a \emph{principal orientation}. These orientations will allow us to choose a canonical connected component of our quasi-cubulation.

\begin{lemma}\label{lem:principalorientation1}
Any principal orientation defines a well-founded orientation.
\end{lemma}

\begin{proof}
Let $x \in X$ be a point and let $\sigma_x$ denote the associated principal orientation. Fixing some partition $\mathcal{P} \in \mathfrak{P}$, we want to prove that the set $\{ \mathcal{Q} \in \mathfrak{P} \mid \sigma_x(\mathcal{Q}) < \sigma_x(\mathcal{P}) \}$ is finite. Because $\# \mathcal{P} \geq 2$ and $\emptyset \notin \mathcal{P}$, there exists a point $y \in X$ which does not belong to $\sigma_x(\mathcal{P})$. Notice that, if $\mathcal{Q}$ belongs to our set, then $x \in \sigma_x(\mathcal{Q})$ but $y \notin \sigma_x(\mathcal{Q})$. Therefore, $\mathcal{Q}$ separates $x$ and $y$. Since there exist only finitely many partitions of $\mathfrak{P}$ separating two given points of $X$, this concludes the proof.
\end{proof}

\begin{lemma}\label{lem:principalorientation2}
Two principal orientations differ only on finitely many partitions.
\end{lemma}

\begin{proof}
Let $x,y \in X$ be two points and let $\sigma_x,\sigma_y$ denote the associated principal orientations. If $\mathcal{P} \in \mathfrak{P}$ is a partition satisfying $\sigma_x(\mathcal{P}) \neq \sigma_y(\mathcal{P})$, then $x$ and $y$ must belong to different sectors delimited by $\mathcal{P}$, ie., $\mathcal{P}$ separates $x$ and $y$. Because there exist only finitely many partitions separating two given points of $X$, we conclude that $\sigma_x$ and $\sigma_y$ differ on only finitely many partitions.
\end{proof}

\noindent
Thus, there is a canonical choice of a connected component of the quasi-cubulation of the popset associated to a space with partitions $(X, \mathfrak{P})$. Namely, this is the connected component which contains all the principal orientations, which we denote by $C(X, \mathfrak{P})$. It is worth noticing that, if a group $G$ acts on $X$ leaving $\mathfrak{P}$ invariant, then $G$ naturally acts on $C(X,\mathfrak{P})$ by isometries, since the set of principal orientations is $G$-invariant (indeed, if $x \in X$ and $g \in G$, then $g \cdot \sigma_x=\sigma_{g \cdot x}$). Moreover, according to the following lemma, the canonical map $(X, d_{\mathfrak{P}}) \hookrightarrow C(X, \mathfrak{P})$ defined by $x \mapsto \sigma_x$ is a (pseudo-)isometric embedding.

\begin{lemma}\label{lem:distspacepartitions}
For every $x,y \in X$, we have $d_{C(X, \mathfrak{P})}(\sigma_x,\sigma_y) = d_{\mathfrak{P}}(x,y)$.
\end{lemma}

\begin{proof}
The partitions on which $\sigma_x$ and $\sigma_y$ differ are precisely the partitions separating $x$ and $y$. Therefore, our equality follows from Corollary \ref{cor:orientationdist}. 
\end{proof}

\noindent
Important examples of spaces with partitions, as already mentionned, are quasi-median graphs themselves. Our next result characterize their quasi-cubulations.

\begin{prop}\label{prop:quasicubulatingqm}
Let $X$ be a quasi-median graph. The quasi-cubulation of $X$, viewed as a space with partitions, is isometric to $X$. 
\end{prop}

\noindent
Our proposition will be essentially a consequence of the following lemma.

\begin{lemma}\label{lem:anyorientationprincipal}
Let $(X, \mathfrak{P})$ be the space with partitions canonically associated to a quasi-median graph $X$. Any orientation $\sigma \in C(X, \mathfrak{P})$ is principal.
\end{lemma}

\begin{proof}
For convenience, we will identify $\mathfrak{P}$ with the set of the hyperplanes of $X$. It is worth noticing that two hyperplanes are transverse as elements of $\mathfrak{P}$ if and only if they are transverse as hyperplanes of $X$. Fix an arbitrary vertex $x \in X$, and let $\sigma_x$ denote its associated principal orientation. Because $\sigma \in C(X, \mathfrak{P})$, there exist finitely many hyperplanes $J_1, \ldots, J_n$ of $X$ on which $\sigma$ and $\sigma_x$ differ. Let $y$ denote the projection of $x$ onto $C=\bigcap\limits_{i=1}^n \sigma(J_i)$ (which is non empty because $\sigma$ is an orientation, thanks to Helly's property \ref{prop:Helly}). Notice that $J_1, \ldots, J_n$ separate $x$ and $y$. We claim that $\sigma= \sigma_y$. 

\medskip \noindent
Let $J$ be a hyperplane of $X$. If $J$ does not separate $x$ and $y$, then $\sigma(J)=\sigma_x(J)= \sigma_y(J)$. Otherwise, $J$ separates $x$ from $C$ according to Lemma \ref{lem:projseparate}. On the other hand, because $\sigma$ is an orientation, $\sigma(J) \cap \sigma(J_i) \neq \emptyset$ for every $1 \leq i \leq n$, so that Helly's property \ref{prop:Helly} implies that $\sigma(J) \cap C \neq \emptyset$. As a consequence, $\sigma(J)$ must be the sector delimited by $J$ which contains $C$; hence $\sigma(J)=\sigma_y(J)$. 

\medskip \noindent
Thus, $\sigma = \sigma_y$. We have proved that $\sigma$ is a principal orientation.
\end{proof}

\begin{proof}[Proof of Proposition \ref{prop:quasicubulatingqm}.]
Suppose that $X$ is a quasi-median graph and that $\mathfrak{P}$ is the collection of the sector-decompositions its the hyperplanes. Consider the map 
$$\Sigma : \left\{ \begin{array}{ccc} X & \longrightarrow & C(X, \mathfrak{P}) \\ x & \longmapsto & \sigma_x \end{array} \right.$$
which sends each vertex of $X$ to the associated principal orientation. Because, for every vertices $x,y \in X$, the distance between $\sigma_x$ and $\sigma_y$ in $C(X, \mathfrak{P})$ is equal to the number of hyperplanes separating $x$ and $y$ in $X$, which is equal to the distance between $x$ and $y$ in $X$, we know that $\Sigma$ defines an isometric embedding $X \hookrightarrow C(X, \mathfrak{P})$. Moreover, the surjectivity of $\Sigma$ follows from Lemma \ref{lem:anyorientationprincipal}, so $X$ and $C(X, \mathfrak{P})$ are isometric.
\end{proof}

\noindent
Proposition \ref{prop:quasicubulatingqm} will be useful to prove some results on quasi-median graphs. Let us mention a first application.

\begin{cor}\label{cor:finitesubgraph}
Let $X$ be a quasi-median graph and $Y \subset X$ a subgraph. Then $Y$ is finite if and only if the set of the sectors separating two vertices of $Y$ is finite.
\end{cor}

\begin{proof}
The implication is clear, so suppose the set of the sectors separating two vertices of $Y$ is finite. Fix some vertex $x \in Y$. For every $y \in Y$, let $D(y)$ denote the set of the sectors containing $y$ but not $x$. Set $\mathcal{D}= \bigcup\limits_{y \in Y} D(y)$. Notice that any sector of $\mathcal{D}$ separates two vertices of $Y$, hence $\# \mathcal{D}<+ \infty$ because $Y$ is finite and that there exist only finitely many sectors separating two given vertices of $Y$. On the other hand, the map
$$\left\{ \begin{array}{ccc} Y & \longrightarrow & 2^{\mathcal{D}} \\ y & \longmapsto & D(y) \end{array} \right.$$
is injective, since $\sigma_y=[\sigma_x,D(y)]$ for every $y \in Y$. Therefore, the finiteness of $\mathcal{D}$ implies the finiteness of $Y$.
\end{proof}

\begin{remark}\label{rem:wallspaces}
\emph{Pocsets} are popsets where all the walls have cardinality two, and \emph{spaces with walls} are spaces with partitions whose partitions all have cardinality two. In these cases, our quasi-cubulation coincides with the usual cubulations, as introduced respectively in \cite{SageevCAT(0)} and \cite{NicaWallspaces, ChatterjiNibloWallspaces}, and so produces a median graph. Alternatively, we can notice that the quasi-median graph produced by quasi-cubulating a pocset is triangle-free. Indeed, as a consequence of Lemma \ref{lem:qctriangle}, the edges of a triangle are labelled by a common wall, so that this wall necessarily delimits at least three sectors. Therefore, it follows from Corollary \ref{cor:whenmedian} that the quasi-median graph we obtain turns out to be a median graph.
\end{remark}

\subsection{Gated hulls}

\noindent
In this section, we notice that gated subgraphs of quasi-median graphs behave like combinatorially convex subcomplexes of CAT(0) cube complexes; see for instance \cite[Paragraph 2.3]{MR2413337}.

\begin{definition}
Let $X$ be a graph and $S \subset X$ a set of vertices. The \emph{gated hull}\index{Gated hulls} of $S$ is the smallest gated subgraph of $X$ containing $S$; alternatively, this is the intersection of all the gated subgraphs of $X$ which contain $S$.
\end{definition}

\noindent
The main result of this section is the following proposition.

\begin{prop}\label{prop:gatedhullhyp}
Let $X$ be a quasi-median graph and $S \subset X$ a subset. The hyperplanes of the gated hull $Y$ of $S$ are precisely the restrictions to $Y$ of the hyperplanes of $X$ separating two vertices of $S$. Moreover, two hyperplanes of $Y$ are transverse if and only if their extensions are transverse in $X$. 
\end{prop}

\noindent
We begin by proving the following lemma, which is well-known for CAT(0) cube complexes.

\begin{lemma}\label{lem:gatedhull}
The gated hull of a subset $S \subset X$ is the intersection of all the sectors containing $S$.
\end{lemma}

\begin{proof}
Let $Y$ denote the intersection of all the sectors containing $S$. Because the intersection of gated subgraphs is gated, it follows from Corollary \ref{cor:sectorgated} that $Y$ is gated. Now, we want to prove that, for any gated subgraph $Z$ containing $S$, necessarily $Y \subset Z$. Suppose by contradiction that there exists a vertex $x \in Y \backslash Z$. In particular, there exists a hyperplane $J$ separating $x$ from its projection onto $Z$, and we deduce from Lemma \ref{lem:projseparate} that $J$ separates $x$ from $Z$. A fortiori, since $S \subset Z$, $J$ separates $x$ and $S$, producing a sector containing $S$ but not $x \in Y$. This contradicts the definition of $Y$. 
\end{proof}

\begin{proof}[Proof of Proposition \ref{prop:gatedhullhyp}.] It is clear that a hyperplane of $Y$ extends to a hyperplane of $X$. However, we need to verify that two distinct hyperplanes of $Y$ extends to two distinct hyperplanes of $X$. More precisely, we want to prove that, if $e,f$ are two edges of $Y$ which are dual to the same hyperplane $J$ in $X$, then they are dual to the same hyperplane in $Y$. In fact, using the isomorphism given by Lemma \ref{lem:carrierproduct}, it is clear that the gated hull of $e \cup f$ in $N(J)$, which must be included into $Y$, contains a sequence of edges between $e$ and $f$ such that two consecutive edges either are two opposite sides of the same square or belong to the same triangle. A fortiori, $e$ and $f$ are dual to the same hyperplane in $Y$.

\medskip \noindent
Now, in order to conclude the proof of the first assertion in our proposition, it is sufficient to show that a hyperplane of $X$ separates two vertices of $Y$ if and only if it separates two vertices of $S$. If $J$ is a hyperplane which does not separate two vertices of $S$, then $S$ must be included into a sector delimited by $J$, and it follows from Lemma \ref{lem:gatedhull} that $J$ is disjoint from $Y$. Conversely, since $S \subset Y$, it is obvious that any hyperplane separating two vertices of $S$ separates two vertices of $Y$. 

\medskip \noindent
Let $J_1,J_2$ be two hyperplanes of $Y$, which we identify with their extensions to $X$ for convenience. Clearly, if they are transverse in $Y$ then they are transverse in $X$. Conversely, we suppose that $J_1,J_2$ are transverse in $X$ and we want to prove that they are transverse in $Y$ as well. 

\medskip \noindent
Fix two edges $(a,b), (c,d) \subset Y$ dual to $J_1, J_2$ respectively. Let $a',b'$ denote the projections onto $N(J_2)$ of $a,b$ respectively. Since $(a,b)$ is not dual to $J_2$, the vertices $a,b$ belong to the same sector delimited by $J_2$ according to Lemma \ref{lem:samesector}, so that $a',b'$ belong to the same connected component $\partial_1$ of $\partial J_2$. Let $\partial_2$ be a connected component of $\partial J_2$ which contains either $c$ or $d$ and which is different from $\partial_1$. Let $a'',b''$ denote the projections onto $\partial_2$ of $a,b$ respectively. Notice that, according to Corollary \ref{cor:projnested}, $a'',b''$ are also the projections onto $\partial_2$ of $a',b'$ respectively, so that $a'$ and $b'$ must be adjacent to $a''$ and $b''$ respectively. Then, because $J_1$ separates $a$ and $b$ and intersects $N(J_2)$, Proposition \ref{prop:2projseparate} implies that $a' \neq b'$ and $a'' \neq b''$; on the other hand, as a consequence of Corollary \ref{cor:projectionlip}, $a'$ and $b'$, and $a''$ and $b''$, must be adjacent. Therefore, the vertices $a',a'',b',b''$ define some square $Q$ whose hyperplanes are $J_1$ and $J_2$. Now, since $a'$ (resp. $b'$) belongs to a geodesic between $a$ and $c$ (resp. between $b$ and $c$), the convexity of $Y$ implies $a' \in Y$ (resp. $b' \in Y$). Similarly, $a''$ (resp. $b''$) belongs to a geodesic between $a$ and $\{c,d \} \cap \partial_2$ (resp. between $b$ and $\{ c,d \} \cap \partial_2$), hence $a'' \in Y$ (resp. $b'' \in Y$). Therefore, $Q \subset Y$, and we conclude that $J_1$ and $J_2$ are transverse in $Y$.
\end{proof}

\noindent
In CAT(0) cube complexes, the convex hull of a finite set is finite. Although we mentionned that gated subgraphs in quasi-median graphs are the analogue of convex subcomplexes in CAT(0) cube complexes, this assertion does not hold in the world of quasi-median graphs, ie., the gated hull of a finite set may be infinite. Indeed, the gated hull of two vertices which belong to an infinite clique will be the whole infinite clique. Nevertheless, some finiteness property can be deduced from the previous proposition.

\begin{definition}\label{def:cubicallyfinite}
A quasi-median graph is \emph{cubically finite}\index{Cubically finite subgraphs} if it contains finitely many hyperplanes.
\end{definition}

\noindent
The terminology is justified by the observation that a CAT(0) cube complex is finite if and only if it contains finitely many hyperplanes. Corollary \ref{cor:cubicallyfinite} will give an equivalent characterization of cubically finite quasi-median graphs.

\begin{cor}\label{cor:finitehull}
The gated hull of a finite subset is cubically finite.
\end{cor}

\begin{proof}
Let $S \subset X$ be a finite subset, and let $Y$ denote its gated hull. We deduce from Proposition \ref{prop:gatedhullhyp} that there exists a bijection between the hyperplanes of $Y$ and the hyperplanes of $X$ which separate at least two vertices of $S$. Because $S$ is finite and that only finitely many hyperplanes separate two given vertices, we conclude that $Y$ contains finitely many hyperplanes, ie., it is cubically finite.
\end{proof}

\subsection{Prisms}

\noindent
Recall that a \emph{weak Cartesian product} of graphs is a connected component of the Cartesian product of these graphs.

\begin{definition}\label{def:prism}
Let $X$ be a graph. An induced subgraph $Y \subset X$ is a \emph{prism}\index{Prisms} if it is a weak Cartesian product of cliques of $X$. In particular, $X$ is a \emph{prism} if it is isomorphic to a weak Cartesian product of complete graphs.
\end{definition}

\noindent
In this section, we want to show that cubes in CAT(0) cube complexes are replaced with prisms in quasi-median graphs. Our first main result is the following proposition.

\begin{prop}\label{prop:transversehypcube}
Let $X$ be a quasi-median graph. If $J_1, \ldots, J_n$ is a collection of pairwise transverse hyperplanes, then $X$ contains a prism whose dual hyperplanes are precisely $J_1, \ldots, J_n$. 
\end{prop}

\noindent
We begin by proving the following lemma.

\begin{lemma}\label{lem:prismhyp}
Let $X$ be a quasi-median graph whose hyperplanes are pairwise transverse. Let $\mathfrak{J}$ denote its collection of hyperplanes, and for every $J \in \mathfrak{J}$, let $K(J)$ be the complete graph whose vertices are the sectors delimited by $J$. Then $X$ is isometric to a connected component of the Cartesian product $\prod\limits_{J \in \mathfrak{J}} K(J)$. In particular, $X$ is a prism.
\end{lemma}

\begin{proof}
According to Proposition \ref{prop:quasicubulatingqm}, we can identify $X$ with the quasi-cubulation $C(X, \mathfrak{P})$ of the canonical space of partitions $(X, \mathfrak{P})$ associated to $X$. In particular, the vertices of $\prod\limits_{J \in \mathfrak{J}} K(J)$ can be naturally thought of maps associating to any hyperplane one of its sectors. If $x \in X$ is a base vertex, we denote by $P$ the connected component of $\prod\limits_{J \in \mathfrak{J}} K(J)$ which contains $\sigma_x$, ie., the connected component containing all the principal orientations, since two principal orientations differ only on finitely many hyperplanes according to Lemma \ref{lem:principalorientation2}. Now, we defined
$$\Psi : \left\{ \begin{array}{ccc} X & \longrightarrow & P \\ y & \longmapsto & \sigma_y \end{array} \right..$$
Notice that $\Psi$ defines an isometric embedding, since, for every vertices $y,z \in X$,
$$\begin{array}{lcl} d( \Psi(y), \Psi(z)) & = & d(\sigma_y, \sigma_z) \\ & = & \# \{ \text{hyperplanes on which $\sigma_y$ and $\sigma_z$ differ} \} \\ & = & \# \{ \text{hyperplanes separating $y$ and $z$} \} \\ & = & d(y,z) \end{array}$$
Now, we want to prove that $\Psi$ is surjective. Let $p \in P$ be a vertex. Because $p$ and $\sigma_x$ belong to the same connected componenent of $\prod\limits_{J \in \mathfrak{J}} K(J)$, $p$ differ from $\sigma_x$ on finitely many coordinates $p_1, \ldots, p_k$. Set
$$\sigma = [\sigma_x,p_1, \ldots, p_k].$$
We claim that $\sigma$ defines an orientation. Notice that, if $\mu$ is an orientation and $D$ a sector, because the hyperplanes of $X$ are pairwise transverse, no two elements of $(\mu(\mathfrak{P}), \subset)$ are comparable, so that any element turns out to be minimal; It follows from Lemma \ref{lem:orientationmin} that $[\mu,D]$ always defines an orientation. By applying this observation successively to $\sigma_x$, to $[\sigma_x,p_1]$, to $[\sigma_x,p_1,p_2]$, and so on, it follows that $\sigma$ defines an orientation.

\medskip \noindent
Thus, we deduce from Lemma \ref{lem:anyorientationprincipal} that $\sigma=\sigma_y$ for some vertex $y \in X$. By construction, we have $\Psi(y)=p$, which concludes the proof.
\end{proof}

\begin{proof}[Proof of Proposition \ref{prop:transversehypcube}.]
Because $J_1, \ldots, J_n$ are pairwise transverse, we deduce from Helly's property \ref{prop:Helly} that $\bigcap\limits_{i=1}^n N(J_i) \neq \emptyset$. Let $x \in \bigcap\limits_{i=1}^n N(J_i)$ be a vertex. For every $1 \leq i \leq n$, fix some sector $D_i$ delimited by $J_i$ which does not contain $x$, and let $x_i \in D_i$ be a vertex adjacent to $x$ (such a vertex exists because $x \in N(J_i)$). Finally, let $C$ denote the gated hull of $\{x,x_1, \ldots, x_n\}$. We claim that $C$ is the prism we are looking for. In fact, we are going to prove the following statement:

\begin{fact}\label{fact:spanningprism}
If $J_1, \ldots, J_n$ is a collection of pairwise transverse hyperplanes, $x \in \bigcap\limits_{i=1}^n N(J_i)$ a vertex, and, for every $1 \leq i \leq n$, $x_i$ a vertex adjacent to $x$ separated from it by $J_i$, then the gated hull of $\{x,x_1, \ldots, x_n\}$ is a prism whose dual hyperplanes are $J_1, \ldots, J_n$.
\end{fact}

\noindent
According to Proposition \ref{prop:gatedhullhyp}, the hyperplanes of $C$ are naturally identified to the hyperplanes of $X$ separating $x$ and $x_i$ or $x_i$ and $x_j$ for some $1 \leq i,j \leq n$. In particular, $J_1, \ldots, J_n$ are hyperplanes of $C$. Now, notice that, if for some $1 \leq i,j \leq n$ the vertices $x_i$ and $x_j$ are adjacent, then $J_i=J_j$, which is impossible; so $d(x_i,x_j)=2$. As a consequence, if $i \neq j$, the concatenation $(x_i,x) \cup (x,x_j)$ is a geodesic, so that any hyperplane separating $x_i$ and $x_j$ must separate either $x$ and $x_i$ or $x$ and $x_j$. We deduce that the hyperplanes of $C$ are precisely $J_1, \ldots, J_n$. 

\medskip \noindent
Moreover, we know from Proposition \ref{prop:gatedhullhyp} that two hyperplanes are transverse in a gated subgraph if and only if they are transverse in the whole graph. Therefore, the hyperplanes of $C$ are pairwise transverse, and it follows from Lemma \ref{lem:prismhyp} that $C$ is a prism on its own right. We conclude that $C$ is a prism of $X$ because, as $C$ is a gated subgraph, it has to be a induced and a clique in $C$ must be a clique in $X$.
\end{proof}

\begin{remark}
Proposition \ref{prop:transversehypcube} does not hold for infinite collections of pairwise transverse hyperplanes. For instance, there exist CAT(0) cube complexes containing infinite collections of pairwise transverse hyperplanes but containing no infinite cubes. See \cite[Figure6]{HruskaWise} for a simple explicit example, or \cite[Section A.1]{article3} for examples admitting interesting group actions.
\end{remark}

\noindent
It is worth noticing that, during the proof of Lemma \ref{lem:prismhyp}, the hypothesis on the hyperplanes of $X$ was not used to prove that $\Psi$ is an isometric embedding. Therefore, any quasi-median graph embeds isometrically into a prism. Let us give an alternative proof of this observation, which will be useful later.

\begin{lemma}\label{lem:embedinfiniteprism}
Let $X$ be a quasi-median graph. Let $\mathfrak{J}$ denote the collection of the hyperplanes of $X$, and, for every $J \in \mathfrak{J}$, fix a clique $C(J)$ dual to $J$. The map 
$$x \mapsto \left( \mathrm{proj}_{C(J)}(x) \right)_{J \in \mathfrak{J}}$$
defines an isometric embedding $X \hookrightarrow \prod\limits_{J \in \mathfrak{J}} C(J)$. 
\end{lemma}

\begin{proof}
Let $x,y \in X$ be two vertices. For every $J \in \mathfrak{J}$, we know that $\mathrm{proj}_{C(J)}(x) \neq \mathrm{proj}_{C(J)}(y)$ if and only if $J$ separates $x$ and $y$. Therefore, the number of coordinates on which $\left( \mathrm{proj}_{C(J)}(x) \right)_{J \in \mathfrak{J}}$ and $\left( \mathrm{proj}_{C(J)}(y) \right)_{J \in \mathfrak{J}}$ differ is equal to the number of hyperplanes separating $x$ and $y$. This precisely means that our map is an isometric embedding, concluding the proof.
\end{proof}

\noindent
In particular, Proposition \ref{prop:transversehypcube} justifies the following definition.

\begin{definition}
The \emph{cubical dimension}\index{Cubical dimension} of a quasi-median graph $X$, denoted by $\dim_{\square}X$, is the maximal number of pairwise intersecting hyperplanes.
\end{definition}

\noindent
We conclude this section by our second and last main result.

\begin{prop}\label{prop:maxprism}
Let $X$ be a quasi-median graph of finite cubical dimension. The map $\mathfrak{J} \mapsto \bigcap\limits_{J \in \mathfrak{J}} N(J)$ defines an $\mathrm{Aut}(X)$-equivariant bijection between the maximal collections of pairwise transverse hyperplanes and the maximal prisms of $X$.
\end{prop}

\noindent
To prove this proposition, the following lemma will be needed.

\begin{lemma}\label{lem:prismgated}
In a quasi-median graph, a prism is gated.
\end{lemma}

\begin{proof}
Let $X$ be a quasi-median graph. Let $\mathcal{S}$ be a collection of cliques of $X$ and $C \subset X$ a prism which is a weak Cartesian product of the cliques in $\mathcal{S}$. 

\medskip \noindent
First, we claim that any edge $e$ of $C$ belongs to a clique of $C$ which is a clique of $X$. This is clear if there exists some $S \in \mathcal{S}$ such that $e \subset S$. Otherwise, $C$ contains a subgraph isomorphic to $S \times [0,1]$, for some $S \in \mathcal{S}$, such that $S \times \{0 \}=S$ and $e \subset S \times \{1 \}$. It follows from Lemma \ref{lem:parallelcliques} that, because $S$ is a clique in $X$, necessarily $S \times \{ 1 \}$ has to be a clique in $X$ as well. This proves our claim.

\medskip \noindent
Now, we can prove that $C$ contains its triangles. Indeed, if $(a,b,c)$ is a triangle with $(a,b) \subset C$, then, according to our previous claim, there exists some clique $Y$ of $C$ which is a clique of $X$ and which contains $(a,b)$. Because $X$ does not contain induced subgraphs isomorphic to $K_4^-$, necessarily $c \in Y$, hence $(a,b,c) \subset Y \subset C$.

\medskip \noindent
Finally, we prove that $C$ is locally convex, which is sufficient to conclude that $C$ is gated according to Proposition \ref{prop:gated}. Let $(a,b,c,d) \subset X$ be a square with $(a,b) \cup (b,c) \subset C$. If $a$ and $c$ are adjacent, we conclude that $(a,b,c,d) \subset Y$ because we already know that $C$ contains its triangles. From now on, we suppose that $a$ and $c$ are not adjacent. Because $C$ is a weak Cartesian product of complete graphs, two adjacent edges generate either a triangle or a square in $C$, so that there exists a vertex $x \in C$ defining a square $(a,b,c,x) \subset C$. If $d=x$, we conclude that $(a,b,c,d) \subset C$ and we are done. Suppose by contradiction that $d \neq x$. In particular, the vertices $a,b,c,d,x$ define a subgraph isomorphic to $K_{2,3}$. Because $X$ does not contain induced subgraphs isomorphic to $K_4^-$ and that we supposed that $a$ and $c$ are not adjacent, we deduce that neither $b$ and $d$, nor $b$ and $x$, nor $d$ and $x$, are adjacent, so that we find an induced subgraph in $X$ which is isomorphic to $K_{2,3}$.
\end{proof}

\begin{proof}[Proof of Proposition \ref{prop:maxprism}.]
Let $J_1, \ldots, J_n$ a maximal collection of pairwise transverse hyperplanes, and set $C= \bigcap\limits_{i=1}^n N(J_i)$. Notice that, if $J$ is a hyperplane intersecting $C$, then, for every $1 \leq i \leq n$, either $J=J_i$ or $J$ is transverse to $J_i$ since $J$ intersects $N(J_i)$. By maximality of our collection, we deduce that no hyperplane different from $J_1, \ldots, J_n$ intersects $C$; in particular, we deduce from Lemma \ref{lem:prismhyp} that $C$ is a prism. Conversely, it follows from Fact \ref{fact:spanningprism} that $J_1, \ldots, J_n$ intersect $C$, so that the hyperplanes dual to $C$ are precisely $J_1, \ldots, J_n$. Finally, the maximality of our collection $J_1, \ldots, J_n$ implies that $C$ cannot be strictly contained in a larger prism, so that $C$ is indeed a maximal prism of $X$. 

\medskip \noindent
The injectivity of our application is clear, because the hyperplanes dual to the final maximal prism are precisely the initial collection of hyperplanes. 

\medskip \noindent
Now, if $C$ is a maximal prism of $X$, and if $J_1, \ldots, J_n$ are the hyperplanes dual to $C$, we claim that $C= \bigcap\limits_{i=1}^n N(J_i)$, proving that our application is surjective. The inclusion $C \subset \bigcap\limits_{i=1}^n N(J_i)$ is clear. If we prove that $J_1, \ldots, J_n$ defines a maximal collection of pairwise transverse hyperplanes, then the reverse inclusion will be a consequence of the maximality of $C$, since we already noticed that $\bigcap\limits_{i=1}^n N(J_i)$ defines a prism in this case. More precisely, we will prove that, if $C$ is a prism of $X$ and $J_1, \ldots, J_n$ its dual hyperplanes, and if there exists a hyperplane $J$ transverse to $J_1, \ldots, J_n$, then $C$ is included into a larger prism, and in particular it is not maximal.

\medskip \noindent
Let $x \in C$ and $y \in N(J)$ be two vertices minimizing the distance between $C$ and $N(J)$. Because $C$ is gated according to Lemma \ref{lem:prismgated}, we apply Lemma \ref{lem:2min} to deduce that the hyperplanes separating $x$ and $y$ are precisely the hyperplanes separating $C$ and $N(J)$. Let $J_0$ denote the hyperplane dual to the first edge of some fixed geodesic $\gamma$ from $x$ to $y$. In particular, $J_0$ separates $x$ and $y$ so that it must separate $C$ and $N(J)$. On the other hand, the hyperplanes $J_1, \ldots, J_n$ intersect both $C$ and $N(J)$, so we deduce that $J_0$ has to be transverse to $J_1, \ldots, J_n$. For every $1 \leq i \leq n$, let $x_i \in C$ be a vertex adjacent to $x$ and separated from it by $J_i$; and let $x_0$ denote the vertex of $\gamma$ adjacent to $x$, so that $J_0$ is the hyperplane separating $x$ and $x_0$. If $C'$ denote the gated hull of $\{x,x_1, \ldots, x_n\}$ and $C''$ the gated hull of $\{x,x_0,x_1, \ldots, x_n\}$, we have $C \subseteq C' \subsetneq C''$. On the other hand, because $J_0,J_1, \ldots, J_n$ are pairwise transverse, we deduce from Fact \ref{fact:spanningprism} that $C''$ is a prism.
\end{proof}

\begin{cor}\label{cor:cubicallyfinite}
Let $X$ be a quasi-median graph. The following assertions are equivalent:
\begin{itemize}
	\item[(i)] $X$ is cubically finite; 
	\item[(ii)] $\dim_{\square}(X)<+ \infty$ and $X$ contains finitely many maximal prisms.
	\item[(iii)] every prism of $X$ is contained in a maximal prism and $X$ contains finitely many maximal prisms;
\end{itemize}
\end{cor}

\begin{proof}
Suppose that $X$ is cubically finite. Of course, the cubical dimension of $X$ is necessarily finite (and bounded above by the number of hyperplanes of $X$). As an immediate consequence of Proposition \ref{prop:maxprism}, it also follows that $X$ contains finitely many maximal prisms. This proves $(i) \Rightarrow (ii)$. Next, suppose $(ii)$. If $X$ contains a prism which is not contained in a maximal prism, then there must exist an increasing sequence of prisms in $X$, and looking at the hyperplanes dual to this sequence of prisms, we find an infinite collection of pairwise transverse hyperplanes, contradicting the finiteness of the cubical dimension. Therefore, $(ii)$ implies $(iii)$. Finally, suppose $(iii)$. If $X$ is covered by $N$ maximal prisms and if $D$ is the maximal cubical dimension of these prisms, then necessarily $X$ contains at most $N\cdot D$ hyperplanes. A fortiori, $X$ is cubically finite, hence $(iii) \Rightarrow (i)$.
\end{proof}

\subsection{Quasi-medians}\label{section:qm}

\noindent
Since quasi-median graphs may contain triangles, the median point of a triple of vertices is generally not well-defined. In this section, we prove nevertheless that any triple admits a \emph{quasi-median}. In fact, quasi-median graphs were essentially defined in \cite{quasimedian} in terms of quasi-medians, where various equivalent definitions are given, including the definition we use in this paper.

\begin{definition}\label{def:qmtriangle}
Let $X$ be a graph and $x,y,z \in X$ three vertices. A \emph{median triangle}\index{Median triangles} of $(x,y,z)$ is a triple of vertices $(x',y',z')$ satisfying
$$\left\{ \begin{array}{l} d(x,y)=d(x,x')+d(x',y')+d(y',y) \\  d(x,z)=d(x,x')+d(x',z')+d(z',z) \\ d(y,z)=d(y,y')+d(y',z')+d(z',z) \end{array} \right. .$$
The quantity $\min(d(x',y'),d(y',z'),d(z',x'))$ will be referred to as the \emph{size} of the median triangle. A \emph{quasi-median}\index{Quasi-median triangles} of $(x,y,z)$ is an \emph{equilateral} median triangle $(x',y',z')$ (ie., $d(x',y')=d(x',z')=d(y',z')$) of minimal size. 
\end{definition}

\noindent
It is worth noticing that any triple $(x,y,z)$ admits at least one median triangle, namely $(x,y,z)$ itself. Also, if $(x',y',z')$ denotes a median triangle of some triple $(x,y,z)$, then, for any choice of geodesics $[x,x']$, $[x',y']$, $[y',y]$ respectively between $x$ and $x'$, $x'$ and $y'$, $y'$ and $y$, the concatenation $[x,x'] \cup [x',y'] \cup [y',y']$ turns out to be a geodesic. And a similar statement holds for $x,x',z',z$ and $y,y',z',z$.
\begin{figure}
\begin{center}
\includegraphics[scale=0.6]{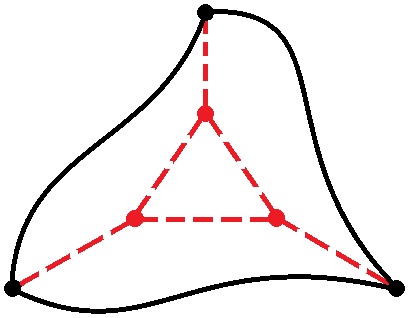}
\end{center}
\caption{A quasi-median triangle.}
\label{figure9}
\end{figure}
\begin{definition}
A graph is \emph{median}\index{Median graphs} if any triple of vertices admits a unique quasi-median of size zero, which we refer to as the \emph{median vertex}.
\end{definition}

\noindent
The main result of this section is the following:

\begin{prop}\label{prop:quasimedian}
In a quasi-median graph $X$, any triple of vertices $x,y,z$ admits a unique quasi-median. Moreover, its size is equal to the number of hyperplanes separating the three vertices $x,y,z$ and its gated hull is a prism of finite cubical dimension.
\end{prop}

\noindent
We begin by proving the following technical lemma.

\begin{lemma}\label{lem:projsectors}
Let $X$ be a quasi-median graph, $\mathcal{S}$ a finite collection of pairwise intersecting sectors, and $x \in X$ a vertex. Any hyperplane separating $x$ from its projection onto $\bigcap \mathcal{S}$ delimits a sector which contains an element of $\mathcal{S}$. 
\end{lemma}

\begin{proof}
Let us define inductively a finite sequence of vertices $(x_i)$ as follows. Set $x_0=x$. Suppose now that $x_i$ is defined for some $i \geq 0$. If $x_i \in \bigcap \mathcal{S}$, then $x_i$ is the last term of our sequence. Otherwise, there exists some $S_i \in \mathcal{S}$ such that $x_i \notin S_i$, and we define $x_{i+1}$ as the projection of $x_i$ onto $S_i$. Let $x_0, \ldots, x_n \in X$ be the sequence defined in this way (and $S_0, \ldots, S_{n-1} \in \mathcal{S}$ the corresponding sequence of sectors). Notice that $x_n \in \bigcap \mathcal{S}$

\medskip \noindent
First of all, let us notice that

\begin{fact}
For every $0 \leq i \leq n-1$ and $1 \leq k \leq n-i$, the vertex $x_{i+k}$ belongs to $S_i$.
\end{fact}

\noindent
We fix $i$ and argue by induction on $k$. Because $x_{i+1}$ is by construction the projection of $x_{i}$ onto $S_i$, of course $x_{i+1}$ belongs to $S_i$. Now, suppose that $x_{i+k}$ belongs to $S_i$. As a consequence of Lemma \ref{cor:projnested}, we have
$$\begin{array}{lcl} x_{i+k+1} & = & \mathrm{proj}_{S_{i+k}}(x_{i+k}) = \mathrm{proj}_{S_{i+k}} \circ \mathrm{proj}_{S_i} (x_{i+k}) \\ \\ & =& \mathrm{proj}_{S_{i+k} \cap S_i} (x_{i+k}) \in S_i \end{array}$$
which concludes the proof of our fact.

\medskip \noindent
For every $0 \leq i \leq n-1$, fix a geodesic $[x_i,x_{i+1}]$. We claim that the concatenation $[x_0,x_1] \cup \cdots \cup [x_{n-1},x_n]$ is a geodesic. Indeed, if $J$ is a hyperplane intersecting $[x_i,x_{i+1}]$ for some $0 \leq i \leq n-1$, then, since $x_{i+1}$ is the projection of $x_i$ onto $S_i$, it follows from Lemma \ref{lem:projseparate} that $J$ separates $x_i$ from $S_i$. On the other hand, we know from our previous fact that $x_j$ belongs to $S_i$ for every $j \geq i+1$, so that $J$ cannot intersect $[x_{i+1},x_{i+2}] \cup \cdots \cup [x_{n-1},x_n]$. This proves that our concatenation turns out to be a geodesic. As a consequence, any hyperplane separating $x_0$ and $x_n$ must separate $x_i$ and $x_{i+1}$ for some $0 \leq i \leq n-1$.

\medskip \noindent
Now, we claim that $x_n$ is the projection of $x=x_0$ onto $\bigcap \mathcal{S}$. Let $p$ denote this projection. Because $x_n \in \bigcap \mathcal{S}$, there exists a geodesic between $x$ and $x_n$ passing through $p$, so that any hyperplane separating $x_n$ and $p$ must separate $x_n$ and $x$, and finally $x_i$ and $x_{i+1}$ for some $0 \leq i \leq n-1$. Because $x_{i+1}$ is the projection of $x_i$ onto $S_i$, it follows from Lemma \ref{lem:projseparate} that $J$ separates $x_i$ from $S_i$. On the other hand, $\bigcap \mathcal{S}$ is included into $S_i$, so it follows that $J$ must be disjoint from $\bigcap \mathcal{S}$, which is impossible since it separates two of its vertices, namely $x_n$ and $p$. As a consequence, there do not exist hyperplanes separating $x_n$ and $p$, ie., $x_n=p$ is the projection of $x$ onto $\bigcap \mathcal{S}$.

\medskip \noindent
We are finally ready to prove our lemma. So let $J$ be a hyperplane separating $x$ and $x_n$. We know that $J$ must separate $x_i$ and $x_{i+1}$ for some $0 \leq i \leq n-1$. And once again because $x_{i+1}$ is the projection of $x_i$ onto $S_i$, we deduce from Lemma \ref{lem:projseparate} that $J$ separates $x_i$ from $S_i$. Therefore, $S_i \in \mathcal{S}$ is included into a sector delimited by $J$, which concludes the proof.
\end{proof}

\noindent
Our previous lemma allows us to prove the following characterisation of median triangles:

\begin{lemma}\label{lem:mediantriangleS}
Let $X$ be a quasi-median graph and $x,y,z \in X$ three vertices. Set $\mathcal{S}$ the set of the sectors $S$ satisfying $|S \cap \{ x,y,z \} | =2$. A triple $(x',y',z')$ is a median triangle of $(x,y,z)$ if and only if there exists some subcollection $\mathcal{S}_0 \subset \mathcal{S}$ such that $x',y',z'$ are respectively the projections of $x,y,z$ onto $\bigcap \mathcal{S}_0$.
\end{lemma}

\begin{proof}
Let  $\mathcal{S}_0 \subset \mathcal{S}$ be a subcollection and let $x',y',z'$ denote respectively the projections of $x,y,z$ onto $\bigcap \mathcal{S}_0$. We want to prove that $(x',y',z')$ is a median triangle of $(x,y,z)$. By symmetry, it is sufficient to prove that $d(x,y)=d(x,x')+d(x',y')+d(y',y)$. Equivalently, fixing some geodesics $[x,x']$, $[x',y']$, $[y',y]$ respectively between $x$ and $x'$, $x'$ and $y'$, $y'$ and $y$, we will show that the concatenation $[x,x'] \cup [x',y'] \cup [y',y]$ is a geodesic. 

\medskip \noindent
Let $J$ be a hyperplane intersecting $[x,x']$. It follows from Lemma \ref{lem:projsectors} that $J$ delimits a sector $S$ such that $S_0 \subset S$ for some $S_0 \in \mathcal{S}_0$. Because $S_0 \in \mathcal{S}_0$, we know that $y \in S_0$, hence $y \in S$. Moreover, $\bigcap \mathcal{S}_0 \subset S_0 \subset S$, so $x',y' \in S$. Thus, the sector $S$ delimited by $J$ contains $x',y',y$, and it follows by convexity of $S$ that $[x',y'] \cup [y',y] \subset S$. A fortiori, $J$ does not intersect $[x',y'] \cup [y',y]$. Similarly, one shows that a hyperplane intersecting $[y,y']$ cannot intersect $[x,x'] \cup [x',y']$. Therefore, $[x,x'] \cup [x',y'] \cup [y',y]$ is a geodesic.

\medskip \noindent
Conversely, suppose that $(x',y',z')$ is a median triangle of $(x,y,z)$. Let $\mathcal{S}_x$ (resp. $\mathcal{S}_y, \mathcal{S}_z$) denote the collection of the sectors which are delimited by the hyperplanes separating $x$ and $x'$ (resp. $y$ and $y'$, $z$ and $z'$) and which contain $x'$ (resp. $y'$, $z'$). Finally, set $\mathcal{S}_0 = \mathcal{S}_x \cup \mathcal{S}_y \cup \mathcal{S}_z$. We want to show that $x',y',z'$ are respectively the projections of $x,y,z$ onto $\bigcap \mathcal{S}_0$. By symmetry, it is sufficient to prove this statement for $x$. 

\medskip \noindent
First of all, fix some geodesics $[x, x'], [x', y'], [x', z'], [y', z'], [y, y'], [z, z']$. If $J$ is a hyperplane separating $x$ and $x'$, we know that, since a geodesic cannot intersect twice a hyperplane, necessarily $J$ must be disjoint from $[x', y']$ (resp. $[x', z']$, $[y',y]$, $[z',z]$), so that $y'$ (resp. $z'$, $y$, $z$) must belong to the same sector delimited by $J$ as $x'$. Otherwise saying, we have proved the following statement, which we record for future use:

\begin{fact}\label{fact:1}
Given a median triangle $(x',y',z')$ of a given triple $(x,y,z)$ and a hyperplane $J$ separating $x$ and $x'$, the vertices $x',y',z',y,z$ all belong to a common sector delimited by $J$. 
\end{fact}

\noindent
A similar statement holds for hyperplanes separating $y$ and $y'$, and for those separating $z$ and $z'$, for the same reason. Consequently, the vertices $x',y',z'$ all belong to $\bigcap \mathcal{S}_0$.

\medskip \noindent
Now, let $p$ denote the projection of $x$ onto $\bigcap \mathcal{S}_0$. Since $x' \in \bigcap \mathcal{S}_0$, there exists a geodesic between $x$ and $x'$ passing through $p$. If $J$ is a hyperplane separating $p$ and $x'$, then $J$ must separate $x$ and $x'$, so that $J$ delimits a sector $S$ which belongs to $\mathcal{S}_x \subset \mathcal{S}_0$ (namely, the sector which contains $x'$). But the inclusion $\bigcap \mathcal{S}_0 \subset S$ implies that $J$ must be disjoint from $\bigcap \mathcal{S}_0$, which is impossible since $J$ separates two vertices of this intersection, namely $x'$ and $p$. Thus, we have proved that no hyperplanes separate $x'$ and $p$, hence $x'=p$. Therefore, $x'$ is the projection of $x$ onto $\bigcap \mathcal{S}_0$, concluding the proof.
\end{proof}

\begin{proof}[Proof of Proposition \ref{prop:quasimedian}.]
Denote by $\mathcal{S}$ the set of the sectors $S$ satisfying $|S \cap \{ x,y,z \} | =2$. We claim that, if $x',y',z'$ denote respectively the projections of $x,y,z$ onto $\bigcap \mathcal{S}$, then $(x',y',z')$ is a quasi-median of $(x,y,z)$. First of all, we know that this is a median triangle according to Lemma \ref{lem:mediantriangleS}. 

\medskip \noindent
Now, we claim that a hyperplane $J$ separates $x'$ and $y'$ if and only if it separates $x,y,z$. Suppose first that $J$ separates $x'$ and $y'$. Because $(x',y',z')$ is a median triangle, this implies that $J$ separates $x$ and $y$. If $z$ belongs to the same sector $S$ delimited by $J$ as $x$ or $y$, then we deduce that $\bigcap \mathcal{S} \subset S$ cannot contain $x'$ or $y'$, which is absurd. Therefore, $J$ separates $x,y,z$. To prove the converse, we will show the following more general assertion:

\begin{fact}\label{fact:sizemedian}
Let $(x'',y'',z'')$ be any median triangle of $(x,y,z)$. If a hyperplane $J$ separates $x,y,z$, then it must separate $x''$ and $y''$.
\end{fact}

\noindent
According to Lemma \ref{lem:mediantriangleS}, there exists some $\mathcal{S}_0 \subset \mathcal{S}$ such that $x'',y'',z''$ are respectively the projections of $x,y,z$ onto $\bigcap \mathcal{S}_0$. Since it follows from Lemma \ref{lem:projsectors} that a hyperplane separating $x$ and $x''$ does not separate $y$ and $z$, we know that $J$ cannot sepate $x$ and $x''$. Similarly, it cannot separate $y$ and $y''$. On the other hand, we know that $J$ separates $x$ and $y$, so we deduce that $J$ must separate $x''$ and $y''$. This concludes the proof of our fact, and also of our claim. 

\medskip \noindent
Similarly, one shows that the hyperplanes separating $x'$ and $z'$, or $y'$ and $z'$, are precisely the hyperplanes separating $x,y,z$. Thus, we have proved that $(x',y',z')$ is an equilateral median triangle of $(x,y,z)$ of size the number of hyperplanes separating $x,y,z$, which we denote by $\mu(x,y,z)$. On the other hand, it follows from Fact \ref{fact:sizemedian} that the size of any median triangle is at least $\mu(x,y,z)$. Thus, we have proved that $(x',y',z')$ is a quasi-median of $(x,y,z)$.

\medskip \noindent
We record the following assertion which we have proved for future use:

\begin{fact}\label{fact:5}
For every triple $(x,y,z)$, its quasi-median $(x',y',z')$ satisfies the following property: a hyperplane separates $x'$ and $y'$ if and only if it separates $x',y',z'$ if and only if it separates $x,y,z$.
\end{fact}

\noindent
Now, we want to prove that it is the unique quasi-median of $(x,y,z)$. So let $(x'',y'',z'')$ be another equilateral median triangle of $(x,y,z)$. According to Lemma \ref{lem:mediantriangleS}, there exists some $\mathcal{S}_0 \subset \mathcal{S}$ such that $x'',y'',z''$ are respectively the projections of $x,y,z$ onto $\bigcap \mathcal{S}_0$. If $\bigcap \mathcal{S}_0 = \bigcap \mathcal{S}$, then $(x'',y'',z'')=(x',y',z')$, so suppose that $\bigcap \mathcal{S} \nsubseteq \bigcap \mathcal{S}_0$. By considering a hyperplane separating a vertex of $\bigcap \mathcal{S}_0 \backslash \bigcap \mathcal{S}$ from its projection onto $\mathcal{S}$, we find a hyperplane $J$ intersecting $\bigcap \mathcal{S}_0$ which is disjoint from $\bigcap \mathcal{S}$. By applying Lemma \ref{lem:projsectors}, it follows that $J$ delimits a sector $S$ which contains an element of $S' \in \mathcal{S}$; say that $y,z \in S'$ and $x \notin S'$. Let $J'$ denote the hyperplane delimiting the sector $S'$. Notice that $y'' \in \bigcap \mathcal{S}_0 \subset S'$, and similarly $z'' \in S'$. On the other hand, we claim that $x'' \notin S'$. Otherwise, $J'$ separates $x$ and $x''$ because $x \notin S'$, and it follows from Lemma \ref{lem:projseparate} that $J'$ separates $x$ from $\bigcap \mathcal{S}_0$. Therefore, it is sufficient to show that $J'$ intersects $\bigcap \mathcal{S}_0$. But $S' \subset S$ implies that
$$S^c \cap \bigcap \mathcal{S}_0 \subset (S')^c \bigcap \mathcal{S}_0,$$
and because $J$ intersects $\bigcap \mathcal{S}_0$, we know that the left-hand side is non empty. Our claim follows. Thus, we have proved that $J'$ separates $x$ and $\{y,z\}$, but does not separate $y$ and $z$. On the other hand, we know from Fact \ref{fact:sizemedian} that any hyperplane separating $x,y,z$ separates $x''$ and $y''$. Therefore, $d(x'',y'') \geq \mu(x,y,z)+1$. A fortiori, the size of our equilateral median triangle must be strictly larger than the size of $(x',y',z')$, so $(x'',y'',z'')$ is not a quasi-median of $(x,y,z)$. This proves that $(x',y',z')$ is the unique quasi-median of $(x,y,z)$.

\medskip \noindent
Finally, let $P$ denote the gated hull of $\{ x',y',z'\}$. As a consequence of Fact \ref{fact:sizemedian}, any hyperplane separating $x, y, z$ must separate $x',y',z'$. Since $d(x',y')=d(y',z')=d(z',x')=\mu(x,y,z)$, we deduce that no other hyperplane can separate two vertices of $\{x',y',z' \}$. Therefore, it follows from Proposition \ref{prop:gatedhullhyp} that the hyperplanes of $P$ correspond to the hyperplanes of $X$ separating $x, y, z$. In particular, they must be pairwise transverse, so that Lemma \ref{lem:prismhyp} implies that $P$ is a prism. This proves the last statement of our proposition.
\end{proof}

\begin{remark}
It can be proved that the convex hull of the quasi-median of any triple of vertices is a product of triangles. Indeed, it follows from the description of convex hulls given in Section \ref{section:convexhull} that the sectors of the convex hull of some set of vertices $S \subset X$, viewed as a quasi-median graph on its own right, correspond to the sectors of $X$ separating two vertices of $S$; as a consequence, the hyperplanes of the convex hull of our quasi-median delimits exactly three sectors, so that it follows from Lemma \ref{lem:prismhyp} that this convex hull is isomorphic to a product of triangles. This statement was originally proved by Mulder in his thesis \cite[(25) p. 149]{Mulder} (see also \cite[Proposition 3]{BandeltChepoi2}). 
\end{remark}

\begin{cor}\label{cor:whenmedian}
A graph is median if and only if it is quasi-median and triangle-free.
\end{cor}

\begin{proof}
Proving that a median graph is quasi-median and triangle-free is left as an exercice. Conversely, if $X$ is a quasi-median graph which is triangle-free, then any hyperplane delimits exacly two sectors. In particular, no hyperplane can separate three vertices at the same time. Therefore, the unique quasi-median given by Proposition \ref{prop:quasimedian} has size zero, ie., it produces a unique median vertex.
\end{proof}

\noindent
Next, let us prove the following application of Corollary \ref{cor:whenmedian}. Recall that, given a graph $X$ and two vertices $x,y \in X$, the \emph{interval}\index{Interval} $I(x,y)$ is defined as the union of all the geodesics between $x$ and $y$.

\begin{prop}\label{prop:intervalmedian}
In a quasi-median graph, an interval is median.
\end{prop}

\noindent
We begin by proving a technical lemma.

\begin{lemma}
Let $X$ be a quasi-median graph, $x,y \in X$ two vertices and $a,b \in I(x,y)$ two adjacent vertices. Then $d(x,a) \neq d(x,b)$.
\end{lemma}

\begin{proof}
Suppose by contradiction that $d(x,a)=d(x,b)=k$. We deduce from the triangle condition that there exists some vertex $p \in X$ adjacent to both $a$ and $b$ such that $d(x,p)=k-1$. Similarly, noticing that $$d(y,a)=d(x,y)-d(x,a)=d(x,y)-d(x,b)=d(y,b),$$ we know that there exists some vertex $q \in X$ adjacent to both $a$ and $b$ such that $d(y,q)=d(y,a)-1=d(x,y)-d(x,a)-1=d(x,y)-k-1$. Now,
$$d(x,p)+d(y,q)+2 = d(x,y) \leq d(x,p)+d(p,q)+d(q,y),$$
hence $d(p,q) \geq 2$. On the other hand $d(p,q) \leq d(p,a)+d(a,q)=2$, so $d(p,q)=2$. As a consequence, the vertices $a,b,p,q$ define an induced subgraph which is isomorphic to $K_4^-$, a contradiction.
\end{proof}

\noindent
The previous lemma has several interesting consequences on the structure of intervals.

\begin{cor}
In a quasi-median graph, an interval is an induced subgraph.
\end{cor}

\begin{proof}
Let $X$ be a quasi-median graph and $x,y \in X$ two vertices. If $a,b \in I(x,y)$ are two vertices adjacent in $X$, we want to prove that they are adjacent in $I(x,y)$. First, the previous lemma implies that $d(x,a) \neq d(x,b)$; say $d(x,a)=d(x,b)-1$. Let $\gamma$ denote the concatenation of the subsegment between $x$ and $a$ of some geodesic between $x$ and $y$ passing through $a$, followed by the edge $(a,b)$, and then followed by the subsegment between $b$ and $y$ of some geodesic between $x$ and $y$ passing through $b$. Now, notice that
$$d(x,y)= d(x,b)+d(b,y)=d(x,a)+1+d(b,y)=d(x,a)+d(a,b)+d(b,y),$$
which is precisely the length of $\gamma$. Therefore, $\gamma$ defines a geodesic between $x$ and $y$ containing the edge $(a,b)$. A fortiori, $(a,b) \subset I(x,y)$.
\end{proof}

\begin{cor}\label{cor:intervaltrianglefree}
In a quasi-median graph, an interval is triangle-free.
\end{cor}

\begin{proof}
Let $X$ be a quasi-median graph, $x,y \in X$ two vertices and $a,b,c \in I(x,y)$ three vertices such that $d(x,b) \leq d(x,a),d(x,c)$ and $b$ is adjacent to both $a$ and $c$. Because $a$ and $b$ are adjacent, and that $d(x,b) \leq d(x,a)$, we deduce that either $d(x,b)=d(x,a)$ or $d(x,b)=d(x,a)-1$; but the former case is impossible according the previous lemma, hence $d(x,b)=d(x,a)-1$. Similarly, we show that $d(x,c)=d(x,b)+1$, so that $d(x,c)=d(x,a)$. It follows from the previous lemma that $a$ and $c$ cannot be adjacent, so that the vertices $a,b,c$ cannot define a triangle.
\end{proof}

\noindent
The following lemma is the last step to prove Proposition \ref{prop:intervalmedian}.

\begin{lemma}
In a quasi-median graph, an interval is a convex subgraph.
\end{lemma}

\begin{proof}
Let $X$ be a quasi-median graph and $x,y \in X$ two vertices. We want to prove that, if $c$ is a vertex which belongs to some geodesic $[a,b]$ between two vertices $a,b \in I(x,y)$, then $c \in I(x,y)$, ie., there exists some geodesic between $x$ and $y$ passing through $c$. More precisely, if we fix two geodesics $[x,c]$ and $[c,y]$, we claim that the concatenation $[x,c] \cup [c,y]$ defines a geodesic. For convenience, we will denote by $[x,a]$ (resp. $[a,y]$) the subsegment between $x$ and $a$ (resp. $a$ and $y$) of some geodesic between $x$ and $y$ passing through $a$; and similarly for $[x,b]$ and $[b,y]$.

\medskip \noindent
Suppose by contradiction that $[x,c] \cup [c,y]$ is not a geodesic, so that there exists some hyperplane $J$ intersecting both $[x,c]$ and $[c,y]$. We will use the following observation:

\begin{fact}
Let $[p,q] \cup [q,r] \cup [r,p]$ be a geodesic triangle and $J$ a hyperplane intersecting $[p,q]$. Then $J$ has to intersect either $[q,r]$ or $[r,p]$.
\end{fact}

\noindent
Indeed, because $J$ intersects the geodesic $[p,q]$, necessarily $J$ separates $[p,q]$. Therefore, if the path $[q,r] \cup [r,p]$ does not meet $J$, it defines a path between $p$ and $q$ which lies in some sector delimited by $J$, which is impossible. Consequently, $J$ has to intersect either $[q,r]$ or $[r,p]$ (or both), proving our fact.

\medskip \noindent
Let $[a,c]$ (resp. $[b,c]$) denote the subsegment of $[a,b]$ between $a$ and $c$ (resp. $b$ and $c$). By applying our fact to the triangle $[a,x] \cup [x,c] \cup [c,a]$, we deduce that $J$ has to intersect either $[x,a]$ or $[a,c]$; and by applying our fact to the triangle $[a,y] \cup [y,c] \cup [c,a]$, we deduce that $J$ has to intersect either $[y,a]$ or $[a,c]$. On the other hand, because $[x,a] \cup [a,y]$ is a geodesic, $J$ intersects it at most once. Therefore, $J$ must intersect $[a,c]$. Similarly, we show that $J$ intersects $[b,c]$. However, $[a,c] \cup [c,b]=[a,b]$ is a geodesic, so $J$ cannot intersect it twice, hence a contradiction.
\end{proof}

\begin{proof}[Proof of Proposition \ref{prop:intervalmedian}.]
Let $X$ be a quasi-median graph and $x,y \in X$ two vertices. Because $I(x,y)$ is a convex subgraph of $X$, the interval $I(x,y)$ is a quasi-median graph on its own right. Moreover, we know thanks to Corollary \ref{cor:intervaltrianglefree} that $I(x,y)$ is triangle-free, so we conclude from Corollary \ref{cor:whenmedian} that $I(x,y)$ is a median graph.
\end{proof}

\begin{remark}
The uniqueness of quasi-medians stated in Proposition \ref{prop:quasimedian} was proved in \cite{quasimedian}, as well as the gatedness of intervals and Proposition \ref{prop:intervalmedian}; and Corollary \ref{cor:whenmedian} was proved in \cite{Mulder}. It is interesting to notice that our langage of hyperplanes turns out to define a common framework to prove all these results. 
\end{remark}

\subsection{Convex hulls}\label{section:convexhull}

\begin{definition}
Let $X$ be a graph and $S \subset X$ a set of vertices. The \emph{convex hull}\index{Convex hull} of $S$ is the smallest convex subgraph of $X$ containing $S$; alternatively, this is the intersection of all the convex subgraphs of $X$ which contain $S$.
\end{definition}

\noindent
In this section, we are interested in characterizing the convex hull of a given set of vertices. In CAT(0) cube complexes, the convex hull of some subspace coincides with the intersection of all the halfspaces containing it. However, according to our dictionnary given by Table \ref{tab:Dictionnary}, the translation of this statement is precisely Lemma \ref{lem:gatedhull}. So we need to introduce another family of specific subgraphs.

\begin{definition}\label{def:multisector}
Let $X$ be a quasi-median graph. A \emph{multisector}\index{Multisectors} $M$ is the subgraph generated by a union $D_1 \cup \cdots \cup D_n$, where $D_1 , \ldots, D_n$ are sectors delimited by a common hyperplane. In particular, if $D$ is a sector delimited by a hyperplane $J$, we refer to the multisector generated by the sectors delimited by $J$ which are different from $D$ as the \emph{cosector}\index{Cosectors} associated to $D$.
\end{definition}

\noindent
As a consequence of the following lemma, a proper multisector is never gated.

\begin{lemma}
Let $J$ be a hyperplane and $D_1, \ldots, D_k$ some sectors delimited by $J$, with $k \geq 2$. Let $M$ denote the multisector defined by $D_1, \ldots, D_k$. The gated hull of $M$ is $N(J) \cup D_1 \cup \cdots \cup D_k$.
\end{lemma}

\begin{proof}
Let $G$ denote the gated hull of $M$, so that we want to prove that 
$$G= N(J) \cup D_1 \cup \cdots \cup D_k.$$
Let $C$ be a clique dual to $J$. Because $k \geq 2$, $M$ contains an edge of $C$, so that $G$ must contain $C$ since any gated subgraph contains its triangles. Therefore, $N(J) \cup D_1 \cup \cdots \cup D_k \subset G$. To conclude, it is sufficient to notice that $N(J) \cup D_1 \cup \cdots \cup D_k$ is gated. This is proved by the following fact, which is slightly more general than what we need here since we allow $k=1$. 

\begin{fact}\label{fact:multisectorgatedhull}
Let $J$ be a hyperplane and $D_1, \ldots, D_k$ some sectors delimited by $J$ with $k \geq 1$. Then $N(J) \cup D_1 \cup \cdots \cup D_k$ is a gated subgraph.
\end{fact}

\noindent
For convenience, set $L= N(J) \cup D_1 \cup \cdots \cup D_k$. Let $x \in X \backslash L$ be a vertex. Because $J$ separates $x$ from $D_1, \ldots, D_k$, any geodesic between $x$ and some vertex of $D_1 \cup \cdots \cup D_k$ must intersect $N(J)$. Therefore, $d(x,L)=d(x,N(J))$. As a consequence, the gate of $x$ in $N(J)$ defines a gate of $x$ in $L$. We deduce that $L$ must be gated.
\end{proof}

\noindent
The main result of this section is the following statement.

\begin{prop}\label{prop:convexhull}
The convex hull of a set of vertices $S$ coincides with the intersection of all the multisectors containing $S$
\end{prop}

\noindent
Our proposition will be essentially a consequence of the following two preliminary lemmas.

\begin{lemma}\label{lem:multisectorconvex}
A multisector is convex.
\end{lemma}

\begin{proof}
Let $M$ be a multisector generated by some sectors $D_1, \ldots, D_k$ delimited by a hyperplane $J$. If $x,y \in M$ are two vertices, we want to prove that any geodesic $\gamma$ between $x$ and $y$ belongs to $M$. Let $1 \leq i,j \leq k$ be such that $x \in D_i$ and $y \in D_j$. If $i=j$, then $\gamma \subset D_i \subset M$ by convexity of $D_i$, and we are done. Now, suppose that $i \neq j$. In particular, $J$ separates $x$ and $y$, so that $\gamma$ must contain an edge $e$ dual to $J$. Let us write $\gamma$ as the concatenation $\gamma_1 \cup e  \cup \gamma_2$, where $\gamma_1$ is the initial segment of $\gamma$ whose last vertex is the initial point of $e$ and $\gamma_2$ the final segment of $\gamma$ starting from the terminating vertex of $e$. We know that, because $\gamma$ is a geodesic, $J$ intersects it once, hence $\gamma_1 \subset D_i$ and $\gamma_2 \subset D_j$; in particular, $\gamma_1, \gamma_2 \subset M$. Moreover, because the endpoints of $e$ belong to $M$, we necessarily have $e \subset M$. A fortiori, we have $\gamma \subset M$, concluding the proof.
\end{proof}

\begin{lemma}\label{lem:multisectorseparating}
Let $C$ be a convex subgraph and $x \notin C$ a vertex. There exists a sector $D$ satisfying $x \in D$ and $D \cap C = \emptyset$.
\end{lemma}

\begin{proof}
Let $y \in C$ be a vertex minimizing the distance ot $x$ in $C$. Let $J$ be a hyperplane separating $x$ and $y$ (such a hyperplane exists since $x \notin C$), and let $D$ denote the sector delimited by $J$ which contains $x$. We claim that $D$ is the sector we are looking for, ie., $D \cap C= \emptyset$. Fix some vertex $z \in C$, and let $(x',y',z')$ be the quasi-median of the triple $(x,y,z)$. Fix some geodesics $[z,z']$, $[y',z']$, $[y,y']$. Notice that, because $C$ is convex, the geodesic $[y,y'] \cup [y',z'] \cup [z',z]$ between $y$ and $z$ must be included into $C$. As a consequence, because $y$ minimizes the distance to $x$ in $C$, necessarily $y=y'$. In particular, the hyperplane $J$, which separates $x$ and $y$, either separates $x$ and $x'$ or it separates $x'$ and $y'$. In the former case, $J$ separates $x$ and $z$, hence $z \notin D$; in the latter case, we deduce from Fact \ref{fact:5} that $J$ separates $x,y,z$, hence $z \notin D$. We conclude that $D \cap C= \emptyset$.
\end{proof}

\begin{proof}[Proof of \ref{prop:convexhull}.]
Let $C$ denote the convex hull of $S$ and $I$ the intersection of all the multisectors containing $S$. We deduce from Lemma \ref{lem:multisectorconvex} that $I$ is convex, hence $C \subset I$. Now, if $x \notin C$, then we deduce from Lemma \ref{lem:multisectorseparating} that there exists some sector $D$ satisfying $x \in D$ and $D \cap C = \emptyset$. In particular, if $M$ denotes the cosector associated to $D$, then $S \subset C \subset M$ and $x \notin M$. A fortiori, $x \notin I$ since $I \subset M$. Thus, we have proved that $X \backslash C \subset X \backslash I$, or equivalently $I \subset C$. We conclude that $I=C$.
\end{proof}

\begin{cor}\label{cor:finiteconvexhull}
The convex hull of a finite subset is finite.
\end{cor}

\begin{proof}
Let $C$ denote the convex hull of a finite set of vertices $S$, and let $D$ be a sector. Three cases may happen. Either $D \cap S= \emptyset$, so that $S$ is included into the cosector associated to $D$, hence $D \cap C= \emptyset$; or $S \subset D$, so that $C \subset D$; or $D$ separates two vertices of $S$. Therefore, a sector separates two vertices of $C$ if and only if it separates two vertices of $S$. On the other hand, because $S$ is finite, there exist only finitely many such sectors, so that we deduce from Corollary \ref{cor:finitesubgraph} that $C$ is finite.
\end{proof}

\subsection{Flat rectangles}

\noindent
In this section, we introduce a particular class of subgraphs we call \emph{flat rectangles}. We used them in \cite{coningoff} to study hyperbolicity in CAT(0) cube complexes. Our goal here is to prove that similar results still hold in quasi-median graphs. In particular, Proposition \ref{prop:fourcycle} and Proposition \ref{prop:qmhyp} for CAT(0) cube complexes correspond respectively to \cite[Theorem 2.13 and Theorem 3.3]{coningoff}. The hyperbolicity of quasi-median graphs was also studied in \cite[Corollary 5]{CDEHV}.

\begin{definition}
A \emph{flat rectangle}\index{Flat rectangles} is an isometric embedding $R:[0,n] \times [0,m] \hookrightarrow X$, where we identify the square complex $[0,n] \times [0,m]$ with its 1-skeleton; if $m=n$, $R$ will be referred to as a \emph{flat square}. If $n,m \leq C$ for some $C \geq 0$, we say that $R$ is \emph{$C$-thin}.
\end{definition}

\noindent
The following lemma provides a useful way to construct flat rectangles.

\begin{definition}
A \emph{quadrangle}\index{Quadrangles} is a quadruple $(a,b,c,d)$ satisfying $b,d \in I(a,c)$ and $a,c \in I(b,d)$. 
\end{definition}

\begin{lemma}\label{lem:productingflatrectangle}
Let $(a,b,c,d)$ be a quadrangle and $[a,b],[b,c]$ two geodesics. There exists a flat rectangle $[0,r] \times [0,s] \hookrightarrow X$ satisyfing $[0,r] \times \{0 \} = [a,b]$, $\{r \} \times [0,s] =[b,c]$ and $(0,s)=d$.
\end{lemma}

\begin{proof}
According to Proposition \ref{prop:quasicubulatingqm}, we can identify $X$ with the quasi-cubulation of its canonical structure of space with partitions. In particular, the geodesic $[a,b]$ produces the sequence of vertices
$$\sigma_a, \ [\sigma_a,A_1], \ [\sigma_a,A_1,A_2], \ldots, [\sigma_a, A_1, \ldots, A_p]= \sigma_b;$$
and similarly, the geodesic $[b,c]$ produces
$$\sigma_b, \ [\sigma_b,B_1], \ [\sigma_b,B_1,B_2], \ldots, [\sigma_b,B_1, \ldots, B_q]= \sigma_c.$$
Let $J_1, \ldots, J_p$ (resp. $H_1, \ldots, H_q$) denote the hyperplanes underlying the sectors $A_1, \ldots, A_p$ (resp. $B_1, \ldots, B_q$). We claim that, for every $1 \leq i \leq p$ and $1 \leq j \leq q$, the hyperplanes $J_i$ and $H_j$ are transverse. For convenience, fix two geodesics $[a,d]$ and $[d,c]$; notice that, because $d \in I(a,c)$, the concatenation $[a,d] \cup [d,c]$ is a geodesic. 

\medskip \noindent
Because $J_i$ intersects a geodesic between $a$ and $c$, namely $[a,b] \cup [b,c]$, necessarily $J_i$ separates $a$ and $c$. In particular, $J_i$ must intersect the geodesic $[a,d] \cup [d,c]$. On the other hand, since $a \in I(b,d)$, the concatenation $[b,a] \cup [a,d]$ is a geodesic, so that $J_i$ cannot intersect it twice. We deduce that $J_i$ intersects $[c,d]$. Therefore, $b,c \in A_i$ and $a,d \notin A_i$. Similarly, we show that $a,b \in B_j$ and $c,d \notin B_j$. In particular, $J_i$ and $H_j$ must be transverse.

\medskip \noindent
As a consequence of Lemma \ref{lem:orientationcommute}, for every $1 \leq i \leq p$ and $1 \leq j \leq q$ we have
$$c= [\sigma_a,A_1, \ldots, A_p, B_1, \ldots, B_q]=[ \sigma_a,A_1 \ldots, A_i,B_1, \ldots, B_q, A_{i+1}, \ldots, A_p],$$
and $\sigma_{ij} = [\sigma,A_1, \ldots, A_i,B_1, \ldots, B_j]$ defines an orientation, ie., a vertex of $X$. Consider the application 
$$R : \left\{ \begin{array}{ccc} [0,p] \times [0,q] & \longrightarrow & X \\ (i,j) & \longmapsto & \sigma_{ij} \end{array} \right.$$
We claim that $R$ is the flat rectangle we are looking for. To conclude, we only have to prove that $R$ is an isometric embedding. Once again according to Lemma \ref{lem:orientationcommute}, we know that 
$$\sigma_{i+r,j+s} = [\sigma_{ij}, A_{i+1}, \ldots,  A_{i+r}, B_{j+1}, \ldots, B_{j+s}].$$
Moreover, we deduce from Lemma \ref{lem:orientationgeod} that
$$A_{i+k} \neq \sigma_a(J_k) = \sigma_{ij}(J_k) \ \text{and} \ B_{j+h} \neq \sigma_b(H_k)= \sigma_a(H_k)= \sigma_{ij}(H_k).$$
Therefore, according to Corollary \ref{cor:orientationdist}, we conclude that $d(\sigma_{ij},\sigma_{i+r,j+s})=r+s$.
\end{proof}

\noindent
Now, we are ready to prove the following proposition, which was fundamental in \cite{coningoff}. A sequence of subgraphs $(Y_1, \ldots, Y_n)$ is called a \emph{cycle of subgraphs} if $Y_i \cap Y_{i+1} \neq \emptyset$ for every $i$ mod $n$. 

\begin{prop}\label{prop:fourcycle}
Let $X$ be a quasi-median graph and $(Y_1,Y_2,Y_3,Y_4)$ a cycle of gated subgraphs. There exists a flat rectangle $[0,a] \times [0,b] \hookrightarrow X$ satisyfing $[0,a] \times \{ 0 \} \subset Y_1$, $[0,a] \times \{ b \} \subset Y_3$, $\{a \} \times [0,b] \subset Y_2$ and $\{0 \} \times [0,b] \subset Y_4$.
\end{prop}

\begin{proof}
Let $a \in Y_1 \cap Y_2$ be a vertex minimizing the distance to $Y_3 \cap Y_4$. Let $b$ (resp. $d,c$) denote its projection onto $Y_3$ (resp. $Y_4$, $Y_3 \cap Y_4$). It is worth noticing that, according to Lemma \ref{lem:projinter2}, we have $b \in Y_2$ and $d \in Y_1$. We claim that $(a,b,c,d)$ defines a quandrangle. For convenience, fix four geodesics $[a,b]$, $[b,c]$, $[c,d]$ and $[d,a]$.

\medskip \noindent
First, $c \in Y_3$ implies $b \in I(a,c)$, and $c \in Y_4$ implies $d \in I(a,c)$. 

\medskip \noindent
Now, notice that, if $J$ is a hyperplane intersecting $[c,d]$, because $c$ is the projection of $d$ onto $Y_3$ as a consequence of Lemma \ref{lem:projinter1}, we deduce from Lemma \ref{lem:projseparate} that $J$ has to be disjoint from $Y_3$. In particular, it cannot intersect $[b,c]$. Similarly, we show that no hyperplane intersecting $[b,c]$ can intersect $[c,d]$. Therefore, the concatenation $[d,c] \cup [c,b]$ is a geodesic according to Proposition \ref{prop:geodesichyp}. As a consequence, $c \in I(b,d)$.

\medskip \noindent
Finally, we want to prove that $a \in I(b,d)$. Because we know from Proposition \ref{prop:intervalmedian} that $I(a,c)$ is a median graph, the triple $(a,b,d)$ admits a median vertex $m$. Noticing that $m \in I(a,b) \subset Y_2$ since $Y_2$ is convex, and similarly $m \in I(a,d) \subset Y_1$ since $Y_1$ is convex, we deduce that $m \in Y_1 \cap Y_2$. On the other hand,
$$d(a,c)=d(a,d)+d(d,c)=d(a,m)+d(m,d)+d(d,c)=d(a,m)+d(m,c),$$
hence $d(m,c) =d(a,c)-d(a,m)$. Since $a$ minimizes the distance to $Y_3 \cap Y_4$ in $Y_1 \cap Y_2$, we deduce that $a=m$, hence $a=m \in I(b,d)$.

\medskip \noindent
Thus, we have proved that $(a,b,c,d)$ defines a quadrangle. According to Lemma \ref{lem:productingflatrectangle}, there exists a flat rectangle $[0,n] \times [0,m] \hookrightarrow X$ such that $(0,0)=a$, $(n,0)=b$, $(n,m)=c$ and $(0,m)=d$. Because $a,b \in Y_2$ and $Y_2$ is convex, necessarily $[0,n] \times \{0 \} \subset Y_2$. Similarly, we show that $\{ n \} \times [0,m] \subset Y_3$, $[0,n] \times \{ m \} \subset Y_4$ and $\{0 \} \times [0,m] \subset Y_1$. 
\end{proof}

\begin{cor}\label{cor:projflat}
Let $Y_1,Y_2 \subset X$ be two gated subgraphs. If $d$ denote the diameter of the projection of $Y_1$ onto $Y_2$, then there exists a flat rectangle $[0,d] \times [0,p] \hookrightarrow X$ satisfying $[0,d] \times \{ 0 \} \subset Y_1$ and $[0,d] \times \{p \} \subset Y_2$.
\end{cor}

\begin{proof}
Let $x,y \in Y_1$ be two vertices such that, if $x',y'$ denote their respective projections onto $Y_2$, then $d(x',y')=d$. Let $J_1, \ldots, J_d$ denote the $d$ hyperplanes separating $x'$ and $y'$. For every $1 \leq i \leq d$, we denote by $J_i^{+}$ (resp. $J_i^-$) the sector delimited by $J_i$ which contains $x'$ (resp. $y'$). Finally, let $C^{\pm}= \bigcap\limits_{i=1}^d J_i^{\pm}$. Notice that $x' \in C^+$ and $y' \in C^-$, so that $Y_2 \cap C^+$ and $Y_2 \cap C^-$ are non empty. Then, as a consequence of Lemma \ref{lem:projseparate}, any hyperplane separating $x$ and $x'$ must be disjoint from $Y_2$, so that no $J_i$ separates $x$ and $x'$, hence $x \in J_i^+$ and finally $x \in C^+$. Similarly, we show that $y \in C^-$. Therefore, the intersection $Y_1 \cap C^+$ and $Y_1 \cap C^-$ are non empty. We have prove that $(Y_1, C^+, Y_2,C^-)$ is a cycle of gated subgraphs, so that Proposition \ref{prop:fourcycle} implies that there exists a flat rectangle $[0,n] \times [0,m] \hookrightarrow X$ satisfying $[0,n] \times \{ 0 \} \subset Y_1$, $[0,n] \times \{m\} \subset Y_2$, $(0,m) \in C^+$ and $(n,m) \in C^-$. Since $J_1, \ldots, J_d$ clearly separate $C^+$ and $C^-$, we have $m \geq d$, so that taking a subrectangle of $[0,n] \times [0,m]$ if necessary produces the flat rectangle we are looking for.
\end{proof}

\noindent
As an application, we are able to characterize hyperbolic quasi-median graphs. First, define a \emph{facing triple} as the data of three hyperplanes such that no one separates the two others, and a \emph{join of hyperplanes} $(\mathcal{H}, \mathcal{V})$ as the data of two collections of hyperplanes $\mathcal{H}, \mathcal{V}$ which do not contain facing triples and such that any hyperplane of $\mathcal{H}$ is transverse to any hyperplane of $\mathcal{V}$; a join $(\mathcal{H}, \mathcal{V})$ is \emph{$K$-thin}, where $K \geq 0$ is some fixed integer, if $\min(\# \mathcal{H}, \# \mathcal{V}) \leq K$.

\begin{prop}\label{prop:qmhyp}
Let $X$ be a quasi-median graph. The following assertions are equivalent:
\begin{itemize}
	\item[(i)] $X$ is hyperbolic;
	\item[(ii)] its join of hyperplanes are uniformly thin;
	\item[(iii)] its flat rectangles are uniformly thin;
	\item[(iv)] its bigons are uniformly thin.
\end{itemize}
\end{prop}

\noindent
Let us begin by proving a preliminary lemma.

\begin{lemma}\label{lem:bigonflatrectangle}
Let $X$ be a quasi-median graph and $x,y \in X$ two vertices. For every vertex $z \in I(x,y)$ and for every geodesic $[x,y]$ between $x$ and $y$, there exists a flat square $R : [0,n] \times [0,n] \hookrightarrow X$ satisfying $z =(0,0)$ and $(n,n) \in [x,y]$.
\end{lemma}

\begin{proof}
Let $w$ be the unique vertex of $[x,y]$ satisfying $d(x,w)=d(x,z)$. Because $I(x,y)$ is median according to Proposition \ref{prop:intervalmedian}, we can introduce the median vertex $a \in X$ of the triple $(w,x,z)$ and the median vertex $b \in X$ of the triple $(w,y,z)$. Notice that $a,b \in I(z,w)$. Moreover, if we fix some geodesics $[x,a]$, $[a,z]$, $[z,b]$ and $[b,y]$, we notice that the concatenation
$$[x,a] \cup [a,z] \cup [z,b] \cup [b,y]$$
defines a geodesic; in particular, $[a,z] \cup [z,b]$ is a geodesic between $a$ and $b$ passing through $z$, hence $z \in I(a,b)$. Similarly, we show that $w \in I(a,b)$. Therefore, there exists a flat rectangle $R$ with $a,b,w,z$ as corners. By proving that
$$d(a,z)=d(a,w)=d(b,z)=d(b,w),$$
we will justify that this flat rectangle is the flat square we are looking for. First, notice that
$$d(x,a)+d(a,z)=d(x,z)=d(x,w)=d(x,a)+d(a,w),$$
hence $d(a,z)=d(a,w)$. Let $\ell_1$ denote this common value. Similarly, because
$$d(y,z)=d(x,y)-d(x,z)=d(x,y)-d(x,w)=d(y,w),$$
we show that $d(b,z)=d(b,w)$. Let $\ell_2$ denote this common value. Next, notice that
$$2 \ell_1= d(z,a)+d(a,w)=d(z,b)+d(b,w)=2 \ell_2,$$
hence $\ell_1= \ell_2$. This concludes the proof.
\end{proof}

\begin{proof}[Proof of Proposition \ref{prop:qmhyp}.]
Suppose that $X$ is $\delta$-hyperbolic for some $\delta \geq 0$ and let $(\mathcal{H}, \mathcal{V})$ be a join of hyperplanes. First notice that the cubical dimension of $X$ must be finite. If it was not the case, $X$ would contain (the 1-skeleton of) an $n$-cube as an isometrically embedded subgraph for every $n \geq 1$ according to Proposition \ref{prop:transversehypcube} and Lemma \ref{lem:prismgated}, but such a cube contains a triangle which is not $(n-1)$-thin. Next, suppose that $\# \mathcal{H}, \# \mathcal{V} \geq \mathrm{Ram}(d)$ for some $d \geq \dim_{\square}(X)$, where $\mathrm{Ram}(\cdot)$ denotes the corresponding Ramsey number (which is defined by the following property: for every $r \geq 1$ and every edge-coloring of the complete graph on at least $\mathrm{Ram}(r)$ vertices with two colors, one can find a monochromatic complete subgraph with at least $r$ vertices). This implies that $\mathcal{H}$ and $\mathcal{V}$ contain respectively some subcollections $\mathcal{H}_0$ and $\mathcal{V}_0$ of $d$ pairwise non transverse hyperplanes (this is a classical argument; see for instance \cite[Lemma 3.7]{coningoff}). Because $\mathcal{H}_0$ contains no facing triple, we can number its elements
$$\mathcal{H}_0= \{ H_1, \ldots, H_d \}$$
so that, for every $2 \leq i \leq d-1$, the hyperplane $H_i$ separates $H_{i-1}$ and $H_{i+1}$. Number similarly the elements of $\mathcal{V}_0$:
$$\mathcal{V}_0= \{ V_1, \ldots, V_d \},$$
i.e., $V_i$ separates $V_{i-1}$ and $V_{i+1}$ for every $2 \leq i \leq d-1$. Now, $(N(V_1),N(H_1),N(V_d),N(H_d))$ defines a cycle of gated subgraphs, so Proposition \ref{prop:fourcycle} yields a flat rectangle $R : [0,n] \times [0,m] \hookrightarrow X$ such that $[0,n] \times \{ 0 \} \subset N(H_1)$, $\{n \} \times [0,m] \subset N(V_d)$, $[0,n] \times \{ m \} \subset N(H_d)$ and $\{0 \} \times [0,m] \subset N(V_1)$. Since $V_1$ and $V_d$ are separated by $d-2$ hyperplanes, it follows that $n \geq d-2$; similarly, we know that $m \geq d-2$. On the other hand, because a flat square $[0,k] \times [0,k] \hookrightarrow X$ contains a geodesic triangle which is not $(k-1)$-thin, we deduce that the flat rectangles of $X$ are all $(\delta+1)$-thin. Therefore, $d \leq \min(n,m)+2 \leq \delta+3.$ Consequently, we have proved that
$$\min ( \# \mathcal{H}, \# \mathcal{V}) \leq \max \left( \mathrm{Ram}(\dim_{\square}(X)), \mathrm{Ram}(\delta+4) \right).$$
This proves the implication $(i) \Rightarrow (ii)$. 

\medskip \noindent
Next, because the hyperplanes dual to a flat rectangle $R : [0,n] \times [0,m] \hookrightarrow X$ defines a join of hyperplanes $(\mathcal{H}, \mathcal{V})$ satisfying $\# \mathcal{H} = m$ and $\# \mathcal{V}=n$, we know that the implication $(ii) \Rightarrow (iii)$ holds.

\medskip \noindent
Now, suppose that $(iii)$ holds, ie., there exists some constant $B$ such that the flat rectangles of $X$ are all $B$-thin. So, if we fix some bigon $(\gamma_1,\gamma_2)$ and some vertex $x \in \gamma_1$, Lemma \ref{lem:bigonflatrectangle} implies that there exists a flat square $R$ such that $x \in R$ and $R \cap \gamma_2$ contains some vertex $y$, and we conclude that
$$d(x, \gamma_2) \leq d(x,y) \leq \mathrm{diam}(R) \leq 2B.$$
This proves the implication $(iii) \Rightarrow (iv)$. Finally, suppose that $(iv)$ holds, ie., there exists some constant $B$ such that any bigon of $X$ is $B$-thin. Let $[x,y] \cup [y,z] \cup [z,x]$ be a geodesic triangle and let $(x',y',z')$ denote the quasi-median of the triple $(x,y,z)$. We fix some geodesics $[x,x']$, $[y,y']$, $[z,z']$, $[x',y']$, $[y',z']$, $[z',x']$. We want to prove that, if $p \in [x,y]$ is a vertex, then $d(p,[y,z] \cup [z,x]) \leq 4B+3$. Because the bigon defined by $[x,y]$ and $[x,x'] \cup [x',y'] \cup [y',y]$ is $B$-thin by assumption, there exists some $q \in [x,x'] \cup [x',y'] \cup [y',y]$ such that $d(p,q) \leq B$. If $q \in [x,x']$ (resp. $q \in [y,y']$), then similarly we find some $q' \in [x,z]$ (resp. $q' \in [y,z]$) such that $d(q,q') \leq B$, hence 
$$d(p, [x,x'] \cup [x',y'] \cup [y',y]) \leq d(p,q') \leq d(p,q)+d(q,q') \leq 2B.$$
From now on, suppose that $q \in [x',y']$. Fix some vertex $q' \in [x',z']$. Again by the same argument, we find some vertex $q'' \in [x,z]$ satisfying $d(q',q'') \leq B$. According to Proposition \ref{prop:quasimedian}, the gated hull of $\{x',y',z'\}$ is a prism $P$, so
$$d(p, [x,x'] \cup [x',y'] \cup [y',y]) \leq d(p,q'') \leq d(p,q)+d(q,q')+d(q',q'') \leq 2B+ \mathrm{diam}(P).$$
On the other hand, a prism of cubical dimension $2n$ contains a bigon which is not $(n-1)$-thin, hence $\dim_{\square} (P) \leq 2B+1$; as a conquence, $\mathrm{diam}(P) \leq 2B+3$. Therefore, 
$$d(p, [x,x'] \cup [x',y'] \cup [y',y]) \leq 2B+2B+3=4B+3.$$
Thus, we have proved that the geodesic triangles of $X$ are all $(4B+3)$-thin, so that $X$ must be hyperbolic. This proves $(iv) \Rightarrow (i)$.
\end{proof}

\subsection{Fixed point theorem}

\noindent
This section is dedicated to the proof of the following result.

\begin{thm}\label{thm:fixedpoint}
A group acting on a quasi-median graph with a bounded orbit stabilises a finite prism.
\end{thm}

\noindent
In the particular case of a finite group acting on a locally finite quasi-median graph, the previous statement follows from \cite[Theorem 4]{bucolic}. Our proof follows the argument given by Roller in \cite[Theorem 11.7]{Roller}, where he shows a similar statement for median algebras (which holds in particular for CAT(0) cube complexes). We begin by proving our theorem for finite quasi-median graphs. 

\begin{prop}\label{prop:fixedpoint}
Let $X$ be a finite quasi-median graph. There exists a prism $P \subset X$ which is stabilised by $\mathrm{Aut}(X)$.
\end{prop}

\begin{proof}
Let $\mathfrak{D}$ denote the sectors $D$ of $X$ satisfying $|D|> \frac{1}{2} |X|$ (by abuse of notation, we denote by $| \cdot|$ the number of vertices of the subgraph which we consider), and let $X^{(1)}$ denote the (finite) intersection $\bigcap\limits_{D \in \mathfrak{D}} D$. Notice that, if $D_1,D_2 \in \mathfrak{D}$, then $D_1 \cap D_2 \neq \emptyset$, since otherwise we would have
$$|X| \geq |D_1|+|D_2| > \frac{1}{2}|X| + \frac{1}{2} |X| = |X|,$$
a contradiction. Therefore, according to Helly's property \ref{prop:Helly}, $X^{(1)}$ is a non empty gated subgraph of $X$, which is clearly $\mathrm{Aut}(X)$-invariant. If $X^{(1)}$ is a single vertex, we are done. Otherwise, as a quasi-median graph on its own right, $X^{(1)}$ must contain some hyperplanes. In particular, if there exists a sector $D$ of $X^{(1)}$ satisfying $|D| > \frac{1}{2} |X^{(1)}|$, then the previous construction can be iterated, and we obtain a new non empty $\mathrm{Aut}(X)$-invariant gated subgraph $X^{(2)}$, and so on. Because $X$ is finite, our sequence $X \supsetneq X^{(1)} \supsetneq X^{(2)} \supsetneq \cdots$ must terminate to some non empty $\mathrm{Aut}(X)$-invariant gated subgraph $Y$. By construction, either $Y$ is a single vertex, so that we are done, or any sector $D$ of $Y$ satisfies $|D| \leq \frac{1}{2} |Y|$. 

\medskip \noindent
In the latter case, we claim that any two hyperplanes of $Y$ are transverse. If not, there exist two nested hyperplanes $\mathcal{P}$ and $\mathcal{Q}$. Let $P$ denote the sector delimited by $\mathcal{P}$ containing $\mathcal{Q}$, and $P'$ the union of the sectors delimited by $\mathcal{P}$ which are different from $P$; similarly, let $Q$ denote the sector delimited by $\mathcal{Q}$ containing $\mathcal{P}$, and $Q'$ the union of the sectors delimited by $\mathcal{Q}$ which are different from $Q$. Notice that $P' \subset Q$ implies $|P'| \leq |Q| \leq \frac{1}{2}|Y|$, hence
$$\frac{1}{2} |Y| \geq |P| = |Y|- |P'| \geq |Y|- \frac{1}{2}|Y|= \frac{1}{2}|Y|.$$
Therefore, $|P|= \frac{1}{2}|Y|$, and finally $|P'|= |Y|-|P|= \frac{1}{2}|Y|$. Similarly, we show that $|Q'|= \frac{1}{2}|Y|$. As a consequence,
$$|Y \backslash (P' \cup Q')| = |Y|- |P'|-|Q'|= |Y| - \frac{1}{2}|Y|  - \frac{1}{2}|Y| =0,$$
which is impossible. Thus, we have proved that the hyperplanes of $Y$ are pairwise transverse. It follows from Lemma \ref{lem:prismhyp} that $Y$ must be a prism.
\end{proof}

\begin{definition}
Let $X$ be a quasi-median graph. We define the \emph{gated topology} on $X^{(0)}$ by taking the set of the sectors of $X$ as a prebasis.
\end{definition}

\noindent
Notice that a sector is open by definition, so that, because two distinct vertices are always separated by at least one hyperplane, the topology is Hausdorff. But a sector is also closed, since it is the complement of the union of the other sectors delimited by the corresponding hyperplane. In particular, the gated topology is totally disconnected and any gated subgraph is closed. 

\begin{lemma}\label{lem:compactball}
With respect to the gated topology, a ball is compact.
\end{lemma}

\noindent
In the sequel, we will use the following notation: if $S$ is a set and $\mathcal{U}$ a collection of subsets, the intersection $\bigcap \mathcal{U}$ denotes $\bigcap\limits_{U \in \mathcal{U}} U$. 

\begin{proof}[Proof of Lemma \ref{lem:compactball}.]
Let $B$ be a ball of radius $r$. Recall that a Hausdorff topological space is compact if and only if, for every collection $\mathcal{F}$ of closed subspaces satisfying $\bigcap \mathcal{F}_s \neq \emptyset$ for every finite subcollection $\mathcal{F}_s \subset \mathcal{F}$, we have $\bigcap \mathcal{F} \neq \emptyset$. Therefore, it is sufficient to prove that, given a collection $\mathcal{H}$ of sectors satisfying $\bigcap \mathcal{H}_s \cap B \neq \emptyset$ for every finite subcollection $\mathcal{H}_s \subset \mathcal{H}$, we have $\bigcap \mathcal{H} \cap B \neq \emptyset$.

\medskip \noindent
If $y \in B$, let $\mathcal{H}_y$ denote $\{H \in \mathcal{H} \mid y \notin H \}$. We claim that $\# \mathcal{H}_y \leq 2r$. Let $H_1, \ldots, H_n \in \mathcal{H}_y$. By our hypothesis, there exists a vertex $z \in B \cap H_1 \cap \cdots \cap H_n$. Clearly, $H_1, \ldots, H_n$ separate $y$ and $z$, hence $\# \mathcal{H}_y \leq d(y,z) \leq 2r$.

\medskip \noindent
Next, we want to prove that $\bigcap \mathcal{H} \neq \emptyset$. Let $x \in B$ be the center of $B$. Since our previous claim states that $\mathcal{H}_x$ is finite, we know that there exists a vertex $y \in B \cap \bigcap \mathcal{H}_x$; notice that any sector of $\mathcal{H}$ contains either $x$ or $y$. Let $\mathcal{H}_0 = \mathcal{H} \backslash ( \mathcal{H}_x \cup \mathcal{H}_y)$ denote the subcollection of the sectors of $\mathcal{H}$ containing both $x$ and $y$. Notice that $x \in \bigcap \mathcal{H}_0 \cap \bigcap \mathcal{H}_y$ and $y \in \bigcap \mathcal{H}_0 \cap \bigcap \mathcal{H}_x$; moreover, because $\mathcal{H}_x$ and $\mathcal{H}_y$ are finite, our hypothesis implies that $B \cap \bigcap \mathcal{H}_x \cap \bigcap \mathcal{H}_y$ is non empty, so that in particular $\bigcap \mathcal{H}_x \cap \bigcap \mathcal{H}_y$ is necessarily non empty. Applying Helly's property \ref{prop:Helly}, we deduce that
$$\bigcap \mathcal{H}=  \bigcap \mathcal{H}_0 \cap \bigcap \mathcal{H}_x \cap \bigcap \mathcal{H}_y \neq \emptyset.$$
Let $z \in \bigcap \mathcal{H}$, and let $(x',y',z')$ denote the quasi-median of the triple $(x,y,z)$. Notice that any sector $H$ of $\mathcal{H}$ contains $z$ and either $x$ or $y$; a fortiori, it must contain either $I(z,x)$ or $I(z,y)$, hence $z' \in H$. Moreover
$$d(x,z')=d(x,x')+d(x',z')=d(x,x')+d(x',y')=d(x,y') \leq d(x,y) \leq r,$$
hence $z' \in B$. Thus, the intersection $\bigcap \mathcal{H} \cap B$ is non empty.
\end{proof}

\begin{lemma}\label{lem:gatedtopology2}
With respect to the gated topology, the set of accumulation points $B(x,r)'$ of a ball $B(x,r)$ centered at $x$ of radius $r$ is included into $B(x,r-1)$. 
\end{lemma}

\begin{proof}
Because $B(x,r)$ is compact, we know that $B(x,r)' \subset B(x,r)$. Therefore, we only have to prove that a vertex $y \in X$ satisfying $d(x,y)=r$ cannot belong to $B(x,r)'$. 

\medskip \noindent
Let $S$ denote the intersection of all the sectors of $X$ which contain $y$ but not $x$. Notice that, because there exist only finitely many sectors separating $x$ and $y$, $S$ is the intersection of finitely many sectors, so that it must be open with respect to the gated topology. Fix some vertex $z \in B(x,r) \cap S$. Notice that no hyperplane separate $x,y,z$ at the same time: otherwise, there would exist a sector containing $y$ but not $x$ and $z$, contradicting our choice of $z$. Therefore, we deduce from Proposition \ref{prop:quasimedian} that the triple $\{x,y,z \}$ admits a median vertex $m$. Similarly, no hyperplane can separate $m$ and $y$, since otherwise there would exist a sector containing $y$ but not $x$ and $z$, hence $m=y$. In particular, $y$ belongs to a geodesic between $x$ and $z$, hence
$$d(x,z)=d(x,y)+d(y,z)=r+d(y,z).$$
On the other hand, $z \in B(x,r)$ implies $d(x,z) \leq r$, hence $y=z$. Thus, we have proved that $B(x,r) \cap S = \{ y \}$. As a consequence, $\{y \}$ is an open subspace of $B(x,r)$, which implies that $y \notin B(x,r)'$. 
\end{proof}

\begin{proof}[Proof of Theorem \ref{thm:fixedpoint}]
Let $G$ be a group acting on a quasi-median graph $X$ with a bounded orbit. Let $B \subset X$ be a bounded $G$-invariant subset, say a bounded orbit of $G$. In particular, there exist some vertex $x \in X$ and some integer $r \geq 0$ such that $B \subset B(x,r)$. We define the sequence of subsets $(B^{(i)})$ by $B^{(0)}=B$ and $B^{(i+1)}= \left( B^{(i)} \right)'$ for every $i \geq 0$. Notice that $B' \subset B(x,r)' \subset B(x,r-1)$ according to Lemma \ref{lem:gatedtopology2}, and we deduce by induction that $B^{(i)} \subset B(x,r-i)$ for every $i \geq 0$. As a consequence, there exists some $i \geq 0$ such that $B^{(i+1)}= \emptyset$. This implies that $B^{(i)}$ must be finite, since otherwise $B^{(i+1)}$ would be non empty by compactness. Thus, $B^{(i)}$ is a finite $G$-invariant subset. Let $Y$ denote its convex hull. Clearly, $Y$ is $G$-invariant, and we deduce from Corollary \ref{cor:finiteconvexhull} that $Y$ is finite as well. Because $Y$ is a quasi-median graph on its own right, we deduce from Proposition \ref{prop:fixedpoint} that $Y$ contains a $G$-invariant prism, which concludes the proof. 
\end{proof}

\subsection{CAT(0)-ness}\label{section:qmCAT0}

\noindent
In this section, we prove that quasi-median graphs can be endowed naturally with a cellular structure making them CAT(0) spaces. See \cite[Appendix]{Leary} and references therein for the link between CAT(0) cube complexes and CAT(0) spaces; and \cite{RetractsChordal} to see how CAT(0) spaces arise similarly from the more general class of (finite) retracts of products of chordal graphs. 

\medskip \noindent
Given a complete graph $K$, consider the Hilbert space $\mathcal{H}(K)$ of the $\ell^2$-functions $K^{(0)} \to [0,1]$. If $\delta_v$, for $v \in K^{(0)}$, denotes the function $u \mapsto \left\{ \begin{array}{cl} 0 & \text{if} \ u \neq v \\ 1 & \text{otherwise} \end{array} \right.$, then the map $v \mapsto \frac{1}{\sqrt{2}} \delta_v$ defines an isometric embedding $K \hookrightarrow \mathcal{H}(K)$. By identifying $K$ with the convex hull of its image in $\mathcal{H}(X)$, we naturally view $K$ as a CAT(0) cellular complex. By Cartesian product, we produce a similar structure on any prism.

\medskip \noindent
A \emph{prism complex} $X$ is a cellular complex obtained by identifying a collection of prisms along their faces by isometries (a \emph{face} of a prism is a prism of smaller dimension). Given two points $x,y \in X$, a chain $\Sigma$ between $x$ and $y$ is a sequence of points $x_1, \ldots, x_n$ such that $x_1=x$, $x_n=y$, and for every $1 \leq i \leq n-1$, $x_i$ and $x_{i+1}$ belong to a same prism $P_i$. (Notice that a chain between two points of $X$ always exists if $X$ is connected: for instance, consider a path in the one-skeleton between two prisms containing our points.) The \emph{length} of our chain is defined by $\ell(\Sigma)= \sum\limits_{i=1}^{n-1} d_{P_i}(x_i,x_{i+1})$. Finally, we introduce a pseudo-distance on $X$ by
$$d : (x,y) \mapsto \inf \{ \ell (\Sigma) \mid \Sigma \ \text{chain between $x$ and $y$} \}.$$
This is the standard way to define a metric on a polyhedral complex, as explained in \cite[Definition 7.4]{MR1744486}. From now on, we will suppose that every prism complex is endowed with this pseudo-metric. By defining a \emph{quasi-median complex} as a prism complex obtained from a quasi-median graph by filling in every clique with a simplex and every 1-skeleton of an $n$-cube with an $n$-cube, we want to prove:

\begin{thm}\label{thm:CAT0}
Quasi-median complexes are CAT(0).
\end{thm}

\noindent
The following lemma will be fundamental in the proof of Theorem \ref{thm:CAT0}.

\begin{lemma}\label{lem:shortening}
Let $X$ be a quasi-median complex and $x,y \in X$ two points. Let $p_x,p_y$ denote the minimal prisms containing $x,y$ respectively. For every chain $\Sigma$ between $x$ and $y$, there exists another chain $\Sigma^*$ between $x$ and $y$ such that $\ell(\Sigma^*) \leq \ell(\Sigma)$ which is included into the gated hull of $\{p_x,p_y \}$.
\end{lemma}

\begin{proof}
Let $H$ denote the gated hull of $\{p_x,p_y\}$ and $\Sigma = (x_1, \ldots, x_n)$. A hyperplane $J$ is \emph{bad relatively to $\Sigma$} if $J$ does not separate $x$ and $y$, and $\Sigma$ is not contained into the sector $J^+$ delimited by $J$ which contains $x$ and $y$. Let $b(\Sigma)$ denote the number of hyperplanes which are bad relatively to $\Sigma$. We want to prove that there exists a chain $\Sigma^*$ between $x$ and $y$ satisfying $b(\Sigma^*)=0$ and $\ell(\Sigma^*) \leq \ell(\Sigma)$. If $b(\Sigma)=0$, it suffices to set $\Sigma^*= \Sigma$. From now on, suppose that $b(\Sigma) \geq 1$, and let $J$ be a bad hyperplane relatively to $\Sigma$ which maximizes the distance to $H$. Let $\partial = N(J) \cap J^+ \subset \partial J$. 

\medskip \noindent
Let $(x_r, \ldots, x_s)$ be a maximal subsegment of $\Sigma$ included into $N(J)$. We claim that the points $x_r,x_s$ belong to $\partial$. 

\medskip \noindent
Suppose by contradiction that $x_r \notin \partial$. Let $P$ be a prism containing $x_{r-1}$ and $x_r$. Because $x_{r-1} \notin N(J)$, $P \nsubseteq N(J)$; but $P \cap N(J) \neq \emptyset$ because $\Sigma$ is a chain. Let $v \in P \cap N(J)$ be a vertex, $J_1, \ldots, J_m$ the hyperplanes dual to $P$, and, for every $1 \leq i \leq m$, $v_i \in P$ a vertex adjacent to $v$ and separated from it by $J_i$. If $J,J_1, \ldots, J_m$ are pairwise transverse, we deduce from Fact \ref{fact:spanningprism} that the gated hull of $\{x,x_1, \ldots, x_m \}$ defines a prism whose dual hyperplanes are $J,J_1, \ldots, J_m$, hence $P \subset N(J)$, a contradiction. Therefore, there exists a hyperplane $J'$ dual to $P$ which is disjoint from $J$. On the other hand, because $x_r \notin \partial$ and $P \nsubseteq N(J)$, the prism $P$, and a fortiori the hyperplane $J'$, must be included into a sector delimited by $J$ which is different from $J^+$, the sector containing $H$. Consequently, we have $d(J',H)>d(J,H)$, contradicting our choice of $J$. Thus, we have proved that necessarily $x_r \in \partial$. Similarly, we show that $x_s \in \partial$.

\medskip \noindent
From Lemma \ref{lem:carrierproduct}, we know that $N(J)$ is naturally isomorphic to the product $F(J) \times C$, where $C$ is a clique dual to $J$ and $F(J)$ the main fiber of $J$; moreover, there exists $c \in C$ such that $\partial = F(J) \times \{c \}$. From $\Sigma$, we construct a new chain $\Sigma'$ by replacing the subsegment $(x_r,\ldots, x_s)$ with its projection onto $\partial$. It is worth noticing that any hyperplane intersecting $\Sigma'$ must intersect $\Sigma$, hence $b(\Sigma') \leq b(\Sigma)$; moreover, $\ell(\Sigma') \leq \ell(\Sigma)$. Finally, if we construct the chain $\Sigma''$ from $\Sigma$ by iterating this process with all the maximal subsegments of $\Sigma$ included into $N(J)$, then $J$ will not be a bad hyperplane relatively to $\Sigma''$, hence $b(\Sigma'') \leq b(\Sigma)-1$. Moreover, $\ell (\Sigma'') \leq \ell(\Sigma)$. 

\medskip \noindent
By iterating this construction, we find a new chain $\Sigma^*$ between $x$ and $y$ such that $b(\Sigma^*)=0$ and $\ell(\Sigma^*) \leq \ell(\Sigma)$. The fact that $b(\Sigma^*)=0$ precisely means that, for every hyperplane $J$ which does not separate $x$ and $y$, $\Sigma^*$ must be included into the sector delimited by $J$ which contains $x$ and $y$. From Lemma \ref{lem:gatedhull}, we deduce that $\Sigma^* \subset H$. 
\end{proof}

\noindent
The next step is to prove Theorem \ref{thm:CAT0} in a special case.

\begin{lemma}\label{lem:finiteCAT0}
Cubically finite quasi-median complexes are complete CAT(0) spaces.
\end{lemma}

\begin{proof}
We argue by induction on the number of hyperplanes. A quasi-median graph without hyperplanes being a single vertex, there is nothing to prove in this case. If the hyperplanes of our quasi-median complex are pairwise transverse, we deduce from Lemma \ref{lem:prismhyp} that it is a prism, so that it must be a complete CAT(0) space as a Cartesian product of complete CAT(0) spaces. From now on, we suppose that our quasi-median complex $X$ contains at least two disjoint hyperplanes $J,J'$. Let $J^+$ denote the sector delimited by $J$ which contains our second hyperplane $J'$. 

\medskip \noindent
Let $Z$ denote the union of $N(J)$ with all the sectors delimited by $J$ which are different from $J^+$, and $\partial$ the connected component of $\partial J$ equal to $N(J) \cap J^+= Z \cap J^+$. Notice that

\begin{claim}\label{claim:gluingdist}
For every vertices $x \in Z$ and $y \in J^+$, 
$$d(x,y)= \inf\limits_{z \in \partial} \left( d(x,z)+d(z,y) \right).$$
\end{claim}

\noindent
Let $\Sigma=(x_1, \ldots, x_n)$ be a chain from $x$ to $y$. Let $x_k$ denote the last element of the sequence $x_1, \ldots, x_n$ which belongs to $Z$. Notice that, since $x_n=y \notin Z$, necessarily $k \leq n-1$, so that $x_{k+1}$ is well-defined, and it belongs to $J^+$. Let $P$ denote the minimal prism containing both $x_{k}$ and $x_{k+1}$. Because $x_{k+1}$ does not belong to $Z$, and that $Z$ is a gated subgraph, there must exist a hyperplane $J'$ separating $x_{k+1}$ from $Z$; a fortiori, $J$ and $J'$ are disjoint. On the other hand, if $x_k$ does not belong to $\partial$, $J$ must separate $x_k$ and $x_{k+1}$ since they belong to distinct sectors delimited by $J$. Thus, the two disjoint hyperplanes $J$ and $J'$ are transverse to the same prism $P$, which is impossible. Therefore, $x_k \in \partial$. We deduce that, if $\Sigma_1=(x_1, \ldots, x_k)$ and $\Sigma_2=(x_k,\ldots, x_n)$ are subchains of $\Sigma$, then
$$\ell(\Sigma)= \ell(\Sigma_1)+ \ell(\Sigma_2) \geq d(x,x_k)+d(x_k,y) \geq \inf\limits_{z \in \partial} \left( d(x,z)+d(z,y) \right).$$
A fortiori, $d(x,y) \geq  \inf\limits_{z \in \partial} \left( d(x,z)+d(z,y) \right)$. The reverse inequality follows from the triangle inequality, which proves our claim.

\medskip \noindent
Next, noticing that $J$ is disjoint from $J^+$ and that $J'$ is disjoint from $Z$, we deduce that the numbers of hyperplanes of $Z$ and $J^+$ are strictly less than the number of hyperplanes of $X$. Notice also that $Z$ is a gated subcomplex of $X$ according to Fact \ref{fact:multisectorgatedhull}, so that it is a quasi-median complex on its own right. Therefore, we can apply our induction hypothesis to deduce that $J^+$ and $Z$ are CAT(0). In order to conclude that $X$ is CAT(0), it is sufficient to show that $\partial$ is a complete convex subspace in both $Z$ and $J^+$, so that the Gluing Theorem \cite[Theorem II.11.1]{MR1744486} applies. (Indeed, Claim \ref{claim:gluingdist} states that the pseudo-distance $d$ coincides with the pseudo-distance defined on the gluing $Z \sqcup_{\partial} J^+$.)

\medskip \noindent
From Lemma \ref{lem:ccpartialJ}, we already know that $\partial$ is a gated subcomplex. In particular, it is a quasi-median complex on its own right, and we deduce that it is a complete CAT(0) space by applying our induction hypothesis since $J$ being disjoint from $\partial$ implies that the number of hyperplanes of $\partial$ is strictly less than the number of hyperplanes of $X$. The convexity of $\partial$ in $Z$ and $J^+$ follows from Lemma \ref{lem:shortening}. Indeed, if we view $\partial$ as a subcomplex of $Z$ (resp. $J^+$), for every points $x,y \in \partial$, the unique geodesic in $Z$ (resp. $J^+$) between $x$ and $y$ must be contained in the gated hull of $\{p_x,p_y \}$, which in its turn must be included into $\partial$ since $\partial$ is gated.
\end{proof}

\begin{cor}\label{cor:qmcgeodesic}
Quasi-median complexes are geodesic metric spaces. 
\end{cor}

\begin{proof}
Let $X$ be a quasi-median complex and $x,y \in X$ two points. Let $Y$ denote the gated hull of $\{p_x,p_y \}$. According to Lemma \ref{lem:shortening}, 
$$\begin{array}{lcl} d_X(x,y) & = & \inf \{ \ell (\Sigma) \mid \Sigma \ \text{chain between $x$ and $y$} \} \\ & = & \inf \{ \ell(\Sigma) \mid \Sigma \subset Y  \ \text{chain between $x$ and $y$} \} \\ & = & d_Y(x,y). \end{array}$$
On the other hand, since $Y$ is a cubically finite quasi-median complex according to Corollary \ref{cor:finitehull}, it follows from Lemma \ref{lem:finiteCAT0} that $Y$ is a CAT(0) space. Thus, the unique geodesic between $x$ and $y$ in $Y$ produces a geodesic in $X$. 
\end{proof}

\begin{proof}[Proof of Theorem \ref{thm:CAT0}.]
Let $X$ be a quasi-median complex and $x,y,z \in X$ three points defining a geodesic triangle $T$. Let $Y$ denote the gated hull of $\{p_x,p_y,p_z\}$. It follows from Lemma \ref{lem:shortening} that $Y$ is convex, hence $T \subset Y$. On the other hand, according to Corollary \ref{cor:finitehull}, $Y$ is a cubically finite quasi-median complex, so that $T$ must satify the CAT(0) inequality since $Y$ is CAT(0) according to Lemma \ref{lem:finiteCAT0}. 
\end{proof}

\subsection{Locally quasi-median prism complexes}\label{section:locallyQM}

\noindent
In \cite[Paragraph 4.2.C]{Gromov1987}, Gromov noticed that a (finite-dimensional) cube complex defines a CAT(0) space if and only if it is simply connected and if the links of its vertices are simplicial flag complexes. Because this condition on the links turns out to be often easy to verify, CAT(0) cube complexes have become a convenient source of CAT(0) spaces. In this section, we are interested in finding a similar criterion for determining whether a prism complex is or not quasi-median. 

\begin{definition}
Let $X$ be a prism complex and $v \in X$ a vertex. The \emph{simplicial part} $\mathrm{link}_{\triangle}(v)$ of the link $\mathrm{link}(v)$ is the subcomplex generated by the edges associated to 2-simplices of $X$; and the \emph{cubical part} $\mathrm{link}_{\square}(v)$ of the link $\mathrm{link}(v)$ is the subcomplex generated by the edges associated to 2-cubes of $X$.
\end{definition}

\noindent
Roughly speaking, our criterion can be thought of as follows: the cubical part of the prism complex must behave like a CAT(0) cube complex; its simplicial part must be a union of simplices such that the intersection between two such simplices is at most a single vertex; and the cubical and simplicial parts must interact in a ``nonpositively curved'' way, meaning that the interior of a prism is always filled in. Before stating the criterion we will be interested in, let us mention that, following \cite{bucolic}, we define a \emph{flag simplicial complex} as a simplicial complex in which every cycle of length three bounds a 2-simplex, and a \emph{flag prism complex} as a prism complex in which every cycle of length three bounds a 2-simplex and every induced cycle of length four bounds a 2-cube. 

\begin{definition}
A prism complex $X$ is \emph{locally quasi-median} if it is flag and if, for every vertex $v \in X$,
\begin{itemize}
	\item $\mathrm{link}_{\square}(v)$ is a simplicial flag complexes;
	\item $\mathrm{link}_{\triangle}(v)$ is a disjoint union of simplices;
	\item a simplex and a cube both containing $v$ and intersecting along (at least) an edge must be contained into a prism.
\end{itemize}
\end{definition}

\noindent
Our main criterion is the following. The difficult part of the proof, namely showing that the 1-skeleton of a simply connected locally quasi-median complex is a weakly modular graph, is contained in \cite{bucolic}.

\begin{thm}
A prism complex is quasi-median if and only if it locally quasi-median and simply connected.
\end{thm}

\begin{proof}
Let $X$ be a prism complex. Suppose that $X$ is quasi-median. It is flag by definition and, because $X$ is a CAT(0) space, according to Theorem \ref{thm:CAT0}, it must be simply connected (and even contractible). Next, let $v \in X$ be a vertex. It follows from Fact \ref{fact:spanningprism} that $\mathrm{link}_{\square}(v)$ is a flag complex, and because the intersection between two distinct cliques is either empty or a single vertex, according to Lemma \ref{lem:cliqueinter}, we also know that $\mathrm{link}_{\triangle}(v)$ is a disjoint union of simplices. Finally, if $C$ is a simplex and $Q$ a cube such that both contain $v$ and such that their intersection contains an edge, then $C$ and $Q$ must be included into the carrier of the hyperplane dual to the clique underlying $C$, so that it follows from Lemma \ref{lem:carrierproduct} that $C$ and $Q$ generate a prism.

\medskip \noindent
Conversely, suppose that $X$ is locally quasi-median and simply connected. Then $X$ is a \emph{bucolic complex}, as defined in \cite{bucolic}, and it follows from \cite[Theorem 1]{bucolic} that the $1$-skeleton $X^{(1)}$ is a weakly modular graph. Moreover, because cycles of length three in $X$ bound 2-simplices and because the simplicial parts of links are disjoint unions of simplices, $K_4^-$ cannot be an induced subgraph of $X^{(1)}$; and similarly, because induced cycles of length four in $X$ bound 2-cubes and because the cubical parts of links are simplicial complexes, $K_{2,3}$ cannot be an induced subgraph of $X^{(1)}$. Therefore, $X^{(1)}$ is a quasi-median graph, and a fortiori $X$ is a quasi-median prism complex.
\end{proof}

\section{Metrizing quasi-median graphs}\label{section:metrizingQM}

In this section, we show that, if each clique $C$ of a quasi-median graph $X$ is endowed with a metric $\delta_C$, in such a way that the collection of all these metrics is \emph{coherent}, then there exists a global metric $\delta$ extending them. Moreover, if a group $G$ acts on $X$ and if our collection of metrics is \emph{$G$-invariant}, then the action $G \curvearrowright X$ induces an isometric action $G \curvearrowright (X, \delta)$. (The existence such collection of metrics will be studied in Section \ref{section:topicalactionsI}.) The main idea is that the global geometry of $(X,\delta)$ reduces to the local geometries of the $(C, \delta_C)$'s.

\subsection{Canonical bijections}\label{section:canonical}

\noindent
Before dealing with metrics on quasi-median graphs, we define a \emph{canonical bijection} $t_{C \to C'} : C \to C'$ between any two cliques $C,C'$ dual to the same hyperplane. This family of bijections will be fundamental not only in this section but in the whole article. 

\begin{definition}
Let $X$ be a quasi-median graph and $C,C'$ two cliques dual to the same hyperplane. The \emph{canonical bijection from $C$ to $C'$}\index{Canonical bijections}, denoted by $t_{C \to C'}$, is the restriction of the projection $\mathrm{proj}_{C'} : X \to C'$ to $C$.
\end{definition}

\noindent
Alternatively, if $J$ denotes the hyperplane dual to $C$ and $C'$, then Lemma \ref{lem:carrierproduct} yields an isomorphism $\Psi : N(J) \to F(J) \times C$, and conjugating $t_{C \to C'}$ by $\Psi$ produces the natural bijection $(C,x) \mapsto (C',x)$. Loosely speaking, $t_{C \to C'}$ translates $C$ along $J$ to $C'$. 

\medskip \noindent
Now, we register some basic facts about the canonical bijections for future use. 

\begin{lemma}\label{lem:projectionandtransfer}
For any two cliques $C,C'$ dual to the same hyperplane, $$\mathrm{proj}_{C'}= t_{C \to C'} \circ \mathrm{proj}_C.$$
\end{lemma}

\noindent
The following observation will be needed to prove Lemma \ref{lem:projectionandtransfer}.

\begin{fact}\label{fact:projandsector}
Let $X$ be a quasi-median graph, $C$ a clique and $x \in X$ a vertex. If $D$ denotes the sector delimited by the hyperplane dual to $C$ which contains $x$, then $D \cap C = \{ \mathrm{proj}_C(x) \}$. 
\end{fact}

\begin{proof}
Let $J$ denote the hyperplane dual to $C$. According to Lemma \ref{lem:projseparate}, any hyperplane separating $x$ and its projection onto $C$ must be disjoint from $C$, so $J$ cannot separate $x$ and $\mathrm{proj}_C(x)$. It follows that $x$ and $\mathrm{proj}_C(x)$ belong to the same sector delimited by $J$, hence $\mathrm{proj}_C(x) \in D \cap C$. On the other hand, any two vertices of $C$ are separated by $J$, so $D \cap C$ contains a single vertex. The conclusion follows.
\end{proof}

\begin{proof}[Proof of Lemma \ref{lem:projectionandtransfer}.]
Let $J$ denote the hyperplane dual to $C,C'$. Fix a vertex $x \in X$ and let $D_x$ denote the sector delimited by $J$ containing $x$. As a consequence of Fact \ref{fact:projandsector}, $x$ and $\mathrm{proj}_C(x)$ belong to $D_x$. Next, by applying Fact \ref{fact:projandsector} twice, we know that
$$\{ \mathrm{proj}_{C'}(x) \} = D_x \cap C' = \{ \mathrm{proj}_{C'} \circ \mathrm{proj}_C (x) \}.$$
Therefore,
$$\mathrm{proj}_{C'}= \mathrm{proj}_{C'} \circ \mathrm{proj}_C = t_{C \to C'} \circ \mathrm{proj}_C,$$
which concludes the proof.
\end{proof}

\begin{lemma}\label{lem:composingbij}
If $C,C',C''$ are three cliques dual to the same hyperplane, then 
$$t_{C' \to C''} \circ t_{C \to C'} = t_{C \to C''}.$$
\end{lemma}

\begin{proof}
By applying Lemma \ref{lem:projectionandtransfer} twice, we deduce that
$$t_{C \to C''} = \mathrm{proj}_{C''} \circ t_{C \to C'} = t_{C' \to C''} \circ \mathrm{proj}_{C'} \circ t_{C \to C'}.$$
On the other hand, it is clear that $\mathrm{proj}_{C'} \circ t_{C \to C'}= t_{C \to C'}$, so that the conclusion follows.
\end{proof}

\begin{cor}
If $C,C'$ are two cliques dual to the same hyperplane, then $$t_{C \to C'}^{-1}= t_{C' \to C}.$$
\end{cor}

\begin{proof}
By applying the previous lemma, we know that
$$t_{C' \to C} \circ t_{C \to C'}= t_{C \to C}= \mathrm{id}_C.$$
The conclusion follows.
\end{proof}

\begin{remark}\label{remark:monodromy}
As a consequence of Lemma \ref{lem:composingbij}, if $C_1, \ldots, C_n$ is a sequence of cliques all dual to the same hyperplanes, where $C_1=C=C_n$ for some clique $C$, then 
$$t_{C_{n-1} \to C_n} \circ t_{C_{n-2} \to C_{n-1}} \circ \cdots \circ t_{C_2 \to C_3} \circ t_{C_1 \to C_2} = \mathrm{id}_{C}.$$
This \emph{monodromy condition} satisfied by the canonical bijections generalizes the fact that, in a CAT(0) cube complex, hyperplanes are two-sided. 
\end{remark}

\subsection{Extending metrics}

\noindent
Given a quasi-median graph $X$, a \emph{system of (pseudo-)metrics}\index{Systems of metrics} is the data, for every clique $C$ of $X$, of a (pseudo-)metric $\delta_C$ defined on $C$. If $C$ and $C'$ are two cliques dual to the same hyperplane, we define a new metric on $C'$ by
$$\delta_{C \to C'} : (x,y) \longmapsto \delta_C \left( t_{C' \to C}(x),t_{C' \to C}(y) \right),$$
where $t_{C \to C'} : C \to C'$ is the canonical bijection defined in Section \ref{section:canonical}. Our system of (pseudo-)metrics is \emph{coherent} if $\delta_{C'}= \delta_{C \to C'}$ for every pair of cliques $C,C'$ dual to the same hyperplane. In this section, our goal is to extend such a collection of (pseudo-)metrics to a global one. In fact, several natural extensions exist, depending on some parameter $p \in [1,+ \infty]$ we fix. 

\medskip \noindent
Given two vertices $x,y \in X$, a \emph{chain} $\Sigma$ between $x$ and $y$ is a sequence of vertices 
$$x_0=x, \ x_1, \ldots, \ x_{n-1}, \ x_n=y$$
such that, for every $0 \leq i \leq n-1$, $x_i$ and $x_{i+1}$ belong to a common prism $P_i=C^i_1 \times \cdots \times C^i_{n(i)}$. The \emph{$\ell^p$-length} of $\Sigma$ is defined by
$$\ell_p(\Sigma)= \sum\limits_{i=0}^{n-1} \delta^p_{P_i}(x_i,x_{i+1}),$$
where $\delta^p_{P_i}$ denotes the $\ell^p$-distance of the product $(C^i_1,\delta_{C^i_1}) \times \cdots \times (C^i_{n(i)}, \delta_{C^i_{n(i)}})$, ie. 
$$\delta_{P_i}^p : \left( (y_k), (z_k) \right) \mapsto \left( \sum\limits_{k=1}^{n(i)} \delta_{C_k^i}(y_k,z_k)^p \right)^{1/p}$$
if $p<+ \infty$ and 
$$\delta_{P_i}^{\infty} : ((y_k),(z_k)) \mapsto \max\limits_{1 \leq k \leq n(i)} \delta_{C_k^i}(y_k,z_k)$$ 
otherwise. It is worth noticing that the quantity $\delta_{P_i}^p(x_i,x_{i+1})$ does not depend on the choice of the prism $P_i$ containing both $x_i$ and $x_{i+1}$, and that the metric $\delta_{P_i}^p$ does not depend on the way we decomposed $P_i$ as a product of cliques, because our system of (pseudo-)metrics is coherent. Finally, we define the pseudo-distance
$$\delta^p(x,y)= \inf \{ \ell_p(\Sigma) \mid \Sigma \ \text{chain between $x$ and $y$} \}.$$
Notice that, for every vertices $x,y \in X$, the quantity $\delta^p(x,y)$ is necessarily finite. Indeed, fix a geodesic $z_1, \ldots, z_n$ from $x$ to $y$. A fortiori, for every $1 \leq i \leq n-1$, the vertices $z_i$ and $z_{i+1}$ are adjacent, so that there exists a unique clique $C_i$ containing them. Therefore, $z_1, \ldots, z_n$ defines a chain $\Sigma$ from $x$ to $y$, hence
$$\delta^p(x,y) \leq \ell_p(\Sigma)= \left( \sum\limits_{i=1}^{n-1} \delta_{C_i}(z_i,z_{i+1})^p \right)^{1/p} < + \infty$$
if $p$ is finite, and
$$\delta^{\infty}(x,y) \leq \ell_{\infty}(\Sigma) = \max\limits_{1 \leq i \leq n-1} \delta_{C_i}(z_i,z_{i+1}) <+ \infty$$
otherwise. Thus, we have defined a global pseudo-distance $\delta^p$. Now, we want to prove that $\delta^p$ turns out to be a metric (if we started from a collection of distances) which extends the $\delta_C$'s. 

\begin{prop}\label{prop:distext}
Let $X$ be a quasi-median graph endowed with a coherent system of metrics. The global pseudo-metric $\delta^p$ is a distance. Moreover, for every gated subgraph $Y$, the inclusion $Y \subset X$ induces an isometric embedding $(Y,\delta^p_Y) \hookrightarrow (X,\delta^p)$ where $\delta_Y^p$ denotes the global metric on $Y$ obtained by extending the (restriction to $Y$ of) our system of metrics.
\end{prop}

\noindent
The second part of the statement is contained in the following lemma, which follows by reproducing word for word the proof of Lemma \ref{lem:shortening}.

\begin{lemma}\label{lem:shorteningchain}
Let $X$ be a quasi-median graph endowed with a coherent system of metrics. Let $x,y \in X$ be two vertices. For every chain $\Sigma$ between $x$ and $y$, there exists another chain $\Sigma^*$ between $x$ and $y$ satisfying $\ell(\Sigma^*) \leq \ell(\Sigma)$ which is contained into the gated hull of $\{x,y \}$.
\end{lemma}

\begin{proof}[Proof of Proposition \ref{prop:distext}.]
Let $x,y \in X$ be two distinct vertices. A fortiori, there must exist a hyperplane $J$ separating them. Set 
$$\delta_J(x,y)= \delta_C(\mathrm{proj}_C(x), \mathrm{proj}_C(y)),$$
where $C$ is a clique dual to $J$. Notice that $\delta_J(x,y)>0$ because $J$ separates $x$ and $y$. Moreover, this quantity does not depend on the choice of $C$. Indeed, if $C$ and $C'$ are two cliques dual to the same hyperplane, then we deduce from Lemma \ref{lem:projectionandtransfer} that
$$\delta_C(p(x),p(y))= \delta_C(t \circ q(x),t \circ q(y)) = \delta_{C \to C'}(q(x),q(y))= \delta_{C'}(q(x),q(y)),$$
where $t= t_{C' \to C} : C' \to C$ is the canonical bijection from $C$ to $C'$, $p : X \to C$ the projection onto $C$ and $q : X \to C'$ the projection onto $C'$. Now, let $\Sigma= (x_1, \ldots, x_n)$ be a chain between $x$ and $y$, and, for every $1 \leq i \leq n-1$, let $P_i$ denote the minimal prism containing both $x_i$ and $x_{i+1}$. Because $J$ separates $x$ and $y$, the set of prisms $\{ P_{m_1}, \ldots, P_{m_r} \}$ intersected by $J$ is non empty (for convenience, suppose that $m_1< \cdots < m_r$). Notice that $x$ and $x_{m_1}$ are not separated by $J$, so that $p(x)=p(x_{m_1})$; similarly, $p(y)=p(x_{m_r+1})$. Therefore,
$$\begin{array}{lcl} \ell_p(\Sigma) & = & \sum\limits_{i=1}^{n-1} \delta_{P_i}^p(x_i,x_{i+1}) \geq \sum\limits_{k=1}^{r-1} \delta_{P_{m_k}}^p(x_{m_k},x_{m_k+1}) \geq \sum\limits_{k=1}^{r-1} \delta_C(p(x_{m_k}),p(x_{m_k+1})) \\ \\ & \geq & \delta_C(p(x_{m_1}),p(x_{m_r+1})) = \delta_C(p(x),p(y)) = \delta_J(x,y) \end{array}$$
Because the quantity $\delta_J(x,y)$ does not depend on the choice of the chain $\Sigma$, it follows that $\delta^p(x,y) \geq \delta_J(x,y)>0$. Thus, we have proved that $\delta^p$ is a distance. The second statement of our proposition is a direct consequence of Lemma \ref{lem:shorteningchain}.
\end{proof}

\noindent
So far, we have given infinitely many ways to extend a coherent system of metrics in order to get a global distance on the vertices of a quasi-median graph. Our next statement shows that, if the cubical dimension of our quasi-median graph is finite, then these global metrics are all Lipschitz-equivalent. When the cubical dimension is infinite, the global geometry may be more interesting with respect to one of the possible metrics. For instance, in the case of CAT(0) cube complexes, we noticed in \cite{coningoff} that the CAT(0) cube complex associated to a $C'(1/4)-T(4)$ polygonal complex is always hyperbolic with respect to the $\ell^{\infty}$-distance, even when it is infinite-dimensional (if so, it cannot be hyperbolic with respect to the $\ell^p$-distance for any $p \in [1,+ \infty)$).

\begin{prop}\label{prop:comparedelta}
Let $X$ be a quasi-median graph endowed with a coherent system of metrics. If $q \in [1,+ \infty]$ satisfies $\frac{1}{p} + \frac{1}{q} = 1$, then
$$\delta^p \leq \delta^1 \leq \dim_{\square}(X)^{1/q} \cdot \delta^p.$$
\end{prop}

\begin{proof}
Let $x,y \in X$ be two vertices. Fix a chain $\Sigma=(x_1, \ldots, x_n)$ between $x$ and $y$, and, for every $1 \leq i \leq n-1$, let $P_i$ denote a prism containing both $x_i$ and $x_{i+1}$. Because $p \geq 1$, we know that $\delta_{P_i}^p \leq \delta_{P_i}^1$ for every $1 \leq i \leq n-1$, hence 
$$\delta^p(x,y) \leq \ell_p(\Sigma) = \sum\limits_{i=1}^{n-1} \delta_{P_i}^p(x_i,x_{i+1}) \leq \sum\limits_{i=1}^{n-1} \delta_{P_i}^1 (x_i,x_{i+1}) = \ell_1(\Sigma).$$
A fortiori, $\delta^p(x,y) \leq \delta^1(x,y)$. On the other hand, it follows from H\"{o}lder's inequality that 
$$\delta_{P_i}^1(x_i,x_{i+1}) \leq \dim_{\square}(P_i)^{1/q} \cdot \delta_{P_i}^p(x_i,x_{i+1}) \leq  \dim_{\square}(X)^{1/q} \cdot \delta_{P_i}^p(x_i,x_{i+1})$$
for every $1 \leq i \leq n-1$, hence 
$$\delta^1(x,y) \leq \ell_1(\Sigma) \leq \dim_{\square}(X)^{1/q} \cdot \ell_p(\Sigma).$$
A fortiori, $\delta^1(x,y) \leq \dim_{\square}(X)^{1/q} \cdot \delta^p(x,y)$. This concludes the proof. 
\end{proof}

\noindent
A natural problem is to determine how the (global) geometry of $(X,\delta^p)$ reduces to the (local) geometries of the cliques $(C,\delta_C)$. So we ask the following question:

\begin{question}
Let $X$ be a quasi-median graph endowed with a coherent system of metrics, and $\mathcal{P}$ a property of metric spaces. If, for every clique $C$, the metric space $(C,\delta_C)$ satisfies the property $\mathcal{P}$, when does $(X,\delta^p)$ satisfies $\mathcal{P}$ as well?
\end{question}

\noindent
In this section and the next ones, we will answer this question for several properties, including being geodesic (Corollary \ref{cor:systgeodesic}), being locally finite (Lemma \ref{lem:whenlocallyfinite}; see also Remark \ref{remark:locallyfinite}), and being CAT(0) (Proposition \ref{prop:CAT0delta}); Proposition \ref{prop:compression} also studies the $\ell^p$-compression of the global metric. 

\medskip \noindent
A good picture to keep in mind is that $(X, \delta^p)$ is a kind of a ``cubical agregate''. A CAT(0) cube complex is a union of cubes, which are products of edges, glued together in a ``non-positively curved'' way. Similarly, $(X,\delta^p)$ is a union of prisms, which are products of the spaces $(C,\delta_C)$, glued together in a ``non-positively curved'' way. From this analogy, it would not be surprising that any property characterizing a ``non-positively curved'' behaviour transfers from local to global. We will motivate this idea by proving the following statement:  

\begin{prop}\label{prop:CAT0delta}
Let $X$ be a quasi-median graph endowed with a coherent system of complete CAT(0) metrics. Then $(X,\delta^2)$ is a CAT(0) space in which gated subgraphs of $X$ are convex.
\end{prop}

\begin{proof}
First, let us prove by induction on the number of hyperplanes that, for any cubically finite gated subgraph $Y \subset X$, the metric space $(Y,\delta^2)$ is a complete CAT(0) space. If $Y$ has no hyperplane, it is a single vertex and there is nothing to prove; if $Y$ has a single hyperplane, it is a single clique and there is nothing to prove. From now on, suppose that $Y$ contains at least two hyperplanes. If the hyperplanes of $Y$ are pairwise transverse, it follows from Lemma \ref{lem:prismhyp} that $Y$ is a prism, so that $(Y,\delta^2)$ is isometric to the $\ell^2$-product of finitely many complete CAT(0) spaces. It follows from \cite[Examples II.1.5]{MR1744486} that $(Y,\delta^2)$ is a complete CAT(0) space. Next, suppose that $Y$ contains two disjoint hyperplanes, say $J_1$ and $J_2$. Let $Y_2$ denote the sector delimited by $J_1$ containing $J_2$ and $Y_1$ the union of all the other sectors with $N(J_1)$. According to Fact \ref{fact:multisectorgatedhull} and Proposition \ref{prop:hypsumup} respectively, $Y_1$ and $Y_2$ are gated. By applying our induction hypothesis, we know that $(Y_1,\delta^2)$ and $(Y_2, \delta^2)$ are complete CAT(0) spaces. Moreover, the intersection $Y_1 \cap Y_2$ is a fiber of $J_1$, which is also a gated subgraph according to Proposition \ref{prop:hypsumup}, so we deduce from our induction hypothesis that $Y_1 \cap Y_2$ is a complete subspace of both $(Y_1,\delta^2)$ and $(Y_2,\delta^2)$; also, note that, as a consequence of Lemma \ref{lem:shorteningchain}, $Y_1 \cap Y_2$ is geodesic, and a fortiori convex since CAT(0) spaces are uniquely geodesic, in both $(Y_1,\delta^2)$ and $(Y_2,\delta^2)$. Next, notice that, by reproducing word for word the proof of Claim \ref{claim:gluingdist}, we get the following general statement:

\begin{fact}\label{fact:gluingdelta}
Let $X$ be a quasi-median graph endowed with a coherent system of metrics, and $Y$ a gated subgraph containing two disjoint hyperplanes $J_1$ and $J_2$.  Let $Y_2$ denote the sector delimited by $J_1$ containing $J_2$, $Y_1$ the union of all the other sectors with $N(J_1)$, and $\partial$ the fiber $N(J_1) \cap Y_2$ of $J_1$. Fix some $p \in [1,+ \infty]$. Then
$$\delta^p(x,y) = \inf\limits_{z \in \partial} \left( \delta^p(x,z) + \delta^p(z,y) \right)$$
for every vertices $x \in Y_1$ and $y \in Y_2$.
\end{fact}

\noindent
In particular, the equality
$$\delta^2(x,y)= \inf\limits_{z \in Y_1 \cap Y_2} \left( \delta^2(x,z)+ \delta^2(z,y) \right)$$
holds for every vertices $x \in Y_1$ and $y \in Y_2$. Therefore, it follows from the Gluing Theorem \cite[Theorem II.11.1]{MR1744486} that $(Y,\delta^2)$ is a complete CAT(0) space. 

\medskip \noindent
Now, we can notice that $(X, \delta^2)$ is a geodesic metric space. Indeed, the interval between two vertices is a cubically finite gated subgraph, and it follows from our previous observation that such an interval must be geodesic. On the other hand, such a geodesic defines a geodesic in $(X,\delta^2)$ as a consequence of Lemma \ref{lem:shorteningchain}. Now, it makes sense to state that gated subgraphs of $X$ are convex in $(X, \delta^2)$, which is also a consequence of Lemma \ref{lem:shorteningchain}.

\medskip \noindent
Finally, let $\Delta=(x,y,z)$ be a geodesic triangle in $(X, \delta^2)$. Since gated subgraphs are convex, necessarily $\Delta$ must be included into the gated hull $H$ of $\{x,y,z\}$. On the other hand, this hull is a cubically finite gated subgraph, so that the CAT(0) inequality must be satisfied in $H$, and a fortiori in $(X, \delta^2)$. This concludes the proof.
\end{proof}

\subsection{More on the global metric $\delta^1$}\label{section:extendingmetrics}

\noindent
According to Proposition \ref{prop:comparedelta}, our global metrics are all Lipschitz equivalent when the cubical dimension of the underlying quasi-median graph is finite, so that it is often possible to choose the metric which is the easiest to handle in the context we are interested in. In this section, our goal is to show that the global metric $\delta^1$ can be described in a more explicit way, justifying the fact it will be our favorite metric most of the time. In all the article, we will set $\delta=\delta^1$ for short.

\medskip \noindent
Fix a quasi-median graph endowed with a coherent system of metrics. For every hyperplane $J$ of $X$, we associate a pseudo-metric $\delta_J$ on $X$ defined by
$$\delta_J : (x,y) \mapsto \delta_C( \mathrm{proj}_C(x), \mathrm{proj}_C(y)),$$
where $C$ is any clique dual to $J$. Notice that $\delta_J$ does not depend on the choice of the clique $C$ since, whenever $C$ and $C'$ are two cliques dual to the same hyperplane, we deduce from Lemma \ref{lem:projectionandtransfer} that
$$\delta_C(p(x),p(y))= \delta_C(t \circ q(x),t \circ q(y)) = \delta_{C \to C'}(q(x),q(y))= \delta_{C'}(q(x),q(y)),$$
where $t= t_{C' \to C} : C' \to C$ is the canonical bijection from $C$ to $C'$, $p : X \to C$ the projection onto $C$ and $q : X \to C'$ the projection onto $C'$. Notice that, for every vertices $x,y \in X$, the quantity $\delta_J(x,y)$ is non zero if and only if $J$ separates $x$ and $y$.

\begin{prop}
Let $X$ be a quasi-median graph endowed with a coherent system of metrics and $x,y \in X$ two vertices. Then
$$\delta(x,y)= \sum\limits_{\text{$J$ hyperplane}} \delta_J(x,y) = \sum\limits_{\text{$J$ separates $x$ and $y$}} \delta_J(x,y).$$
\end{prop}

\begin{proof}
We want to prove by induction on $k$ that, for any gated subgraph $Y$ containing at most $k$ hyperplanes, the equality $$\delta(x,y) = \sum\limits_{\text{$J$ separates $x$ and $y$}} \delta_J(x,y)$$ holds for every $x,y \in Y$. 

\medskip \noindent
Let $Y$ be a gated subgraph with finitely many hyperplanes. If $Y$ contains no hyperplane then $Y$ is a single vertex and there is nothing to prove. And if the hyperplanes of $Y$ are pairwise transverse, then Lemma \ref{lem:prismhyp} implies that $Y$ is a single prism, so that the equality follows. From now on, suppose that $Y$ contains at least two disjoint hyperplanes, say $J_1$ and $J_2$. Let $Y_2$ denote the sector delimited by $J_1$ which contains $J_2$, and $Y_1$ the union of all the other sectors with $N(J_1)$. According to  Fact \ref{fact:multisectorgatedhull} and Proposition \ref{prop:hypsumup} respectively, $Y_1$ and $Y_2$ are gated subgraphs. Moreover, the numbers of hyperplanes of $Y_1$ and $Y_2$ are smaller than the number of hyperplanes of $Y$, so that our induction hypothesis applies, ie., 
$$\delta(x,y) = \sum\limits_{\text{$J$ separates $x$ and $y$}} \delta_J(x,y)$$
for every vertices $x$ and $y$ which both belong to either $Y_1$ or $Y_2$. In particular, we know that the equality we want to prove holds if we consider two vertices of $Y_1$ or two vertices of $Y_2$, so, in order to conclude, it is sufficient to verify our equality for two vertices $x \in Y_1$ and $y \in Y_2$. On the other hand, if we denote by $\partial$ the fiber $N(J_1) \cap Y_2$ of $J_1$, we know from Fact \ref{fact:gluingdelta} that
$$\delta(x,y) = \inf\limits_{z \in \partial} \left( \delta(x,z)+ \delta(z,y) \right).$$
Fix some $z \in \partial$. Let $x',y'$ denote respectively the projections of $x,y$ onto $\partial$, and let $(x'',y'',z'')$ be the quasi-median of the triple $(x',y',z)$ (as defined in Section \ref{section:qm}). Notice that, because $\partial$ is a quasi-median graph on its own right, and because there exists a unique quasi-median in a quasi-median graph according to Proposition \ref{prop:quasimedian}, necessarily $x'',y'',z''$ belong to $\partial$. 

\medskip \noindent
For convenience, let us introduce the following notation. If $a,b \in X$ are two vertices, we denote by $\mathfrak{J}(a,b)$ the set of the hyperplanes separating $a$ and $b$. Notice that, if $c \in I(a,b)$, ie., if $c$ lies on a geodesic between $a$ and $b$, then $\mathfrak{J}(a,b) = \mathfrak{J}(a,c) \sqcup \mathfrak{J}(c,b)$.

\medskip \noindent
Because $x'$ is the projection of $x$ onto $\partial$ and that $z \in \partial$, necessarily $x' \in I(x,z)$. Moreover, by definition of a median triangle, there exists a geodesic between $x'$ and $z$ passing through $x''$ and $z''$. Therefore,
$$\mathfrak{J}(x,z)= \mathfrak{J}(x,x') \sqcup \mathfrak{J}(x',x'') \sqcup \mathfrak{J}(x'',z'') \sqcup \mathfrak{J}(z'',z).$$
Similarly,
$$\mathfrak{J}(z,y)= \mathfrak{J}(z,z'') \sqcup \mathfrak{J}(z'',y'') \sqcup \mathfrak{J}(y'',y') \sqcup \mathfrak{J}(y',y)$$
As a consequence,
$$\begin{array}{lcl} \delta(x,z)+ \delta(z,y) & = & \displaystyle \sum\limits_{J \in \mathfrak{J}(x,z)} \delta_J(x,z)+ \sum\limits_{J \in \mathfrak{J}(z,y)} \delta_J(z,y) \\ \\ & \geq & \displaystyle \sum\limits_{J \in \mathfrak{J}(x,x') \sqcup \mathfrak{J}(x',x'')} \delta_J(x,z) + \sum\limits_{J \in \mathfrak{J}(x'',z'')} \delta_J(x,z) \\ & & \displaystyle + \sum\limits_{J \in \mathfrak{J}(z'',y'')} \delta_J(z,y) + \sum\limits_{J \in \mathfrak{J}(y'',y') \sqcup \mathfrak{J}(y',y)} \delta_J(z,y) \end{array}$$
Notice that, if $J \in \mathfrak{J}(x,x') \sqcup \mathfrak{J}(x',x'')$, then $z$ and $y''$ are not separated by $J$. Indeed, if $J \in \mathfrak{J}(x,x')$, then it follows from Lemma \ref{lem:projseparate} that $J$ separates $x$ from $\partial$, but $z$ and $y''$ belongs to $\partial$, so $J$ cannot separate them; and if $J \in \mathfrak{J}(x',x'')$, our claim follows from Fact \ref{fact:1}. So
$$\sum\limits_{J \in \mathfrak{J}(x,x') \sqcup \mathfrak{J}(x',x'')} \delta_J(x,z)= \sum\limits_{J \in \mathfrak{J}(x,x') \sqcup \mathfrak{J}(x',x'')} \delta_J(x,y''),$$
and similarly, 
$$\sum\limits_{J \in \mathfrak{J}(y'',y') \sqcup \mathfrak{J}(y',y)} \delta_J(z,y) = \sum\limits_{J \in \mathfrak{J}(y'',y') \sqcup \mathfrak{J}(y',y)} \delta_J(y'',y).$$
Next, notice that it follows from Fact \ref{fact:5} that $\mathfrak{J}(x'',z'')= \mathfrak{J}(z'',y'')= \mathfrak{J}(x'',y'')$, so that
$$\begin{array}{lcl} \displaystyle \sum\limits_{J \in \mathfrak{J}(x'',z'')} \delta_J(x,z)+ \sum\limits_{J \in \mathfrak{J}(z'',y'')} \delta_J(z,y) & = & \displaystyle \sum\limits_{J \in \mathfrak{J}(x'',y'')} \left( \delta_J(x,z)+ \delta_J(z,y) \right) \\ \\ & \geq & \displaystyle \sum\limits_{J \in \mathfrak{J}(x'',y'')} \delta_J(x,y)= \sum\limits_{J \in \mathfrak{J}(x'',y'')} \delta_J(x,y'') \end{array}$$
The last equality is justified by the fact that no hyperplane of $\mathfrak{J}(x'',y'')$ separates $y$ and $y''$. Indeed, a hyperplane $J$ separating $y$ and $y''$ separates either $y$ and $y'$ or $y'$ and $y''$. In the former case, it follows from Lemma \ref{lem:projseparate} that $J$ separates $y$ from $\partial$, so that it cannot separate $y$ and $y''$ since $y'' \in \partial$; and in the latter case, the conclusion is a consequence of Fact \ref{fact:1}. By noticing that
$$\mathfrak{J}(x,x') \sqcup \mathfrak{J}(x',x'') \sqcup \mathfrak{J}(x'',y'')= \mathfrak{J}(x,y'') \ \text{and} \ \mathfrak{J}(y'',y') \sqcup \mathfrak{J}(y',y) = \mathfrak{J}(y'',y),$$
we finally get 
$$\begin{array}{lcl} \delta(x,z)+ \delta(z,y) & \geq & \displaystyle \sum\limits_{J \in \mathfrak{J}(x,x') \sqcup \mathfrak{J}(x',x'') \sqcup \mathfrak{J}(x'',y'')} \delta_J(x,y'') + \sum\limits_{J \in \mathfrak{J}(y'',y') \sqcup \mathfrak{J}(y',y)} \delta_J(y'',y) \\ \\ & \geq & \displaystyle \sum\limits_{J \in \mathfrak{J}(x,y'')} \delta_J(x,y'') + \sum\limits_{J \in \mathfrak{J}(y'',y)} \delta_J(y'',y) = \delta(x,y'')+\delta(y'',y) \end{array}$$
where the last equality is justified by the fact that $x,y''$ both belong to $Y_1$ and that $y'',y$ both belong to $Y_2$. Notice that $y'' \in I(x,y)$. Indeed, by definition of a quasi-median triangle, we know that there exists a geodesic $[x',y'] \subset \partial$ between $x'$ and $y'$ passing through $y''$. On the other hand, a hyperplane separating $x$ and $x'$ must separate $x$ from $\partial$ according to Lemma \ref{lem:projseparate}, and a fortiori from $Y_2$, so that it cannot intersect $[x',y']$ or separate $y$ and $y'$; and similarly, a hyperplane separating $y$ and $y'$ cannot intersect $[x',y']$ or separate $x$ and $x'$. Therefore, fixing some geodesics $[x,x']$ and $[y,y']$ respectively between $x$ and $x'$, and $y$ and $y'$, we deduce that the concatenation $[x,x'] \cup [x',y'] \cup [y',y]$ is a geodesic since it cannot intersect a hyperplane twice. So we have proved that $y''$ belongs to some geodesic between $x$ and $y$, which precisely means that $y''$ belongs to the interval $I(x,y)$. A fortiori, $y'' \in \partial \cap I(x,y)$. Thus, we have proved that
$$\delta(x,y) = \inf\limits_{z \in \partial \cap I(x,y)} \left( \delta(x,z)+ \delta(z,y) \right).$$
On the other hand, for every $z \in \partial \cap I(x,y)$,
$$\delta(x,z)+ \delta(z,y)= \sum\limits_{J \in \mathfrak{J}(x,z)} \delta_J(x,z) + \sum\limits_{J \in \mathfrak{J}(z,y)} \delta_J(y,z).$$
Notice that, for every $J \in \mathfrak{J}(x,z)$, the projections of $y$ and $z$ onto $N(J)$ coincide since $J$ does not separate $y$ and $z$, hence $\delta_J(x,z)=\delta_J(x,y)$. Similarly, $\delta_J(z,y)= \delta_J(x,y)$ for every $J \in \mathfrak{J}(z,y)$. Thus,
$$\delta(x,z)+ \delta(z,y) = \sum\limits_{J \in \mathfrak{J}(x,z) \sqcup \mathfrak{J}(z,y)} \delta_J(x,y) = \sum\limits_{J \in \mathfrak{J}(x,y)} \delta_J(x,y).$$
It follows that
$$\delta(x,y) = \sum\limits_{J \in \mathfrak{J}(x,y)} \delta_J(x,y),$$
concluding the proof.
\end{proof}

\begin{ex}
If $\delta_C : (x,y) \mapsto \left\{ \begin{array}{cl} 0 & \text{if} \ x=y \\ 1 & \text{otherwise} \end{array} \right.$ for every clique $C$ of $X$, then $\delta=d$. 
\end{ex}

\noindent
We know from Proposition \ref{prop:distext} that gated subgraphs define isometrically embedded subspaces of $(X,\delta^p)$, for every $p \in [1,+\infty]$. When $p=1$, we are also able to show that they are \emph{convex} in the following meaning:

\begin{definition}
Let $(S,d)$ be a metric space. The \emph{interval} of two points $x,y \in S$ is the set
$$I(x,y)= \{ z \in X \mid d(x,y)=d(x,z)+d(z,y) \}.$$
A subset $R \subset S$ is \emph{convex} if $I(x,y) \subset R$ for every $x,y \in R$.
\end{definition}

\noindent
It is worth noticing that, if $(S,d)$ is a geodesic metric space, then $I(x,y)$ is the union of all the geodesics between $x$ and $y$; and a subset $R \subset S$ is convex if and only if every geodesic between two points of $R$ stays in $R$, which is the usual definition of convexity for geodesic metric spaces.

\begin{prop}\label{prop:convexdelta1}
Let $X$ be a quasi-median graph endowed with a coherent system of metrics. Every gated subgraph of $X$ is convex in $(X,\delta)$. 
\end{prop}

\noindent
This proposition is essentially a consequence of the following lemma, which is an adaptation of  Lemma \ref{lem:shortening}.

\begin{lemma}\label{lem:shorteningquantitatif}
Let $X$ be a quasi-median graph endowed with a coherent system of metrics. Let $x,y,z \in X$ be three vertices such that $z$ does not belong to the gated hull of $\{x,y \}$. There exists some $\epsilon(z)>0$ such that, for every chain $\Sigma$ between $x$ and $y$ passing through $z$, there exists another chain $\Sigma^*$ between $x$ and $y$ satisfying $\ell(\Sigma^*) \leq \ell(\Sigma)- \epsilon(z)$.
\end{lemma}

\begin{proof}
Let $J$ be a hyperplane separating $z$ from the gated hull $H$ of $\{x,y\}$ which maximizes the distance to $H$. Notice that $z$ belongs to $N(J)$ since otherwise there would exist a hyperplane separating $z$ from $N(J)$, contradicting the maximality of $J$. Let $J^+$ denote the sector delimited by $J$ which contains $H$, and set $\partial = N(J) \cap J^+ \subset \partial J$. 

\medskip \noindent
Let $(x_r, \ldots, x_s)$ be the maximal subsegment of $\Sigma$ included into $N(J)$ which contains $z$. We claim that the points $x_r,x_s$ belong to $\partial$. 

\medskip \noindent
Suppose by contradiction that $x_r \notin \partial$. Let $P$ be a prism containing $x_{r-1}$ and $x_r$. Because $x_{r-1} \notin N(J)$, $P \nsubseteq N(J)$; but $P \cap N(J) \neq \emptyset$ because $\Sigma$ is a chain. Let $v \in P \cap N(J)$ be a vertex, $J_1, \ldots, J_m$ the hyperplanes dual to $P$, and, for every $1 \leq i \leq m$, $v_i \in P$ a vertex adjacent to $v$ and separated from it by $J_i$. If $J,J_1, \ldots, J_m$ are pairwise transverse, we deduce from Fact \ref{fact:spanningprism} that the gated hull of $\{x,x_1, \ldots, x_m \}$ defines a prism whose dual hyperplanes are $J,J_1, \ldots, J_m$, hence $P \subset N(J)$, a contradiction. Therefore, there exists a hyperplane $J'$ dual to $P$ which is disjoint from $J$. On the other hand, because $x_r \notin \partial$ and $P \nsubseteq N(J)$, the prism $P$, and a fortiori the hyperplane $J'$, must be included into a sector delimited by $J$ which is different from $J^+$, the sector containing $H$. Consequently, we have $d(J',H)>d(J,H)$, contradicting our choice of $J$. Thus, we have proved that necessarily $x_r \in \partial$. Similarly, we show that $x_s \in \partial$.

\medskip \noindent
Let $\Sigma'$ denote the subchain $(x_r,\ldots,x_s)$. For every $r \leq i \leq s-1$, fix a prism $P_i \subset N(J)$ containing both $x_i$ and $x_{i+1}$. Without loss of generality, we can suppose that $P_i$ contains a clique dual to $J$, so that $P_i$ decomposes a product $C_i \times P_i'$ where $P_i'$ is a prism included into $\partial$ and $C_i$ a clique dual to $J$. Notice that, if $m : X \to \partial$ denotes the projection onto $\partial$, then $(m(x_r),\ldots, m(x_s))$ defines a chain since $m(x_i)$ and $m(x_{i+1})$ both belong to $P_i'$ for every $r \leq i \leq s-1$; let $\Sigma''$ denote this new chain. Now, if we denote by $q : X \to C$ the projection onto some clique $C$ dual to $J$ and $q_i : X \to C_i$ the projection onto $C_i$ for every $r \leq i \leq s$, we have
$$\begin{array}{lcl} \ell(\Sigma') & = & \displaystyle \sum\limits_{i=r}^{s-1} \delta_{P_i}(x_i,x_{i+1}) =  \sum\limits_{i=r}^{s-1} \left( \delta_{P_i'}(m(x_i),m(x_{i+1})) + \delta_{C_i}(q_i(x_i),q_i(x_{i+1})) \right) \\ \\ & = & \displaystyle \ell(\Sigma'') + \sum\limits_{i=r}^{s-1} \delta_C(q(x_i),q(x_{i+1})) \end{array}$$
On the other hand, we know by assumption that there exists  some $r+1 \leq k \leq s-1$ satisfying $x_k=z$, so that
$$\sum\limits_{i=r}^{s-1} \delta_C(q(x_i),q(x_{i+1})) \geq \sum\limits_{i=r}^{k-1} \delta_C(q(x_i),q(x_{i+1})) \geq \delta_C(q(x_r),q(x_k))= \delta_C(q(x),q(z)).$$
Therefore, $\ell(\Sigma') \geq \ell(\Sigma'')+ \epsilon(z)$ where 
$$\epsilon(z)= \min \{ \delta_J(x,z) \mid \text{$J$ separates $z$ from $\{x,y\}$} \}.$$ 
Notice that $\epsilon(z)>0$ since there exists only finitely many hyperplanes separating $z$ from $\{x,y\}$. A fortiori, if $\Sigma^*$ denote the chain obtained from $\Sigma$ by replacing $\Sigma'$ with $\Sigma''$, 
$$\ell(\Sigma^*) = \ell(\Sigma) + \ell(\Sigma'')-\ell(\Sigma') \leq \ell(\Sigma)- \epsilon(z).$$
This concludes the proof.
\end{proof}

\begin{proof}[Proof of Proposition \ref{prop:convexdelta1}.]
Let $Y$ be a gated subgraph and $x,y \in Y$ two vertices. We want to prove that, if $z \in X$ is a vertex which does not belong to $Y$, then $z \notin I(x,y)$. This will prove that $I(x,y)$ must be included into $Y$, concluding the proof of the proposition.

\medskip \noindent
Let $\epsilon(z)$ denote the constant given by Lemma \ref{lem:shorteningquantitatif}. Fix a chain $\Sigma_1$ (resp. $\Sigma_2$) between $x$ and $z$ (resp. between $z$ and $y$) such that $|\ell(\Sigma_1)- \delta(x,z)| \leq \epsilon(z)/4$ (resp. $|\ell(\Sigma_2)- \delta(z,y)| \leq \epsilon(z)/4$). A fortiori, the concatenation of $\Sigma_1$ and $\Sigma_2$ defines a chain $\Sigma$ between $x$ and $y$ passing through $z$. According to Lemma \ref{lem:shorteningquantitatif}, there exists a chain $\Sigma^*$ between $x$ and $y$ such that $\ell(\Sigma^*) \leq \ell(\Sigma)- \epsilon(z)$. We have
$$\begin{array}{lcl} \delta(x,y) & \leq & \ell(\Sigma^*) \leq \ell(\Sigma)-\epsilon(z) \\ \\ & \leq & \ell(\Sigma_1)+ \ell(\Sigma_2)- \epsilon(z) \leq \delta(x,z) + \frac{\epsilon(z)}{4} + \delta(z,y) + \frac{\epsilon(z)}{4} - \epsilon(z) \\ \\ & \leq & \delta(x,z)+\delta(z,y) - \frac{\epsilon(z)}{2} < \delta(x,z) + \delta(z,y) \end{array}$$
Therefore, $z$ does not belong to the interval $I(x,y)$. 
\end{proof}

\noindent
Now, fix a quasi-median graph $X$ endowed with a coherent system of geodesic metrics. A natural problem is to determine whether $(X,\delta)$ is geodesic as well, and if so, to determine (at least some of) its geodesics. For this purpose, let us define a \emph{broken geodesic} $\gamma= \gamma_1 \cup \cdots \cup \gamma_n$ between two vertices $x,y \in X$ as the data of 
\begin{itemize}
	\item a sequence of cliques $C_1, \cdots, C_n$ such that $x \in C_1$, $y \in C_n$ and such that $C_i \cap C_{i+1}$ reduces to a single vertex $x_i$ for every $1 \leq i \leq n-1$;
	\item a sequence of paths $\gamma_1 \subset C_1, \ldots, \gamma_n \subset C_n$ such that $\gamma_i$ is a geodesic in $(C_i, \delta_{C_i})$ between $x_i$ and $x_{i+1}$ for every $1 \leq i \leq n-1$, with the convention that $x_0=x$ and $x_n=y$;
\end{itemize}
such that $x_0,x_1, \ldots, x_{n-1},x_n$ defines a geodesic in $(X,d)$ between $x$ and $y$. 

\medskip \noindent
It is worth noticing that any two vertices of $X$ are linked by a broken geodesic. Indeed, fix two vertices $x,y \in X$ and consider some geodesic $x_0=x, \ x_1, \ldots, x_{n-1}, \ x_n=y$ between $x$ and $y$ in $(X,d)$. For every $0 \leq i \leq n-1$, let $C_i$ denote the unique clique of $X$ containing the edge $(x_i,x_{i+1})$ and $\gamma_i$ a geodesic between $x_i$ and $x_{i+1}$ in $(C_i,\delta_{C_i})$. Then
$$\gamma = \gamma_0 \cup \gamma_1 \cdots \cdots \cup \gamma_{n-1} \cup \gamma_n$$
defines a broken geodesic between $x$ and $y$.

\begin{lemma}\label{lem:brokengeod}
A broken geodesic defines a geodesic in $(X, \delta)$. 
\end{lemma}

\begin{proof}
Let $\gamma = \gamma_0 \cup \gamma_1 \cup \cdots \gamma_n$ be a broken geodesic between two vertices $x,y \in X$. For every $0 \leq i \leq n-1$, let $C_i$ denote the clique containing $\gamma_i$ and $x_i,x_{i+1} \in C_i$ the endpoints of $\gamma_i$. By definition,
$$\varphi : x_0=x, \ x_1, \ldots, x_{n-1}, \ x_n=y$$
is a geodesic between $x$ and $y$ in $(X,d)$. We claim that $\gamma$ defines a geodesic in $(X, \delta)$.

\medskip \noindent
Let $a,b \in \gamma$. If there exists some $0 \leq i \leq n$ such that $a,b \in \gamma_i$, then it is clear that $\delta(a,b)=\delta_{C_i}(a,b)$ is equal to the length of the subsegment of $\gamma$ between $a$ and $b$. So let us suppose that there exist $i< j$ such that $a \in \gamma_i \backslash \{ x_{i+1} \}$ and  $b \in \gamma_j \backslash \{ x_{j} \}$. Set
$$\overline{\varphi} : a,\ x_{i+1}, \ x_{i+2}, \ldots, x_{j-1}, \ x_j, \ b.$$
Notice that $\overline{\varphi}$ defines a geodesic in $(X,d)$. Indeed, the clique containing $(a,x_{i+1})$ is $C_i$, the clique containing $(x_j,b)$ is $C_j$, and finally, for every $i+1 \leq r \leq j-1$, the clique containing $(x_r,x_{r+1})$ is $C_r$, so that no hyperplane intersects $\overline{\varphi}$ twice (otherwise, a hyperplane would intersect $\varphi$ twice, which is impossible). As a consequence, the hyperplanes separating $a$ and $b$ are precisely the hyperplanes dual to the cliques mentionned above, hence
$$\delta(a,b)= \delta_{C_i}(a,x_{i+1})+ \sum\limits_{r=i+1}^{j-1} \delta_{C_r}(x_r,x_{r+1})+ \delta_{C_j}(x_j,b).$$
By noticing that this is precisely the length of the subsegment of $\gamma$ between $a$ and $b$, we conclude that $\gamma$ is a geodesic.
\end{proof}

\noindent
Thus, combined with Proposition \ref{prop:convexdelta1}, we deduce from this lemma that:

\begin{cor}\label{cor:systgeodesic}
Let $X$ be quasi-median graph endowed with a coherent system of geodesic metrics. Then $(X,\delta)$ is a geodesic metric space in which the gated subgraphs are convex.
\end{cor}

\begin{remark}
In the statements of Lemma \ref{lem:brokengeod} and Corollary \ref{cor:systgeodesic}, we consider geodesic metric spaces which are either all \emph{discrete} (ie., with geodesics defined on $\mathbb{N}$) or all \emph{continuous} (ie., with geodesics defined on $\mathbb{R}$). More precisely, Lemma \ref{lem:brokengeod} proves that a broken geodesic constructed from discrete (resp. continuous) geodesics defines a discrete (resp. continuous) geodesic, and Corollary \ref{cor:systgeodesic} states that the global metric associated to a system of discrete (resp. continuous) geodesic spaces is discrete (resp. continuous) geodesic. 
\end{remark}

\noindent
In order to motivate the idea that studying $(X,\delta^1)$ may be useful to deduce properties of $(X,\delta^p)$ when the cubical dimension of $X$ is finite, let us notice that, as a consequence of Proposition \ref{prop:comparedelta}, $(X,\delta^1)$ is proper if and only if so is $(X,\delta^p)$. This observation allows us to prove the following criterion, which is however purely technical and will be used only for the proof of Theorem \ref{thm:producingCAT0groups}. 

\begin{lemma}\label{lem:proper}
Let $Y$ be a quasi-median graph of finite cubical dimension endowed with a coherent system of complete and proper CAT(0) metrics. Suppose that $Y$ contains a convex subgraph $X$ satisfying the following conditions:
\begin{itemize}
	\item every vertex of $X$ is contained into finitely many cliques of $X$;
	\item every vertex of $Y$ belongs to a clique containing an edge of $X$;
	\item there exists some constant $K >0$ such that, for any clique $C$ and for any distinct vertices $x,y \in X$, we have $\delta_C(x,y) \geq K$.
\end{itemize}
Then $(Y,\delta^p)$ is proper.
\end{lemma}

\begin{proof}
According to Proposition \ref{prop:comparedelta}, it is sufficient to prove that $(Y, \delta^1)$ is proper. We first want to prove that any ball of $(Y,\delta^1)$ contains only finitely many points of $X$. Notice that, because any vertex of $Y$ is adjacent to some vertex of $X$, we may suppose without loss of generality that our ball is centered at some point of $X$.

\medskip \noindent
Let $B(x,R)$ be a ball of radius $R$ (with respect to $\delta^1$) centered at some point $x \in X$. Suppose that $y \in B(x,R) \cap X$. By considering some broken geodesic between $x$ and $y$, provided by Lemma \ref{lem:brokengeod}, we find a geodesic $x_1, \ldots, x_r$ between $x$ and $y$ in $Y$ (with respect to its combinatorial metric) such that, if $C_i$ denotes the clique containing $x_i$ and $x_{i+1}$ for every $1 \leq i \leq r-1$, then $\delta_{C_i}(x_i,x_{i+1}) \geq K$ (notice that $x_i \in X$ for every $1 \leq i \leq r$ by convexity of $X$). In particular,
$$R \geq \delta^1(x,y) = \sum\limits_{i=1}^{r-1} \delta_{C_i}(x_i,x_{i+1}) \geq Kr,$$
hence $r \leq R/K$. Thinking of $y$ as a variable point, because there exist only finitely many cliques of $X$ containing $x$, one has only finitely many choices on $C_1$, and because $(C_1, \delta_{C_1})$ is proper and that $\delta_{C_1}(x,x_1) \geq K$, one has only finitely many choices on $x_1$ as a point of $C_1$. Therefore, there are finitely many possible choices on $x_1$. Similarly, $x_1$ being fixed, one has only finitely many choices on $C_2$ since there exist only finitely many cliques of $X$ containing $x_1$, and, because $(C_2, \delta_{C_2})$ is proper and that $\delta_{C_2}(x_1,x_2)$, on has only finitely many choices on $x_2$ as a point of $X_2$. Therefore, there are finitely many possible choices on $x_2$. By iterating the argument at most $R/K$ times, we conclude that there exist only finitely many possible $y$'s, ie., $B(x,R) \cap X$ is finite.

\medskip \noindent
Now, we are ready to conclude the proof of our lemma. So, fixing some bounded sequence $(y_n)$ of vertices (with respect to $\delta^1$), we want to find a converging subsequence (with respect to $\delta^1$). For every $n \geq 0$, let $C_n$ be a clique of $X$ such that the unique clique of $Y$ containing $C_n$ contains $y_n$. If $B$ denotes a ball containing all the $y_n$'s, because we proved that $B \cap X$ is finite necessarily $(y_n)$ eventually lies inside some $C_m$ up to taking a subsequence. Since $(C_m,\delta_{C_m})$ is proper, we conclude that this subsequence must contain a converging subsequence, which concludes the proof.
\end{proof}

\noindent
For future use, we record the first claim of the previous proof:

\begin{claim}\label{claim:Xdiscrete}
Under the hypotheses of Lemma \ref{lem:proper}, any ball of $(Y,\delta^p)$ contains finitely many vertices of $X$.
\end{claim}

\subsection{Actions of groups}

\noindent
In this section, we are interested in group actions. First notice that, if a group $G$ acts on a quasi-median graph $X$, endowed with a system of (pseudo-)metrics which is coherent and \emph{$G$-invariant} (ie., $\delta_{gC}(gx,gy)=\delta_C(x,y)$ for every clique $C$ and every vertices $x,y \in C$), then $G$ acts on $(X,\delta)$. Indeed, if we fix two vertices $x,y \in X$ and an element $g \in G$, then, by denoting $J_1, \ldots, J_n$ the hyperplanes separating $x$ and $y$ and by fixing some clique $C_i$ dual to $J_i$ for every $1 \leq i \leq n$, then
$$\begin{array}{lcl} \delta(gx,gy) & = & \sum\limits_{i=1}^n \delta_{gC_i}(\mathrm{proj}_{gC_i}(gx), \mathrm{proj}_{gC_i}(gy)) = \sum\limits_{i=1}^n \delta_{gC_i}( g \cdot \mathrm{proj}_{C_i}(x), g \cdot \mathrm{proj}_{C_i}(y)) \\ \\ & = & \sum\limits_{i=1}^n \delta_{C_i}( \mathrm{proj}_{C_i}(x), \mathrm{proj}_{C_i}(y)) = \delta(x,y) \end{array}$$

\medskip \noindent
Now, we would like to determine when the induced action $G \curvearrowright (X,\delta)$ is metrically proper or geometric. We give some criteria  below. 

\begin{lemma}\label{lem:whencobounded}
Let $G$ be a group acting on a quasi-median graph $X$ endowed with a $G$-invariant coherent system of metrics. Suppose that
\begin{itemize}
	\item $X$ contains finitely many $G$-orbits of cliques;
	\item for every clique $C$, the action $\mathrm{stab}(C) \curvearrowright (C,\delta_C)$ is cobounded.
\end{itemize}
Then $G$ acts coboundedly on $(X,\delta)$.
\end{lemma}

\begin{proof}
Let $C_1, \ldots, C_n$ be a collection of cliques so that any clique of $X$ is a translate of some $C_i$. For every $1 \leq i \leq n$, fix a vertex $x_i \in C_i$ and let $A_i$ be the diameter of a fundamental domain containing $x_i$ for the action $\mathrm{stab}(C_i) \curvearrowright (C_i,\delta_{C_i})$. Set $A= \max\limits_{1 \leq i \leq n} A_i$ and $B=\max\limits_{2 \leq i \leq n} \delta(x_1,x_i)$. For every $x \in X$, there exists some $g \in G$ and $1 \leq i \leq n$ such that $gx \in C_i$. Then, there exists some $h \in \mathrm{stab}(C_i)$ such that $\delta(x_i,hgx) \leq A$. Therefore,
$$\delta(x_1,x) \leq \delta(x_1,x_i) + \delta(x_i,hgx) \leq A+B.$$
As a consequence, the ball of radius $A+B$ centered at $x_1$ defines a fundamental domain for the action $G \curvearrowright (X, \delta)$. A fortiori, this action is cobounded. 
\end{proof}

\begin{lemma}\label{lem:whenmetricallyproper}
Let $G$ be a group acting on a quasi-median graph $X$ endowed with a $G$-invariant coherent system of metrics. Suppose that
\begin{itemize}
	\item any vertex of $X$ belongs to finitely many cliques;
	\item vertex-stabilisers are finite;
	\item for every clique $C$, the space $(C, \delta_C)$ is locally finite;
	\item suppose that the metrics $\delta_C$ are \emph{uniformly discrete}, ie., there exists some constant $K>0$ such that, for any clique $C$ and any points $x,y \in C$, $\delta_C(x,y) \geq K$.
\end{itemize}
Then $G$ acts metrically properly on $(X,\delta)$.
\end{lemma}

\noindent
Essentially, this result will be a consequence of the following lemma:

\begin{lemma}\label{lem:whenlocallyfinite}
Let $X$ be a quasi-median graph endowed with a coherent system of metrics. Suppose that
\begin{itemize}
	\item any vertex of $X$ belongs to finitely many cliques;
	\item for every clique $C$, the space $(C, \delta_C)$ is locally finite;
	\item suppose that the metrics $\delta_C$ are uniformly discrete, ie., there exists some constant $K>0$ such that, for any clique $C$ and any points $x \neq y \in C$, $\delta_C(x,y) \geq K$.
\end{itemize}
Then $(X, \delta)$ is locally finite.
\end{lemma}

\begin{proof}
Up to rescaling all the $\delta_C$'s, we may suppose without loss of generality that $K=1$. As a consequence, $d \leq K \cdot \delta = \delta$.

\medskip \noindent
We want to prove that, for every $x \in X$ and every $R \geq 0$, the ball $B_{\delta}(x,R)$ is finite; noticing that
$$B_{\delta}(x,R+1) \subset \bigcup\limits_{y \in B_{\delta}(x,R)} B_{\delta}(y,1),$$
it is sufficient to prove that $B_{\delta}(x,1)$ is finite for every $x \in X$. 

\medskip \noindent
Notice that, for every $y \in B_{\delta}(x,1)$, we have $d(x,y) \leq \delta(x,y) \leq 1$, ie., $x$ and $y$ belong to a common clique. Let $C_1, \ldots, C_n$ denote the cliques containing $x$. Notice that, for every $1 \leq i \leq n$ and every $y \in C_i$, we have $\delta(x,y)= \delta_{C_i} (x,y)$. Consequently,
$$B_{\delta}(x,1) \subset \bigcup\limits_{i=1}^n B_{\delta_{C_i}}(x,1),$$
so that we deduce that $B_{\delta}(x,1)$ is finite since the balls $B_{\delta_{C_i}}(x,1)$ are themselves finite by local finiteness of the spaces $(C_i, \delta_{C_i})$. 
\end{proof}

\begin{remark}\label{remark:locallyfinite}
In the statement of Lemma \ref{lem:whenlocallyfinite}, the second condition is of course necessary. However, the first and third conditions are not necessary. Nevertheless, the three conditions turn out to define a necessary and sufficient criterion for the local finiteness of $(X,\delta)$ if a group acts on $X$ with finitely many orbits of cliques. 
\end{remark}

\begin{proof}[Proof of Lemma \ref{lem:whenmetricallyproper}.]
Let $x \in X$ and $R \geq 0$. For every $y \in B_{\delta}(x,R)$, fix some $g_y \in G$ such that $g_y \cdot x=y$; if no such element of $G$ exists, set $g_y=1$. Then, set $F= \{ g_y \mid y \in B_{\delta}(x,R) \}$. Notice that
$$\# F \leq \# \mathcal{B}_{\delta}(x,R) <+ \infty$$
according to the previous lemma. Now, if $g \in G$ satisfies $d_{\delta}(x,gx) \leq R$, then there exists some $f \in F$ such that $gx=fx$, hence $g \in F \cdot \mathrm{stab}(x)$. Thus, we have proved that
$$\{ g \in G \mid d_{\delta}(x,gx) \leq R \} \subset F \cdot \mathrm{stab}(x).$$
On the other hand, we know that $F$ is finite and $\mathrm{stab}(x)$ is finite by hypothesis. Therefore, $G$ acts metrically properly on $(X, \delta)$.
\end{proof}

\noindent
By combining Lemma \ref{lem:whencobounded} and Lemma \ref{lem:whenmetricallyproper}, we deduce the following statement.

\begin{cor}\label{cor:mwhengeometric}
Let $G$ be a group acting on a quasi-median graph $X$ endowed with a $G$-invariant coherent system of discrete metrics. Suppose that:
\begin{itemize}
	\item any vertex of $X$ belongs to finitely many cliques;
	\item vertex-stabilisers are finite;
	\item $X$ contains finitely many $G$-orbits of cliques;
	\item for every clique $C$, the action $\mathrm{stab}(C) \curvearrowright (C,\delta_C)$ is geometric.
\end{itemize}
Then $G$ acts geometrically on $(X,\delta)$. 
\end{cor}

\subsection{$\ell^p$-compression}\label{section:compression}

\noindent
Roughly speaking, the $\ell^p$-compression of a metric space $X$ is a real number in $[0,1]$ which quantifies how much it is necessary to deform the geometry of $X$ in order to embed it in some $L^p$-space; a precise definition is given below. In this section, our goal is to find a lower bound on the $\ell^p$-compression of the global metric associated to some coherent system of metrics on a quasi-median graph with respect to the $\ell^p$-compressions of the local metrics. 

\begin{definition}
Let $f : X \to Y$ be a Lipschitz map between two metric spaces. The \emph{compression} of $f$, denoted by $\alpha(f)$, is the supremum of the $\alpha$'s such that there exists some constant $C >0$ so that the inequality 
$$C \cdot d(x,y)^{\alpha} \leq d(f(x),f(y))$$
holds for every $x,y \in X$. Given a metric space $X$ and some $p \geq 1$, the \emph{$\ell^p$-compression} of $X$ is the supremum of the $\alpha(f)$'s where $f$ is a Lipschitz map from $X$ to a $L^p$-space. In particular, if $G$ is a finitely generated group, its $\ell^p$-compression\index{Compression ($\ell^p$-)} is defined as the $\ell^p$-compression of $G$ endowed with the word metric associated to a finite generating set; notice that the choice of this generating set does not modify the compression since any two word metrics (with respect to two finite generating sets) are Lipschitz-equivalent. 
\end{definition}

\noindent
It is worth noticing that, if $X$ is a uniformly discrete metric space (ie., there exists some $B>0$ such that $d(x,y) \geq B$ for every distinct points $x,y \in X$), it is possible to define alternatively the $\ell^p$-compression of $X$ as the supremum of the $\alpha$'s such that there exist a Lipschitz map $f : X \to L$ to some $L^p$-space such that 
$$C\cdot d(x,y)^{\alpha}-D \leq d(f(x),f(y))$$
for some $C,D>0$ and for all $x,y \in X$; see for instance \cite[Lemma 2.1]{Dreesen1}. As a consequence, the $\ell^p$-compression of a finitely generated group turns out to be a quasi-isometric invariant.   

\begin{definition}
Let $X$ be a quasi-median graph endowed with a coherent system of metrics. A \emph{system of Lipschitz $\ell^p$-maps} is the data, for each clique $C$ of $X$, of one Lipschitz map $f_C : (C,\delta_C) \to L_C$ from the clique $C$ to some $L^p$-space $L_C$. Such a system is \emph{coherent} if
\begin{itemize}
	\item $L_C=L_{C'}$ if $C$ and $C$' are two cliques dual to the same hyperplane;
	\item for every two cliques $C,C'$ of $X$ dual to the same hyperplane and for every two vertices $x,y \in C$, $$f_C(x,y)= f_{C'}( t_{C \to C'}(x), t_{C\to C'}(y))$$ where $t_{C \to C'} : C \to C'$ is the canonical bijection.
\end{itemize}
\end{definition}

\noindent
The main result of this section is the following.

\begin{prop}\label{prop:compression}
Let $X$ be a quasi-median graph endowed with a coherent system of uniformly discrete metrics $\{ (C,\delta_C) \mid C \ \text{clique} \}$. Suppose that the collection $\{ (C, \delta_C) \mid \text{$C$ clique} \}$ contains finitely many isometry classes. Then
$$\alpha_p(X,\delta) \geq \min \left( \frac{1}{p}, \inf\limits_{C ~ \text{clique}} \alpha_p(C,\delta_C)\right).$$
\end{prop}

\begin{proof}
First of all, notice that rescaling all the $\delta_C$'s by a common positive constant does not modify the compressions $\alpha_p(X,\delta)$ and $\alpha_p(C,\delta_C)$. Therefore, since our collection of metrics is uniformly discrete, we can suppose without loss of generality that $\delta_C(x,y) \geq 1$ for every clique $C$ and every distinct vertices $x,y \in C$. Next, if $\min\limits_{\text{$C$ clique}} \alpha_p(C,\delta_C)=0$, there is nothing to prove, so we suppose that $\min\limits_{\text{$C$ clique}} \alpha_p(C,\delta_C)>0$.

\medskip \noindent
Fix some collection $\mathcal{Q}$ of cliques of $X$ such that each hyperplane contains exactly one clique of $\mathcal{Q}$, and some $0 < \epsilon < \min\limits_{\text{$C$ clique}} \alpha_p(C,\delta_C)$. For every clique $C \in \mathcal{Q}$, fix a Lipschitz map $f_C$ from $(C, \delta_C)$ to some $L^p$-space $L_C$ such that $\alpha(f_C) \geq \alpha_p(C,\delta_C) - \epsilon$. Because $\{ (C, \delta_C) \mid C \in \mathcal{Q} \}$ contains only finitely many isometry classes, we can choose the maps $f_C$ so that, for every positive $\eta < \min\limits_{C \in \mathcal{Q}} \alpha_p(C,\delta_C)$, there exist constants $A_{\eta}, B>0$ (which do not depend on $C$) such that
$$A_{\eta} \cdot \delta_C(x,y)^{\alpha(f_C)-\eta} \leq \| f_C(x)- f_C(y) \| \leq B \cdot \delta_C(x,y)$$
for every $C \in \mathcal{Q}$ and every $x,y \in C$. Now, if $C$ is an arbitrary clique of $X$, the hyperplane dual to $C$ must contain a unique clique $Q$ of $\mathcal{Q}$. We set $L_C=L_Q$ and $f_C = f_Q \circ t_{C \to Q}$, in order to define a system of Lipschitz $\ell^p$-maps, which is coherent by construction; notice that, because our system of metrics is coherent, the canonical bijection $t_{C \to Q} : (C,\delta_C) \to (Q,\delta_Q)$ defines an isometry, so that, for every positive $\eta < \min\limits_{C \in \mathcal{Q}} \alpha_p(C,\delta_C)$, the inequalities
$$A_{\eta} \cdot \delta_C(x,y)^{\alpha(f_Q)-\eta} \leq \| f_C(x)- f_C(y) \| \leq B \cdot \delta_C(x,y)$$
hold for every clique $C$ labelled by $Q \in \mathcal{Q}$ and every vertices $x,y \in C$. In particular,
$$\alpha(f_C) = \alpha(f_Q) \geq \alpha_p(Q,\delta_Q)- \epsilon = \alpha_p(C, \delta_C)- \epsilon.$$
For convenience, let $p_Q : X \to Q$ denote the projection onto $Q$ for every $Q \in \mathcal{Q}$. Define
$$f : \left\{ \begin{array}{ccc} X & \to & L= \bigoplus\limits_{C \in \mathcal{Q}}^{\ell^p} L_C \\ x & \mapsto & ( f_C(p_C(x))-f_C(p_C(x_0)))_C \end{array} \right.,$$
where $x_0 \in X$ is a basepoint we fix. 

\begin{claim}\label{claim:compressionf}
$f$ is well-defined, Lipschitz and satisfies
$$\alpha(f) \geq \min \left( \frac{1}{p}, \min\limits_{C ~ \text{clique}} \alpha(f_C)\right).$$
\end{claim}

\noindent
First, let us verify that $f$ is well-defined. For every vertex $x \in X$, only finitely many hyperplanes separate $x$ and $x_0$, so that $p_C(x)$ and $p_C(x_0)$ differ only for finitely many $C \in \mathcal{Q}$; a fortiori, $f(x)$ belongs to $L$. Next, in order to justify that $f$ is a Lipschitz map, notice that, for every vertices $x,y \in X$, we have
$$\begin{array}{lcl} \| f(x)- f(y) \| & = & \displaystyle  \left( \sum\limits_{C \in \mathcal{Q}} \| f_C(p_C(x)) - f_C(p_C(y)) \|^p \right)^{1/p} \\ \\ & \leq & \displaystyle  B \cdot \left( \sum\limits_{C \in \mathcal{Q}} \delta_C(p_C(x),p_C(y))^{p} \right)^{1/p} \\ \\ & \leq & \displaystyle B \cdot \sum\limits_{C \in \mathcal{Q}} \delta_C(p_C(x),p_C(y)) = B \cdot \delta(x,y) \end{array}$$
The final step is to understand the compression of $f$. For convenience, we set $c=  \min \left( \frac{1}{p}, \min\limits_{C ~ \text{clique}} \alpha(f_C) - \eta\right)$ for some fixed positive $\eta < \min\limits_{\text{$C$ clique}} \alpha_C$. For every vertices $x,y \in X$, 
$$\begin{array}{lcl} \| f(x)-f(y) \|^{1/c} & = & \displaystyle \left( \sum\limits_{C \in \mathcal{Q}} \| f_C(p_C(x))- f_C(p_C(y)) \|^p \right)^{1/pc} \\ \\ & \geq & \displaystyle  A_{\eta}^{1/c} \cdot \left( \sum\limits_{C \in \mathcal{Q}} \delta_C(p_C(x), p_C(y))^{p(\alpha(f_C)-\eta)} \right)^{1/pc} \\ \\ & \geq & \displaystyle A_{\eta}^{1/c} \cdot  \left( \sum\limits_{C \in \mathcal{Q}} \delta_C(p_C(x), p_C(y))^{pc} \right)^{1/pc} \\ \\ & \geq & \displaystyle A_{\eta}^{1/c} \cdot \sum\limits_{Q \in \mathcal{Q}} \delta_C(p_C(x),p_C(y)) = A_{\eta}^{1/c} \cdot \delta(x,y) \end{array}$$
Thus, we have shown that $\alpha(f) \geq  \min \left( \frac{1}{p}, \min\limits_{C ~ \text{clique}} \alpha(f_C) - \eta\right)$, and because the inequality holds for every sufficently small $\eta>0$, this concludes the proof of our claim. Because 
$$\alpha_p(X, \delta) \geq \alpha(f) \geq \min \left( \frac{1}{p}, \min\limits_{C~\text{clique}} \alpha(f_C) \right) \geq \min \left( \frac{1}{p}, \min\limits_{C~\text{clique}} \alpha(C,\delta_C) - \epsilon \right)$$
is true for every $\epsilon>0$, our proposition follows.
\end{proof}

\section{Cubulating quasi-median graphs}\label{section:cubulatingQM}

\noindent
In this section, we show that, if each clique $C$ of a quasi-median graph $X$ is endowed with a collection of \emph{walls} $\mathcal{W}(C)$, in such a way that the family of all these \emph{wallspaces} is \emph{coherent}, then there exists a global collection of walls $\mathcal{HW}$ extending them. Moreover, if a group $G$ acts on $X$ and if our collection of wallspaces is \emph{$G$-invariant}, then the action $G \curvearrowright X$ induces an action $G \curvearrowright (X, \mathcal{HW})$. (The existence such collection of wallspaces will be studied in Section \ref{section:topicalactionsI}.) The main point is that the global combinatorics of $(X,\mathcal{HW})$ reduces to the local combinatorics of the $(C, \mathcal{W}(C))$'s. In the final subsection, we mention how this formalism can be adapted to some generalisations of spaces with walls, namely \emph{measured wallspaces} and \emph{spaces with labelled partitions}. 

\medskip \noindent
The rest of this introduction is used to fixed the definitions and notation related to spaces with walls.

\medskip \noindent 
Given a set $X$, a \emph{wall} is a partition $W=\{ h,h^c \}$ where $h \notin \{ \emptyset, X \}$; we will refer to $h$ and its complementary $h^c$ as the \emph{halfspaces} of $W$. If $x,y \in X$ are two points, we say that $W$ \emph{separates} $x$ and $y$ if they belong to distinct halfspaces of $W$. Finally, given a collection of walls $\mathcal{W}$, we say that $(X, \mathcal{W})$ is a \emph{space with walls}\index{Spaces with walls} if any two points of $X$ are separated by only finitely many walls. Notice that this condition allows us to define a pseudo-distance $d_{\mathcal{W}}$ on $X$ by counting the number of walls separating two given points. 

\begin{figure}
\begin{center}
\includegraphics[scale=0.6]{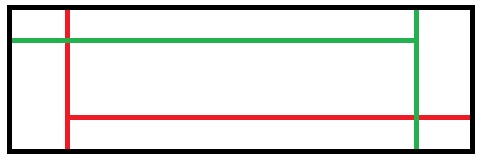}
\end{center}
\caption{Transverse partitions delimiting two disjoint sectors.}
\label{figure6}
\end{figure}

\medskip \noindent
In particular, spaces with walls are examples of spaces with partitions as defined in Section \ref{section:spaceswithpartitions} (see Remark \ref{rem:wallspaces}). In our definition of spaces with walls, we allow \emph{duplicates}, ie., two walls may be \emph{indistinguishable}\index{Indistinguishable walls} in the sense that they induce the same partition on $X$; compare with \cite[Remark 2.2]{HruskaWise}. We cubulate spaces with walls by quasi-cubulating them as spaces with partitions, as explained in Section \ref{section:spaceswithpartitions}. The main difference between our cubulation and the cubulation for instance described in \cite{HruskaWise} is that we define when two walls are nested and then we say that two walls which are not nested are transverse, whereas usually one says that two walls $\{h_1,h_1^c \}$ and $\{h_2,h_2^c \}$ are transverse if the four intersections
$$h_1 \cap h_2, \ h_1 \cap h_2^c, \ h_1^c \cap h_2, \ h_1^c \cap h_2^c$$
are non empty and then one says that two walls which are not transverse are nested. Observe that if $\mathcal{P}_1$ and $\mathcal{P}_2$ are two partitions such that $A_1 \cap A_2 \neq \emptyset$ for every $A_1 \in \mathcal{P}_1$ and $A_2 \in \mathcal{P}_2$, then $\mathcal{P}_1$ and $\mathcal{P}_2$ are transverse, so that the usual definition of transversality implies ours. However, the converse does not hold in general, as illustrated by Figure \ref{figure6}.

\medskip \noindent
Of course, given a space with walls $(X, \mathcal{W})$ and the CAT(0) cube complex $C(X, \mathcal{W})$ obtained by cubulating it, we are interested in reading the geometry of the CAT(0) cube complex $C(X, \mathcal{W})$ directly from $(X, \mathcal{W})$. For this purpose, we will essentially use Theorem \ref{thm:quasicubulation} in order to link the combinatorics of the hyperplanes of $C(X, \mathcal{W})$ with the combinatorics of the walls of $(X, \mathcal{W})$. For instance, we define the \emph{dimension} of $(X, \mathcal{W})$ as the maximal cardinality of a collection of pairwise transverse walls, and we show:

\begin{lemma}\label{lem:dimspacewalls}
Let $(X, \mathcal{W})$ be a space with walls and let $C(X, \mathcal{W})$ denote the CAT(0) cube complex obtained by cubulating $(X, \mathcal{W})$. Then $\dim(X, \mathcal{W})= \dim C(X, \mathcal{W})$. 
\end{lemma}

\begin{proof}
The dimension of $C(X, \mathcal{W})$ is equal to the maximal cardinality of a collection of pairwise transverse hyperplanes of $C(X, \mathcal{W})$ (see for instance Proposition \ref{prop:transversehypcube}). Therefore, the bijection between the walls of $\mathcal{W}$ and the hyperplanes of $C(X, \mathcal{W})$ given by Theorem \ref{thm:quasicubulation} provides the conclusion.
\end{proof}

\noindent
In the sequel, given a group $G$ acting on a space with walls $(X, \mathcal{W})$, we will say that the action $G \curvearrowright (X, \mathcal{W})$ satisfies some property $\mathcal{P}$ if the associated action of $G$ on the CAT(0) cube complex $C(X,\mathcal{W})$ obtained by cubulating $(X, \mathcal{W})$ satisfies $\mathcal{P}$ as well. In particular,

\begin{lemma}
A group $G$ acts metrically properly on a space with walls $(X,\mathcal{W})$ if and only if it acts metrically properly on the pseudo-metric space $(X,d_{\mathcal{W}})$.
\end{lemma}

\begin{proof}
Let $C(X, \mathcal{W})$ denote the CAT(0) cube complex obtained by cubulating $(X, \mathcal{W})$. Fix some vertex $x \in X$ and let $\sigma_x$ denote the principal orientation associated. In particular, $\sigma_x$ is a vertex of $C(X, \mathcal{W})$. Notice that, using Lemma \ref{lem:distspacepartitions}, we know that
$$d_{C(X, \mathcal{W})}(\sigma_x, g \cdot \sigma_x)= d_{C(X, \mathcal{W})}(\sigma_x, \sigma_{g \cdot x})= d_{\mathcal{W}}(x,g \cdot x)$$
for every $g \in G$. As a consequence, the action $G \curvearrowright C(X, \mathcal{W})$ is metrically proper if and only if the action $G \curvearrowright (X,d_{\mathcal{W}})$ is metrically proper as well.
\end{proof}

\begin{lemma}
Let $(X, \mathcal{W})$ be a finite dimensional space with walls. A group $G$ acts cocompactly on $(X, \mathcal{W})$ if and only if $X$ contains finitely many $G$-orbits of maximal collections of pairwise transverse walls.
\end{lemma}

\begin{proof}
Let $C(X, \mathcal{W})$ denote the CAT(0) cube complex obtained by cubulating $(X, \mathcal{W})$. Notice that, according to Lemma \ref{lem:dimspacewalls}, $C(X, \mathcal{W})$ is finite dimensional. As a consequence, Proposition \ref{prop:maxprism} applies. Therefore, $G$ acts cocompactly on $C(X, \mathcal{W})$ if and only if $C(X, \mathcal{W})$ has finitely many $G$-orbits of maximal cubes if and only if $C(X, \mathcal{W})$ contains finitely many $G$-orbits of maximal collections of pairwise transverse hyperplanes. According to Theorem \ref{thm:quasicubulation}, the latter condition is equivalent to the following statement: $X$ contains finitely many $G$-orbits of maximal collections of pairwise transverse walls.
\end{proof}

\noindent
Such a result can also be formulated for \emph{special} actions, as defined by Haglund and Wise in \cite{MR2377497}. Let us adapt the definition for popsets (and, in particular, for spaces with partitions, including quasi-median graphs).

\begin{definition}
Let $G$ be a group acting on a popset $(X,<, \mathfrak{P})$. The action is \emph{special} if none of the following \emph{pathological configurations} happen:
\begin{itemize}
	\item there exist a partition $\mathcal{P} \in \mathfrak{P}$ and an element $g \in G$ such that $\mathcal{P}$ and $g \mathcal{P}$ are transverse;
	\item there exist a partition $\mathcal{P} \in \mathfrak{P}$ and an element $g \in G$ such that $\mathcal{P}$ and $g \mathcal{P}$ are tangent;
	\item there exist three partitions $\mathcal{P}_1, \mathcal{P}_2, \mathcal{P}_3 \in \mathfrak{P}$ and an element $g \in G$ such that $\mathcal{P}_1$ is transverse to $\mathcal{P}_2$, $\mathcal{P}_3$ is tangent to $\mathcal{P}_2$, and $g \mathcal{P}_1= \mathcal{P}_3$. 
\end{itemize}
\end{definition}

\begin{remark}
We emphasize that a group acting (geometrically and) specially on some CAT(0) cube complex may not be a \emph{special group} as defined in \cite{MR2377497}. For instance, it is proved in \cite{HsuWiseNonRF} that there exists a nonpositively-curved square of finite groups whose fundamental group is not residually finite. Such a group acts geometrically and specially on some CAT(0) square complex but is not (virtually) special. Nevertheless, a group acting geometrically and virtually specially on a CAT(0) cube complex turns out to be virtually special if it is residually finite, or more generally if it is virtually torsion-free. 
\end{remark}

\begin{remark}
Usually, when defining special actions on CAT(0) cube complexes, the following additional condition is also required: if $J$ is a hyperplane delimiting two halfspaces $J^+$ and $J^-$, then $gJ^+ \neq J^-$ for every $g \in \mathrm{stab}(J)$. This condition does not appear in the definition of specialness given in \cite[Definition 3.2]{MR2377497}. Nevertheless, it is worth noticing that, when our group turns out to be residually finite or virtually torsion-free (ie., when our group turns out to be special, which is the case which interests us), up to taking a finite-index subgroup this condition is always satisfied \cite[Proposition 3.10]{MR2377497}. 
\end{remark}

\noindent
It clearly follows from Theorem \ref{thm:quasicubulation} that the action on a popset is special if and only if the action on the associated quasi-median graph is special as well. In particular, acting specially on a space with walls implies acting specially on the dual CAT(0) cube complex.

\begin{definition}
Let $G$ be a group acting on a space with walls $(X, \mathcal{W})$. The \emph{obstruction to the specialness} $\mathrm{Obs}(G \curvearrowright (X, \mathcal{W}))$ is the set of the elements of $G$ which appear in a pathological configuration.
\end{definition}

\noindent
We conclude this introduction by the following easy observation, which will be useful later.

\begin{fact}\label{fact:obs}
Let $G$ be a group acting on a space with walls $(X, \mathcal{W})$. Suppose that $G$ contains a subgroup $H$ such that the induced action $H \curvearrowright (X, \mathcal{W})$ is special. If $F$ is a set of representatives of $G/H$ different from $H$, then
$$\mathrm{Obs}(G \curvearrowright (X, \mathcal{W})) \subset FH.$$
\end{fact}

\begin{proof}
Since $G= FH \sqcup H$ and $\mathrm{Obs}(G \curvearrowright (X, \mathcal{W})) \cap H = \emptyset$, because the action $H \curvearrowright (X, \mathcal{W})$ is special, the conclusion follows.
\end{proof}

\subsection{Extending walls}\label{section:extendingwall}

\noindent
Let $X$ be a quasi-median graph. If $W= \{ h,h^c \}$ is a wall defined on some clique $C$, we extend it as a wall on $X$ by
$$\overline{W}= \left\{ \bigcup\limits_{x \in h} [C,x], \bigcup\limits_{x \in h^c} [C,x] \right\} = \left\{ \mathrm{proj}_C^{-1}(h), \mathrm{proj}_C^{-1}(h^c) \right\}.$$
Recall that, if $x \in C$, $[C,x]$ denotes the sector delimited by the hyperplane dual to $C$ which contains $x$. More generally, if $\mathcal{W}(C)$ is a collection of walls of $C$, we introduce
$$\overline{\mathcal{W}}(C)= \left\{ \overline{W} \mid W \in \mathcal{W}(C) \right\}.$$
If $C'$ is a second clique such that $C$ and $C'$ are dual to the same hyperplane, we are able to transfer a collection of walls of $C$ to $C'$ by defining
$$\mathcal{W}(C \to C')= \left\{ t(W)=\{ t(h),t(h^c) \} \mid W= \{ h,h^c \} \in \mathcal{W}(C) \right\},$$
where $t=t_{C \to C'}$ is the canonical bijection from $C$ to $C'$ as defined in Section \ref{section:canonical}. It is worth noticing that this operation does not modify the extended collection of walls:

\begin{fact}\label{fact:transfert}
For any two cliques dual to the same hyperplane, the equality 
$$\overline{\mathcal{W}}(C\to C')= \overline{\mathcal{W}}(C)$$
holds, where $\overline{\mathcal{W}}(C \to C') = \{ \overline{W} \mid W \in \mathcal{W}(C \to C') \}$. 
\end{fact}

\begin{proof}
Let $M \in \overline{\mathcal{W}}(C)$ be a wall. So there exists some $W= \{h,h^c \} \in \mathcal{W}(C)$ such that $M= \overline{W}$. We deduce from Lemma \ref{lem:projectionandtransfer} that
$$M= \left\{ \mathrm{proj}_C^{-1}(h), \mathrm{proj}_C^{-1}(h^c) \right\} = \left\{ \mathrm{proj}_{C'}^{-1}(t(h)), \mathrm{proj}_{C'}^{-1}(t(h^c)) \right\} \in \overline{\mathcal{W}}(C \to C'),$$
since $\{t(h),t(h^c) \} \in \mathcal{W}(C \to C')$. Conversely, consider a wall $M \in \overline{\mathcal{W}}(C \to C')$. So there exists some $W = \{ t(h),t(h^c) \} \in \mathcal{W}(C \to C')$, where $\{ h,h^c \} \in \mathcal{W}(C)$, such that $M= \overline{W}$. Once again, we deduce from Lemma \ref{lem:projectionandtransfer} that
$$M = \left\{ \mathrm{proj}_{C'}^{-1}(t(h)), \mathrm{proj}_{C'}^{-1}(t(h^c)) \right\} = \left\{ \mathrm{proj}_C(h), \mathrm{proj}_C(h^c) \right\} \in \overline{\mathcal{W}}(C),$$
since $\{ h,h^c \} \in \mathcal{W}(C)$. This concludes the proof.
\end{proof}

\begin{definition}
Let $X$ be a quasi-median graph. A \emph{system of wallspaces}\index{Systems of wallspaces} is the data of a set of walls $\mathcal{W}(C)$ for each clique $C$, such that $(C, \mathcal{W}(C))$ defines a space with walls. It is \emph{coherent} if $\mathcal{W}(C')= \mathcal{W}(C \to C')$ for every cliques $C,C'$ dual to the same hyperplane.
\end{definition}

\noindent
From now on, we fix a quasi-median graph $X$ endowed with a coherent system of wallspaces. Notice that, if $C$ and $C'$ are two cliques dual to the same hyperplane, then, because our system of wallspaces is coherent and thanks to Fact \ref{fact:transfert}, we have
$$\overline{\mathcal{W}}(C')= \overline{\mathcal{W}}(C \to C') = \overline{\mathcal{W}}(C).$$
As a consequence, it makes sense to set, for every hyperplane $J$ of $X$, $\mathcal{W}(J)= \overline{\mathcal{W}}(C)$, where $C$ is any clique dual to $J$. Thus, the coherent system of wallspaces of $X$ defines a collection of walls, which we will refer to as \emph{hyperplane-walls}\index{Hyperplane-walls}, denoted by
$$\mathcal{HW}=  \bigcup\limits_{\text{$J$ hyperplane of $X$}} \overline{\mathcal{W}}(J).$$
Notice that every clique is endowed with the pseudo-metric associated to its collection of walls. Our first result essentially states that the pseudo-distance associated to $\mathcal{HW}$ coincides with the global pseudo-metric obtained from this system of pseudo-metrics as explained in Section \ref{section:extendingmetrics}.

\begin{lemma}\label{lem:distanceHW}
Let $x,y \in X$ be two vertices. Let $J_1, \ldots, J_n$ denote the hyperplanes separating $x$ and $y$, and, for every $1 \leq i \leq n$, fix a clique $C_i$ dual to $J_i$. If $p_i : X \to C_i$ denotes the projection onto $C_i$ for every $1 \leq i \leq n$, then 
$$d_{\mathcal{HW}}(x,y)= \sum\limits_{i=1}^n d_{\mathcal{W}(C_i)}(p_i(x),p_i(y)).$$
\end{lemma}

\begin{proof}
For every $1 \leq i \leq n$, let $q_i : X \to N(J_i)$ denote the projection onto $N(J_i)$. Because a hyperplane-wall is a union of sectors delimited by some fixed hyperplane, the underlying hyperplane of a hyperplane-wall separating $x$ and $y$ must separate $x$ and $y$ as well, hence
$$d_{\mathcal{HW}}(x,y)= \sum\limits_{J \ \text{hyperplane}} d_{\mathcal{W}(J)}(x,y)= \sum\limits_{i=1}^n d_{\mathcal{W}(J_i)}(x,y).$$
Moreover, for every $1 \leq i \leq n$, the quantity $d_{\mathcal{W}(J_i)}(x,y)$ depends only on the sectors containing $x$ and $y$, hence $d_{\mathcal{W}(J_i)}(x,y)=d_{\mathcal{W}(J_i)}(q_i(x),q_i(y))$. Noticing that $p_i(q_i(x))=p_i(x)$ and similarly $p_i(q_i(y))=p_i(y)$, according to Corollary \ref{cor:projnested}, it is sufficient to prove the following fact in order to deduce that $d_{\mathcal{W}(J_i)}(x,y)= d_{\mathcal{W}(C_i)}(p_i(x),p_i(y)),$
which concludes the proof of our lemma.

\begin{fact}\label{fact:pseudodistcoh}
Let $J$ be a hyperplane, $x,y \in N(J)$ two vertices and $C$ a clique dual to $J$. If $p : X \to C$ denotes the projection onto $C$, then $d_{\mathcal{W}(J)}(x,y)=d_{\mathcal{W}(C)}(p(x),p(y)).$
\end{fact}

\noindent
Let $M \in \mathcal{W}(J)$. Because $\mathcal{W}(J)= \overline{\mathcal{W}}(C)$, there exists a wall $W= \{ h,h^c\} \in \mathcal{W}(C)$ such that $M= \overline{W} = \left\{ p^{-1}(h), p^{-1}(h^c) \right\}$. Then $M$ separates $x$ and $y$ if and only if $W$ separates $p(x)$ and $p(y)$. Therefore $d_{\mathcal{W}(J)}(x,y)=d_{\mathcal{W}(C)}(p(x),p(y)).$
\end{proof}

\noindent
The second step is to understand when two hyperplane-walls are transverse or tangent. This is the aim of our next two lemmas.

\begin{lemma}\label{lem:transverseHW}
Let $M_1 \in \mathcal{W}(J_1)$ and $M_2 \in \mathcal{W}(J_2)$ be two hyperplane-walls. Then $M_1$ and $M_2$ are transverse if and only if either $J_1$ and $J_2$ are transverse, or $J_1=J_2$ and there exist two transverse walls $W_1,W_2\in \mathcal{W}(C)$, where $C$ is a clique dual to $J_1=J_2$, such that $M_1= \overline{W_1}$ and $M_2= \overline{W_2}$. 
\end{lemma}

\begin{proof}
If $J_1$ and $J_2$ are nested, then, because the halfspaces of $M_1$ and $M_2$ are unions of sectors delimited by $J_1$ and $J_2$ respectively, some halfspace delimited by $M_1$ or $M_2$ must be strictly included into another halfspace delimited by the other wall. A fortiori, $M_1$ and $M_2$ are nested.

\medskip \noindent
Suppose that $J_1$ and $J_2$ are transverse. According to Proposition \ref{prop:transversehypcube}, there exists a prism $P$, which is a product of two cliques $C_1, C_2$ dual to $J_1,J_2$ respectively. Let $u$ be the single vertex of the intersection $C_1 \cap C_2$, and, for $i=1,2$, let $u_i \in C_i$ be a vertex such that $M_i$ separates $u$ and $u_i$. Let $w$ be the vertex $(u_1,u_2)$ of $P= C_1 \times C_2$. Moreover, $M_1$ separates $u_2$ and $w$, because $u_1$ and $w$ belong to the same sector delimited by $J_1$; similarly, $M_2$ separates $u_1$ and $w$. Thus, if $M_i^+$ denotes the halfspace delimited by $M_i$ which contains $u$, and $M_i^-$ its complement, then 
$$u \in M_1^+ \cap M_2^+, \ u_1 \in M_1^- \cap M_2^+, \ u_2 \in M_1^+ \cap M_2^-, \ w \in M_1^- \cap M_2^-.$$
Therefore, $M_1$ and $M_2$ are transverse.

\medskip \noindent
Finally, suppose that $J_1=J_2$. So, given a clique $C$ dual to $J_1=J_2$, there exist two walls $W_1, W_2 \in \mathcal{W}(C)$ such that $M_1= \overline{W}_1$ and $M_2= \overline{W}_2$. Say $W_1= \{ h_1, h_1^c \}$ and $W_2= \{ h_2, h_2^c \}$, so that, if $p : X \to C$ denotes the projection onto $C$, we have $M_1= \{ p^{-1}(h_1), p^{-1}(h_1^c) \}$ and $M_2= \{ p^{-1}(h_2),p^{-1}(h_2^c) \}$. Notice that the inclusion $h_1 \subset h_2$ is equivalent to $p^{-1}(h_1) \subset p^{-1}(h_2)$, and similarly the inclusion $h_2 \subset h_1$ is equivalent to $p^{-1}(h_2) \subset p^{-1}(h_1)$. As a consequence, $M_1$ and $M_2$ are nested if and only if $W_1$ and $W_2$ are nested as well; a fortiori, $M_1$ and $M_2$ are transverse if and only if $W_1$ and $W_2$ are transverse as well. 
\end{proof}

\begin{cor}\label{cor:dimHW}
Let $X$ be a quasi-median graph endowed with a coherent system of wallspaces. Then $\dim (X, \mathcal{HW}) \leq \dim_{\square}(X) \cdot \sup\limits_{C \ \mathrm{clique}} \dim (C, \mathcal{W}(C))$.
\end{cor}

\begin{proof}
Let $M_1, \ldots, M_n$ be a collection of pairwise transverse hyperplane-walls. Let $J_1, \ldots, J_m$ denote the collection of pairwise distinct hyperplanes of $X$ underlying the walls $M_1, \ldots, M_n$. We deduce from Lemma \ref{lem:transverseHW} that $J_1, \ldots, J_m$ are pairwise transverse, hence $m \leq \dim_{\square}(X)$. Moreover, if $J_j$ is the underlying hyperplane of $M_{i_1}, \ldots, M_{i_k}$, then, according to Lemma \ref{lem:transverseHW}, there exist pairwise transverse walls $W_{i_1}, \ldots, W_{i_k} \in \mathcal{W}(C_j)$, where $C_j$ is a clique dual to $J_j$, such that $M_{i_r}= \overline{W_{i_r}}$ for every $1 \leq r \leq k$. Therefore, $k \leq \dim (C_j, \mathcal{W}(C_j))$. The conclusion follows.
\end{proof}

\begin{lemma}\label{lem:tangentHW}
Suppose that for any clique $C$ and any two vertices $x,y \in C$, there exists a wall in $\mathcal{W}(C)$ separating $x$ and $y$. If two hyperplane-walls $M_1 \in \mathcal{W}(J_1)$ and $M_2 \in \mathcal{W}(J_2)$ are tangent then either $J_1$ and $J_2$ are tangent or $J_1=J_2$ and there exist two tangent walls $W_1,W_2 \in \mathcal{W}(C)$, where $C$ is a clique dual to $J_1=J_2$, such that $M_1= \overline{W}_1$ and $M_2= \overline{W}_2$. 
\end{lemma}

\begin{proof}
First, notice that $J_1$ and $J_2$ cannot be transverse, since otherwise $M_1$ and $M_2$ would be transverse as well according to Lemma \ref{lem:transverseHW}. Therefore, either $J_1$ and $J_2$ are nested or $J_1=J_2$. If $J_1$ and $J_2$ are nested but not tangent, there exists a hyperplane $J$ separating $J_1$ and $J_2$, and taking a wall $M \in \mathcal{W}(J)$ separating $\mathrm{proj}_{N(J)}(N(J_1))$ and $\mathrm{proj}_{N(J)}(N(J_2))$ produces a hyperplane-wall separating $M_1$ and $M_2$. Therefore, if $J_1$ and $J_2$ are nested, they must be tangent. Finally, suppose that $J_1=J_2$, so that there exist two walls $W_1,W_2 \in \mathcal{W}(C)$, where $C$ is a clique dual to $J_1=J_2$, such that $M_1= \overline{W}_1$ and $M_2= \overline{W}_2$. Clearly, if there exists a wall $W \in \mathcal{W}(C)$ separating $W_1$ and $W_2$, then $\overline{W}$ separates $M_1$ and $M_2$, so $W_1$ and $W_2$ must be tangent.
\end{proof}

\noindent
For instance, setting $\mathcal{W}(C)= \{ \{ \{x \} , \{x\}^c \} \mid x \in C \}$ for every clique $C$ produces a coherent system of wallspaces. Notice that, if $C$ has cardinality two, then we introduced two indistinguishable walls; for convenience, we identify these two walls. We denote by $\mathcal{SW}$ the resulting collection of walls on $X$, which we refer to as \emph{sector-walls}\index{Sector-walls}. Notice that a sector-wall is the data of a sector and its corresponding cosector (as defined in Section \ref{section:convexhull}). This provides a systematic way to cubulate quasi-median graphs, allowing us to show that admitting a ``nice'' action on a quasi-median graph is equivalent to admitting a ``nice'' action on a CAT(0) cube complex. This idea is made precise by the following proposition.

\begin{prop}\label{prop:quasimedianimplycubical}
Let $G$ be a group acting on a quasi-median graph $X$. Let $C(X, \mathcal{SW})$ denote the CAT(0) cube complex obtained by (quasi-)cubulating $(X, \mathcal{SW})$. Suppose that one of the following assertions holds:
\begin{itemize}
	\item $G \curvearrowright X$ is not elliptic;
	\item $G \curvearrowright X$ is metrically proper;
	\item $G \curvearrowright X$ is geometric.
\end{itemize}
Then the action $G \curvearrowright C(X, \mathcal{SW})$ satisfies the corresponding property.
\end{prop}

\noindent
Essentially, Proposition \ref{prop:quasimedianimplycubical} will be a consequence of the various observations made by the next preliminary lemmas.

\begin{lemma}\label{lem:distanceSW}
For every vertices $x,y \in X$, $d(x,y) \leq d_{\mathcal{SW}}(x,y) \leq 2d(x,y)$. 
\end{lemma}

\begin{proof}
By definition of $\mathcal{SW}$, $d_{\mathcal{SW}}(x,y)$ is equal to the number of hyperplanes separating $x$ and $y$ counted with multiplicity, depending on whether it delimits more than two sectors. More precisely, if we set 
$$n(J)= \left\{ \begin{array}{cl} 1 & \text{if $J$ delimits two sectors} \\ 2 & \text{if $J$ delimits at least three sectors} \end{array} \right.$$
and if $J_1, \ldots, J_n$ are the hyperplanes separating $x$ and $y$, then 
$$d_{\mathcal{SW}}(x,y)= \sum\limits_{i=1}^n n(J_i).$$
Since $1 \leq n(J_i) \leq 2$ for every $1 \leq i \leq n$, and $n=d(x,y)$, the conclusion follows.
\end{proof}

\begin{cor}\label{cor:cubulationelliptic}
Let $X$ be a quasi-median graph and let $C(X, \mathcal{SW})$ denote the CAT(0) cube complex obtained by (quasi-)cubulating $(X, \mathcal{SW})$. For any $H \leq \mathrm{Aut}(X)$, the induced action $H \curvearrowright C(X, \mathcal{SW})$ is elliptic if and only if $H$ stabilises a prism of $X$. 
\end{cor}

\begin{proof}
Let $x \in X$ be a vertex, and let $\sigma_x$ denote the associated principal orientation, thought of as a vertex of $C(X, \mathcal{SW})$. Noticing that, for every $g \in H$,
$$d_{C(X, \mathcal{SW})}(\sigma_x,g \cdot \sigma_x)= d_{C(X, \mathcal{SW})}(\sigma_x,\sigma_{g \cdot x})= d_{\mathcal{SW}}(x,g \cdot x),$$
we deduce from Lemma \ref{lem:distanceSW} that the action $H \curvearrowright C(X, \mathcal{SW})$ has a bounded orbit if and only if the action $H \curvearrowright X$ has a bounded orbit as well. On the other hand, we know from Theorem \ref{thm:fixedpoint} that $H \curvearrowright X$ has a bounded orbit if and only if $H$ stabilises a prism.
\end{proof}

\begin{lemma}\label{lem:transverseSW}
Two sector-walls are transverse if and only if their underlying hyperplanes are transverse.
\end{lemma}

\begin{proof}
For every clique and every vertices $x,y \in C$, the walls
$$\{ \{x \} , \{x \}^c \}, \ \{ \{y \}, \{y \}^c \} \in \mathcal{W}(C)$$
are nested, so that our lemma follows from Lemma \ref{lem:transverseHW}.
\end{proof}

\noindent
For our next and last preliminary lemma, we need the following definition:

\begin{definition}
The \emph{simplicial dimension} of a quasi-median graph $X$, denoted by $\dim_{\triangle} X$, is the maximal cardinality of a clique of $X$. The \emph{dimension} of $X$ is $\dim(X)= \max ( \dim_{\square} X, \dim_{\triangle} X)$.
\end{definition}

\begin{lemma}\label{lem:cubulationcube}
Let $X$ be a quasi-median graph of finite dimension and let $C(X, \mathcal{SW})$ denote the CAT(0) cube complex obtained by (quasi-)cubulating $(X, \mathcal{SW})$. There exists an $\mathrm{Aut}(X)$-invariant surjective map from the maximal cubes of $C(X, \mathcal{SW})$ to the maximal prisms of $X$ such that the pre-image of a maximal prism of $C(X, \mathcal{SW})$ has cardinality at most $\dim_{\square}(X) \cdot \dim_{\triangle} (X)$.
\end{lemma}

\begin{proof}
Let $C$ be a maximal cube of $C(X, \mathcal{SW})$. Because $\dim C(X, \mathcal{SW})= \dim_{\square}X<+ \infty$, this cube is uniquely determined by the hyperplanes dual to it, which produces a maximal collection of pairwise transverse hyperplanes (see Proposition \ref{prop:maxprism}). According to Theorem \ref{thm:quasicubulation}, to such a collection corresponds a maximal collection of pairwise transverse sector-walls of $X$, producing a maximal collection of pairwise transverse hyperplanes of $X$ according to Lemma \ref{lem:transverseSW}, and finally a maximal prism $P(C)$ according to Proposition \ref{prop:maxprism}. Notice that, according to the various results we used, the map $C \mapsto P(C)$ is $\mathrm{Aut}(X)$-invariant. 

\medskip \noindent
By construction, two maximal cubes $C,C'$ of $C(X, \mathcal{SW})$ satisfy $P(C)=P(C')$ if and only if the two collections of pairwise transverse sector-walls produce the same collection of pairwise transverse hyperplanes. Because there exist at most $\dim_{\square}(X)$ pairwise transverse hyperplanes and that each hyperplane delimits at most $\dim_{\triangle}(X)$ sectors, we deduce that the cardinality of a pre-image of $C \mapsto P(C)$ is at most $\dim_{\square}(X) \cdot \dim_{\triangle}(X)$.
\end{proof}

\begin{proof}[Proof of Proposition \ref{prop:quasimedianimplycubical}.]
The first point is a direct consequence of Corollary \ref{cor:cubulationelliptic}, and the second one a direct consequence of Lemma \ref{lem:distanceSW}. Now, suppose that $G$ acts geometrically on $X$. We already know that the action $G \curvearrowright C(X, \mathcal{SW})$ is metrically proper according to the previous point. Next, since a group acts geometrically on $X$, necessarily $X$ must be uniformly locally finite. In particular, $X$ is finite dimensional, so that Lemma \ref{lem:cubulationcube} applies. Because the action $G \curvearrowright X$ has finitely many orbits of maximal prisms, we deduce that the action $G \curvearrowright C(X, \mathcal{SW})$ has finitely many orbits of maximal cubes. Consequently, $G \curvearrowright C(X, \mathcal{SW})$ is geometric.
\end{proof}

\subsection{Actions of groups}\label{section:localglobalwall}

\noindent
In this section, we are interested in group actions. First notice that, if a group $G$ acts on a quasi-median graph $X$, endowed with a system of wallspaces which is coherent and \emph{$G$-invariant} (ie., $g\mathcal{W}(C)= \mathcal{W}(gC)$ for every clique $C$ and every $g \in G$), then $G$ acts on the space with walls $(X,\mathcal{HW})$. Now, we would like to determine when this induced action $G \curvearrowright (X,\mathcal{HW})$ is metrically proper, cocompact, virtually special, etc. We give some criteria  below. 

\begin{prop}\label{prop:producingproperaction}
Let $G$ be a group acting on a quasi-median graph $X$ endowed with a $G$-invariant coherent system of wallspaces. Suppose that
\begin{itemize}
	\item any vertex of $X$ belongs to finitely many cliques;
	\item any vertex-stabiliser is finite;
	\item for every clique $C$, any two points of $C$ are separated by a wall of $\mathcal{W}(C)$;
	\item for every clique $C$, the space with walls $(C, \mathcal{W}(C))$ is locally finite.
\end{itemize}
Then the action $G \curvearrowright (X, \mathcal{HW})$ is metrically proper.
\end{prop}

\begin{proof}
Let us consider the system of pseudo-metrics $(C, d_{\mathcal{W}(C)})$. Notice that, because any two points of a clique are sepated by a wall, this is a system of metrics. As a consequence of Fact \ref{fact:pseudodistcoh}, this system is coherent, and according to Lemma \ref{lem:distanceHW}, $d_{\mathcal{HW}}$ is the global metric extending this system. Thus, our proposition is a direct consequence of Lemma \ref{lem:whenmetricallyproper}. 
\end{proof}

\noindent
For the next proposition, we need to introduce some notation. If $P=C_1 \times \cdots C_n$ is a prism of $X$, we denote by $\mathcal{W}(P)$ the collection of walls $\mathcal{W}(C_1) \times \cdots \times \mathcal{W}(C_n)$. 

\begin{prop}\label{prop:producingcocompactaction}
Let $G$ be a group acting on a quasi-median graph $X$ endowed with a $G$-invariant coherent system of wallspaces. Suppose that
\begin{itemize}
	\item the cubical dimension of $X$ is finite;
	\item for every clique $C$, $\mathcal{W}(C) \neq \emptyset$ and $\dim (C, \mathcal{W}(C))<+ \infty$;
	\item $X$ contains finitely many $G$-orbits of prisms;
	\item for every maximal prism $P$, the action $\mathrm{stab}(P) \curvearrowright (P, \mathcal{W}(P))$ is cocompact.
\end{itemize}
Then the action $G \curvearrowright (X, \mathcal{HW})$ is cocompact.
\end{prop}

\begin{proof}
Define an \emph{$M$-collection} as the data of a maximal prism $P$ and, for each hyperplane $J$ dual to it, a maximal collection of pairwise transverse walls of $\mathcal{W}(J) \subset \mathcal{W}(P)$. We claim that any maximal collection of pairwise transverse hyperplane-walls naturally defines an $M$-collection. It is worth noticing that $\dim (X, \mathcal{HW})$ is finite according to Corollary \ref{cor:dimHW}, so that such a collection must be finite.

\medskip \noindent
Let $M_1, \ldots, M_n \in \mathcal{HW}$ be a maximal collection of pairwise transverse hyperplane-walls. Let $J_1, \ldots, J_m$ denote the associated collection of pairwise distinct hyperplanes of $X$; according to Lemma \ref{lem:transverseHW}, these hyperplanes are pairwise transverse. If there exists a hyperplane $J'$ transverse to $J_1, \ldots, J_m$, then, according to Lemma \ref{lem:transverseHW}, any wall $M' \in \mathcal{W}(J')$ (there exists such a wall since $\mathcal{W}(J')$ is non empty by assumption) will be transverse to $M_1, \ldots, M_n$, contradicting the maximality of our collection. Therefore, $J_1, \ldots, J_m$ is a maximal collection of pairwise transverse hyperplanes. Let $P$ denote the maximal prism which is associated to this collection by Proposition \ref{prop:maxprism}. Fix some $1 \leq j \leq n$, and let $M_{i_1}, \ldots, M_{i_k}$ be the hyperplane-walls of our collection with $J_j$ as underlying hyperplane. If there existed a hyperplane-wall $M \in \mathcal{W}(J_j)$ transverse to $M_{i_1}, \ldots, M_{i_k}$ then, according to Lemma \ref{lem:transverseHW}, $M$ would be transverse to $M_1, \ldots, M_n$, contradicting the maximality of our collection. Therefore, $M_{i_1}, \ldots, M_{i_k}$ is a maximal collection of pairwise transverse hyperplane-walls in $\mathcal{W}(J_j)$. Thus, $P$ and $M_1, \ldots, M_n$ naturally defines an $M$-collection.

\medskip \noindent
On the other hand, we know that there exist finitely many $G$-orbits of (maximal) prisms of $X$ and, for each maximal prism $P$, the action $\mathrm{stab}(P) \curvearrowright (P, \mathcal{W}(P))$ is cocompact, so there must exist only finitely many $G$-orbits of $M$-collections. A fortiori, there exist only finitely many $G$-orbits of maximal collections of pairwise transverse hyperplane-walls in $\mathcal{HW}$, ie., the action $G \curvearrowright (X, \mathcal{HW})$ is cocompact.
\end{proof}

\noindent
Before stating our next proposition, recall that a subgroup $H \leq G$ is a \emph{retract}\index{Retracts (subgroups)} if there exists an epimorphism $r : G \twoheadrightarrow H$ (called a \emph{retraction}) satisfying $r_{|H}= \mathrm{Id}_H$. 

\begin{prop}\label{prop:producingspecialaction}
Let $G$ be a group acting on a quasi-median graph $X$ endowed with a $G$-invariant coherent system of wallspaces. Suppose that
\begin{itemize}
	\item for every clique $C$, two distinct vertices of $C$ are separated by some wall of $\mathcal{W}(C)$;
	\item the action $G \curvearrowright X$ is special.
\end{itemize}
Then $\mathrm{Obs}(G \curvearrowright (X, \mathcal{HW})) = \bigcup\limits_{J \ \mathrm{hyperplane}} \mathrm{Obs}(\mathrm{stab}(J) \curvearrowright (N(J), \mathcal{W}(J)))$. If moreover
\begin{itemize}
	\item for every hyperplane $J$, $\mathrm{stab}(J)$ is a finitely-generated retract of $G$;
	\item $X$ has finitely many $G$-orbits of hyperplanes;
	\item for every hyperplane $J$, $\mathrm{stab}(J) \curvearrowright (N(J), \mathcal{W}(J))$ is virtually special;
\end{itemize}
then $G \curvearrowright (X, \mathcal{HW})$ is virtually special.
\end{prop}

\noindent
We recall that $\mathcal{W}(J)$ denotes the collection of extended walls $\overline{\mathcal{W}}(C)$, where $C$ is any clique dual to $J$. (This collection does not depend on the choice of $C$ since our system is coherent.) 

\begin{proof}[Proof of Proposition \ref{prop:producingspecialaction}.] 
Let $g \in \mathrm{Obs}(G \curvearrowright (X, \mathcal{HW}))$. For convenience, let $\mathrm{Obs}(J)$ denote the obstruction $\mathrm{Obs}(\mathrm{stab}(J) \curvearrowright (N(J),\mathcal{W}(J)))$ for every hyperplane $J$ of $X$. Three cases may happen.

\medskip \noindent
Suppose that there exists a wall $M \in \mathcal{HW}$ such that $M$ and $gM$ are transverse. Let $J$ denote the underlying hyperplane of $M$. According to Lemma \ref{lem:transverseHW}, either $gJ=J$ or $J$ and $gJ$ are transverse. The latter case is impossible since $G \curvearrowright X$ is special. Therefore, $g \in \mathrm{stab}(J)$. We conclude that $g \in \mathrm{Obs}(J)$. 

\medskip \noindent
Suppose that there exists a wall $M \in \mathcal{HW}$ such that $M$ and $gM$ are tangent. Let $J$ denote the underlying hyperplane of $M$. Notice that $J$ and $gJ$ cannot be tangent since $G \curvearrowright X$ is special. Therefore, we deduce from Lemma \ref{lem:tangentHW} that $g \in \mathrm{Obs}(J)$. 

\medskip \noindent
Suppose finally that there exist three walls $M_1,M_2,M_3 \in \mathcal{HW}$ such that $M_1$ and $M_2$ are tangent, $M_2$ and $M_3$ are transverse, and $gM_1=M_3$. Let $J_1,J_2,J_3$ denote the underlying hyperplanes of $M_1,M_2,M_3$ respectively; notice that $gJ_1=J_3$. We distinguish five cases.
\begin{itemize}
	\item If $J_1,J_2,J_3$ are pairwise distinct, then $J_1$ and $J_2$ must be tangent and $J_2$ and $J_3$ transverse, which is impossible since $G \curvearrowright X$ is special.
	\item If $J_1=J_2 \neq J_3$, then $gJ_2 = gJ_1=J_3$ is transverse to $J_2$, contradicting the fact that $G \curvearrowright X$ is special.
	\item If $J_2=J_3 \neq J_1$, then $g^{-1}J_2 =  g^{-1}J_3=J_1$ is tangent to $J_2$, contradicting the fact that $G \curvearrowright X$ is special.
	\item If $J_1=J_3 \neq J_2$, then $J_2$ is both transverse and tangent to $J_1=J_3$, which is impossible.
	\item If $J_1=J_2=J_3$, let $J$ denote this common hyperplane. In particular, $gJ=gJ_1=J_3=J$ hence $g \in \mathrm{stab}(J)$. We conclude that $g \in \mathrm{Obs}(J)$.
\end{itemize}
We have proved that $\mathrm{Obs}(G \curvearrowright (X, \mathcal{HW})) \subset \bigcup\limits_{J \ \mathrm{hyperplane}} \mathrm{Obs}(\mathrm{stab}(J) \curvearrowright (N(J), \mathcal{W}(J)))$. The reverse inclusion is clear, concluding the proof of the first assertion of our proposition.

\medskip \noindent
Now, suppose that hyperplane-stabilisers are finitely generated retracts, that $X$ has finitely many $G$-orbits of hyperplanes, and that the actions $\mathrm{stab}(J) \curvearrowright (N(J), \mathcal{W}(J))$ are all virtually special.

\medskip \noindent
Let $J_1, \ldots, J_n$ be a set of representatives of the hyperplanes of $X$ under the action of $G$. Let $1 \leq i \leq n$. Because the action $\mathrm{stab}(J_i) \curvearrowright (N(J_i),\mathcal{W}_i)$ is virtually special, and that $\mathrm{stab}(J_i)$ is finitely generated, there exists a normal finite-index subgroup $H_i \leq \mathrm{stab}(J_i)$ such that the induced action $H_i \curvearrowright (N(J_i), \mathcal{W}(J_i))$ is special. In particular, if $F_i$ denotes a set of representatives of $\mathrm{stab}(J_i)/H_i$ different from $H_i$, then
$$\mathrm{Obs}(\mathrm{stab}(J_i) \curvearrowright (N(J_i), \mathcal{W}(J_i))) \subset F_i H_i$$
according to Fact \ref{fact:obs}. Fixing a retraction $r_i : G \to \mathrm{stab}(J_i)$, we introduce
$$K_i= \mathrm{ker} \left( G \overset{r_i}{\twoheadrightarrow} \mathrm{stab}(J_i) \twoheadrightarrow \mathrm{stab}(J_i) / H_i \right).$$
Notice that $K_i$ is a normal finite-index subgroup of $G$ satisfying $H_i \leq K_i$ and $K_i \cap F_i= \emptyset$. Set $K= \bigcap\limits_{i=1}^n K_i$. This defines a normal finite-index subgroup of $G$. 

\medskip \noindent
We claim that the induced action $K \curvearrowright (X, \mathcal{HW})$ is special. Let
$$g \in \mathrm{Obs}(G \curvearrowright (X, \mathcal{HW})) \subset \bigcup\limits_{k \in G} \bigcup\limits_{i=1}^n k F_iH_i k^{-1},$$
say $g=kfhk^{-1}$ for some $1 \leq i \leq n$, $k \in G$, $f \in F_i$ and $h \in H_i$. Then we deduce from $g \in K$ that
$$f \in k^{-1}Kkh^{-1} = Kh^{-1} \subset K_ih^{-1}=K_i,$$
contradicting the fact that $K_i \cap F_i = \emptyset$. Therefore, $K \cap \mathrm{Obs}(G \curvearrowright (X, \mathcal{HW})) = \emptyset$.
\end{proof}

\begin{prop}\label{prop:producingproperaction2}
Let $G$ be a group acting on a quasi-median graph $X$ endowed with a $G$-invariant coherent system of wallspaces. Suppose that
\begin{itemize}
	\item for every clique $C$, any two vertices of $C$ are separated by a wall of $\mathcal{W}(C)$;
	\item for every prism $P$ of $X$, the action $\mathrm{stab}(P) \curvearrowright (P, \mathcal{W}(P))$ is properly discontinuous;
	\item vertex-stabilisers are finite.
\end{itemize}
Then the action $G \curvearrowright (X, \mathcal{HW})$ is properly discontinuous. 
\end{prop}

\noindent
Recall that the action $G \curvearrowright (X, \mathcal{HW})$ is properly discontinuous if the action of $G$ on the CAT(0) cube complex obtained by cubulating $(X, \mathcal{HW})$ is properly discontinuous. Since an action on a CAT(0) cube complex is properly discontinuous if and only if its vertex-stabilisers are finite, Proposition \ref{prop:producingproperaction2} precisely means that the stabiliser of any orientation of $\mathcal{HW}$ must be finite. Of course, the first step toward the proof of our proposition is to understand the orientations of $\mathcal{HW}$. 

\medskip \noindent
Fix a prism $P$ of $X$ which is a Cartesian product of cliques $C_1 \times \cdots \times C_n$, and, for every $1 \leq i \leq n$, let $\sigma_i$ be a non principal orientation of $\mathcal{W}(C_i)$. In particular, $\sigma_P=\prod\limits_{i=1}^n \sigma_i$ defines an orientation of $\mathcal{W}(P)$. Now, we extend $\sigma_P$ on $\mathcal{HW}$ in the following way: if $W \in \mathcal{W}(J)$ for some hyperplane $J$ which is disjoint from $P$, set $\sigma_P(W)$ as the unique halfspace of $W$ containing $P$. It is not difficult to verify that $\sigma_P$ defines an orientation. An orientation arising in this way will be called \emph{semiprincipal}. 

\begin{lemma}\label{lem:dichotomyorientations}
Let $X$ be a quasi-median graph endowed with a coherent system of wallspaces. Suppose that, for every clique $C$ of $X$, any two vertices of $C$ are separated by a wall of $\mathcal{W}(C)$. Then any orientation of $\mathcal{HW}$ is either principal or semiprincipal. 
\end{lemma}

\begin{proof}
Fix some orientation $\sigma$ of $\mathcal{HW}$. Let $\sigma_x$ be a principal orientation minimizing the distance to $\sigma$ in $C(X, \mathcal{HW})$, ie., the number of walls on which $\sigma_x$ and $\sigma$ differ, and let $\mathfrak{I}$ denote the collection of walls on which they differ.

\medskip \noindent
We claim that, for any hyperplane $J$ associated to a wall $W \in \mathfrak{I}$, we have $x \in N(J)$. So fix a wall $W \in \mathfrak{I}$ and let $J$ denote its underlying hyperplane. Let $y$ denote the projection of $x$ onto $N(J)$. Of course, $\sigma_x(W')=\sigma_y(W')$ for any wall $W' \in \mathcal{HW}$ which does not separate $x$ and $y$. Now, suppose that $W' \in \mathcal{HW}$ is a wall separating $x$ and $y$. In particular, its underlying hyperplane $J'$ separates $x$ and $y$. As a consequence of Lemma \ref{lem:projseparate}, $J'$ separates $x$ and $N(J)$, so that $N(J)$ must be included into a single sector $S$ delimited by $J'$; notice that $y \in S$. On the other hand, because $\sigma_x(W) \neq \sigma(W)$, we know that $\sigma(W)$ does not contain $x$, so $\sigma(W) \subset S$. By noticing that $\sigma(W) \subset S \subset \sigma_y(W')$, we conclude that 
$$\sigma(W')= \sigma_y(W') \neq \sigma_x(W').$$
Therefore, since $\sigma_x$ is a principal orientation minimizing the number of walls on which $\sigma_x$ and $\sigma$ differ, we deduce that there exist no wall of $\mathcal{HW}$ separating $x$ and $y$. On the other hand, a hyperplane $H$ separating $x$ and $y$ would produce a wall separating $x$ and $y$ by taking a wall of $\mathcal{W}(H)$ separating the projections of $x$ and $y$ onto $N(H)$ (such a wall would exist by assumption), so we deduce that no hyperplane separates $x$ and $y$, ie., $x=y$. A fortiori, $x \in N(J)$. 

\medskip \noindent
We claim that any two hyperplanes underlying two walls of $\mathfrak{I}$ are transverse. Let $J_1,J_2$ be two non transverse hyperplanes such that $x \in N(J_1) \cap N(J_2)$ and such that $J_1$ is the underlying hyperplane of some $W_1 \in \mathfrak{I}$. Fix a wall $W_2 \in \mathcal{W}(J_2)$. Because $W_1 \in \mathfrak{I}$, we know that $\sigma(W_1) \neq \sigma_x(W_1)$, hence $\sigma(W_1) \subset \sigma_x(W_2)$. This implies that $\sigma(W_2)= \sigma_x(W_2)$; a fortiori, $W_2 \notin \mathfrak{I}$. This proves our claim.

\medskip \noindent
Thus, if $\mathfrak{H}$ denotes the set of the hyperplanes underlying some wall of $\mathfrak{I}$, we have proved that $\mathfrak{H}$ defines a (finite) collection of pairwise transverse hyperplanes satisfying $x \in \bigcap\limits_{J \in \mathfrak{H}} N(J)$. It follows from Fact \ref{fact:spanningprism} that there exists a prism $P=C_1 \times \cdots \times C_n$ such that $x \in P$ and such that $\mathfrak{H}$ is precisely the set of the hyperplanes dual to it. 

\medskip \noindent
Now, we claim that, for every $1 \leq i \leq n$, the restriction $\sigma_{|\mathcal{W}(C_i)}$ defines an orientation of $\mathcal{W}(C_i)$ which is not principal. Indeed, suppose by contradiction that there exists some $1 \leq i \leq n$ such that the restriction $\sigma_{|\mathcal{W}(C_i)}$ is a principal orientation of $\mathcal{W}(C_i)$; let $y \in C_i$ denote the corresponding vertex. Because $x$ and $y$ both belong to $C_i$, the orientations $\sigma_x$ and $\sigma_y$ may only differ on walls with underlying hyperplane $J_i$, the hyperplane dual to $C_i$. On the other hand, $\sigma_y$ and $\sigma$ coincide on every wall with underlying hyperplane $J_i$. It follows that the number of walls on which $\sigma_y$ and $\sigma$ differ is stricly smaller than the number of walls on which $\sigma_x$ and $\sigma$ differ, contradicting our choice of $\sigma_x$. This proves our claim. 

\medskip \noindent
Let $\sigma_P$ denote the semiprincipal orientation associated to $\prod\limits_{i=1}^n \sigma_{|\mathcal{W}(C_i)}$. We claim that $\sigma= \sigma_P$. Let $W \in \mathcal{HW}$ be a wall and let $J$ denote its underlying hyperplane. If $W \in \mathcal{W}(C_i)$ for some $1 \leq i \leq n$ (ie., if $J$ intersects the prism $P$), then by construction we know that $\sigma(W)= \sigma_P(W)$. Next, if $J$ is disjoint from $P$, $\sigma_P(W)$ must be the halfspace delimited by $W$ which contains $P$ (or equivalently, $x \in P$), hence
$$\sigma_P(W)= \sigma_x(W)= \sigma(W)$$
because $J \notin \mathfrak{H}$ so that $W \notin \mathfrak{I}$. A fortiori, $\sigma$ is a semiprincipal orientation.
\end{proof}

\begin{proof}[Proof of Proposition \ref{prop:producingproperaction2}.]
In order to prove that the action $G \curvearrowright (X, \mathcal{HW})$ is properly discontinuous, we will prove that the stabiliser of any orientation $\sigma$ of $\mathcal{HW}$ is finite. According to Lemma \ref{lem:dichotomyorientations}, we know that $\sigma$ is either principal or semiprincipal. In the former case, say $\sigma=\sigma_x$ for some vertex $x \in X$, we have $\mathrm{stab}(\sigma_x)= \mathrm{stab}(x)$, which is finite by assumption. Next, suppose that $\sigma$ is semiprincipal. So there exists a prism $P=C_1 \times \cdots \times C_n$ such that, for every $1 \leq i \leq n$, the restriction $\sigma_{\mathcal{W}(C_i)}$ defines an orientation which is not principal, and such that, for any wall $W \in \mathcal{HW}$ whose underlying hyperplane does not intersect $P$, $\sigma(W)$ is the halfspace delimited by $W$ which contain $P$. Let $\mathfrak{H}$ denote the set of the hyperplanes $J$ such that the restriction $\sigma_{|\mathcal{W}(J)}$ is principal, and let $S(J)$ denote the sector delimited by $J$ which contains $P$; notice that, if we fix some clique $C$ dual to $J$, then $S(J)$ is also the sector $[C,x]$ where $x \in C$ is defined by $\sigma_{|\mathcal{W}(C)} = \sigma_x$. The intersection $C= \bigcap\limits_{J \in \mathfrak{H}} S(J)$ must be $\mathrm{stab}(\sigma)$-invariant since the collection $\mathfrak{H}$ and the sectors $S(J)$ are defined only from $\sigma$. We claim that $C=P$.

\medskip \noindent
Let $J \in \mathfrak{H}$. Fix a clique $C$ dual to $J$ and let $x \in C$ be the vertex defined by $\sigma_{|\mathcal{W}(C)} = \sigma_x$, so that $S(J)= [C,x]$. Suppose by contradiction that $x$ is different from the projection $y$ of $P$ onto $C$. Then there exists a wall $W \in \mathcal{W}(C)$ separating $x$ and $y$. Clearly, $\sigma(W)=\sigma_x(W)$ is disjoint from $P$ since $y \notin \sigma(W)$. On the other hand, by definition of $P$, we know that $P \subset \sigma(W')$ for every $W' \in \mathcal{W}(C)$ since $J$ does not intersect $P$, whence a contradiction. Therefore, $x=y$ holds, and we deduce that $S(J)$ is the sector delimited by $J$ containing $P$. A fortiori, $P \subset C$. Conversely, let $x \in X$ be a vertex satisfying $x \notin P$. Let $y \in P$ denote the projection of $x$ onto $P$ (which is well-defined since prisms are gated according to Lemma \ref{lem:prismgated}). Because $x \notin P$, necessarily $x \neq y$ so that there exists a hyperplane $J$ separating $x$ and $y$. Notice that $J$ does not intersect $P$ according to Lemma \ref{lem:projseparate}, so $J \in \mathfrak{H}$. Moreover, we saw that $S(J)$ is the sector delimited by $J$ which contains $P$, so $x \notin S(J)$ and $C \subset S(J)$, hence $x \notin C$. A fortiori, $C \subset P$, which concludes our claim.

\medskip \noindent
Thus, we have proved that the prism $P$ is $\mathrm{stab}(\sigma)$-invariant. On the other hand, we know by assumption that the action $\mathrm{stab}(P) \curvearrowright (P, \mathcal{W}(P))$ is properly discontinuous, so that the stabiliser (in $\mathrm{stab}(P)$) of the orientation $\prod\limits_{i=1}^n \sigma_{|\mathcal{W}(C_i)}$ of $\mathcal{W}(P)$ is finite. We conclude that the stabiliser (in $G$) of $\sigma$ must be finite as well. 
\end{proof}

\subsection{From actions to walls}\label{section:actionwall}

\noindent
In order to apply the different propositions proved in the previous section, we need to find a collection of walls on each clique of our quasi-median graph. In the cases we will be interested in, the induced action of each clique-stabiliser on the corresponding clique often turns out to be transitive and free, so that finding walls on a clique is equivalent to finding walls on its stabiliser. For this purpose, one possibility is to consider an action on a CAT(0) cube complex and then to pullback the walls defined on the cube complex to our group. In fact, the construction presented below works in full generality for actions on quasi-median graphs, but, as noticed by Proposition \ref{prop:quasimedianimplycubical}, we do not really loss generality.

\medskip \noindent
Let $G$ be a group acting on a CAT(0) cube complex $X$. Fix a base vertex $x_0 \in X$. Without loss of generality, we suppose that the combinatorial convex hull of the orbit $G \cdot x_0$ is the whole $X$. To any hyperplane $J$ of $X$, delimiting two halfspaces $J^+$ and $J^-$, we associate the wall on $G$
$$\mathfrak{m}(J)= \left\{ \{ g \in G \mid g \cdot x_0 \in J^+ \}, \{ g \in G \mid g \cdot x_0 \in J^- \} \right\};$$
and we introduce the collection of walls 
$$\mathcal{M}(G \curvearrowright X)= \left\{ \mathfrak{m}(J) \mid J \ \text{hyperplane of} \ X \right\}.$$ 
Notice that, for every hyperplane $J$ and every $g \in G$, we have $g \cdot \mathcal{M}(J)= \mathcal{M} ( g \cdot J)$, so that $\mathcal{M}(G \curvearrowright X)$ defines a $G$-invariant collection of walls. Moreover, for every hyperplane $J$ of $X$ and every $g,h \in G$, the wall $\mathfrak{m}(J)$ separates $g$ and $h$ if and only if $J$ separates $g \cdot x_0$ and $h \cdot x_0$, so that
$$d_{\mathcal{M}(G \curvearrowright X)}(g,h)= d_X(g \cdot x_0, h \cdot x_0).$$
As a consequence, $(G,\mathcal{M}=\mathcal{M}(G \curvearrowright X))$ defines a space with walls. But does the properties of the action $G \curvearrowright X$ transfer to $G \curvearrowright (G, \mathcal{M})$? The next two lemmas provide positive answers for the metrical properness and the cocompactness.

\begin{lemma}\label{lem:transferproper}
If $G \curvearrowright X$ is metrically proper, then so is $G \curvearrowright (G, \mathcal{M})$.
\end{lemma}

\begin{proof}
Because the action $G \curvearrowright X$ is metrically proper, we deduce that for every $R \geq 0$ 
$$\# \{ g \in G \mid d_{\mathcal{M}}(1,g) \leq R \} = \# \{ g \in G \mid d_X(x_0,g \cdot x_0) \leq R \} <+ \infty.$$
Therefore, the action $G \curvearrowright (G, \mathcal{W})$ is metrically proper as well.
\end{proof}

\begin{lemma}\label{lem:transfercocompact}
If $G \curvearrowright X$ is cocompact and $X$ locally finite, then $G \curvearrowright (G, \mathcal{M})$ is cocompact.
\end{lemma}

\begin{proof}
Let $Y \subset X$ be a compact fundamental domain containing $x_0$. Because $X$ is locally finite, the $1$-neighborhood $Y^{+1}$ of $Y$ must be finite. Let $J_1, \ldots, J_p$ denote the hyperplanes intersecting $Y^{+1}$. For every $1 \leq i \leq p$, there exists some $g_i \in G$ such that $J_i$ separates $x_0$ and $g_i \cdot x_0$. Set $R= \max \{ d(x_0,g_i \cdot x_0) \mid 1 \leq i \leq p \}$. 

\medskip \noindent
Let $\mathfrak{m}(H_1), \ldots, \mathfrak{m}(H_n) \in \mathcal{M}$ be a collection of pairwise transverse walls. Up to reordering these walls, suppose that, for some $1 \leq s \leq n$, $H_1, \ldots, H_s$ is a maximal subcollection of $H_1, \ldots, H_n$ of pairwise transverse hyperplanes. If we fix some $x \in \bigcap\limits_{i=1}^s N(H_i)$, we can suppose, up to a translation, that $x \in Y$, so that $H_1, \ldots, H_s \in \{ J_1, \ldots, J_p \}$. 

\medskip \noindent
Now, fix some $s+1 \leq j \leq n$. We know that $H_j$ must be disjoint from some $H_k$, $1 \leq k \leq n$. Let $H_j^+$ denote the halfspace delimited by $H_j$ which contains $H_k$, and similarly, let $H_k^+$ denote the halfspace delimited by $H_k$ which contains $H_j$. If $G \cdot x_0$ intersected $H_j^+ \cap H_k^+$, then $\mathfrak{m}(H_j)$ and $\mathfrak{m}(H_k)$ would be nested. So, because we supposed $\mathfrak{m}(H_j)$ and $\mathfrak{m}(H_k)$ transverse, necessarily $G \cdot x_0$ does not intersect $H_j^+ \cap H_k^+$, so that $H_j$ and $H_k$ are indistinguishable. As a consequence, there exists some $1 \leq \ell \leq p$ such that $H_j$ and $J_{\ell}=H_k$ are indistinguishable. In particular, $H_j$ must separate $x_0$ and $g_{\ell} \cdot x_0$, so that $H_j$ intersects the ball $B(x_0,R) \subset Y^{+R}$.

\medskip \noindent
Thus, we have proved that $H_i$ intersects $Y^{+R}$ for every $1 \leq i \leq n$. Notice that, because $X$ is locally finite, $Y^{+R}$ must be finite. In particular, the cardinality $n$ of our collection is bounded above by the (finite) number of hyperplanes intersecting $Y^{+R}$, hence

\begin{fact}\label{fact:transferdim}
If $G \curvearrowright X$ is cocompact and $X$ locally finite, then $\dim (G, \mathcal{M})<+ \infty$.
\end{fact}

\noindent
But we have also proved that there exist only finitely many $G$-orbits of (maximal) collections of pairwise transverse walls of $\mathcal{M}$, ie., $G \curvearrowright (G, \mathcal{M})$ is cocompact.
\end{proof}

\noindent
As a direct consequence of Lemma \ref{lem:transferproper} and Lemma \ref{lem:transfercocompact},

\begin{cor}
If $G \curvearrowright X$ is geometric, then so is $G \curvearrowright (G, \mathcal{M})$.
\end{cor}

\noindent
It is worth noticing that the dimension of $(G, \mathcal{M})$ does not depend only on $X$. To be precise, if given a hyperplane $J$ of $X$ we denote by $\iota(J)$ the number of hyperplanes indistinguishable to $J$, we have 
$$\dim(G, \mathcal{M}) = \max \left\{ \sum\limits_{i=1}^n \iota(J_i) \mid J_1, \ldots, J_n \ \text{pairwise transverse} \right\}.$$
Nevertheless, we know that the dimension is finite whenever the cube complex $X$ is finite dimensional. In fact, for geometric actions, it is possible to introduce another space with walls with a better control on the dimension. If $G$ acts on a CAT(0) cube complex $X$, we denote by $\mathcal{N}= \mathcal{N}(G \curvearrowright X)$ the collection of walls obtained from $\mathcal{M}$ by identifying two indistinguishable walls. Explicitely, the walls of $\mathcal{N}$ are the $\mathfrak{m}(J)$'s where $\mathfrak{m}(J_1)= \mathfrak{m}(J_2)$ if $J_1$ and $J_2$ are indistinguishable.

\begin{lemma}\label{lem:transfergeometricaction}
If the action $G\curvearrowright X$ is geometric, then so is $G \curvearrowright (G, \mathcal{N})$. Moreover, $\dim(G, \mathcal{N}) \leq \dim X$. 
\end{lemma}

\begin{proof}
First, we claim that there exists a constant $R \geq 0$ such that any two indistinguishable hyperplanes $J_1,J_2$ satisfies $d(N(J_1),N(J_2)) \leq R$. 

\medskip \noindent
Let $\{H_1, \ldots, H_n \}$ be a set of representatives for the action of $G$ on the hyperplanes of $X$. For every $1 \leq i \leq n$, fix some $g_i \in G$ such that $H_i$ separates $x_0$ and $g_i \cdot x_0$. Set $R= \max \{ d(x_0,g_i \cdot x_0) \mid 1 \leq i \leq n\}$. Now, let $J_1,J_2$ be two indistinguishable hyperplanes. There exists some $g \in G$ such that $g \cdot J_1 \in  \{ H_1, \ldots, H_n\}$, say $g \cdot J_1=H_1$. Therefore, the hyperplanes $H_1$ and $g \cdot J_2$ are necessarily indistinguishable, so that they both separate $x_0$ and $g_1 \cdot x_0$, hence
$$d(N(J_1),N(J_2))=d( N(H_1),N(gJ_2)) \leq d(x_0,g_1 \cdot x_0) \leq R.$$
This proves our claim.

\medskip \noindent
Now, fixing some $x \in X$ and $K \geq \dim(X)$, we want to prove that the set
$$F = \{ g \in G \mid d_{\mathcal{N}}(x,g \cdot x) \leq K \}$$
is finite. Let $g \in G$ satisfy $d(x,g \cdot x) \geq \mathrm{Ram}(K(R+3)+1)$, where $\mathrm{Ram}(\cdot)$ denotes the Ramsey number. If so, we know that $x$ and $g \cdot x$ are separated by at least $\mathrm{Ram}(K(R+3)+1)$ hyperplanes, so that $x$ and $g \cdot x$ must be separated by at least $K(R+3)+1$ disjoint hyperplanes, say $J_0, \ldots, J_{K(R+3)}$. Because, for every $0 \leq i \leq K-1$, the hyperplanes $J_{i(R+3)}$ and $J_{(i+1)(R+3)}$ are separated by at least $R+1$ hyperplanes, we deduce from our previous observation that $J_{i(R+3)}$ and $J_{(i+1)(R+3)}$ cannot be indistinguishable. Thus, $\mathfrak{m}(J_0), \mathfrak{m}(J_{R+3}), \ldots, \mathfrak{m}(J_{K(R+3)})$ define $K+1$ pairwise distinct walls of $\mathcal{N}$ separating $x$ and $g \cdot x$, hence $g \notin F$. Therefore,
$$F \subset \{ g \in G \mid d(x,g \cdot x) \leq \mathrm{Ram}(K(R+3)+1) \},$$
and we conclude that $F$ is finite because the action $G \curvearrowright X$ is metrically proper. A fortiori, the action $G \curvearrowright (X, \mathcal{N})$ is metrically proper as well. 

\medskip \noindent
Finally, we claim that the action $G \curvearrowright (X, \mathcal{N})$ is cocompact, ie., there exist only finitely many orbits of (maximal) collections of pairwise transverse walls of $\mathcal{N}$. Notice that if $\mathfrak{m}(J_1)$ and $\mathfrak{m}(J_2)$ are transverse (as two elements of $\mathcal{N}$) then $J_1$ and $J_2$ must be transverse. Indeed, if $J_1$ and $J_2$ are two disjoint hyperplanes of $X$, and if we denote by $J_1^+$ (resp. $J_2^+$) the halfspace delimited by $J_1$ (resp. $J_2$) containing $J_2$ (resp. $J_1$), then either $J_1^+ \cap J_2^+ \cap G \cdot x_0= \emptyset$, so that $J_1$ and $J_2$ are indistinguishable (and $\mathfrak{m}(J_1)= \mathfrak{m}(J_2)$ in $\mathcal{N}$), or $\mathfrak{m}(J_1)$ and $\mathfrak{m}(J_2)$ are nested. Therefore, any collection of pairwise transverse walls of $\mathcal{N}$ defines a collection of pairwise transverse hyperplanes of $X$. Because we already know that $X$ contains only finitely many orbits of collections of pairwise transverse hyperplanes, the conclusion follows.
\end{proof}

\begin{question}\label{question:special}
If $G \curvearrowright X$ is geometric and virtually special, is $G \curvearrowright (G, \mathcal{N})$ virtually special as well?
\end{question}

\noindent
In order to apply the propositions proved in the previous section, we also want any two distinct elements of $G$ to be separated by some wall of $\mathcal{M}$ or $\mathcal{N}$. The following lemma gives a sufficient condition which implies this property, and moreover Lemma \ref{lem:modifycubing} below essentially states that we can always assume that this condition is satisfied. 

\begin{lemma}\label{lem:trivialstabiliser}
If the stabiliser of $x_0$ is trivial, then two distinct elements of $G$ are separated by some $\mathfrak{m}(J)$. 
\end{lemma}

\begin{proof}
Let $g,h \in G$ be two distinct elements. Because the stabiliser of $x_0$ is trivial, necessarily $g \cdot x_0 \neq h \cdot x_0$. Thus, if $J$ denotes a hyperplane separating $g \cdot x_0$ and $h \cdot x_0$, then $\mathfrak{m}(J)$ must separate $g$ and $h$.
\end{proof}

\begin{lemma}\label{lem:modifycubing}
Let $G$ be a group acting on a CAT(0) cube complex $X_0$. Then $G$ acts on a CAT(0) cube complex $X$ containing $X_0$ so that the action $G \curvearrowright X_0$ extends to an action $G \curvearrowright X$ and $X$ contains a vertex whose stabiliser is trivial. Moreover, the action $G \curvearrowright X$ is properly discontinuous (resp. metrically proper, cocompact) if and only if the action $G \curvearrowright X_0$ is properly discontinuous (resp. metrically proper, cocompact) as well; also, if $G$ is a finitely generated residually finite group and if $G \curvearrowright X_0$ is properly discontinous and virtually special, then $G \curvearrowright X$ is virtually special. 
\end{lemma}

\begin{proof}
Let $x_0 \in X_0$ be a base vertex and let $\Omega$ denote its $G$-orbit. Let $X$ be the cube complex constructed from $X_0$ by adding one vertex $(x,g)$ for every $x \in \Omega$ and $g \in \mathrm{stab}(x)$, and one edge between $x$ and $(x,g)$ for every $x \in \Omega$ and $g \in \mathrm{stab}(x)$. It is clear that $X$ satisfies the Gromov-link condition and that it is simply connected, so it is a CAT(0) cube complex. 

\medskip \noindent
Now, we extend the action $G \curvearrowright X_0$ to an action $G \curvearrowright X$. For every $x \in \Omega$, fix some $h_x \in G$ such that $h_x \cdot x_0=x$. For every $g,k \in G$ and $x \in \Omega$, define
$$g \cdot (x,k)= (gx,gk h_x h_{gx}^{-1});$$
notice that
$$gkh_x h_{gx}^{-1} \cdot gx= gkh_x \cdot x_0 = gk \cdot x = g \cdot x,$$
so that $gkh_xh_{gx}^{-1} \in \mathrm{stab}(gx)$. Moreover, 
$$\begin{array}{lcl} g_1 \cdot ( g_2 \cdot (x,k)) & = & g_1 \cdot (g_2x,g_2kh_xh_{g_2x}^{-1}) \\ \\ & = & (g_1g_2x, g_1 \cdot g_2kh_xh_{g_2x}^{-1} \cdot h_{g_2x} h_{g_1g_2x}^{-1}) \\ \\ & = & (g_1g_2x, g_1g_2 k h_x h_{g_1g_2x}^{-1}) = g_1g_2 \cdot (x,k) \end{array}$$
so we have defined a group action $G \curvearrowright X$, which extends $G \curvearrowright X_0$ by construction. 

\medskip \noindent
Fixing some $x \in \Omega$, we claim that the vertex $(x,1) \in X$ has trivial stabiliser. Indeed, if $g \in G$ fixes $(x,1)$, then $(x,1)= g \cdot (x,1)= (gx, gh_xh_{gx}^{-1})$. As a consequence, $gx=x$, ie., $g \in \mathrm{stab}(x)$, so that $h_{gx}=h_x$. Therefore, our relation becomes $(x,1)=(x,g)$, hence $g=1$. 

\medskip \noindent
This proves the first assertion of our lemma. Next, it is clear that the action $G \curvearrowright X$ is properly discontinuous (resp. metrically proper, cocompact) if and only if the action $G \curvearrowright X_0$ is properly discontinuous (resp. metrically proper, cocompact) as well. Finally, suppose that $G$ is a finitely generated residually finite group and that $G \curvearrowright X_0$ is properly discontinous and virtually special. Let $H \leq G$ be a finite-index subgroup acting specially on $X_0$. Then
$$\mathrm{Obs}(H \curvearrowright X) \subset \bigcup\limits_{x \in \Omega} \mathrm{stab}(x) = \bigcup\limits_{g \in G} g \cdot \mathrm{stab}(x_0) \cdot g^{-1}.$$
Because $G$ acts properly discontinuously on $X_0$, $\mathrm{stab}(x)$ must be finite. Therefore, $G$ contains a normal finite-index $K$ such that $K \cap \mathrm{stab}(x_0)= \{ 1 \}$. A fortiori, $K \cap \mathrm{Obs}(H \curvearrowright X)= \emptyset$, so that $H \cap K$ defines a finite-index subgroup of $G$ acting specially on $X$. Consequently, the action $G \curvearrowright X$ is virtually special.
\end{proof}

\subsection{Generalizations of spaces with walls}\label{section:generalizedwalls}

\noindent
In \cite{measuredwallspaces}, Chelix, Martin and Valette introduces \emph{spaces with measured walls}, which may be thought of as ``continuous spaces with walls''; they may also be compared to real trees generalizing simplicial trees.

\begin{definition}
A \emph{space with measured walls}\index{Spaces with measured walls} $(X, \mathcal{W}, \mathcal{B}, \mu)$ is the data of a collection of walls $\mathcal{W}$ on a set $X$, a $\sigma$-algebra $\mathcal{B}$ on $\mathcal{W}$, and a measure $\mu$ on $(\mathcal{W}, \mathcal{B})$, so that, for every $x,y \in X$, the set $\mathcal{W}(x \mid y)$ of the walls separating $x$ and $y$ belongs to $\mathcal{B}$ and has finite measure.
\end{definition}

\noindent
In particular, the definition allows us to introduce the pseudo-distance $$d_{\mathcal{W}} : (x,y) \mapsto \mu \left( \mathcal{W}(x \mid y) \right).$$ 
Spaces with measured walls were essentially introduced because of their relations with \emph{Kazhdan's property (T)} and \emph{a-T-menability} (also known as \emph{Haagerup's property}). Recall that a (discrete) group is \emph{a-T-menable}\index{A-T-menable groups} if it admits a proper action on a Hilbert space by affine isometries. In this article, we are not interested in Kazhdan's property (T) since, according to Proposition \ref{prop:quasimedianimplycubical}, any group acting without global fixed point on a quasi-median graph acts without global fixed point on a CAT(0) cube complex, and so does not satisfy Kazhdan's property (T) according to \cite{MR1459140}. The relations between a-T-menability and spaces with measured walls were strenghtened by Chatterji, Dru\c{t}u and Haglund in \cite[Theorem 1.3]{medianviewpoint} when proving that

\begin{thm}\emph{\cite{medianviewpoint}}
A (discrete) group is a-T-menable if and only if it acts metrically properly on a space with measured walls.
\end{thm}

\noindent
It is possible to generalize Proposition \ref{prop:producingproperaction} to spaces with measured walls. To do so, given a quasi-median graph $X$, define a \emph{system of measured wallspaces} by the data of a space with measured walls $(C, \mathcal{W}(C), \mathcal{B}(C), \mu_C)$ for every clique $C \subset X$; it is \emph{coherent} if the underlying system of wallspaces is coherent on its own right. Let us fix a coherent system of measured wallspaces on $X$. Recall that, for every hyperplane $J$, we set $\mathcal{W}(J)= \overline{\mathcal{W}} (C)$ where $C$ is any clique dual to $J$, so that $\mathcal{B}(C)$ naturally induces a $\sigma$-algebra $\mathcal{B}(J)$ on $\mathcal{W}(J)$, and $\mu_C$ a measure $\mu_J$ on $\mathcal{B}(J)$. On our collection of walls $\mathcal{W} = \mathcal{SW} \sqcup \mathcal{HW}$, we consider the $\sigma$-algebra 
$$\mathcal{B}= \mathcal{B}_0 \sqcup \coprod\limits_{J \ \text{hyperplane of} \ X} \mathcal{B}(J),$$
where $\mathcal{B}_0$ is the discrete $\sigma$-algebra defined on $\mathcal{SW}$. Then, we introduce the measure $\mu$ on $\mathcal{W}$ defined by
$$\mu \left( A \sqcup \coprod\limits_{J \ \text{hyperplane of} \ X} A_J \right) = \#A + \sum\limits_{J \ \text{hyperplane of} \ X} \mu_J(A_J),$$
where $A_J \in \mathcal{B}(J)$ for every $J$ and $A \in \mathcal{B}_0$. Using Lemma \ref{lem:distanceSW} and following word for word the proof of Lemma \ref{lem:distanceHW}, we show that:

\begin{lemma}\label{lem:distanceMHW}
Let $x,y \in X$ be two vertices. Let $J_1, \ldots, J_n$ denote the hyperplanes separating $x$ and $y$, and, for every $1 \leq i \leq n$, fix a clique $C_i$ dual to $J_i$. If $p_i : X \to C_i$ denotes the projection onto $C_i$ for every $1 \leq i \leq n$, then 
$$d(x,y) \leq \mu~ \mathcal{W}(x \mid y) \leq 2d(x,y)+ \sum\limits_{i=1}^n \mu_{C_i} \mathcal{W}_{C_i}(p_i(x) \mid p_i(y)).$$
\end{lemma}

\noindent
As a consequence, $(X, \mathcal{M},\mathcal{B}, \mu)$ defines a space with measured walls. Moreover, by following word for word the proof of Proposition \ref{prop:producingproperaction}, we show:

\begin{prop}\label{prop:producingproperactionmhw}
Let $G$ be a group acting on a quasi-median graph $X$ endowed with a $G$-invariant coherent system of measured wallspaces. Suppose that
\begin{itemize}
	\item any vertex of $X$ belongs to finitely many cliques;
	\item any vertex-stabiliser is finite;
	\item for every clique $C$, the pseudo-metric space $(C,d_{\mathcal{W}(C)})$ is locally finite;
	\item there exists some constant $K>0$ such that, for every clique $C$ and every distinct vertices $x,y \in C$, $\mu_{C} \mathcal{W}_C( x \mid y) \geq K$.
\end{itemize}
Then the action $G \curvearrowright (X, \mathcal{W}, \mathcal{B},\mu)$ is metrically proper.
\end{prop}

\noindent
In \cite{labelledpartitions}, Arnt introduced a generalization of spaces with measured walls, called \emph{spaces with labelled partitions}, motivated by generalizations of a-T-menability. 

\begin{definition}
A \emph{space with labelled partitions}\index{Spaces with labelled partitions} $(X, \mathcal{P}, F(\mathcal{P}))$ is the data of a set $X$, a family $\mathcal{P}$ of labelling functions from $X$ to $\mathbb{K}$ (where $\mathbb{K}= \mathbb{R}$ or $\mathbb{C}$), and $F(\mathcal{P})= \mathcal{F}(\mathcal{P})/ \mathcal{F}(\mathcal{P})_0$ a Banach space, where $(\mathcal{F}(\mathcal{P}), \| \cdot \|)$ is a semi-normed space of $\mathbb{K}$-valued functions on $\mathcal{P}$ and $\mathcal{F}(\mathcal{P})_0 = \{ \xi \in \mathcal{F}(\mathcal{P}) \mid \| \xi \| =0 \}$, such that
$$c(x,y) : \left\{ \begin{array}{ccc} \mathcal{P} & \longrightarrow & \mathbb{K} \\ f & \longmapsto & f(x)-f(y) \end{array} \right. \in \mathcal{F}(\mathcal{P})$$
for every $x,y \in X$. 
\end{definition}

\noindent
The definition allows us to define a pseudo-distance on $X$ by
$$d_{\mathcal{P}} : (x,y) \mapsto \| c(x,y) \|.$$

\begin{thm}\emph{\cite{labelledpartitions}}
Let $G$ be a (discrete) group $G$. If $G$ acts properly by affine isometries on a Banach space $B$ then there exists a structure $(G, \mathcal{P}, F(\mathcal{P}))$ of space with labelled partitions on $G$ such that $G$ acts properly on $(G, \mathcal{P},F(\mathcal{P}))$ by left-multiplication. Moreover, there exists a linear isometric embedding $F(\mathcal{P}) \hookrightarrow B$. Conversely, if $G$ acts properly on a space with labelled partitions $(X, \mathcal{P},F(\mathcal{P}))$ then there exists a proper isometric affine action of $G$ on a Banach space $B$, which is a closed subspace of $F(\mathcal{P}$. 
\end{thm}

\noindent
It is possible to generalize Proposition \ref{prop:producingproperaction} to spaces with labelled partitions. To do so, given a quasi-median graph $X$, define a \emph{system of spaces with labelled partitions} by the data of a space with labelled partitions $(C, \mathcal{P}_C, F(\mathcal{P}_C))$ for each clique $C \subset X$. Fixing some clique $C$, any function $f : C \to \mathbb{K}$ naturally extends to $\bar{f} : X \to \mathbb{K}$ by defining $\bar{f}(x)= f( \mathrm{proj}_C(x))$ for every $x \in X$. By extension, we define
$$\overline{\mathcal{P}}_C = \{ \bar{f} : X \to \mathbb{K} \mid f \in \mathcal{P}_C \} \ \text{and} \ \mathcal{F}(\overline{\mathcal{P}}_C) = \{ \bar{\xi} : \bar{f} \mapsto \xi(f) \mid \xi \in \mathcal{F}(\mathcal{P}_C) \}.$$
Next, let $C,C'$ be two cliques dual to the same hyperplane, and let $t=t_{C \to C'}$ denote the canonical bijection $C \to C'$. We define
$$\mathcal{P}_{C \to C'} = \{ f \circ t^{-1} \mid f \in \mathcal{P}_C \} \ \text{and} \ \mathcal{F}( \mathcal{P}_{C \to C'}) = \{ f \mapsto \xi (f \circ t) \mid \xi \in \mathcal{F}(\mathcal{P}_C) \}.$$
Using the map $\xi \mapsto \left( f \mapsto \xi(f \circ t) \right)$, we transfer the structure of semi-normed vector space of $\mathcal{F}(\mathcal{P}_C)$ to $\mathcal{F}(\mathcal{P}_{C \to C'})$. As a consequence, $F(\mathcal{P}_{C \to C'})= \mathcal{F}(\mathcal{P}_{C \to C'}) / \mathcal{F}(\mathcal{P}_{C \to C'})_0$ is naturally isomorphic to $F(\mathcal{P}_C)$, so that it defines a Banach space. Moreover, if $x,y \in C'$ and
$$c(x,y) : \left\{ \begin{array}{ccc} \mathcal{P}_{C \to C'} & \longrightarrow & \mathbb{K} \\ f & \longmapsto & f(x)-f(y) \end{array} \right.$$
then $c(x,y)(f)=f \circ t (t^{-1}(x)) -  f \circ t (t^{-1}(y)) = \kappa(t^{-1}(x),t^{-1}(y))(f \circ t)$ for every $f \in \mathcal{P}_{C \to C'}$, where
$$\kappa(a,b) : \left\{ \begin{array}{ccc} \mathcal{P}_C & \longrightarrow & \mathbb{K} \\ f & \longmapsto & f(a)-f(b) \end{array} \right. \in \mathcal{F}(\mathcal{P}_C)$$
for every $a,b \in C$. Consequently, $c(x,y) \in \mathcal{F}(\mathcal{P}_{C \to C'})$. 

\medskip \noindent
Thus, we have proved that $(C', \mathcal{P}_{C \to C'}, F(\mathcal{P}_{C \to C'}))$ defines a space with labelled partitions. Our system of spaces with labelled partitions is \emph{coherent} if $\mathcal{P}_{C'}= \mathcal{P}_{C \to C'}$ and $\mathcal{F}(\mathcal{P}_{C \to C'})= \mathcal{F}(\mathcal{P}_{C'})$ for every cliques $C,C'$ dual to the same hyperplane. From now on, we suppose that this is the case.

\begin{lemma}
If $C,C'$ are two cliques dual to the same hyperplane, then $\overline{\mathcal{P}}_C = \overline{\mathcal{P}}_{C'}$ and $\mathcal{F}(\overline{\mathcal{P}}_{C'}) = \mathcal{F}( \overline{\mathcal{P}}_C)$.
\end{lemma}

\begin{proof}
Let $f \in \mathcal{P}_C$. Then $f \circ t^{-1} \in \mathcal{P}_{C \to C'}$ and, for every $x \in X$,
$$\bar{f}(x)=f( \mathrm{proj}_C(x)) = f( t^{-1} ( \mathrm{proj}_{C'}(x)))= \overline{f \circ t^{-1}}(x),$$
where the second equality is justified by Lemma \ref{lem:projectionandtransfer}, hence $\bar{f} = \overline{f \circ t^{-1}} \in \overline{\mathcal{P}}_{C \to C'}$. Conversely, if $f \in \mathcal{P}_{C \to C'}$, then $f \circ t \in \mathcal{P}_C$ and, for every $x \in X$,
$$\bar{f}(x)=f( \mathrm{proj}_{C'}(x)) = f( t( \mathrm{proj}_C(x))) = \overline{f \circ t}(x),$$
hence $\bar{f} = \overline{f \circ t} \in \overline{\mathcal{P}}_C$. Thus, we have proved that $\overline{\mathcal{P}}_C = \overline{\mathcal{P}}_{C \to C'}$. On the other hand, because our system of spaces with labelled partitions is coherent, we know that $\mathcal{P}_{C'}= \mathcal{P}_{C \to C'}$, hence $\overline{\mathcal{P}}_C= \overline{\mathcal{P}}_{C \to C'} = \overline{\mathcal{P}}_{C'}$.

\medskip \noindent
Next, because $\mathcal{P}_{C'} = \mathcal{P}_{C \to C'}$, we have
$$\mathcal{F}( \overline{\mathcal{P}}_{C'}) = \left\{ \left\{ \begin{array}{ccc} \overline{\mathcal{P}}_{C'} & \to & \mathbb{K} \\ \bar{f} & \mapsto & \xi(f) \end{array} \right. \left| \xi \in \mathcal{F}( \mathcal{P}_{C'}) \right. \right\} = \left\{ \left\{ \begin{array}{ccc} \overline{\mathcal{P}}_{C'} & \to & \mathbb{K} \\ \bar{f} & \mapsto & \zeta(f \circ t) \end{array} \right. \left| \zeta \in \mathcal{F}( \mathcal{P}_{C}) \right. \right\}.$$
Now, for every $f \in \mathcal{P}_{C'}$, $\bar{f}= \overline{(f \circ t) \circ t^{-1}}$ where $f \circ t \in \mathcal{P}_C$; because we proved in the previous paragraph that $\bar{g} = \overline{g \circ t^{-1}}$ for every $g \in \mathcal{P}_C$, we deduce that $\bar{f} = \overline{f \circ t}$. Therefore, 
$$\mathcal{F}(\overline{\mathcal{P}}_{C'}) = \left\{ \left\{ \begin{array}{ccc} \mathcal{P}_C & \to & \mathbb{K} \\ \bar{g} & \mapsto & \zeta(g) \end{array} \right. \left| \zeta \in \mathcal{F}(\mathcal{P}_C) \right. \right\} = \mathcal{F}(\overline{\mathcal{P}}_C),$$
which concludes the proof of our lemma.
\end{proof}

\noindent
If $J$ is a hyperplane of $X$, this allows us to define $\mathcal{P}_J = \overline{\mathcal{P}}_C$ and $\mathcal{F}(\mathcal{P}_J)= \mathcal{F}(\overline{\mathcal{P}}_C)$ for every clique $C$ dual to $J$. Set 
$$\mathcal{P} = \mathcal{P}_0 \sqcup \coprod\limits_{J \ \text{hyperplane of} \ X} \mathcal{P}_J,$$
where $\mathcal{P}_0$ is the set of all the characteristic functions $\chi_D$ of the sectors $D \in\mathfrak{D}$ of $X$. Then, define 
$$\mathcal{F}(\mathcal{P})= \ell^p(\mathfrak{D}) \oplus \bigoplus\limits_{J \ \text{hyperplane of} \ X}^{\ell^p} \mathcal{F}(\mathcal{P}_J).$$
To conclude that $(X, \mathcal{P}, F(\mathcal{P}))$ defines a space with labelled partitions, it is sufficient to prove the following lemma by following the lines of the proof of Lemma \ref{lem:distanceHW} (we leave the proof as an exercice): 

\begin{lemma}
Let $x,y \in X$ be two vertices. Let $J_1, \ldots, J_n$ denote the hyperplanes separating $x$ and $y$, and, for every $1 \leq i \leq n$, fix a clique $C_i$ dual to $J_i$. If $p_i : X \to C_i$ denotes the projection onto $C_i$ for every $1 \leq i \leq n$, then 
$$c(x,y) = \sum\limits_{D \ \text{separates $x$ and $y$}} \chi_D + \sum\limits_{i=1}^n \kappa_i(p_i(x),p_i(y)),$$
where $\kappa_i(a,b) : f \mapsto f(a)-f(b) \in \mathcal{F}(\mathcal{P}_{C_i})$ for every $a,b \in C_i$.
\end{lemma}

\noindent
In particular, notice that 
$$d_{\mathcal{P}}(x,y) = \| c(x,y)\| = 2d(x,y) + \left( \sum\limits_{i=1}^n d_i(p_i(x),p_i(y))^p \right)^{1/p},$$
where $d_i$ denotes the pseudo-distance on $C_i$. Thus, by reproducing the proof of Proposition \ref{prop:producingproperaction}, we show:

\begin{prop}\label{prop:producingproperactionlp}
Let $G$ be a group acting on a quasi-median graph $X$ endowed with a $G$-invariant coherent system of spaces with labelled partitions. Suppose that
\begin{itemize}
	\item any vertex of $X$ belongs to finitely many cliques;
	\item any vertex-stabiliser is finite;
	\item for every clique $C$, the pseudo-metric space $(C, d_{\mathcal{P}(C)})$ is locally finite;
	\item there exists some constant $K>0$ such that, for every clique $C$ and every distinct vertices $x,y \in C$, $d_{\mathcal{P}(C)}(x,y) \geq K$.
\end{itemize}
Then the action $G \curvearrowright (X, \mathcal{P}, F(\mathcal{P}))$ is metrically proper.
\end{prop}

\noindent
We emphasize the fact that, in the conclusion of the previous proposition, the space $F(\mathcal{P})$ is obtained as a $\ell^p$-sums of local $F(\mathcal{P}_C)$'s with some $\ell^p$-space, where $p \geq 1$ is an arbitrary integer. Consequently, as a consequence of the statement below, our construction applies naturally to prove the \emph{a-$L^p$-menability} of some groups:

\begin{cor}\emph{\cite{labelledpartitions}}
Let $G$ be a (discrete) group and $p \geq 1$ with $p \notin 2 \mathbb{Z} \backslash \{ 2\}$. Then $G$ is \emph{a-$L^p$-menable}\index{A-$L^p$-menable groups} if and only if $G$ acts properly on a space with labelled partitions $(X, \mathcal{P},F(\mathcal{P})$ where $F(\mathcal{P})$ is isometrically isomorphic to a closed subspace of an $L^p$-space. 
\end{cor}

\noindent
Although in the following we restrict ourself to a-$L^p$-menability, we emphasize that our arguments hold for other classes of Banach spaces.

\section{Topical actions on quasi-median graphs I}\label{section:topicalactionsI}

\noindent
In this section, we show that, if a group $G$ acts \emph{topically-transitively} on a quasi-median graph $X$, then we can choose metrics or collections of walls on clique-stabilisers to produce coherent and $G$-invariant systems of metrics or wallspaces. Next, this construction is used to prove combination theorems of the form: Fix a group property $\mathcal{P}$ and suppose that a given group $G$ acts topically-transitively on some quasi-median graph $X$ satisfying some finiteness properties (depending on the property $\mathcal{P}$ we are looking at). If clique-stabilisers satisfy $\mathcal{P}$, then so does $G$.

\subsection{Creating invariant coherent systems I}

\noindent
In this section, we are interested in the following problem. Fix a group $G$ acting on a quasi-median graph $X$ and a collection of cliques $\mathcal{C}$ such that every $G$-orbit of hyperplanes intersects $\mathcal{C}$ along a single clique. Roughly speaking, $\mathcal{C}$ is a collection of cliques of reference. On each clique $C \in \mathcal{C}$, fix a metric $\delta_C$ or a collection of walls $\mathcal{W}(C)$. Is it possible to extend these data to a coherent and $G$-invariant system of metrics or wallspaces? The main result of this section, namely Theorem \ref{thm:extensionunique}, completely answers this question.

\medskip \noindent
First, we need to introduce some vocabulary. Let $G$ be a group acting on a quasi-median graph $X$. Fix some clique $C$ of $X$ and let $J$ denote the hyperplane dual to $C$. For every $x \in X$, recall that $[C,x]$ denotes the sector delimited by $J$ which contains $x$. Define
$$\rho_C : \left\{ \begin{array}{ccc} \mathrm{stab}(J) & \to & \mathrm{Bij}(C) \\ g & \mapsto & (x \mapsto [C,gx] \cap C) \end{array} \right.$$
Otherwise saying, we look at the permutation of the vertices of $C$ corresponding to the permutation of the sectors delimited by $J$ which is induced by $g \in \mathrm{stab}(J)$. If $C'$ is a clique dual to the same hyperplane as $C$, we denote by $\rho_{C' \to C} : \mathrm{stab}(C') \to \mathrm{Bij}(C)$ the restriction of $\rho_C$ to $\mathrm{stab}(C')$. 

\begin{thm}\label{thm:extensionunique}
Let $G$ be a group acting on a quasi-median graph $X$. Fix some collection of cliques $\mathcal{C}$ such that every $G$-orbit of hyperplanes intersects it along a single clique. Suppose that each clique of $\mathcal{C}$ is endowed with a metric $\delta_C$ (resp. a structure of wallspace $\mathcal{W}(C)$). There exists a $G$-invariant coherent system of metrics (resp. wallspaces) extending $\{ (C,\delta_C) \mid C \in \mathcal{C} \}$ (resp. $\{ (C, \mathcal{W}(C)) \mid C \in \mathcal{C} \}$) if and only if, for every $C \in \mathcal{C}$, the metric $\delta_C$ (resp. the collection of walls $\mathcal{W}(C)$) is $\mathrm{Im} (\rho_C)$-invariant. Moreover, the extension, when it exists, is unique. 
\end{thm}

\noindent
We begin by proving a few preliminary results about the maps $\rho_C$. 

\begin{lemma}\label{lem:rho2}
Let $X$ be a quasi-median graph and $C,C'$ two cliques dual to a common hyperplane $J$. The equality $$\rho_{C'}(g) \cdot t_{C \to C'}(x) = t_{C \to C'} ( \rho_C(g) \cdot x )$$
holds for every $g \in \mathrm{stab}(J)$ and every $x \in C$.
\end{lemma}

\begin{proof}
$\rho_{C'}(g) \cdot t_{C\to C'}(x)$ is the unique vertex of $C'$ contained in the sector delimited by $J$ containing $g \cdot t_{C \to C'}(x)$. Since $x$ and $t_{C \to C'}(x)$ belong to the same sector delimited by $J$, a fortiori $\rho_{C'}(g) \cdot t_{C \to C'}(x)$ is also the unique vertex of $C'$ contained in the sector delimited by $J$ containing $g \cdot x$. Next, since $\rho_{C'}(g) \cdot t_{C \to C'}(x)$ and $t_{C' \to C} \left( \rho_{C'}(g) \cdot t_{C \to C'}(x) \right)$ belong to the same sector delimited by $J$, we deduce that $t_{C' \to C} \left( \rho_{C'}(g) \cdot t_{C \to C'}(x) \right)$ is the unique vertex of $C$ contained in the sector delimited by $J$ containing $g \cdot x$, hence $\rho_C(g) \cdot x = t_{C' \to C} \left( \rho_{C'}(g) \cdot t_{C \to C'}(x)\right)$. We conclude that $\rho_{C'}(g) \cdot t_{C \to C'}(x) = t_{C \to C'} ( \rho_C(g) \cdot x )$.
\end{proof}

\begin{lemma}\label{lem:rho?}
Let $G$ be a group acting on a quasi-median graph $X$ and $C$ a clique of $X$. Denote by $J$ the hyperplane dual to $C$. The equality
$$\rho_C(g) \cdot x = t_{gC \to C}(g \cdot x)$$
holds for every vertex $x \in C$ and every element $g \in \mathrm{stab}(J)$.
\end{lemma}

\begin{proof}
Fix a vertex $x \in C$ and an element $g \in \mathrm{stab}(J)$. By definition of $\rho_C$, $\rho_C(g) \cdot x$ is the unique vertex of $C$ which belongs to the same fiber of $J$ as $g \cdot x$. This description also defines $t_{gC \to C}(g \cdot x)$, so the conclusion follows. 
\end{proof}

\begin{cor}\label{cor:rho1}
Let $G$ be a group acting on some quasi-median graph $X$ and $J$ a hyperplane dual to some clique $C$. The equality 
$$k \cdot t_{kC \to C}(x) = \rho_{kC}(k) \cdot x$$
holds for every $k \in \mathrm{stab}(J)$ and $x \in kC$.
\end{cor}

\begin{proof}
Fix some $k \in \mathrm{stab}(J)$ and $x \in kC$. By applying the previous lemma, we deduce that
$$\rho_{kC}(k) \cdot x = t_{k^2C \to kC}(k \cdot x) = k \cdot t_{kC \to C}(x),$$
and the conlusion follows. 
\end{proof}

\noindent
The maps constructed by our next preliminary result will allow us to link an arbitrary clique of our quasi-median graph to the clique of $\mathcal{C}$ which \emph{labels} it, in the following meaning. 

\begin{definition}
Let $G$ be a group acting on a quasi-median graph $X$ and $\mathcal{C}$ a collection of cliques such that any $G$-orbit of hyperplanes intersects it along a single clique. A clique $C$ is \emph{labelled} by $Q \in \mathcal{C}$ if $Q$ is the unique clique of $\mathcal{C}$ dual to the same hyperplane as one of the translates of $C$.
\end{definition}

\noindent
We will use these maps to transfer structures (metrics or collections of walls) from the cliques of $\mathcal{C}$ to all the cliques of the quasi-median graph.

\begin{prop}\label{prop:topicalI0}
Let $G$ be a group acting on a quasi-median graph $X$. Fix a collection of cliques $\mathcal{C}$ such that every $G$-orbit of hyperplanes intersects $\mathcal{C}$ along a single clique. There exist
\begin{itemize}
	\item a morphism $\psi_{C} : \mathrm{stab}(Q) \to \mathrm{Im}(\rho_C)$;
	\item a bijection $\Phi_{C} : Q \to C$;
\end{itemize}
for every clique $C$ labelled by some $Q \in \mathcal{C}$, such that
\begin{itemize}
	\item $\Phi_C$ is $\psi_C$-equivariant;
	\item if $C,C'$ are two cliques dual to the same hyperplane, $\Phi_{C'}= t_{C \to C'} \circ \Phi_C$;
	\item for every clique $C$ and every $g \in G$, there exists some $h \in \mathrm{Im}(\rho_Q)$ such that $g \cdot \Phi_C(x)= \Phi_{gC}(h \cdot x)$ for every $x \in Q$.
\end{itemize}
Moreover, the morphism $\psi_{C}$ is, up to conjugacy, the restriction $\rho_{Q \to C'}$ to $\mathrm{stab}(Q)$ of $\rho_{C'}$, where $C'$ is a $G$-translate of $C$. 
\end{prop}

\begin{proof}
For every hyperplane $J$ of $X$, fix some $h(J) \in G$ such that $h(J) \cdot J$ contains a clique of $\mathcal{C}$. By extension, if $C$ is a clique of $J$, we set $h(C)=h(J)$. Now, define
$$\Phi_C : \left\{ \begin{array}{ccc} Q & \to & C \\ x & \mapsto & h(C)^{-1} \cdot t_{Q \to h(C)C}(x) \end{array} \right. \ \text{and} \ \psi_C : \left\{ \begin{array}{ccc} \mathrm{stab}(Q) & \to & \mathrm{Im}(\rho_C) \\ g & \mapsto & h(C)^{-1} \rho_{h(C)C}(g) h(C) \end{array} \right.$$
Notice that, for every $g \in \mathrm{stab}(Q)$ and $x \in Q$, Lemma \ref{lem:rho2} implies that
$$\begin{array}{lcl} \Phi_C(gx) & = & h(C)^{-1} t_{Q \to h(C)C}(gx) = h(C)^{-1} \rho_{h(C)C}(g) t_{Q \to h(C)C}(x) \\ \\ & = & \psi_C(g)h(C)^{-1} t_{Q \to h(C)C}(x) = \psi_C(g) \cdot \Phi_C(x) \end{array}$$
Therefore, $\Phi_C$ is $\psi_C$-equivariant. Next, for every $g \in G$ and $x \in C$, by applying Corollary \ref{cor:rho1} proved below and Lemma \ref{lem:rho2}, we deduce that
$$\begin{array}{lcl} g \Phi_C(x) & = & gh(C)^{-1} t_{Q \to h(C)C}(x) \\ \\ & = & h(gC)^{-1} \left( h(gC)gh(C)^{-1} \right) t_{h(gC)gC \to h(C)C} \circ t_{Q \to h(gC)gC}(x) \\ \\ & = & h(gC)^{-1} \rho_{h(gC)gC}(h(gC)gh(C)^{-1}) t_{Q \to h(gC)gC}(x) \\ \\ & = & h(gC)^{-1} t_{Q \to h(gC)gC} \left( \rho_Q(h(gC)gh(C)^{-1}) \cdot x \right) \\ \\ & = & \Phi_{gC} \left( \rho_Q \left( h(gC)gh(C)^{-1} \right) \cdot x \right) \end{array}$$
Thus, we have proved the third point of our proposition. Finally, if $C'$ is a clique dual to the same hyperplane as $C$, then, by noticing that $h(C)=h(J)=h(C')$, we deduce that
$$\begin{array}{lcl} \Phi_{C'}(x) & = & h(C')^{-1} t_{Q \to h(C')C'}(x) \\ \\ & = & h(J)^{-1} t_{h(J)C \to h(J)C'} \circ t_{Q \to h(J)C} (x) \\ \\ & = & t_{C \to C'} \left( h(C)^{-1} t_{Q \to h(C)C}(x) \right) = t_{C \to C'} \circ \Phi_C(x) \end{array}$$
for every $x \in C'$, hence $\Phi_{C'}= t_{C \to C'} \circ \Phi_C$. 
\end{proof}

\noindent
Now, we are ready to prove the existence claimed by Theorem \ref{thm:extensionunique}. 

\begin{cor}\label{cor:strongsystwall0}
Let $G$ be a group acting on a quasi-median graph $X$. Fix a collection of cliques $\mathcal{C}$ such that every $G$-orbit of hyperplanes intersects $\mathcal{C}$ along a single clique, and, for every $C \in \mathcal{C}$, let $\mathcal{W}(C)$ be an $\mathrm{Im}(\rho_C)$-invariant collection of walls making $(C, \mathcal{W}(C))$ a wallspace. Then there exists a $G$-invariant coherent system of wallspaces extending $\{ (C, \mathcal{W}(C)) \mid C \in \mathcal{C} \}$ such that, for every clique $C$, labelled by some $Q \in \mathcal{C}$, the spaces $(C, \mathcal{W}(C))$ and $(Q, \mathcal{W}(Q))$ are isomorphic. \\
Moreover, under the additional assumptions that $\mathrm{fix}(C)= \{1\}$ and that $\rho_{Q \to C'} : \mathrm{stab}(Q) \to \mathrm{Im}(\rho_{C'})$ is an isomorphism for every translate $C'$ of $C$ which is dual to the same hyperplane as $Q$, the actions $\mathrm{stab}(C) \curvearrowright (C, \mathcal{W}(C))$ and $\mathrm{stab}(Q) \curvearrowright (Q, \mathcal{W}(Q))$ turn out to be isomorphic.
\end{cor}

\begin{proof}
Let $\Phi_C$ and $\psi_C$ be the maps given by Proposition \ref{prop:topicalI0}. We set
$$\mathcal{W}(C) = \Phi_C ~ \mathcal{W}(Q)$$
for every clique $C$ labelled by $Q \in \mathcal{C}$. If $C'$ is a clique dual to the same hyperplane as $C$, then 
$$\mathcal{W}(C')= \Phi_{C'} \mathcal{W}(Q)= t_{C \to C'} \circ \Phi_C \mathcal{W}(Q)= t_{C \to C'} \mathcal{W}(C) = \mathcal{W}(C \to C').$$
Therefore, our system is coherent. Next, if $g \in G$, then
$$g \mathcal{W}(C)= g \Phi_C \mathcal{W}(Q)= \Phi_{gC} \left( h \cdot \mathcal{W}(Q) \right)= \Phi_{gC} \mathcal{W}(Q) = \mathcal{W}(gC),$$
where $h \in \mathrm{Im}(\rho_Q)$ is the element given by the third point of Proposition \ref{prop:topicalI0}. Notice that the third equality is justified by the fact that $\mathcal{W}(Q)$ is $\mathrm{Im}(\rho_Q)$-invariant. Consequently, our system is also $G$-invariant. 

\medskip \noindent
Finally, we claim that, if $\mathrm{fix}(C)= \{1\}$ and if $\rho_{Q \to C'}$ is an isomorphism for every translate $C'$ of $C$ which is dual to the same hyperplane as $Q$, then $\rho_{C \to C}$ defines an isomorphism $\mathrm{stab}(C) \to \rho_C(\mathrm{stab}(C))$ and $\Phi_C$ is a $\left( \rho_{C \to C}^{-1} \circ \psi_C \right)$-equivariant isomorphism $(Q, \mathcal{W}(\mathcal{Q})) \to (C, \mathcal{W}(C))$. The first assertion follows from the observation that $\mathrm{ker}(\rho_{C \to C}) = \mathrm{fix}(C)$. Next, we know that $\Phi_C$ is $\rho_{C \to C}^{-1} \circ \psi_C$-equivariant because $\Phi_C$ is already $\psi_C$-equivariant according to Proposition \ref{prop:topicalI0} and that $\psi_C(g)$ and $\rho_{C \to C}^{-1} \circ \psi_C(g)$ induce the same permutation on $C$ for every $g \in \mathrm{stab}(Q)$. Finally, $\rho_{C \to C}^{-1} \circ \psi_C$ turns out to be an isomorphism since $\psi_C$ is itself an isomorphism as a consequence of our second assumption combined with the last assertion of Proposition \ref{prop:topicalI0}. A fortiori, the actions $\mathrm{stab}(C) \curvearrowright (C, \mathcal{W}(C))$ and $\mathrm{stab}(Q) \curvearrowright (Q, \mathcal{W}(Q))$ are isomorphic.
\end{proof}

\begin{cor}\label{cor:strongsystdist0}
Let $G$ be a group acting on a quasi-median graph $X$. Fix a collection of cliques $\mathcal{C}$ such that every $G$-orbit of hyperplanes intersects $\mathcal{C}$ along a single clique, and, for every $C \in \mathcal{C}$, let $\delta_C$ be a $\mathrm{Im}(\rho_C)$-invariant metric on $C$. Then there exists a $G$-invariant coherent system of metrics $\{ (C, \delta_C) \mid C \in \mathcal{C} \}$ such that, for every clique $C$, labelled by some $Q \in \mathcal{C}$, the spaces $(C, \delta_C)$ and $(Q, \delta_Q)$ are isometric. \\
Moreover, under the additional assumptions that $\mathrm{fix}(C)= \{1\}$ and that $\rho_{Q \to C'} : \mathrm{stab}(Q) \to \mathrm{Im}(\rho_{C'})$ is an isomorphism for every translate $C'$ of $C$ which is dual to the same hyperplane as $Q$, the actions $\mathrm{stab}(C) \curvearrowright (C, \delta_C)$ and $\mathrm{stab}(Q) \curvearrowright (Q, \delta_Q)$ turn out to be isomorphic.
\end{cor}

\begin{proof}
Let $\Phi_C$ and $\psi_C$ be the maps given by Proposition \ref{prop:topicalI0}. We define
$$\delta_C : (x,y) \mapsto \delta_Q(\Phi_C^{-1}(x), \Phi_C^{-1}(y))$$
for every clique $C$ labelled by $Q \in \mathcal{C}$. If $C'$ is a clique dual to the same hyperplane as $C$, then  
$$\begin{array}{lcl} \delta_{C'}(x,y) & = & \delta_{Q}( \Phi_{C'}^{-1}(x), \Phi_{C'}^{-1}(y)) = \delta_{Q}( \Phi_C^{-1} \circ t_{C' \to C} (x), \Phi_C^{-1} \circ t_{C' \to C}(y)) \\ \\ & = & \delta_C( t_{C' \to C}(x), t_{C' \to C}(y)) = \delta_{C \to C'}(x,y) \end{array}$$
for every $x,y \in C'$. Therefore, our system is coherent. Next, if $g \in G$, then for every $x,y \in C$ we have
$$\begin{array}{lcl} \delta_{gC}(gx,gy) & = & \delta_Q \left( \Phi_{gC}^{-1}(gx), \Phi_{gC}^{-1}(gy) \right) = \delta_Q \left( h \cdot \Phi_C^{-1}(x), h \cdot \Phi_{C}^{-1}(y) \right) \\ \\ & = & \delta_Q \left(\Phi_C^{-1}(x),\Phi_{C}^{-1}(y) \right) = \delta_C(x,y) \end{array}$$
where $h \in \mathrm{Im}(\rho_Q)$ is the element given by the third point of Proposition \ref{prop:topicalI0}. Notice that the third equality is justified by the fact that $\delta_Q$ is $\mathrm{Im}(\rho_Q)$-invariant. Consequently, our system is also $G$-invariant. 

\medskip \noindent
Finally, we claim that, if $\mathrm{fix}(C)= \{1\}$ and if $\rho_{Q \to C'}$ is an isomorphism for every translate $C'$ of $C$ which is dual to the same hyperplane as $Q$, then $\rho_{C \to C}$ defines an isomorphism $\mathrm{stab}(C) \to \rho_C(\mathrm{stab}(C))$ and $\Phi_C$ is a $\rho_{C \to C}^{-1} \circ \psi_C$-equivariant isomorphism $(Q, \delta_Q) \to (C, \delta_C)$. The first assertion follows from the observation that $\mathrm{ker}(\rho_{C \to C}) = \mathrm{fix}(C)$. Next, we know that $\Phi_C$ is $\left( \rho_{C \to C}^{-1} \circ \psi_C \right)$-equivariant because $\Phi_C$ is already $\psi_C$-equivariant according to Proposition \ref{prop:topicalI0} and that $\psi_C(g)$ and $\rho_{C \to C}^{-1} \circ \psi_C(g)$ induce the same permutation on $C$ for every $g \in \mathrm{stab}(Q)$. Finally, $\rho_{C \to C}^{-1} \circ \psi_C$ turns out to be an isomorphism since $\psi_C$ is itself an isomorphism as a consequence of our second assumption combined with the last assertion of Proposition \ref{prop:topicalI0}. A fortiori, the actions $\mathrm{stab}(C) \curvearrowright (C, \delta_C)$ and $\mathrm{stab}(Q) \curvearrowright (Q, \delta_Q)$ are isomorphic.
\end{proof}

\begin{proof}[Proof of Theorem \ref{thm:extensionunique}.]
The existence claimed by our theorem is proved by the previous two corollaries. Now, suppose that there exists a $G$-invariant coherent system of metrics $\{ (C,\delta_C) \mid \text{$C$ clique} \}$ extending $\{ (C,\delta_C) \mid C \in \mathcal{C} \}$. Fixing a clique $C \in \mathcal{C}$, we claim that $\delta_C$ is $\mathrm{Im}(\rho_C)$-invariant, ie., 
$$\delta_C( \rho_C(g) \cdot x , \rho_C(g) \cdot y) = \delta_C(x,y)$$
for every $x,y \in C$ and $g \in \mathrm{stab}(J)$ where $J$ denotes the hyperplane dual to $C$. Indeed, it follows from the Lemma \ref{lem:rho?} and the fact that our system of metrics is $G$-invariant and coherent that
$$\begin{array}{lcl} \delta_C( \rho_C(g) \cdot x , \rho_C(g) \cdot y) & =& \delta_C \left( t_{gC \to C}(g \cdot x), t_{gC \to C} (g \cdot y) \right) \\ \\ & = & \delta_{gC}(g \cdot x , g \cdot y) = \delta_C(x,y) \end{array}$$
Next, we claim that $\{ (C,\delta_C) \mid \text{$C$ clique} \}$ is uniquely determined by $\{ (C, \delta_C) \mid C \in \mathcal{C} \}$. So let $C$ be an arbitrary clique of $X$ and $x,y \in C$ two vertices. By definition of $\mathcal{C}$, there exists some $g \in G$ and $Q \in \mathcal{C}$ such that $Q$ and $gC$ are dual to the same hyperplane. Now, because our system is $G$-invariant and coherent, we know that
$$\delta_C(x,y)= \delta_{gC}(gx,gy)= \delta_Q(t_{gC \to Q}(gx),t_{gC \to Q}(gy)).$$
Therefore, the metric $\delta_C$ is completely determined by $\delta_Q$.

\medskip \noindent
We argue similarly for wallspaces. Suppose that there exists a $G$-invariant coherent system of wallsapces $\{ (C, \mathcal{W}(C)) \mid \text{$C$ clique} \}$ extending $\{ (C, \mathcal{W}(C)) \mid C \in \mathcal{C} \}$. Fixing a clique $C \in \mathcal{C}$, we claim that $\mathcal{W}(C)$ is $\mathrm{Im}(\rho_C)$-invariant. Indeed,  it follows from the Lemma \ref{lem:rho?} and the fact that our system of metrics is $G$-invariant and coherent that
$$\rho_C(g) \mathcal{W}(C) = t_{gC \to C} \left( g \mathcal{W}(C) \right) = t_{gC \to C} \mathcal{W}(gC) = \mathcal{W}(C)$$
for every $g \in \mathrm{stab}(J)$ where $J$ denotes the hyperplane dual to $C$. Next, we claim that $\{ (C, \mathcal{W}(C)) \mid \text{$C$ clique} \}$ is uniquely determined by $\{ (C, \mathcal{W}(C)) \mid C \in \mathcal{C} \}$. So let $C$ be an arbitrary clique of $X$. By definition of $\mathcal{C}$, there exists some $g \in G$ and $Q \in \mathcal{C}$ such that $Q$ and $gC$ are dual to the same hyperplane. Now, because our system is $G$-invariant and coherent, we know that
$$\mathcal{W}(C) = g^{-1} \mathcal{W}(gC)= g^{-1} \cdot t_{Q \to gC} \left( \mathcal{W}(Q) \right).$$
A fortiori, $\mathcal{W}(C)$ is completely determined by $\mathcal{W}(Q)$. 
\end{proof}

\subsection{Topical and topical-transitive actions}

\noindent
Given a group acting on a quasi-median graph, we have shown in the previous section that $\mathrm{Im}(\rho_C)$-invariance is a necessary and sufficient condition to be able to extend some partial system of metrics or wallspaces to a global coherent and invariant sytem. In this section, our goal is to notice that, if our group acts in a specific way, then this condition may be weakened into a $\mathrm{stab}(C)$-invariance. The first kind of actions we are interested in is:

\begin{definition}
Let $G$ be a group acting on a quasi-median graph $X$. The action $G \curvearrowright X$ is \emph{$\mathcal{C}$-topical}\index{Topical actions} if $\mathcal{C}$ is a collection of cliques such that every $G$-orbit of hyperplanes intersects it on a single clique, and if $\mathrm{Im}(\rho_C) \subset \rho_C(\mathrm{stab}(C))$ for every $C \in \mathcal{C}$.
\end{definition}

\noindent
It is worth noticing that, if $G$ is a group acting $\mathcal{C}$-topically on a quasi-median graph $X$ and if $\delta_C$ is a $\mathrm{stab}(C)$-invariant metric defined on some clique $C$ of $X$, then $\delta_C$ is automatically $\mathrm{Im}(\rho_C)$-invariant. Of course, the same statement holds if metrics are replaced with collections of walls. As a consequence, Corollaries \ref{cor:strongsystwall0} and \ref{cor:strongsystdist0} immediately imply:

\begin{prop}\label{prop:strongsystwall}
Let $G$ be a group with a $\mathcal{C}$-topical action on a quasi-median graph $X$. For every $C \in \mathcal{C}$, let $\mathcal{W}(C)$ be a $\mathrm{stab}(C)$-invariant collection of walls such that $(C, \mathcal{W}(C))$ is a space with walls. Then there exists a $G$-invariant coherent system of wallspaces extending $\{ (C, \mathcal{W}(C)) \mid C \in \mathcal{C} \}$ such that, for every clique $C$, labelled by some $Q \in \mathcal{C}$, the spaces $(C, \mathcal{W}(C))$ and $(Q, \mathcal{W}(Q))$ are isomorphic. \\
Moreover, under the additional assumptions that $\mathrm{fix}(C)=\{1\}$ and that $\rho_{Q \to C'}$ is an isomorphism for every translate $C'$ of $C$ which is dual to the same hyperplane as $Q$, the actions $\mathrm{stab}(C) \curvearrowright (C, \mathcal{W}(C))$ and $\mathrm{stab}(Q) \curvearrowright (Q, \mathcal{W}(Q))$ turn out to be isomorphic.
\end{prop}

\begin{prop}\label{prop:strongsystdist}
Let $G$ be a group with a $\mathcal{C}$-topical action on a quasi-median graph $X$. For every $C \in \mathcal{C}$, let $\delta_C$ be a $\mathrm{stab}(C)$-invariant metric on $C$. Then there exists a $G$-invariant coherent system of metrics $\{ (C, \delta_C) \mid C \in \mathcal{C} \}$ such that, for every clique $C$, labelled by some $Q \in \mathcal{C}$, the spaces $(C, \delta_C)$ and $(Q, \delta_Q)$ are isometric.\\
Moreover, under the additional assumptions that $\mathrm{fix}(C)=\{1\}$ and that $\rho_{Q \to C'}$ is an isomorphism for every translate $C'$ of $C$ which is dual to the same hyperplane as $Q$, the actions $\mathrm{stab}(C) \curvearrowright (C, \delta_C)$ and $\mathrm{stab}(Q) \curvearrowright (Q, \delta_Q)$ turn out to be isomorphic.
\end{prop}

\noindent
The second kind of actions we are interested in, namely \emph{topical-transitive} actions, is the main type of actions which we will study. (In fact, topical actions will only be used in Section  \ref{section:topicalactionsII}, by noticing that, when a group acts topically-transitively on some quasi-median graph, then the induced action on a quasi-median graph obtained by inflating the hyperplanes of the previous graph turns out to be topical.) 

\begin{definition}
Let $G$ be a group acting on a quasi-median graph $X$. The action is \emph{topical-transitive}\index{Topical-transitive actions} if $\mathrm{Im}(\rho_C) \subset \rho_C(\mathrm{stab}(C))$ for every clique $C$ of $X$ (ie., for every hyperplane $J$, every clique $C$ dual to $J$ and every $g \in \mathrm{stab}(J)$, there exists some $h \in \mathrm{stab}(C)$ such that $g$ and $h$ induce the same permutation on the set of fibers delimited by $J$) and if, for every clique $C$ of $X$,
\begin{itemize}
	\item either $C$ is finite and $\mathrm{stab}(C)= \mathrm{fix}(C)$;
	\item or $\mathrm{stab}(C) \curvearrowright C$ is free and transitive on the vertices.
\end{itemize}
If $\mathcal{C}$ denotes a collection of cliques such that any $G$-orbit of hyperplanes intersects it along a single clique, then we decompose it as the disjoint union $\mathcal{C}_1 \sqcup \mathcal{C}_2$ where $\mathcal{C}_1$ denotes the subcollection of cliques $C \in \mathcal{C}$ such that $\mathrm{stab}(C) \curvearrowright C$ is transitive and free on the vertices, and where $\mathcal{C}_2$ denotes the subcollection of cliques $C \in \mathcal{C}$ such that $C$ is finite and satisfies $\mathrm{stab}(C)= \mathrm{fix}(C)$. 
\end{definition}

\noindent
The next lemma will allow us to simplify the previous two propositions. 

\begin{lemma}\label{lem:rho3}
Let $G$ be a group acting topically-transitively on a quasi-median graph $X$ and let $C,C'$ be two cliques dual to the same hyperplane $J$. If $\mathrm{fix}(C)= \{1 \}$, then the restriction $\rho_{C \to C'} : \mathrm{stab}(C) \to \mathrm{Im}(\rho_{C'})$ of $\rho_{C'}$ is an isomorphism.
\end{lemma}

\begin{proof}
Let $h \in \mathrm{stab}(J)$. Because $G$ acts topically-transitively on $X$, there exists some $g \in \mathrm{stab}(C)$ which induces the same permutation on the set of fibers of $J$ as $h$. Consequently, $\rho_{C \to C'}(g) = \rho_{C'}(h)$. Therefore, $\rho_{C \to C'}$ is surjective. Next, notice that $\mathrm{ker}(\rho_{C\to C'})$ is the set of elements of $\mathrm{stab}(C)$ inducing a trivial permutation on the fibers of $J$, ie., $\mathrm{ker}(\rho_{C \to C'})= \mathrm{fix}(C)$. The conclusion follows.
\end{proof}

\noindent
Thus, Corollaries \ref{cor:strongsystwall0} and \ref{cor:strongsystdist0} immediately imply:

\begin{prop}\label{prop:wtopicalI}
Let $G$ be a group with a topical-transitive action on a quasi-median graph $X$. Fix a collection of cliques $\mathcal{C}$ such that any $G$-orbit of hyperplanes intersects it along a single hyperplane. For every $C \in \mathcal{C}$, let $\mathcal{W}(C)$ be a $\mathrm{stab}(C)$-invariant collection of walls making $(C, \mathcal{W}(C))$ a wallspace. Then there exists a $G$-invariant coherent system of wallspaces extending $\{ (C, \mathcal{W}(C)) \mid C \in \mathcal{C} \}$ such that, for every clique $C$, labelled by some $Q \in \mathcal{C}$ satisfying $\mathrm{fix}(Q)= \{1\}$, the actions $\mathrm{stab}(C) \curvearrowright (C, \mathcal{W}(C))$ and $\mathrm{stab}(Q) \curvearrowright (Q, \mathcal{W}(Q))$ are isomorphic.
\end{prop}

\begin{prop}\label{prop:mtopicalI}
Let $G$ be a group with a topical-transitive action on a quasi-median graph $X$. Fix a collection of cliques $\mathcal{C}$ such that any $G$-orbit of hyperplanes intersects it along a single hyperplane. For every $C \in \mathcal{C}$, let $\delta_C$ be a $\mathrm{stab}(C)$-invariant metric on $C$. Then there exists a unique $G$-invariant coherent system of metrics extending $\{ (C, \delta_C) \mid C \in \mathcal{C} \}$ such that, for every clique $C$, labelled by some $Q \in \mathcal{C}$ satisfying $\mathrm{fix}(Q)= \{1\}$, the actions $\mathrm{stab}(C) \curvearrowright (C, \delta_{C})$ and $\mathrm{stab}(Q) \curvearrowright (Q, \delta_Q)$ are isomorphic.
\end{prop}

\subsection{Relatively hyperbolic groups acting on quasi-median graphs}

\noindent
In this section, we are interested in determining when a group acting topically-transitively on a quasi-median graph is hyperbolic. In this situation, conditions on the geometry of the quasi-median graph are necessary. Indeed, the direct product of two infinite groups acts topically-transitively on a prism, but is never hyperbolic. In fact, our conclusion will be stronger: we prove that a group acting topically-transitively on a suitable hyperbolic quasi-median graph will be hyperbolic relatively to its clique-stabilisers. In the sequel, we use the definition of relatively hyperbolic groups given by Bowditch in \cite{relativelyhyperbolic}. 

\begin{definition}\label{def:relativehyp}
A finitely generated group $G$ is \emph{hyperbolic relatively to a collection of subgroups $\mathcal{H}=\{ H_1, \ldots, H_n \}$}\index{Relatively hyperbolic groups} if $G$ acts by isometries on a graph $\Gamma$ such that:
\begin{itemize}
	\item $\Gamma$ is hyperbolic;
	\item $\Gamma$ contains finitely many orbits of edges;
	\item edge-stabilisers are finite;
	\item each vertex-stabilizer is either finite or is conjugated to some $H_i$;
	\item any $H_i$ stabilises a vertex;
	\item $\Gamma$ is \emph{fine}, ie., any edge belongs only to finitely many simple loops (or \emph{cycle}) of any given length.
\end{itemize}
A subgroup conjugated to some $H_i$ is \emph{peripheral}. $G$ is just said \emph{relatively hyperbolic} if it is relatively hyperbolic with respect to a finite collection of proper subgroups.
\end{definition}

\noindent
The main criterion of this section is the following:

\begin{thm}\label{thm:qmrelativelyhyp}
Let $G$ be a group acting topically-transitively on a hyperbolic quasi-median graph $X$. Suppose that
\begin{itemize}
	\item every vertex belongs to finitely many cliques;
	\item vertex-stabilisers are finite and clique-stabilisers are finitely generated;
	\item $X$ contains finitely many $G$-orbits of cliques;
	\item the hyperplanes dual to two infinite cliques cannot be transverse;
	\item carriers of hyperplanes dual to infinite cliques are cubically finite.
\end{itemize}
Then $G$ is hyperbolic relatively to the stabilisers of hyperplanes dual to infinite cliques.
\end{thm}

\noindent
In order to apply our definition of relatively hyperbolic groups, we need to find a fine hyperbolic graph on which $G$ acts. A classical construction for this purpose is to \emph{cone-off} the space on which our group originally acts. 

\begin{definition}
Let $X$ be a graph and $\mathcal{Q}$ be a collection of subgraphs.
\begin{itemize}
	\item The \emph{cone-off of $X$ over $\mathcal{Q}$}\index{Cone-offs} is the graph obtained from $X$ by adding an edge between two vertices whenever they both belong to a common subgraph of $\mathcal{Q}$.
	\item The \emph{usual cone-off of $X$ over $\mathcal{Q}$} is the graph obtained from $X$ by adding a vertex for each subgraph $Q \in \mathcal{Q}$ and linking by an edge $Q$ to each vertex belonging to $Q$.
\end{itemize}
\end{definition}

\noindent
Depending on the context, considering either the cone-off or the usual cone-off of a graph can be more or less convenient. Nevertheless, it is worth noticing that the inclusion from the set of vertices of the cone-off into the set of vertices of the usual cone-off defines a quasi-isometry. As a consequence, we know that if one of these two cone-off's is hyperbolic, then it follows that the second one must be hyperbolic as well. 

\medskip \noindent
Notice that a discrete geodesic metric space $(M,d)$ naturally defines a graph. Indeed, the graph whose vertices are the elements of $M$ and whose edges link two elements $x,y \in M$ satisfying $d(x,y)=1$ is naturally isometric to $(M,d)$. As a consequence, given a quasi-median graph $X$ endowed with a coherent system of discrete geodesic metrics, the metric space $(X,\delta)$ can be thought of as a graph. Our next proposition provides a criterion to determine when the cone-off of this graph turns out to be fine (we refer to Section \ref{section:extendingwall} for the definition of \emph{sector-walls}). This generalises \cite[Theorem 5.7]{coningoff}.

\begin{prop}\label{prop:criterionfine}
Let $X$ be a quasi-median graph endowed with a coherent system of discrete geodesic metrics, and $\mathcal{Q}$ a collection of gated subgraphs of $X$. Denote by $\delta$ the global distance extending our system of metrics. Suppose that
\begin{itemize}
	\item $(X,\delta)$ is a locally finite graph;
	\item $\mathcal{Q}$ is locally finite, ie., finitely many subspaces of $\mathcal{Q}$ contain a given edge of $(X,\delta)$;
	\item there exists a constant $C \geq 0$ such that, for any distinct $Q_1,Q_2 \in \mathcal{Q}$, at most $C$ sector-walls intersect both $Q_1$ and $Q_2$.
\end{itemize}
Then the usual cone-off $Y$ of the graph $(X,\delta)$ over $\mathcal{Q}$ is fine.
\end{prop}

\noindent
We recall from Section \ref{section:extendingwall} that a sector-wall is the data of a sector and its complement. A given sector-wall \emph{intersects} a subgraph $Q$ if it separates two vertices of $Q$. 

\begin{proof}
Let $e \in Y$ be an edge and fix one of its endpoints $a \in X$. To a given cycle $\gamma \subset Y$ of length $n$ containing $e$, we associate a loop $\overline{\gamma} \subset X$ containing $a$ in the following way. The cycle $\gamma$ passes through a sequence of cone-vertices $C_1, \ldots, C_k$ (which we identify with subspaces of $\mathcal{Q}$). For every $1 \leq j \leq k$, let $x_j,y_j \in X$ be the two vertices of $\gamma \cap C_j \cap X$, and choose a broken geodesic $[x_j,y_j]$ between $x_j$ and $y_j$. Also, for every $1 \leq j \leq k-1$, choose a broken geodesic $[y_j,x_{j+1}]$ between $y_j$ and $x_{j+1}$. Finally, choose two broken geodesics $[a,x_1]$ and $[y_n,a]$, respectively between $a$ and $x_1$, and $y_n$ and $a$. Now, we set
$$\overline{\gamma} = [a,x_1] \cup [x_1,y_1] \cup \cdots \cup [x_k,y_k] \cup [y_k,a].$$
Fix some $1 \leq j \leq k$. For every $i \neq j$, denote by $\mathfrak{J}_i$ the collection of the hyperplanes separating $x_j$ and $y_j$ and crossing the subgraph $C_i$; also, denote by $\mathfrak{H}$ the collection of the hyperplanes separating $x_j$ and $y_j$ which do not belong to any $\mathfrak{J}_i$. One has
$$\delta(x_j,y_j) \leq \sum\limits_{i \neq j} \sum\limits_{J \in \mathfrak{J}_i} \delta_J(x_j,y_j) + \sum\limits_{J \in \mathfrak{H}} \delta_J(x_j,y_j).$$
Fix some $i \neq j$. It is worth noticing that, since $\gamma$ is a cycle (ie., a simple loop), the cone-vertices $C_1, \ldots, C_k$ must be pairwise distinct, hence $C_i \neq C_j$. Moreover, because $C_i$ and $C_j$ are gated subgraphs of $X$, for every $J \in \mathfrak{J}_i$, they must both contain an entire clique dual to $J$, so that any sector-wall induced by $J$ intersects necessarily $C_i$ and $C_j$. Notice that $J$ induces $\# C(J)$ sector-walls, where $C(J)$ is an arbitrary clique dual to $J$, hence
$$\sum\limits_{J \in \mathfrak{J}_i} \# C(J) \leq C.$$
Consequently,
$$\sum\limits_{i \neq j} \sum\limits_{J \in \mathfrak{J}_i} \delta_J(x_j,x_j) \leq \sum\limits_{i \neq j} \sum\limits_{J \in \mathfrak{J}_i} \mathrm{diam}(C(J),\delta) \leq \sum\limits_{i \neq j} \sum\limits_{J \in \mathfrak{J}_i} \# C(J) \leq (k-1)C.$$
Now, fix some $J \in \mathfrak{H}$. If $\gamma_j$ (resp. $\overline{\gamma}_j$) denotes the path from $y_j$ to $x_j$ obtained from $\gamma$ (resp. $\overline{\gamma}$) by removing the subpath between $x_j$ and $y_j$, let $\mathcal{E}(J)$ be the collection of the edges of $\gamma_j \cap X$ which are dual to $J$. One has
$$\delta_J(x_j,y_j) \leq \sum\limits_{(a,b) \in \overline{\gamma}_j} \delta_J(a,b) \ \text{where} \ \delta_J(a,b)= \left\{ \begin{array}{cl} 1 & \text{if $(a,b)$ is dual to $J$} \\  0 & \text{otherwise} \end{array} \right. .$$
Since we know that $J$ is disjoint from $C_1, \ldots, C_k$, it follows that
$$\delta_J(x_j,y_j) \leq \sum\limits_{(a,b) \in \mathcal{E}(J)} \delta_J(a,b) = \# \mathcal{E}(J).$$
Finally, because $\mathcal{E}(J) \cap \mathcal{E}(J') = \emptyset$ for every two distinct hyperplanes $J,J' \in \mathfrak{H}$ and that $\gamma$ contains only $n$ edges, we deduce that
$$\sum\limits_{J \in \mathfrak{H}} \delta_J(x_j,y_j) \leq \sum\limits_{J \in \mathfrak{H}} \# \mathcal{E}(J) \leq n.$$
Our conclusion is that
$$\delta(x_j,y_j) \leq \sum\limits_{i \neq j} \sum\limits_{J \in \mathfrak{J}_i} \delta_J(x_j,y_j) + \sum\limits_{J \in \mathfrak{J}} \delta_J(x_j,y_j) \leq (k-1)C+n \leq (C+1)n.$$
Therefore, we get
$$\begin{array}{lcl} \mathrm{diam}_{\delta} ( \overline{\gamma}) & \leq & \displaystyle \delta(a,x_1)+ \sum\limits_{j=1}^k \delta(x_j,y_j) + \sum\limits_{j=1}^{k-1} \delta(y_j,x_{j+1})+ \delta(a,x_1) \\ \\ & \leq & k(C+1)n + n \leq (C+2)n^2 \end{array}$$
We have proved that $\overline{\gamma}$ is included into the ball $B= B_{\delta}(a, (C+2)n^2)$. Let $\dot{B}$ denote the cone-off of this ball over $\{Q \cap B \mid Q \in \mathcal{Q} \}$, so that $\overline{\gamma} \subset B$ implies $\gamma \subset \dot{B}$. It is worth noticing that $\dot{B}$ depends only on the edge $e$, the constant $C$ and the integer $n$, so $\dot{B}$ contains all the cycles of length $n$ passing through $e$. Since $\dot{B}$ is finite by the local finiteness of $(X, \delta)$ and $\mathcal{Q}$, we conclude that there exist only finitely many such cycles. This proves that $Y$ is finite. 
\end{proof}

\begin{proof}[Proof of Theorem \ref{thm:qmrelativelyhyp}.]
Let $\mathcal{C}$ be a collection of cliques such that any $G$-orbit of hyperplanes intersects it along a single clique, and let $\mathcal{C}= \mathcal{C}_1 \sqcup \mathcal{C}_2$ denote the associated decomposition of $\mathcal{C}$. Fix a base point $x_0 \in X$ and set $x_0(C)= \mathrm{proj}_C(x_0)$ for every clique $C$.

\medskip \noindent
If $C \in \mathcal{C}_2$, set $\delta_C : (x,y) \mapsto \left\{ \begin{array}{cl} 1 & \text{if $x \neq y$} \\ 0  & \text{otherwise} \end{array} \right.$. If $C \in \mathcal{C}_1$, use the bijection $\mathrm{stab}(C) \to C$ defined by $g \mapsto g \cdot x_0(C)$ to transfer a word metric of $\mathrm{stab}(C)$ (which is finitely generated) on $C$; let $\delta_C$ denote this metric. According to Proposition \ref{prop:mtopicalI}, $\{(C, \delta_C) \mid C \in \mathcal{C}\}$ extends as a $G$-invariant coherent system of metrics; observe that these metrics are discrete and geodesic according to Corollary \ref{cor:systgeodesic}. 

\medskip \noindent
Let $\mathcal{Q}$ denote the collection of subspaces of $X$ containing the finite cliques and the carriers of the hyperplanes dual to infinite cliques. Because $\mathcal{Q}$ is $G$-invariant, $G$ acts on the usual cone-off $Y$ of $X$ over $\mathcal{Q}$. Notice that 
\begin{itemize}
	\item $Y$ contains finitely many orbits of edges, since $X$ contains finitely many orbits of cliques and that the action $\mathrm{stab}(C) \curvearrowright C$ is transitive for every infinite clique $C$ of $X$;
	\item edge-stabilisers of $Y$ are trivial since vertex-stabilisers of $X$ are trivial themselves;
	\item vertex-stabilisers are either trivial or the stabiliser of a hyperplane of $X$ dual to an infinite clique.
\end{itemize}
In order to conclude that $G$ is hyperpolic relatively the stabilisers of hyperplanes of $X$ dual to infinite cliques, it is sufficient to show that $Y$ is a fine hyperbolic graph. Let us begin by proving that $Y$ is hyperbolic. In fact, because the cone-off and the usual cone-off of a space are quasi-isometric, it is sufficient to prove that the cone-off $Z$ of $(X, \delta)$ over $\mathcal{Q}$ is hyperbolic. By definition, $Z$ is obtained from $(X, \delta)$ by adding an edge between two vertices whenever they belong either to the same finite clique of $X$ or to a same hyperplane of $X$ which is dual to an infinite clique. A fortiori, $Z$ is also obtained from the quasi-median graph $X$ by adding an edge between two vertices whenever they belong to a same hyperplane which is dual to an infinite clique. We know by assumption that $X$ is hyperbolic, so we will use the following criterion, mentionned in \cite[Lemma 5.5]{arXiv:1304.1246}:

\begin{fact}
Let $Z$ be a graph obtained from another graph $X$ by adding edges. Suppose that $X$ is hyperbolic and that there exists some constant $M>0$ such that, for any two vertices $x,y \in X$ adjacent in $Z$ and any geodesic $\gamma$ in $X$ between $x$ and $y$, the diameter of $\gamma$ in $Z$ is at most $M$. Then $Z$ is hyperbolic.
\end{fact}

\noindent
In our context, if $x,y \in X$ are adjacent in $Z$, either they are adjacent in $X$, so that the edge between $x$ and $y$ is the unique geodesic in $X$ between $x$ and $y$, which clearly has diameter at most one in $Z$; or $x$ and $y$ belong to a hyperplane $J$ which is dual to an infinite clique. In the latter case, any geodesic between $x$ and $y$ must be included into $N(J)$, and a fortiori must have diameter at most one in $Z$. Consequently, the previous criterion applies, and we conclude that $Z$ is hyperbolic.

\medskip \noindent
Finally, we want to prove that $Y$ is a fine graph by applying Proposition \ref{prop:criterionfine}. First notice that $\mathcal{Q}$ is a collection of $\delta$-convex subspaces as a consequence of Corollary \ref{cor:systgeodesic}; that $(X, \delta)$ is locally finite according to Lemma \ref{lem:whenlocallyfinite}; and that $\mathcal{Q}$ is locally finite since an edge of $X$ belongs to a single clique of $X$ and only to finitely many hyperplanes of $X$ (otherwise, the endpoints of our edge would belong to infinitely many cliques of $X$). 

\medskip \noindent
Because $X$ contains finitely many orbits of cliques, there exists a constant $A_1$ such that any finite clique of $X$ has cardinality at most $A_1$. Also, because a hyperplane dual to an infinite clique must be transverse to finitely many hyperplanes and that $X$ contains finitely many orbits of hyperplanes, there exists a constant $A_2$ such that any hyperplane dual to an infinite clique is transverse to at most $A_2$ hyperplanes. Set $A= \max (A_1,A_2)$. Now, let $Q_1,Q_2 \in \mathcal{Q}$. We claim that at most $A^2$ sector-walls intersect both $Q_1$ and $Q_2$. We distinguish two cases. Firstly, if $Q_1$ is a finite clique, then the underlying hyperplane of any sector-wall intersecting both $Q_1$ and $Q_2$ must be the hyperplane dual to $Q_1$, which delimits at most $A$ sector-walls. Secondly, if $Q_1$ is the carrier of a hyperplane dual to an infinite clique, then the underlying hyperplane of any sector-wall intersecting both $Q_1$ and $Q_2$ must intersect $Q_1$; we already know that there exist at most $A$ such hyperplanes. Moreover, by assumption each of these hyperplanes must be dual to a finite clique, and so it delimits at most $A$ sector-walls. Therefore, at most $A^2$ sector-walls intersect both $Q_1$ and $Q_2$. This concludes the proof of our theorem.
\end{proof}

\noindent
Because a group which is hyperbolic relatively to hyperbolic groups must be hyperbolic itself, we deduce the following statement from Theorem \ref{thm:qmrelativelyhyp}. 

\begin{cor}\label{cor:beinghyperbolic}
Let $G$ be a group acting topically-transitively on a hyperbolic quasi-median graph $X$. Suppose that
\begin{itemize}
	\item every vertex belongs to finitely many cliques;
	\item vertex-stabilisers are finite;
	\item $X$ contains finitely many $G$-orbits of cliques;
	\item the hyperplanes dual to two infinite cliques cannot be transverse;
	\item carriers of hyperplanes dual to infinite cliques are cubically finite.
\end{itemize}
If $\mathrm{stab}(C)$ is hyperbolic for every clique $C$, then $G$ is hyperbolic as well.
\end{cor}

\begin{proof}
Theorem \ref{thm:qmrelativelyhyp} implies that $G$ is hyperbolic relatively to stabilisers of the hyperplanes of $X$ dual to infinite cliques. By assumption, the carrier of such a hyperplane must be cubically finite; in particuler, it contains only finitely many cliques, which implies that its stabilisers contains a clique-stabiliser as a finite-index subgroup. If clique-stabilisers are hyperbolic, it follows that $G$ is hyperbolic relatively to hyperbolic groups, so that $G$ must be hyperbolic (see for instance \cite[Corollary 2.41]{OsinRelativeHyp}).
\end{proof}

\subsection{Cubulating groups acting on quasi-median graphs I}\label{section:cubulatingI}

\noindent
In this section, we combine Proposition \ref{prop:wtopicalI} with some of the criteria proved in Section \ref{section:localglobalwall} in order to show that, if a group acts suitably on a quasi-median graph and if clique-stabilisers have good cubical properties, then these properties extend to the whole group.

\begin{prop}\label{prop:CAT0metricallyproper}
Let $G$ be a group acting topically-transitively on a quasi-median graph $X$. Suppose
\begin{itemize}
	\item any vertex of $X$ belongs to finitely many cliques;
	\item any vertex-stabiliser is finite.
\end{itemize}
If clique-stabilisers act metrically properly on CAT(0) cube complexes, then so does $G$.
\end{prop}

\begin{proof}
Let $\mathcal{C}$ be a collection of cliques such that any $G$-orbit of hyperplanes intersects it along a single clique, and let $\mathcal{C}= \mathcal{C}_1 \sqcup \mathcal{C}_2$ denote the associated decomposition of $\mathcal{C}$. Let $C \in \mathcal{C}$. If $C \in \mathcal{C}_2$, set $\mathcal{W}(C)= \{  \{ \{x \} \cup \{x \}^c \} \mid x \in C \}$. Otherwise, if $C \in \mathcal{C}_1$, fixing a base point $x_0(C) \in C$, $\psi_C : g \mapsto g \cdot x_0(G)$ defines a $\mathrm{stab}(G)$-equivariant bijection from $\mathrm{stab}(C)$ onto $C$ since the action $\mathrm{stab}(C) \curvearrowright C$ is transitive and free. By assumption, there exists a CAT(0) cube complex $X(C)$ on which $\mathrm{stab}(C)$ acts metrically properly; according to Lemma \ref{lem:modifycubing}, we may suppose without loss of generality that $X(C)$ contains a vertex $y_0(C)$ whose stabiliser is trivial. Set $\mathcal{W}(C) = \psi_C ~ \mathcal{M}( \mathrm{stab}(C) \curvearrowright X(C))$. 

\medskip \noindent
Proposition \ref{prop:wtopicalI} allows us to extend $\{ (C, \mathcal{W}(C)) \mid C \in \mathcal{C} \}$ to a $G$-invariant coherent system of wallspaces. Now, we want to apply Proposition \ref{prop:producingproperaction}.

\medskip \noindent
If $C \in \mathcal{C}$, then any two vertices of $C$ are separated by some wall of $\mathcal{W}(C)$: if $C \in \mathcal{C}_2$, this is clear from the definition of $\mathcal{W}(C)$, and if $C \in \mathcal{C}_1$, this is a consequence of Lemma \ref{lem:trivialstabiliser}. Next, we claim that $(C, \mathcal{W}(C))$ is locally finite. Once again, this is clear if $C \in \mathcal{C}_2$. If $C \in \mathcal{C}_1$, because the action $\mathrm{stab}(C) \curvearrowright C$ is transitive and free, we deduce that for every $x \in C$ and $R \geq 0$,
$$\# \{ y \in C \mid d_{\mathcal{W}(C)}(x,y) \leq R \}= \{ g \in \mathrm{stab}(C) \mid d_{\mathcal{W}(C)}(x,gx) \leq R \},$$
which is finite since we know from Lemma \ref{lem:transferproper} that the action $\mathrm{stab}(C) \curvearrowright (C, \mathcal{W}(C))$ is metrically proper. A fortiori, $(C, \mathcal{W}(C))$ is locally finite.

\medskip \noindent
On the other hand, according to Proposition \ref{prop:wtopicalI}, if $C'$ is an arbitrary clique of $X$, then there exists some $C \in \mathcal{C}$ such that the spaces with walls $(C,\mathcal{W}(C))$ and $(C',\mathcal{W} (C'))$ are isomorphic. As a consequence, $(C', \mathcal{W}(C'))$ is locally finite, and two vertices of $C'$ are separated by some wall of $\mathcal{W}(C')$. Therefore, by applying Proposition \ref{prop:producingproperaction}, we conclude that $G$ acts metrically properly on the CAT(0) cube complex obtained by cubulating $(X, \mathcal{HW})$.
\end{proof}

\begin{prop}\label{prop:cubulatinggeometrically}
Let $G$ be a group acting topically-transitively on a quasi-median graph $X$. Suppose that
\begin{itemize}
	\item any vertex of $X$ belongs to finitely many cliques;
	\item any vertex-stabiliser is finite;
	\item the cubical dimension of $X$ is finite;
	\item $X$ contains finitely many orbits of prisms;
	\item for every maximal prism $P=C_1 \times \cdots \times C_n$, $\mathrm{stab}(P)= \mathrm{stab}(C_1) \times \cdots \times \mathrm{stab}(C_n)$.
\end{itemize}
If, for every clique $C$, $\mathrm{stab}(C)$ acts geometrically on a CAT(0) cube complex of dimension $d(C)$, then $G$ acts geometrically on a CAT(0) cube complex of dimension at most $\mathrm{dim}_{\square}(X) \cdot \max \{ d(C) \mid C \ \text{clique}\}$.
\end{prop}

\begin{proof}
Let $\mathcal{C}$ be a collection of cliques such that any $G$-orbit of hyperplanes intersects it along a single clique, and let $\mathcal{C}= \mathcal{C}_1 \sqcup \mathcal{C}_2$ denote the associated decomposition of $\mathcal{C}$. Let $C \in \mathcal{C}$. If $C \in \mathcal{C}_2$, set $\mathcal{W}(C)= \{  \{ \{x \} \cup \{x \}^c \} \mid x \in C \}$. Otherwise, if $C \in \mathcal{C}_1$, fixing a base point $x_0(C) \in C$, $\psi_C : g \mapsto g \cdot x_0(G)$ defines a $\mathrm{stab}(G)$-equivariant bijection from $\mathrm{stab}(C)$ onto $C$ since the action $\mathrm{stab}(C) \curvearrowright C$ is transitive and free. By assumption, there exists a CAT(0) cube complex $X(C)$ on which $\mathrm{stab}(C)$ acts geometrically; according to Lemma \ref{lem:modifycubing}, we may suppose without loss of generality that $X(C)$ contains a vertex $y_0(C)$ whose stabiliser is trivial. Set $\mathcal{W}(C) = \psi_C ~ \mathcal{N}( \mathrm{stab}(C) \curvearrowright X(C))$. 

\medskip \noindent
Proposition \ref{prop:wtopicalI} allows us to extend $\{ (C, \mathcal{W}(C)) \mid C \in \mathcal{C} \}$ to a $G$-invariant coherent system of wallspaces. Now, we want to apply Proposition \ref{prop:producingproperaction} and Proposition \ref{prop:producingcocompactaction} in order to deduce that the action $G \curvearrowright (X, \mathcal{HW})$ is metrically proper and cocompact. We already know from the proof of the previous proposition that the hypotheses of Proposition \ref{prop:producingproperaction} hold. Similarly, since it is clear that $\mathcal{W}(C) \neq \emptyset$ and $\dim(C, \mathcal{W}(C)) < + \infty$ for every $C \in \mathcal{C}$, we deduce that the same assertion holds for every clique.

\medskip \noindent
Finally, if $P= C_1 \times \cdots \times C_n$ is a prism, since $\mathcal{W}(P)= \mathcal{W}(C_1) \times \cdots \times \mathcal{W}(C_n)$ and $\mathrm{stab}(P)= \mathrm{stab}(C_1) \times \cdots \times \mathrm{stab}(C_n)$, we deduce that the action $\mathrm{stab}(P) \curvearrowright (P, \mathcal{W}(P))$ is cocompact since the action $\mathrm{stab}(C) \curvearrowright (C,\mathcal{W}(C))$ is cocompact for every clique $C$. Indeed, we know that this statement is true for $C \in \mathcal{C}_2$, and for $C\in \mathcal{C}_1$ according to Lemma \ref{lem:transfercocompact}, and otherwise we know that the action $\mathrm{stab}(C) \curvearrowright (C,\mathcal{W}(C))$ is isomorphic to the action $\mathrm{stab}(C') \curvearrowright (C',\mathcal{W}(C'))$ for some $C' \in \mathcal{C}$, according to Proposition \ref{prop:wtopicalI}, so that its cocompactness follows as well.

\medskip \noindent
We conclude that $G$ acts geometrically on the CAT(0) cube complex obtained by cubulating $(X, \mathcal{HW})$. 

\medskip \noindent
The assertion on the dimension of the cube complex we obtain follows from Corollary \ref{cor:dimHW}, combined with the observation that $\dim(C', \mathcal{W}(C'))= \dim(C, \mathcal{W}(C))$ for some $C \in \mathcal{C}$, since $(C', \mathcal{W}(C'))$ is isomorphic to $(C, \mathcal{W}(C))$ for some $C \in \mathcal{C}$ according to Proposition \ref{prop:wtopicalI}, and $\dim(C, \mathcal{W}(C))=1$ if $C \in \mathcal{C}_2$ and $\dim(C, \mathcal{W}(C)) \leq d(C)$ if $C \in \mathcal{C}_1$ according to Lemma \ref{lem:transfergeometricaction}.
\end{proof}

\subsection[A-T-menability and a-B-menability]{A-T-menability and a-$\mathcal{B}$-menability}

\noindent
In this section, we mention how to extend Proposition \ref{prop:CAT0metricallyproper} to construct metrically proper actions on spaces with measured walls and spaces with labelled partitions. We begin with spaces with measured walls. The first step is to extend Proposition \ref{prop:mwtopicalI}. By following closely its proof, we show that:

\begin{prop}\label{prop:mwtopicalI}
Let $G$ be a group with a topical-transitive action on a quasi-median graph $X$ and $\mathcal{C}$ a collection of cliques such that any $G$-orbit of hyperplanes intersects it along a single clique. For every $C \in \mathcal{C}$, let $\mathcal{W}(C)$ be a $\mathrm{stab}(C)$-invariant collection of walls, a $\sigma$-algebra $\mathcal{B}(C)$ on $\mathcal{W}(C)$, and a measure $\mu_C$, such that $(C, \mathcal{W}(C), \mathcal{B}(C),\mu_C)$ is a space with measured walls. Then there exists a $G$-invariant coherent system of measured wallspaces extending $\{ (C, \mathcal{W}(C), \mathcal{B}(C), \mu_C) \mid C \in \mathcal{C} \}$ such that, for every clique $C'$, labelled by some $C \in \mathcal{C}$, the actions $\mathrm{stab}(C') \curvearrowright (C', \mathcal{W}(C'), \mathcal{B}(C'), \mu_{C'})$ and $\mathrm{stab}(C) \curvearrowright (C, \mathcal{W}(C), \mathcal{B}(C), \mu_C)$ are isomorphic.
\end{prop}

\noindent
Afterwards, Proposition \ref{prop:CAT0metricallyproper} naturally extends as

\begin{prop}\label{aTmenablegroups}
Let $G$ be a group acting topically-transitively on a quasi-median graph $X$. Suppose
\begin{itemize}
	\item any vertex of $X$ belongs to finitely many cliques;
	\item any vertex-stabiliser is finite.
\end{itemize}
If clique-stabilisers are a-T-menable, then so is $G$.
\end{prop}

\begin{proof}[Sketch of proof.]
Let $\mathcal{C}$ be a collection of cliques such that any $G$-orbit of hyperplanes intersects it along a single clique, and let $\mathcal{C}= \mathcal{C}_1 \sqcup \mathcal{C}_2$ denote the associated decomposition of $\mathcal{C}$. Let $C \in \mathcal{C}$. If $C \in \mathcal{C}_2$, set $\mathcal{W}(C)= \{  \{ \{x \} \cup \{x \}^c \} \mid x \in C \}$ so that $(C, \mathcal{W}(C))$ is a discrete space with measured walls. Otherwise, if $C\in \mathcal{C}_1$, we know from \cite{medianviewpoint} that $\mathrm{stab}(C)$ acts metrically properly on a space with measured walls $X(C)$; by reproducing the construction of Lemma \ref{lem:modifycubing}, we may suppose without loss of generality that $X(C)$ contains a point $x_0(C)$ whose stabiliser is trivial. Next, we pullback the walls from $X(C)$ to $\mathrm{stab}(C)$ thanks to the map $g \mapsto g \cdot x_0(C)$, making $\mathrm{stab}(C)$ a space with measured walls (see also \cite[Lemma 3.9]{medianviewpoint}). Finally, use Proposition \ref{prop:mwtopicalI} in order to extend these local spaces with measured walls to obtain a space with measured walls $(X, \mathcal{W}, \mathcal{B}, \mu)$, and conclude by using Proposition \ref{prop:producingproperactionmhw} that $G$ acts metrically properly on $(X, \mathcal{W}, \mathcal{B}, \mu)$. Therefore, $G$ must be a-T-menable.
\end{proof}

\noindent
In exactly the same way, we show the following combination result by using spaces with labelled partitions (see Section \ref{section:generalizedwalls}).

\begin{prop}\label{aBmenablegroups}
Let $G$ be a group acting topically-transitively on a quasi-median graph $X$ and $p \geq 1$ such that $p \notin 2 \mathbb{Z}$. Suppose that
\begin{itemize}
	\item any vertex of $X$ belongs to finitely many cliques;
	\item any vertex-stabiliser is finite.
\end{itemize}
If clique-stabilisers are a-$L^p$-menable, then so is $G$.
\end{prop}

\subsection{Creating coherent invariant systems II}

\noindent
In this section, our goal is to show that the maps given by Proposition \ref{prop:topicalI0} can be constructed in a natural way. The control we get on these maps will be fundamental in Sections \ref{section:equicompressiongeneral} and \ref{section:topicalactionsII}.

\medskip \noindent 
From now on, fix a group $G$ acting topically-transitively on a quasi-median graph $X$ and a collection of cliques $\mathcal{C}$ such that every $G$-orbit of hyperplanes intersects it along a single clique. Let $\mathcal{C}=\mathcal{C}_1 \sqcup \mathcal{C}_2$ denote the corresponding decomposition of $\mathcal{C}$. Given a clique $C$ and an element $g \in G$, we decompose canonically $g$ as a product $g=p_C(g) \cdot s_C(g)$ where $s_C(g) \in \mathrm{stab}(C)$ and $p_C(g) \cdot C= gC$. The idea we have in mind is roughly the following: $g$ sends $C$ onto $gC$ first by ``rotating'' $C$ via the action of $s_C(g) \in \mathrm{stab}(C)$, and then by translating $C$ to $gC$ without ``rotating'' the clique. 

\begin{definition}\label{def:topicalps}
Let $G$ be a group acting $\mathcal{C}$-transitively on a quasi-median graph $X$. Fix a base point $x_0 \in X$, and let us set $x_0(C)= \mathrm{proj}_C(x_0)$ for every clique $C$. Consider an element $g \in G$ and a clique $C$. If $C$ is labelled by $\mathcal{C}_1$, let $p_C(g)$ denote the unique element of $g \cdot \mathrm{stab}(C)$ satisfying $p_C(g) \cdot x_0(C)=x_0(gC)$; set $s_C(g)= p_C(g)^{-1}g$. If $C$ is labelled by $\mathcal{C}_2$, set $p_C(g)=g$ and $s_C(g)=1$. It is worth noticing that $p_C(p_C(g))=p_C(g)$. 
\end{definition}

\noindent
Let us study how the operations $p_C(\cdot)$ and $s_C(\cdot)$ behave. 

\begin{claim}\label{claim:pC}
For every clique $C$ and every $g_1,g_2 \in G$, $p_C(g_1g_2)=p_{g_2C}(g_1) p_C(g_2)$.
\end{claim}

\begin{proof}
If $C$ is labelled by a clique of $\mathcal{C}_2$, then
$$p_C(g_1g_2)=g_1g_2= p_{g_2C}(g_1) p_C(g_2).$$
Otherwise, if $C$ is labelled by a clique of $\mathcal{C}_1$, notice that
$$p_{g_2C}(g_1) p_C(g_2) \cdot x_0(C)= p_{g_2C}(g_1) \cdot x_0(g_2C) = x_0(g_1g_2C),$$
so $p_C(g_1g_2) = p_{g_2C}(g_1) p_C(g_2)$.
\end{proof}

\begin{cor}\label{cor:sCgh}
For every clique $C$ and every $g_1,g_2 \in G$, $$s_C(g_1g_2)=p_C(g_2)^{-1} s_{g_2C}(g_1)p_C(g_2) \cdot s_C(g_2).$$
\end{cor}

\begin{proof}
As a consequence of Claim \ref{claim:pC}, 
$$\begin{array}{lcl} s_C(g_1g_2) & = & p_C(g_1g_2)^{-1}g_1g_2 = p_C(g_2)^{-1}p_{g_2C}(g_1)^{-1}g_1g_2 \\ \\ & = & p_C(g_2)^{-1} s_{g_2C}(g_1) g_2 = p_C(g_2)^{-1} s_{g_2C}(g_1)p_C(g_2) s_C(g_2) \end{array}$$
This concludes the proof.
\end{proof}

\begin{claim}\label{claim:wellphi}
Let $C,Q$ be two cliques, with $Q \in \mathcal{C}$, and $g,h \in G$ two elements satisfying $p_C(g)=g$, $p_C(h)=h$ such that the cliques $gC$, $Q$ and $hC$ are dual to the same hyperplane. Then $t_{gC \to Q}(gx)= t_{hC \to Q}(hx)$ for every $x \in C$.
\end{claim}

\begin{proof}
Let $J$ denote the hyperplane dual to $Q,gC,hC$. Because $hg^{-1} \cdot gC=hC$, necessarily $hg^{-1} \in \mathrm{stab}(J)$. We want to prove that $hg^{-1} \in \ker(\rho_Q)$. If $Q \in \mathcal{C}_2$, this is clear since $\mathrm{stab}(Q)= \mathrm{fix}(Q)$. From now on, suppose that $Q \in \mathcal{C}_1$. Notice that, because $p_C(g)=g$ and $p_C(h)=h$,
$$hg^{-1} \cdot x_0(gC) = h \cdot x_0(C)= x_0(hC).$$
Therefore, $hg^{-1}$ stabilises the sector delimited by $J$ which contains $x_0(gC)$ and $x_0(hC)$ (which is also the sector containing $x_0$). A fortiori, as a permutation of $Q$, $\rho_Q(hg^{-1})$ must fix $x_0(Q)$. On the other hand, since the action is topical-transitive, $\rho_Q(hg^{-1})$ induces the same permutation on $Q$ as some elements of $\mathrm{stab}(Q)$, so that, because $\mathrm{stab}(Q)$ acts freely on $Q$, it follows that $\rho_Q(hg^{-1}) \in \mathrm{fix}(Q) = \{1\}$.

\medskip \noindent
Thus, we have proved that $hg^{-1}$ stabilises every sector delimited by $J$, so that, for every vertex $x \in C$, $gx$ and $hx=hg^{-1} \cdot gx$ must belong to the same sector delimited by $J$; in particular, $t_{gC \to Q}(gx)= t_{hC \to Q}(hx)$.
\end{proof}

\begin{claim}\label{claim:wellvarphi}
Let $C,Q$ be two cliques, with $Q \in \mathcal{C}$, and $h,k \in G$ two elements satisfying $p_C(h)=h$, $p_C(k)=k$ such that the cliques $hC$, $Q$ and $kC$ are dual to the same hyperplane. Then $\rho_Q(hgh^{-1})=\rho_Q(kgk^{-1})$ for every $g \in \mathrm{stab}(C)$. 
\end{claim}

\begin{proof}
If $Q \in \mathcal{C}_2$, there is nothing to prove since $\mathrm{Im}(\rho_Q)$ is reduced to the identity. From now on, suppose that $Q \in \mathcal{C}_1$. As in the previous proof, we notice that $hk^{-1} \in \mathrm{stab}(J)$, where $J$ is the hyperplane dual to $Q,hC,kC$, and we show that $\rho_Q(hk^{-1})= 1$, ie., $hk^{-1}$ stabilises every sector delimited by $J$. In particular, $kg \cdot x_0(C)$ and $hk \cdot x_0(C)= hk^{-1} \cdot kg \cdot x_0(C)$ belong to the same sector delimited by $J$. On the other hand, if $S$ denotes the sector delimited by $J$ which contains $x_0$, then $S$ also contains $x_0(hC)$ and $x_0(kC)$; moreover, since $p_C(h)=h$ and $p_C(k)=k$, we know that $kgk^{-1} \cdot x_0(kC)=kg \cdot x_0(C)$ and $hkh^{-1} \cdot x_0(hC)= hk \cdot x_0(C)$. As a consequence, we deduce that $kgk^{-1} \cdot S= hgh^{-1} \cdot S$, which also implies that
$$\rho_Q(kgk^{-1}) \cdot x_0(Q)= \rho_Q(hkh^{-1}) \cdot x_0(Q).$$
Since $\mathrm{stab}(Q)$ acts freely on $Q$ and that $\rho_Q(kgk^{-1})= \rho_{kC \to Q}(kgk^{-1})$ and $\rho_Q(hkh^{-1})= \rho_{hC \to Q}(hkh^{-1})$ induce the same permutations as elements of $\mathrm{stab}(Q)$, since the action is topical-transitive, we conclude that $\rho_Q(hgh^{-1})=\rho_Q(kgk^{-1})$.
\end{proof}

\noindent
Finally, we are ready to introduce the maps which will replace the maps given by Proposition \ref{prop:topicalI0}. 

\begin{definition}\label{def:mapsfortopical}
Let $G$ be a group acting topically-transitively on a quasi-median graph $X$ and $\mathcal{C}$ a collection of cliques such that every $G$-orbit of hyperplanes intersects it along a single clique. Let $C$ be a clique labelled by some $Q \in \mathcal{C}$, and fix some $g \in G$ satisfying $p_C(g)=g$ such that $gC$ and $Q$ are dual to the same hyperplane. Now, we define
$$\phi_C : \left\{ \begin{array}{ccc} C & \to & Q \\ x  & \mapsto & t_{gC \to Q}(gx) \end{array} \right. \ \text{and} \ \varphi_C : \left\{ \begin{array}{ccc} \mathrm{stab}(C) & \to & \mathrm{Im}(\rho_Q) \\ h & \mapsto & \rho_Q(ghg^{-1}) \end{array} \right. .$$
\end{definition}

\noindent
It is worth noticing that our maps do not depend on the choice of $g$ according to Claim \ref{claim:wellphi} and Claim \ref{claim:wellvarphi}. Moreover, $\phi_C$ is a bijection, and it follows from Lemma \ref{lem:rho3} that $\varphi_C$ induces an isomorphism $\mathrm{stab}(C) \to \mathrm{stab}(Q)$ if $Q \in \mathcal{C}_1$. 

\begin{claim}\label{claim:divers}
The following assertions hold:
\begin{itemize}
	\item[(i)] if $C,C'$ are two cliques dual to the same hyperplane, $\phi_{C'} \circ t_{C \to C'}= \phi_C$;
	\item[(ii)] for every clique $C$, $x \in C$ and $g \in G$, $\phi_{gC}(gx)= \varphi_C(s_C(g)) \cdot \phi_C(x)$;
	\item[(iii)] for every clique $C$, $g \in G$ and $h \in \mathrm{stab}(gC)$, $\varphi_C(p_C(g)^{-1}h p_C(g)) = \varphi_{gC}(h)$.
\end{itemize}
\end{claim}

\begin{proof}
We first prove point $(i)$. So let $C,C'$ be two cliques dual to the same hyperplane, say labelled by $Q \in \mathcal{C}$. Fix two elements $h(C),h(C') \in G$ satisfying $p_C(h(C))=h(C)$ and $p_C(h(C'))=h(C')$ such that $Q$, $h(C) \cdot C$ and $h(C') \cdot C'$ are dual to the same hyperplane. For every $x \in C$, we have
$$\begin{array}{lcl} \phi_{C'}(t_{C \to C'}(x)) & = & t_{h(C')C' \to Q} \left( h(C') \cdot t_{C \to C'}(x) \right) \\ \\ & = & t_{h(C')C' \to Q } \circ t_{h(C')C \to h(C')C'} (h(C') \cdot x) \\ \\ & = & t_{h(C')C \to Q} (h(C') \cdot x) \end{array}$$
Finally, we deduce from Claim \ref{claim:wellphi} that
$$\phi_{C'}(t_{C \to C'}(x))= t_{h(C)C \to Q}(h(C) \cdot x) = \phi_C(x).$$
Now, focus on point $(ii)$. So let $C$ be a clique, $x \in C$ a vertex and $g \in G$. Assume that $Q \in \mathcal{C}$ labels $C$ and $gC$. Fix two elements $h(C),h(gC) \in G$ satisyfing $p_C(h(C))=h(C)$ and $p_{gC}(h(gC))=h(gC)$ such that $Q$, $h(C) \cdot C$ and $h(gC) \cdot gC$ are dual to the same hyperplane. As an application of Claim \ref{claim:wellphi},
$$\phi_C(x)= t_{h(C)C \to Q}(h(C) \cdot x) = t_{h(gC)gC \to Q} \left( p_C(h(gC)g) \cdot x \right).$$
Therefore, by applying Lemma \ref{lem:rho2}, we deduce that
$$\begin{array}{lcl} \varphi_C(s_C(g)) \cdot \phi_C(x) & = & \rho_Q( h(C)s_C(g)h(C)^{-1}) \cdot t_{h(gC)gC \to Q} ( p_C(h(gC)g) \cdot x) \\ \\ & = & t_{h(gC)gC \to Q} \left( \rho_{h(gC)gC} \left( h(C) s_C(g) h(C)^{-1} \right) p_C(h(gC)g) \cdot x \right) \end{array}$$
Thanks to Claim \ref{claim:wellvarphi}, notice that
$$\begin{array}{lcl} \rho_{h(gC)gC} \left( h(C) s_C(g) h(C)^{-1} \right) & = & \rho_{h(gC)gC} \left( p_C(h(gC)g) s_C(g) p_C(h(gC)g)^{-1} \right) \\ \\ & = & p_C(h(gC)g) s_C(g) p_C(h(gC)g)^{-1} \end{array}$$
since $p_C(h(gC)g) s_C(g) p_C(h(gC)g)^{-1} \in \mathrm{stab}( h(gC)gC)$. Consequently,
$$\varphi_C(s_C(g)) \cdot \phi_C(x) = t_{h(gC)gC \to Q} \left( p_C(h(gC)g) s_C(g) \cdot x \right).$$
On the other hand, by applying Claim \ref{claim:pC}, we know that
$$p_C(h(gC)g)s_C(g)= p_{gC}(h(gC)) p_C(g)s_C(g)=h(gC)g,$$
hence
$$\varphi_C(s_C(g)) \cdot \phi_C(x) = t_{h(gC)gC \to Q} \left( h(gC)g \cdot x \right) = \phi_{gC}(gx).$$
Finally, let us prove point $(iii)$. Let $C$ be a clique, $g \in G$ an element satisfying $p_C(g)=g$, and $h \in \mathrm{stab}(gC)$. Assume that $Q \in \mathcal{C}$ labels $C$ and $gC$. Fix some $h(C) \in G$ satisfying $p_C(h(C))=h(C)$ such that $h(C) \cdot C$ and $Q$ are dual to the same hyperplane. If $Q \in \mathcal{C}_2$, there is nothing to prove since $\mathrm{Im}(\rho_Q)$ is trivial, so suppose that $Q \in \mathcal{C}_1$. Notice that $h(C)g^{-1} \cdot gC = h(C) \cdot C$ and $Q$ are dual to the same hyperplane and that, using Claim \ref{claim:pC},
$$p_{gC}(h(C)g^{-1})= p_C(h(C))p_{gC}(g^{-1}) = h(C)g^{-1},$$
where the equality $p_{gC}(g^{-1})=g^{-1}$ is justified by the equality $g^{-1}x_0(gC)=x_0(C)$, since $p_C(g)=g$ implies that $g \cdot x_0(C) = x_0(gC)$. Therefore, 
$$\varphi_C(g^{-1}hg) = \rho_Q(h(C)g^{-1}hgh(C)^{-1}) = \varphi_{gC}(h).$$
This concludes the proof of our claim.
\end{proof}

\noindent
We conclude this section by reproving Propositions \ref{prop:wtopicalI} and \ref{prop:mtopicalI} by using our new maps. 

\begin{proof}[Proof of Proposition \ref{prop:wtopicalI}.]
For every clique $C$ labelled by $Q\in \mathcal{C}$, set $\mathcal{W}(C)= \phi_C^{-1} \mathcal{W}(Q)$. Notice that, by applying Claim \ref{claim:divers}, we deduce that for every $g \in G$
$$g \mathcal{W}(C)= g \phi_C^{-1} \mathcal{W}(Q) = \phi_{gC}^{-1} \left( \varphi_C(s_C(g)) \cdot \mathcal{W}(Q) \right) = \mathcal{W}(gC),$$
since $\varphi_C(s_C(g)) \in \mathrm{stab}(Q)$ implies $\varphi_C(s_C(g)) \cdot \mathcal{W}(Q) = \mathcal{W}(Q)$. Moreover, if $C'$ is a second clique dual to the same hyperplane as $C$, then Claim \ref{claim:divers} implies that
$$\mathcal{W}(C')= \phi_{C'}^{-1} \mathcal{W}(Q)= t_{C'\to C} \circ \phi_C^{-1} \mathcal{W}(Q) = t_{C' \to C} \mathcal{W}(C)= \mathcal{W}(C \to C').$$
Thus, our system is both $G$-invariant and coherent. Finally, it follows from Claim \ref{claim:divers} that $\phi_C$ defines a $\varphi_C$-equivariant isomorphism from $(C, \mathcal{W}(C))$ onto $(Q, \mathcal{W}(Q))$ if $Q \in \mathcal{C}_1$. 
\end{proof}

\begin{proof}[Proof of Proposition \ref{prop:mtopicalI}.]
For every clique $C$ labelled by $Q\in \mathcal{C}$, let us define
$$\delta_C : (x,y) \mapsto \delta_Q \left( \phi_C(x), \phi_C(y) \right).$$
Fix some $g\in G$. Because $\delta_Q$ is $\mathrm{stab}(Q)$-invariant, and thanks to Claim \ref{claim:divers}, we deduce that, for every $x,y \in C$,
$$\begin{array}{lcl} \delta_{gC}(gx,gy) & = & \delta_Q \left( \phi_{gC}(gx), \phi_{gC}(gy) \right) \\ \\ & = & \delta_Q \left( \varphi_C(s_C(g)) \cdot \phi_C(x) , \varphi_C(s_C(g)) \cdot \phi_C(y) \right) \\ \\ & = & \delta_Q \left( \phi_C(x), \phi_C(y) \right) = \delta_C(x,y) \end{array}$$
Next, if $C'$ is a second clique dual to the same hyperplane as $C$, then Claim \ref{claim:divers} implies that
$$\begin{array}{lcl} \delta_{C'}(x,y) & = & \delta_Q\left( \phi_{C'}(x), \phi_{C'}(y) \right) \\ \\ & = & \delta_Q \left( \phi_C \circ t_{C' \to C}(x), \phi_C \circ t_{C' \to C} (y) \right) \\ \\ & = & \delta_{C} \left( t_{C' \to C}(x), t_{C' \to C}(y) \right) = \delta_{C \to C'}(x,y) \end{array} $$
for every $x,y \in C'$. Thus, our system is both $G$-invariant and coherent. Finally, it follows from Claim \ref{claim:divers} that 

\begin{fact}\label{fact:phiCequiisom}
The map $\phi_C$ defines a $\varphi_C$-equivariant isometry from $(C, \delta_C)$ onto $(Q, \delta_Q)$ if $Q \in \mathcal{C}_1$. 
\end{fact}

\noindent
This concludes the proof.
\end{proof}

\subsection{Equivariant $\ell^p$-compression}\label{section:equicompressiongeneral}

\noindent
In Section \ref{section:compression}, we studied the $\ell^p$-compression of the global metric on a quasi-median graph associated to a coherent system of metrics. Now, we suppose that some group acts on our quasi-median graph, and we are interested in the \emph{equivariant $\ell^p$-compression}, that is, in coarsely embedding our metric space into some $L^p$-space in an equivariant way. More precisely:

\begin{definition}
Let $G$ be a group acting on a metric space $X$ by isometries, and $p \geq 1$. The \emph{equivariant $\ell^p$-compression}\index{Equivariant compression ($\ell^p$-)} of $X$, denoted by $\alpha_p^*(X)$, is the supremum of the $\alpha(f)$'s where $G$ acts on some $L^p$-space $L$ by affine isometries and $f : X \to L$ is a Lipschitz map which is $G$-equivariant. (Recall that $G$ acts on $L$ by affine isometries if there exist a linear representation $g \mapsto \pi_g$ and a map $b : G \to L$, satisfying $b(gh)=\pi_g(b(h))+b(g)$ for every $g,h \in G$, such that $g \cdot x = \pi_g(x)+b(g)$ for every $x \in X$ and $g \in G$.) In particular, the \emph{equivariant $\ell^p$-compressio}n of a finitely generated group is the equivariant $\ell^p$-compression of the metric space obtained from this group when endowed with the word length associated to some finite generating set and with the left-multiplication. 
\end{definition}

\noindent
The following preliminary lemma may be useful as it allows us to suppose that a given equivariant coarse embedding vanishes at some basepoint.

\begin{lemma}\label{lem:translation}
Let $G$ be a group acting on some metric space $X$ and on some $L^p$-space $L$ by affine isometries, say $g \mapsto \pi_g(\cdot)+b(g)$, and let $f : X \to L$ be a $G$-equivariant map. For every $x_0 \in L$, the map $g \mapsto \pi_g( \cdot) + b(g)+x_0-\pi_g(x_0)$ defines an action $G \curvearrowright L$ by affine isometries such that $f'=f( \cdot)+x_0$ is a $G$-equivariant map $X \to L$. 
\end{lemma}

\begin{proof}
Set $b' : g \mapsto b(g)+x_0- \pi_g(x_0)$. For every $g,h \in G$, 
$$\begin{array}{lcl} b'(gh) & = & \pi_g(b(h))+b(g) +x_0- \pi_{gh}(x_0) = \pi_g(b'(h))+b(g)+x_0 - \pi_g(x_0) \\ \\ & = & \pi_g(b'(h))+b'(g) \end{array}$$
Thus, $g \mapsto \pi_g( \cdot)+b'(g)$ defines an action $G \curvearrowright L$ by affine isometries. Next, for every $x \in L$ and $g \in G$,
$$\begin{array}{lcl} f'(gx) & = & f(gx)+x_0 = \pi_g(f(x))+b(g)+x_0 = \pi_g(f'(x))+b(g)+x_0-\pi_g(x_0) \\ \\ & = & \pi_g(f'(x))+b'(g) \end{array}$$
so $f'$ is indeed $G$-equivariant with respect to our new action. 
\end{proof}

\noindent
The main result of this section is the following proposition.

\begin{prop}\label{prop:equicompression}
Let $G$ be a group acting topically-transitively on a quasi-median graph $X$. Let $\mathcal{C}$ be a collection of cliques such that any $G$-orbit of hyperplanes intersects it along a single clique, and let $\mathcal{C}= \mathcal{C}_1 \sqcup \mathcal{C}_2$ denote the associated decomposition of $\mathcal{C}$. Suppose that 
\begin{itemize}
	\item $\mathrm{stab}(C)$ is finitely generated for every $C \in \mathcal{C}_1$ and finite for every $C \in \mathcal{C}_2$;
	\item $\mathcal{C}_1$ is finite and the cardinalities of the cliques of $\mathcal{C}_2$ are uniformly bounded.
\end{itemize}
There exists a $G$-invariant coherent system of metrics on $X$ such that 
$$\alpha_p^*(X,\delta) \geq \min \left( \frac{1}{p}, \min\limits_{C \in \mathcal{C}_1} \alpha_p^*( \mathrm{stab}(C)) \right).$$
As a consequence, if some orbit map $G \to (X,\delta)$ is Lipschitz and has compression $\alpha$,
$$\alpha_p^*(G) \geq \alpha \cdot \min \left( \frac{1}{p}, \min\limits_{C \in \mathcal{C}_1} \alpha_p^*( \mathrm{stab}(C)) \right).$$
\end{prop}

\begin{proof}
If $\min\limits_{C \in \mathcal{C}_1} \alpha_p^*( \mathrm{stab}(C)) =0$, there is nothing to prove, so we suppose that this quantity is positive. Fix a basepoint $x_0 \in X$, and, for every clique $C$, let $x_0(C)$ denote the projection of $x_0$ onto $C$. Choose a finite generating set for every $\mathrm{stab}(C)$ where $C \in \mathcal{C}_1$, and define a metric $\delta_C$ on $C$ by transfering the associated word length from $\mathrm{stab}(C)$ to $C$ thanks to the orbit map $g \mapsto g \cdot x_0(C)$; if $C \in \mathcal{C}_2$, set $\delta_C : (x,y) \mapsto \left\{ \begin{array}{cl} 1 & \text{if $x \neq y$} \\ 0  & \text{otherwise} \end{array} \right.$. Finally, apply Proposition \ref{prop:mtopicalI} to extend this collection of metrics to a $G$-invariant coherent system of metrics.

\medskip \noindent
Fix some $\epsilon >0$. For every $C \in \mathcal{C}$, let $f_C : (C,\delta_C) \to L_C$ be a $\mathrm{stab}(C)$-equivariant Lipschitz map to some $L^p$-space satisfying $\alpha(f_C) \geq \alpha_p^*(C,\delta_C)- \epsilon$. Moreover, thanks to our assumptions on $\mathcal{C}$, we can choose the maps $f_C$ so that, for every $0< \eta < \min\limits_{C\in \mathcal{C}} \alpha_p^*(C,\delta_C)$, there exist positive constants $A_{\eta}$ and $B$ (which does not depend on $C$) such that the inequalities
$$A_{\eta} \cdot \delta_C(x,y)^{\alpha_p(f_C)- \eta} \leq \| f_C(x)- f_C(y) \| \leq B \cdot \delta_C(x,y)$$
hold for every $C \in \mathcal{C}$ and every $x,y \in C$. 
Without loss of generality, suppose that an action of $\mathrm{stab}(C)$ on $L_C$ is trivial if $C\in \mathcal{C}_2$. Notice that, according to Lemma \ref{lem:translation}, we can suppose without loss of generality that $f_C(x_0(C))=0$ for every $C \in \mathcal{C}$ (and this process does not modify the constants of our previous inequalities). 

\medskip \noindent
Now, let $C$ be an arbitrary clique of $X$, say labelled by $Q \in \mathcal{C}$. Set $L_C=L_Q$, $f_C= f_Q \circ \phi_C$ and define an action $\mathrm{stab}(C) \curvearrowright L_C$ in the following way: if $Q \in \mathcal{C}_2$, take the trivial action; and if $Q \in \mathcal{C}_1$, define the action from $\mathrm{stab}(Q) \curvearrowright L_Q$ by using the isomorphism $\varphi_C : \mathrm{stab}(C) \to \mathrm{stab}(Q)$, ie., $g \cdot x = \varphi_C(g) \cdot x$ for every $g \in \mathrm{stab}(C)$ and $x \in L_C=L_Q$. (Recall that the functions $\phi_C$ and $\varphi_C$ were defined in Definition \ref{def:mapsfortopical}.) Observe that, since the maps $\phi_C$ define isometries $(C,\delta_C) \to (Q, \delta_Q)$ according to Fact \ref{fact:phiCequiisom}, we know that, for every $0< \eta < \min\limits_{C\in \mathcal{C}} \alpha_p^*(C,\delta_C)$, the inequalities
$$A_{\eta} \cdot \delta_C(x,y)^{\alpha_p(f_Q)- \eta} \leq \| f_C(x)- f_C(y) \| \leq B \cdot \delta_C(x,y)$$
hold for every clique $C$ labelled by $Q \in \mathcal{C}$ and every $x,y \in C$. In particular, $\alpha_p(f_C)= \alpha_p(f_Q)$. Next, notice that $f_C$ is $\mathrm{stab}(C)$-equivariant. Indeed, this is clear if $Q \in \mathcal{C}_2$, and if $Q \in \mathcal{C}_1$, we deduce from the point $(ii)$ of Claim \ref{claim:divers} that, for every $x \in L_C$ and $g \in \mathrm{stab}(C)$, we have
$$f_C(gx)=f_Q \circ \phi_C(gx) = f_Q( \varphi_C(g) \cdot \phi_C(x))= \varphi_C(g) \cdot f_Q \circ \phi_C(x)= g \cdot f_C(x),$$
by noticing that $s_C(g)=g$ since $g \in \mathrm{stab}(C)$. Notice also that, for every clique $C$, 
$$f_C(x_0(C))= f_Q \circ \phi_C(x_0(C))= f_Q(x_0(Q))=0.$$
Moreover, once again as a consequence of Claim \ref{claim:divers}, our system of  Lipschitz $\ell^p$-maps is coherent. Indeed, if $C,C'$ are two cliques dual to the same hyperplane, then they are labelled by the same clique of $\mathcal{C}$ say $Q$, so that $L_C=L_{C'}$, and
$$f_C= f_Q \circ \phi_C = f_Q \circ \phi_{C'} \circ t_{C \to C'} = f_{C'} \circ t_{C \to C'}.$$
Let us prove a few claims about the functions we have defined.

\begin{claim}\label{claim:compression1}
Let $C,C'$ be two cliques dual to the same hyperplane $J$.
\begin{itemize}
	\item[(i)] For every $g \in G$ and $h \in L_C=L_{C'}$, $s_C(g) \cdot h = s_{C'}(g) \cdot h$ holds.
	\item[(ii)] The equality  $f_C \circ \mathrm{proj}_C = f_{C'} \circ \mathrm{proj}_{C'}$ holds.
\end{itemize}
\end{claim}

\noindent
Point $(i)$ is clear if $C$ and $C'$ are labelled by $\mathcal{C}_2$ because the actions $\mathrm{stab}(C) \curvearrowright L_C$ and $\mathrm{stab}(C') \curvearrowright L_{C'}$ are trivial; and otherwise, it follows from the observation that, as a consequence of Lemma \ref{lem:varphisC} below, we have the equality $\varphi_C(s_C(g)) = \varphi_{C'}(s_{C'}(g))$, ie., the images of $s_C(g)$ and $s_{C'}(g)$ into the group of affine isometries of $L_C=L_{C'}$ coincide. Next, using Claim \ref{claim:divers}, we find that
$$f_C \circ \mathrm{proj}_C= f_Q \circ \phi_C \circ t_{C \to C'} \circ \mathrm{proj}_{C'} = f_Q \circ \phi_{C'} \circ \mathrm{proj}_{C'}= f_{C'} \circ \mathrm{proj}_{C'},$$
where $Q$ is the clique of $\mathcal{C}$ labelling both $C$ and $C'$, whence Point $(ii)$. 

\begin{claim}\label{claim:compression2}
Let $C$ be a clique and $g \in G$. The following statements hold:
\begin{itemize}
	\item[(i)] $L_C=L_{gC}$;
	\item[(ii)] for every $h \in \mathrm{stab}(C)$, the image of $h$ in $\mathrm{Isom}(L_C)$ is the same as the image of $p_C(g)hp_C(g)^{-1}$ in $\mathrm{Isom}(L_{gC})= \mathrm{Isom}(L_C)$; 
	\item[(iii)] for every $x \in X$, $f_{gC}(gx)=f_C(s_C(g)x)$.
\end{itemize}
\end{claim}

\noindent
Point $(i)$ follows by noticing that $C$ and $gC$ are labelled by the same clique of $\mathcal{C}$, say $Q$, so $L_C=L_Q=L_{gC}$. Next, Point $(ii)$ is clear if $C$ is labelled by $\mathcal{C}_2$ because the images we consider are both trivial; and otherwise it follows from Claim \ref{claim:divers} stating that $\varphi_C(h)= \varphi_{gC} \left( p_C(g) h p_C(g)^{-1} \right)$. Finally, notice that we get
$$\phi_{gC}(gx)= \varphi_C(s_C(g)) \cdot \phi_C(x) = \phi_C(s_C(g) x)$$
by applying Claim \ref{claim:divers} twice, so that 
$$f_{gC}(gx)= f_Q \circ \phi_{gC}(gx) = f_Q \circ \phi_C(s_C(g) x) = f_C(s_C(g)x),$$  
whence Point $(iii)$, concluding the proof of our claim. 

\medskip \noindent
Now we are ready to extend our collection of Lipschitz maps to a global embedding. For convenience, we denote by $\mathfrak{J}$ the set of all the hyperplanes of $X$. 
For every hyperplane $J \in \mathfrak{J}$, choose a clique $C(J)$. Set $L_J=L_{C(J)}$ and $f_J= f_{C(J)} \circ \mathrm{proj}_{C(J)} : X \to L_J$. Now, define
$$f : \left\{ \begin{array}{ccc} X & \to & L= \bigoplus\limits_{J \in \mathfrak{J}}^{\ell^p} L_J \\ x & \mapsto & \left( f_J(x)\right)_J \end{array} \right. .$$
Since $f_J(x_0) =0$ for every $J \in \mathfrak{J}$, our map $f$ is the same as the map constructed in the proof of Proposition \ref{prop:compression}, so that Claim \ref{claim:compressionf} applies. Therefore, the map $f$ is well-defined, Lipschitz and satisfies
$$ \alpha(f) \geq \min \left( \frac{1}{p}, \min\limits_{C~\text{clique}} \alpha(f_C) \right) \geq \min \left( \frac{1}{p}, \min\limits_{C~\text{clique}} \alpha_p^*(C, \delta_C)- \epsilon \right).$$
Therefore, if we prove that $G$ acts on $L$ by affine isometries so that $f$ is $G$-invariant, it will follow that
$$\alpha^*_p(X, \delta) \geq \alpha(f) \geq  \min \left( \frac{1}{p}, \min\limits_{C~\text{clique}} \alpha_p^*(C, \delta_C)- \epsilon \right),$$
and because this is true for every $\epsilon >0$, we will conclude that
$$\alpha_p^*(X, \delta) \geq \min \left( \frac{1}{p}, \min\limits_{\text{$C$ clique}} \alpha_p^*(C, \delta_C) \right).$$
Let us begin by defining an action of $G$ on $L$. If $g \in G$ and $(h_J)_J \in L$, set
$$g \cdot (h_J)_J= \left( s_{C(J)}(g) \cdot h_J \right)_{gJ}.$$
(Recall that $L_{C(J)}= L_{gC(J)}=L_{C(gJ)}$ for every $J \in \mathfrak{J}$, so the expression makes sense.) Notice that Lemma \ref{lem:allbut} below implies that $g$ is identity on all but finitely many coordinates, so $g \cdot (h_J)_J$ belongs to $L$. Moreover, $g$ is an affine isometry of $L$ since it is a affine isometry on each coordinate. We claim that this defines an action $G \curvearrowright L$. Let $g_1,g_2 \in G$ and $(h_J)_J \in L$. By using Corollary \ref{cor:sCgh} and the point $(ii)$ of Claim \ref{claim:compression2}, we find that
$$\begin{array}{lcl} g_1g_2 \cdot (h_J)_J & = & \left(s_{C(J)}(g_1g_2) \cdot h_J \right)_{g_1g_2J} \\ \\ & = & \left( p_{C(J)}(g_2)^{-1}s_{g_2C(J)}(g_1) p_{C(J)}(g_2) \cdot s_{C(J)}(g_2) \cdot h_J \right)_{g_1g_2J} \\ \\ & = & \left( s_{g_2C(J)}(g_1) s_{C(J)}(g_2) \cdot h_J \right)_{g_1g_2J} \end{array}$$
Moreover, it follows from the point $(i)$ of Claim \ref{claim:compression1} that
$$s_{g_2C(J)}(g_1) s_{C(J)}(g_2) \cdot h_J= s_{C(g_2J)}(g_1) s_{C(J)}(g_2) \cdot h_J $$
for every $J \in \mathfrak{J}$, so we deduce that
$$g_1g_2 \cdot (h_J)_J = \left(  s_{C(g_2J)}(g_1) s_{C(J)}(g_2) \cdot h_J \right)_{g_1g_2J}= g_1 \cdot \left( s_{C(J)}(g_2) \cdot h_J \right)_{g_2J}= g_1 \cdot \left( g_2 \cdot \left(h_J \right)_J \right).$$
Thus, we have defined an action $G \curvearrowright L$. Now, we want to prove that $f$ is $G$-equivariant. Let $x \in X$ be any vertex and $g \in G$. For any hyperplane $J$, by using the point $(iii)$ of Claim \ref{claim:compression2} and the point $(ii)$ of Claim \ref{claim:compression1}, we find that
$$\begin{array}{lcl} f_J(gx) & = & f_{C(J)}( \mathrm{proj}_{C(J)} (gx)) = f_{gg^{-1}C(J)} \left( g \cdot \mathrm{proj}_{g^{-1}C(J)}(x) \right) \\ \\ & = & f_{g^{-1}C(J)} \left( s_{g^{-1}C(J)}(g) \cdot \mathrm{proj}_{g^{-1}C(J)}(x) \right) \\ \\ & = & s_{g^{-1}C(J)}(g) \cdot f_{g^{-1}C(J)} \left( \mathrm{proj}_{g^{-1}C(J)} (x) \right) \\ \\ & =& s_{g^{-1}C(J)}(g) \cdot f_{C(g^{-1}J)} \left( \mathrm{proj}_{C(g^{-1}J)} (x) \right)= s_{g^{-1}C(J)}(g) \cdot f_{g^{-1}J}(x) \end{array}$$
Finally, we deduce from the point $(i)$ of Claim \ref{claim:compression1} that
$$f_J(gx)=s_{C(g^{-1}J)}(g) \cdot f_{g^{-1}J}(x).$$
From this observation, it follows that,
$$ \begin{array}{lcl} f(gx) & = & \left( f_J(gx) \right)_J = \left( s_{C(g^{-1}J)}(g) \cdot f_{g^{-1}J}(x)  \right)_{gg^{-1}J} \\ \\ & = & g \cdot \left( f_{g^{-1}J}(x) \right)_{g^{-1}J} =  g \cdot \left( f_J(x) \right)_J=g \cdot f(x)\end{array}$$
Thus, we have proved that
$$\alpha_p^*(X, \delta) \geq \min \left( \frac{1}{p}, \min\limits_{\text{$C$ clique}} \alpha_p^*(C, \delta_C) \right).$$
On the other hand, notice that, if $C$ is a clique labelled by $\mathcal{C}_2$, then $\alpha_p^*(C,\delta_C) = 1$ because both $C$ and $\mathrm{stab}(C)$ are finite. Next, if $C$ is a clique labelled by some $Q \in \mathcal{C}_1$, then $\phi_C$ defines a $\varphi_C$-equivariant isometry $(C,\delta_C) \to (Q, \delta_Q)$ according to Fact \ref{fact:phiCequiisom}, so that $\alpha_p^*(C,\delta_C)= \alpha_p^*(Q,\delta_Q)$; and because $g \mapsto g \cdot x_0(C)$ induces an equivariant isometry $\mathrm{stab}(Q) \to (Q, \delta_Q)$ by construction, we know that $\alpha_p^*(Q, \delta_Q)= \alpha_p^* (\mathrm{stab}(Q))$. Hence $\alpha_p^*(C,\delta_C) = \alpha_p^*(\mathrm{stab}(Q))$. Therefore,
$$\min\limits_{\text{$C$ clique}} \alpha_p^*(C, \delta_C) = \min\limits_{C \in \mathcal{C}} \alpha_p^*( \mathrm{stab}(C)) = \min\limits_{C \in \mathcal{C}_1} \alpha_p^*( \mathrm{stab}(C)),$$
which concludes the proof.
\end{proof}

\begin{lemma}\label{lem:varphisC}
Let $G$ be a group acting topically-transitively on a quasi-median graph $X$ and $\mathcal{C}$ a collection of cliques such that any $G$-orbit of hyperplanes intersects it along a single clique. Let $C,C' \subset X$ be two cliques dual to the same hyperplane. For every $h \in G$, $\varphi_C(s_C(h))= \varphi_{C'}(s_{C'}(h))$. 
\end{lemma} 

\begin{proof}
Let $\mathcal{C}= \mathcal{C}_1 \sqcup \mathcal{C}_2$ denote the associated decomposition of $\mathcal{C}$. If $C$ and $C'$ are labelled by a clique of $\mathcal{C}_2$, there is nothing to prove, since $\varphi_C$ and $\varphi_{C'}$ are trivial. So suppose that $C$ and $C'$ are labelled by a clique $Q \in \mathcal{C}_1$. For convenience, let $a=h^{-1} \cdot x_0(hC)$ and $b=h^{-1}x_0(hC')$. By applying the points $(i)$ and $(ii)$ of Claim \ref{claim:divers}, we find that
$$\begin{array}{lcl} \varphi_C(s_C(h)) \cdot \phi_C(a) & = & \phi_{hC}(ha) = \phi_{hC}(x_0(hC)) \\ \\ & = & \phi_{hC} \left( t_{hC' \to hC} (x_0(hC')) \right) \\ \\ & = & \phi_{hC'} (x_0(hC')) = \phi_{hC'}(hb) \\ \\ & = & \varphi_{C'}(s_{C'}(h)) \cdot \phi_{C'}(b) \end{array}$$
On the other hand, since $x_0(hC)$ and $x_0(hC')$ belong to the same sector delimited by $hJ$, where $J$ denotes the hyperplane dual to $C$ and $C'$, a fortiori $a$ and $b$ belong to the same sector delimited by $J$, hence $b=t_{C \to C'}(a)$. As a consequence,
$$\phi_{C'}(b)= \phi_{C'} \left( t_{C \to C'} (a) \right) = \phi_C(a),$$
using the point $(i)$ of Claim \ref{claim:divers}. Thus, we have proved that
$$\varphi_C(s_C(h)) \cdot \phi_C(a) = \varphi_{C'} (s_{C'}(h)) \cdot \phi_C(a).$$
Since $\mathrm{stab}(Q)$ acts freely on $Q$, we conclude that $\varphi_C(s_C(h))= \varphi_{C'}(s_{C'}(h))$.
\end{proof}

\begin{lemma}\label{lem:allbut}
Let $G$ be a group acting topically-transitively on a quasi-median graph $X$ and $\mathcal{C}$ a collection of cliques such that any $G$-orbit of hyperplanes intersects it along a single clique. Fix a collection $\mathcal{Q}$ of cliques such that any hyperplane of $X$ contains exactly one clique of $\mathcal{Q}$. For every $g \in G$, $s_C(g)=1$ for all but finitely many $C \in \mathcal{Q}$. 
\end{lemma}

\begin{proof}
Let $\mathcal{C}= \mathcal{C}_1 \sqcup \mathcal{C}_2$ denote the associated decomposition of $\mathcal{C}$. Notice that $\mathrm{proj}_{gC}(x_0)=x_0(gC)$ and $\mathrm{proj}_{gC}(gx_0) = g \cdot \mathrm{proj}_C(x_0)=g x_0(C)$. But 
$$x_0(gC)= p_C(g) x_0(C) \neq p_C(g)s_C(g) x_0(C) = g x_0(C)$$
if $s_C(g) \neq 1$. Indeed, $s_C(g) \neq 1$ implies that $C$ is necessarily labelled by some clique of $\mathcal{C}_1$, which implies that $\mathrm{stab}(C)$ acts freely on the vertices of $C$, hence $s_C(g) x_0(C) \neq x_0(C)$. A fortiori, since $x_0$ and $gx_0$ have different projections onto $gC$, the hyperplane $gJ$ must separate $x_0$ and $gx_0$ where $J$ denotes the hyperplane dual to $C$. Equivalently, $J$ separates $x_0$ and $g^{-1}x_0$. Therefore, $s_C(g)=1$ for every clique $C$ whose dual hyperplane does not separate $x_0$ and $g^{-1}x_0$. Because there exist only finitely many hyperplanes separating $x_0$ and $g^{-1}x_0$, the conclusion follows.
\end{proof}

\begin{remark}\label{remark:whichsystem}
In view of the uniqueness provided by Theorem \ref{thm:extensionunique}, notice that Proposition \ref{prop:equicompression} applies to any system of metrics defined in the following way: Fix a basepoint $x_0$, and set $x_0(C)= \mathrm{proj}_C(x_0)$ for every clique $C$. Endow any clique $C \in \mathcal{C}_2$ with the discrete metric $(x,y) \mapsto \left\{ \begin{array}{cl} 1 & \text{if $x \neq y$} \\ 0  & \text{otherwise} \end{array} \right.$, and any clique $C \in \mathcal{C}_1$ with the word metric of $\mathrm{stab}(C)$ with respect to some finite generating set thanks to the map $g \mapsto g \cdot x_0(C)$.
\end{remark}

\begin{remark}
Under the hypotheses of Proposition \ref{prop:equicompression}, we can similarly prove the following lower bound on the non equivariant $\ell^p$-compression of $(X,\delta)$:
$$\alpha_p(X,\delta) \geq \min \left( \frac{1}{p}, \inf\limits_{C \in \mathcal{C}} \alpha_p( \mathrm{stab}(C)) \right).$$
However, if the cubical dimension of $X$ is finite, we think that this bound can be improved by removing the term $1/p$. See Question \ref{question:compression}.
\end{remark}

\section{Inflating the hyperplanes of a quasi-median graph}\label{section:inflatingintro}

\noindent
This section is dedicated to a geometric construction, performed on quasi-median graphs, which will be fundamental in the next section.

\subsection{Gluing quasi-median graphs}

\noindent
The construction we focus on in this section is probably well-known from the experts, at least for finite graphs. Using the terminology introduced in \cite{quasimedian}, we prove that the \emph{gated amalgam} of two quasi-median graphs is quasi-median, although we formulate it in a more geometric way, as a gluing construction.

\medskip \noindent
In order to fix the notation, given three graphs $A,B,C$ and two embeddings $\phi : C \hookrightarrow A$ and $\psi : C \hookrightarrow B$, we define the \emph{gluing} $A \underset{C}{\star} B$ as the quotient
$$(A \sqcup B) / ( a \sim b \ \text{if there exists $c \in C$ such that $a=\phi(c)$ and $b=\psi(c)$}).$$
For convenience, $A,B,C$ will be identified with their images in $A \underset{C}{\star} B$. In particular, $A \cap B=C$. 

\begin{prop}\label{prop:gluing}
Let $A_1,A_2$ be two quasi-median graphs and $C$ a graph with two isometric embedding $\varphi_1 : C \hookrightarrow A_1$, $\varphi_2 : C \to A_2$. Suppose that $\varphi_1(C)$ (resp. $\varphi_2(C)$) is gated in $A_1$ (resp. $A_2$). Then $A_1 \underset{C}{\star} A_2$ is a quasi-median graph. Moreover, $A_1$ and $A_2$ are gated in $A_1 \underset{C}{\star} A_2$.
\end{prop}

\begin{proof}
For convenience, set $X= A_1 \underset{C}{\star} A_2$. We first claim that $A_1$ is gated in $X$. Let $x \in X$. If $x \in A_1$, then clearly $x$ is gate of $x$ in $A_1$. Now, suppose that $x \in A_2$, and let $x'$ denote the gate of $x$ in $\varphi_1(C)=A_1 \cap A_2$. If $y \in A_1$, fix a geodesic $[x,y]$ between $x$ and $y$. Let $x''$ denote some vertex of $[x,y]$ which belongs to $A_1 \cap A_2$. By definition of $x'$, there exists a geodesic between $x$ and $x''$ passing through $x'$. By replacing the subsegment $[x,x''] \subset [x,y]$ with this new geodesic, we produce a geodesic between $x$ and $y$ passing through $x'$. Consequently, $x'$ is a gate of $x$ in $A_1$. Thus, we have proved that $A_1$ is gated in $X$. Similarly, we show that $A_2$ is gated in $X$. 

\medskip \noindent
Now, we want to verify the triangle condition. So let $p,q \in X$ be two adjacent vertices and $u \in X$ a third vertex satisfying $d(u,p)=d(u,q)=k$. If either $u,p,q \in A_1$ or $u,p,q \in A_2$, there is nothing to prove since we already know that $A_1$ and $A_2$ are quasi-median (and because, as a consequence of the previous paragraph, $A_1$ and $A_2$ are isometrically embedded subgraphs of $X$). Otherwise, noticing that no vertex of $A_1 \backslash (A_1 \cap A_2)$ is adjacent to a vertex of $A_2 \backslash (A_1 \cap A_2)$, we have either $u \in A_1$ and $p,q \in A_2$, or $u \in A_2$ and $p,q \in A_1$. Say we are in the former case, the latter being symmetric. Let $u'$ denote the gate of $u$ in $A_2$. Then $d(u',p)=d(u',q)=k-d(u,u')$. Since $A_2$ satisfies the triangle condition, there must exist some vertex $r \in A_2$ adjacent to both $p$ and $q$ such that $d(u',r)=k-d(u,u')-1$. Noticing that $d(u,r)=d(u,u')+d(u',r)=k-1$, we conclude that $r$ is the vertex we are looking for. 

\medskip \noindent
Then, we want to verify the quadrangle condition. First, let us prove the following fact.

\begin{fact}\label{fact:ingluing}
If $p,q,r,u \in X$ are four vertices such that $r$ is adjacent to both $p$ and $q$ and that $d(u,p)=d(u,q)=k$ and $d(u,r)=k+1$, then either $p,q,r \in A_1$ or $p,q,r \in A_2$.
\end{fact}

\noindent
Suppose that $p \in A_1$ and $u,q \in A_2$. If $q \in A_1 \cap A_2$, then we deduce that $p,q,r \in A_1$ since $A_1$ is gated, and we are done. So we suppose that $q \notin A_1$. Similarly, if $p \in A_1 \cap A_2$, then necessarily $p,q,r \in A_2$ since $A_2$ is gated, so that we suppose that $p \notin A_2$. As a consequence, $r \in A_1 \cap A_2$. Let $u'$ denote the gate of $u$ in $A_1$. Noticing that
$$\left\{ \begin{array}{l} d(u',r)=d(u,r)-d(u,u')=k+1-d(u,u') \\ d(u',p)= d(u,p)-d(u,u')=k-d(u,u') \end{array} \right.$$
we deduce that $d(u',r)=d(u',p)+1=d(u',p)+d(p,r)$. Therefore, $p$ belongs to a geodesic between $u',r \in A_1 \cap A_2$. Since $A_2$ is convex, we get a contradiction because we supposed $p \notin A_2$. This concludes the proof of our fact.

\medskip \noindent
Now, in order to verify the quadrangle condition in $X$, fix four vertices $u,p,q,r \in X$ such that $r$ is adjacent to both $p$ and $q$, and that $d(u,p)=d(u,q)=k$ and $d(u,r)=k+1$. If either $u,p,q,r \in A_1$ or $u,p,q,r \in A_2$, then there is nothing to prove. Otherwise, according to the previous fact, either $p,q,r \in A_1$ and $u \in A_2$, or $p,q,r \in A_2$ and $u \in A_1$. Say we are in the latter case, the former being symmetric. We argue similarly as for the triangle condition: if $u'$ denotes the gate of $u$ in $A_2$, then we apply the quadrangle condition in $A_2$ to $u',p,q,r$, producing a vertex $s$ adjacent to both $p$ and $q$, which is precisely the vertex we are looking for.

\medskip \noindent
Finally, suppose that by contradiction that $X$ contains an induced subgraph $K$ isomophic to $K_{2,3}$. It follows from Fact \ref{fact:ingluing} that some two consecutive edges of $K$ both belong to either $A_1$ or $A_2$, say $A_1$. Because $A_1$ is convex, we deduce that $K \subset A_1$, contradicting that $A_1$ is quasi-median. Similarly, if we suppose by contradiction that $X$ contains an induced subgraph $K$ isomorphic to $K_4^-$, then we deduce from Fact \ref{fact:ingluing} that some two consecutive edges of $K$ both belong to either $A_1$ or $A_2$, say $A_1$, so that, since $A_1$ contains its triangles, $K$ must be included into $A_1$. We get a contradiction with the fact that $A_1$ is quasi-median.
\end{proof}

\subsection{Inflating one hyperplane}

\noindent
In this section, we prove that it is possible to add new vertices to the cliques dual to a given hyperplane and to get a new quasi-median graph.

\begin{definition}
Let $X$ be a quasi-median graph, $J$ one of its hyperplanes and $S(J)$ an arbitrary set. For every clique $C$ dual to $J$, fix a copy of $S(J)$ denoted by $S(C)$, and for every pair of parallel cliques $C$ and $C'$, fix a bijection $f_{C \to C'}  : S(C) \to S(C')$, so that the following \emph{monodromy condition} is satisfied: for every sequence of successively parallel cliques $C_1, \ldots, C_n$ where $C_1=C_n=C$, 
$$f_{C_{n-1} \to C_n} \circ f_{C_{n-2} \to C_{n-1}} \circ \cdots \circ f_{C_2 \to C_3} \circ f_{C_1 \to C_2} = \mathrm{id}_{S(C)};$$ 
for instance, set $f_{C \to C'}= \mathrm{id}_{S(J)}$. We define a new graph $Y$ from $X$ by adding the vertices $\bigsqcup S(C)$, edges between any two vertices of $C \cup S(C)$, and an edge between $x$ and $f_{C \to C'}(x)$ for every cliques $C,C'$ which are opposite in some prism and every vertex $x \in C$. We refer to $Y$ as a graph obtained from $X$ by \emph{inflating the hyperplane $J$}. 
\end{definition}

\begin{figure}
\begin{center}
\includegraphics[scale=0.6]{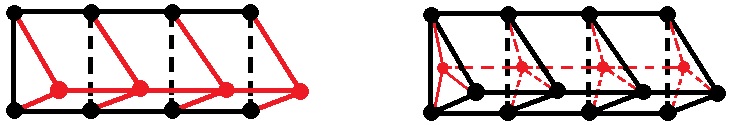}
\end{center}
\caption{Examples of inflations.}
\label{figure10}
\end{figure}

\noindent
It is worth noticing that the choices of our bijections $f_{C \to C'}$ does not disturbe $Y$ up to isometry (provided the monodromy condition is satisfied). Therefore, the isometry class of a graph obtained from a quasi-median graph by inflating one of its hyperplanes depends only on the cardinality of our set $S(J)$. Figure \ref{figure10} gives examples of inflations (adding just one vertex for each clique) of $K_2 \times P_3$ and $K_3 \times P_3$. 

\begin{prop}\label{prop:inflatingonehyp}
Let $Y$ be a graph obtained from a quasi-median graph $X$ by inflating one of its hyperplanes $J$. Then $Y$ is a quasi-median graph in which $X$ is convex, and the map $J \mapsto J \cap X$ defines a bijection from the hyperplanes of $Y$ onto those of $X$. 
\end{prop}

\begin{proof}
For every clique $C$ of $X$ which is dual to $J$, let $C_+$ denote the complete subgraph generated by $C \cup S(C)$. Let $Z$ denote the subgraph of $Y$ generated by
$$\{ C_+ \mid C \ \text{clique of $X$ dual to $J$} \}.$$
Finally, let $\mathcal{S}$ denote the collection of all the sectors delimited by $J$ in $X$. For convenience, we index the sectors of $\mathcal{S}$ by some ordinal $\xi$, ie., $\mathcal{S}= \{ S_{\alpha} \mid 0 < \alpha \leq \xi \}$. For every $\alpha \leq \xi$, set $F_{\alpha}= S_{\alpha} \cap Z$. Notice that, because $S_{\alpha} \subset X$, we also have $F_{\alpha}=S_{\alpha} \cap Z \cap X=S_{\alpha} \cap N(J)$. 

\medskip \noindent
Fixing some clique $C$ of $X$, we can define a projection $p : \bigcup\limits_{C' \ \text{clique}} S(C') \to S(C)$ in the following way: if $x \in S(C')$, set $$p(x)=f_{C_{n-1} \to C_n} \circ \cdots \circ f_{C_1 \to C_2}(x),$$ where $C_1, \ldots, C_n$ is a sequence of successively parallel cliques between $C_1=C'$ and $C_n=C$. The monodromy condition ensures that $p(x)$ does not depend on our choice of the sequence of cliques. Notice that the restriction of $p$ to some $S(C')$ defines a bijection $S(C') \to S(C)$. Moreover, Lemma \ref{lem:carrierproduct} gives an isometric embedding $\Psi : N(J) \hookrightarrow F(J) \times C_+$. This allows us to define
$$\Phi : \left\{ \begin{array}{ccc} Z & \longrightarrow & F(J) \times C_+ \\ x & \longmapsto & \left\{ \begin{array}{cl} \Psi(x) & \text{if} \ x \in N(J) \\ (C',p(x)) & \text{if} \ x \in S(C') \end{array} \right. \end{array} \right.$$

\begin{claim}\label{claim:Phi}
The map $\Phi : Z \to C_+ \times F(J)$ defines an isomorphism. Moreover, for every $\alpha \leq \xi$, there exists a vertex $v_{\alpha} \in C \subset C_+$ such that $F_{\alpha}= \Phi^{-1}(F(J) \times \{v_{\alpha} \})$.
\end{claim}

\noindent
First, we prove that $\Phi$ is injective. Let $x,y \in Z$ be two distinct vertices such that $x \in (C_1)_+$ and $y \in(C_2)_+$ for some cliques $C_1,C_2$ of $X$. If $C_1 \neq C_2$, then $\Phi(x)$ and $\Phi(y)$ differ on their first coordinates, so that $\Phi(x) \neq \Phi(y)$. So suppose that $C_1=C_2$. If $x \in C_1$ and $y \in C_2$, then $\Phi(x) =\Psi(x) \neq \Psi(y)=\Phi(y)$ since we already know that $\Psi$ is injective. If $x \in C_1$ and $y \notin C_2$, then the second coordinate of $\Phi(x)$ belong to $C$ whereas the second coordinate of $\Phi(y)$ belongs to $S(C)$, hence $\Phi(x) \neq \Phi(y)$. If $x \notin C_1$ and $y \in C_2$, the situation is symmetric to the previous one. Finally, if $x \notin C_1$ and $y \notin C_2$, $x \neq y$ implies $p(x) \neq p(y)$, hence $\Phi(x) \neq \Phi(y)$. Thus, $\Phi$ is injective.

\medskip \noindent
Now, we want to verify that $\Phi$ is surjective. Let $(C',v) \in F(J) \times C_+$. If $v \in C$, then there exists some $x \in N(J) \subset Z$ such that $\Phi(x)=\Psi(x)=(C',v)$. Otherwise, if $v \in S(C)$, there exists some $x \in S(C')$ such that $p(x)=v$, because $p$ induces a bijection $S(C') \to S(C)$, so that $\Phi(x)=(C',v)$. Thus, $\Phi$ is surjective.

\medskip \noindent
Then, suppose that $x,y \in Z$ are two adjacent vertices. In order to show that $\Phi(x)$ and $\Phi(y)$ are also adjacent, we distinguish three cases:
\begin{itemize}
	\item if $x,y \in N(J)$, then $\Phi(x)= \Psi(x)$ and $\Phi(y)= \Psi(y)$ must be adjacent;
	\item if $x \in S(C')$ and $y \in C'_+$ for some clique $C'$, then $\Phi(x)$ and $\Phi(y)$ belong to the clique $\{C'\} \times C_+$;
	\item if $x \in S(C_1)$ and $y=f_{C_1 \to C_2}(x) \in S(C_2)$ for some cliques $C_1,C_2$, then $\Phi(x)=(C_1,p)$ and $\Phi(y)=(C_2,p)$, where $p(x)=p=p(y)$, so that $\Phi(x)$ and $\Phi(y)$ are adjacent since $C_1$ and $C_2$ are parallel.
\end{itemize}
All the other possible cases are symmetric to one of the three above. Thus, $\Phi$ sends adjacent vertices to adjacent vertices. Conversely, suppose that $\Phi(x)$ and $\Phi(y)$ are adjacent in $F(J) \times C_+$. We distinguish three cases:
\begin{itemize}
	\item if $\Phi(x)$ and $\Phi(y)$ belong to the same clique $\{C' \} \times C_+$, then $x$ and $y$ belong to the complete subgraph $C_+'$, so that $x$ and $y$ must be adjacent;
	\item if $\Phi(x)=(C_1,v)$ and $\Phi(y)=(C_2,v)$, where $C_1$ and $C_2$ are two parallel cliques of $X$ and $v \in C$, then we already know that $x$ and $y$ must be adjacent since $\Phi(x)=\Psi(x)$ and $\Phi(y)= \Psi(y)$;
	\item if $\Phi(x)=(C_1,v)$ and $\Phi(y)=(C_2,v)$, where $C_1$ and $C_2$ are two parallel cliques of $X$ and $v \in S(C)$, then $v=p(y)=p \circ f_{C_1 \to C_2}(x)$, hence $y=f_{C_1 \to C_2}(x)$ because $p$ is injective on $C_2$, so that $x$ and $y$ must be adjacent.
\end{itemize}
All the other possible cases are symmetric to one of the three above. We conclude that $\Phi : Z \to F(J) \times C_+$ defines an isomorphism.

\medskip \noindent
Finally, according Lemma \ref{lem:ccpartialJ}, there exists some $v_{\alpha} \in C$ such that 
$$F_{\alpha}= \Psi^{-1} (F(J) \times \{ v_{\alpha} \}) = \Phi^{-1}( F(J) \times \{v_{\alpha} \}).$$
This concludes the proof of our first claim.

\medskip \noindent
Now, we define by transfinite induction the following increasing sequence of subgraphs
\begin{itemize}
	\item $Z^{(0)}=Z$;
	\item $Z^{(\alpha+1)} = Z^{(\alpha)} \cup S_{\alpha+1}$ for every $\alpha< \xi$;
	\item $Z^{(\lambda)}= \bigcup\limits_{\alpha< \lambda} Z^{(\alpha)}$ for every limit ordinal $\lambda \leq \xi$.
\end{itemize}
Notice that $Y=Z^{(\xi)}$. Indeed, fixing an edge $e \subset Y$, three cases may happen. Either $e$ has an endpoint which does not belong to $X$, so that $e \subset Z$; or $e \subset X$ is dual to $J$, so that $e \subset N(J) \subset Z$; or $e \subset X$ is not dual to $J$, so that $e \subset S_{\alpha} \subset Z^{(\alpha)}$ for some $\alpha \leq \xi$.

\begin{claim}\label{claim:amalgam}
For every $\alpha< \xi$, $Z^{(\alpha+1)} = Z^{(\alpha)} \underset{F_{\alpha+1}}{\star} S_{\alpha+1}$.
\end{claim}

\noindent
First, let us prove by transfinite induction that, $\beta$ being a fixed successor ordinal, $$Z^{(\alpha)} \cap S_{\beta} = \left\{ \begin{array}{cl} F_{\beta} & \text{if} \ \alpha < \beta \\ S_{\beta} & \text{otherwise} \end{array} \right.$$ for every $\alpha \leq \xi$. Of course, $Z^{(0)} \cap S_{\beta}= Z \cap S_{\beta} = F_{\beta}$ just by definition of $F_{\beta}$. Then, if we suppose that our assertion is true for $Z^{(\alpha)}$, we have
$$\begin{array}{lcl} Z^{(\alpha+1)} \cap S_{\beta} & = & (Z^{(\alpha)} \cup S_{\alpha+1}) \cap S_{\beta} = (Z^{(\alpha)} \cap S_{\beta}) \cup (S_{\alpha+1} \cap S_{\beta}) \\ \\ & = & \left\{ \begin{array}{cl} F_{\beta} & \text{if} \ \alpha < \beta \\ S_{\beta} & \text{otherwise} \end{array} \right. \bigcup \left\{ \begin{array}{cl} \emptyset & \text{if} \ \alpha+1 \neq \beta \\ S_{\beta} & \text{otherwise} \end{array} \right. = \left\{ \begin{array}{cl} F_{\beta} & \text{if} \ \alpha+1 < \beta \\ S_{\beta} & \text{otherwise} \end{array} \right. ; \end{array} $$
and if $\lambda$ is a limit ordinal such that our assertion is true for every $\alpha< \lambda$,
$$Z^{(\lambda)} \cap S_{\beta} = \bigcup\limits_{\alpha< \lambda} \left( Z^{(\alpha)} \cap S_{\beta} \right) = \left\{ \begin{array}{cl} F_{\beta} & \text{if} \ \lambda \leq \beta \\ S_{\beta} & \text{otherwise} \end{array} \right. ;$$
notice that, because $\beta$ is not a limit ordinal, $\lambda \leq \beta$ is equivalent to $\lambda< \beta$. This concludes the proof of our assertion. As a consequence, we know that $Z^{(\alpha)} \cap S_{\alpha+1}= F_{\alpha}$ for every $\alpha < \xi$. 

\medskip \noindent
To conclude the proof of our claim, it is sufficient to show that no vertex of $S_{\alpha+1} \backslash F_{\alpha+1}$ is adjacent to a vertex of $Z^{(\alpha)} \backslash F_{\alpha+1}$. We already know that no vertex of $S_{\alpha+1} \backslash F_{\alpha+1} \subset X$ is adjacent to a vertex of $X \backslash S_{\alpha+1} \supset X \cap \left( Z^{(\alpha)} \backslash F_{\alpha+1} \right)$ since the projection (in $X$) onto $N(J)$ of any vertex of $S_{\alpha+1}$ belongs to $F_{\alpha+1}$. Moreover, the vertices of $N(J)$ are the only vertices of $X$ adjacent to new vertices, so that, since $N(J) \cap S_{\alpha+1}= F_{\alpha+1}$, no vertex of $S_{\alpha+1} \backslash F_{\alpha+1}$ is adjacent to a vertex of $Y \backslash X \supset X \backslash \left( Z^{(\alpha)} \backslash F_{\alpha+1} \right)$. This concludes the proof.

\medskip \noindent
Now, we are able to prove by transfinite induction that, for every $\alpha \leq \xi$, $Z^{(\alpha)}$ is a quasi-median graph in which $Z^{(\beta)}$ is gated for every $\beta< \alpha$. First, as a consequence of Claim \ref{claim:Phi} we know that $Z^{(0)}=Z$ is isomorphic to $F(J) \times C_+$; on the other hand, $F(J)$ is a quasi-median graph according to Proposition \ref{prop:hypsumup}. We deduce that $Z$ is a quasi-median graph as a Cartesian product of two quasi-median graphs. 

\medskip \noindent
Next, suppose that our assertion holds for $Z^{(\alpha)}$. According to Claim \ref{claim:amalgam}, we know that $Z^{(\alpha+1)} = Z^{(\alpha)} \underset{F_{\alpha+1}}{\star} S_{\alpha}$. Moreover, we deduce from Claim \ref{claim:Phi} that $F_{\alpha+1}= \Phi^{-1} (F(J) \times \{ v_{\alpha+1} \})$ for some vertex $v_{\alpha+1} \in C_+$, so that $F_{\alpha+1}$ must be gated in $Z$. On the other hand, by our induction hypothesis, $Z$ is also gated in $Z^{(\alpha)}$ so $F_{\alpha+1}$ must be gated in $Z^{(\alpha)}$. Finally, we deduce from Proposition \ref{prop:hypsumup} that $F_{\alpha+1}$ is gated in $S_{\alpha+1}$, so that Proposition \ref{prop:gluing} implies that $Z^{(\alpha+1)}$ is a quasi-median graph in which $Z^{(\alpha)}$ is gated. Moreover, we know by our induction hypothesis that, for $\beta< \alpha$, $Z^{(\beta)}$ is gated in $Z^{(\alpha)}$, so $Z^{(\beta)}$ is gated in $Z^{(\alpha+1)}$ as well.

\medskip \noindent
Finally, suppose that $\lambda$ is a limit ordinal such that our assertion holds for every $\alpha < \lambda$. Then $Z^{(\lambda)}$ is an increasing sequence of quasi-median graphs $\bigcup\limits_{\alpha< \lambda} Z^{(\alpha)}$ such that, for every $\beta < \gamma < \lambda$, $Z^{(\beta)}$ is gated in $Z^{(\gamma)}$. We deduce that $Z^{(\lambda)}$ is a quasi-median graph in which $Z^{(\beta)}$ is gated for every $\beta< \lambda$. This concludes the proof of our assertion.

\medskip \noindent
As a consequence, $Y=Z^{(\xi)}$ is a quasi-median graph. This concludes the proof of the first assertion of our proposition. Now, let us show that $X$ is convex in $Y$.

\medskip \noindent
Let $x,y \in X \subset Y$ be two vertices, and let $x_1=x,x_2, \ldots, x_{n-1},x_n=y$ be a path in $Y$ between $x$ and $y$ which is not included into $X$. In particular, there must exist a subpath $x_p,x_{p+1}, \ldots, x_{p+r-1},x_{p+r}$ such that $r \geq 2$, $x_p,x_{p+r} \in X$ and $x_{p+1}, \ldots, x_{p+r-1} \notin X$. Notice that necessarily $x_p, \ldots, x_{p+r} \in Z$. Let $\Phi : Z \to F(J) \times C_+$ be the isomorphism given by Claim \ref{claim:Phi}. For every $0 \leq i \leq r$, set $\Phi(x_{p+i})=(C_i,v_i)$. Because $x_p \in X$, we deduce that $v_0 \in C_0$. This allows us to define the sequence of vertices $w_0, \ldots, w_r$ by: $w_0=v_0$ and, for every $1 \leq i \leq $, $w_{i+1} \in C_{i+1}$ is the vertex opposite to $w_i \in C_i$. On the other hand, we deduce from $x_p \in X$ and $x_{p+1} \notin X$ that $C_0=C_1$, $v_0 \in C_0$ and $v_1 \in S(C_0)$; in particular $w_0=w_1$. Similarly, we deduce from $x_{p+r} \in X$ and $x_{p+r-1} \notin X$ that $C_r=C_{r-1}$, $v_r \in C_r$ and $v_{r-1} \in S(C_r)$, so that $w_{r-1}=w_r$. Consequently, the sequence of vertices
$$\Gamma : (C_0,w_0), \ (C_2,w_2), \ (C_3,w_3), \ldots, (C_{r-1},w_{r-1}), \ (C_r,v_r)$$
defines a path of length $r-1$ in $F(J) \times C_+$ between $\Phi(x_p)$ and $\Phi(x_{p+r})$. Thus, by replacing the subsegment $x_p, \ldots, x_{p+r}$ of our initial path with $\Phi^{-1}(\Gamma)$ produces a shorter path between $x$ and $y$. This proves that any path between two vertices of $X \subset Y$ which is not included into $X$ can be shortened. As a consequence, $X$ is convex in $Y$.

\medskip \noindent
To conclude the proof of our proposition, we now focus on the hyperplanes of $Y$. According to Lemma \ref{lem:hypconvexsub} below, the map $J \mapsto J \cap X$ defines a bijection from the hyperplanes of $Y$ intersecting $X$ onto the hyperplanes of $X$. Therefore, it is sufficient to show that any hyperplane of $Y$ intersects $X$. For this purpose, we will prove that any clique of $Y$ either is some $C_+$ or is parallel to a clique of $X$, which is sufficient to conclude.

\medskip \noindent
Let $Q$ be a clique of $Y$. If the vertices of $Q$ all belong to $X$, then $Q \subset X$ and there is nothing to prove. So suppose that $Q$ contains a vertex $q \in Y \backslash X$. From the construction of $Y$, it is clear that all the vertices of $Y$ adjacent to a vertex of $Y \backslash X$ must be contained into $Z$, hence $Q \subset Z$. Let $C$ be the unique clique of $X$ such that $q \in C_+$. By identifying $Z$ with $F(J) \times C_+$ thanks to Claim \ref{claim:Phi}, two cases may happen: either $Q=C_+$, or $Q= K \times \{q \}$ for some clique $K$ of $F(J)$ containing the vertex $C \in F(J)$. In the latter case, fixing some $x \in C$, $Q$ is parallel to $K \times \{x \}$, which is a clique of $X$. This concludes the proof.
\end{proof}

\begin{lemma}\label{lem:hypconvexsub}
Let $Y$ be a quasi-median graph and $X \subset Y$ a convex subgraph. The map $J \mapsto J \cap X$ defines a bijection from the hyperplanes of $Y$ intersecting $X$ onto the hyperplanes of $X$.
\end{lemma}

\begin{proof}
For convenience, if $C$ is a clique of $X$, let $C_+$ denote the unique clique of $Y$ containing $X$. 

\medskip \noindent
We first need to verify that the map $J \mapsto J \cap X$ is well-defined, ie., that the image $J \cap X$ is indeed a hyperplane of $X$. For this purpose, it is sufficient to notice that two cliques $C,C'$ of $X$ are dual to the same hyperplane of $X$ if and only if $C_+,C'_+$ are dual to the same hyperplane of $Y$. The implication is clear since $X$ is a subgraph of $Y$. Conversely, suppose that $C_+$ and $C_+'$ are dual to the same hyperplane of $Y$. Let $Z$ denote the union of all the geodesics between $C$ and $C'$. Because $X$ is convex, $Z \subset X$; and because the carrier $N(J)$ is gated, $Z \subset N(J)$. By identifying $N(J)$ with $F(J) \times C_+$ thanks to Lemma \ref{lem:carrierproduct}, we deduce that $Z$ contains a sequence of successively parallel cliques of $X$ between $C$ and $C'$. A fortiori, $C$ and $C'$ are dual to the same hyperplane of $X$.

\medskip \noindent
Thus, we have proved that $J \mapsto J \cap X$ defines a map from the hyperplanes of $Y$ intersecting $X$ to the hyperplanes of $X$. This map is clearly injective since the intersection between two distinct hyperplanes of $Y$ contains at most a single vertex. Moreover, any hyperplane of $X$ extends to a hyperplane of $Y$, so that this map is also surjective.
\end{proof}

\begin{figure}
\begin{center}
\includegraphics[scale=0.5]{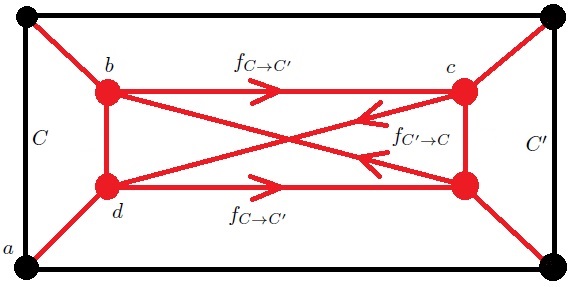}
\end{center}
\caption{Inflation of the square not satisfying the monodromy condition.}
\label{figure11}
\end{figure}

\begin{remark}
Proposition \ref{prop:inflatingonehyp} ensures that the monodromy condition is sufficient to produce a quasi-median graph. On the other hand, Figure \ref{figure11} provides an example of an inflation of the square $K_{2,2}$ which does not satisfy the monodromy condition; the graph which is produced is not quasi-median since it contains an induced subgraph isomorphic to $K_4^-$ (generated by the vertices $a,b,c,d$). In fact, if you suppose that the graph you obtain is quasi-median, then the map $f_{C \to C'}$ is precisely the restriction of the canonical bijection $C_+ \to C_+'$ to the added vertices $S(C)$. On the other hand, the family of the canonical bijections $t_{C \to C'}$ of any quasi-median graph turns out to satisfy the monodromy condition; see Remark \ref{remark:monodromy}. Therefore, the monodromy condition is necessary to produce a quasi-median graph.
\end{remark}

\subsection{Inflating more hyperplanes}\label{section:inflatingclosure}

\noindent
In this section, our aim is to show how to inflate all the hyperplanes of a quasi-median graph simultaneously.

\medskip \noindent
Let $X$ be a quasi-median graph. For every hyperplane $J$ of $X$, let $S(J)$ be an arbitrary set. Order the hyperplanes of $X$ by some ordinal $\xi$, so that $\{J_{\alpha} \mid \alpha \leq \xi \}$ is the set of the hyperplanes of $X$. Define an increasing sequence of quasi-median graphs $(X^{\alpha})$ by transfinite induction in the following way:
\begin{itemize}
	\item First, set $X^1=X$.
	\item Let $\alpha$ be an ordinal. Suppose that $X^{\gamma}$ is convex in $X^{\delta}$ for every $\gamma \leq \delta \leq \alpha$, and that the map $J \mapsto J \cap X$ defines a bijection from the hyperplanes of $X^{\alpha}$ onto those of $X$. Set $X^{\alpha+1}$ as the graph obtained from $X^{\alpha}$ by inflating the hyperplane extending $J_{\alpha}$; this is a quasi-median graph according to Proposition \ref{prop:inflatingonehyp}. Notice that, for every $\beta \leq \alpha$, $X^{\beta}$ is convex in $X^{\alpha+1}$, since $X^{\beta}$ is convex in $X^{\alpha}$ by assumption and that $X^{\alpha}$ is convex in $X^{\alpha+1}$ according to Proposition \ref{prop:inflatingonehyp}. Moreover, combining our assumption about the hyperplanes of $X^{\alpha}$ with Proposition \ref{prop:inflatingonehyp}, we know that the map $$J \mapsto J \cap X^{\alpha} \mapsto J \cap X^{\alpha} \cap X = J \cap X$$ defines a bijection from the hyperplanes of $X^{\alpha+1}$ onto those of $X$.
	\item Let $\lambda$ be a limit ordinal. Suppose that $X^{\alpha}$ is convex in $X^{\beta}$ for every $\alpha \leq \beta< \lambda$, and that, for every $\alpha< \lambda$, the map $J \mapsto J \cap X$ defines a bijection from the hyperplanes of $X^{\alpha}$ onto those of $X$. Set $X^{\lambda}= \bigcup\limits_{\alpha< \lambda} X^{\alpha}$. It is clear that, for every $\alpha< \lambda$, $X^{\alpha}$ is convex in $X^{\lambda}$, since any geodesic between two vertices of $X^{\alpha}$ must be contained in $X^{\beta}$ for some $\beta <\lambda$ sufficiently large and that $X^{\alpha}$ is convex in $X^{\gamma}$ for every $\gamma \geq \alpha$. In particular, since the property of being quasi-median can be read on finitely many vertices, we deduce that $X^{\lambda}$ is quasi-median since the $X^{\alpha}$'s are quasi-median themselves. Finally, according to Lemma \ref{lem:hypconvexsub}, and because we know that $X=X^1$ is convex in $X^{\lambda}$, it is sufficient to show that every hyperplane of $X^{\lambda}$ intersects $X$ in order to deduce that the map $J \mapsto J \cap X$ defines a bijection from the hyperplanes of $X^{\lambda}$ and those of $X$. Let $J$ be a hyperplane of $X^{\lambda}$. There exists some $\alpha < \lambda$ such that $X^{\alpha}$ contains an edge of some clique of $X^{\lambda}$ dual to $J$. In particular, $J$ intersects $X^{\alpha}$, so that $J \cap X^{\alpha}$ defines a hyperplane of $X^{\alpha}$ according to Lemma \ref{lem:hypconvexsub}. We conclude from our initial assumption about the hyperplanes of $X^{\alpha}$ that $J \cap X^{\alpha} \cap X= J \cap X$ is a hyperplane of $X$. A fortiori, $J$ intersects $X$.
\end{itemize}
The graph $X^{\xi}$ is obtained from $X$ by \emph{inflating its hyperplanes}\index{Inflating hyperplanes}. An immediate consequence of the construction is:

\begin{fact}\label{bigfact:inflating}
Let $Y$ be a graph obtained from a quasi-median graph $X$ by inflating its hyperplanes. Then $Y$ is a quasi-median graph containing $X$ as a convex subgraph, such that the map $J \mapsto J \cap X$ defines a bijection from the hyperplanes of $Y$ onto those of $X$.
\end{fact}

\noindent
Now, our goal is to show that, up to isometry, the graph we obtain does not depend on the ordering of the hyperplanes of $X$ we chose. This is a consequence of Proposition \ref{prop:inflatingclosure} below.

\begin{definition}
Let $X$ be a quasi-median graph and $Y \subset X$ a subgraph. A prism $P$ of $X$ is \emph{admissible relatively to $Y$} if, for every clique $C$ of $P$, $C \cap Y$ is a clique of $Y$. The \emph{closure}\index{Closures of subgraphs} of $Y$ in $X$ is defined as the union of all the prisms of $X$ which are admissible relatively to $Y$. 
\end{definition}

\begin{prop}\label{prop:inflatingclosure}
Let $X$ be a quasi-median graph. For every hyperplane $J$, fix an arbitrary set $S(J)$. Let $Y_1$ denote the corresponding graph we obtain by inflating the hyperplanes of $X$. For every hyperplane $J$ of $X$, let $C(J)$ denote a clique dual to $J$. If $\mathfrak{J}$ denotes the collection of all the hyperplanes of $X$, we have canonical embeddings
$$X \hookrightarrow \prod\limits_{J \in \mathfrak{J}} C(J) \hookrightarrow \prod\limits_{J \in \mathfrak{J}} C(J) \cup S(J).$$
Let $Y_2$ denote the closure of $X$ in $\prod\limits_{J \in \mathfrak{J}} C(J) \cup S(J)$. Then $Y_1$ and $Y_2$ are isometric. 
\end{prop}

\begin{proof}
The map $x \mapsto \left( \mathrm{proj}_{C(J)}^X (x): J \in \mathfrak{J} \right)$ induces an embedding $X \hookrightarrow \prod\limits_{J \in \mathfrak{J}} C(J)$ according to Lemma \ref{lem:embedinfiniteprism}. Now, for every hyperplane $J \in \mathfrak{J}$, let $J_+$ denote the unique hyperplane of $Y_1$ extending $J$ and $C(J_+)$ the unique clique of $Y_1$ containing $C(J)$; in particular, $C(J_+)=C(J) \cup S(J)$. Similarly, the map $y \mapsto \left( \mathrm{proj}_{C(J_+)}^{Y_2}(y) : J \in \mathfrak{J} \right)$ induces an embedding 
$$Y_1 \hookrightarrow \prod\limits_{J \in \mathfrak{J}} C(J_+)= \prod\limits_{J \in \mathfrak{J}} C(J) \cup S(J).$$ 
Lemma \ref{lem:projectionsurgraph} proved below precisely means that the following diagram is commutative:
\begin{displaymath}
\xymatrix{ X \ar[d] \ar[rr] & & \prod\limits_{J \in \mathfrak{J}} C(J) \ar[rr] & & \prod\limits_{J \in \mathfrak{J}} C(J) \cup S(J) \\ Y_1 \ar[urrrr] & & & & }
\end{displaymath}
From now on, we will suppose that $Y_1$ and $Y_2$ are two subgraphs of $Z=\prod\limits_{J \in \mathfrak{J}} C(J) \cup S(J)$ containing $X$. Our goal is to prove that $Y_1=Y_2$.

\medskip \noindent
Let $P=C_1 \times \cdots \times C_n$ be a prism of $Z$ which is admissible relatively to $X$, ie., $C_i \cap X$ is a clique of $X$ for every $1 \leq i \leq n$. Fix some $1 \leq i \leq n$. If $J_i$ denotes the hyperplane of $X$ dual to $C_i \cap X$, then $C_i= (C_i \cap X) \cup S(J_i)$. We know that $S(J_i) \subset Y_1$, since $S(J_i)$ are the new vertices added when inflating $J_i$, and $C_i \cap X \subset X \subset Y_1$, so $C_i \subset Y_1$. A fortiori, $P \subset Y_1$. Thus, we deduce that any prism of $Z$ which is admissible relatively to $X$ must be included into $Y_1$, whence $Y_2 \subset Y_1$.

\medskip \noindent
Conversely, let $p \in Y_1$ be a vertex. If $p \in X$, a fortiori $p \in Y_2$ since $X \subset Y_2$. Otherwise, there exists a clique $C$ of $X$ such that $p$ belongs to the clique $C_+$ of $Y_1$ obtained from $C$ after inflating the hyperplane $J$ dual to $C$. Clearly, as a prism of $Z$, $C_+= C \cup S(J)$ is admissible relatively to $X$, hence $C_+ \subset Y_2$, and finally $p \in Y_2$. This proves that $Y_1 \subset Y_2$, concluding the proof of our proposition.
\end{proof}

\begin{lemma}\label{lem:projectionsurgraph}
Let $Y$ be a graph obtained from a quasi-median graph $X$ by inflating its hyperplanes. Fix $C$ some clique of $X$ and let $C_+$ denote the unique clique of $Y$ containing $C$. If $p : X \to C$ is the projection of $X$ onto $C$ and $q : Y \to C_+$ the projection of $Y$ onto $C_+$, then $p=q_{|X}$. 
\end{lemma}

\noindent
We begin by proving the lemma when $Y$ is obtained from $X$ by inflating only one of its hyperplanes.

\begin{lemma}
Let $Y$ be a graph obtained from a quasi-median graph $X$ by inflating one of its hyperplanes. Fix $C$ some clique of $X$ and let $C_+$ denote the unique clique of $Y$ containing $C$. If $p : X \to C$ is the projection of $X$ onto $C$ and $q : Y \to C_+$ the projection of $Y$ onto $C_+$, then $p=q_{|X}$. 
\end{lemma}

\begin{proof}
Let $x \in X$ and $p \in C_+ \backslash C$ be two vertices. We want to prove that $d(x,C)<d(x,p)$. This assertion implies the conclusion of our lemma because $X$ is isometrically embedded into $Y$.

\medskip \noindent
For convenience, let $J_+$ denote the unique hyperplane of $Y$ extending $J$. Fix some geodesic $[x,p]$ between $x$ and $p$ and let $x'$ be the first vertex of $[x,p]$ which belongs to $N(J_+)$. In particular, $x'$ is adjacent to a vertex which does not belong to $N(J_+)$, since otherwise it would not be the first vertex of $[x,p]$ belonging to $N(J_+)$. By noticing that any vertex of $Y$ adjacent to some vertex of $Y \backslash X$ must belong to $N(J_+)$, we conclude that $x' \in X$. According to Claim \ref{claim:Phi}, $N(J_+)$ is naturally isomorphic to $F(J) \times C_+$. Let $F= F(J) \times \{ x'' \}$ denote the fiber containing $x'$; notice that $x'' \in X$. The structure of $N(J_+)$ as a cartesian product implies that $d(x',p)=d(x',x'')+d(x'',p)=d(x',x'')+1$. Thus,
$$d(x,p)=d(x,x')+d(x',x'')+1 \geq d(x,x'')+1 \geq d(x,C)+1,$$
concluding the proof.
\end{proof}

\begin{proof}[Proof of Lemma \ref{lem:projectionsurgraph}.]
Index the hyperplanes of $X$ by some ordinal $\xi$, say $\{ J_{\alpha} \mid \alpha \leq \xi \}$, so that, if $(X^{\alpha})$ denotes the sequence of graphs obtained from $X$ by successively inflating its hyperplanes following the order induced by $\xi$, then $Y= X^{\xi}$. Say that the hyperplane dual to $C$ is $J_{\alpha}$. It is clear that a clique remains a clique after inflating a hyperplane which is not dual to it, so $C$ is a clique of $X^{\alpha}$ and $C_+$ is the unique clique of $X^{\alpha+1}$ containing $C$. By applying the previous lemma, we know that
$$\mathrm{proj}_{C_+}^{X^{\alpha+1}} \left|_{X^{\alpha}} \right. = \mathrm{proj}_C^{X^{\alpha}}.$$
On the other hand, the following observation is clear, since any convex subgraph is isometrically embedded:

\begin{fact}
Let $B$ be a quasi-median graph and $A$ a convex subgraph. If $C$ is a clique of $B$ included into $A$, then $\mathrm{proj}_C^A = \mathrm{proj}_C^B \left|_{A} \right.$.
\end{fact}

\noindent
As a consequence, by applying this fact twice, we deduce that 
$$\mathrm{proj}_C^X = \mathrm{proj}_C^{X^{\alpha}} \left|_{X} \right. \ \text{and} \ \mathrm{proj}_{C_+}^{X^{\alpha+1}} = \mathrm{proj}_{C_+}^Y \left|_{X^{\alpha+1}} \right. .$$
Therefore, since $X=X^1$ is included into every $X^{\beta}$ for $\beta \leq \xi$,
$$\mathrm{proj}_{C_+}^Y \left|_{X} \right. = \mathrm{proj}_{C_+}^{X^{\alpha+1}} \left|_{X} \right. = \mathrm{proj}_C^{X^{\alpha}} \left|_{X} \right. = \mathrm{proj}_C^X.$$
This concludes the proof.
\end{proof}

\noindent
We conclude this section by describing the maximal prisms of a graph obtained by inflating the hyperplanes of a quasi-median graph (of finite cubical dimension): they are precisely the extensions of the maximal prisms of the initial graph. (It is worth noticing that this statement does not hold without the assumption of maximality of the prisms; see for instance Figure \ref{figure10}.)

\begin{lemma}\label{lem:inflatingprism}
Let $Y$ be a graph obtained from a quasi-median graph $X$ of finite cubical dimension by inflating its hyperplanes. A prism $P$ of $Y$ is maximal if and only if there exist $C^1, \ldots, C^m$ cliques of $X$ such that $C^1 \times \cdots \times C^m$ is a maximal prism of $X$ and $P= C^1_+ \times \cdots \times C^m_+$. 
\end{lemma}

\begin{proof}
Suppose that $P$ is a maximal prism of $Y$. Let $J_1, \ldots, J_m$ denote the hyperplanes dual to $P$; according to Proposition \ref{prop:maxprism}, this is a maximal collection of pairwise transverse hyperplanes of $Y$. As a consequence of Fact \ref{bigfact:inflating}, the hyperplanes $J_1 \cap X, \ldots, J_m \cap X$ define a maximal collection of pairwise transverse hyperplanes of $X$, so that it follows from Proposition \ref{prop:maxprism} that $\bigcap\limits_{i=1}^m N(J_i \cap X)$ defines a maximal prism $C^1 \times \cdots \times C^m$ of $X$. By applying Proposition \ref{prop:maxprism} another time, we deduce that the prism $C^1_+ \times \cdots C^m_+$ is maximal in $Y$ since the hyperplanes dual to it, namely $J_1, \ldots, J_m$, define a maximal collection of pairwise transverse hyperplanes of $Y$. On the other hand, we know that
$$C^1_+ \times \cdots \times C_+^m \subset \bigcap\limits_{i=1}^m N(J_i) = P,$$
hence $P= C_+^1 \times \cdots C^m_+$. Conversely, suppose that it is possible to write $P=C^1_+ \times \cdots \times C^m_+$ where $C^1, \ldots, C^m$ is a collection of cliques of $X$ such that $C^1 \times \cdots \times C^m$ is a maximal prism of $X$. Let $J_1, \ldots, J_m$ denote the hyperplanes of $X$ dual to $X$, and $J_1^+, \ldots, J_m^+$ their respective extensions to $Y$. As a consequence of Proposition \ref{prop:maxprism} and Fact \ref{bigfact:inflating}, we deduce that $\{J_1, \ldots, J_m\}$ is a maximal collection of pairwise transverse hyperplanes of $X$, and then that $J_1^+, \ldots, J_m^+$ define a maximal collection of pairwise transverse hyperplanes of $Y$. On the other hand, $J_1^+, \ldots, J_m^+$ are precisely the hyperplanes dual to the cliques $C_+^1, \ldots, C_+^m$, or equivalently to $P$. Therefore, Proposition \ref{prop:maxprism} implies that $P$ is a maximal prism of $Y$.
\end{proof}

\section{Topical actions on quasi-median graphs II}\label{section:topicalactionsII}

\noindent
In this section, we prove other combination theorems for groups acting topically-transitively on quasi-median graphs. We essentially follow the arguments used in Section \ref{section:topicalactionsI}, which we apply to other contexts thanks to the construction detailed in Section \ref{section:inflatingintro}.

\subsection{Inflating topical actions}

\noindent
Roughly speaking, our aim in this section is to show that, if a group $G$ acts on a quasi-median graph $X$ and if each clique-stabiliser acts on a given set, then it is possible to inflate the hyperplanes of $X$ in order to add new vertices to the cliques so that we get an action of $G$ on a new quasi-median graph $Y$ such that the action on a clique of $Y$ by its stabiliser is the same as the action of the associated clique-stabiliser of $X$ on its set. More precisely, we prove:

\begin{prop}\label{prop:inflatinghypaction}
Let $G$ be a group with a topical-transitive action on a quasi-median graph $X$. Let $\mathcal{C}$ be a collection of cliques such that any $G$-orbit of hyperplanes intersects it along a single clique, and let $\mathcal{C}= \mathcal{C}_1 \sqcup \mathcal{C}_2$ denote the associated decomposition of $\mathcal{C}$. For every clique $C\in \mathcal{C}_1$, fix a set $Y(C)$ on which $\mathrm{stab}(C)$ acts and suppose that $Y(C)$ contains a point of trivial stabiliser; for every clique $C \in \mathcal{C}_2$, fix a set $Y(C)$ satisfying $\# Y(C) \geq \# C$ in which $\mathrm{stab}(C)$ acts trivially. Then there exists a quasi-median graph $Y$, obtained from $X$ by inflating its hyperplanes, so that, if $C_+$ denotes the unique clique of $Y$ containing $C$ for every $C \in \mathcal{C}$, then 
\begin{itemize}
	\item the action $G \curvearrowright X$ extends to a $\mathcal{C}_+$-topical action $G \curvearrowright Y$;
	\item for every clique $C \in \mathcal{C}$, there exists a $\mathrm{stab}(C)$-equivariant bijection $C_+ \to Y(C)$;
	\item if $C,C'$ are two cliques of $X$ dual to the same hyperplane, with $C \in \mathcal{C}_1$, then $\rho_{C'_+}$ induces an isomorphism $\mathrm{stab}(C_+) \to \mathrm{stab}(C'_+)$.
\end{itemize}
\end{prop}

\noindent
The first step of the proof is to embed $X$ into the appropriate product $Z=\prod\limits_{J \in \mathfrak{J}} Y(J)$ and to show that the action $G \curvearrowright X$ extends to an action $G \curvearrowright Z$. This construction is of independent interest, and will be useful later, so we prove it separately.

\begin{construction}\label{construction:inflating}
Let $G$ be a group with a topical-transitive action on a quasi-median graph $X$. Let $\mathcal{C}$ be a collection of cliques such that any $G$-orbit of hyperplanes intersects it along a single clique, and let $\mathcal{C}= \mathcal{C}_1 \sqcup \mathcal{C}_2$ denote the associated decomposition of $\mathcal{C}$. For every clique $C\in \mathcal{C}_1$, fix a set $Y(C)$ on which $\mathrm{stab}(C)$ acts and suppose that $Y(C)$ contains a point of trivial stabiliser; for every clique $C \in \mathcal{C}_2$, fix a set $Y(C)$ satisfying $\# Y(C) \geq \# C$ on which $\mathrm{stab}(C)$ acts trivially. For every hyperplane $J \in \mathfrak{J}$, labelled by some $C \in \mathcal{C}$, set $Y(J)=Y(C)$. Then $X$ naturally embeds in $Z= \prod\limits_{J \in \mathfrak{J}} Y(J)$, and the action $G \curvearrowright X$ extends to an action $G \curvearrowright Z$.
\end{construction}

\noindent
For every clique $C \in \mathcal{C}_1$, fix a point $y_0(C) \in Y(C)$ whose stabiliser is trivial. For every clique $C \in \mathcal{C}_2$, fix an injection $f_C : C \hookrightarrow Y(C)$. Let $x_0 \in X$ be a base point. For every clique $C$, set $x_0(C)= \mathrm{proj}_C(x_0)$. 

\medskip \noindent
For every hyperplane $J$ of $X$, fix a clique $C(J)$ dual to it; by convention, if $J$ contains a clique of $\mathcal{C}$, we require $C(J) \in \mathcal{C}$. If $J$ is labelled by some $C \in \mathcal{C}$, set $S(J)= Y(C) \backslash f_C(C)$ if $C \in \mathcal{C}_2$ and $S(J)= Y(C) \backslash \mathrm{stab}(C) \cdot y_0(C)$ if $C \in \mathcal{C}_1$. 

\medskip \noindent
Fix some hyperplane $J \in \mathfrak{J}$ labelled by some $C \in \mathcal{C}$. For convenience, set $Y(J)=Y(C)$. Define a bijection $\phi_J : C(J) \cup S(J) \to Y(J)$ in the following way. If $C \in \mathcal{C}_2$, set $\phi_{J| C(J)}= f_C \circ \phi_{C(J)}$, where $\phi_{C(J)} : C(J) \to C=Y(J)$ is the bijection given by Definition \ref{def:mapsfortopical}, and next identify $\phi_{J | S(J)}$ with the inclusion $S(J) \subset Y(C)$. Suppose that $C \in \mathcal{C}_1$ and let $x \in C(J) \cup S(J)$. If $x \in C(J)$, then there exists a unique $g \in \mathrm{stab}(C)$ such that $\phi_{C(J)}(x)= g \cdot x_0(C)$; set $\phi_J(x)=g \cdot y_0(C)$. If $x \in S(J)$, set $\phi_J(x)= x$.

\medskip \noindent
According to Lemma \ref{lem:embedinfiniteprism}, we have natural embeddings
\begin{displaymath}
\xymatrix{ X \ar[rr] & & \prod\limits_{J \in \mathfrak{J}} C(J) \ar[rr] & & \prod\limits_{J \in \mathfrak{J}} C(J) \cup S(J) \ar[rr]_{\underset{J \in \mathfrak{J}}{\times} \phi_J} & & \prod\limits_{J \in \mathfrak{J}} Y(J)=Z}
\end{displaymath}
From now on, we see $X$ as a subgraph of $Z$. For every $g \in G$ and $(x_J)_J \in Z$, define
$$g \cdot (x_J)_J = \left( \varphi_{C(J)} ( s_{C(J)}(g) ) \cdot x_J \right)_{gJ} ;$$
recall that $s_{C}(\cdot)$ and $\varphi_C( \cdot)$ are defined respectively by Definition \ref{def:topicalps} and Definition \ref{def:mapsfortopical}. We want to prove that this defines an action $G \curvearrowright Z$ extending $G \curvearrowright X$.

\begin{claim}
For every hyperplane $J$ and every $g_1,g_2 \in G$, $$\varphi_{C(J)} \left( s_{C(J)}(g_1g_2) \right) =\varphi_{C(g_2J)} \left( s_{c(g_2J)}(g_1) \right) \cdot \varphi_{C(J)} \left( s_{C(J)} (g_2) \right).$$
\end{claim}

\begin{proof}
By applying Corollary \ref{cor:sCgh} and the point $(iii)$ of Claim \ref{claim:divers}, we find that
$$\begin{array}{lcl} \varphi_{C}(s_C(g_1g_2)) & = & \varphi_C \left( p_C(g_2)^{-1}s_{g_2C}(g_1) p_C(g_2) \right) \cdot \varphi_C( s_C(g_2)) \\ \\ & = & \varphi_{g_2C}(s_{g_2C}(g_1)) \cdot \varphi_C(s_C(g_2)) \end{array}$$
Finally, by noticing that the cliques $C(g_2J)$ and $g_2C(J)$ are dual to the same hyperplane, thanks to Lemma \ref{lem:varphisC} we conclude that
$$\begin{array}{lcl} \varphi_{C(J)} \left( s_{C(J)}(g_1g_2) \right) & = & \varphi_{g_2C(J)}(s_{g_2C(J)}(g_1)) \cdot \varphi_{C(J)} \left( s_{C(J)}(g_2) \right) \\ \\ & = & \varphi_{C(g_2J)} \left( s_{C(g_2J)}(g_1) \right) \cdot \varphi_{C(J)} \left( s_{C(J)}(g_2) \right) \end{array}$$
\end{proof}

\noindent
Now we are ready to verify that we have defined an action on $Z$. Indeed, as a consequence of the previous claim, 
$$\begin{array}{lcl} g_1 \cdot \left( g_2 \cdot (x_J)_J \right) & = & g_1 \cdot \left( \varphi_{C(J)} \left( s_{C(J)} (g_2) \right) \cdot x_J \right)_{g_2J} \\ \\ & = & g_1 \cdot \left( \varphi_{C(g_2^{-1}J)} \left( s_{C(g_2^{-1}J)}(g_2) \right) \cdot x_{g_2^{-1}J} \right)_J \\ \\ & = & \left( \varphi_{C(J)} \left( s_{C(J)}(g_1) \right) \cdot \varphi_{C(g_2^{-1}J)} \left( s_{C(g_2^{-1}J)} (g_2) \right) \cdot x_{g_2^{-1}J} \right)_{g_1J} \\ \\ & = & \left( \varphi_{C(g_2J)} \left( s_{C(g_2J)}(g_1) \right) \cdot \varphi_{C(J)}  \left( s_{C(J)} (g_2) \right) \cdot x_J \right)_{g_1g_2J} \\ \\ & = & \left( \varphi_{C(J)} \left( s_{C(J)}(g_1g_2) \right) \cdot x_J \right)_{g_1g_2J} = (g_1g_2) \cdot (x_J)_J \end{array}$$
Thus, $G$ acts on $Z$. Now, we want to prove that this action extends $G \curvearrowright X$.

\medskip \noindent
Let us make explicit the embedding $X \hookrightarrow Z$. We have
$$\left\{ \begin{array}{ccccccc} X & \to & \prod\limits_{J \in \mathfrak{J}} C(J) & \to & \prod\limits_{J \in \mathfrak{J}} C(J) \cup S(J) & \to & \prod\limits_{J \in \mathfrak{J}} Y(J) \\ x & \mapsto & \left( \mathrm{proj}_{C(J)}(x) \right)_J & \mapsto & \left( \mathrm{proj}_{C(J)}(x) \right)_J & \mapsto & \left( \phi_J \left( \mathrm{proj}_{C(J)}(x) \right) \right)_J \end{array} \right.$$
Therefore, we need to prove that for every $x \in X$
$$\phi_{gJ} \left( \mathrm{proj}_{C(gJ)}(gx) \right) = \varphi_{C(J)} \left( s_{C(J)}(g) \right) \cdot \phi_{J}( \mathrm{proj}_{C(J)}(x))$$
Let $Q \in \mathcal{C}$ label $C(J),C(gJ),gC(J)$. Suppose first that $Q \in \mathcal{C}_1$. Then, by definition, the restriction of $\phi_J : C(J) \cup S(J) \to Y(Q)$ to $C(J)$ is the concatenation $L_Q \circ \phi_{C(J)}$, where
$$L_Q : \left\{ \begin{array}{ccc} Q & \to & Y(Q) \\ h \cdot x_0(Q) & \mapsto & h \cdot y_0(Q) \end{array} \right.$$
Similarly, the restriction of $\phi_{gJ} : C(gJ) \cup S(gJ) \to Y(Q)$ to $C(gJ)$ is the concatenation $L_Q \circ \phi_{C(gJ)}$. Notice that $L_Q(hx) = h \cdot L_Q(x)$ for every $x \in Q$ and $h \in \mathrm{stab}(Q)$. Now, from
$$\begin{array}{lcl} \mathrm{proj}_{C(gJ)}(gx) & = & t_{gC(J) \to C(gJ)} \left( \mathrm{proj}_{gC(J)}(gx) \right) \\ \\ & = & t_{gC(J) \to C(gJ)} \left( g \cdot \mathrm{proj}_{C(J)}(x) \right) \end{array}$$
together with the points $(i)$ and $(ii)$ of Claim \ref{claim:divers}, we deduce that
$$\begin{array}{lcl} \phi_{gJ} \left( \mathrm{proj}_{C(gJ)} (gx) \right) & = & L_Q \circ \phi_{C(gJ)} \left( \mathrm{proj}_{C(gJ)}(gx) \right) \\ \\ & = & L_Q \circ \phi_{C(gJ)} \left( t_{gC(J) \to C(gJ)} \left( g \cdot \mathrm{proj}_{C(J)}(x) \right) \right) \\ \\ & = & L_Q \circ \phi_{gC(J)} \left( g \cdot \mathrm{proj}_{C(J)}(x) \right) \\ \\ & = & L_Q \left( \varphi_{C(J)} (s_{C(J)}(g)) \cdot \phi_{C(J)} \left( \mathrm{proj}_{C(J)}(x) \right) \right) \\ \\ & = & \varphi_{C(J)}(s_{C(J)}(g)) \cdot \phi_J \left( \mathrm{proj}_{C(J)}( x) \right) \end{array}$$
This is the equality we are looking for. Next, suppose that $Q \in \mathcal{C}_2$. Similarly, we have
$$\begin{array}{lcl} \phi_{gJ} \left( \mathrm{proj}_{C(gJ)} (gx) \right) & = & f_Q \circ \phi_{C(gJ)} \left( \mathrm{proj}_{C(gJ)}(gx) \right) \\ \\ & = &  f_Q \circ \phi_{C(gJ)} \left( t_{gC(J) \to C(gJ)} \left( g \cdot \mathrm{proj}_{C(J)}(x) \right) \right) \\ \\ & = & f_Q \circ \phi_{gC(J)} \left( g \cdot \mathrm{proj}_{C(J)}(x) \right) \\ \\ & = & f_Q \left( \varphi_{C(J)} (s_{C(J)}(g)) \cdot \phi_{C(J)} \left( \mathrm{proj}_{C(J)}(x) \right) \right) \\ \\ & = & f_Q \circ \phi_{C(J)} \left( \mathrm{proj}_{C(J)}(x) \right) \\ \\ & = & \phi_J \left( \mathrm{proj}_{C(J)}( x) \right) \end{array}$$
since $\varphi_{C(J)}$ is trivial. Thus, we have proved that the action $G \curvearrowright Z$ extends $G \curvearrowright X$.

\begin{proof}[Proof of Proposition \ref{prop:inflatinghypaction}.]
From the sets $S(J)$ of our previous construction, we deduce a quasi-median graph $Y$ obtained from $X$ by inflating its hyperplanes. For every $J \in \mathfrak{J}$, let $J_+$ denote the unique hyperplane of $Y$ extending $J$. Let $C(J_+)$ denote the unique clique of $Y$ containing $C(J)$, so $C(J_+)=C(J) \cup S(J)$. According to Lemma \ref{lem:projectionsurgraph}, we have a commutative diagram
\begin{displaymath}
\xymatrix{ X \ar[d] \ar[rr] & & \prod\limits_{J \in \mathfrak{J}} C(J) \ar[rr] & & \prod\limits_{J \in \mathfrak{J}} C(J) \cup S(J) \ar[rr]_{\underset{J \in \mathfrak{J}}{\times} \phi_J} & & \prod\limits_{J \in \mathfrak{J}} Y(J)=Z \\ Y \ar[urrrr] & & & & & & }
\end{displaymath}
Here, each $Y(J)$ is thought of as a complete graph. From now on, we see $X$ and $Y$ as subgraphs of $Z$. Previously, we defined an action $G \curvearrowright Z$ extending $G \curvearrowright X$. In particular, $X$ is $G$-invariant, so that its closure in $Z$, namely $Y$, must be $G$-invariant as well. Consequently, we have defined an action $G \curvearrowright Y$.

\medskip \noindent
Now, we want to prove the second point of our statement. Let $C \in \mathcal{C}$ be a clique and let $J$ denote its dual hyperplane. Notice that $C=C(J)$, so that $C_+=C(J) \cup S(J)$ and $Y(C)=Y(J)$. We claim that our bijection $\phi_J : C(J) \cup S(J) \to Y(J)$ is $\mathrm{stab}(C)$-equivariant. If $C \in \mathcal{C}_2$, there is nothing to prove since $\mathrm{stab}(C)$ acts trivially on both spaces. Suppose that $C \in \mathcal{C}_1$. Let $x \in C(J) \cup S(J)$ and $g \in \mathrm{stab}(C)$. If $x \in S(J)$ then $gx \in S(J)$ as well, so that $\phi_J(g \cdot x)=g \cdot x=g\cdot \phi_J(x)$. Otherwise, if $x \in C(J)$, notice that $\phi_{C(J)}$ is the identity since $C(J)=C$; therefore, if $h \in \mathrm{stab}(C)$ is the unique element satisfying $x= h \cdot x_0(C)$, then $\phi_J(g \cdot x) = \phi_J(gh \cdot x_0(C)) = gh \cdot y_0(C)= g \cdot \phi_J(x)$. This concludes the proof of our second point.

\medskip \noindent
Now, we want to prove that the action $G \curvearrowright Y$ is $\mathcal{C}_+$-topical. Let $C \in \mathcal{C}$. In particular, if $J$ denotes the hyperplane dual to $C$, then $C=C(J)$. We have
$$\left\{ \begin{array}{ccccccc} X & \to & \prod\limits_{J \in \mathfrak{J}} C(J) & \to & \prod\limits_{J \in \mathfrak{J}} C(J_+) & \to & \prod\limits_{J \in \mathfrak{J}} Y(J) \\ C & \mapsto & C(J) & \mapsto & C(J) \subset C(J_+) & \mapsto & Y(J) \end{array} \right.$$
Therefore, we need to prove that the image of $\rho_{Y(J)}$ is included into $\rho_{Y(J)}(\mathrm{stab}(Y(J)))$. Notice that, according to Fact \ref{bigfact:inflating}, there exists an equivariant bijection between the set of hyperplanes of $X$ and the set of hyperplanes of $Y$, so that $\mathrm{stab}(J_+) = \mathrm{stab}(J)$. As a consequence, if $g \in \mathrm{stab}(J_+)$, then 
$$\rho_{Y(J)}(g) = \left( x \mapsto \varphi_{C(J)} \left( s_{C(J)}(g) \right) \cdot x \right),$$
so $\rho_{Y(J)}(g)$ induces the same permutation on $Y(J)$ as $s_{C(J)}(g) \in \mathrm{stab}(C(J))= \mathrm{stab}(Y(J))$. Thus, we have verified that our action is indeed $\mathcal{C}_+$-topical.

\medskip \noindent
Let us conclude the proof of our proposition by showing its third point. Notice that, if $J$ denotes the hyperplane dual to $C$, then $C=C(J)$. Since we know that the action $G \curvearrowright Y$ is $\mathcal{C}_+$-topical, the image of $\rho_{C(J_+)}$ must be included into $\mathrm{stab}(C(J_+)) \subset \mathrm{Bij}(C(J_+))$. In fact, since $C(J_+)=C(J) \cup S(J)$ and that $C(J)$ is stable under the action of $\mathrm{stab}(C(J_+)) = \mathrm{stab}(C(J))$ on $C(J)$, the image of $\rho_{C(J_+)}$ must be included into the subgroup $\mathrm{Bij}_{C(J)}(C(J_+))$ of $\mathrm{Bij}(C(J_+))$ leaving $C(J)$ invariant. Moreover, since $\mathrm{stab}(C(J))$ acts freely on $C(J)$, the natural surjection $\pi : \mathrm{Bij}_{C(J)}(C(J_+)) \twoheadrightarrow \mathrm{Bij}(C(J))$ is injective on $\mathrm{stab}(C(J))$. As a consequence, if we are able to prove that $\rho_{C(J_+)} \circ \pi$ turns out to induce an isomorphism $\mathrm{stab}(C') \to \mathrm{stab}(C(J)) \subset \mathrm{Bij}(C(J))$, then it will follow that that $\rho_{C(J_+)}$ induces an isomorphism $\mathrm{stab}(C') \to \mathrm{stab}(C(J)) \subset \mathrm{Bij}(C(J_+))$. On the other hand, since $G \curvearrowright Y$ extends the action $G \curvearrowright X$, and noticing that $\mathrm{stab}(J_+)= \mathrm{stab}(J)$, we deduce that $\rho_{C(J_+)} \circ \pi = \rho_{C(J)}$. According to Lemma \ref{lem:rho3}, we know that $\rho_{C(J)}$ is an isomorphism, so the conclusion follows.
\end{proof}

\subsection{Cubulating groups acting on quasi-median graphs II}\label{section:cubulatingII}

\noindent
This section completes Section \ref{section:cubulatingI} by showing how to create properly discontinous and virtually cocompact special actions on CAT(0) cube complexes from a topical-transitive action on a quasi-median graph.

\begin{prop}\label{prop:properlydiscontinuous}
Let $G$ be a group acting topically-transitively on a quasi-median graph $X$ with finite vertex-stabilisers. Let $\mathcal{C}$ be a collection of cliques such that any $G$-orbit of hyperplanes intersects it along a single clique, and let $\mathcal{C}= \mathcal{C}_1 \sqcup \mathcal{C}_2$ denote the associated decomposition of $\mathcal{C}$. Suppose that every clique of $\mathcal{C}_2$ has cardinality two. If clique-stabilisers act properly discontinuously on a CAT(0) cube complex, then so does $G$.
\end{prop}

\begin{proof}
For every $C \in \mathcal{C}_1$, let $Y(C)$ be a CAT(0) cube complex on which $\mathrm{stab}(C)$ acts properly discontinuously; according to Lemma \ref{lem:modifycubing}, we may suppose without loss of generality that $Y(C)$ contains a vertex with trivial stabiliser. For every $C \in \mathcal{C}_2$, set $Y(C)=C$. Let $Y$ be the quasi-median graph provided by Proposition \ref{prop:inflatinghypaction} on which $G$ acts $\mathcal{C}_+$-topically. 

\medskip \noindent
Let $C \in \mathcal{C}$. If $C \in \mathcal{C}_2$, then $C_+=C$ is an edge $(x,y)$, and we set $\mathcal{W}(C_+)= \{ \{x\}, \{y\} \}$. Otherwise, if $C \in \mathcal{C}_1$, according to Proposition \ref{prop:inflatinghypaction}, there exists a $\mathrm{stab}(C)$-equivariant bijection $\varphi_C : C_+ \to Y(C)$, and we set $\mathcal{W}(C_+)= \varphi_C^{-1} \mathcal{W}(Y(C))$. According to Proposition \ref{prop:strongsystwall}, our collection $\{ (C_+, \mathcal{W}(C_+)) \mid C \in \mathcal{C} \}$ extends to a $G$-invariant coherent system of wallspaces. We claim that the induced action $G \curvearrowright C(Y, \mathcal{HW})$ is properly discontinuous, where $C(Y, \mathcal{HW})$ is the CAT(0) cube complex obtained by cubulating $(Y, \mathcal{HW})$. 

\medskip \noindent
So we need to verify that the stabiliser of any orientation of $\mathcal{HW}$ is finite. According to Lemma \ref{lem:dichotomyorientations}, such an orientation is either principal or semiprincipal. On the other hand, it is clear that, for every $C \in \mathcal{C}$, any orientation of $\mathcal{W}(C_+)$ must be principal, and this statement extends to all the cliques of $Y$ as a consequence of Proposition \ref{prop:strongsystwall}. Therefore, the orientations of $\mathcal{HW}$ are principal, so that it is sufficient to show that $G$ acts on $Y$ with finite vertex-stabiliser. 

\medskip \noindent
Let $y \in Y$ be a vertex. If $y \in X$, then $\mathrm{stab}(y)$ is finite by assumption. Otherwise, because $Y$ is obtained from $X$ by inflating its hyperplanes, there exists a unique clique $C$ of $X$ such that $y \in C_+$, where $C_+$ denotes the unique clique of $Y$ containing $C$. As a consequence, $\mathrm{stab}(y) \subset \mathrm{stab}(C_+)$. If $C$ is labelled by $\mathcal{C}_2$, then $\mathrm{stab}(C_+)= \mathrm{stab}(C)=\mathrm{fix}(C)$ is finite since vertex-stabilisers of $X$ are finite, so that $\mathrm{stab}(y)$ must be finite as well. Next, suppose that $C$ is labelled by $\mathcal{C}_1$. According to Proposition \ref{prop:strongsystwall}, the action $\mathrm{stab}(C_+) \curvearrowright C_+$ is isomorphic to the action $\mathrm{stab}(Q) \curvearrowright Q_+$ for some $Q \in \mathcal{C}_1$ (the morphism $\psi_C$ in the statement of Proposition \ref{prop:strongsystwall} is an isomorphism because the morphism $\rho_{Q \to C'}$ to which it is conjugated turns out to be an isomorphism according to Proposition \ref{prop:inflatinghypaction}); on the other hand, the action $\mathrm{stab}(Q) \curvearrowright Q_+$ is isomorphic to the action $\mathrm{stab}(Q) \curvearrowright Y(Q)$ according to Proposition \ref{prop:inflatinghypaction}. A fortiori, vertex-stabilisers of $C_+$ with respect to the action $\mathrm{stab}(C_+) \curvearrowright C_+$ are finite, so that $\mathrm{stab}(y)$ must be finite.
\end{proof}

\begin{prop}\label{prop:cocompactspecial}
Let $G$ be a group acting topically-transitively on a quasi-median graph $X$. Suppose that
\begin{itemize}
	\item any vertex of $X$ belongs to finitely many cliques;
	\item any vertex-stabiliser is finite;
	\item the cubical dimension of $X$ is finite;
	\item $X$ contains finitely many $G$-orbits of prisms and hyperplanes;
	\item for every maximal prism $P=C_1 \times \cdots C_n$, $\mathrm{stab}(P)= \mathrm{stab}(C_1) \times \cdots \times \mathrm{stab}(C_n)$;
	\item the action $G \curvearrowright X$ is special;
	\item for every hyperplane $J$, $\mathrm{stab}(J)$ is a retract of $G$;
\end{itemize}
If clique-stabilisers act geometrically and virtually specially on CAT(0) cube complexes, then so does $G$.
\end{prop}

\begin{proof}
Let $\mathcal{C}$ be a collection of cliques such that any $G$-orbit of hyperplanes intersects it along a single clique, and let $\mathcal{C}= \mathcal{C}_1 \sqcup \mathcal{C}_2$ denote the associated decomposition of $\mathcal{C}$. For every $C \in \mathcal{C}_1$, let $Y(C)$ be a CAT(0) cube complex on which $\mathrm{stab}(C)$ acts geometrically and virtually specially; according to Lemma \ref{lem:modifycubing}, we may suppose without loss of generality that $Y(C)$ contains a vertex with trivial stabiliser. For every $C \in \mathcal{C}_2$, set $Y(C)=C$. According to Proposition \ref{prop:inflatinghypaction}, there exists a quasi-median graph $Y$ obtained from $X$ by inflating its hyperplanes so that the action $G \curvearrowright X$ extends to a $\mathcal{C}_+$-topical action $G \curvearrowright Y$. Moreover, $C_+=C$ for every $C \in \mathcal{C}_2$, and, for every $C \in \mathcal{C}_1$, there exists a $\mathrm{stab}(C)$-equivariant bijection $\varphi_C : C_+ \to Y(C)$. 

\medskip \noindent
If $C \in \mathcal{C}_2$, set $\mathcal{W}(C_+)= \{ \{ \{x \} \cup \{x \}^c \} \mid x \in C_+ \}$. If $C \in \mathcal{C}_1$, set $\mathcal{W}(C_+)= \varphi_C^{-1} \mathcal{W}(Y(C))$. According to Proposition \ref{prop:strongsystwall}, our collection $\{ (C_+, \mathcal{W}(C_+)) \mid C \in \mathcal{C} \}$ extends to a coherent $G$-invariant system of wallspaces; moreover, given any clique $C$ of $Y$, there exists some $Q \in \mathcal{C}$ such that the spaces with walls $(C, \mathcal{W}(C))$ and $(Y(Q),\mathcal{W}(Y(Q)))$ are isomorphic. We also deduce from the combination of the third point of Proposition \ref{prop:inflatinghypaction} with the last assertion of Proposition \ref{prop:strongsystwall}, that, for every clique $C$ of $X$ labelled by some $Q \in \mathcal{C}_1$, the actions $\mathrm{stab}(C_+) \curvearrowright (C_+,\mathcal{W}(C_+))$ and $\mathrm{stab}(Q) \curvearrowright Y(Q)$ are isomorphic. 

\medskip \noindent
Let $Z$ denote the CAT(0) cube complex obtained by cubulating the space with walls $(X, \mathcal{HW})$. According to Proposition \ref{prop:producingproperaction}, the action $G \curvearrowright Z$ is metrically proper. As a consequence of Lemma \ref{lem:inflatingprism}, if $P_+$ is a maximal prism of $Y$, then $P_+=C^1_+ \times \cdots \times C^n_+$ for some prism $P=C^1 \times \cdots \times C^n$ of $X$. In particular, 
$$\mathrm{stab}(P_+)= \mathrm{stab}(P)= \mathrm{stab}(C^1) \times \cdots \times \mathrm{stab}(C^n) = \mathrm{stab}(C^1_+) \times \cdots \times \mathrm{stab}(C^n_+).$$
On the other hand, we know that the actions $\mathrm{stab}(C^i_+) \curvearrowright (C^i_+, \mathcal{W}(C^i_+))$ are cocompact, for every $1 \leq i \leq n$, so the action $\mathrm{stab}(P_+) \curvearrowright (P_+, \mathcal{W}(P_+))$ must be cocompact as well. By applying Proposition \ref{prop:producingcocompactaction}, we deduce that $G$ acts cocompactly on $Z$. 

\medskip \noindent
Finally, fix $J$ a hyperplane of $Y$, $C$ a clique dual to $J$ and $g \in \mathrm{stab}(J)$. In particular, $\mathcal{W}(J)= \overline{\mathcal{W}}(C)$. Notice that:
\begin{itemize}
	\item if there exists some $W \in \mathcal{W}(C)$ such that $\overline{W}$ and $g \cdot \overline{W}$ are transverse, then $W$ and $\rho_C(g) \cdot W$ are transverse;
	\item if there exists some $W \in \mathcal{W}(C)$ such that $\overline{W}$ and $g \cdot \overline{W}$ are tangent, then $W$ and $\rho_C(g) \cdot W$ are tangent;
	\item if there exist $W_1,W_2 \in \mathcal{W}(C)$ such that $\overline{W}_1$ and $\overline{W}_2$ are transverse and $g \cdot \overline{W}_1$ and $\overline{W}_2$ are tangent, then $W_1$ and $W_2$ are transverse and $\rho_C(g) \cdot W_1$ and $W_2$ are tangent.
\end{itemize}
Therefore, if $H$ denotes a finite-index subgroup of $\mathrm{stab}(C)$ such that the induced action $H \curvearrowright (C, \mathcal{W}(C))$ is special, then $\rho_C^{-1}(H)$ defines a finite-index subgroup of $\mathrm{stab}(J)$ acting specially on $(N(J), \mathcal{W}(J))$. A fortiori, the action $\mathrm{stab}(J) \curvearrowright (N(J), \mathcal{W}(J))$ is virtually special. It follows from Proposition \ref{prop:producingspecialaction} that $G$ acts virtually specially on $Z$.
\end{proof}

\begin{remark}
In the statements of Proposition \ref{prop:properlydiscontinuous} and Proposition \ref{prop:cocompactspecial}, the dimension of the CAT(0) cube complex which is produced is bounded above by 
$$\dim_{\square}X \cdot \sup \{ d(C) \mid C \in \mathcal{C}_1 \},$$
where $d(C)$ denote the CAT(0) cube complex on which $\mathrm{stab}(C)$ acts. This follows from Corollary \ref{cor:dimHW}, Proposition \ref{prop:strongsystwall}, and from the observation that inflating the hyperplanes of a quasi-median graph does not modify its cubical dimension. 
\end{remark}

\subsection{CAT(0) groups acting on quasi-median graphs}

\noindent
This section is dedicated to the proof of the following combination theorem, which combines Proposition \ref{prop:inflatinghypaction} with Proposition \ref{prop:CAT0delta} in order to produce \emph{CAT(0) groups}. Following \cite{MR1744486}, a \emph{CAT(0) group}\index{CAT(0) groups} is defined as a group acting metrically properly and cocompactly on some CAT(0) space. Notice that such a CAT(0) space is necessarily proper.

\begin{thm}\label{thm:producingCAT0groups}
Let $G$ be a group acting topically-transitively on a quasi-median graph $X$. Suppose that
\begin{itemize}
	\item any vertex of $X$ belongs to finitely many cliques;
	\item any vertex-stabiliser is finite;
	\item the cubical dimension of $X$ is finite;
	\item $X$ contains finitely many $G$-orbits of prisms;
	\item for every maximal prism $P=C_1 \times \cdots C_n$, $\mathrm{stab}(P)= \mathrm{stab}(C_1) \times \cdots \times \mathrm{stab}(C_n)$.
\end{itemize}
If clique-stabilisers are CAT(0), then so is $G$.
\end{thm}

\noindent
This theorem will be an easy consequence of the following general statement together with Proposition \ref{prop:CAT0delta}.

\begin{prop}\label{prop:producinggeometricactions}
Let $G$ be a group acting topically-transitively on a quasi-median graph $X$. Let $\mathcal{C}$ be a collection of cliques such that any $G$-orbit of hyperplanes intersects it along a single clique, and let $\mathcal{C}= \mathcal{C}_1 \sqcup \mathcal{C}_2$ denote the associated decomposition of $\mathcal{C}$. Suppose that
\begin{itemize}
	\item any vertex of $X$ belongs to finitely many cliques;
	\item any vertex-stabiliser is finite;
	\item the cubical dimension of $X$ is finite;
	\item $X$ contains finitely many $G$-orbits of prisms;
	\item for every maximal prism $P=C_1 \times \cdots C_n$, $\mathrm{stab}(P)= \mathrm{stab}(C_1) \times \cdots \times \mathrm{stab}(C_n)$;
	\item for every $C \in \mathcal{C}_1$, $\mathrm{stab}(C)$ acts geometrically on a proper metric space $Y(C)$ containing a point $x_0(C)$ whose stabiliser is trivial; and for every $C \in \mathcal{C}_2$, $\mathrm{stab}(C)$ acts trivially on a compact metric space $Y(C)$ satisfying $\# Y(C) \geq \# C^{(0)}$.
\end{itemize}
There exists a quasi-median graph $Y$ endowed with a coherent system of metrics such that $G$ acts metrically properly and cocompactly on $(Y, \delta^p)$ for every $p \in [1,+ \infty]$ and such that, for every clique $Q$ of $Y$, $(Q,\delta_Q)$ is isometric to $Y(C)$ for some $C \in \mathcal{C}$. 
\end{prop}

\begin{proof}
Let $G$ act on the quasi-median graph $Y$ provided by Proposition \ref{prop:inflatinghypaction}.

\medskip \noindent
For every clique $C \in \mathcal{C}$, we use the bijection $C_+ \to Y(C)$ provided by Proposition \ref{prop:inflatinghypaction} in order to transfer the distance of $Y(C)$ to $C_+$. According to Proposition \ref{prop:strongsystdist}, our collection of metric spaces $\{(C_+,\delta_{C_+}) \mid C \in \mathcal{C} \}$ extends to a $G$-invariant coherent system of metrics $\{(C, \delta_C) \mid C \ \text{clique} \}$ of $Y$. Moreover, we know from Proposition \ref{prop:strongsystdist}, and the combination of the third point of Proposition \ref{prop:inflatinghypaction} with the last assertion of Proposition \ref{prop:strongsystdist}, that each action $\mathrm{stab}(C_+) \curvearrowright (C_+,\delta_{C_+})$, where $C$ is labelled by $\mathcal{C}_1$, must be isomorphic to an action $\mathrm{stab}(Q) \curvearrowright Y(Q)$ for some $Q \in \mathcal{C}_1$. On the other hand, for every clique $C$ labelled by $\mathcal{C}_2$, we know that $\mathrm{stab}(C_+)= \mathrm{stab}(C) = \mathrm{fix}(C)$ is finite since vertex-stabilisers of $X$ are finite, and it follows from Proposition \ref{prop:strongsystdist} that $(C_+,\delta_{C_+})$ is compact. As a consequence,

\begin{fact}
For every clique $C$ of $X$, the action $\mathrm{stab}(C_+) \curvearrowright (C_+,\delta_{C_+})$ is metrically proper and cocompact.
\end{fact}

\noindent
We claim that the action $G \curvearrowright (Y,\delta^p)$ is cocompact. 

\medskip \noindent
First, notice that this action is cobounded. Indeed, we know that $X$ contains finitely many orbits of prisms, so we deduce from Lemma \ref{lem:inflatingprism} that $Y$ also contains finitely many orbits of prisms. Let $P_1, \ldots, P_n$ be a set of representatives for the action of $G$ on the maximal prisms of $Y$. Fixing some $1 \leq i \leq n$, the prism $P_i$ decomposes as $C^1_+ \times \cdots \times C^{k(i)}_+$ for some cliques $C^1, \ldots, C^{k(i)}$ of $X$. Moreover, 
$$\mathrm{stab}(P_i)= \mathrm{stab}(P_i \cap X)= \mathrm{stab}(C^1) \times \cdots \times \mathrm{stab}(C^{k(i)})= \mathrm{stab}(C^1_+) \times \cdots \times \mathrm{stab}(C^{k(i)}_+).$$
Because the actions $\mathrm{stab}(C^j_+) \curvearrowright (C^j_+,\delta_{C^j_+})$ are cocompact, we deduce that the global action $\mathrm{stab}(P_i) \curvearrowright (P_i, \delta^2)$ is cocompact as well. Therefore, if $F_i \subset P_i$ denotes a compact fundamental domain for this action, then $\bigcup\limits_{i=1}^n F_i$ defines a bounded fundamental domain for the action $G \curvearrowright (Y,\delta^p)$. 

\medskip \noindent
Finally, we want to apply Lemma \ref{lem:proper} in order to prove that $(Y,\delta^p)$ is proper, which will conclude the proof of our claim. Set
$$K_1 = \inf \{ d_{Y(C)}(x_0(C),g \cdot x_0(C)) \mid C \in \mathcal{C}_1, \ g \in \mathrm{stab}(C) \backslash \{ 1 \} \}$$
and
$$K_2 = \inf \{ \delta_{C_+} (x,y) \mid x \neq y \in C, \ C \in \mathcal{C}_2 \},$$
and finally $K= \min(K_1,K_2)$. 
Notice that, because the actions $\mathrm{stab}(C) \curvearrowright Y(C)$ are metrically proper and that the $x_0(C)$'s have trivial stabilisers, necessarily $K_1>0$. Moreover, since every clique of $\mathcal{C}_2$ is finite, necessarily $K_2 >0$. Therefore, $K>0$. 

\medskip \noindent
Let $x,y \in X \subset Y$ be two distinct vertices belonging to a common clique; let us write this clique as $C_+$ for some clique $C$ of $X$. Let $Q \in \mathcal{C}$ and $\Phi_C : Q_+ \to C_+$ be the clique and the bijection provided by Proposition \ref{prop:topicalI0} (which are used in Proposition \ref{prop:strongsystdist} in order to produce our system of metrics, so that $\Phi_C$ induces an isometry $(Q_+,\delta_{Q_+}) \to (C_+,\delta_{C_+})$), and $\Psi_Q : Q_+ \to Y(Q)$ the bijection provided by Proposition \ref{prop:inflatinghypaction} (which we used at the begining of our proof in order to define $\delta_{Q_+}$, so that $\Psi_Q$ induces an isometry $(Q_+,\delta_{Q_+}) \to (Y(Q),d_{Y(Q)})$ by construction). Because $x,y \in X$, necessarily $\Phi_C^{-1}(x), \Phi_C^{-1}(y) \in X$ since $X$ is convex into $Y$. Now, two cases may happen. First, suppose that $Q \in \mathcal{C}_1$. We deduce from the definition of $\Psi_C$ that there must exist some $g,h \in \mathrm{stab}(Q)$ such that $\Psi_C \circ \Phi_C^{-1}(x)=g \cdot x_0(Q)$ and $\Psi_C \circ \Phi_C^{-1}(y)=h \cdot x_0(Q)$. Notice that, because $x \neq y$ and because the stabiliser of $x_0(Q)$ is trivial, necessarily $g \neq h$. By construction, we deduce
$$\begin{array}{lcl} \delta_{C_+}(x,y) & = & \delta_{Q_+}( \Phi_C^{-1}(x), \Phi_C^{-1}(y)) = d_{Y(Q)}( \Psi_C \circ \Phi_C^{-1}(x), \Psi_C \circ \Phi_C^{-1}(y)) \\ \\ & = & d_{Y(Q)}(g \cdot x_0(Q),h \cdot x_0(Q)) \geq K_1 \end{array}$$
Next, suppose that $Q \in \mathcal{C}_2$. Then
$$\delta_{C_+}(x,y)= \delta_{Q_+} ( \Phi_C^{-1}(x), \Phi_C^{-1}(y)) \geq K_2.$$
Thus, we have proved that $\delta_{C_+}(x,y) \geq K$. 

\medskip \noindent
Consequently, Lemma \ref{lem:proper} applies, so that $(Y, \delta^p)$ is proper.

\medskip \noindent
We claim that the action $G \curvearrowright (Y, \delta^p)$ is metrically proper. Let $x \in X \subset Y$ denote the base point used in Construction \ref{construction:inflating} and fix some constant $R \geq 0$. We want to prove that the set
$$F = \{ g \in G \mid \delta^p(x,g \cdot x) \leq R \}$$
is finite. Notice that, according to Claim \ref{claim:Xdiscrete}, the ball $B(x,R)$ (with respect to $\delta^p$) intersects finitely many cliques of $X$, say $C_1, \cdots, C_m \subset X$. For every $1 \leq i \leq m$, let $x_i \in X$ denote the projection of $x$ onto $C_i$. Without loss of generality, we may suppose that $x \in C_1$, so that $x_1=x$. Fixing some $1 \leq i \leq m$, set $p_i=p_{C_1}(g)$ if there exists some $g \in G$ such that $gC_1=C_i$ (notice that $p_{C_1}(\cdot)$ is defined with respect to the action $G \curvearrowright X$); otherwise, set $p_i=1$. It is worth noticing that there are only finitely many possibilities for the choice of $p_i$, because

\begin{fact}
If $g,h \in G$ satisfy $gC_1=C_i$ and $hC_1=C_i$, then either $p_{C_1}(h)=p_{C_1}(g)$ if $C_1$ is labelled by $\mathcal{C}_1$, or $p_{C_1}(h)\in p_{C_1}(g) \cdot \mathrm{stab}(C_1)$ otherwise. 
\end{fact}

\noindent
Indeed, if $C_1$ is labelled by $\mathcal{C}_2$, then $p_{C_1}(h)^{-1}p_{C_1}(g)=h^{-1}g \in \mathrm{stab}(C_1)$; and if $C_1$ is labelled by $\mathcal{C}_1$, the equality follows from 
$$p_{C_1}(h)^{-1}p_{C_1}(g) \cdot x_1= p_{C_1}(h)^{-1} \cdot x_i = x_1,$$
because $\mathrm{stab}(C_1)$ acts freely on the vertices of $C_1$. This proves the fact. For convenience, set $P_i= \{ p_i \}$ if $C_1$ is labelled by $\mathcal{C}_1$ and $P_i= \mathrm{stab}(C_1)$ otherwise. Notice that $P_i$ is necessarily finite. 

\medskip \noindent
Next, set 
$$S= \{ s \in \mathrm{stab}(C_1) \mid \delta_{C_1}(x,s x) \leq R \cdot \dim_{\square}(X) \}.$$ 
Now, let $g \in F$. First, there exists some $1 \leq i \leq m$ such that $g C_1=C_i$, so $p_{C_1}(g) \in P_i$. In the sequel, we will denote $\delta_{C_i}$ for the restriction of $\delta_{(C_i)_+}$ to $C_i$. Next, as a consequence of Proposition \ref{prop:comparedelta} and the observation that $x_i,gx$ belong to the same clique, we get
$$\begin{array}{lcl} \delta_{C_1}(x,s_{C_1}(g)x) & = & \delta_{C_i}(p_{C_1}(g) x, p_{C_1}(g)s_{C_1}(g)x)= \delta_{C_i}(x_i,gx) \\ \\ & \leq & \delta^1(x,x_i)+ \delta_{C_i}(x_i,gx) = \delta^1(x,gx) \\ \\ & \leq & \dim_{\square}(X) \cdot \delta^p(x,gx) \leq \dim_{\square}(X) \cdot R \end{array}$$
Therefore, $g= p_{C_1}(g)s_{C_1}(g) \in P_iS$. We conclude that $F$ is finite since $F \subset \bigcup\limits_{i=1}^m P_iS$, where $S$ is finite because the action $\mathrm{stab}(C_1) \curvearrowright (C_1, \delta_{C_1})$ is metrically proper.

\medskip \noindent
Thus, we have proved that $G$ acts properly discontinuously and cocompactly on $(X,\delta^p)$. 
\end{proof}

\noindent
Now, we are ready to prove Theorem \ref{thm:producingCAT0groups}.

\begin{proof}[Proof of Theorem \ref{thm:producingCAT0groups}.]
Let $\mathcal{C}$ be a collection of cliques such that any $G$-orbit of hyperplanes intersects it along a single clique, and let $\mathcal{C}= \mathcal{C}_1 \sqcup \mathcal{C}_2$ denote the associated decomposition of $\mathcal{C}$. For every $C \in \mathcal{C}_1$, let $Y(C)$ be a proper CAT(0) space on which $\mathrm{stab}(C)$ acts geometrically; we may suppose without loss of generality that $Y(C)$ contains a vertex $x_0(C)$ with trivial stabiliser (follow the construction used in the proof of Lemma \ref{lem:modifycubing}). For every $C \in \mathcal{C}_2$, let $Y(C)$ be any infinite compact CAT(0) space, say the segment $[0,1]$. Now, let $(Y,\delta^2)$ be the space provided by Proposition \ref{prop:producinggeometricactions} on which $G$ acts metrically properly and cocompactly. According to Proposition \ref{prop:CAT0delta}, $(Y,\delta^2)$ is a CAT(0) metric space. 

\medskip \noindent
Thus, we have proved that $G$ acts metrically properly and cocompactly on some CAT(0) space, ie., $G$ is a CAT(0) group.
\end{proof}

\section{Application to graph products}\label{section:appli1}

\noindent
In this introduction, we fix the main definitions and notation about graph products which we will use in the sequel.

\medskip \noindent
Given a \emph{simplicial graph} $\Gamma$ (ie., a graph containing neither multiple edges nor loops) and a collection of groups $\mathcal{G}=\{ G_u \mid u \in V(\Gamma) \}$ indexed by the vertices of $\Gamma$, we define the \emph{graph product}\index{Graph products} $\Gamma \mathcal{G}$ as the quotient $$\left( \underset{u \in V(\Gamma)}{\ast} G_u \right) / \langle \langle [g,h]=1, g \in G_u, h \in G_v \ \text{if} \ (u,v) \in E(\Gamma) \rangle \rangle.$$
Every group of $\mathcal{G}$, called a \emph{vertex-group}, naturally embeds into $\Gamma \mathcal{G}$. For convenience, we will identify each vertex-group with its image into the graph product.

\noindent
A \emph{word} in $\Gamma \mathcal{G}$ is a product $g_1 \cdots g_n$ for some $n \geq 0$ and, for every $1 \leq i \leq n$, $g_i \in G$ for some $G \in \mathcal{G}$; the $g_i$'s are the \emph{syllables} of the word, and $n$ is the \emph{length} of the word. Clearly, the following operations on a word does not modify the element of $\Gamma \mathcal{G}$ it represents:
\begin{description}
	\item[Cancellation:] delete the syllable $g_i=1$;
	\item[Amalgamation:] if $g_i,g_{i+1} \in G$ for some $G \in \mathcal{G}$, replace the two syllables $g_i$ and $g_{i+1}$ by the single syllable $g_ig_{i+1} \in G$;
	\item[Shuffling:] if $g_i$ and $g_{i+1}$ belong to two adjacent vertex-groups, switch them.
\end{description}
A word is \emph{reduced} if its length cannot be shortened by applying these elementary moves, and \emph{semi-reduced} if it can be reduced without any cancellation. In practice, if $g=g_1 \cdots g_n$ is a reduced word and $h$ is a syllable, then a reduction of the product $gh$ is given by
\begin{center}
$\left\{ \begin{array}{ll} g_1 \cdots g_n & \text{if}~h=1 \\ g_1 \cdots g_{i-1} \cdot g_{i+1} \cdots g_n & \text{if}~g_i~ \text{shuffles to the end and}~g_i=h^{-1} \\ g_1 \cdots g_{i-1} \cdot g_{i+1} \cdots g_n \cdot (g_ih) & \text{if}~g_i~ \text{shuffles to the end and}~g_i \neq h^{-1} \end{array} \right.$
\end{center}
In particular, every element of $\Gamma \mathcal{G}$ can be represented by a reduced word, and this word is unique up to the shuffling operation; in fact, we will often identify an element of $\Gamma \mathcal{G}$ with any of the reduced words representing it. This allows us to define the \emph{length} of an element $g \in \Gamma \mathcal{G}$, denoted by $|g|$, as the length of any reduced word representing $g$; and its \emph{support}, denoted by $\supp(g)$, as the set of vertices of $\Gamma$ which corresponds exactly to the vertex-groups containing the syllables of $g$. For more information, we refer to \cite{GreenGP} (see also \cite{HsuWise}). The following definition will also be useful:

\begin{definition}\label{def:headtail}
Let $\Gamma$ be a simplicial graph, $\mathcal{G}$ a collection of groups labelled by $V(\Gamma)$, and $g \in \Gamma \mathcal{G}$ an element. The \emph{head} of $g$, denoted by $\mathrm{head}(g)$, is the set of syllables of some reduced word representing $g$ which appear as the first syllable of some reduced word representing $g$. Similarly, the \emph{tail} of $g$, denoted by $\mathrm{tail}(g)$, is the set of syllables of some reduced word representing $g$ which appear as the last syllable of some reduced word representing $g$.
\end{definition}

\noindent
It is worth noticing that, since any two reduced words representing a given element differ only by a sequence of shufflings, the previous sets of syllables does not depend on the the reduced word we consider.

\medskip \noindent
Finally, we give some definitions on simplicial graphs. Given a vertex $v$, 
\begin{itemize}
	\item the \emph{link} of $v$, denoted by $\link(v)$, is the subgraph generated by the vertices which are adjacent to $v$;
	\item the \emph{star} of $v$, denoted by $\st(v)$, is the subgraph generated by $v$ and $\link(v)$.
\end{itemize}
Given a graph $\Gamma$, a subgraph $\Lambda$ is \emph{induced} if two vertices of $\Lambda$ are adjacent in $\Lambda$ if and only if they are adjacent in $\Gamma$. Subgraphs of particular interest will be \emph{join subgraphs}: $\Lambda \leq \Gamma$ is a \emph{join subgraph} if it contains itself two disjoint subgraphs $\Lambda_1, \Lambda_2$ such that every vertex of $\Lambda$ lies in either $\Lambda_1$ or $\Lambda_2$, and every vertex of $\Lambda_1$ is adjacent to any vertex of $\Lambda_2$.

\subsection{Cubical-like geometry}\label{section:GPgeometry}

\noindent
In this section, we associate to any graph product a quasi-median graph, and we study its geometry.

\medskip \noindent
Let $\Gamma$ be a simplicial graph and $\mathcal{G}$ be a collection of non trivial groups labelled by $V(\Gamma)$. We denote by $\X$ the Cayley graph of the graph product $\Gamma \mathcal{G}$ with respect to the generating set $\bigsqcup\limits_{G \in \mathcal{G}} G \backslash \{ 1 \}$. Explicitely, $\X$ is the graph whose vertices are the elements of $\Gamma \mathcal{G}$ and whose edges link two elements $g,h \in \Gamma \mathcal{G}$ whenever there exists some $s \in \bigsqcup\limits_{G \in \mathcal{G}} G \backslash \{ 1 \}$ such that $g=hs$ or $h=gs$. In the sequel, we will often identify $\Gamma \mathcal{G}$ with the set of vertices of $\X$.

\begin{prop}\label{prop:Xquasimedian}
$\X$ is a quasi-median graph.
\end{prop}

\noindent
We first need to understand the geodesics of $\X$.

\begin{lemma}\label{lem:geodesic}
Let $g,h \in \Gamma \mathcal{G}$. Write $g^{-1}h$ as a reduced word $u_1 \cdots u_n$. Then the sequence of vertices $$g,gu_1,gu_1u_2, \ldots, gu_1 \cdots u_n=h$$ defines a geodesic between $g$ and $h$ in $X(\Gamma, \mathcal{G})$. Conversely, any geodesic between $g$ and $h$ is labelled by a reduced word representing $g^{-1}h$.
\end{lemma}

\begin{proof}
Any path between $g$ and $h$ produces a sequence of vertices
$$g,gu_1,gu_1u_2 \ldots, gu_1 \cdots u_n=h,$$
where $u_1, \ldots, u_n$ are syllables such that $u_1 \cdots u_n=g^{-1}h$. The elements $g$ and $h$ being fixed, such a path has minimal length if and only if $u_1 \cdots u_n$ is a word of minimal length representing $g^{-1}h$, ie., it is a reduced word.
\end{proof}

\noindent
As a direct consequence of the previous lemma, we deduce:

\begin{cor}\label{cor:GPdist}
Let $g,h \in \Gamma, \mathcal{G}$. The distance between $g$ and $h$ in $X(\Gamma, \mathcal{G})$ is $|g^{-1}h|$. 
\end{cor}

\noindent
Before proving Proposition \ref{prop:Xquasimedian}, we demonstrate the following preliminary lemma.

\begin{lemma}\label{lem:Xtriangle}
The edges of a triangle in $\X$ are labelled by the same vertex-group.
\end{lemma}

\begin{proof}
Let $a,b,c \in \X$ be three pairwise adjacent vertices. Up to translating by $a^{-1}$, we may suppose without loss of generality that $a=1$. As a consequence, $b$ and $c$ belong to two vertex-group $G_u$ and $G_v$ respectively, which are the vertex-groups labelling the edges $(1,b)$ and $(1,c)$ respectively. On the other hand, $|b^{-1}c|=d(b,c)=1$ so $b$ and $c$ must belong to the same vertex-group, ie., $u=v$. We conclude that the vertex-group $G_u=G_v$ labels the three edges $(1,b)$, $(1,c)$ and $(b,c)$.
\end{proof}

\begin{proof}[Proof of Proposition \ref{prop:Xquasimedian}.]
First, we verify the triangle condition. So let $a,p,q \in \X$ be three vertices such that $p$ and $q$ are adjacent and $d(a,p)=d(a,q)=n$. Up to translating by $a^{-1} \in \Gamma \mathcal{G}$, we may suppose without loss of generality that $a=1$. In particular, $|p|=d(1,p)=n$ so that there exists a reduced word $p_1 \cdots p_n$ of length $n$ representing $p$. Then, because $p$ and $q$ are adjacent, there exists a non trivial syllable $s \in G_u$ such that $q=ps$. On the other hand, $|q|=d(1,q)=n$ whereas $ps$ is a word of length $n+1$. Therefore, there must exist a syllable $p_i \in G_u \backslash \{ s^{-1} \}$ which shuffles to the end in $p_1 \cdots p_n$ such that $p_1 \cdots p_{i-1}p_{i+1} \cdots p_n (p_is)$ is a reduced word representing $q$. Setting $c=p_1 \cdots p_{i-1}p_{i+1} \cdots p_n$, we have $d(1,c)=n-1$, $d(c,p)= |c^{-1}p|=|p_i|=1$ and $d(c,q)= |c^{-1}q|=|p_is|=1$. Thus, $c$ is the vertex we are looking for.

\medskip \noindent
Next, we verify the quadrangle condition. So let $u,p,q,r \in \X$ be three vertices such that $p,q$ are both adjacent to $r$, $d(u,p)=d(u,q)=n$ and $d(u,r)=n+1$. Up to translating by $u^{-1}$, we may suppose without loss of generality that $u=1$. In particular, $|p|=d(1,p)=n$ so there exists a reduced word $p_1 \cdots p_n$ of length $n$ representing $p$. Because $r$ is adjacent to $p$, there exists a non trivial syllable $a \in G_u$ such that $r=pa$; similarly, because $q$ is adjacent to $r$, there exists a non trivial syllable $b \in G_v$ such that $q=rb=pab$. Notice that, $p_1 \cdots p_n a$ being a word of length $n+1$ representing $r$, we deduce from $d(1,r)=n+1$ that this product is reduced. Since we already know that the triangle condition holds, we may suppose without loss of generality that $p$ and $q$ are not adjacent, hence $|ab|=|p^{-1}q|= d(p,q)=2$, so that $a$ and $b$ necessarily belong to distinct vertex-groups, ie., $u \neq v$. On the other hand, the word $p_1 \cdots p_n ab$ representing $q$ cannot be reduced since it has length $n+2$ whereas $d(1,q)=n$. Since we already know that $p_1 \cdots p_na$ is reduced, there must exist some syllable $p_i=b^{-1}$ which shuffles to the end in $p_1 \cdots p_na$, so that $p_1 \cdots p_{i-1} p_{i+1} \cdots p_na$ is a reduced word representing $q$. Setting $c=p_1 \cdots p_{i-1}p_{i+1} \cdots p_n$, we have $d(1,c)=n-1$, $d(c,p)=|c^{-1}p|=|p_i|=1$ and $d(c,q)=|c^{-1}q|= |a|=1$. Therefore, $c$ is the vertex we are looking for.

\medskip \noindent
Now, notice that, if $a,p \in \X$ are two vertices satisfying $d(a,p)=2$, then there exists at most two geodesics between $a$ and $p$. Indeed, according to Lemma \ref{lem:geodesic}, there exists a bijection between the geodesics between $a$ and $p$ and the reduced words representing $a^{-1}p$. Since $|a^{-1}p|=d(a,p)=2$, we conclude. As a consequence, $\X$ cannot contain an induced subgraph isomorphic to $K_{2,3}$. 

\medskip \noindent
Finally, let $(a,b,c,d)$ be a square in $\X$, where $b$ and $d$ are adjacent. According to Lemma \ref{lem:Xtriangle}, the edges of the triangle $a,b,d$ are labelled by the same vertex-group of $\Gamma$, and similarly, the edges of the triangle $b,c,d$ are labelled by the same vertex-group of $\Gamma$. A fortiori, the edges $(a,b)$ and $(b,c)$ are labelled by the same vertex-group of $\Gamma$, so that $a$ and $c$ must be adjacent in $\X$. This proves that $\X$ does not contain any induced subgraph isomorphic to $K_4^-$.
\end{proof}

\paragraph{Cliques and Prisms of $\X$.}
Now, we know that $\X$ is a quasi-median graph. Therefore, in order to understand its geometry, it is necessary to understand first what are its cliques, prisms, and hyperplanes. We begin with the cliques of $\X$.

\begin{lemma}\label{lem:Xclique}
The cliques of $\X$ are the subgraphs $gG_u$, where $g \in \Gamma \mathcal{G}$ and $u \in V(\Gamma)$.
\end{lemma}

\begin{proof}
Let $C$ be a clique. According to Lemma \ref{lem:Xtriangle}, all the edges of $C$ are labelled by the same vertex-group, say $G_u$. Thus, if we choose a vertex $g \in C$, then $C \subset gG_u$. On the other hand, since $G_u \backslash \{ 1 \}$ is included into the generating set used to define the Cayley graph $\X$, we know that $gG_u$ is a complete subgraph, hence $C=gG_u$. Conversely, we want to prove that, if $g \in \Gamma \mathcal{G}$ and $u \in V(\Gamma)$, then $gG_u$ defines a clique of $\X$. We know that $gG_u$ is a complete subgraph, so there exists a clique $C$ such that $gG_u \subset C$. On the other hand, as a consequence of Lemma \ref{lem:Xtriangle}, all the edges of $C$ must be labelled by the same vertex-group, which is necessarily $G_u$. Therefore, $C \subset gG_u$, so that $gG_u=C$ is a clique of $\X$. 
\end{proof}

\begin{cor}\label{cor:Xprisms}
The prisms of $\X$ are the translates of the subgraphs $\bigoplus\limits_{v \in \Lambda} G_v$, where $\Lambda$ is a complete subgraph of $\Gamma$. As a consequence, $\dim_{\square}(X)= \mathrm{clique}(\Gamma)$.
\end{cor}

\begin{proof}
It is clear that, if $\Lambda$ is a complete subgraph of $\Gamma$, then $P= \bigoplus\limits_{v \in \Lambda} G_v$ is a prism of $\X$ since we know from the previous lemma that each $G_v$ is a clique.

\medskip \noindent
Conversely, let $P$ be a prism of $\X$. Fix a vertex $x \in P$ and a collection of cliques $\mathcal{C}$ containing $x$ such that $P$ is the weak Cartesian product of the cliques of $\mathcal{C}$; up to translating $P$, we will suppose that $x=1$. Let $\Lambda$ denote the set of vertices of $\Gamma$ associated to the vertex-groups labelling the cliques of $\mathcal{C}$. As a consequence of the previous lemma, $\mathcal{C}= \{ G_u \mid u \in \Lambda \}$; and, as a consequence of Lemma \ref{lem:Xsquare}, we deduce from the fact that two distinct cliques $G_u$ and $G_v$, where $u,v \in \Lambda$, generate a prism that $u$ and $v$ must be two adjacent vertices of $\Gamma$, so that $\Lambda$ defines a complete subgraph of $\Gamma$. By definition, $P$ is the connected component of the Cartesian product of the cliques of $\mathcal{C}$ which contains $1$, which is precisely $P= \bigoplus\limits_{v \in \Lambda} G_v$.
\end{proof}

\paragraph{Hyperplanes of $\X$.}
According to Lemma \ref{lem:Xtriangle}, the edges of a triangle are labelled by the same vertex-group. Moreover, according to our next lemma, two opposite edges of some square in $\X$ are also labelled by the same vertex-group. As a consequence, all the edges of a given hyperplane of $\X$ are labelled by the same vertex-group. Thus, if $J_u$ denotes the hyperplane dual to the clique $G_u$ for every vertex $u \in V(\Gamma)$, then the hyperplanes of $\X$ are the $gJ_v$'s where $g \in \Gamma \mathcal{G}$ and $v \in V(\Gamma)$. 

\begin{lemma}\label{lem:Xsquare}
For any induced square $C \subset \X$, there exist an element $g \in \Gamma \mathcal{G}$, two adjacent vertices $u,v \in V(\Gamma)$, and two syllables $a \in G_u \backslash \{ 1 \}$, $b \in G_v \backslash \{ 1 \}$ such that $g$, $ga$, $gb$ and $gab=gba$ are the vertices of $C$. 
\end{lemma}

\begin{proof}
Let $g \in C$ be one vertex of our square. The two vertices of $C$ which are adjacent to $g$ can be written as $ga$ and $gb$ for some syllables $a \in G_u$ and $b \in G_v$. Let us write the last vertex of $C$ as $gac$ for some $c \in G_w$. Because our square is induced, we have $|b^{-1}ac|=d(gb,gac)=2$, so word $b^{-1}ac$ cannot be reduced. On the other hand, $|ac|=d(g,gac)=2$ so the word $ac$ is reduced; and $u \neq v$, since otherwise the vertices $ga$ and $gb$ would be adjacent, contradicting the fact that our square is induced, so the word $b^{-1}a$ is reduced. Finally, we conclude that $b^{-1}=c^{-1}$ shuffles to the end in $b^{-1}a$. Thus, $u$ and $v$ must be adjacent, and $gac=gab=gba$.
\end{proof}

\noindent
Since any hyperplane of $\X$ is a translate of some $J_u$, we only need to understand such a hyperplane in order to understand all the hyperplanes of $\X$. This is the purpose of the next lemma.

\begin{lemma}\label{lem:hyperplane}
Let $u \in V(\Gamma)$. An edge $e$ of $X(\Gamma, \mathcal{G})$ belongs to $J_u$ if and only if $e=(ah,ag)$ for some $a \in \langle \link(u) \rangle$ and $g,h \in G_u$ where $g \neq h$. 
\end{lemma}

\begin{figure}
\begin{center}
\includegraphics[scale=0.7]{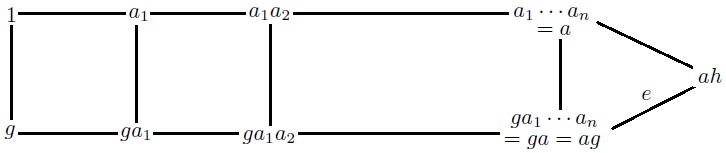}
\end{center}
\caption{The edge $(ag,ah)$ belongs to $J_u$.}
\label{hyperplane}
\end{figure}

\begin{proof}
Let $e_0$ denote the edge $(1,g)$ for some fixed $g \in G_u \backslash \{ 1 \}$. If an edge $e$ belongs to $J_u$, then there exists a sequence of edges $e_0,e_1, \ldots, e_{n-1},e_n=e$ such that, for every $1 \leq i \leq n-1$, the edges $e_i$ and $e_{i+1}$ either are parallel sides of some square or belong to the same clique. We want to show by induction over $n$ that there exist some $a \in \langle \link(u) \rangle$ and some distinct $g,h \in G_u$ such that $e=(ah,ag)$. 

\medskip \noindent
If $n=0$, then $e=e_0=(1,g)$, and there is nothing to prove. If $n \geq 1$, then our induction hypothesis implies that $e_{n-1}=(ah,ag)$ for some $a \in \langle \link(u) \rangle$ and $g,h \in G_u$ with $g \neq h$. Now, two cases may happen: $e_{n-1}$ and $e=e_n$ either are parallel sides of some square or belong to the same clique. In the latter case, because the edge $e_{n-1}$ is labelled by an element of $G_u$, we deduce that there exist some $b,c \in G_u$ such that $e=(ahb,agc)$, and we are done. In the former case, there exist a vertex $v \in V(\Gamma)$, adjacent to $u$, and some $k \in G_v \backslash \{ 1 \}$ such that $e=(ahk,agk)$; since $u$ and $v$ are adjacent in $\Gamma$, $k$ commutes with $h$ and $g$, hence $e=((ak)h,(ak)g)$ with $ak \in \langle \link(u) \rangle$. The conclusion follows.

\medskip \noindent
Conversely, suppose that $e=(ah,ag)$ for some $a \in \langle \link(u) \rangle$ and $g,h \in G_u$ with $g \neq h$. If we write $a$ as a reduced word $a_1 \cdots a_n$, then $$(1,g),(a_1,ga_1),(a_1a_2,ga_1a_2) \ldots, (a_1 \cdots a_n,ga_1 \cdots a_n)=(a,ga)=(a,ag)$$ defines a sequence of edges such that two successive edges are parallel sides of some square; see Figure \ref{hyperplane}. Then, the edges $(a,ag)$ and $(ah,ag)$ belong to the same simplex, so we conclude that $(1,g)$ and $(ah,ag)$ are dual to the same hyperplane; this precisely means that $(ah,ag) \in J_u$.
\end{proof}
\begin{cor}\label{cor:hyperplane}
For every $u \in V(\Gamma)$, $N(J_u)= \langle \st(u) \rangle$. 
\end{cor}

\begin{proof}
If $x \in N(J_u)$, then there exists an edge $e \in J_u$ with $x$ as an endpoint. According to the previous lemma, there must exist some $a \in \langle \link(u) \rangle$ and $h \in G_u$ such that $x=ah$. A fortiori, $x \in \langle \st(u) \rangle$. Conversely, suppose that $x \in \langle \mathrm{star}(u) \rangle$. Since $\langle \st(u) \rangle= \langle \link(u) \rangle \times G_u$, we can write $x=ah$ for some $a \in \langle \link(u) \rangle$ and $h \in G_u$. If we choose some $g \in G_u \backslash \{ h \}$, we deduce from the previous lemma that $e=(ah,ag)$ is an edge of $J_u$ with $x$ as an endpoint. Therefore, $x \in N(J_u)$.
\end{proof}

\noindent
Now, we want to understand how the vertex-groups labelling two hyperplanes behave when these two hyperplanes are transverse or tangent. As a consequence, we will be able to deduce the following proposition.

\begin{prop}\label{prop:Xspecial}
The action $\Gamma \mathcal{G} \curvearrowright \X$ is special.
\end{prop}

\noindent
Our statement describing the behaviour of vertex-groups labelling hyperplanes is the following:

\begin{lemma}\label{lem:labelhyp}
If two hyperplanes of $\X$ are transverse, then they are labelled by two adjacent vertex-groups; if they are tangent, then the vertex-groups which label them are distinct and non adjacent. 
\end{lemma}

\begin{proof}
If $e$ and $f$ are two adjacent edges, then $e=(g,gh)$ and $f=(g,gk)$ for some $g,h,k \in \Gamma \mathcal{G}$ where $h$ and $k$ are two syllables which belong to two vertex-groups $G_u$ and $G_v$ respectively. Let $J$ and $H$ denote respectively the hyperplanes dual to $e$ and $f$. Then three cases may happen:
\begin{itemize}
	\item $u=v$, so that $e$ and $f$ belong to the same clique, hence $J=H$;
	\item $(u,v) \in E(\Gamma)$, so that $e$ and $f$ generate a square in $\X$, hence $J$ and $H$ are transverse;
	\item $u \neq v$ and $(u,v) \notin E(\Gamma)$, so that $e$ and $f$ do not generate a square in $\X$, hence $J$ and $H$ tangent.
\end{itemize}
This proves our lemma.
\end{proof}

\begin{proof}[Proof of Proposition \ref{prop:Xspecial}.]
Let $J$ be a hyperplane and $g \in \Gamma \mathcal{G}$. Since $J$ and $gJ$ are labelled by the same vertex-group, it follows from the previous lemma that $J$ and $gJ$ cannot be transverse nor tangent. 

\medskip \noindent
Next, let $J_1,J_2,J_3$ be three hyperplanes such that $J_1$ is tangent to $J_2$, and $J_3$ transverse to $J_2$. Let $G_1,G_2,G_3$ denote the vertex-groups labelling the hyperplanes $J_1,J_2,J_3$ respectively. It follows from the previous lemma that $G_3$ is adjacent to $G_2$ but $G_1$ cannot be adjacent to $G_2$. Consequently, there does not exist some $g \in \Gamma \mathcal{G}$ satisfying $gJ_1=J_3$.
\end{proof}

\paragraph{Flat rectangles of $\X$.} 
We conclude this section by describing the flat rectangles of $\X$.

\begin{lemma}\label{lem:flatrectangle}
An induced subgraph $R \subset \X$ is a flat rectangle if and only if there exist a join subgraph $\Lambda_1 \ast \Lambda_2 \leq \Gamma$ and syllables $g_1, \ldots, g_n \in \langle \Lambda_1 \rangle$, $h_1, \ldots,h_m \in \langle \Lambda_2 \rangle$ such that the products $g_1 \cdots g_n$ and $h_1 \cdots h_m$ are reduced and such that $R$ is generated by the vertices $$ \{ k g_1 \cdots g_i h_1 \cdots h_j \mid 0 \leq i \leq n, 0 \leq j \leq m \}$$ for some $k \in \Gamma \mathcal{G}$.
\end{lemma}

\begin{proof}
Let $R$ be a flat rectangle. For convenience, we will identify it with the 1-skeleton of the Euclidean rectangle $[0,m] \times [0,n]$ to which it is isometric. Up to translating by an element of $\Gamma \mathcal{G}$, and so replacing $R$ with $k^{-1} R$ for some $k \in \Gamma \mathcal{G}$, we may suppose without loss of generality that the vertices $(0,0)$ is the identity $1$. Let $g$ (resp. $h$) denote the vertex $(0,n)$ (resp. $(m,0)$). Because $R$ is isometrically embedded, $\{0 \} \times [0,n]$ is a geodesic between $1$ and $g$, so we can write $g$ as reduced word $g_1 \cdots g_n$ such that the vertex $(0,i)$ corresponds to $g_1 \cdots g_i$ for every $0 \leq i \leq n$. Similarly, let $h_1 \cdots h_m$ be a reduced word representing $h$ such that the vertex $(i,0)$ corresponds to $h_1 \cdots h_i$ for every $0 \leq i \leq m$. 

\medskip \noindent
For every $1 \leq i \leq n$, let $u_i \in V(\Gamma)$ be such that $g_i \in G_{u_i}$; similarly, for every $1 \leq i \leq m$, let $v_i \in V(\Gamma)$ be such that $h_i \in G_{v_i}$. Now, let $\Lambda$ and $\Xi$ denote the subgraphs respectively generated by $\{ u_1, \ldots, u_n \}$ and $\{v_1, \ldots, v_m \}$. Notice that, for every $1 \leq i \leq n$ and $1 \leq j \leq m$, the hyperplanes $g_1 \cdots g_{i-1}J_{u_i}$ and $h_1 \cdots h_{j-1} J_{v_j}$, which are respectively dual to the ediges $[i-1,i] \times \{ 0 \}$ and $\{0\} \times [j-1,j]$, are transverse; so it follows from Lemma \ref{lem:labelhyp} that $u_i$ and $v_j$ are two adjacent vertices. Therefore, the subgraph generated by $\Lambda$ and $\Xi$ is a join subgraph $\Lambda \ast \Xi \leq \Gamma$. 

\medskip \noindent
Now, we claim that, if $w$ denotes the vertex $(i,j)$, then $w=h_1 \cdots h_i g_1 \cdots g_j$. Indeed, the initial segment of length $i$ of the geodesic $([0,i] \times \{ 0 \} ) \cup ( \{i\} \times [0,j])$ is labelled by $h_1 \cdots h_i$, and similarly, the initial segment of length $j$ of the geodesic $( \{i\} \times [0,j]) \cup ([0,i] \times \{ 0 \})$ is labelled by $g_1 \dots g_j$, so we deduce that $h_1 \cdots h_i$ and $g_1 \cdots g_j$ are two prefixes of $w$ whose supports are disjoint. Therefore, $h_1 \cdots h_i g_1 \cdots g_j$ must be itself a prefix of $w$. On the other hand, $|w|=d(1,w)=i+j$, hence $w=h_1 \cdots h_i g_1 \cdots g_j$. 

\medskip \noindent
Thus, the flat rectangle $R$ is generated by the vertices $$ \{ k g_1 \cdots g_i h_1 \cdots h_j \mid 0 \leq i \leq n, 0 \leq j \leq m \}.$$ Conversely, suppose that there exist a join subgraph $\Lambda_1 \ast \Lambda_2 \leq \Gamma$ and syllables $g_1, \ldots, g_n \in \langle \Lambda_1 \rangle$, $h_1, \ldots,h_m \in \langle \Lambda_2 \rangle$ such that the products $g_1 \cdots g_n$ and $h_1 \cdots h_m$ are reduced and such that $R$ is generated by the vertices $$ \{ k g_1 \cdots g_i h_1 \cdots h_j \mid 0 \leq i \leq n, 0 \leq j \leq m \}$$ for some $k \in \Gamma \mathcal{G}$. From our description of the geodesics in $\X$, we notice that $R$ is a geodesic subgraph, so that it is in particular isometrically embedded; furthermore, we deduce that the map which sends $k g_1 \cdots g_i h_1 \cdots h_j$ to $(i,j)$, for every $1 \leq i \leq n$ and $1 \leq j \leq m$, defines an isomorphism $R \to [0,n] \times [0,m]$.
\end{proof}

\subsection{Properties of graph products}

\noindent
So far, we have proved that graph products act on specific quasi-median graphs and we have studied these actions and the geometry of these quasi-median graphs. Now, we are ready to apply the different combination results proved in Section \ref{section:topicalactionsI} and Section \ref{section:topicalactionsII}; all their hypotheses are verified by the following statement.

\begin{prop}\label{prop:GPsumup}
Let $\Gamma$ be a simplicial graph and $\mathcal{G}$ a collection of groups labelled by $V(\Gamma)$. 
The following assertions hold:
\begin{itemize}
	\item[(i)] the action $\Gamma \mathcal{G} \curvearrowright \X$ is topical-transitive;
	\item[(ii)] any vertex of $\X$ belongs to $\# V(\Gamma)$ cliques;
	\item[(iii)] vertex-stabilisers are trivial;
	\item[(iv)] $\dim_{\square} \X = \mathrm{clique}(\Gamma)$;
	\item[(v)] $\X$ contains $\# V(\Gamma)$ $\Gamma \mathcal{G}$-orbits of cliques;
	\item[(vi)] for every prism $P=C_1 \times \cdots \times C_n$, $\mathrm{stab}(P)= \mathrm{stab}(C_1) \times \cdots \times \mathrm{stab}(C_n)$;
	\item[(vii)] the action $\Gamma \mathcal{G} \curvearrowright \X$ is special;
	\item[(viii)] for every hyperplane $J$, $\mathrm{stab}(J)$ is a rectract of $\Gamma \mathcal{G}$.
\end{itemize}
\end{prop}

\noindent
We first need to describe the projection onto a clique of $\X$. Given a vertex-group $G \in \mathcal{G}$, we define a map $\pi_G : \Gamma \mathcal{G} \to G$ by setting, for every $g \in \Gamma \mathcal{G}$, $\pi_G(g)$ as the syllable of $\mathrm{head}(g)$ which belongs to $G$, if such a syllable exists, and $1$ otherwise.

\begin{lemma}\label{lem:projXclique}
Let $C$ be the clique of $\X$ corresponding to a vertex-group $G \in \mathcal{G}$, and $g \in \Gamma \mathcal{G}$ a vertex of $\X$. Then $\pi_G(g)$ is the projection of $g$ onto $C$. 
\end{lemma}

\begin{proof}
If $\pi_G(g)=1$, then, for every $s \in G$, the product $sg$ is reduced so that $|sg|>|g|$. Therefore, $1$ is the vertex of $C$ minimizing the distance to $g$, ie., $1= \pi_G(g)$ is the projection of $g$ onto $C$. 

\medskip \noindent
If $\pi_G(g)=g_1$, write $g$ as a reduced product $g_1g_2$. By noticing that $(sg_1)g_2$ is a reduced product, we deduce that the length $|sg|=|sg_1g_2|$, where $s \in G$, is minimised when $s=g_1^{-1}$; a fortiori, $d(s^{-1},g)=|sg|$ is minimised when $s= g_1^{-1}$, ie., $g_1=\pi_G(g)$ is the projection of $g$ onto $C$. 
\end{proof}

\begin{proof}[Proof of Proposition \ref{prop:GPsumup}.]
Point $(ii)$ and Point $(v)$ follow from Lemma \ref{lem:Xclique}, Point $(iii)$ is clear, Point $(iv)$ and Point $(vi)$ follows from Corollary \ref{cor:Xprisms}, and Point $(vii)$ follows from Proposition \ref{prop:Xspecial}. 

\medskip \noindent
Next, according to Corollary \ref{cor:hyperplane}, it is sufficient to prove that, for every vertex $v \in V(\Gamma)$, the subgroup $\langle \mathrm{star}(v) \rangle$ is a retract of $\Gamma \mathcal{G}$ in order to deduce Point $(viii)$. The map
$$g \mapsto \left\{ \begin{array}{cl} g & \text{for every} \ g \in G_u, \ u \in \mathrm{star}(v) \\ 1 & \text{for every} \ g \in G_u, \ u \notin \mathrm{star}(v) \end{array} \right. $$
induces a required retraction $\Gamma \mathcal{G} \twoheadrightarrow \langle \mathrm{star}(v) \rangle$.

\medskip \noindent
Finally, let us prove Point $(i)$. 
Fix a clique $C$ of $\X$. Up to translating by an element of $\Gamma \mathcal{G}$, we may suppose that $C$ corresponds to some vertex-group $G_u$, where $u \in V(\Gamma)$. According to Lemma \ref{lem:Xclique}, the action $\mathrm{stab}(C) \curvearrowright C$ corresponds to the action of a subgroup on itself by left-multiplication, so it must be free and transitive. Now, we claim that $\mathrm{Im}(\rho_C) = \rho_C( \mathrm{stab}(C))$. The inclusion $\rho_C(\mathrm{stab}(C)) \subset \mathrm{Im}(\rho_C)$ is clear. Conversely, let $g \in \mathrm{stab}(J_u)$ and $h \in N(J_u)$. It follows from Corollary \ref{cor:hyperplane} that 
$$\mathrm{stab}(J_u)= \langle \mathrm{star}(u) \rangle = \langle \mathrm{link}(u) \rangle \times G_u \subset \Gamma \mathcal{G}$$
and that
$$N(J_u) =  \langle \mathrm{star}(u) \rangle = \langle \mathrm{link}(u) \rangle \times G_u \subset \X.$$
Therefore, we can write $g=g_1g_2$ and $h=h_1h_2$ where $g_1,h_1 \in \langle \mathrm{link}(u) \rangle$ and $g_2,h_2 \in G_u$. By noticing that
$$\pi_{G_u}(gg_2^{-1} \cdot h) = \pi_{G_u}( g_1h_1 \cdot h_2)= h_2 = \pi_{G_u} (h_1h_2)= \pi_{G_u}(h),$$
we deduce from Lemma \ref{lem:projXclique} that $gg_2^{-1} \cdot h$ and $h$ belong to the same sector delimited by $J_u$. Since this is true for every $h \in N(J_u)$, we know that $gg_2^{-1}=g_1$ stabilises each sector delimited by $J_u$. Therefore,
$$\rho_C(g)= \rho_C(g_1^{-1}g) = \rho_C(g_2),$$
so that $g_2 \in G_u = \mathrm{stab}(C)$ induces the same permutation on $C$ as $g$. This proves our claim. Thus, we have proved that the action $\Gamma \mathcal{G} \curvearrowright \X$ is topical-transitive.
\end{proof}

\noindent
Let us begin by studying cubical properties of graph products. First, Proposition \ref{prop:properlydiscontinuous} gives:

\begin{thm}\label{thm:GPccproper}
A graph product of groups acting properly discontinously on CAT(0) cube complexes acts properly discontinuously on a CAT(0) cube complex.
\end{thm}

\noindent
It is worth noticing that, in this statement, no assumption is made on the graph along which the graph product is defined. Next, when this graph is finite, we are able to prove more thanks to Proposition \ref{prop:CAT0metricallyproper}, Proposition \ref{prop:cubulatinggeometrically} and Proposition \ref{prop:cocompactspecial}. (Notice that Haglund already proved in \cite[Theorem J]{MR2413337} that a graph product (along a finite graph) of special groups is special.)

\begin{thm}\label{thm:GPcc}
Let $\Gamma$ be a finite simplicial graph and $\mathcal{G}$ a collection of groups labelled by $V(\Gamma)$. Suppose that the groups of $\mathcal{G}$ act metrically properly (resp. geometrically, geometrically and virtually specially) on CAT(0) cube complexes, then the graph product $\Gamma \mathcal{G}$ acts metrically properly (resp. geometrically, geometrically and virtually specially) on some CAT(0) cube complex.
\end{thm}

\noindent
In fact, when the simplicial graph is countable, we are also able to construct metrically proper actions by weighting the local metrics, or equivalently by duplicating the walls, as shown in the proof of the next proposition. 

\begin{prop}\label{prop:GPmetricallyproper}
Let $\Gamma$ be a countable simplicial graph and $\mathcal{G}$ a collection of groups labelled by $V(\Gamma)$. Suppose that the groups of $\mathcal{G}$ act metrically properly on CAT(0) cube complexes (resp. spaces with measured walls, spaces with labelled partitions). Then the graph product $\Gamma \mathcal{G}$ acts metrically properly on some CAT(0) cube complex (resp. space with measured walls, space with labelled partitions).
\end{prop}

\begin{proof}[Sketch of proof.]
We will write the proof only for spaces with measured walls. The other cases are similar. For convenience, set $X=\X$. Recall from Proposition \ref{prop:GPsumup} that the action $\Gamma \mathcal{G} \curvearrowright X$ is topical-transitive.

\medskip \noindent
Let $\{u_1,u_2, \ldots\}$ denote the vertices of $\Gamma$. For convenience, the vertex-group $G_{u_i}$ will be denoted by $G_i$. For every $i \geq 1$, since $G_i$ is a-T-menable by assumption, there exists a space with measured walls $(X_i, \mathcal{W}_i, \mathcal{B}_i, \mu_i)$ on which $G_i$ acts metrically properly; by reproducing the construction of Lemma \ref{lem:modifycubing}, we may suppose without loss of generality that $X_i$ contains a point $x_i$ whose stabiliser is trivial. Since the action of $G_i$ on $(X_i, \mathcal{W}_i, \mathcal{B}_i, \mu_i)$ is metrically proper and that the stabiliser of $x_i$ is trivial, there must exists some constant $K_i>0$ such that $d_{\mathcal{W}_i}(x_i,gx_i) \geq K_i$ for every $g \in G_i \backslash \{ 1 \}$. Set $n(i)= \lfloor i/K_i \rfloor +1$ for every $i \geq 1$. Notice that, if $\mathcal{W}_i^{n(i)}$ denotes the disjoint union of $n(i)$ copies of $\mathcal{W}_i$, then we get naturally a new space with measured walls $(X_i, \mathcal{W}_i^{n(i)}, \mathcal{B}_i^{n(i)}, \mu_i^{n(i)})$ on which $G_i$ also acts metrically properly. Next, we pullback the walls from $X_i$ to $G_i$ thanks to the map $g \mapsto g \cdot x_i$, making a space with measured walls $(G_i, \mathcal{M}_i, \mathcal{D}_i, \nu_i)$ (see also \cite[Lemma 3.9]{medianviewpoint}). Notice that, for every disctinct $g,h \in G_i$,
$$d_{\mathcal{M}_i}(g,h) = d_{\mathcal{W}_i^{n(i)}} (gx_i,hx_i) = n(i) \cdot d_{\mathcal{W}_i}(x_i,g^{-1}hx_i) \geq n(i)K_i \geq i.$$
Finally, use Proposition \ref{prop:mwtopicalI} in order to extend these local spaces with measured walls to obtain a space with measured walls $(X, \mathcal{W}, \mathcal{B}, \mu)$. 

\medskip \noindent
We claim that $\Gamma \mathcal{G}$ acts metrically properly on $(X, \mathcal{W}, \mathcal{B}, \mu)$. Following the proof of Lemma \ref{lem:whenmetricallyproper}, it is sufficient to show that $(X, \mathcal{W}, \mathcal{B}, \mu)$ is locally finite (with respect to $d_{\mathcal{W}}$). So fixing some vertex $x \in X$ and some constant $R \geq 0$, we want to prove that the ball $B_{\mathcal{W}}(x,R)$ is finite. Following the proof of Lemma \ref{lem:whenlocallyfinite}, we notice that it is sufficient to show that the ball $B_{\mathcal{W}}(x,1)$ is finite and that
$$B_{\mathcal{W}}(x,1) \subset \bigcup\limits_{i \geq 1} B_{\mathcal{W}(C_i)}(x,1),$$
where $C_i$ denotes the unique clique labelled by $G_i$ containing $x$ for every $i \geq 1$. Noticing that $C_i \cap B_{\mathcal{W}(C_i)} (x,1) = \{x \}$ for every $i \geq 2$, we deduce the infinite union of the previous expression turns out to be finite. This concludes the proof.
\end{proof}

\noindent
Since discrete a-T-menable (resp. a-$L^p$-menable) groups are precisely the groups which act metrically properly on some space with measured walls (resp. specific space with labelled partitions), we deduce the following statement (first proved in \cite{AntolinDreesen}, without restriction on $p$ for the a-$L^p$-menability).

\begin{cor}\label{cor:aTaB}
Graph products of discrete a-T-menable (resp. a-$L^p$-menable) groups along countable simplicial graphs are a-T-menable (resp. a-$L^p$-menable, where $p \notin 2 \mathbb{Z}$).
\end{cor}

\noindent
Next, Theorem \ref{thm:producingCAT0groups} produces the following statement.

\begin{thm}\label{thm:GPCAT0}
Graph products of CAT(0) groups along finite graphs are CAT(0).
\end{thm}

\noindent
Now, let us state and prove the following sufficient condition for two graph products to be quasi-isometric. Our result should be compared with \cite{GPcommensurability}, where a sufficient condition for two graph products to be commensurable is proved. More precicely, we will prove that some graph products are \emph{Lipschitz-equivalent}. 

\begin{definition}
Given some $K >0$, two metric spaces $X,Y$ are \emph{$K$-Lipschitz equivalent}\index{Lipschitz equivalence} if there exists a bijection $f : X \to Y$ such that, for every $x,y \in X$, 
$$\frac{1}{K} \cdot d(x,y) \leq d(f(x),f(y)) \leq K \cdot d(x,y).$$
We say that two metric spaces are \emph{Lipschitz equivalent} if they are $K$-Lipschitz equivalent for some $K>0$.
\end{definition}

\noindent
It is worth noticing that, according to \cite{QIvsLIP}, being Lipschitz equivalent is strictly stronger than being quasi-isometric; on the other hand, it is proved in \cite{WhyteLip} that two non amenable finitely generated groups are quasi-isometric if and only if they are Lipschitz equivalent. Our criterion is the following:

\begin{thm}\label{thm:GPlip}
Let $\Gamma$ be a finite simplicial graph and $\mathcal{G}, \mathcal{H}$ two collections of finitely generated groups labelled by $V(\Gamma)$. Suppose that, for every $v \in V(\Gamma)$, the vertex-groups $G_v$ and $H_v$ are Lipschitz equivalent. Then the graph products $\Gamma \mathcal{G}$ and $\Gamma \mathcal{H}$ are Lipschitz equivalent.
\end{thm}

\noindent
We begin by proving the following general preliminary statement:

\begin{lemma}\label{lem:Lipequi}
Let $X$ be a quasi-median graph and $\{ (C,\delta_C) \mid C\}$, $\{ (C,\eta_C) \mid C \}$ two coherent systems of metrics. Suppose that there exist some $K>0$ such that, for every clique $C$, the identity $(C,\delta_C) \to (C,\eta_C)$ is a $K$-Lipschitz equivalence. Then the global metrics $\delta$ and $\eta$ are $K$-Lipschitz equivalent.
\end{lemma}

\begin{proof}
Let $x,y \in X$ be two vertices. Let $J_1, \ldots, J_n$ denotes the hyperplanes separating $x$ and $y$, fix a clique $C_i$ dual to $J_i$ for every $1 \leq i \leq n$, and let $p_i : X \to C_i$ denote the projection onto $C_i$. Then
$$\delta(x,y)= \sum\limits_{i=1}^n \delta_{C_i}(p_i(x),p_i(y)) \leq K \sum\limits_{i=1}^n \eta_{C_i}(p_i(x),p_i(y)) = K \cdot \eta(x,y),$$
and similarly, $\eta(x,y) \leq K \cdot \delta(x,y)$. Therefore, $\delta$ and $\eta$ are $K$-Lipschitz equivalent.
\end{proof}

\begin{proof}[Proof of Theorem \ref{thm:GPlip}.]
For every $v \in V(\Gamma)$, let $\varphi_v : G_v \to H_v$ be a Lipschitz equivalence, where we fixed two word metrics $\delta_v^1$ and $\delta_v^2$ on $G_v$ and $H_v$ with respect to finite generating sets $S_v$ and $R_v$ respectively; up to post-composing $\varphi_v$ with a translation of $H_v$, we may suppose without loss of generality that $\varphi_v(1)=1$. From this collection of maps, we define the following bijection:
$$\varphi : \left\{ \begin{array}{ccc} \bigsqcup\limits_{v \in V(\Gamma)} G_v & \to & \bigsqcup\limits_{v \in V(\Gamma)} H_v \\ g & \mapsto & \varphi_v(g) \ \text{if $g \in G_v$} \end{array} \right.$$
Finally, let $\Phi$ be the map which sends a word $w(\ell_1, \ldots, \ell_r)$, where $\ell_i \in \bigsqcup\limits_{v \in V(\Gamma)} G_v$ for every $1 \leq i \leq r$, to the word $w(\varphi(\ell_1),\ldots, \varphi(\ell_r))$. Because $\Phi$ sends two words equal in $\Gamma \mathcal{G}$ to two words equal in $\Gamma \mathcal{H}$, it induces a bijection $\Gamma \mathcal{G} \to \Gamma \mathcal{H}$. Moreover, because $\Phi$ sends a reduced word to a reduced word of the same length, we deduce that $\Phi$ induces an isomorphism $\X \to X(\Gamma, \mathcal{H})$. 

\medskip \noindent
Now, for every clique $C$ of $\X$, write $C=gG_u$ for some $g \in \Gamma \mathcal{G}$ and $u \in V(\Gamma)$ such that $\mathrm{tail}(g)$ does not contain a syllable which belongs to $G_u$ (notice that there exists a unique such element satisfying $C=gG_u$), and use the bijection $G_u \to C$ defined by $h \mapsto gh$ to transfer the metric $\delta_u^1$ from $G_u$ to $C$. We claim that our system of metrics $\{ (C, \delta_C^1) \mid C \ \text{clique} \}$ is coherent and $\Gamma \mathcal{G}$-invariant. 

\medskip \noindent
Let $C$ be a clique, $x,y \in C$ two vertices and $g \in \Gamma \mathcal{G}$. Write $C=hG_u$ for some $u \in V(\Gamma)$ and $h \in \Gamma \mathcal{G}$ such that the tail of $h$ does not contain a syllable belonging to $G_u$. Next, write $g=g_1g_2$ where $g_2$ is the syllable of $\mathrm{tail}(g)$ which belongs to $G_u$ and which commutes with $h$, if such a syllable exists, or $g_2=1$ otherwise. Notice that $gh=g_1hg_2$ where the tail of $g_1h$ does not contain a syllable belonging to $G_u$. As a consequence, writting $x=ha$ and $y=hb$ for some $a,b \in G_u$, we have
$$\delta_{gC}^1(gx,gy)= \delta_u^1(g_2a,g_2b).$$
On the other hand, $\delta_u$ is $G_u$-left-invariant, hence
$$\delta_{gC}^1(gx,gy)= \delta_u^1(a,b)= \delta_C^1(x,y).$$
Thus, we have proved that our system of metrics is $\Gamma \mathcal{G}$-invariant. Next, let $C,C'$ be two cliques dual to the same hyperplane. Because we already know that our system of metrics is $\Gamma \mathcal{G}$-invariant, we may suppose without loss of generality that $C=G_u$ for some $u \in V(\Gamma)$. Write $C'=gG_u$ where $g \in \Gamma \mathcal{G}$ is such that its tail does not contain a syllable belonging to $G_u$. (Recall that, since $C$ and $C'$ are dual to the same hyperplane, they are labelled by the same vertex-group.) Let $gx,gy \in C'$ be two vertices. By noticing that $\mathrm{supp}(g) \subset \mathrm{link}(u)$ as a consequence of Lemma \ref{lem:hyperplane}, we deduce from Lemma \ref{lem:projXclique} that
$$t_{C' \to C}(gx)= \pi_G(gx)=x \ \text{and} \ t_{C' \to C}(gy)= \pi_G(gy) = y,$$
hence
$$\delta_{C'}^1(gx,gy) = \delta_u(x,y)= \delta_C( t_{C' \to C}(gx), t_{C' \to C}(gy)).$$
Thus, we have proved that our system of metrics is coherent, concluding the proof of our claim.

\medskip \noindent
The global metric $\delta^1$ extending the $\delta_u^1$'s turns out to be very natural:

\begin{claim}\label{claim:Cayley}
$(\X, \delta^1)$ is isometric to the Cayley graph $\mathrm{Cayl} \left( \Gamma \mathcal{G}, S= \bigsqcup\limits_{u \in V(\Gamma)} S_u \right)$, where the $S_u$'s are the finite generating sets we used to defined $\delta^1$. 
\end{claim}

\noindent
Let $g,h \in \X$ be two vertices. Write $g^{-1}h$ as a reduced product $\ell_1 \cdots \ell_n$, where, for every $1 \leq i \leq n$, $\ell_i \in G_{u_i}$ for some $u_i \in V(\Gamma)$. According to Lemma \ref{lem:geodesic},
$$g, \ g \ell_1, \ g\ell_1 \ell_2, \ldots, \ g \ell_1 \cdots \ell_n$$
defines a geodesic in $\X$ from $g$ to $h$. By noticing that, for every $0 \leq i \leq n-1$, $g \ell_1 \cdots \ell_i$ and $g \ell_1 \cdots \ell_{i+1}$ belongs to the clique $g \ell_1 \cdots \ell_i G_{u_{i+1}}$, we deduce that
$$\delta^1(g,h)= \sum\limits_{i=0}^{n-1} \delta_{g \ell_1 \cdots \ell_i G_{u_{i+1}}} (g \ell_1 \cdots \ell_i, g \ell_1 \cdots \ell_{i+1}) = \sum\limits_{i=0}^{n-1} \delta_{G_{u_{i+1}}} (1,\ell_{i+1}) = \sum\limits_{i=1}^n \mathrm{lg}_{S_u}(\ell_i).$$
On the other hand, it follows from the normal form of graph products described at the beginning of Section \ref{section:appli1} that 
$$\sum\limits_{i=1}^n \mathrm{lg}_{S_u}(\ell_i) = \mathrm{lg}_S(g^{-1}h)= d_{\mathrm{Cayl}(\Gamma \mathcal{G},S)}(g,h),$$
concluding the proof of our claim.

\medskip \noindent
Similarly, from the $\delta_u^2$'s, we define a coherent and $\Gamma \mathcal{H}$-invariant system of metrics $\{ (C,\delta_C^2) \mid C \ \text{clique} \}$ on $X(\Gamma , \mathcal{H})$ such that $(X(\Gamma, \mathcal{H}),\delta^2)$ is isometric to the Cayley graph $\mathrm{Cayl} \left( \Gamma \mathcal{H}, \bigsqcup\limits_{u \in V(\Gamma)} R_u \right)$. 

\medskip \noindent
By construction of $\Phi$, for every clique $C$ of $\X$, the identity $(C,\delta_C^1) \to (C,\Phi^{-1} \delta_{\Phi(C)}^)2$ is biLipschitz, so that we deduce from Lemma \ref{lem:Lipequi} that $\Phi$ induces a Lipschitz equivalence between $(\X,\delta^1)$ and $(X(\Gamma, \mathcal{H}), \delta^2)$. A fortiori, the graph products $\Gamma \mathcal{G}$ and $\Gamma \mathcal{H}$ have two Lipschitz equivalent Cayley graphs, which concludes the proof.
\end{proof}

\noindent
It is worth noticing that, in the previous proof, we showed that the quasi-median graphs $X(\Gamma, \mathcal{H})$ and $X(\Lambda, \mathcal{G})$ are isometric if $\Gamma = \Lambda$ and if $\# G_v = \# H_v$ for every $v \in V(\Gamma)$. 

\begin{fact}\label{fact:isometricX}
Let $\Gamma$ be a simplicial graph and $\mathcal{G}, \mathcal{H}$ two collections of groups labelled by $V(\Gamma)$. Suppose that, for every vertex $v \in V(\Gamma)$, $G_v$ and $H_v$ have the same cardinality. The quasi-median graphs $\X$ and $X(\Gamma,\mathcal{H})$ are isometric.
\end{fact}

\noindent
A natural problem is to determine whether or not the converse holds. We suspect that the answer is positive (except maybe in a few particular cases). 

\begin{question}
When are the quasi-median graphs $X(\Gamma, \mathcal{H})$ and $X(\Lambda, \mathcal{G})$ isometric?
\end{question}

\noindent
Finally, let us apply Proposition \ref{prop:equicompression} to find a lower bound on the equivariant $\ell^p$-compression of graph products. The inequality we obtain was previously stated and proved in \cite[Corollary 4.4]{AntolinDreesen}. 

\begin{thm}\label{thm:GPcompression}
Let $\Gamma$ be a finite simplicial graph and $\mathcal{G}$ a collection of finitely generated groups indexed by $V(\Gamma)$. For every $p \geq 1$,
$$\alpha^*_p(\Gamma \mathcal{G}) \geq \min \left( \frac{1}{p}, \min\limits_{G \in \mathcal{G}} \alpha_p^*(G) \right).$$
\end{thm}

\begin{proof}
According to Proposition \ref{prop:GPsumup}, $\Gamma \mathcal{G}$ acts $\mathcal{C}$-topically-transitively on $\X$, where $\mathcal{C}$ denotes the collection of cliques corresponding to the vertex-groups. Now, endow every clique $C \in \mathcal{C}$ with the word metric associated to some finite generating set of the vertex-group labelling $C$, and extend this data to a system of metrics on $\X$ thanks to Proposition \ref{prop:mtopicalI}. Proposition \ref{prop:equicompression} applies to such a system (see Remark \ref{remark:whichsystem}, where our base point is $1$), so
$$\alpha_p^*( \X, \delta) \geq \min \left( \frac{1}{p}, \min\limits_{G \in \mathcal{G}} \alpha_p^*(G) \right).$$
On the other hand, it follows from Claim \ref{claim:Cayley} that $(\X, \delta)$ is isometric to the Cayley graph of $\Gamma \mathcal{G}$ with respect to some finite generating set, hence $\alpha_p^*(\Gamma \mathcal{G}) = \alpha_p^*( \X, \delta)$. The conclusion follows.
\end{proof}

\begin{remark}\label{rem:compressionGP}
It is worth noticing that the inequality given by Theorem \ref{thm:GPcompression} often turns out to be an equality. Indeed, it is clear that $\alpha_p^*(\Gamma \mathcal{G}) \leq \alpha_p^*(G)$ for every $G \in \mathcal{G}$ since a vertex-group is quasi-isometrically embedded; and a graph product often contains a quasi-isometrically embedded non abelian free subgroup (see \cite[Lemma 4.6]{AntolinDreesen}) so that $\alpha_p^*(\Gamma \mathcal{G}) \leq \alpha_p^*(\mathbb{F}_2)= \max(1/2,1/p)$. See \cite[Theorem 4.7]{AntolinDreesen} for a precise statement.
\end{remark}

\begin{remark}
By applying Proposition \ref{prop:compression}, we are also able to show that
$$\alpha_p( \Gamma \mathcal{G}) \geq \min \left( \frac{1}{p}, \min\limits_{G \in \mathcal{G}} \alpha_p(G) \right).$$
However, this lower bound is not optimal in general. In fact, it is proved in \cite[Corollary 5.8]{AntolinDreesen} that $\alpha_p(\Gamma \mathcal{G}) = \min\limits_{G \in \mathcal{G}} \alpha_p(G)$. 
\end{remark}

\noindent
As a conclusion of this section, it is worth noticing that, in the case of graph products of finite groups (along a finite simplicial graph), the actions on the associated quasi-median graphs are geometric. Although it follows from Proposition \ref{prop:quasimedianimplycubical} that such a graph product must also act geometrically (and virtually specially) on some CAT(0) cube complex, we really expect that the explicit structure of the quasi-median graphs involved will lead to considering graphs of finite groups as a source of convenient and interesting examples. Motivations can already be found in the study of right-angled Coxeter groups, see for instance \cite{BHSC,DaniThomas,AsDimRACG,HypRACGlargedim,GPsurfacesub,divCCC,HypBoundRACG,RACGboundSier} and \cite{DavisCovEucl,Davis1989,DavisJanusLafont,DavisOkun} for applications to topology and geometry. To stress out this idea, in addition of our combination results applied to this situation, let us mention the following result, which generalizes \cite[Theorem 1.1]{RACGproductoftrees} proved for right-angled Coxeter groups (with a different method, although a cubical proof following the lines of the proof below is possible).

\begin{thm}\label{thm:embedxtrees}
Let $\Gamma$ be a finite simplicial graph and $\mathcal{G}$ a collection of finite groups labelled by $V(\Gamma)$. There exits an equivariant quasi-isometric embedding from the graph product $\Gamma \mathcal{G}$ into the Cartesian product of $\chi(\Gamma)$ trees, where $\chi(\Gamma)$ denotes the chromatic number of $\Gamma$.
\end{thm}

\noindent
This theorem is an immediate consequence of the proposition proved below, which we believe to be well-known from specialists in the context of CAT(0) cube complexes. In this statement, we use the following definition: the \emph{crossing graph} $\Delta X$ of a CAT(0) cube complex $X$ (or more generally, of a quasi-median graph) is the graph whose vertices are the hyperplanes of $X$ and whose edges link two transverse hyperplanes. It is worth noticing that the crossing graph may not be connected (in fact, it is disconnected if and only if there exist cut vertices in the cube complex). 

\begin{prop}\label{prop:embedxtrees}
Let $G$ be a group acting geometrically on a CAT(0) cube complex $X$. Suppose that, for every hyperplane $J$ and every element $g \in G$, the hyperplanes $J$ and $gJ$ are not tranverse. There exists an equivariant quasi-isometric embedding from $X$ into the Cartesian product of $\chi(\Delta X/G)$ trees.
\end{prop}

\begin{proof}
A coloring of $\Delta X/G$ with $n=\chi(\Delta X/G)$ colors induces a coloring of $\Delta X$ with $n$ colors. Indeed, if $J_1$ and $J_2$ are two adjacent vertices of $\Delta X$, ie., two transverse hyperplanes of $X$, then by assumption no element of $G$ sends $J_1$ to $J_2$, so that their images in $\Delta X/G$ are distinct, and we conclude that $J_1$ and $J_2$ have different colors in $\Delta X$. For convenience, let $\mathfrak{J}$ denote the set of hyperplanes of $X$, and, for every $1 \leq k \leq n$, $\mathfrak{J}_k$ the subset of the hyperplanes colored using the $k$-th color. 

\medskip \noindent
Fix some $1 \leq k \leq n$, and let $T_k$ denote the CAT(0) cube complex obtained by cubulating the space with walls $(X, \mathfrak{J}_k)$. Notice that $G$ acts on $T_k$ because $\mathfrak{J}_k$ is $G$-invariant and that there exists a natural equivariant map $\phi_k : X \to T_k$, defined by sending every vertex of $X$ to its corresponding principal orientation, such that $d_{T_k}(\phi_k(x),\phi_k(y))$ is equal to the number of the hyperplanes separating $x$ and $y$ which belong to $\mathfrak{J}_k$, say $m_k(x,y)$, for every vertices $x,y \in X$. Moreover, because no two hyperplanes of $\mathfrak{J}_k$ are transverse, the CAT(0) cube complex $T_k$ does not contain two transverse hyperplane, ie., it is a tree.

\medskip \noindent
Set $T= T_1 \times \cdots \times T_n$ and $\phi = \phi_1 \times \cdots \times \phi_n : X \to T$. Of course, $G$ acts naturally on $T$. For every vertices $x,y \in X$, we have
$$d_T(\phi(x),\phi(y))= \sum\limits_{k=1}^n d_{T_k}(\phi_k(x),\phi_k(y)) = \sum\limits_{k=1}^n m_k(x,y) = d_X(x,y),$$
so $\phi$ is an isometric embedding from $X$ into $T$. Next, since the $\phi_k$'s are all $G$-invariant, we deduce that for every $g \in G$ and every $x,y \in X$,
$$\phi(g \cdot x)=(\phi_k(g \cdot x))= (g \cdot \phi_k(x))= g \cdot (\phi_k(x))= g \cdot \phi(x).$$
Therefore, the map $\phi$ defines an equivariant isometric embedding $X \hookrightarrow T$. On the other hand, we know from the Milnor-Svarc lemma that, for some fixed vertex $x_0 \in X$, the map $\varphi : g \mapsto g \cdot x_0$ defines a quasi-isometry $G \to X$. The map $\phi \circ \varphi : G \to T$ is the equivariant quasi-isometric embedding we are looking for.
\end{proof}

\begin{proof}[Proof of Theorem \ref{thm:embedxtrees}.]
According to Proposition \ref{prop:quasimedianimplycubical}, $\Gamma \mathcal{G}$ acts geometrically on the CAT(0) cube complex $C(X, \mathcal{SW})$ obtained by cubulating the space with walls $(X, \mathcal{SW})$. Moreover, according to Theorem \ref{thm:quasicubulation}, there exists an equivariant bijection from the hyperplanes of $\X$ and those of $C(X, \mathcal{SW})$ which respects transversality. As a consequence, since $\Gamma \mathcal{G}$ does not send a hyperplane of $\X$ to a hyperplane transverse to it, according to Lemma \ref{lem:labelhyp}, we deduce that a similar statement holds for the action of $\Gamma \mathcal{G}$ on $C(X, \mathcal{SW})$. Therefore, the only point to check in order to deduce Theorem \ref{thm:embedxtrees} from Proposition \ref{prop:embedxtrees} is that the quotient $\Delta / \Gamma \mathcal{G}$, where $\Delta$ denotes the crossing graph of $\X$, is isomorphic to $\Gamma$. Because any hyperplane of $\X$ is a translate of some $J_u$ where $u \in V(\Gamma)$, and that the hyperplanes $J_u$ and $J_v$ do not belong to the same $\Gamma \mathcal{G}$-orbit for every distinct vertices $u,v \in V(\Gamma)$, the map $\Delta \to \Gamma$ defined by $gJ_u \mapsto u$ induces a well-defined map $\Delta/ \Gamma \mathcal{G} \to \Gamma$. Moreover, we know from Lemma \ref{lem:labelhyp} that, for every vertices $u,v \in V(\Gamma)$, the hyperplanes $J_u$ and $J_v$ are transverse if and only if $u$ and $v$ are adjacent in $V(\Gamma)$. Thus, our map $\Delta / \Gamma \mathcal{G} \to \Gamma$ defines an isomorphism. 
\end{proof}

\subsection{When is a graph product hyperbolic?}

\noindent
In \cite{MeierGP}, Meier determined precisely when a graph product is hyperbolic. More precisely, he proved:

\begin{thm}\label{Meier}
Let $\Gamma$ be a finite simplicial graph and $\mathcal{G}$ a collection of non trivial groups indexed by $V(\Gamma)$. The graph product $\Gamma \mathcal{G}$ is hyperbolic if and only if 
\begin{itemize}
	\item every group of $\mathcal{G}$ is hyperbolic;
	\item no two infinite vertex-groups are adjacent in $\Gamma$;
	\item two vertex-groups adjacent to a common infinite vertex-group must be adjacent;
	\item the graph $\Gamma$ is square-free.
\end{itemize}
\end{thm}

\noindent
We will give two proofs of this result. The first one follows Meier's argument but replaces the CAT(0) cube complex he constructed with the quasi-median graph $\X$. In our opinion, this is the most simple and natural argument. The second proof is an application of Theorem \ref{thm:qmrelativelyhyp}. In this case, the conclusion is in fact stronger: you can remove the assumption that vertex-groups are hyperbolic and prove that the graph product is hyperbolic relatively to groups commensurable to vertex-groups. This observation will be generalized in the next section.

\begin{proof}[First proof of Theorem \ref{Meier}.]
Suppose that $\Gamma \mathcal{G}$ is hyperbolic. Notice that, fixing some vertex-group $G \in \mathcal{G}$, the map $\Gamma \mathcal{G} \to G$ defined by sending any vertex-group different from $G$ to the identity produces a retraction onto $G$. Because a retract in a hyperbolic group must be hyperbolic, we deduce that vertex-groups are necessarily hyperbolic. This proves the first condition. Next, if $G_1,G_2$ are two adjacent vertex-groups, the subgroup $\langle G_1,G_2 \rangle$ is isomorphic to $G_1 \times G_2$, so, since an infinite hyperbolic group contains an infinite order element and that a hyperbolic group cannot contain a subgroup isomorphic to $\mathbb{Z}^2$, we deduce that $G_1$ and $G_2$ cannot be both infinite. This proves the second condition. Now, if $G_1,G_2$ are two non adjacent vertex-groups both adjacent to a third vertex-group $G_3$, then the subgroup $\langle G_1,G_2,G_3 \rangle$ is isomorphic to $G_3 \times (G_1 \ast G_2)$. Once again, because a hyperbolic group does not contain any subgroup isomorphic to $\mathbb{Z}^2$, $G_3$ cannot be infinite, proving the third condition. Finally, a square is a bipartite complete graph $K_{2,2}$, so if $\Gamma$ contains an induced square then the subgroup it spans is isomorphic to $(G_1 \ast G_2) \times (G_3 \ast G_4)$ for some vertex-groups $G_1,G_2,G_3,G_4 \in \mathcal{G}$. But such a subgroup cannot exist since the hyperbolicity of $\Gamma \mathcal{G}$ implies that it does not contain any subgroup isomorphic to $\mathbb{Z}^2$.

\medskip \noindent
Conversely, suppose that our three conditions hold. Suppose that there exists an infinite vertex-group $G_u$ where $u \in V(\Gamma)$. Let $\Gamma_0$ denote the subgraph $\Gamma \backslash \{ u \}$. Then the graph product $\Gamma \mathcal{G}$ splits as the amalgamated product
$$\Gamma \mathcal{G}= \Gamma_0 \mathcal{G} \underset{\langle \mathrm{link}(u) \rangle}{\ast} \langle \mathrm{star}(u) \rangle.$$
Notice that, because $G_u$ is infinite, the vertex-groups adjacent to $G_u$ must be finite and pairwise adjacent, so that $\langle \mathrm{link}(u) \rangle$ is the direct product of finite groups. As a consequence, $\langle \mathrm{star}(u) \rangle = G_u \times \langle \mathrm{link}(u) \rangle$ contains $G_u$ as a finite-index subgroup; in particular, it must be hyperbolic. 

\medskip \noindent
By iterating the argument, we find a hierarchy for $\Gamma \mathcal{G}$. Precisely, there exists a sequence of groups $H_1, \ldots, H_n$ such that $H_1 = \Gamma_f \mathcal{G}$ where $\Gamma_f$ is the subgraph of $\Gamma$ generated by the vertices whose associated groups are finite, $H_n= \Gamma \mathcal{G}$, and, for every $1 \leq i \leq n-1$, $H_{i+1}$ decomposes as an amalgamated product $H_i \underset{F}{\ast} K$ where $F$ is a finite group and $K$ is hyperbolic. Because amalgamating two hyperbolic groups along a finite subgroup produces a hyperbolic group, we deduce that it is sufficient to show that $\Gamma_f \mathcal{G}$ is hyperbolic to prove that $\Gamma \mathcal{G}$ is hyperbolic.

\medskip \noindent
Notice that $\Gamma_f$ is square-free since $\Gamma$ is square-free itself. On the other hand, because any vertex-group of $\Gamma_f \mathcal{G}$ is finite, $X(\Gamma_f,\mathcal{G})$ is a Cayley graph of $\Gamma_f \mathcal{G}$ with respect to a finite generating set. Therefore, the conclusion follows from the following observation:

\begin{fact}\label{fact:hypiffsquarefree}
Let $\Gamma$ be simplicial graph and $\mathcal{G}$ a collection of non trivial groups indexed by the vertices of $\Gamma$. The graph $\X$ is hyperbolic if and only if $\Gamma$ is square-free and $\mathrm{clique}(\Gamma)<+ \infty$.
\end{fact}

\noindent
If $\dim_{\square}(\Gamma)= \mathrm{clique}(\Gamma) = + \infty$, then $\X$ cannot be hyperbolic, because the $1$-skeleton of $(2n+1)$-cube contains a triangle which is not $n$-thin. Next, if $\Gamma$ contains an induced square $(a,b,c,d)$, then, choosing some non trivial elements $g \in G_a$, $h \in G_b$, $k \in G_c$, $\ell \in G_d$, we find an $n$-thick flat square
$$\{ (gk)^i (hk)^j \mid 0 \leq i,j \leq n \} \subset \X$$
for arbitrarily large $n$. A fortiori, $\X$ cannot be hyperbolic.

\medskip \noindent
Conversely, suppose that $\X$ is not hyperbolic but that $\mathrm{clique}(\Gamma)<+ \infty$. According to Proposition \ref{prop:qmhyp}, there exists an arbitrarily thick flat rectangle in $\X$. In particular, there exists a flat square $R : [0,n] \times [0,n] \hookrightarrow \X$ where $n> \mathrm{clique}(\Gamma)$. According to Lemma \ref{lem:flatrectangle}, there exist a join subgraph $\Lambda_1 \ast \Lambda_2 \subset \Gamma$, an element $k \in \Gamma \mathcal{G}$ and syllables $g_1, \ldots , g_n \in \Lambda_1 \mathcal{G}$, $h_1, \ldots , h_n \in \Lambda_2 \mathcal{G}$ such that the products $g=g_1 \cdots g_n$ and $h=h_1 \cdots h_n$ are reduced, and such that (the image of) $R$ coincides with
$$\{ k \cdot g_1 \cdots g_i \cdot h_1 \cdots h_j \mid 0 \leq i,j \leq n \}.$$
Notice that, if the $g_i$'s pairwise commute, we deduce from the fact that the product $g_1 \cdots g_r$ is reduced the vertex-groups associated to the $g_i$'s must be pairwise distinct; of course, they are also pairwise adjacent, hence $n \leq \mathrm{clique}(\Gamma)$, which contradicts our assumption on $n$. Therefore, $\mathrm{supp}(g)$ must contain at least two non adjacent vertices, say $u_1,u_2 \in \Lambda_1$. Similarly, $\mathrm{supp}(h)$ must contain at least two non adjacent vertices, say $v_1, v_2 \in \Lambda_2$. Because any vertex of $\Lambda_1$ is adjacent to any vertex of $\Lambda_2$, we deduce that the vertices $u_1,u_2,v_1,v_2$ generate an induced square of $\Gamma$. This concludes the proof.
\end{proof}

\begin{proof}[Second proof of Theorem \ref{Meier}.]
We already proved in the previous proof that the hypotheses of our statement are necessary. Conversely, suppose that they are satisfied, and let us prove that $\Gamma \mathcal{G}$ is hyperbolic.

\medskip \noindent
Notice that, because two infinite vertex-groups cannot be adjacent, two hyperplanes of $\X$ dual to infinite cliques cannot be transverse: this is a consequence of Lemma \ref{lem:labelhyp}. Next, let $J$ be a hyperplane of $\X$ dual to an infinite clique; say $J$ is labelled by the vertex $u \in V(\Gamma)$ such that the vertex-group $G_u$ is infinite. Up to translating by an element of $\Gamma \mathcal{G}$, we may suppose without loss of generality that $J=J_u$. According to Corollary \ref{cor:hyperplane}, the carrier of $J$ is $\langle \mathrm{star}(u) \rangle = \langle \mathrm{link}(u) \rangle \times G_u$. On the other hand, because the link of $u$ is complete by assumption, $\langle \mathrm{link}(u) \rangle$ is just a prism, so that $N(J)$ is just a prism as well. In particular, $N(J)$ is cubically finite. As a consequence of Corollary \ref{cor:beinghyperbolic}, which applies thanks to Proposition \ref{prop:GPsumup} and our previous observations, we deduce that $\Gamma \mathcal{G}$ is hyperbolic.
\end{proof}

\subsection{When is a graph product relatively hyperbolic?}

\noindent
In this section, our goal is to characterize relatively hyperbolic graph products. This fills in a gap between the characterisation of hyperbolic graph products, due to Meier \cite{MeierGP} and discussed in the previous section, and the characterisation of acylindrically hyperbolic graph products, due to Minasyan and Osin \cite{arXiv:1310.6289} (see also the next section for another point of view on this subject). We essentially follow the proof of the characterisation of relatively hyperbolic right-angled Coxeter groups we wrote in \cite{coningoff}. Before stating our main result, we need to introduce some definitions.

\medskip \noindent
Given a finite simplicial graph $\Gamma$ and a collection of groups $\mathcal{G}$ labelled by $V(\Gamma)$, we will say that a subgraph $\Lambda \leq \Gamma$ is \emph{vast} if the subgroup of $\Gamma \mathcal{G}$ generated by the vertex-groups corresponding to the vertices of $\Lambda$, ie., $\Lambda \mathcal{G}$, is infinite; otherwise, $\Lambda$ is said \emph{narrow}. Notice that a subgraph is narrow if and only if it is complete and all the vertex-groups labellings its vertices are finite. A join $\Lambda_1 \ast \Lambda_2 \leq \Gamma$ is \emph{large} if both $\Lambda_1$ and $\Lambda_2$ are vast. 

\begin{definition}
Let $\Gamma$ be a finite simplicial graph and $\mathcal{G}$ a collection of groups labelled by $V(\Gamma)$. For every subgraph $\Lambda \subset \Gamma$, let $\mathrm{cp}(\Lambda)$ denote the subgraph of $\Gamma$ generated by $\Lambda$ and the vertices $v \in \Gamma$ such that $\mathrm{link}(v) \cap \Lambda$ is vast. Now, define the collection of subgraphs $\mathfrak{J}^n(\Gamma)$ of $\Gamma$ by induction in the following way:
\begin{itemize}
	\item $\mathfrak{J}^0(\Gamma)$ is the collection of all the large joins in $\Gamma$;
	\item if $C_1, \ldots, C_k$ denote the connected components of the graph whose set of vertices is $\mathfrak{J}^n(\Gamma)$ and whose edges link two subgraphs with vast intersection, we set $\mathfrak{J}^{n+1}(\Gamma) = \left( \mathrm{cp} \left( \bigcup\limits_{\Lambda \in C_1} \Lambda \right), \ldots, \mathrm{cp} \left( \bigcup\limits_{\Lambda \in C_k} \Lambda \right) \right)$.
\end{itemize}
Because $\Gamma$ is finite, the sequence $(\mathfrak{J}^n(\Gamma))$ must eventually be constant and equal to some collection $\mathfrak{J}^{\infty}(\Gamma)$. Finally, let $\mathfrak{J}(\Gamma)$ denote the collection of subgraphs of $\Gamma$ obtained from $\mathfrak{J}^{\infty}(\Gamma)$ by adding the singletons corresponding to the vertices of $\Gamma \backslash \bigcup \mathfrak{J}^{\infty}(\Gamma)$.  
\end{definition}

\begin{thm}\label{thm:relativehyp}
Let $\Gamma$ be a finite simplicial graph not reduced to a single vertex and $\mathcal{G}$ a collection of finitely generated groups labelled by $V(\Gamma)$. The graph product $\Gamma \mathcal{G}$ is relatively hyperbolic if and only if $\mathfrak{J}(\Gamma) \neq \{ \Gamma \}$. If so, $\Gamma \mathcal{G}$ is hyperbolic relatively to $\{ \Lambda \mathcal{G} \mid \Lambda \in \mathfrak{J}(\Gamma) \}$. 
\end{thm}

\noindent
Following \cite{coningoff}, we first introduce \emph{join decompositions} of $\Gamma$, generalizing the decomposition $\mathfrak{J}(\Gamma)$ we constructed above, and prove that they lead to some relative hyperbolicity of $\Gamma \mathcal{G}$.

\begin{definition}
Let $\Gamma$ be a finite simplicial graph and $\mathcal{G}$ a collection of groups labelled by $V(\Gamma)$. A \emph{join decomposition} of $\Gamma$ is a collection of subgraphs $(\Gamma_1, \ldots , \Gamma_n)$, with $n \geq 0$, such that:
\begin{itemize}
	\item any large join of $\Gamma$ is included into $\Gamma_i$ for some $1 \leq i \leq n$;
	\item $\Gamma_i \cap \Gamma_j$ is narrow for every $1 \leq i < j \leq n$;
	\item for every vertex $v \in \Gamma$, $link(v) \cap \Gamma_i$ vast implies $v \in \Gamma_i$;
	\item $\Gamma= \Gamma_1 \cup \cdots \cup \Gamma_n$.
\end{itemize}
\end{definition}

\noindent
There exists at least one join decomposition: the \emph{trivial decomposition} $\{ \Gamma \}$.

\begin{prop}\label{prop:relativehyp}
Let $\Gamma$ be a finite simplicial graph and $\mathcal{G}$ a collection of finitely generated groups labelled by $V(\Gamma)$. If $\{ \Gamma_1, \ldots, \Gamma_n \}$ is a join decomposition of $\Gamma$, then the graph product $\Gamma \mathcal{G}$ is hyperbolic relative to $\{ \Gamma_1 \mathcal{G} , \ldots, \Gamma_n \mathcal{G} \}$. 
\end{prop}

\noindent
In order to prove this relative hyperbolicity, we will apply Definition \ref{def:relativehyp}, ie., we will construct a fine hyperbolic graph acted upon by our graph product. We already proved a criterion to get fine cone-offs, namely Proposition \ref{prop:criterionfine}, and now, before proving Proposition \ref{prop:relativehyp}, we need a criterion to get hyperbolic cone-offs. 

\begin{prop}\label{prop:hypconeoff}
Let $X$ be a quasi-median graph of finite cubical dimension and $\mathcal{Q}$ a collection of convex subgraphs. Let $Y$ denote the cone-off of $X$ over $\mathcal{Q}$. If there exist constants $L,C \geq 0$ such that any $L$-thick flat rectangle of $X$ has diameter at most $C$ in $Y$, then $Y$ is $\delta$-hyperbolic for some $\delta$ depending only on $L$, $C$ and $\dim_{\square}(X)$.
\end{prop}

\begin{proof}
For every vertices $x,y \in Y$, let $\eta(x,y)$ denote the subgraph generated in $Y$ by the interval $I(x,y)$ in $X$. Now, we want to apply the following criterion, due to Bowditch \cite[Proposition 3.1]{Bowditchcriterion}:

\begin{prop}
Let $T$ be a graph and $D \geq 0$. Suppose that a connected subgraph $\eta(x, y)$, containing $x$ and $y$, is associated to any pair of vertices $\{ x, y \} \in T^2$ such that:
\begin{itemize}
	\item for any vertices $x,y \in T$, $d(x,y) \leq 1$ implies $\mathrm{diam}~ \eta(x,y) \leq D$;
	\item for any vertices $x,y,z \in T$, we have $\eta(x,y) \subset ( \eta(x,z) \cup \eta(z,y) )^{+D}$.
\end{itemize}
Then $T$ is $\delta$-hyperbolic for some $\delta$ depending only on $D$.
\end{prop}

\noindent
Let us verify the first condition. So let $x,y \in Y$ be two vertices satisfying $d_Y(x,y) \leq 1$. If $x=y$ then $\eta(x,y)= \{x \}$, hence $\mathrm{diam}_Y \eta(x,y)=0$. Otherwise, two cases may happen. Either $x$ and $y$ are adjacent in $X$, hence $\mathrm{diam}_Y \eta(x,y) \leq 1$; or there exists some $Q \in \mathcal{Q}$ such that $x,y \in Q$, so that $\eta(x,y) \subset Q$ because $Q$ is a convex subgraph of $X$, and we conclude that $\mathrm{diam}_Y \eta(x,y) \leq 1$. 

\medskip \noindent
Now, we focus on the second condition. Let $x,y,z \in Y$ be three vertices, and fix some vertex $p \in \eta(x,y)$, ie., there exists a geodesic $[x,y]$ (in $X$) between $x$ and $y$ passing through $p$. We want to prove that $p$ is contained in the $(\max(2L,C) + \dim_{\square}(X))$-neighborhood of $\eta(x,z) \cup \eta(y,z)$. Let $(x',y',z')$ denote the quasi-median of the triple $(x,y,z)$, and fix six geodesics $[x,x']$, $[y,y']$, $[z,z']$, $[x',y']$, $[y',z']$ and $[x',z']$ (in $X$). Notice that the two geodesics $[x,y]$ and $[x,x'] \cup [x',y'] \cup [y',y]$ define a bigon in $X$. By applying Lemma \ref{lem:bigonflatrectangle}, we know that there exists a flat square $R : [0,n] \times [0,n] \hookrightarrow X$ such that $(0,0)=p$ and $p'=(n,n) \in [x,x'] \cup [x',y'] \cup [y',y]$. Notice that, if $n \geq L$, necessarily (the image of) $R$ has diameter at most $C$ in $Y$, hence $d_Y(p,p') \leq \max(2L,C)$. If $p' \in [x,x'] \cup [y',y]$, then
$$d_Y(p, \eta(x,y) \cup \eta(y,z)) \leq d_Y(p,p') \leq \max(2L,C)$$
since $[x,x'] \cup [y',y] \subset \eta(x,z) \cup \eta(z,y)$. From now on, suppose that $p' \in [x',y']$. According to Proposition \ref{prop:quasimedian}, the gated hull of $\{x',y',z'\}$ is a prism $P$. By noticing that
$$\mathrm{diam}_Y(P) \leq \mathrm{diam}_X(P) \leq \dim_{\square}(X),$$
we deduce that there exists $p'' \in [x',z']$ such that $d_Y(p',p'') \leq \dim_{\square}(X)$. Since $p'' \in I(x,z) \subset \eta(x,z)$, we conclude that
$$d_Y(p, \eta(x,z)) \leq d_Y(p,p'') \leq d_Y(p,p')+d_Y(p',p'') \leq \max(2L,C) + \dim_{\square}(X).$$
This concludes the proof.
\end{proof}

\begin{proof}[Proof of Proposition \ref{prop:relativehyp}.]
According to Proposition \ref{prop:GPsumup}, the action $\Gamma \mathcal{G} \curvearrowright \X$ is topical-transitive. Fix a collection of cliques $\mathcal{C}$ such that any $G$-orbit of hyperplanes intersects it along a single clique, and recall that $\mathcal{C}= \mathcal{C}_1$. For every $C \in \mathcal{C}$, let $\delta_C$ be a word metric on the associated vertex-group (which is finitely generated by assumption). According to Proposition \ref{prop:mtopicalI}, our collection $\{(C,\delta_C) \mid C \in \mathcal{C} \}$ extends to a coherent $\Gamma \mathcal{G}$-invariant system of metrics. We will think of $(\X,\delta)$ as a graph (see the remark preceding Proposition \ref{prop:criterionfine}), which we denote by $X$ for convenience. Setting $\mathcal{Q}= \{ g \cdot \Gamma_i \mathcal{G} \mid g \in \Gamma \mathcal{G}, 1 \leq i \leq n\}$, let $Y$ denote the usual cone-off of $X$ over $\mathcal{Q}$. We claim that the hypotheses of Definition \ref{def:relativehyp} are satisfied. Because the stabiliser of any vertex of $X$ is trivial, the only point to verify is that $Y$ is a fine hyperbolic graph. This will be done thanks to Proposotion \ref{prop:criterionfine} and Proposition \ref{prop:hypconeoff}.

\medskip \noindent
In order to show that $Y$ is hyperbolic, it is sufficient to prove that the cone-off $Z$ of $X$ over $\mathcal{Q}$ is hyperbolic, since taking either the cone-off or the usual cone-off does not disturb the quasi-isometry class of the graph we obtain. Because our collection of subgraphs of $\Gamma$ covers $V(\Gamma)$, $Z$ can also be obtained from $\X$ by coning-off $\mathcal{Q}$. On the other hand, if $R$ is a $(\mathrm{clique}(\Gamma)+1)$-thick flat rectangle of $\X$, it follows from Lemma \ref{lem:flatrectangle} that there exists a large join $\Lambda= \Lambda_1 \ast \Lambda_2 \subset \Gamma$ such that $R \subset g \cdot \Lambda \mathcal{G}$ for some $g \in \Gamma \mathcal{G}$. By definition of join decompositions, there must exist some $1 \leq i \leq n$ such that $\Lambda \subset \Gamma_i$. From
$$R \subset g \cdot \Lambda \subset g \cdot \Gamma_i \in \mathcal{Q},$$
we deduce that (the image of) $R$ in $Z$ has diameter at most $1$. Consequently, $Z$ (and a fortiori $Y$) is hyperbolic according to Proposition \ref{prop:hypconeoff}.

\medskip \noindent
Next, let us verify the hypotheses of Proposition \ref{prop:criterionfine}.
\begin{itemize}
	\item Every $Q \in \mathcal{Q}$ is convex in $X$. This follows from Corollary \ref{cor:systgeodesic} and from the observation that every $Q \in \mathcal{Q}$ is gated in $\X$ (as a consequence of Proposition \ref{prop:gated}).
	\item $X$ is locally finite. This follows from Lemma \ref{lem:whenlocallyfinite} because every vertex of $\X$ belongs to at most $\# V(\Gamma)$ cliques, the local metrics are locally finite since they are word metrics defined on finitely generated groups, and for the same reason the local metrics are uniformly discrete.
	\item $\mathcal{Q}$ is locally finite in $X$. Indeed, elements of $\mathcal{Q}$ are cosets of $\Gamma_i \mathcal{G}$'s, and because two cosets of the same subgroup must be either equal or disjoint, a collection of pairwise intersecting elements of $\mathcal{Q}$ in $X$ is necessarily of the form $\{g_i \cdot \Gamma_i \mathcal{G} \mid i \in I \}$ where $I \subset \{ 1, \ldots, n \}$ and $g_i \in \Gamma \mathcal{G}$ for all $i \in I$. Therefore, at most $n$ elements of $\mathcal{Q}$ may pairwise intersect.
\end{itemize}
Finally, in order to conclude the proof, we will show that, given any two distinct elements $Q_1, Q_2 \in \mathcal{Q}$, at most $F \cdot \mathrm{clique}(\Gamma)$ sector-walls intersect both $Q_1$ and $Q_2$ (in $\X$), where $F= \max \{ \# G \mid G \in \mathcal{G} \ \text{finite} \}$. 

\medskip \noindent
First, we claim that a hyperplane (of $\X$) intersecting both $Q_1$ and $Q_2$ cannot be labelled by an infinite vertex-group. Suppose by contradiction that such a hyperplane exists, say a hyperplane $J$ which is labelled by the infinite vertex-group $G_u$. Because the edges dual to $J$ are all labelled by $u$, necessarily $Q_1$ and $Q_2$ contain edges labelled by $u$, so that, if we write $Q_1 = g_1 \cdot \Lambda_1$ and $Q_2=g_2 \cdot \Lambda_2$ for some $g_1,g_2 \in \Gamma \mathcal{G}$ and $\Lambda_1, \Lambda_2 \in \{ \Gamma_1, \ldots, \Gamma_n \}$, then $u \in V(\Lambda_1) \cap V(\Lambda_2)$. Since $G_u$ is infinite, it follows that the intersection $\Lambda_1 \cap \Lambda_2$ is vast, hence $\Lambda_1= \Lambda_2$; let $\Lambda$ denote this common subgraph. Notice that, because $u \in V(\Lambda)$ and because $G_u$ is infinite, necessarily $\mathrm{star}(u) \subset \Lambda$. As a consequence of Corollary \ref{cor:hyperplane}, we deduce that
$$N(J) = g_1 \cdot \langle \mathrm{star}(u) \rangle \subset g_1 \cdot \Lambda \mathcal{G}=Q_1,$$
and similarly $N(J) \subset Q_2$. Therefore, $N(J) \subset Q_1 \cap Q_2$. Finally, since two cosets of the same subgroup must be either equal or disjoint, we conclude that $Q_1=Q_2$ holds, contradicting our initial assumption. This concludes the proof of our first claim.

\medskip \noindent
Next, we claim that at most $\mathrm{clique}(\Gamma)$ hyperplanes (of $\X$) intersect both $Q_1$ and $Q_2$. Suppose by contradiction that at least $\mathrm{clique}(\Gamma)+1$ hyperplanes intersect both $Q_1$ and $Q_2$. By combining Corollary \ref{cor:projflat} and Lemma \ref{lem:flatrectangle}, we deduce that there exist a join subgraph $\Xi_1 \ast \Xi_2 \subset \Gamma$, an element $k \in \Gamma \mathcal{G}$ and syllables $g_1, \ldots , g_r \in \Xi_1 \mathcal{G}$, $h_1, \ldots , h_s \in \Xi_2 \mathcal{G}$ such that the products $g=g_1 \cdots g_r$ and $h=h_1 \cdots h_s$ are reduced, such that $kg_1 \cdots g_i \in Q_1$ and $khg_1 \cdots g_i \in Q_2$ for every $1 \leq i \leq r$, and such that $r \geq \mathrm{clique}(\Gamma)+1$. Notice that, if the $g_i$'s pairwise commute, we deduce from the fact that the product $g_1 \cdots g_r$ is reduced that the vertex-groups associated to the $g_i$'s must be pairwise distinct; of course, they are also pairwise adjacent, hence $r \leq \mathrm{clique}(\Gamma)$, which contradicts our assumption on $r$. Therefore, $\mathrm{supp}(g)$ must contain two non adjacent vertices. On the other hand, because $g$ labels paths in $Q_1$ and $Q_2$, if we write $Q_1 = g_1 \cdot \Lambda_1$ and $Q_2=g_2 \cdot \Lambda_2$ for some $g_1,g_2 \in \Gamma \mathcal{G}$ and $\Lambda_1, \Lambda_2 \in \{ \Gamma_1, \ldots, \Gamma_n \}$, then $\mathrm{supp}(g) \subset V(\Lambda_1) \cap V(\Lambda_2)$ must hold. A fortiori, the intersection $\Lambda_1 \cap \Lambda_2$ is vast, hence $\Lambda_1= \Lambda_2$; let $\Lambda$ denote this common subgraph. Now, for every $1 \leq i \leq s$, notice that 
$$\mathrm{link}(h_i) \cap \Lambda \supset \Xi_1 \cap \Lambda \supset \mathrm{supp}(g);$$
therefore, because $\mathrm{supp}(g)$ contains two non adjacent vertices, it follows that the intersection $\mathrm{link}(h_i) \cap \Lambda$ is vast, hence $h_i \in \Lambda$. A fortiori, $h$ must belong to $\Lambda$. As a consequence, $khg \in Q_1 \cap Q_2$, and because two cosets of the same subgroup must be either equal or disjoint, we conclude that $Q_1=Q_2$, contradicting our initial assumption. This concludes the proof of our second claim.

\medskip \noindent
Finally, we are able to bound the number of sector-walls intersecting both $Q_1$ and $Q_2$. From our previous two claims, we know that there exist at most $\mathrm{clique}(\Gamma)$ hyperplanes intersecting both $Q_1$ and $Q_2$, and that each of these hyperplanes are labelled by a finite group; a fortiori, such a hyperplane delimits at most $F$ sector-walls. Therefore, at most $F \cdot \mathrm{clique}(\Gamma)$ sector-walls may intersect both $Q_1$ and $Q_2$.
\end{proof}

\noindent
It is clear from its definition that $\mathfrak{J}(\Gamma)$ is a join decomposition. Furthermore, according to the next lemma, it may be thought of as the minimal join decomposition of $\Gamma$.

\begin{proof}[Proof of Theorem \ref{thm:relativehyp}.]
Suppose that $\mathfrak{J}(\Gamma) \neq \{ \Gamma \}$. Because $\Gamma$ is not reduced to a single vertex, necessarily $\Lambda \mathcal{G}$ is a proper subgroup of $\Gamma \mathcal{G}$ for every $\Lambda \in \mathfrak{J}(\Gamma) \backslash J^{\infty}(\Gamma)$. Next, from the construction of $\mathfrak{J}^{\infty}(\Gamma)$, it follows by induction that any element of $\mathfrak{J}^{\infty}(\Gamma)$ is vast. Thus, if $\Lambda, \Gamma \in \mathfrak{J}^{\infty}(\Gamma)$, because the intersection $\Lambda \cap \Gamma= \Lambda$ is vast, necessarily $\Lambda= \Gamma$; as a consequence, if $\Gamma \in \mathfrak{J}^{\infty}(\Gamma)$ then $\mathfrak{J}^{\infty}(\Gamma)= \{ \Gamma \}$, so that $\mathfrak{J}(\Gamma)= \{ \Gamma \}$. Since we supposed that $\mathfrak{J}(\Gamma) \neq \{ \Gamma \}$, we conclude that $\Lambda \mathcal{G}$ is a proper subgroup of $\Gamma\mathcal{G}$ for every $\Lambda \in \mathfrak{J}^{\infty}(\Gamma)$. Consequently, if $\mathfrak{J}(\Gamma) \neq \{ \Gamma \}$, we deduce from Proposition \ref{prop:relativehyp} that $\Gamma \mathcal{G}$ is hyperbolic relatively to a finite collection of proper subgroups.

\medskip \noindent
Conversely, suppose that $\mathfrak{J}(\Gamma)= \{ \Gamma \}$. Because $\Gamma$ is not reduced to a single vertex, this implies that $\mathfrak{J}^{\infty}(\Gamma)= \{ \Gamma \}$. We want to prove that, if $\Gamma \mathcal{G}$ is hyperbolic relatively to some finite collection of subgroups $\mathcal{H}$, then $\Gamma \mathcal{G} \in \mathcal{H}$, which implies that $\Gamma \mathcal{G}$ is not relatively hyperbolic. This is a direct consequence of the following claim.

\medskip \noindent
Independently of the assumption $\mathfrak{J}(\Gamma)= \{ \Gamma \}$, we claim that, whenever $\Gamma \mathcal{G}$ is hyperbolic relatively to a finite collection of subgroups $\mathcal{H}$, for every $n \geq 0$ and every $\Lambda \in \mathfrak{J}^n(\Gamma)$ there exists some $H \in \mathcal{H}$ such that $\Lambda \mathcal{G} \subset H$. Because any subgroup isomorphic to a direct product of two infinite groups has to be included into a peripheral subgroup (see for instance \cite[Theorems 4.16 and 4.19]{OsinRelativeHyp}), the statement holds for $n=0$. Now suppose that this statement holds for some $n \geq 0$, and let $\Lambda \in \mathfrak{J}^{n+1}(\Gamma)$. The subgraph $\Lambda$ corresponds to a connected component $C=\{\Lambda_1, \ldots, \Lambda_k\}$ of the graph whose set of vertices is $\mathfrak{J}^n(\Gamma)$ and whose edges link two subgraphs with a vast intersection, ie., $\Lambda= \mathrm{cp}(\Lambda_1 \cup \cdots \cup \Lambda_k)$. By our induction hypothesis, for every $1 \leq i \leq k$, the subgroup $\langle \Lambda_i \rangle$ is included into some peripheral subgroup $H_i$. Notice that, for every $1 \leq i <j \leq k$, $H_i \cap H_j$ contains $\langle \Lambda_i \rangle \cap \langle \Lambda_j \rangle = \langle \Lambda_i \cap \Lambda_j \rangle$, which is infinite, hence $H_i=H_j$ because peripheral subgroups define an almost malnormal collection (see for instance \cite[Theorems 1.4 and 1.5]{OsinRelativeHyp}). Therefore, the $\langle \Lambda_i \rangle$'s are all included into the same peripheral subgroup $H$. Then, for every vertex $v \in \mathrm{cp}(\Lambda_1 \cup \cdots \cup \Lambda_k) \backslash (\Lambda_1 \cup \cdots \cup \Lambda_k)$, it follows from the definition of $\mathrm{cp}(\cdot)$ that the intersection 
$$\langle \Lambda_1 \cup \cdots \cup \Lambda_k \rangle \cap \langle \Lambda_1 \cup \cdots \cup \Lambda_k \rangle^v$$ 
is infinite, so that we deduce from
$$H \cap H^v \supset \langle \Lambda_1 \cup \cdots \cup \Lambda_k \rangle \cap \langle \Lambda_1 \cup \cdots \cup \Lambda_k \rangle^v$$
that $H \cap H^v$ is infinite, and finally that $v \in H$ since once again peripheral subgroups define an almost malnormal collection. Therefore,
$$\langle \Lambda \rangle = \langle \mathrm{cp}(\Lambda_1 \cup \cdots \cup \Lambda_k) \rangle \subset H.$$
This concludes the proof of our claim.
\end{proof}

\begin{remark}
The collection of peripheral subgroups provided by Theorem \ref{thm:relativehyp} is not minimal in general. For instance, if $\Gamma$ is the graph containing two vertices and no edges, and if $\mathcal{G}$ are two free groups, then $\Gamma \mathcal{G}$ is the free product of two free groups. Thus, $\{ \{1 \} \}$ is the minimal collection of peripheral subgroups of $\Gamma \mathcal{G}$ whereas we deduce by applying Theorem \ref{thm:relativehyp} that $\Gamma \mathcal{G}$ is hyperbolic relatively to the two free factors. Nevertheless, if we write $$\mathfrak{J}(\Gamma) = \mathfrak{J}^{\infty}(\Gamma) \sqcup \{ \{ v_i \} \mid 1 \leq i \leq m \},$$ and if $\{H_1^i,\ldots, H_{n(i)}^i \}$ is a minimal collection of peripheral subgroups of $G_{v_i}$ for every $1 \leq i \leq m$, then $$\{ \Lambda \mathcal{G} \mid \Lambda \in \mathfrak{J}^{\infty} \} \sqcup \bigsqcup\limits_{i=1}^m \{ H_1^i, \ldots, H_{n(i)}^i \}$$ is a minimal collection of peripheral subgroups of the graph product $\Gamma \mathcal{G}$. 
\end{remark}

\subsection{Embedding graph products}\label{section:GPcurvegraph}

\noindent
Motivated by the curve graph of a surface, several hyperbolic graphs were introduced to study various classes of groups, including, in the context of CAT(0) cube complexes, the contact graph \cite{MR3217625}, the contracting graph \cite{coningoff} and the crossing graph \cite{Roller} (see also \cite{SplittingObstruction}); and in other contexts, \cite{extensiongraph} for right-angled Artin groups and \cite{GarsideCurve} for Garside groups. In \cite{embeddingRAAG}, Kim and Koberda define the \emph{extension graph} of a right-angled Artin group, and in \cite{extensiongraph}, they motivate the analogy with the curve graph of a surface through a variety of results. This graph turns out to be precisely the \emph{crossing graph} of the quasi-median graph associated to the right-angled Artin group viewed as a graph product, ie., the graph whose vertices are the hyperplanes of the quasi-median graph and whose edges link two transverse hyperplanes. Therefore, Kim and Koberda's analogy can be extended to arbitrary graph products, in which the underlying quasi-median graph gets a crucial role: it is the surface on which the mapping class group acts. We refer to Section \ref{section:open} in which we elaborate this idea.

\medskip \noindent
The main result of this section, Theorem \ref{thm:embeddinggraphproduct} below, generalizes \cite[Theorem 1.3]{embeddingRAAG} and stresses out our analogy.

\begin{definition}
Let $X$ be a quasi-median graph. The \emph{crossing graph}\index{Crossing graphs} of $X$ is the graph, denoted by $\Delta X$, whose vertices are the hyperplanes of $X$ and whose edges link two transverse hyperplanes. 
\end{definition}

\begin{definition}
Let $G$ be a group acting on a quasi-median graph $X$ and $J$ a hyperplane of $X$. The \emph{rotative stabiliser}\index{Rotative stabilisers} of $J$ is
$$\mathrm{stab}_{\circlearrowleft}(J)= \bigcap \{ \mathrm{stab}(C) \mid \text{$C$ clique of $J$} \}.$$
\end{definition}

\begin{thm}\label{thm:embeddinggraphproduct}
Let $G$ be a group acting on a quasi-median graph $X$ with trivial vertex-stabilisers. For every hyperplane $J$, choose a residually finite subgroup $H_J$ of its rotative stabiliser. If $\Gamma$ is a finite subgraph of the crossing graph of $X$, there exists a collection $\mathcal{G}$ of finite-index subgroups of our $H_J$'s such that the graph product $\Gamma \mathcal{G}$ embeds into $G$. 
\end{thm}

\noindent
The idea of the proof is to play ping-pong with the rotative stabilisers of the hyperplanes corresponding to the vertices of our finite subgraph. The statement below generalizes the ping-pong lemma known for right-angled Artin groups to arbitrary graph products (see \cite{pingponglemmas} and references therein). The proof is completely similar.

\begin{prop}\label{prop:pingpong}
Let $G$ be a group acting on a set $X$, $\Gamma$ a simplicial graph and $\mathcal{H}= \{H_v \mid v \in V(\Gamma) \}$ a collection of subgroups. Suppose that $\bigcup\limits_{v \in V(\Gamma)} H_v$ generates $G$ and that $g$ and $h$ commute for every $g \in H_u$ and $h \in H_v$ if $u$ and $v$ are two adjacent vertices of $\Gamma$. Next, suppose that there exist a collection $\{ X_v \mid v \in V(\Gamma) \}$ of subsets of $X$ and a point $x_0 \in X \backslash \bigcup\limits_{v \in V(\Gamma)} X_v$ satisfying:
\begin{itemize}
	\item if $u,v \in V(\Gamma)$ are adjacent, then $g \cdot X_u \subset X_u$ for every $g \in H_v \backslash \{ 1 \}$;
	\item if $u,v \in V(\Gamma)$ are not adjacent and distinct, then $g \cdot X_u \subset X_v$ for every $g \in H_v \backslash \{ 1 \}$;
	\item for every $u \in V(\Gamma)$ and $g \in H_u \backslash \{ 1 \}$, $g \cdot x_0 \in X_u$.
\end{itemize}
Then $G$ is isomorphic to the graph product $\Gamma \mathcal{H}$. 
\end{prop}

\begin{proof}
By our assumptions on the subgroups of $\mathcal{H}$, we deduce that there exists a natural surjective morphism $\Gamma \mathcal{H} \twoheadrightarrow G$. In order to show that this morphism is also injective, we want to prove that, for any non empty reduced word $w$ of $\Gamma \mathcal{H}$, thought of as an element of $G$, $w \cdot x_0 \in X_u$ where $u$ is a vertex of $\Gamma$ which belongs to the support of the head of $w$; notice that, since $x_0 \notin \bigcup\limits_{v \in V(\Gamma)} X_v$ by assumption, this implies that $w \cdot x_0 \neq x_0$. This is sufficient to conclude the proof of our proposition. In particular, this will prove the following fact:

\begin{fact}\label{fact:pingpong}
For every non trivial $g \in G$, we know that $g \cdot x_0 \in \bigcup\limits_{u \in V(\Gamma)} X_u$.
\end{fact}

\noindent
We argue by induction on the length of $w$. If $w$ has length one, then $w \in H_u \backslash \{ 1 \}$ for some $u \in V(\Gamma)$. Our third assumption implies $w \cdot x_0 \in X_u$. Next, suppose that $w$ has length at least two. Write $w=gw'$ where $g$ is the first syllable of $w$, and $w'$ the rest of the word. Say $g \in H_u \backslash \{ 1 \}$. We know from our induction hypothesis that $w' \cdot x_0 \in X_v$ where $v$ is a vertex of $\Gamma$ which belongs to the support of the head of $w'$. Notice that $u \neq v$ since otherwise the word $gw'$ would not be reduced. Two cases may happen: either $u$ and $v$ are not adjacent, so that our second assumption implies that $w \cdot x_0 \in g \cdot X_v \subset X_u$; or $u$ and $v$ are adjacent, so that our first assumption implies that $w \cdot x_0 \in g \cdot X_v \subset X_v$. It is worth noticing that, in the former case, $u$ clearly belongs to the support of the head of $w$ since $g$ belongs to the head of $w$; in the latter case, if we write $w'=hw''$ where $h$ is a syllable of the head of $w'$ which belongs to $H_v$, then
$$w = gw'= g hw'' = hgw'',$$
so $h$ also belongs to the head of $w$, and a fortiori $v$ belongs to the support of the head of $w$. This concludes the proof.
\end{proof}

\noindent
In order to apply Proposition \ref{prop:pingpong} in the proof of Theorem \ref{thm:embeddinggraphproduct}, the following two preliminary lemmas will be necessary.

\begin{lemma}\label{lem:rotativestabcommute}
Let $G$ be a group acting on a quasi-median graph $X$. Suppose the action free on the vertices. If $J_1$ and $J_2$ are two transverse hyperplanes, then $g$ and $h$ commute for every $g \in \mathrm{stab}_{\circlearrowleft}(J_1)$ and $h \in \mathrm{stab}_{\circlearrowleft}(J_2)$.
\end{lemma}

\begin{proof}
According to Proposition \ref{prop:transversehypcube}, there exists a prism $P=C_1 \times C_2$ such that the cliques $C_1,C_2$ are dual to the hyperplanes $J_1,J_2$ respectively. Let $x \in X$ satisfy $C_1 \cap C_2 = \{ x \}$. Notice that, because $g \in \mathrm{stab}_{\circlearrowleft}(J_1) \subset \mathrm{stab}(C_1)$, necessarily $g \cdot x \in C_1$; similarly, $h \cdot x \in C_2$. In the same way, $g$ stabilises $C_1 \times \{ g \cdot x \}$, hence $gh \cdot x \in C_1 \times \{ g \cdot x \}$. Of course, since $x$ and $h \cdot x$ are adjacent vertices, $g \cdot x$ and $gh \cdot x$ must be adjacent as well, so that $x$, $g \cdot x$, $h \cdot x$ and $gh \cdot x$ define a square in $X$. Moreover, this square is induced since otherwise we would be able to deduce that $J_1 = J_2$. The argument implies that $x$, $h \cdot x$, $g \cdot x$ and $hg \cdot x$ define an induced square in $X$. If $hg \cdot x \neq gh \cdot x$, these two squares define a subgraph isomorphic to $K_{2,3}$. Because $X$ is quasi-median, this subgraph cannot be induced, so, since our two squares are induced, the vertices $gh \cdot x$ and $hg \cdot x$ must be adjacent. We get a contradiction by noticing that the four vertices $h \cdot x$, $gh \cdot x$, $g \cdot x$ and $hg \cdot x$ define an induced subgraph isomorphic to $K_4^-$. Therefore, $hg \cdot x = gh \cdot x$, or equivalently $[g,h] \cdot x =x$. Because vertex-stabilisers are trivial, we conclude that $[g,h]=1$, i.e., $g$ and $h$ commute.
\end{proof}

\begin{lemma}\label{lem:rotativestab}
Let $G$ be a group acting on a quasi-median graph $X$ with trivial vertex-stabilisers. Let $J_1,J_2$ be two transverse hyperplanes and $S$ a sector delimited by $J_1$. For every $g \in \mathrm{stab}_{\circlearrowleft}(J_2)$, $g \cdot S =S$. 
\end{lemma}

\begin{proof}
According to Proposition \ref{prop:transversehypcube}, there exists a prism $P=C_1 \times C_2$ such that the cliques $C_1,C_2$ are dual to the hyperplanes $J_1,J_2$ respectively. Let $x \in C_1$ such that $S= [C,x]$. Without loss of generality, we may suppose that $C_1 \cap C_2 = \{ x \}$ (otherwise, replace $C_1$ and $C_2$ with $C_1 \times \{ x \}$ and $\{x \} \times C_2$ respectively). Fix some vertex $y \in C_1 \backslash \{ x \}$. Because $g$ stabilises $C_2$ and $\{ y \} \times C_2$, necessarily $gx \in C_2$ and $gy \in \{ y \} \times C_2$. Moreover, because $x$ and $y$ are linked by an edge, we deduce that $gy \in C_1 \times \{ gx \}$, so that $C_1 \times \{ gx \}$ contains both $gx$ and $gy$; as a consequence, $gC_1= C_1 \times \{ gx \}$. In particular, this implies that $g \in \mathrm{stab}(J_1)$ since $C_1$ and $gC_1$ are parallel. Therefore, in order to conclude that $gS=S$, it is sufficient to show that $gS \subset S$.

\medskip \noindent
Let $y \in S=[C,x]$, i.e., $\mathrm{proj}_{C_1}(y)=x$. Then
$$gx= g \cdot \mathrm{proj}_{C_1}(y)= \mathrm{proj}_{gC_1}(gy),$$
hence $gy \in [gC_1,gx]$. On the other hand, $C_1= C_1 \times \{x\}$ and $gC_1= C_1 \times \{ gx \}$, so $x$ and $gx$ belong to the same sector delimited by the hyperplane dual to both $C_1$ and $gC_1$, namely $J_1$, so $[gC_1,gx]=[C_1,x]=S$. Therefore, $gy \in S$. This concludes the proof.
\end{proof}

\begin{proof}[Proof of Theorem \ref{thm:embeddinggraphproduct}.]
Fix a base point $x_0 \in X$. Let $\Gamma$ be a finite subgraph of $\Delta X$. If $J \in V(\Gamma)$ is a hyperplane of $X$, let $S_1, \ldots, S_n$ denote the sectors it delimits which contain $x_0$ or another hyperplane of $X$ which belongs to $V(\Gamma)$. For every $1 \leq i,j \leq n$, let $h_{ij} \in H_J$ be an element satisfying $h_{ij} \cdot S_i = S_j$; notice that there exists at most one such element because $H_J$ acts freely on the set of the sectors delimited by $J$. If no such element exists, set $h_{ij}=1$. Because $H_J$ is residually finite, there exists some finite-index subgroup $K_J \leq H_J$ satisfying $K_J \cap F_J \subset \{1 \}$ where $F_J = \{ h_{ij} \mid 1 \leq i,j \leq n\}$. Finally, let $X_J$ denote the union of all the sectors delimited by $J$ different from $S_1, \ldots, S_n$.

\medskip \noindent
Now, we want to prove that $\Gamma$, $\mathcal{G}= \{ K_J \mid J \in V(\Gamma) \}$, $x_0$ and the $X_J$'s satisfy the hypotheses of Proposition \ref{prop:pingpong} in order to conclude that the subgroup $\langle K_J \mid J \in V(\Gamma) \rangle$ is isomorphic to the graph product $\Gamma \mathcal{G}$. 

\medskip \noindent
From the definition of the $X_J$'s, it is clear that $x \notin \bigcup\limits_{J \in V(\Gamma)} X_J$. Next, it follows from Lemma \ref{lem:rotativestabcommute} that, if $J_1,J_2 \in V(\Gamma)$ are two adjacent vertices, then any element of $K_{J_1}$ commutes with any element of $K_{J_2}$; moreover, as a consequence of Lemma \ref{lem:rotativestab}, any element of $K_{J_1}$ stabilises any sector delimited by $J_2$, so it necessarily stabilises $X_{J_2}$. 

\medskip \noindent
Let $J \in V(\Gamma)$ and $g \in K_J \backslash \{ 1 \}$. Denote by $S$ the sector delimited by $J$ which contains $x_0$. By definition of $K_J$, we know that $g \notin F_J$, so the sector $g \cdot S$ cannot contain $x_0$ or a hyperplane of $V(\Gamma)$, hence $g \cdot x_0 \in g \cdot S \subset X_J$.

\medskip \noindent
Finally, let $J_1,J_2 \in V(\Gamma)$ be two transverse hyperplanes and $g \in K_{J_1} \backslash \{1 \}$. We defined $X_{J_2}$ as the union of all the sectors delimited by $J_2$ which do not contain $x_0$ or another hyperplane of $V(\Gamma)$. Therefore, if $S_2$ is one of these sectors, then the sector $S_1$ delimited by $J_1$ which contains $J_2$ must contain $S_2$ as well. But $g \in K_{J_1} \backslash \{1 \}$ implies $g \notin F_J$, so the sector $g \cdot S_1$ cannot contain $x_0$ or a hyperplane of $V(\Gamma)$, hence $g \cdot S_1 \subset X_{J_1}$. We conclude that $g \cdot S_2 \subset g \cdot S_1 \subset X_{J_1}$. A fortiori, $g \cdot X_{J_2} \subset X_{J_1}$.
\end{proof}

\noindent
It is worth noticing that the conclusion of Theorem \ref{thm:embeddinggraphproduct} may be empty. Indeed, a group may act on a quasi-median graph with large hyperplane-stabilisers but with trivial rotative stabilisers of hyperplanes, so that the graph product we produce is trivial.  

\begin{cor}\label{cor:embedRAAG}
Let $G$ be a group acting on a quasi-median graph $X$ with trivial vertex-stabilisers. Suppose that the rotative stabiliser of every hyperplanes contains an infinite cyclic group. If $\Gamma$ is a finite induced subgraph of $\Delta X$, then the right-angled Artin group $A(\Gamma)$ embeds into $G$. 
\end{cor}

\noindent
As an application, let us mention the following statement:

\begin{cor}
Let $\Gamma$ be a simplicial graph and $\mathcal{G}$ a collection of groups labelled by $V(\Gamma)$. Suppose that every group of $\mathcal{G}$ contains an infinite-order element. If $\mathrm{diam}(\Gamma) \geq 3$ then, for every finite forest $F$, the right-angled Artin group $A(F)$ embeds into the graph product $\Gamma \mathcal{G}$.
\end{cor}

\begin{proof}
According to Fact \ref{fact:isometricX}, $\X$ is isomorphic to $X(\Gamma,\mathcal{Z})$, where $\mathcal{Z}$ contains only infinite cyclic groups, so the crossing graph of $\X$ is isomorphic to Kim and Koberda's extension graph $\Gamma^e$. According to \cite[Proposition 5.2]{embeddingRAAG}, if $\Gamma$ contains a path with four vertices $P_4$ as an induced subgraph, the conclusion holds thanks to Corollary \ref{cor:embedRAAG}. Since we know that $\mathrm{diam}(\Gamma) \geq 3$, this concludes the proof.
\end{proof}

\noindent
The above discussion suggests our next question, in which, for every $n \geq 0$, we denote by $P_n$ the segment of length $n$, $C_n$ the cycle of length $n$ and $E_n$ the join of a vertex with $n$ pairwise non adjacent vertices.

\begin{question}
Let $n \geq 0$ be an integer, $\Gamma$ a simplicial graph and $\mathcal{G}$ a collection of groups labelled by $V(\Gamma)$. When does the crossing graph of $\X$ contain the path $P_n$, the cycle $C_n$ or the star $E_n$?
\end{question}

\noindent
For instance, it can be proved that the crossing graph of $\X$ contains a cycle of length at least five if and only if so does $\Gamma$.

\section{Application to wreath products}\label{section:appli3}

\noindent
Recall that, given two groups $G,H$ and an $H$-set $\Omega$, the \emph{permutational wreath product}\index{Permutational wreath products} of $G$ and $H$ with respect to $\Omega$ is defined as the semidirect product
$$G \wr_{\Omega} H = \left( \bigoplus\limits_{p \in \Omega} G \right) \rtimes H,$$
where $H$ acts on the direct sum by permuting the coordinates. More precisely, an element $(h,\varphi)$ of $G \wr_{\Omega} H$ is given by an element $h \in H$ and a map $\varphi : \Omega \to G$ with finite support, denoted by $\varphi \in G^{(\Omega)}$; and the product is defined by
$$(h_1, \varphi_1) \cdot (h_2, \varphi_2) = (h_1h_2, \varphi_1( \cdot) \varphi_2(h_1^{-1} \cdot))$$
for every $h_1,h_2 \in H$ and $\varphi_1,\varphi_2 \in G^{(\Omega)}$. If $\Omega=H$, on which $H$ acts by left-multiplication, we denote by $G \wr H$ the associated permutational wreath product. 

\medskip \noindent
In this section, given two groups $G,H$, a CAT(0) cube complex $X$ on which $H$ acts and a union of $H$-orbits $\Omega \subset X$, we construct an action of $G \wr_{\Omega} H$ on a quasi-median graph $\mathfrak{W}$, we call the \emph{graph of wreaths}, which turns out to be topical-transitive if points of $\Omega$ have trivial stabilisers. Therefore, our applications apply to permutational wreath products of the form
$$\left( \underset{\text{$n$ direct sums}}{\underbrace{\bigoplus\limits_{h \in H} G \oplus \cdots \oplus \bigoplus\limits_{h \in H} G}} \right) \rtimes H,$$
and in particular to wreath products $G \wr H$. Recall that such a wreath product $G \wr H$ has a nice interpretation as a lamplighter group: Fixing two generating sets $S$ and $R$ of $G$ and $H$ respectively, $S \cup R$ generates $G \wr H$ (identifying $G$ with its copy in the direct sum labelled by the neutral element of $H$). An element of $G \wr H$ can be described as a copy $\Gamma$ of the Cayley graph of $H$ constructed from $R$ whose vertices are labelled by elements of $G$, in such a way that all but finitely many vertices are labelled by $1$, together with an arrow labelling some vertex of $\Gamma$. Formally, the labelled graph $\Gamma$ encodes a map of $G^{(H)}$ and the arrow an element of $H$. Now, right-multiplicating the corresponding element of $G \wr H$ by an element of $s \in S$ multiplies the element of $G$ labelling the vertex where the arrow is by $s$, and right-multiplicating by an element of $r \in R$ moves the arrow from a vertex $v$ to its neighbor $vr$. 

\begin{figure}
\begin{center}
\includegraphics[scale=0.9]{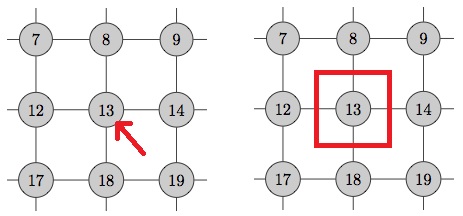}
\end{center}
\caption{Replacing the arrow with a rectangle.}
\label{figure1}
\end{figure}

\medskip \noindent
Let us describe the construction of our quasi-median graph for the wreath product $G \wr \mathbb{Z}^2$. According to the previous paragraph, an element of $G \wr \mathbb{Z}^2$, thought of as a lamplighter group, can be described by an infinite grid whose vertices are labelled by elements of $G$, such that all but finitely many vertices are labelled by $1$, together with an arrow labelling some vertex. See Figure \ref{figure1}.

\begin{figure}
\begin{center}
\includegraphics[scale=0.6]{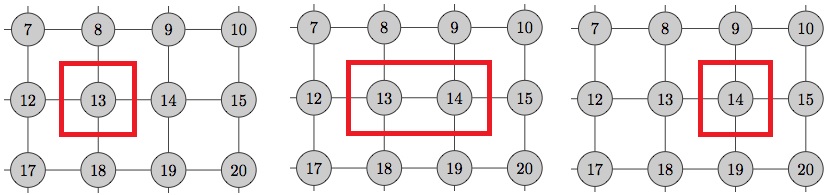}
\end{center}
\caption{Passing from a vertex to an adjacent vertex by elementary moves.}
\label{figure17}
\end{figure}

\medskip \noindent
Essentially, our construction lies on the following idea: replace the arrow of the previous description with a rectangle (whose corners have their coordinates in $\frac{1}{2} \mathbb{Z}$) containing a single vertex of the grid (see Figure \ref{figure1}), and, instead of moving the arrow from a vertex to one of its neighbors, move the sides of the rectangle independently. For instance, in order to move the rectangle to a vertex to one of its neighbors, three moves are necessary; see Figure \ref{figure17}. More formally, we define a \emph{wreath} $(R, \varphi)$ as the data of a rectangle $R$ and a map $\varphi : \mathbb{Z}^2 \to G$ with finite support. Now, our elementary moves on a given wreath $(R, \varphi)$ are the following: modify the label of a vertex which belongs to (the interior of) $R$, or translate one side of $R$ by a unit vector. Among the wreaths, we recover the group $G \wr \mathbb{Z}^2$ as the wreaths whose rectangles contain a single vertex of the grid. Moreover, we have a natural action of $G \wr \mathbb{Z}^2$ of the set of wreaths extending the left-multiplication:
$$( \{p \}, \psi) \cdot (R, \varphi) = \left( R+p, \psi(\cdot)+ \varphi(\cdot -p) \right).$$
Now, we define the \emph{graph of wreaths} $\mathfrak{W}$ as the graph whose vertices are the wreaths and whose edges link two wreaths such that one can be obtained from another by one of our elementary moves. The wreath product $G \wr \mathbb{Z}^2$ acts on $\mathfrak{W}$ via the action described above. 

\medskip \noindent
Interestingly, if $G= \mathbb{Z}/2 \mathbb{Z}$, then the graph of wreaths turns out to be a median graph (or equivalently, the $1$-skeleton of a CAT(0) cube complex); but this graph is no longer median whenever $G$ has cardinality at least three. Nevertheless, it turns out to be quasi-median, and this is the graph we will be interested in.

\subsection{Graphs of finite subcomplexes of CAT(0) cube complexes}

\noindent
Let $X$ be a CAT(0) cube complex and let $\mathcal{FS}(X)$ denote the set of the non empty finite convex subcomplexes of $X$. We define a metric on $\mathcal{FS}(X)$ by
$$d_{\mathcal{FS}(X)} : (C_1,C_2) \mapsto \# \left( \mathcal{D}(C_1) \oplus \mathcal{D}(C_2) \right)$$
where $\oplus$ denotes the symmetric difference and $\mathcal{D}(C)$ the collection of the halfspaces containing $C$. Notice that $\mathcal{D}(C_1) \oplus \mathcal{D}(C_2)$ is finite for every $C_1,C_2 \in \mathcal{FS}(X)$ because a hyperplane delimiting one halfspace of this collection either intersects $C_1$ or $C_2$, or separates $C_1$ and $C_2$, so $d_{\mathcal{FS}(X)}$ takes only finite values. Let us verify that this is indeed a distance. Fix three non empty finite convex subcomplexes $C_1,C_2,C_3$. Since $\mathcal{D}(C_1) \oplus \mathcal{D}(C_2)= \mathcal{D}(C_2) \oplus \mathcal{D}(C_1)$, it follows that $d_{\mathcal{FS}(X)} (C_1,C_2) = d_{\mathcal{FS}(X)}(C_2,C_1)$. Moreover, we deduce from the inclusion
$$\mathcal{D}(C_1) \oplus \mathcal{D}(C_3) \subset \left( \mathcal{D}(C_1) \oplus \mathcal{D}(C_2) \right) \cup \left( \mathcal{D}(C_2) \oplus \mathcal{D}(C_3) \right)$$
that $d_{\mathcal{FS}(X)}(C_1,C_3) \leq d_{\mathcal{FS}(X)}(C_1,C_2)+ d_{\mathcal{FS}(X)}(C_2,C_3)$. Finally, if $d_{\mathcal{FS}(X)}(C_1,C_2)=0$, then necessarily $\mathcal{D}(C_1)= \mathcal{D}(C_2)$, hence
$$C_1= \bigcap\limits_{D \in \mathcal{D}(C_1)} D = \bigcap\limits_{D \in \mathcal{D}(C_2)} D = C_2$$
since a convex subcomplex is equal to the intersection of all the halfspaces which contain it. Notice that the metric $d_{\mathcal{FS}(X)}$ can be described alternatively in the following way (where, for every subcomplex $Y \subset X$, $\mathcal{H}(Y)$ denotes the set of the hyperplanes intersecting the gated hull of $Y$):

\begin{fact}\label{fact:dFSalter}
For every $C_1,C_2 \in \mathcal{FS}(X)$, 
$$d_{\mathcal{FS}(X)}(C_1,C_2)= 2 \cdot \# \mathcal{H}(C_1 \cup C_2) - \# \mathcal{H}(C_1) - \# \mathcal{H}(C_2).$$
\end{fact}

\noindent
Set $F = 2 \cdot \chi_{\mathcal{H}(C_1 \cup C_2)} - \chi_{\mathcal{H}(C_1)} - \chi_{\mathcal{H}(C_2)}$ and fix some hyperplane $J$. We distinguish five possible cases:
\begin{itemize}
	\item If $J \notin \mathcal{H}(C_1 \cup C_2)$, then a halfspace delimited by $J$ either contains both $C_1$ and $C_2$ or is disjoint from from $C_1$ and $C_2$. A fortiori, $\mathcal{D}(C_1) \oplus \mathcal{D}(C_2)$ contains no halfspace delimited by $J$. Notice that $F(J)=0$.
	\item If $J \in \mathcal{H}(C_1) \backslash \mathcal{H}(C_2)$, then one halfspace delimited by $J$ contains $C_2$ and not $C_1$, and the other does not contain neither $C_1$ nor $C_2$. A fortiori, $\mathcal{D}(C_1) \oplus \mathcal{D}(C_2)$ contains precisely one halfspace delimited by $J$. Notice that $F(J)=1$.
	\item If $J \in \mathcal{H}(C_2) \backslash \mathcal{H}(C_1)$, the situation is symmetric: $\mathcal{D}(C_1) \oplus \mathcal{D}(C_2)$ contains precisely one halfspace delimited by $J$ and $F(J)=1$.
	\item If $J \in \mathcal{H}(C_1) \cap \mathcal{H}(C_2)$, then none of the two halfspaces delimited by $J$ belongs to $\mathcal{D}(C_1) \oplus \mathcal{D}(C_2)$. Notice that $F(J)=0$. 
	\item If $J$ separates $C_1$ and $C_2$, then one halfspace delimited by $J$ contains $C_1$ but not $C_2$, and the other contains $C_2$ but not $C_1$. A fortiori, $\mathcal{D}(C_1) \oplus \mathcal{D}(C_2)$ contains the two halfspaces delimited by $J$. Notice that $F(J)=2$.
\end{itemize}
Thus, we have proved that, for every hyperplane $J$, $\mathcal{D}(C_1) \oplus \mathcal{D}(C_2)$ contains exactly $F(J)$ haflspaces delimited by $J$. Therefore, 
$$d_{\mathcal{FS}(X)}(C_1,C_2) = \# \mathcal{D}(C_1) \oplus \mathcal{D}(C_2) = \sum\limits_{J \in \mathfrak{J}} F(J),$$
where $\mathfrak{J}$ denotes the set of all the hyperplanes of $X$. This concludes the proof of our fact.

\medskip \noindent
Our next result states that the metric space we have defined turns out to be geodesic. 

\begin{prop}\label{prop:FSgeodesic}
The metric space $(\mathcal{FS}(X),d_{\mathcal{FS}(X)})$ is geodesic.
\end{prop}

\noindent
First of all, we need to describe the intervals of our metric space $(\mathcal{FS}(X), d_{\mathcal{FS}(X)})$. 

\begin{lemma}\label{lem:FSinterval}
Let $C_1,C_2,Q \in \mathcal{FS}(X)$ be three subcomplexes. The equality
$$d_{\mathcal{FS}(X)}(C_1,C_2)= d_{\mathcal{FS}(X)}(C_1,Q)+ d_{\mathcal{FS}(X)}(Q,C_2)$$
holds if and only if $Q$ is contained into the convex hull of $C_1 \cup C_2$ and if the inclusion $\mathcal{D}(Q) \subset \mathcal{D}(C_1) \cup \mathcal{D}(C_2)$ is satisfied.
\end{lemma}

\begin{proof}
Recall that, for any set $S$ and any subsets $A,B,C \subset S$, the equalities
\begin{itemize}
	\item $(A \oplus C) \cup (B \oplus C) = (A \oplus B) \cup (C \backslash (A \cup B)) \cup ((A \cap B) \backslash C)$;
	\item $(A \oplus C) \cap (B \oplus C) = (C \backslash (A \cup B)) \cup ((A \cap B ) \backslash C)$;
\end{itemize}
hold. Because we know that
$$\mathcal{D}(C_1) \oplus \mathcal{D}(C_2) \subset \left( \mathcal{D}(C_1) \oplus \mathcal{D}(Q) \right) \cup \left( \mathcal{D}(Q) \oplus \mathcal{D}(C_2) \right),$$
the equality we are interested in holds if and only 
$$\mathcal{D}(C_1) \oplus \mathcal{D}(C_2) = \left( \mathcal{D}(C_1) \oplus \mathcal{D}(Q) \right) \sqcup \left( \mathcal{D}(Q) \oplus \mathcal{D}(C_2) \right),$$
which is also equivalent, according to our two previous equalities, to
$$\mathcal{D}(Q) \backslash \left( \mathcal{D}(C_1) \cup \mathcal{D}(C_2) \right) = \left( \mathcal{D}(C_1) \cap \mathcal{D}(C_2) \right) \backslash \mathcal{D}(Q) = \emptyset.$$
Notice that $\mathcal{D}(C_1) \cap \mathcal{D}(C_2) \subset \mathcal{D}(Q)$ precisely means that any halfspace containing both $C_1$ and $C_2$ must contain $Q$ as well, ie., $Q$ is contained in the convex hull of $C_1 \cup C_2$ since the convex subcomplex is equal to the intersection of all the halfspaces which contain it. This concludes the proof.
\end{proof}

\begin{proof}[Proof of Proposition \ref{prop:FSgeodesic}.]
Let $C_1,C_2 \in \mathcal{FS}(X)$ be two subcomplexes. We want to prove that there exists a geodesic between $C_1$ and $C_2$ in $\mathcal{FS}(X)$. To do that, we argue by induction on $d_{\mathcal{FS}(X)}(C_1,C_2)$. If $d_{\mathcal{FS}(X)}(C_1,C_2)=0$ there is nothing to prove, so we suppose that $d_{\mathcal{FS}(X)}(C_1,C_2) \geq 1$. Let $Q$ denote the convex hull of $C_1$ and $C_2$. It follows from Lemma \ref{lem:FSinterval} that $Q$ belongs to the interval between $C_1$ and $C_2$, ie., 
$$d_{\mathcal{FS}(X)}(C_1,C_2)= d_{\mathcal{FS}(X)}(C_1,Q)+ d_{\mathcal{FS}(X)}(Q,C_2),$$
so that we are able to conclude by applying twice our induction hypothesis if $Q$ is different from both $C_1$ and $C_2$. From now on, suppose that $C_2=Q$ (the case $C_1=Q$ being symmetric), so that $C_1 \subset C_2$. Let $x \in C_2$ be a vertex which does not belong to $C_1$ but which is adjacent to one of its vertices, say $x'$, and let $Q'$ denote the convex hull of $C_1 \cup \{ x\}$. By noticing that the hyperplane separating $x$ and $x'$ is the unique hyperplane intersecting $Q'$ but not $C_1$, it follows that $d_{\mathcal{FS}(X)}(C_1,Q')=1$. On the other hand, it follows from Lemma \ref{lem:FSinterval} that
$$d_{\mathcal{FS}(X)}(C_1,C_2)= d_{\mathcal{FS}(X)}(C_1,Q')+ d_{\mathcal{FS}(X)}(Q',C_2).$$
In particular, $d_{\mathcal{FS}(X)}(Q',C_2) = d_{\mathcal{FS}(X)}(C_1,C_2)-1$, so that our induction hypothesis applies to $Q'$ and $C_2$, ie., there exists a geodesic in $\mathcal{FS}(X)$ between $Q'$ and $C_2$. By adding the point $C_1$ to this path, we find a geodesic in $\mathcal{FS}(X)$ between $C_1$ and $C_2$. This concludes the proof.
\end{proof}

\noindent
As a consequence of Proposition \ref{prop:FSgeodesic}, we deduce that $\mathcal{FS}(X)$ can be naturally thought of as a graph by linking two subcomplexes at distance one appart with respect to $d_{\mathcal{FS}(X)}$, and if so, the associated graph metric coincides with $d_{\mathcal{FS}(X)}$. 

\begin{definition}
Let $X$ be a CAT(0) cube complex. The \emph{graph of finite subcomplexes} of $X$, also denoted by $\mathcal{FS}(X)$, is the graph whose vertices are the non empty finite convex subcomplexes of $X$ and whose edges link two subcomplexes at $d_{\mathcal{FS}(X)}$-distance one appart. 
\end{definition}

\noindent
For convenience, we introduce the following notation. If $C_1,C_2 \in \mathcal{FS}(X)$ are two adjacent subcomplexes, then there exists some halfspace $H$ of $X$ such that $\mathcal{D}(C_2)= \mathcal{D}(C_1) \oplus \{ H \}$. If so, we denote by $[C_1,H]$ the subcomplex $C_2$. (Notice that, if $C_3$ is another neighbor of $C_1$ in $\mathcal{FS}(X)$, then $\mathcal{D}(C_3) = \mathcal{D}(C_1) \oplus \{ H'\}$ for some halfspace $H' \neq H$ since 
$$\bigcap\limits_{K \in \mathcal{D}(C_1) \oplus \{ H' \}} K = C_3 \neq C_2= \bigcap\limits_{K \in \mathcal{D}(C_1) \oplus \{ H \}} K.$$
Therefore, our notation cannot be ambiguous.) If $C \in \mathcal{FS}(X)$ is a subcomplex and $H_1,H_2, \ldots$ a sequence of halfspaces, we write $[C,H_1,H_2]$ instead of $[[C,H_1],H_2]$, and $[C,H_1,H_2,H_3]$ instead of $[[C,H_1,H_2],H_3]$, and so on. 

\begin{remark}
If $C \in \mathcal{FS}(X)$ is a subcomplex and $H \notin \mathcal{D}(C)$ a halfspace of $X$ such that $[C,H]$ is well-defined, then $[C,H]$ is a minimal subcomplex of $X$ containing strictly $C$. Alternatively, $[C,H]$ is the convex hull of the union of $C$ together with a vertex which does not belong to $H$ (and a fortiori to $C$) but which is adjacent to some vertex of $C$. This stresses out the idea that $[C,H]$ is one the subcomplexes of $\mathcal{FS}(X)$ which are ``nearest'' to $C$. Geometrically, $[C,H]$ is the union of $C$ together with all the cubes intersecting both $C$ and the hyperplane underlying $H$. Similarly, if $H \in \mathcal{D}(C)$ is such that $[C,H]$ is well-defined, then, because $C = [[C,H],H]$ with $H \notin [C,H]$, the previous discussion implies that $[C,H]$ is a maximal subcomplex of $X$ contained properly into $C$, and it can be obtained from $C$ by removing the interiors of all the cubes intersecting the hyperplane underlying $H$. 
\end{remark}

\noindent
We conclude this section by noticing that the graphs we have associated to  CAT(0) cube complexes are themselves CAT(0) cube complexes. 

\begin{prop}\label{prop:FSmedian}
The graph of finite subcomplexes $\mathcal{FS}(X)$ of any CAT(0) cube complex $X$ is a median graph. Moreover, $\mathcal{FS}(X)$ is locally finite if and only if $X$ is locally finite as well.
\end{prop}

\begin{proof}
Let $X$ be a CAT(0) cube complex and $C_1,C_2,C_3 \in \mathcal{FS}(X)$ three subcomplexes. Let $I$ denote the intersection of intervals
$$I(C_1,C_2) \cap I(C_2,C_3) \cap I(C_3,C_1).$$
According to Lemma \ref{lem:FSinterval}, any subcomplex which belongs to $I$ must be included into the intersection $C$ of the convex hulls of the unions of $C_1$ and $C_2$, $C_2$ and $C_3$, and $C_3$ and $C_1$, ie.,
$$\bigcap\limits_{1 \leq i<j \leq 3} \{ \text{halfspaces containing both $C_i$ and $C_j$} \}.$$
First, notice that $C$ is a finite convex subcomplex of $X$, since it is an intersection of such subcomplexes. Next, $C$ is non empty, as a consequence of the following claim. 

\begin{claim}\label{claim:FSmedian}
Let $x_1 \in C_1$, $x_2 \in C_2$ and $x_3 \in C_3$ be three vertices. The median point $m=m(x_1,x_2,x_3)$ belongs to $C$.
\end{claim}

\noindent
For every $1 \leq i < j \leq 3$, $m$ belongs to a geodesic between $x_i$ and $x_j$. A fortiori, $m$ belongs to the convex hull of $C_i \cup C_j$. Therefore, $m$ belongs to the intersection of the convex hulls between $C_1$ and $C_2$, $C_2$ and $C_3$, and $C_3$ and $C_1$, which is precisely $C$. This proves our claim.

\medskip \noindent
Thus, we have proved that $C \in \mathcal{FS}(X)$. We claim that $C$ is the unique median point of the triple $\{C_1,C_2,C_3 \}$, ie., $I= \{ C \}$. 

\medskip \noindent
First, let $Q \in \mathcal{FS}(X)$ be a subcomplex strictly included into $C$. So there exists a vertex $x \in C \backslash Q$. Because $Q$ is convex, there exists a hyperplane $J$ separating $x$ and $Q$; let $D$ denote the halfspace delimited by $J$ which contains $Q$. A fortiori, $J$ intersects $C$, so $D$ canoot contain two subcomplexes among $\{ C_1,C_2,C_3 \}$. Say that $D$ does not contain neither $C_1$ nor $C_2$. It follows from Lemma \ref{lem:FSinterval} that $Q$ does not belong to the interval $I(C_1,C_2)$ in $\mathcal{FS}(X)$. Therefore, $C$ is the unique candidate for a median point, ie., $I \subset \{ C \}$.

\medskip \noindent
Now, in order to conclude that $C$ is a median point, according to Lemma \ref{lem:FSinterval} it is sufficient to show that, for every $1 \leq i <j \leq 3$, any halfspace $D$ containing $C$ must contain either $C_i$ or $C_j$. Suppose by contradiction that such a $D$ does not contain neither $C_i$ nor $C_j$. So there exist vertices $x_i \in C_i$ and $x_j \in C_j$ which do not belong to $D$. If $\{1,2,3 \} \backslash \{ i,j \} = \{ k\}$, fix a vertex $x_k \in C_k$. According to Claim \ref{claim:FSmedian}, the median point $m=m(x_1,x_2,x_3)$ belongs to $C$. On the other hand, $x_i$ and $x_j$ belong to $D^c$, which is convex, so $m \in D^c$. In particular, we deduce that $D^c \cap C \neq \emptyset$, which contradicts the fact that $C$ is contained into $D$. Thus, we have proved that $\mathcal{FS}(X)$ is a median graph.

\medskip \noindent
Now, suppose that $X$ is not locally finite, ie., there exists at least one vertex $x \in X$ with infinitely many neighbors, say $x_1,x_2, \ldots \in X$. By noticing that, for every $i \geq 1$, the edge $[x,x_i]$ is a neighbor of $\{x \}$ in $\mathcal{FS}(X)$, we deduce that $\mathcal{FS}(X)$ is not locally finite. Now, suppose that $X$ is locally finite. Let $C \in \mathcal{FS}(X)$ be a subcomplex, and $[C,H] \in \mathcal{FS}(X)$ one of its neighbors, where $H$ is a halfspace of $X$ delimited by some hyperplane $J$. Because $C$ is finite, only finitely many hyperplanes intersect $C$, so we may suppose that $J$ is disjoint from $C$; in particular, this implies that $C \subset [C,H]$. Fix two vertices $a \in C$ and $b \in N(J)$ minimising the distance between $C$ and $N(J)$. According to Lemma \ref{lem:2min}, any hyperplane separating $a$ and $b$ must separate $C$ and $N(J)$, so, if $a \neq b$, we deduce that the distance between $C$ and $[C,H]$ is at least two, contradicting our assumption. Therefore, $a$ must be equal to $b$, so that $J$ is dual to an edge with a vertex of $C$ as one of its endpoints. Consequently, the number of neighbors of $C$ in $\mathcal{FS}(X)$ is at most
$$2 \left( \# \mathcal{H}(C) + \sum\limits_{v \in C} \mathrm{deg}(v) \right) < + \infty.$$
This concludes the proof.
\end{proof}

\subsection{Graph of wreaths}\label{section:graphofwreaths}

\noindent
Fix two groups $G$ and $H$, and suppose that $H$ acts on a CAT(0) cube complex $X$. Let $\Omega \subset X$ denote an $H$-invariant collection of vertices (ie., a union of $H$-orbits).  

\begin{definition}
A \emph{wreath} $(C, \varphi)$ is the data of a non empty finite convex subcomplex $C \subset X$ and a function $\varphi : \Omega \to G$ with finite support, ie., such that $\varphi(p)=1$ for all but finitely many $p \in \Omega$ (written $\varphi \in G^{(\Omega)}$). 
\end{definition}

\begin{definition}
The \emph{graph of wreaths} $\mathfrak{W}$ is the graph whose vertices are the wreaths and whose edges link two wreaths $(C, \mathfrak{g})$ and $(Q, \mathfrak{h})$ either if $\mathfrak{g}= \mathfrak{h}$ and $Q=[C,H]$ for some halfspace $H$, or if $C=Q$ and $\mathfrak{g}$ and $\mathfrak{h}$ differ on a single point of $C \cap \Omega$. 
\end{definition}

\noindent
It is worth noticing that an edge of the graph of wreaths is naturally labelled either by a halfspace (and its underlying hyperplane) or by a vertex of $\Omega$. 

\medskip \noindent
Notice also that the permutational wreath product $G \wr_{\Omega} H$ can be naturally thought of as a set of vertices of the graph of wreaths: those with a first coordinate which is a singleton. Furthermore, the action of $G \wr_{\Omega} H$ on itself by left-multiplication extends to an action on the set of all the wreaths in the following way:
$$( h, \psi) \cdot (C, \varphi) = (hC, \psi( \cdot ) \varphi (h^{-1} \cdot ) )$$
for every $(h, \psi) \in G \wr_{\Omega} H$ and $(C, \varphi) \in \mathfrak{W}$. Observe that this action preserves adjacency, so that we have defined an action of $G \wr_{\Omega} H$ on the graph of wreaths $\mathfrak{W}$, which turns out to be a quasi-median graph according to our following result.

\begin{prop}\label{prop:wreathqm}
The graph of wreaths $\mathfrak{W}$ is a quasi-median graph.
\end{prop}

\noindent
In all the other applications of our paper, there is a natural basepoint in the quasi-median graph we consider, and, up to an isometry, it may be supposed that a given vertex is this basepoint. Often, this freedom simplifies the proofs. However, our action on the graph of wreaths is not vertex-transitive, so that such a simplification may be not possible. On the other hand, thanks to the next lemma, we understand completely the orbits under the action of $G \wr_{\Omega} H$.

\begin{lemma}\label{lem:Worbit}
Two wreaths $(C_1, \varphi_1), (C_2, \varphi_2) \in \mathfrak{W}$ belong to the same $(G \wr_{\Omega} H)$-orbit if and only if there exists some $h \in H$ such that $hC_1=C_2$.
\end{lemma}

\begin{proof}
If $(C_1, \varphi_1)$ and $(C_2, \varphi_2)$ belong to the same orbit, there exists some $(h, \psi) \in G \wr_{\Omega} H$ such that
$$(C_2, \varphi_2)= ( h, \psi) \cdot (C_1, \varphi_1 ) =  (hC_1, \psi( \cdot) \varphi_1( h^{-1} \cdot)).$$
In particular, $hC_1=C_2$. Conversely, suppose that there exists some $h \in H$ such that $hC_1=C_2$. Then
$$\left( h, \varphi_2( \cdot) \varphi_1(h^{-1} \cdot)^{-1} \right) \cdot (C_1, \varphi_1)= (C_2, \varphi_2),$$
ie., $(C_1, \varphi_1)$ and $(C_2, \varphi_2)$ belong to the same orbit.
\end{proof}

\noindent
The first step towards the proof or Proposition \ref{prop:wreathqm} is to understand the metric of the graph of wreaths. This is the goal of the next preliminary lemma, in which, given some subcomplex $Y$ of $X$, we denote by $\mathcal{H}(Y)$ the set of the hyperplanes of the convex hull of $Y$; alternatively, $\mathcal{H}(Y)$ is the collection of the hyperplanes separating at least two vertices of $Y$.

\begin{lemma}\label{lem:wrdist}
The distance in $\mathfrak{W}$ between two wreaths $(C_1, \varphi_1)$ and $(C_2, \varphi_2)$ is 
$$2 \# \mathcal{H}(C_1 \cup C_2 \cup \mathrm{supp}(\varphi_1^{-1} \varphi_2)) - \# \mathcal{H}(C_1) - \# \mathcal{H}(C_2) + \# \mathrm{supp}(\varphi_1^{-1} \varphi_2).$$
\end{lemma}

\begin{proof}
Let $Q$ denote the convex hull of $C_1 \cup C_2 \cup \mathrm{supp}(\varphi_1^{-1} \varphi_2)$. We have
$$d((C_1, \varphi_1),(C_2, \varphi_2)) \leq d((C_1, \varphi_1),(Q, \varphi_1))+d((Q, \varphi_1),(Q, \varphi_2)) + d((Q,\varphi_2),(C_2, \varphi_2)).$$
Clearly, the distance between $(C_1, \varphi_1)$ and $(Q, \varphi_1)$ (resp. $(Q, \varphi_2)$ and $(C_2, \varphi_2)$) is bounded above by the distance between $C_1$ and $Q$ (resp. $Q$ and $C_2$) in the graph $\mathcal{FS}(X)$, so it follows from Fact \ref{fact:dFSalter} that
$$d((C_1, \varphi_1),(Q, \varphi_1)) \leq \# \mathcal{H}(Q)- \# \mathcal{H}(C_1) \ \text{and} \ d((Q, \varphi_2),(C_2, \varphi_2)) \leq \# \mathcal{H}(Q)- \# \mathcal{H}(C_2).$$
Finally, since $\mathrm{supp}(\varphi_1^{-1} \varphi_2) \subset Q$, we know that
$$d((Q, \varphi_1), (Q,\varphi_2)) \leq \# \mathrm{supp}(\varphi_1^{-1} \varphi_2).$$
Therefore, our distance between  $(C_1, \varphi_1)$ and $(C_2, \varphi_2)$ is bounded above by the quantity we are considering. Conversely, let $\gamma$ be a path in $\mathfrak{W}$ between  $(C_1, \varphi_1)$ and $(C_2, \varphi_2)$, and let 
$$Q_1=C_1, Q_2, \ldots, Q_{r-1}, Q_r= C_2$$ 
denote the path in $\mathcal{FS}(X)$ which is the projection of $\gamma$ (ie., consider the sequence of the first coordinates of vertices of $\gamma$). Write
$$F=2 \cdot \chi_{\mathcal{H}(C_1 \cup C_2 \cup \mathrm{supp}(\varphi_1^{-1} \varphi_2))} - \chi_{\mathcal{H}(C_1)} - \chi_{\mathcal{H}(C_2)}$$
and fix a hyperplane $J \in \mathcal{H}(C_1 \cup C_2 \cup \mathrm{supp}(\varphi_1^{-1} \varphi_2))$. We distinguish five cases.
\begin{itemize}
	\item If $J \in \mathcal{H}(C_1) \cap \mathcal{H}(C_2)$, notice that $F(J)=0$.
	\item Suppose that $J \in \mathcal{H}(C_1) \backslash \mathcal{H}(C_2)$. There must exist some $1 \leq i \leq r-1$ such that $Q_{i+1}=[Q_i,H]$ where $H$ is the halfspace delimited by $J$ containing $C_2$, so that at least one edge of $\gamma$ is labelled by $J$. Notice that $F(J)=1$.
	\item Suppose that $J \in \mathcal{H}(C_2) \backslash \mathcal{H}(C_1)$. By symmetry, at least one edge of $\gamma$ is labelled by $J$, and $F(J)=1$. 
\end{itemize}
From now on, suppose that $J \notin \mathcal{H}(C_1) \cup \mathcal{H}(C_2)$. Notice that $F(J)=2$.
\begin{itemize}
	\item Suppose that $J$ separates $C_1$ and $C_2$. Let $J^-$ (resp. $J^+$) denote the halfspace delimited by $J$ which contains $C_1$ (resp. $C_2$). Let $2 \leq p \leq r$ be the first index such that $J^- \notin Q_p$ (which exists since $J^- \notin \mathcal{D}(Q_r)$). Because $J^- \in Q_{p-1}$, necessarily $Q_p=[Q_{p-1},J^-]$. Notice that $J^-$ and $J^+$ do not belong to $\mathcal{D}(Q_p)$. On the other hand, $J^+ \in \mathcal{D}(C_2)$, so there must exist some $p+1 \leq q \leq r$ such that $J^+ \in \mathcal{H}(Q_q)$. Without loss of generality, suppose that $q$ is the first such index, so that $Q_q=[Q_{q-1},J^+]$. As a consequence, there exist at least two edges of $\gamma$ labelled by $J$.
	\item Suppose that $J$ separates $C_1 \cup C_2$ and some point $v \in \mathrm{supp}(\varphi_1^{-1} \varphi_2)$. Let $J^-$ denote the halfspace delimited by $J$ containing $v$ and $J^+$ the halfspace containing $C_1 \cup C_2$. There must exist some $2 \leq p \leq r$ such that $v \in Q_p$. Thus, $J^+$ belongs to $\mathcal{D}(C_1)$ but not to $\mathcal{D}(Q_p)$. Let $2 \leq a \leq p$ be the first index satisfying $J^+ \notin \mathcal{H}(Q_a)$. In particular, $Q_a=[Q_{a-1},J^+]$. On the other hand, we know that $J^+$ belongs to $\mathcal{D}(Q_r)$, so there must exist some index $a+1 \leq b \leq r$ such that $J^+ \notin \mathcal{D}(Q_b)$. Without loss of generality, suppose that $b$ is the first such index, so that $Q_b=[Q_{b-1},J^+]$. As a consequence, at least two edges of $\gamma$ are labelled by $J$.
\end{itemize}
Thus, we have proved that, for every hyperplane $J \in \mathcal{H}(C_1 \cup C_2 \cup \mathrm{supp}(\varphi_1^{-1} \varphi_2))$, there exist at least $F(J)$ edges of $\gamma$ labelled by $J$. A fortiori, the number of edges of $\gamma$ is at least $\sum\limits_{J \in \mathfrak{J}} F(J)$, where $\mathfrak{J}$ denotes the set of all the hyperplanes of $X$. On the other hand, we also know that, for every $p \in \mathrm{supp}(\varphi_1^{-1} \varphi_2)$, there must be some edge of $\gamma$ labelled by $p$. Therefore, 
$$\mathrm{length}(\gamma) \geq \sum\limits_{J \in \mathfrak{J}} F(J)+ \# \mathrm{supp}(\varphi_1^{-1} \varphi_2)$$
which is precisely the quantity we are looking for. This concludes the proof.
\end{proof}

\noindent
A consequence of our previous lemma is that the distance between any two wreaths is necessarily finite, hence:

\begin{cor}
The graph of wreaths $\mathfrak{W}$ is connected.
\end{cor}

\noindent
Just by enumerating all the possible cases, we deduce from Lemma \ref{lem:wrdist} the following corollary.

\begin{cor}\label{cor:distadj}
Let $(C,1),(Q, \psi),(Q',\psi')$ be three wreaths. Suppose that $(Q,\psi)$ and $(Q',\psi')$ are adjacent in the graph of wreaths, and let $k$ denote the distance between $(C,1)$ and $(Q,\psi)$. Then the distance between $(C,1)$ and $(Q',\psi')$ is equal to 
\begin{itemize}
	\item $k+1$ if $Q=Q'$ and $\mathrm{supp}(\psi) \subsetneq \mathrm{supp}(\psi')$;
	\item $k$ if $Q=Q'$ and $\mathrm{supp}(\psi)= \mathrm{supp}(\psi')$;
	\item $k-1$ if $Q=Q'$ and $\mathrm{supp}(\psi') \subsetneq \mathrm{supp}(\psi)$;
	\item $k+1$ if $\psi= \psi'$ and $Q'=[Q,H]$ where $H$ is a halfspace containing $Q$, $C$ and $\mathrm{supp}(\psi)$;
	\item $k-1$ if $\psi= \psi'$ and $Q'= [Q,H]$ where $H$ is a halfspace containing both $C$ and $\mathrm{supp}(\psi)$ which is delimited by a hyperplane $J \in \mathcal{H}(Q)$;
	\item $k+1$ if $\psi= \psi'$ and $Q'= [Q,H]$ where $H$ is a halfspace delimited by some hyperplane $J \in \mathcal{H}(Q \cup C \cup \mathrm{supp}(\psi)) \cap \mathcal{H}(Q)$;
	\item $k-1$ if $\psi= \psi'$ and $Q'= [Q,H]$ where $H$ is a halfspace delimited by some hyperplane $J$ satisfying $J \in \mathcal{H}(Q \cup C \cup \mathrm{supp}(\psi))$ and $J \notin \mathcal{H}(Q)$.
\end{itemize}
\end{cor}

\noindent
Another consequence of Lemma \ref{lem:wrdist} is that the graph of wreaths $\mathfrak{W}$ can be decomposed as a union of gated subgraphs, called \emph{leaves}, which are all isometric to the graph of finite subcomplexes $\mathcal{FS}(X)$. 

\begin{definition}
Given some map $\varphi \in G^{(\Omega)}$, the \emph{leaf} $\mathfrak{W}(\varphi)$ is the subgraph of $\mathfrak{W}$ generated by the vertices $\{ (C, \varphi) \mid C \in \mathcal{FS}(X) \}$.
\end{definition}

\begin{lemma}\label{lem:Wleaf}
Let $\varphi \in G^{(\Omega)}$ be a map. The leaf $\mathfrak{W}(\varphi)$ is a gated subgraph of $\mathfrak{W}$ and the map $C \mapsto (C, \varphi)$ defines an isometry $\mathcal{FS}(X) \to \mathfrak{W}(\varphi)$. Morever, for every wreath $(C, \phi) \in \mathfrak{W}$, the projection of $(C, \phi)$ onto $\mathfrak{W}(\varphi)$ is $(Q, \varphi)$ where $Q$ denotes the convex hull of $C \cup \mathrm{supp}(\varphi^{-1}\phi)$. 
\end{lemma}

\begin{proof}
Let $(C, \phi) \in \mathfrak{W}$ be a wreath and let $Q$ denote the convex hull of $C \cup \mathrm{supp}(\varphi^{-1}\phi)$. As a consequence of Lemma \ref{lem:FSinterval}, notice that, for every $(A, \varphi) \in \mathfrak{W}(\varphi)$, we have
$$d((C, \phi),(Q, \varphi))+d((Q, \varphi),(A, \varphi)) = d((C, \phi),(A, \varphi)).$$
This means that $(Q, \varphi)$ is a gate for $(C, \phi)$. Therefore, the leaf $\mathfrak{W}(\varphi)$ is gated and $(Q,\varphi)$ is the projection of $(C, \phi)$ onto $\mathfrak{W}(\varphi)$. Finally, the fact that the map $C \mapsto (C, \varphi)$ defines an isometry $\mathcal{FS}(X) \to \mathfrak{W}(\varphi)$ follows from Lemma \ref{lem:wrdist} and Fact \ref{fact:dFSalter}.
\end{proof}

\noindent
Now we are ready to prove Proposition \ref{prop:wreathqm}.

\begin{proof}[Proof of Proposition \ref{prop:wreathqm}.]
We begin by verifying the triangle condition. Let $w_a=( A, \zeta)$, $w_b = (B, \varphi)$ and $w_c=(C, \psi)$ be three wreaths such that $w_b$ and $w_c$ are adjacent in $\mathfrak{W}$, and such that the distances $d(w_a,w_b)$ and $d(w_a,w_c)$ are equal, say to $k$. According to Lemma \ref{lem:Worbit}, we can suppose without loss of generality that $\zeta=1$. It follows from Corollary \ref{cor:distadj} that necessarily $B=C$ and $\mathrm{supp}(\varphi) = \mathrm{supp}(\psi)$. Moreover, because $w_b$ and $w_c$ are adjacent (and, in particular, distinct), we know that $\varphi$ and $\psi$ differ on a single point $p \in B=C$. Notice that $\varphi(p) \neq 1$ and $\psi(p) \neq 1$. Therefore, setting
$$\xi : q \mapsto \left\{ \begin{array}{cl} \varphi(q)=\psi(q) & \text{if} \ q \neq p \\ 1 & \text{otherwise} \end{array} \right.,$$
the wreath $m=(B,\xi)$ is adjacent to both $w_b$ and $w_c$. Moreover, it follows from Lemma \ref{lem:wrdist} that
$$d(w_a,m)= d(w_a,w_b) + \# \mathrm{supp}(\varphi) - \# \mathrm{supp}(\xi)= k-1.$$
Consequently, $m$ is the wreath we are looking for. 

\medskip \noindent
Now, we want to verify the quandrangle condition. Let $w_a=(A,\xi)$, $w_b=(B,\phi)$, $w_c=(C, \psi)$ and $w_d=(D, \varphi)$ be four wreaths such that $w_c$ is adjacent to both $w_b$ and $w_d$ in $\mathfrak{W}$, such that the distances $d(w_a,w_b)$ and $d(w_a,w_d)$ are equal, say to $k$, and such that $d(w_a,w_c)=k+1$. According to Lemma \ref{lem:Worbit}, we can suppose without loss of generality that $\xi=1$. We want to prove that there exists a fifth vertex $m \in \mathfrak{W}$ which is adjacent to both $w_b$ and $w_d$ and satisfying $d(w_a,m)=k-1$. First, suppose that the edge between $w_b$ and $w_c$ is labelled by a point of $\Omega$. Following Corollary \ref{cor:distadj}, we distinguish three cases.
\begin{itemize}
	\item Suppose that $B=C$, that $\phi$ and $\psi$ differ on a single vertex $p \in B$ where $\phi(p)=1$, that $C=D$, and that $\psi$ and $\varphi$ differ on a single vertex $q \in C$ where $\varphi(q)=1$. Notice that $p \neq q$, since otherwise we would have $w_b=w_d$. Let $M$ denote the subcomplex $B=C=D$, and $\zeta$ the map which is equal to $\psi$ (and $\varphi$) outside $\{p,q\}$ and to the identity otherwise. Notice that $\zeta$ and $\phi$ differ only on $q$, and that $\zeta$ and $\varphi$ differ only on $p$, so $m=(M,\zeta)$ is adjacent to both $w_b$ and $w_d$. Moreover, from the equality $\mathrm{supp}(\zeta) = \mathrm{supp}(\phi) \backslash \{ q \}$, we deduce thanks to Lemma \ref{lem:wrdist} that $d(w_a,m)=k-1$.
	\item Suppose that $B=C$, that $\phi$ and $\psi$ differ on a single vertex $p \in B$ where $\phi(p) =1$, that $\psi= \varphi$, and that $D=[C,H]$ where $H$ is a halfspace delimited by some hyperplane $J$ satisfying $J \in \mathcal{H}(C \cup A \cup \mathrm{supp}(\psi))$ and $J \notin \mathcal{H}(C)$. Notice that $C \subset D$, so that $p$ belongs to $D$. This implies that $m=(D,\phi)$ is adjacent to $w_d$. We also notice that $m$ is adjacent to $w_b$ since $D=[B,H]$. Finally, because $\mathrm{supp}(\phi)= \mathrm{supp}(\psi) \sqcup \{p \}$ and $p \in D$, we deduce that $\mathcal{H}(A \cup D \cup \mathrm{supp}(\phi)) = \mathcal{H}(A \cup D \cup \mathrm{supp}(\psi))$, so that $$d(w_a,m)= d(w_a,w_d)+ \# \mathrm{supp}(\psi)- \# \mathrm{supp}(\phi)=k-1$$
according to Lemma \ref{lem:wrdist}. 
	\item Suppose that $B=C$, that $\phi$ and $\psi$ differ on a single vertex $p \in B$ where $\phi(p) =1$, that $\psi= \varphi$, and that $D=[C,H]$ where $H$ is a halfspace containing both $A$ and $\mathrm{supp}(\psi)$ which is delimited by some hyperplane $J \in \mathcal{H}(C)$. Notice that $p \in D$, so that $m=(D, \phi)$ si adjacent to $w_d$. Moreover, since $D=[B,H]$, $m$ is also adjacent to $w_b$. Finally, because $\mathrm{supp}(\varphi)= \mathrm{supp}(\phi) \sqcup \{p \}$ and $p \in D$, we deduce that  $\mathcal{H}(A \cup D \cup \mathrm{supp}(\phi)) = \mathcal{H}(A \cup D \cup \mathrm{supp}(\psi))$, so that $$d(w_a,m)= d(w_a,w_d)+ \# \mathrm{supp}(\psi)- \# \mathrm{supp}(\phi)=k-1$$
according to Lemma \ref{lem:wrdist}.  
\end{itemize}
Next, because the cases where the edge between $w_c$ and $w_d$ is labelled by a point of $\Omega$ are symmetric to the previous ones, let us suppose that the edges between $w_b$ and $w_c$, and $w_c$ and $w_d$, are not labelled by points of $\Omega$. In this situation, the maps $\phi$, $\psi$ and $ \varphi$ are all equal, say to some $\zeta \in G^{(\Omega)}$. In particular, the wreaths $w_b$, $w_c$ and $w_d$ belong to the same leaf $\mathfrak{W}(\zeta)$. According to Lemma \ref{lem:Wleaf}, if $Q$ denotes the convex hull of $A \cup \mathrm{supp}(\zeta)$ then $q=(Q,\zeta)$ is the projection of $w_a$ onto $\mathfrak{W}(\zeta)$. On the other hand, we know from Proposition \ref{prop:FSmedian} and Lemma \ref{lem:Wleaf} that the subgraph $\mathfrak{W}(\zeta)$ is median, so the median point $m=(M,\zeta)$ of $(Q,\zeta)$, $w_b$ and $w_d$ is well-defined. Notice that 
$$\begin{array}{lcl} d(m,w_b) & = & d(q,w_b)-d(q,m)=d(w_a,w_b)-d(w_a,q)-d(q,m) \\ \\ & = & d(w_a,w_d)-d(w_a,q)-d(q,m)= d(q,w_d)-d(q,m) \\ \\ & = & d(m,w_d) \end{array}$$
On the other hand, $d(w_b,w_d)=d(w_b,m)+d(m,w_d)$ and $1 \leq d(w_b,w_d) \leq 2$. It follows that $d(m,w_b)=d(m,w_d)=1$, ie., $m$ is adjacent to both $w_b$ and $w_d$. Moreover, 
$$\begin{array}{lcl} d(w_a,m) & = & d(w_a,q)+d(q,m) = d(w_a,q)+d(q,w_b)-d(w_b,m) \\ \\ & = & d(w_a,w_b)-d(w_b,m)=k-1 \end{array}$$
This concludes the proof of the quadrangle condition.

\medskip \noindent
Finally, we want to prove that $\mathfrak{W}$ contains no induced subgraphs isomorphic to $K_4^-$ or $K_{3,2}$. Let $Y \subset \mathfrak{W}$ be a subgraph isomorphic to $K_4^-$. If $(A,\varphi)$ and $(B,\psi)$ are the two non adjacent vertices of $Y$, it follows from Lemma \ref{lem:Wtriangle} below that $A=B$ and that $\varphi$ and $\psi$ differ on a single vertex $p \in A=B$. We deduce that $(A, \varphi)$ and $(B,\psi)$ must be adjacent in $\mathfrak{W}$, so that $Y$ is not an induced subgraph. Now, suppose that $Y$ is isomorphic to $K_{3,2}$. It follows from Lemma \ref{lem:WK32} below that $Y$ cannot be an induced subgraph. This concludes the proof of our proposition. 
\end{proof}

\begin{lemma}\label{lem:Wtriangle}
Let $(A,\phi)$, $(B, \varphi)$ and $(C, \psi)$ be three pairwise adjacent wreaths. Then $A=B=C$ and there exists some vertex $p \in A=B=C$ such that the maps $\phi$, $\varphi$ and $\psi$ differ only at $p$.
\end{lemma}

\begin{proof}
According to Lemma \ref{lem:Worbit}, there exists an isometry $g \in G \wr_{\Omega} H$ sending $(A,\phi)$ to $(A,1)$. Because the distances from $(B,\varphi)$ and $(C,\psi)$ to $(A,\phi)$ are equal, we deduce from Corollary \ref{cor:distadj} that the edge between $g \cdot (B,\varphi)$ and $g \cdot (C,\psi)$ is labelled by some $p \in \Omega$. Because the action of $G \wr_{\Omega} H$ on the graph of wreaths preserves the types of the labellings (ie., either a halfspace or a point of $\Omega$), we deduce that the edge between $(B,\varphi)$ and $(C,\psi)$ is labelled by $p$ as well. By symmetry, we also know that the edges between $(A, \phi)$ and $(B,\varphi)$, and $(A,\phi)$ and $(C,\psi)$, are labelled by points $q,r \in \Omega$ respectively. A fortiori, $A=B=C$. Moreover, since $(A,\phi)$, $(B, \varphi)$ and $(C, \psi)$ are pairwise adjacent, necessarily $p,q,r \in A=B=C$. Next, we know that $\varphi$ differ from $\phi$ only at $q$, and from $\psi$ only at $p$. Therefore, $\phi$ and $\psi$ may differ only on $\{p,q \}$. On the other hand, we know that $\phi$ and $\psi$ differ only at $r$, hence $r \in \{p,q\}$. By symmetry, we also know that $p \in \{ q,r \}$ and $q \in \{p,r\}$. If $r=p$, then we deduce from $q \in \{ p,r\}$ that $p=q=r$; and if $r=q$, we deduce from $p \in \{q,r\}$ that $p=q=r$. Therefore, $\phi$, $\varphi$ and $\psi$ differ at a single point $p=q=r$. This concludes the proof. 
\end{proof}

\begin{lemma}\label{lem:WK32}
Between two wreaths at distance two appart in the graph of wreaths $\mathfrak{W}$, there exist at most two geodesics.
\end{lemma}

\begin{proof}
Let $(A, \phi)$ and $(B,\varphi)$ be two wreaths at distance two appart in $\mathfrak{W}$. According to Lemma \ref{lem:Worbit}, we can suppose without loss of generality that $\phi=1$. Notice that, because $(A,1)$ and $(B,\varphi)$ are linked by a path of length two, necessarily 
\begin{itemize}
	\item[(i)] either $\# \mathrm{supp}(\varphi)=2$ and $B=A$;
	\item[(ii)] or $\# \mathrm{supp}(\varphi)=1$ and $B=[A,H]$ for some halfspace $H$;
	\item[(iii)] or $\mathrm{supp}(\varphi)= \emptyset$ and $B=[A,H_1,H_2]$ for some halfspaces $H_1,H_2$. 
\end{itemize}
In case $(i)$, if we denote by $\varphi_1$ (resp. $\varphi_2$) the map equal to $\varphi(p)$ at $p$ (resp. $\varphi(q)$ at $q$) and $1$ otherwise, then
$$(A,1), \ (A, \varphi_1), \ (A, \varphi_1 \varphi_2) \ \text{and} \ (A,1), \ (A, \varphi_2), \ (A,\varphi_1 \varphi_2)$$
are the only possible geodesics from $(A,1)$ to $(B,\varphi)$ in $\mathfrak{W}$, if they are well-defined. In case $(ii)$, 
$$(A,1), \ (A, \varphi), \ ([A,H],\varphi) \ \text{and} \ (A,1), \ ([A,H],1), \ ([A,H], \varphi)$$
are the only possible geodesics from $(A,1)$ to $(B,\varphi)$ in $\mathfrak{W}$, if they are well-defined. And finally, in case $(iii)$, 
$$(A,1), \ ([A,H_1],1), \ ([A,H_1,H_2],1) \ \text{and} \ (A,1), \ ([A,H_2],1), \ ([A,H_2,H_1],1)$$
are the only possible geodesics from $(A,1)$ to $(B,\varphi)$ in $\mathfrak{W}$, if they are well-defined. Thus, there exist at most two geodesics between $(A,1)$ and $(B,\varphi)$. 
\end{proof}

\paragraph{Cliques of $\mathfrak{W}$.} For convenience, let us introduce the following notation. If $\varphi \in G^{(\Omega)}$, $p \in \Omega$ and $g \in G$, set
$$\tau_p^g \varphi : q \mapsto \left\{ \begin{array}{cl} \varphi(q) & \text{if $q \neq p$} \\ g & \text{if $q=p$} \end{array} \right. .$$
The following lemma classifies completely the cliques of $\mathfrak{W}$. As a consequence, a clique of $\mathfrak{W}$ is naturally labelled by a halfspace (and its underlying hyperplane) or by a point of $\Omega$. 

\begin{lemma}\label{lem:Wclique}
A clique of $\mathfrak{W}$ is either an edge labelled by a halfspace or the subgraph generated by a collection of vertices $\left\{ (C, \tau_p^g\varphi) \mid  g \in G \right\}$ where $C \in \mathcal{FS}(X)$, $p \in C$ and $\varphi \in G^{(\Omega)}$ are fixed.
\end{lemma}

\begin{proof}
An edge labelled by a halfspace is a complete subgraph of $\mathfrak{W}$, and according to Lemma \ref{lem:Wtriangle}, it does not belong to any triangle. Therefore, it is a clique. Next, let $Q$ be a clique containing an edge labelled by a point $p \in \Omega$, say linking two wreaths $(C,\varphi)$ and $(C, \tau_p^g \varphi)$. Notice that necessarily $p \in C$. If $(A,\psi)$ belongs to $Q$, then $(A,\psi)$, $(C,\varphi)$ and $(C,\tau_p^g \varphi)$ define a triangle in $\mathfrak{W}$, so that it follows from Lemma \ref{lem:Wtriangle} that $A=C$ and $\psi= \tau_p^h \varphi$ for some $h \in G$. Therefore, $Q \subset \left\{ (C, \tau_p^g\varphi) \mid  g \in G \right\}$. On the other hand, $\left\{ (C, \tau_p^g\varphi) \mid  g \in G \right\}$ is a complete subgraph of $\mathfrak{W}$, hence $Q= \left\{ (C, \tau_p^g\varphi) \mid  g \in G \right\}$ by maximality. This concludes the proof.
\end{proof}

\begin{cor}\label{cor:Wcliquestab}
Let $Q= \left\{ (C, \tau_p^g \varphi) \mid g \in G \right\}$ be a clique. Then
$$\mathrm{stab}(Q) = \left\{  \left(  h , \tau_p^k ( \varphi(\cdot) \varphi(h^{-1} \cdot)^{-1}) \right) \mid h \in \mathrm{stab}(p) \cap \mathrm{stab}(C), k \in G \right\}.$$
\end{cor}

\begin{proof}
Let $(h,\psi) \in G \wr_{\Omega} H$. Clearly, $(h,\psi)$ belongs to $\mathrm{stab}(Q)$ if and only if, for every $g \in G$, there exists some $j \in G$ such that 
$$(h, \psi) \cdot (C, \tau_p^g \varphi) = (C, \tau_p^j \varphi),$$
ie., $hC=C$ and $\psi( \cdot) (\tau_p^g \varphi)(h^{-1} \cdot) = \tau_p^j\varphi( \cdot).$

\medskip \noindent
Notice that, if $hp \neq p$, then the previous equality evaluating at $hp$ provides $\psi(hp)g= \varphi(hp)$. Therefore, if $h$ does not belong to $\mathrm{stab}(p)$ but $(h, \psi)$ belongs to $\mathrm{stab}(Q)$, then, for every $g \in G$, there must exist some $j \in G$ such that $g = \psi(hp)^{-1} \varphi(hp)$, which is impossible since the right-hand side of the equality does not depend on $g$. 

\medskip \noindent
Consequently, $(h,\psi)$ belongs to $\mathrm{stab}(Q)$ if and only if $hC=C$, $hp=p$ and, for every $g \in G$, there exists some $j \in G$ such that $\psi(p) \cdot g = j$ and $\psi(q) \cdot \varphi(h^{-1}q)= \varphi(q)$ for every $q \neq p$. Notice that the second equality determines $\psi$ on $\Omega \backslash \{ p \}$, and that the first equality, once $\psi(p)$ fixed, determined $j$. Therefore, $(h,\psi)$ belongs to $\mathrm{stab}(Q)$ if and only if $h \in \mathrm{stab}(C) \cap \mathrm{stab}(p)$ and $\psi( \cdot ) = \varphi( \cdot ) \varphi(h^{-1} \cdot)^{-1}$ on $\Omega \backslash \{ p \}$. This concludes the proof.
\end{proof}

\paragraph{Hyperplanes of $\mathfrak{W}$.} As a consequence of Lemma \ref{lem:Wclique}, the edges of a given clique have the same label. It also follows from the lemma below that two edges which are the opposite sides of some square have the same label as well. Therefore, all the edges of a given hyperplane have the same label, so that a hyperplane is naturally labelled by a halfspace (and its underlying hyperplane) or a point of $\Omega$. 

\begin{lemma}\label{lem:Wsquare}
 Let $K$ be a square of $\mathfrak{W}$ with $(A,1)$ as one of its vertices. Three cases may occur:
\begin{itemize}
	\item either there exist distinct $p,q \in A \cap \Omega$ such that the vertices of $K$ are $(A,1)$, $(A,\varphi_1)$, $(A, \varphi_2)$ and $(A,\varphi_1 \varphi_2)$ for some $\varphi_1, \varphi_2 \in G^{(\Omega)}$ satisfying $\mathrm{supp}(\varphi_1)=\{p\}$ and $\mathrm{supp}(\varphi_2)= \{q \}$;
	\item or there exist some point $p \in A \cap \Omega$ and some halfspace $H$ of $X$ such that the vertices of $K$ are $(A,1)$, $(A, \varphi)$, $([A,H],1)$ and $([A,H],\varphi)$;
	\item or there exist two halfspaces $H_1,H_2$ of $X$ such that the vertices of $K$ are $(A,1)$, $([A,H_1],1)$, $([A,H_2],1)$ and $([A,H_1,H_2],1)=([A,H_2,H_1],1)$.
\end{itemize}
\end{lemma}

\begin{proof}
$(B,\varphi)$. Notice that, because $(A,1)$ and $(B,\varphi)$ are linked by a path of length two, necessarily 
\begin{itemize}
	\item[(i)] either $\# \mathrm{supp}(\varphi)=2$ and $B=A$;
	\item[(ii)] or $\# \mathrm{supp}(\varphi)=1$ and $B=[A,H]$ for some halfspace $H$;
	\item[(iii)] or $\mathrm{supp}(\varphi)= \emptyset$ and $B=[A,H_1,H_2]$ for some halfspaces $H_1,H_2$. 
\end{itemize}
In case $(i)$, if we denote by $\varphi_1$ (resp. $\varphi_2$) the map equal to $\varphi(p)$ at $p$ (resp. $\varphi(q)$ at $q$) and $1$ otherwise, then
$$(A,1), \ (A, \varphi_1), \ (A, \varphi_1 \varphi_2) \ \text{and} \ (A,1), \ (A, \varphi_2), \ (A,\varphi_1 \varphi_2)$$
are the only possible paths from $(A,1)$ to $(B,\varphi)$ in $\mathfrak{W}$, if they are well-defined. On the other hand, we know that there exist at least two such paths, so they must be well-defined, and this provides the vertices of our square $K$. In case $(ii)$, 
$$(A,1), \ (A, \varphi), \ ([A,H],\varphi) \ \text{and} \ (A,1), \ ([A,H],1), \ ([A,H], \varphi)$$
are the only possible paths from $(A,1)$ to $(B,\varphi)$ in $\mathfrak{W}$, if they are well-defined. On the other hand, we know that there exist at least two such paths, so they must be well-defined, and this provides the vertices of our square $K$. Finally, in case $(iii)$, 
$$(A,1), \ ([A,H_1],1), \ ([A,H_1,H_2],1) \ \text{and} \ (A,1), \ ([A,H_2],1), \ ([A,H_2,H_1],1)$$
are the only possible paths from $(A,1)$ to $(B,\varphi)$ in $\mathfrak{W}$, if they are well-defined. On the other hand, we know that there exist at least two such paths, so they must be well-defined, and this provides the vertices of our square $K$.
\end{proof}

\noindent
The next proposition describes completely the hyperplanes of $\mathfrak{W}$ which are labelled by points of $\Omega$. Notice that, as a consequence, for every $p \in \Omega$ there exists a unique hyperplane labelled by $p$. 

\begin{prop}\label{prop:WhypOmega}
Let $J$ be a hyperplane of $\mathfrak{W}$ labelled by some vertex $p \in \Omega$. An edge linking two wreaths $(C_1, \varphi_1)$ and $(C_2, \varphi_2)$ is dual to $J$ if and only if $C_1=C_2$, $p \in C_1=C_2$, and $\varphi_1$ and $\varphi_2$ differ only at $p$. In particular, $N(J)$ is the subgraph of $\mathfrak{W}$ generated by $\{ (C, \varphi) \mid p \in C \}$. Moreover, the fibers of $J$ are the subgraphs generated by $\{ (C, \varphi) \mid p \in C, \ \varphi(p)=g\}$ where $g \in G$. 
\end{prop}

\begin{proof}
Fix some clique $Q$ dual to $J$, say $Q= \left\{ (C, \tau_p^g \varphi) \mid g \in G \right\}$ where $\varphi \in G^{(\Omega)}$ is a map and $C \in \mathcal{FS}(X)$ a subcomplex satisfying $p \in C$. 

\medskip \noindent
We claim that, for every $(A,\psi) \in \mathfrak{W}$, the projection of $(A,\psi)$ onto $Q$ is $(C, \tau_p^{\psi(p)} \varphi)$. Indeed, it follows from Lemma \ref{lem:wrdist} that, for every $g \in G$, the distance in $\mathfrak{W}$ between $(A, \psi)$ and $(C, \tau_p^g \varphi)$ is equal to 
$$2 \cdot \# \mathcal{H} \left( A \cup C \cup \mathrm{supp}(\psi^{-1} \tau_p^g \varphi) \right) - \# \mathcal{H}(A) - \# \mathcal{H}(C) + \# \mathrm{supp}( \psi^{-1} \tau_p^g \varphi).$$
But $\mathcal{H} \left( A \cup C \cup \mathrm{supp}(\psi^{-1} \tau_p^g \varphi) \right)$ does not depend on $g$ since $p$ belongs to $C$, so the distance is minimal precisely when $\mathrm{supp}(\psi^{-1} \tau_p^g \varphi)$ is minimal, which happens when $g= \psi(p)$. This proves our claim.

\medskip \noindent
As a consequence, two wreaths $(A_1, \varphi_1)$ and $(A_2, \varphi_2)$ are separated by $J$ if and only if $\varphi_1(p) \neq \varphi_2(p)$. In particular, if $(A_1, \varphi_1)$ and $(A_2, \varphi_2)$ are adjacent, they are separated by $J$ (ie., the edge linking them belongs to $J$) if and only $A_1=A_2$ and $\varphi_1, \varphi_2$ differ only at $p$. The additional condition $p \in A_1=A_2$ is automatic, since otherwise $(A_1, \varphi_1)$ and $(A_2, \varphi_2)$ would not be adjacent. Consequently, a wreath $(A, \psi)$ is an endpoint of an edge dual to $J$ (ie., belongs to $N(J)$) if and only if it is a neighbor of $(A, \tau_p^g \psi)$ for every $g \in G$, which is equivalent to $p \in A$. A fortiori, $N(J)$ is the subgraph of $\mathfrak{W}$ generated by $\left\{ (A, \psi) \mid p \in A \right\}$, and it follows from our description of the projection onto $Q$ that the fibers of $J$ are the subgraphs generated by $\left\{ (A, \psi) \mid p \in A, \psi(p)=g \right\}$ where $g \in G$. 
\end{proof}

\begin{cor}\label{cor:Wstabhyp}
Let $J$ be a hyperplane of $\mathfrak{W}$ labelled by some vertex $p \in \Omega$. Then 
$$\mathrm{stab}(J)= \left\{ ( h, \psi) \in G \wr_{\Omega} H \mid h \in \mathrm{stab}(p) \right\}.$$
\end{cor}

\begin{proof}
Fix some $(h,\psi) \in G \wr_{\Omega} H$. As a consequence of Proposition \ref{prop:WhypOmega}, $(h, \psi)$ belongs to $\mathrm{stab}(J)$ if and only if, for every subcomplex $C \in \mathcal{FS}(X)$ containing $p$ and every $\varphi \in G^{(\Omega)}$, there exist a subcomplex $C' \in \mathcal{FS}(X)$ containing $p$ and a map $\phi \in G^{(\Omega)}$ such that $( h, \psi) \cdot (C, \varphi) = (C', \phi).$ Because $( h, \psi) \cdot (C, \varphi) = (hC, \psi(\cdot) \varphi(h^{-1} \cdot))$, we deduce that $(h, \psi)$ belongs to $\mathrm{stab}(J)$ if and only if $p \in hC$ for every subcomplex $C \in \mathcal{FS}(X)$ containing $p$. Since $\{p \}$ is such a subcomplex, we deduce that this condition is also equivalent to $h \in \mathrm{stab}(p)$.
\end{proof}

\noindent
About hyperplanes of $\mathfrak{W}$ labelled by halfspaces, a precise description is not necessary for our purpose. We only prove the following lemma:

\begin{lemma}\label{lem:Wstabfiber}
Let $J$ be a hyperplane of $\mathfrak{W}$ labelled by some halfspace $H$ of $X$. Then $J$ has two fibers, and any element of $\mathrm{stab}(J)$ stabilises them.
\end{lemma}

\begin{proof}
It follows from Lemma \ref{lem:Wclique} that a clique labelled by some hyperplane of $X$ has cardinality two, so $J$ has necessarily two fibers. If $(A,\phi)$ and $([A,H],\phi)$ are the endpoints of some edge $e$ dual to $J$, we can naturally orient $e$ from $(A,\phi)$ to $([A,H],\phi)$ if $\# \mathcal{H}(A)< \# \mathcal{H}([A,H],\phi)$, and from $([A,H],\phi)$ to $(A,\phi)$ otherwise. It is clear that $\mathrm{stab}(J)$ acts on the edges dual to $J$ by preserving their orientations. On the other hand, it follows from Lemma \ref{lem:Wsquare} that two opposite edges in a square have the same orientation, so we can denote the two fibers of $J$ as $F^-$ and $F^+$ so that the orientation of any edge dual to $J$ correspond to going from $F^-$ to $F^+$. Therefore, $\mathrm{stab}(J)$ must stabilise both $F^-$ and $F^+$. 
\end{proof}

\subsection{Cubulating wreath products}

\noindent
So far, we have constructed actions of some wreath products on quasi-median graphs. In this section, our goal is to prove the these actions are topical-transitive and to apply the results of Sections \ref{section:topicalactionsI} and \ref{section:topicalactionsII}.

\medskip \noindent
Fix two groups $G$ and $H$, and suppose that $H$ acts on a CAT(0) cube complex $X$. Let $\Omega \subset X$ denote an $H$-invariant collection of vertices, and fix $\Gamma \subset \Omega$ a set of representatives. Let $\mathcal{C}_1$ denote the collection of the cliques 
$$Q_p= \left\{ ( \{p \}, \tau_p^g 1) \mid g \in G \right\}, \ p \in \Gamma,$$
and $\mathcal{C}_2$ any collection of cliques labelled by halfspaces such that an orbit of hyperplane intersects $\mathcal{C}_2$ in at most one clique. Set $\mathcal{C}= \mathcal{C}_1 \sqcup \mathcal{C}_2$.

\begin{prop}\label{prop:wreathtopical}
Suppose that the stabiliser of every point of $\Omega$ is trivial. The permutational wreath product $G \wr_{\Omega} H$ acts topically-transitively on $\mathfrak{W}$. 
\end{prop}

\begin{proof}
Let $Q$ be a clique of $\mathfrak{W}$ and let $J$ denote the hyperplane dual to it. If $Q$ is labelled by a halfspace, it follows from Lemma \ref{lem:Wstabfiber} that $\mathrm{stab}(Q)=\mathrm{fix}(Q)$. It also follows that $\mathrm{stab}(J)$ stabilises the fibers of $J$, so the image of $\mathrm{stab}(J)$ under $\rho_Q$ is trivial, and a fortiori included into $\rho_Q(\mathrm{stab}(Q))$. 

\medskip \noindent
Otherwise, we know from Lemma \ref{lem:Wclique} that $Q$ can be written as $\left\{ (C,\tau^g_p \varphi) \mid g \in G \right\}$ for some $C \in \mathcal{FS}(X)$, $p \in C$ and $\varphi \in G^{(\Omega)}$, and it follows from Corollary \ref{cor:Wcliquestab} that $\mathrm{stab}(Q)= \left\{ \left( 1, \tau_p^k(1) \right) \mid k \in G \right\}$ since $\mathrm{stab}(p)$ is trivial by assumption. Consequently, $\mathrm{stab}(Q)$ acts freely and transitively on the vertices of $Q$. Next, according to Proposition \ref{prop:WhypOmega}, the fibers of $J$ are the subgraphs generated by $\left\{ (C, \varphi ) \mid \varphi(p)=g \right\}$ where $g \in G$; and according to Corollary \ref{cor:Wstabhyp}, the stabiliser of $J$ is $\left\{ (h,\psi) \mid h \in \mathrm{stab}(p) \right\}$. Consequently, the fibers of $J$ are labelled by $G$, and the action of an element $(h,\psi) \in \mathrm{stab}(J)$ on the set of the fibers of $J$ corresponds to the left-multiplication by $\psi(p)$. As a consequence, the element $\left( 1, \tau_p^{\psi(p)} 1 \right) \in \mathrm{stab}(Q)$ induces the same permutation on the set of the fibers of $J$ as $(h,\psi)$. Therefore, the image of $\mathrm{stab}(J)$ under $\rho_Q$ is included into $\rho_Q(\mathrm{stab}(C))$. 
\end{proof}

\noindent
We need two last preliminary lemmas before giving applications to wreath products.

\begin{lemma}\label{lem:Wvertexstab}
For every wreath $(C, \varphi)$,
$$\mathrm{stab}((C, \varphi)) = \left\{ \left( h, \varphi(\cdot) \varphi(h^{-1} \cdot)^{-1} \right) \mid h \in \mathrm{stab}(C) \right\}.$$
\end{lemma}

\begin{proof}
Let $(C,\varphi) \in \mathfrak{W}$ be a wreath and $(h, \psi) \in G \wr_{\Omega} H$. Then $(h,\psi)$ belongs to the stabiliser of $(C,\varphi)$ if and only if
$$\left( hC, \psi(\cdot) \varphi(h^{-1} \cdot) \right) = (h,\psi) \cdot (C, \varphi)= (C, \varphi),$$
ie., $hC=C$ and $\psi(\cdot) = \varphi(\cdot) \varphi(h^{-1} \cdot)^{-1}$. Thus, 
$$\mathrm{stab}((C, \varphi)) = \left\{ \left( h, \varphi(\cdot) \varphi(h^{-1} \cdot)^{-1} \right) \mid h \in \mathrm{stab}(C) \right\},$$
which concludes the proof.
\end{proof}

\begin{lemma}\label{lem:Wfiniteness}
The following assertions hold:
\begin{itemize}
	\item if $X$ is locally finite, any vertex of $\mathfrak{W}$ belongs to finitely many cliques;
	\item if $H \curvearrowright X$ is properly discontinuous, vertex-stabilisers of $\mathfrak{W}$ are finite.
\end{itemize}
\end{lemma}

\begin{proof}
Let $(C, \varphi) \in \mathfrak{W}$ be a wreath. The cliques of $\mathfrak{W}$ labelled by points of $\Omega$ which contain $(C, \varphi)$ are the
$$\left\{ (C, \tau_p^g \varphi) \mid g \in G \right\},$$
where $p \in C$. Thus, there exist at most $\# C$ such cliques. Next, the cliques of $\mathfrak{W}$ labelled by halfspaces which contain $(C,\varphi)$ are the edges $[(C,\varphi),(C', \varphi)]$, where $C'$ is a neighbor of $C$ in the graph of finite subcomplexes $\mathcal{FS}(X)$. As a consequence of Proposition \ref{prop:FSmedian}, there exist at most finitely many such cliques if $X$ is locally finite. This concludes the proof of the first statement of our lemma. The second statement follows directly from Lemma \ref{lem:Wvertexstab}.
\end{proof}

\noindent
Thanks to Proposition \ref{prop:wreathtopical} and Lemma \ref{lem:Wfiniteness} above, Propositions \ref{prop:CAT0metricallyproper}, \ref{aTmenablegroups} and \ref{prop:properlydiscontinuous} apply, hence:

\begin{thm}\label{thm:Wfinal}
Suppose that $H \curvearrowright X$ is properly discontinuous and that every point of $\Omega$ has trivial stabiliser. The following statements hold:
\begin{itemize}
	\item if $X$ is locally finite and if $G$ acts metrically properly on a CAT(0) cube complex, then $G \wr_{\Omega} H$ acts metrically properly on some CAT(0) cube complex;
	\item if $G$ acts properly discontinuously on a CAT(0) cube complex, then so does $G \wr_{\Omega} H$;
	\item if $X$ is locally finite and if $G$ is a-T-menable (resp. a-$L^p$-menable, where $p \notin 2 \mathbb{Z}$), then $G \wr_{\Omega} H$ is a-T-menable (resp. a-$L^p$-menable).
\end{itemize}
\end{thm}

\noindent
As a consequence, we are able to reprove \cite[Theorem 5.C.3]{CornulierCommensurated}.

\begin{cor}\label{cor:Wproper}
A wreath product of groups acting properly discontinuously on CAT(0) cube complexes acts properly discontinuously on a CAT(0) cube complex.
\end{cor}

\begin{proof}
Let $H$ be a group acting properly discontinuously on some CAT(0) cube complex $X$. According to Lemma \ref{lem:modifycubing}, we can suppose without loss of generality that $X$ contains a vertex $x_0$ with trivial stabiliser. Set $\Omega = H \cdot x_0$. If $G$ is another group acting properly discontinuously on some CAT(0) cube complex, we deduce from the previous theorem that the wreath product $G \wr_{\Omega} H = G \wr H$ acts properly discontinuously on some CAT(0) cube complex as well.
\end{proof}

\begin{remark}\label{remark:explicitW}
The previous proof is not quite explicit, because we applied general constructions and results from previous sections, but interestingly the CAT(0) cube complex we obtain can be constructed directly as follows. Let $G,H$ be two groups acting on two CAT(0) cube complexes $Y,X$ respectively. Suppose that $X$ contains a vertex $x_0$ with trivial stabiliser, and set $\Omega = H \cdot x_0$. Define $W$ as the graph whose vertices are the pairs $(C, \varphi)$, where $C \in \mathcal{FS}(X)$ is a subcomplex of $X$ and $\varphi : \Omega \to Y^{(0)}$ a map with finite support; and whose edges link two pairs $(C_1, \varphi_1)$ and $(C_2, \varphi_2)$ if either $C_1=C_2$ and $\varphi_1, \varphi_2$ differ on a single point $p \in \Omega \cap C_1$ in such a way that $\varphi_1(p)$ and $\varphi_2(p)$ are two adjacent vertices of $Y$, or $C_1,C_2$ are adjacent in $\mathcal{FS}(X)$ and $\varphi_1=\varphi_2$. Then $W$ is the median graph (or equivalently, the CAT(0) cube complex) which is constructed in the proof of Corollary \ref{cor:Wproper}. We make this construction explicit in \cite{aTmW}.
\end{remark}

\begin{remark}
It also follows from Theorem \ref{thm:Wfinal} that, if $G$ is a (discrete) a-T-menable group and $H$ a group acting properly discontinuously on some CAT(0) cube complex, then the wreath product $G \wr H$ is a-T-menable. However, a much more general conclusion already holds according to \cite{CSVgeneral}, since it is proved there that a-T-menability is stable under any wreath products. (This generalises our previous observation because, according to \cite{Haagerup}, a group acting properly discontinuously on some CAT(0) cube complex is necessarily a-T-menable.) In fact, by modifying the construction of the graph of wreaths, it is possible to construct a median space when starting from two median spaces instead of CAT(0) cube complexes, which allows to reprove that the wreath product of two (discrete) a-T-menable groups is a-T-menable thanks to \cite{medianviewpoint}. See \cite{aTmW} for more details.
\end{remark}

\noindent
From now on, fix a group $H$ acting properly discontinuously on some locally finite CAT(0) cube complex $X$, the orbit $\Omega$ of some vertex $x_0 \in X$ which has trivial stabiliser, and a finite group $F$. In the graph of wreaths $\mathfrak{W}$ associated to the wreath product $F \wr H$, there are two classes of hyperplanes, which are both $(F \wr H)$-invariant: those labelled by halfspaces of $X$ and those labelled by vertices of $\Omega$. Let $\mathfrak{W}_c$ denote the quasi-median graph obtained by quasi-cubulating the space of partitions defined by $\mathfrak{W}$ and its collection of hyperplanes labelled by halfspaces of $X$ (see Section \ref{section:spaceswithpartitions} for definitions). Roughly speaking, $\mathfrak{W}_c$ is obtained from $\mathfrak{W}$ by cutting along the hyperplanes labelled by vertices of $\Omega$ and next identifying all the fibers of a given hyperplane together; in the context of CAT(0) cube complexes, such an operation is referred to as a \emph{restriction quotient} in \cite{MR2827012}. Because all the triangles of $\mathfrak{W}$ are dual to hyperplanes labelled by vertices of $\Omega$, according to Lemma \ref{lem:Wtriangle}, we know that the hyperplanes of $\mathfrak{W}_c$ delimit precisely two sectors (since so do the hyperplanes of $\mathfrak{W}$ labelled by halfspaces of $X$), so that $\mathfrak{W}_c$ must be triangle-free. As a consequence of Corollary \ref{cor:whenmedian}, it follows that $\mathfrak{W}_c$ is a median graph, or equivalently, a CAT(0) cube complex. This observation is independent of the assumptions that $F$ is finite and $X$ locally finite. Now we claim that the induced action of $F \wr H$ on $\mathfrak{W}_c$ is metrically proper. According to Lemma \ref{lem:distspacepartitions}, the distance in $\mathfrak{W}_c$ between two wreaths $(C_1, \varphi_1)$ and $(C_2, \varphi_2)$ is equal to the number of hyperplanes of $\mathfrak{W}$ labelled by halfspaces of $X$ which separate them, so that it follows from the proof of Lemma \ref{lem:wrdist} that
$$d_{\mathfrak{W}_c} \left( (C_1, \varphi_1),(C_2, \varphi_2) \right) = 2 \cdot \# \mathcal{H} \left( C_1 \cup C_2 \cup \mathrm{supp}(\varphi_1^{-1} \varphi_2)  \right) - \# \mathcal{H}(C_1) - \# \mathcal{H}(C_2).$$
In particular, for every $w=( g, \varphi) \in F \wr H$,
$$d_{\mathfrak{W}_c} \left( w \cdot (\{x_0 \},1), (\{x_0\},1) \right) = 2 \cdot \# \mathcal{H} \left( \{ x_0, gx_0\} \cup \mathrm{supp}(\varphi) \right).$$
Therefore, if $d_{\mathfrak{W}_c} \left( w \cdot (\{x_0 \},1), (\{x_0\},1) \right) \leq R$ for some $R \geq 0$, then $\{gx_0 \}$ and $\mathrm{supp}(\varphi)$ must be included into the ball $B(x_0,R)$, which is finite since $X$ is locally finite. So there are finitely many choices on $gx_0$ and $\mathrm{supp}(\varphi)$. Because $H$ acts properly discontinuously on $X$, we deduce that there exist finitely many choices on $g$; and because $F$ is finite, we deduce that there exist finitely many choices on $\varphi$. Thus, we have proved that the action of $F \wr H$ on $\mathfrak{W}_c$ is metrically proper. 

\medskip \noindent
The point is that $\mathfrak{W}_c$ is smaller than the cube complex provided by Theorem \ref{thm:Wfinal}, and interestingly, it turns out to be finite-dimensional in some cases. 

\begin{lemma}
The dimension of $\mathfrak{W}_c$ is equal to the dimension of $\mathcal{FS}(X)$. 
\end{lemma}

\begin{proof}
Suppose that $\mathfrak{W}_c$ contains $n$ pairwise transverse hyperplanes. According to Theorem \ref{thm:quasicubulation}, $\mathfrak{W}$ must contain $n$ pairwise transverse hyperplanes labelled by halfspaces of $X$, and the intersection of the neighborhouds of these hyperplanes must contain a prism $P$ of cubical dimension $n$ according to Proposition \ref{prop:transversehypcube}. Because the edges of $P$ are labelled by halfspaces of $X$, necessarily $P$ must be contained into some leaf $\mathfrak{W}(\varphi)$. But $\mathfrak{W}(\varphi)$ is isomorphic to $\mathcal{FS}(X)$ according to Lemma \ref{lem:Wleaf}. This proves the inequality
$$\dim \mathfrak{W}_c \leq \dim \mathcal{FS}(X).$$
The reverse inequality is clear because the leaves of $\mathfrak{W}$ isometrically embed into $\mathfrak{W}_c$ as no hyperplane of $\mathcal{W}$ labelled by a vertex of $\Omega$ intersects a leaf.
\end{proof}

\noindent
For instance, $\mathfrak{W}_c$ is finite-dimensional if $X= \mathbb{R}^n$ for some $n \geq 1$. Indeed, because a non empty finite convex subcomplex of $\mathbb{R}^n$ is a rectangular cuboid, an element of $\mathcal{FS}( \mathbb{R}^n)$ is determined by $2n$ integers $(x_1,y_1, \ldots, x_n,y_n)$ satisfying $x_i \leq y_i$ for every $1 \leq i \leq n$. As a consequence, $\mathcal{FS}(\mathbb{R}^n)$ is the subcomplex of $\mathbb{R}^{2n}$ generated by the vertices
$$\{ (x_1,y_1,\ldots, x_n,y_n) \mid x_i \leq y_i \ \text{for} \ 1 \leq i \leq n \}.$$
It follows from the previous lemma that $\dim \mathfrak{W}_c = 2n$. Thus, we have proved:

\begin{prop}\label{prop:FZnfinitedimension}
For every interger $n \geq 1$ and every finite group $F$, the wreath product $F \wr \mathbb{Z}^n$ acts metrically properly on a CAT(0) cube complex of dimension $2n$.
\end{prop}

\noindent
Since a finite-dimensional CAT(0) cube complex has $\ell^2$-compression one \cite{propertyA}, it would be interesting to know how much the orbit is distorded in this action. This is the purpose of the following lemma. 

\begin{lemma}\label{lem:Wcorbitcompression}
For its action on $\mathfrak{W}_c$, the wreath product $F \wr \mathbb{Z}^n$ has an orbit with compression at least $1/(2n-1)$. 
\end{lemma}

\begin{proof}
Our basepoint is $w_0= (\{x_0 \}, 1)$ where $x_0=0$. Fix some non trivial element $(g,\varphi) \in F \wr \mathbb{Z}^n$. As a finite generating set of $F \wr \mathbb{Z}^n$, we consider the union $F$ with the standard generating set of $\mathbb{Z}^n$. With respect to the associated word metric,
$$d_{F \wr \mathbb{Z}^n}((g,\varphi),1) = TS(x_0,\mathrm{supp}(\varphi), gx_0) + \# \mathrm{supp}(\varphi),$$
where $TS(x,F,y)$ denotes the length of a shortest path in $\mathbb{Z}^n$ starting from $x$, ending at $y$, and passing through each point of $F$, where $x,y \in \mathbb{Z}^n$ are vertices and $F \subset \mathbb{Z}^n$ a finite collection of vertices. First of all, notice that 
$$\# \mathrm{supp}(\varphi) \leq TS(x_0, \mathrm{supp}(\varphi), g \cdot x_0) +1.$$
Indeed, by indexing $\mathrm{supp}(\varphi)$ suitably, say as  $\{ f_1, \ldots, f_r\}$, one has
$$ \begin{array}{lcl} TS(x_0, \mathrm{supp}(\varphi), g \cdot x_0) & = &  \displaystyle d(x_0,f_1)+ \sum\limits_{i=1}^{r-1} d(f_i,f_{i+1}) + d(f_r,g \cdot x_0) \\ \\ & \geq & \displaystyle  \sum\limits_{i=1}^{r-1} d(f_i,f_{i+1}) \geq r-1= \# \mathrm{supp}(\varphi) -1 \end{array}$$
Next, Lemma \ref{lem:TSeuclid} will show later (and independently of our lemma) that, for every $\alpha < n/(2n-1)$, there exists some constant $C >0$ such that the inequality
$$2 \cdot \# \mathcal{H}( F \cup \{ x,y\}) + \# F \geq C \cdot TS(x,F,y)^{\alpha}$$
holds for every $x,y \in \mathbb{Z}^n$ and for every finite subset $F \subset \mathbb{Z}^n$. Therefore, fixing some $\alpha < n/(2n-1)$, there exists a constant $C>0$, depending only on $n$ and $\alpha$, such that
$$\begin{array}{lcl} d_{F \wr \mathbb{Z}^n}((g,\varphi),1) & = & TS(x_0, \mathrm{supp}(\varphi),g \cdot x_0) + \# \mathrm{supp}(\varphi) \\ \\ & \leq & 2 \cdot TS(x_0, \mathrm{supp}(\varphi), g \cdot x_0) +1 \\ \\ & \leq & C \cdot \left( \# \mathcal{H}( \mathrm{supp}(\varphi) \cup \{ x_0,g \cdot x_0 \})+ \# \mathrm{supp}(\varphi) \right)^{1/\alpha} +1 \end{array}$$
Now, notice that, for every rectangular cuboid $R$ of $\mathbb{Z}^n$, the inequality $V \leq H^n$ holds, where $V$ and $H$ denote respectively the numbers of vertices and hyperplanes of $R$. Consequently,
$$\# \mathrm{supp}(\varphi) \leq \# \mathcal{H}(\mathrm{supp}(\varphi))^n \leq \# \mathcal{H}( \mathrm{supp}(\varphi) \cup \{x_0, g \cdot x_0 \})^n,$$
which implies that
$$ \begin{array}{lcl} d_{F \wr \mathbb{Z}^n}( (g,\varphi),1) & \leq & 2^{1/\alpha} \cdot C \cdot \# \mathcal{H}( \mathrm{supp}(\varphi) \cup \{ x_0, g \cdot x_0 \})^{n/\alpha} +1 \\ \\ & \leq & 2^{1/\alpha} \cdot C \cdot d_{\mathfrak{W}_c} ( (g,\varphi) \cdot w_0,w_0) ^{n/\alpha}+1 \\ \\ & \leq & \left( 1+ 2^{1/\alpha} C \right) \cdot d_{\mathfrak{W}_c}( (g, \varphi) \cdot w_0,w_0)^{n/\alpha} \end{array}$$
This proves our lemma.
\end{proof}

\noindent
As a consequence, we are able to reprove a particular case of \cite[Corollary 1.2]{TesseraWP}:

\begin{cor}
For every finite group $F$, $\alpha_2 \left( F \wr \mathbb{Z} \right)=1$. 
\end{cor}

\noindent
Notice that the lower bound $\alpha_2 \left( F \wr \mathbb{Z}^2 \right) \geq 1/3$ we find is not optimal, since it is proved in \cite[Remark 3.4]{NaorPeres} that $\alpha_2(F \wr \mathbb{Z}^2)=1/2$.

\medskip \noindent
Finally, we conclude with a last consequence of our construction.

\begin{prop}
For every finite group $F$, the wreath product $F \wr \mathbb{Z}$ quasi-isometrically embeds into a product of two simplicial trees.
\end{prop}

\begin{proof}
According to Lemma \ref{lem:Wcorbitcompression}, $F\wr \mathbb{Z}$ quasi-isometrically embeds into the CAT(0) cube complex $\mathfrak{W}_c$, so it is sufficient to embed isometrically $\mathfrak{W}_c$ into a product of two simplicial trees. According to Theorem \ref{thm:quasicubulation}, the hyperplanes of $\mathfrak{W}_c$ are in bijection with the hyperplanes of $\mathfrak{W}$ which are labelled by halfspaces of $X$. Passing through such a hyperplane in $\mathfrak{W}$ corresponds to translating one side of the segment representing the first coordinate of a wreath; if this translation is made on the left-side (resp. right-side) of the segment, we say that the associated hyperplane is \emph{left} (resp. \emph{right}). Let $T_l$ (resp. $T_r$) denote the CAT(0) cube complex obtained by cubulating $\mathfrak{W}_c$ (or equivalently $\mathfrak{W}$) with respect to the collection of left hyperplanes (resp. right hyperplanes). By noticing thanks to the description of the squares of $\mathfrak{W}$ provided by Lemma \ref{lem:Wsquare} that two left hyperplanes, or two right hyperplanes, cannot be transverse, we deduce that $T_l$ and $T_r$ are trees. From the canonical maps $\eta_l : \mathfrak{W}_c \to T_l$ and $\eta_r : \mathfrak{W}_c \to T_r$, sending a vertex to the principal orientation it defines, we get a map $\eta = \eta_l \times \eta_r : \mathfrak{W}_c \to T_l \times T_r$. If, for every vertices $x,y \in \mathfrak{W}_c$, we denote by $h_l$ (resp. $h_r$) the number of left hyperplanes (resp. right hyperplanes) separating $x$ and $y$, we have
$$\begin{array}{lcl} d_{T_l \times T_r} (\eta(x),\eta(y)) & = & d_{T_l}(\eta_l(x), \eta_l(y)) + d_{T_r} (\eta_r(x),\eta_r(y)) \\ \\ & = & h_l(x,y)+h_r(x,y) = d_{\mathfrak{W}_c}(x,y) \end{array}$$
where we used Lemma \ref{lem:distspacepartitions} to determine the distances in $T_l$ and $T_r$. This concludes the proof.
\end{proof}

\subsection{Equivariant $\ell^p$-compression}

\noindent
Fix two groups $G$ and $H$ with two finite generating sets, and suppose that $H$ acts on a CAT(0) cube complex $X$ with a vertex $x_0$ of trivial stabiliser. Let $\mathfrak{W}$ be the graph of wreaths associated to $\Omega = H \cdot x_0$; notice that $G \wr_{\Omega} H = G \wr H$. In this section, we want to apply Proposition \ref{prop:equicompression} to study the equivariant $\ell^p$-compression of the wreath product $G \wr H$. The first step is to choose a convenient system of metrics on the graph of wreaths $\mathfrak{W}$. Let $Q$ be a clique of $\mathfrak{W}$. If $Q$ is labelled by a halfspace, set $\delta_Q : (x,y) \mapsto \left\{ \begin{array}{cl} 1 & \text{if $x \neq y$} \\ 0  & \text{otherwise} \end{array} \right.$. Otherwise, if $Q$ is labelled by some point $p \in \Omega$, we set $\delta_Q\left( (C,\varphi),(C, \phi) \right) = d_G \left( \varphi(p), \psi(p) \right)$ for every $(C,\varphi),(C, \phi) \in Q$.

\medskip \noindent
We claim that this system of metrics is coherent and invariant. First, let $Q$ be a clique of $\mathfrak{W}$, $x,y \in Q$ two vertices and $w=(k,\psi) \in G \wr H$. If $Q$ is labelled by a halfspace, it is clear that $\delta_{wQ}(w \cdot x , w \cdot y) = \delta_Q(x,y)$ because $\delta_Q$ and $\delta_{wQ}$ are two discrete metrics. So suppose that $Q$ is labelled by some point $p \in \Omega$. Notice that $wQ$ is labelled by $kp \in \Omega$. Writting $x=(C, \varphi)$ and $y=(C,\psi)$, we have
$$\begin{array}{lcl} \delta_{wQ}(w \cdot x , w \cdot y) & = & \delta_{wQ} \left( (kC, \psi(\cdot) \varphi(k^{-1} \cdot)), (kC, \psi(\cdot) \phi(k^{-1} \cdot)) \right) \\ \\ & = & d_G \left( \psi(kp) \varphi(p), \psi(kp) \phi(p) \right) = d_G \left( \varphi(p), \psi(p) \right) = \delta_Q(x,y) \end{array}$$
Therefore, our system is $(G \wr H)$-invariant. Next, let $Q_1$ and $Q_2$ be two cliques of $\mathfrak{W}$ dual to the same hyperplane $J$ and $x,y \in Q_1$ two vertices. For short, let $t$ denote the canonical bijection $t_{Q_1 \to Q_2} : Q_1 \to Q_2$. If $J$ is labelled by some halfspace of $X$, then clearly $\delta_{Q_2}(t(x),t(y)) = \delta_{Q_1}(x,y)$ since $\delta_{Q_1}$ and $\delta_{Q_2}$ are two discrete metrics. So suppose that $J$ is labelled by some $p \in \Omega$. As a consequence, there exist two subcomplexes $C_1,C_2 \in \mathcal{FS}(X)$ containing $p$ and two maps $\varphi_1,\varphi_2 \in G^{(\Omega)}$ such that 
$$Q_1 = \left\{ (C_1, \tau_p^g \varphi_1) \mid g \in G \right\} \ \text{and} \ Q_2 = \left\{ (C_2, \tau_p^g \varphi_2) \mid g \in G \right\}.$$
Write $x=(C_1, \tau_p^g \varphi_1)$ and $y=(C_1, \tau_p^h \varphi_1)$ where $g,h \in G$. According to Proposition \ref{prop:WhypOmega}, $x$ and $(C_2, \tau_p^g \varphi_2) \in Q_2$ belong to the same fiber of $J$, hence $t(x)= (C_2, \tau_p^g \varphi_2)$; similarly, $t(y) = (C_2, \tau_p^h \varphi_2)$. Thus,
$$\delta_{Q_2}(t(x),t(y)) = d_G(g,h) = \delta_{Q_1}(x,y).$$
Therefore, our system is also coherent, concluding the proof of our claim. 

\medskip \noindent
Let $\delta$ be the global metric associated to our system of metrics. According to Remark \ref{remark:whichsystem} (where our basepoint is $(\{x_0 \},1)$), Proposition \ref{prop:equicompression} applies to $(\mathfrak{W},\delta)$. However, the canonical orbit map $G \wr H \to (\mathfrak{W},\delta)$, which is associated to the basepoint $(\{ x_0 \}, 1)$, may not be a quasi-isometric embedding, so we need to study its compression. This is done by Proposition \ref{prop:wrcompressionW} below. By combining this result with Proposition \ref{prop:equicompression}, we get the following statement, which is the main result of this section:

\begin{thm}\label{thm:Wequicompression}
Let $G,H$ be two groups. Suppose that $H$ acts on a CAT(0) cube complex $X$ with an orbit map which as compression $\alpha$. Then
$$\alpha_p^*(G \wr H) \geq \alpha \cdot TS(X) \cdot \min \left( \frac{1}{p} , \alpha_p^*(G) \right).$$
\end{thm}

\noindent
The constant $TS(X)$ used in this statement is defined as follows:

\begin{definition}
Let $X$ be a CAT(0) cube complex, $x,y \in X$ two vertices and $F \subset X$ a finite collection of vertices. We denote by $TS(x,F,y)$ the length of the shortest path in $X$ starting from $x$, ending at $y$, and passing through each vertex of $F$; and by $TC(x,F,y)$ twice the number of hyperplanes of the convex hull of $F \cup \{ x,y \}$. The \emph{travelling salesman constant} of $X$, denoted by $TS(X)$, is the supremum of the $\alpha$'s such that there exists some constant $C>0$ for which the inequality
$$TC(x,F,y)+ \# F \geq C \cdot TS(x,F,y)^{\alpha}$$
holds for every vertices $x,y \in X$ and every finite collection of vertices $F \subset X$.
\end{definition}

\begin{remark}
The quantity $TC(x,F,y)$ is linked to a kind of travelling salesman problem in $\mathcal{FS}(X)$. Indeed, $TC(x,F,y)$ turns out to be also the length of the shortest path 
$$C_1= \{ x \}, \ C_2, \ldots, \ C_{n-1}, \ C_n= \{ y \}$$
in the graph of finite subcomplexes $\mathcal{FS}(X)$ such that, for every $z \in F$, there exists some $1 \leq i \leq n-1$ satisfying $z \in C_i$.
\end{remark}

\noindent
The following proposition is the main technical result of this section. In its statement, we endow $G \wr H$ with the generating set which is the union of the two finite generating sets of $G$ and $H$ that we fixed at the beginning of the section. Notice that the word metric associated to this generating set is
$$d_{G \wr H}((g, \varphi),(h,\psi)) = TS_H(g, \mathrm{supp}(\varphi^{-1} \psi),h) + \sum\limits_{k \in \mathrm{supp}(\varphi^{-1} \psi)} d_G(\varphi(k), \psi(k)),$$
where $TS_H(g,F,h)$ denotes the length of the shortest path in the Cayley graph of $H$ starting from $g$, ending at $h$, and passing through each point of $F$. 

\begin{prop}\label{prop:wrcompressionW}
Suppose that there exist a vertex $x_0 \in X$ and some constants $C_1,C_2>0$, and $\alpha \in (0,1]$ such that
$$C_1 \cdot d_H(h_1,h_2)^{\alpha} \leq d_X(h_1x_0,h_2x_0) \leq C_2 \cdot d_H(h_1,h_2)$$
for every $h_1,h_2 \in H$. Suppose that there exist some constants $C_3>0$ and $\beta \in (0,1]$ such that
$$TC(x,F,y) + \# F \geq C_3 \cdot TS(x,F,y)^{\beta}$$
for every $x,y \in X$ and every finite collection of vertices $F \subset X$. Without loss of generality, we assume that $C_1,C_3 \leq 1$ and $C_2 \geq 1/2$. The inequalities
$$\min \left( C_1, (C_1/3)^{\beta}C_3 \right) \cdot \left( d_{G \wr H} \right)^{\alpha \beta} \leq \delta \leq 2C_2 \cdot d_{G \wr H}$$
hold. In particular, the compression of the canonical orbit map $G \wr H \to (\mathfrak{W},\delta)$ has compression at least $\alpha \cdot TS(X)$. 
\end{prop}

\noindent
From now on, we identify $\Omega$ with $H$ (via the map $h \mapsto h \cdot x_0$), so that the support $\mathrm{supp}(\varphi)$ of any map $\varphi \in G^{(\Omega)}$ will be thought of either as a subset of $\Omega$ or as a subset of $H$. We begin by showing an explicit formula for our global metric $\delta$.

\begin{lemma}\label{lem:Wdelta}
For every $(h,\varphi) \in G \wr H$,
$$\delta \left( (\{x_0\}, 1),( \{ hx_0 \},\varphi) \right)= TC(x_0, \mathrm{supp}(\varphi), hx_0) + \sum\limits_{p \in \mathrm{supp}(\varphi)} \mathrm{lg}_G(\varphi(p)),$$
where $\mathrm{lg}_G$ denotes the word length associated to the finite generating set of $G$ we fixed.
\end{lemma}

\begin{proof}
Let $C_1, \ldots, C_n$ be a geodesic in $\mathcal{FS}(X)$ from $\{x_0\}$ to the convex hull $Q$ of $\mathrm{supp}(\varphi) \cup \{x_0,hx_0 \}$, and $D_1, \ldots, D_m$ a geodesic in $\mathcal{FS}(X)$ from $Q$ to $\{ hx_0 \}$. Write $\mathrm{supp}(\varphi)= \{ p_1, \ldots, p_r \}$, and define a sequence of maps $\varphi_0, \ldots, \varphi_r \in G^{(\Omega)}$ by induction in the following way: $\varphi_0=1$ and $\varphi_{i+1} = \tau_{p_{i+1}}^{\varphi(p_{i+1})} \varphi_i$ for every $0 \leq i \leq r-1$. We claim that the path $\gamma$ in $\mathfrak{W}$ defined as
$$(C_1,1) ,\ldots, \left( (C_n,1)=(Q,\varphi_0) \right), \ldots, \left( (Q,\varphi_r)=(D_1, \varphi) \right), \ldots, (D_m,\varphi)$$
turns out to be a geodesic in $\mathfrak{W}$ from $(\{x_0 \},1)$ to $(\{hx_0 \}, \varphi)$. Indeed, we deduce from Fact \ref{fact:dFSalter} and Lemma \ref{lem:wrdist} that
$$\begin{array}{lcl} \mathrm{length}(\gamma) & = & d_{\mathcal{FS}(X)}(\{x_0\},Q)+ d_{\mathcal{FS}(X)}(Q, \{hx_0\}) + \# \mathrm{supp}(\varphi) \\ \\ & = & 2 \cdot \# \mathcal{H}(Q) + \# \mathrm{supp}( \varphi) = d_{\mathfrak{W}} \left( (\{x_0\}, 1),( \{ hx_0 \},\varphi) \right) \end{array}$$
Notice that, for every $1 \leq i \leq n-1$, the wreaths $(C_i,1)$ and $(C_{i+1},1)$ belong to a common clique in which they are at distance one appart with respect to the local metric of this clique; similarly, for every $1 \leq i \leq m-1$, the wreaths $(D_i,\varphi)$ and $(D_{i+1},\varphi)$ belong to a common clique in which they are at distance one appart with respect to the local metric of this clique. (This is because any clique labelled by a halfspace of $X$ is endowed with a discrete metric.) Finally, for every $ 0 \leq i \leq r-1$, the wreaths $(Q,\varphi_i)$ and $(Q, \varphi_{i+1})$ belong to a common clique labelled by $p_{i+1}$ in which they are at distance $d_G(\varphi_i(p_{i+1}), \varphi_{i+1}(p_{i+1})) = \mathrm{lg}_G(\varphi(p_{i+1}))$ appart with respect to the local metric of this clique. Thus, 
$$\delta \left( (\{x_0\}, 1),( \{ hx_0 \},\varphi) \right)= 2 \cdot \# \mathcal{H}(Q) + \sum\limits_{p \in \mathrm{supp}(\varphi)} \mathrm{lg}_G(\varphi(p)),$$
concluding the proof of our lemma.
\end{proof}

\begin{proof}[Proof of Proposition \ref{prop:wrcompressionW}.]
Recall that if for every $g,h \in H$ and every finite set $F \subset H$ we denote by $TS_H(g,F,h)$ the length of the shortest path in the Cayley graph of $H$ starting from $g$, ending at $h$, and passing through each point of $F$, then the word metric associated to our finite generating set of $G\wr H$ is 
$$d_{G \wr H} \left( (g,\varphi), (h,\psi) \right) = TS_H(g, \mathrm{supp}(\varphi^{-1} \psi), h) + \sum\limits_{k \in \mathrm{supp}(\varphi^{-1}\psi)} d_G\left( \varphi(k), \psi(k) \right)$$
for every $(g, \varphi), (h,\psi) \in G \wr H$. To avoid any ambiguity, we also denote $TS$ by $TS_X$. 

\medskip \noindent
Let $(h,\varphi) \in G\wr H$. Suppose that $(h,\varphi) \in G$. This amounts to say that $h=1$ and $\mathrm{supp}(\varphi)=\{1\}$. By applying Lemma \ref{lem:Wdelta}, we find that
$$\begin{array}{lcl} \delta(1,(\{hx_0\}, \varphi)) & = & \displaystyle TC(x_0, \mathrm{supp}(\varphi),hx_0) + \sum\limits_{p \in \mathrm{supp}(\varphi)} \mathrm{lg}_G(\varphi(p)) \\ \\  & = & \mathrm{lg}_G( \varphi(1)) = d_{G \wr H} (1,(h,\varphi)) \end{array}$$
Because $2C_2 \geq 1$ and $(C_1/3)^{\beta}C_3 \leq 1$, it follows that the inequalities claimed in our proposition hold. Next, suppose that $(h,\varphi) \in H$. This amounts to say that $\mathrm{supp}(\varphi)$ is empty. We have
$$\begin{array}{lcl} \delta(1,(\{hx_0\}, \varphi)) & = & \displaystyle TC(x_0, \mathrm{supp}(\varphi),hx_0) + \sum\limits_{p \in \mathrm{supp}(\varphi)} \mathrm{lg}_G(\varphi(p)) \\ \\  & = & 2 \cdot d_X(x_0,hx_0) \end{array}$$
On the other hand, $d_{G\wr H}(1,(h,\varphi))= d_H(1,h)$. Therefore, 
$$ \delta(1,(\{hx_0\}, \varphi)) \leq 2C_2 \cdot d_H(1,h)= 2C_2 \cdot d_{G\wr H}(1,(h,\varphi))$$
and
$$ \delta(1,(\{hx_0\}, \varphi)) \geq C_1 \cdot d_H(1,h)^{\alpha} = C_1 \cdot d_{G\wr H}(1,(h,\varphi))^{\alpha} \geq C_1 \cdot d_{G \wr H}(1,(h, \varphi))^{\alpha \beta}.$$
Thus, the inequalities claimed in our proposition hold.

\medskip \noindent
From now on, suppose that $(h,\varphi)$ belongs neither to $G$ nor to $H$. In particular, this implies that $\mathrm{supp}(\varphi)$ is non empty, so that we are able to index $\mathrm{supp}(\varphi)$ as $\{h_1, \ldots, h_r \}$ so that
$$TS_H(1, \mathrm{supp}(\varphi),h) = d_H(1,h_1)+ \sum\limits_{i=1}^{r-1} d_H(h_i,h_{i+1})+d_H(h_r,h).$$
Notice that, necessarily,
$$TS_X(x_0, \mathrm{supp}(\varphi),hx_0) \leq d_X(x_0,h_1x_0)+ \sum\limits_{i=1}^{r-1} d_X(h_ix_0,h_{i+1}x_0)+d_X(h_rx_0,hx_0).$$
Therefore, thanks to Lemma \ref{lem:Wdelta} and to Claim \ref{claim:wrcW1} proved below, we find that
$$\begin{array}{lcl} \delta(1,(\{hx_0\}, \varphi)) & = & \displaystyle TC(x_0, \mathrm{supp}(\varphi),hx_0) + \sum\limits_{p \in \mathrm{supp}(\varphi)} \mathrm{lg}_G(\varphi(p)) \\ \\ & \leq & \displaystyle 2 \cdot TS_X(x_0, \mathrm{supp}(\varphi),hx_0) + \sum\limits_{p \in \mathrm{supp}(\varphi)} \mathrm{lg}_G(\varphi(p)) \\ \\ & \leq & \displaystyle 2 C_2 \cdot \left( d_H(1,h_1)+ \sum\limits_{i=1}^{r-1} d_H(h_i,h_{i+1}) + d_H(h_i,h) \right) \\ & & \displaystyle + \sum\limits_{p \in \mathrm{supp}(\varphi)} \mathrm{lg}_G(\varphi(p)) \\ \\ & \leq & 2C_2 \cdot d_{G \wr H} (1,(h,\varphi)) \end{array}$$
This proves the second inequality of our proposition. In order to prove the first inequality, we begin with the following observation:

\begin{fact}
For every $x,y \in X$ and every finite collection of vertices $F \subset X$, the inequality 
$$TC(x,F,y)+ \# F \geq 3^{-\beta}C_3 \left( TS(x,F,y)+ \# F  \right)^{\beta}$$
holds unless $x=y$ and $F = \{ x \}$. 
\end{fact}

\noindent
First of all, notice that $TS(x,F,y) \geq \frac{1}{2} \# F$. If $\# F=0$, there is nothing to prove. If $\# F=1$, the inequality also holds because our assumption implies that $TC(x,F,y) \neq 0$. Now, suppose that $\# F \geq 2$. Number the elements of $F$ as $\{f_1, \ldots, f_n\}$ so that
$$TS(x,F,y)= d(x,f_1)+ \sum\limits_{i=1}^{n-1} d(f_i,f_{i+1}) + d(f_n,y).$$
This equality implies that $TS(x,F,y) \geq n-1= \# F -1 \geq \frac{1}{2} \# F$. As a consequence,
$$\left( TS(x,F,y)+ \# F \right)^{\beta} \leq 3^{\beta} TS(x,F,y)^{\beta} \leq \frac{3^{\beta}}{C_3} \left( TC(x,F,y)+ \# F \right),$$
which proves our fact.

\medskip \noindent
Notice that, because $(h,\varphi)$ does not belong to $G$, the equality $\mathrm{supp}(\varphi)= \{1=h \}$ does not hold, so that the previous fact apply. We write
$$\begin{array}{lcl} \delta(1,(\{hx_0\},\varphi)) & = & \displaystyle TC(x_0,\mathrm{supp}(\varphi),hx_0)+  \sum\limits_{p \in \mathrm{supp}(\varphi)} \mathrm{lg}_G(\varphi(p)) \\ \\  & = & \displaystyle TC(x_0,\mathrm{supp}(\varphi),hx_0) + \# \mathrm{supp}(\varphi)+  \sum\limits_{p \in \mathrm{supp}(\varphi)} \left( \mathrm{lg}_G(\varphi(p))-1 \right)  \\ \\ & \geq & \displaystyle 3^{-\beta}C_3 \cdot \left( TS_X(x_0, \mathrm{supp}(\varphi),hx_0) + \# \mathrm{supp}(\varphi) \right)^{\beta} +  \sum\limits_{p \in \mathrm{supp}(\varphi)} \left( \mathrm{lg}_G(\varphi(p))-1 \right) \\ \\ & \geq & \displaystyle 3^{-\beta} C_3 \cdot \left( \left( TS_X(x_0, \mathrm{supp}(\varphi),hx_0) + \# \mathrm{supp}(\varphi) \right)^{\beta} +  \sum\limits_{p \in \mathrm{supp}(\varphi)} \left( \mathrm{lg}_G(\varphi(p)) -1 \right)^{\beta} \right) \\ \\ & \geq &  \displaystyle 3^{-\beta} C_3 \cdot \left( TS_X(x_0, \mathrm{supp}(\varphi),hx_0) + \# \mathrm{supp}(\varphi) + \sum\limits_{p \in \mathrm{supp}(\varphi)} \left( \mathrm{lg}_G(\varphi(p))-1 \right) \right)^{\beta} \\ \\ & \geq & \displaystyle 3^{-\beta} C_3 \cdot \left( TS_X(x_0, \mathrm{supp}(\varphi),hx_0) + \sum\limits_{p \in \mathrm{supp}(\varphi)} \mathrm{lg}_G(\varphi(p)) \right)^{\beta} \\ \\ & \geq & 3^{-\beta}C_3 \cdot \eta(1,(h,\varphi))^{\beta}
\end{array}$$
where $\eta(1,(h,\varphi))$ denotes $TS_X(x_0, \mathrm{supp}(\varphi),hx_0) +  \sum\limits_{p \in \mathrm{supp}(\varphi)} \mathrm{lg}_G(\varphi(p))$. Now, index $\mathrm{supp}(\varphi)$ as $\{h_1, \ldots, h_r\}$ so that
$$TS_X(x_0, \mathrm{supp}(\varphi),hx_0)= d_X(x_0,h_1x_0)+ \sum\limits_{i=1}^{r-1} d_X(h_ix_0,h_{i+1}x_0) + d_X(h_rx_0,hx_0).$$
Notice that, necessarily,
$$TS_H(1, \mathrm{supp}(\varphi),h) \leq d_H(1,h_1)+ \sum\limits_{i=1}^{r-1}d_H(h_i,h_{i+1})+ d(h_r,h).$$
Therefore, 
$$\begin{array}{lcl} TS_X(x_0, \mathrm{supp}(\varphi),hx_0) & \geq & \displaystyle C_1 \cdot \left(  d_H(1,h_1)^{\alpha}+ \sum\limits_{i=1}^{r-1}d_H(h_i,h_{i+1})^{\alpha}+ d(h_r,h)^{\alpha} \right) \\ \\ & \geq &  \displaystyle C_1 \cdot \left(  d_H(1,h_1)+ \sum\limits_{i=1}^{r-1}d_H(h_i,h_{i+1})+ d(h_r,h) \right)^{\alpha} \\ \\ & \geq & \displaystyle C_1 \cdot TS_H(1, \mathrm{supp}(\varphi),h)^{\alpha}
\end{array}$$
We conclude that
$$\begin{array}{lcl} \eta(1,(h,\varphi)) & = & \displaystyle TS_X(x_0, \mathrm{supp}(\varphi), hx_0)+ \sum\limits_{p \in \mathrm{supp}(\varphi)} \mathrm{lg}_G(\varphi(p)) \\ \\ & \geq & \displaystyle C_1 \cdot TS_H(x_0, \mathrm{supp}(\varphi),hx_0)^{\alpha} + \sum\limits_{p \in \mathrm{supp}(\varphi)} \mathrm{lg}_G(\varphi(p)) \\ \\ & \geq & \displaystyle C_1 \cdot \left( TS_H(x_0, \mathrm{supp}(\varphi),hx_0)^{\alpha}+ \sum\limits_{p \in \mathrm{supp}(\varphi)} \mathrm{lg}_G(\varphi(p))^{\alpha} \right) \\ \\ & \geq & \displaystyle C_1 \cdot \left( TS_H(x_0, \mathrm{supp}(\varphi),hx_0)+ \sum\limits_{p \in \mathrm{supp}(\varphi)} \mathrm{lg}_G(\varphi(p)) \right)^{\alpha} \\ \\ & \geq & C_1 \cdot \left( d_{G \wr H}(1,(h,\varphi)) \right)^{\alpha}
\end{array}$$
Thus,
$$\delta(1,(\{hx_0\},\varphi)) \geq 3^{-\beta}C_3 \cdot \eta(1,(h,\varphi))^{\beta} \geq (C_1/3)^{\beta}C_3 \cdot d_{G \wr H}(1,(h,\varphi))^{\alpha \beta},$$
concluding the proof.
\end{proof}

\begin{claim}\label{claim:wrcW1}
For every $x,y \in X$ and every finite collection of vertices $F\subset X$, $$TC(x,F,y) \leq 2 \cdot TS_X(x,F,y).$$ 
\end{claim}

\begin{proof}
Choose a path $\gamma$ in $X$ starting from $x$, ending at $y$, and passing through each vertex of $F$. If $J$ is a hyperplane of the convex hull of $F \cup \{ x , y \}$, then $J$ must separate two vertices of $F \cup \{ x, y \}$. Because $\gamma$ passes through these two vertices, we deduce that $J$ intersects $\gamma$. Therefore, the length of $\gamma$ is at least the number of hyperplanes of the convex hull of $F \cup \{x,y\}$. 
\end{proof}

\noindent
In order to apply Theorem \ref{thm:Wequicompression}, we need to determine the constant $TS(X)$ of the CAT(0) cube complex $X$ on which $H$ acts. This raises natural questions: how to compute $TS(X)$? what are the possible values of $TS(X)$? is there a link between $TS(X)$ and the geometry of $X$? The first thing we can say is that this constant necessarily belongs to the interval $[1/2,1]$. Below, we will see that the values $1/2$ and $1$ can be attained, so this interval is the smallest possible. 

\begin{lemma}\label{lem:TS12}
For any unbounded CAT(0) cube complex $X$, $TS(X)$ belongs to $[1/2,1]$. 
\end{lemma}

\begin{proof}
Let $x,y \in X$ be two vertices and $F \subset X$ a finite collection of vertices. If $F= \{ f_1, \ldots, f_r \}$, then
$$\begin{array}{lcl} TS(x,F,y) & \leq & \displaystyle d(x,f_1)+ \sum\limits_{i=1}^{r-1} d(f_i,f_{i+1}) + d(f_r,y) \leq (r+1) \cdot \mathrm{diam}(F \cup \{x,y \}) \\ \\ & \leq & 2r \cdot \# \mathcal{H}(F \cup \{x,y \}) = \# F \cdot TC(x,F,y) \\ \\ & \leq & \left( TC(x,F,y) + \# F \right)^2 \end{array}$$
Therefore, $TS(X) \geq 1/2$. On the other hand, a path starting from $x$, ending at $y$, and passing through each point of $F$, must intersect any hyperplane separating two vertices of $F \cup \{x,y\}$, hence $TS(x,F,y) \geq TC(x,F,y)/2$. Now, because $X$ is unbounded, we can choose $F= \emptyset$ and $x,y$ at distance $n$ appart for some $n$ arbitrarily large; if so, notice that $TC(x,F,y)=2n$. Thus, if
$$TC(x,F,y) + \# F \geq C \cdot TS(x,F,y)^{\alpha}$$
for some $\alpha, C>0$, then $2n \geq \frac{C}{2^{\alpha}} \cdot n^{\alpha}$. Since this inequality must hold for every $n \geq 1$, we deduce that $\alpha \leq 1$. A fortiori, $TS(X) \leq 1$. 
\end{proof}

\noindent
In view of the lower bound given by Theorem \ref{thm:Wequicompression}, we are interested in understanding when $TS(X)$ is as large as possible, ie., when $TS(X)=1$. This is the purpose of our next proposition. 

\begin{prop}\label{prop:TS1}
Let $X$ be a uniformly locally finite CAT(0) cube complex. The equality $TS(X)=1$ holds if and only if $X$ is hyperbolic. 
\end{prop}

\noindent
The computation of our constant for hyperbolic CAT(0) cube complexes is essentially a consequence of the following lemma. The proof was written jointly with J\'er\'emie Chalopin and Victor Chepo\"{i}. 

\begin{lemma}\label{lem:VleqH}
Let $X$ be a finite $\delta$-hyperbolic CAT(0) cube complex in which each vertex has at most $N$ neighbors. Then
$$\# V(X) \leq N^{2\delta \cdot \left(1+ 2 \cdot \mathrm{Ram}(2\delta+1) \right)} \cdot \# \mathcal{H}(X) +1,$$
where $V(X)$ and $\mathcal{H}(X)$ denote the set of the vertices and hyperplanes of $X$ respectively, and $\mathrm{Ram}(\cdot)$ the Ramsey number.
\end{lemma}

\begin{proof}
First of all, recall some terminology from \cite{coningoff}. A \emph{facing triple} is the data of three hyperplanes such that none separates the other two. A \emph{join of hyperplanes} $(\mathcal{H}, \mathcal{V})$ is the data of two collections of hyperplanes which do not contain facing triples such that any hyperplane of $\mathcal{H}$ is transverse to any hyperplane of $\mathcal{V}$. Such a join is \emph{$D$-thin} if $\min \left( \# \mathcal{H}, \# \mathcal{V} \right) \leq D$. According to the proof of \cite[Theorem 3.3]{coningoff}, it follows from the $\delta$-hyperbolicity of $X$ that any join of hyperplanes of $X$ must be $\mathrm{Ram} \left( \max \left( \delta, \dim(X) \right)+1 \right)$-thin. On the other hand, because an $2n$-cube contains a triangle which is not $(n+1)$-thin, we deduce that $\dim(X) \leq 2 \delta$. Therefore, setting $D= \mathrm{Ram}(2 \delta+1)$, the join of hyperplanes of $X$ are all $D$-thin. In this proof, we will be also interested in the distance $d_{\infty}$ on the vertices of $X$, which is the metric associated to the graph obtained from $X^{(1)}$ by adding an edge between any two vertices which belong to a common cube. Alternatively, the $d_{\infty}$-distance between two vertices $x,y \in X$ is the cardinality of a maximal collection of pairwise disjoint hyperplanes separating $x$ and $y$. 

\medskip \noindent
Now, fix a basepoint $x_0 \in X$. For every hyperplane $J$ of $X$, denote by $P(J)$ the set of the vertices of $N(J)$ minimizing the $d_{\infty}$-distance to $x_0$. We begin by proving two preliminary facts.

\begin{fact}\label{fait1}
For every hyperplane $J$, $P(J)$ is included into a ball of radius $2dD$.
\end{fact}

\noindent
Let $z \in N(J)$ denote the projection of $x_0$ onto $N(J)$. Let $H_1, \ldots, H_s$ be a maximal collection of pairwise disjoint hyperplanes separating $x_0$ and $z$, ie., $s=d_{\infty}(x_0,z)$. Say that $H_i$ separates $H_{i-1}$ and $H_{i+1}$ for every $2 \leq i \leq s-1$, and that $H_1$ separates $H_s$ and $z$. Since any hyperplane separating $x_0$ and $z$ must separate $x_0$ from $N(J)$, we know that $H_1, \ldots, H_s$ separate $x_0$ and $N(J)$.  Let $y \in P(J)$ be an arbitrary vertex and $\mathcal{V}$ a maximal collection of pairwise disjoint hyperplanes separating $z$ and $y$, ie., $\# \mathcal{V}=d_{\infty}(z,y)$. Notice that $H_1, \ldots, H_s$ separate $x_0$ and $y$ because $z$ belongs to the interval between $x_0$ and $y$. Set $\mathcal{V}= \mathcal{V}_1 \sqcup \mathcal{V}_2$ where $\mathcal{V}_1$ is the set of the hyperplanes of $\mathcal{V}$ which are transverse to $H_{D+1}$ and $\mathcal{V}_2$ its complementary in $\mathcal{V}$. By noticing that $\{H_1, \ldots, H_{D+1} \}$  and $\mathcal{V}_1$ define a join of hyperplanes, we deduce that $\# \mathcal{V}_1 \leq D$. On the other hand, since $y$ minimizes the $d_{\infty}$-distance to $x_0$ in $N(J)$, necessarily 
$$s=d_{\infty}(x_0,z) \geq d_{\infty}(x_0,y),$$
and because $\mathcal{V}_2 \sqcup \{ H_{D+1}, \ldots, H_s \}$ is collection of pairwise disjoint hyperplanes separating $x_0$ and $y$, necessarily
$$d_{\infty}(x_0,y) \geq s-D + \# \mathcal{V}_2,$$
hence $\# \mathcal{V}_2 \leq D$. Finally, we conclude that
$$d(y,z) \leq d \cdot d_{\infty}(y,z) = d \cdot \left( \# \mathcal{V}_1+ \# \mathcal{V}_2 \right) \leq 2dD.$$
Thus, we have proved that $P(J)$ is included into the ball of radius $2dD$ centered at $z$.

\begin{fact}\label{fait2}
Let $J$ be a hyperplane separating two adjacent vertices $y,z \in N(J)$ which satisfy $d_{\infty}(x_0,z)= d_{\infty}(x_0,y)+1$. Then $y \in P(J)$. 
\end{fact}

\noindent
Let $J_1, \ldots, J_r$ be a maximal collection of pairwise disjoint hyperplanes separating $x_0$ and $z$, ie., $r= d_{\infty}(x_0,z)$. Since $d_{\infty}(x_0,y) < d_{\infty}(x_0,z)$, necessarily $J_r$ separates $y$ and $z$, ie., $J_r=J$. As a consequence, the hyperplanes $J_1, \ldots, J_{r-1}$ must separate $x_0$ from $N(J)$, so that 
$$d_{\infty}(x_0,w) \geq r-1 = d_{\infty}(x_0,y)$$
for every $w \in N(J)$. Equivalently, $y$ minimizes the $d_{\infty}$-distance to $x_0$ in $N(J)$, ie., $y \in P(J)$.

\medskip \noindent
Now, we are ready to prove our lemma. Set the following collection of edges:
$$\mathcal{E}= \{ (x,y) \in E(X) \mid d_{\infty}(x_0,y)= d_{\infty}(x_0,x)+1 \}.$$
If $(x,y) \in \mathcal{E}$, we refer to $y$ as the top endpoint of $(x,y)$. If $(x,y),(a,b) \in \mathcal{E}$ are two edges dual to the same hyperplane $J$, we deduce from Fact \ref{fait2} that $a,x \in P(J)$. On the other hand, if we denote by $B(k)$ the maximal cardinality of a ball of radius $k$ in $X$, we know that the cardinality of $P(J)$ is at most $B(2dD)$ according to Fact \ref{fait1}, so $J$ intersects at most $B(1)B(2dD)$ edges of $\mathcal{E}$. We deduce that
$$\# \mathcal{E} \leq B(1)B(2dD) \cdot \# \mathcal{H}(X).$$
Now, fix a vertex $x \in X$ which is different from $x_0$. Let $x_0,x_1, \ldots, x_n=x$ be a geodesic between $x_0$ and $x$ with respect to $d_{\infty}$, and let $J_1, \ldots, J_n$ be a maximal collection of pairwise disjoint hyperplanes separating $x_0$ and $x=x_n$. Because $d_{\infty}(x_0,x_{n-1}) < d_{\infty}(x_0,x_n)$, necessarily $J_n$ must separate $x_{n-1}$ and $x_n$. As a consequence, since $x_{n-1}$ and $x_n$ belong to a common cube $C$, the hyperplane $J_n$ must be dual to some edge of $C$. Let $a \in C$ denote the vertex of $C$ which is adjacent to $x_{n-1}$ and which is separated from $x_{n-1}$ by $J_n$. On the one hand,
$$d_{\infty}(x_0,a) \leq d_{\infty}(x_0,x_{n-1}) + d_{\infty}(x_{n-1},a) \leq n-1+1=n;$$
on the other hand, because $J_1, \ldots, J_n$ separate $x_0$ and $a$, necessarily $d_{\infty}(x_0,a) \geq n$. Consequently, $d_{\infty}(x_0,a)=n$, so $(x_{n-1},a) \in \mathcal{E}$. Thus, we have proved that any vertex of $X$ which is different from $x_0$ is at distance at most $d-1$ from the top endpoint of an edge of $\mathcal{E}$. Therefore,
$$\# V(X) \leq 2 B(d-1) \cdot \# \mathcal{E} \leq B(d-1)B(1)B(2dD) \cdot \# \mathcal{H}(X) +1,$$
which concludes the proof because $B(k) \leq N^k$ for every $k \geq 0$. 
\end{proof}

\noindent
In the previous lemma, it is clear the hyperbolicity of the cube complex is necessary. For instance, as $n$ tends to infinity, the number of vertices of the square complex $[0,n]^2$ is quadratic whereas its number of hyperplanes is linear. A uniform control on the number of neighbors of the vertices is also required, as shown by the following example (which was communicated to us by Victor Chepo\"{i}). 

\begin{ex}
Given a triangle-free simplicial graph $\Gamma$, let $X(\Gamma)$ denote the graph whose vertex set is $\{ \emptyset \} \cup V(\Gamma) \cup E(\Gamma)$ and whose edges link an edge of $E(\Gamma)$ to its endpoints (which belong to $V(\Gamma)$) and $\emptyset$ to any vertex of $V(\Gamma)$. This defines a median graph (or, equivalently, a CAT(0) square complex). Notice that $X(\Gamma)$ is $6$-hyperbolic (since it has diameter at most six), and it contains $\# V(\Gamma) + \# E(\Gamma)+1$ vertices and $\# V(\Gamma)$ hyperplanes. However, there exist triangle-free graphs whose number of edges is superlinear in the number of vertices, for instance the bipartite complete graph $K_{n,n}$. It is worth noticing that we have no control on the degree of the vertex $\emptyset$.
\end{ex}

\noindent
In order to prove the implication of the statement of our proposition, we need to show that $TS(\mathbb{R}^2) < 1$, since a non-hyperbolic CAT(0) cube complex must contain an arbitrarily large piece of $\mathbb{R}^2$ as an isometrically embedded subcomplex (see for instance \cite[Corollary 5]{CDEHV} or \cite[Theorem 3.3]{coningoff}). In fact, we are able to prove more.

\begin{lemma}\label{lem:TSeuclid}
For every $d \geq 1$, $\displaystyle TS\left( \mathbb{R}^d \right) = \frac{d}{2d-1}.$
\end{lemma}

\begin{proof}
The equality $TS(\mathbb{R})=1$ is an easy exercice left to the reader. From now on, suppose that $d \geq 2$. First, consider the set
$$F = \left\{ (k_1n^{d-1}, \ldots, k_dn^{d-1}) \mid 0 \leq k_i \leq n \right\} \subset \mathbb{R}^d$$
and fix two arbitrary vertices $x,y \in F$. Notice that $\# F = (n+1)^d$ and that the convex hull of $F \cup \{ x,y\}=F$ is $[0,n^d]^d$; as a consequence, $TC(x,F,y)=2dn^d$. By noticing that any two vertices of $F$ are at distance at least $n^{d-1}$ appart, we deduce that $TS(x,F,y) \geq n^{d-1} \left( \# F -1 \right)$. Therefore,
$$\frac{\ln \left( TS(x,F,y) \right)}{ \ln \left( TC(x,F,y) + \#F \right)} \geq \frac{\ln \left( n^{d-1} \left( (n+1)^d-1 \right) \right)}{\ln \left( 2dn^d + (n+1)^d \right)} \underset{n \to + \infty}{\longrightarrow} \frac{2d-1}{d}.$$
On the other hand, if $\alpha \in (0,1]$ and $C>0$ are such that
$$C \cdot TS(x,F,y)^{\alpha} \leq TC(x,F,y)+ \# F,$$
then 
$$\frac{\ln \left( TS(x,F,y) \right)}{ \ln \left( TC(x,F,y) + \#F \right)} \leq \frac{\ln \left( TS(x,F,y) \right)}{ \ln(C)+ \alpha \ln \left( TS(x,F,y) \right)} \underset{n \to + \infty}{\longrightarrow} \frac{1}{\alpha},$$
so that $\alpha \leq \frac{d}{2d-1}$. Thus, we have proved that $TS(\mathbb{R}^d) \leq \frac{d}{2d-1}$. 

\medskip \noindent
Now, let $x,y \in \mathbb{R}^d$ be two vertices and $F \subset \mathbb{R}^d$ a finite collection of vertices. The convex hull of $F \cup \{x,y\}$ is a rectangular cuboid, which we can write, up to a translation, as $P=\prod\limits_{i=1}^d [0,a_i]$. Set $p= \sum\limits_{i=1}^d a_i$, and, for every $2 \leq i \leq d$, let $c_i$ satisfy $c_i^d \leq a_i < (c_i+1)^d$; notice that $c_i \leq a_i^{1/d} \leq p^{1/d}$. In particular, for every $2 \leq i \leq d$, we can write $a_i= c_i^d + \epsilon_i$ where 
$$\epsilon_i \leq (c_i+1)^d -c_i^d = \sum\limits_{k=0}^{d-1} \left( \begin{matrix} d \\ k \end{matrix} \right) c_i^k \leq (d+1)! \cdot c_i^{d-1} \leq (d+1)! \cdot p^{d-1}.$$ 
Decomposing each $[0,a_i]$ as the union 
$$\left[ 0,c_i^{d-1} \right] \cup \left[ c_i^{d-1},2 c_i^{d-1} \right] \cup \cdots \cup \left[ (c_i-1)c_i^{d-1},c_i^{d} \right] \cup \left[ c_i^{d},c_i^{d}+ \epsilon_i \right]$$
leads to a decomposition of $Q=\prod\limits_{i=2}^d [0,a_i]$ as a union of $\prod\limits_{i=2}^d (c_i+1) \leq (p^{1/d}+1)^{d-1}$ rectangular cuboids of diameter at most $(d+1)! \sum\limits_{i=2}^d c_i^{d-1} \leq (d+2)! \cdot p^{(d-1)/d}$. Index this collection of cuboids as $Q_1, \ldots, Q_r$ in such a way that $0$ belongs to $Q_1$ and that, for every $1 \leq i \leq r-1$, $Q_i$ and $Q_{i+1}$ are adjacent. As a consequence, we get a decomposition of $P$ as the union of the cuboids $P_i = [0,a_1] \times Q_i$, $1 \leq i \leq r$. For every $1 \leq i \leq r$, let $N_i$ denote the number of points of $F$ which lie in $P_i$, and for convenience set $N= \sum\limits_{i=1}^r N_i$, which is the cardinality of $\#F$. 

\medskip \noindent
Next, we construct a path $\gamma = \gamma_0 \cup \gamma_1 \cup \cdots \cup \gamma_r \cup \gamma_{r+1}$ in $P$ in the following way:
\begin{itemize}
	\item $\gamma_0$ is a geodesic from $x$ to $0$.
	\item Fix some $1 \leq i \leq r$, and suppose that $\gamma_{i-1}$ is defined. Index the $N_i$ points of $F \cap P_i$ as $(x_1,y_1), \ldots, (x_{N_i}, y_{N_i})$ in such a way that the sequence $(x_i)$ is non decreasing, and let $\gamma_i$ be a concatenation of geodesics starting from the endpoint of $\gamma_{i-1}$, passing through $(x_1,y_1)$, next through $(x_2,y_2)$, and so on until $(x_{N_i}, y_{N_i})$, and finally following a geodesic to $P_i \cap P_{i+1} \cap \{0 \} \times Q$ (if $i=r$, do not add this final segment). 
	\item $\gamma_{r+1}$ is a geodesic from the endpoint of $\gamma_r$ to $y$. 
\end{itemize}
The lengths of $\gamma_0$ and $\gamma_{r+1}$ are at most $\mathrm{diam}(P)=p$, and the length of each $\gamma_i$ is at most
$$\sum\limits_{k=1}^{N_i} (x_{k+1}-x_k) + \sum\limits_{k=1}^{N_i} |y_{k+1}-y_k| +2a_1 \leq 3a_1 + N_i \cdot (d+2)! \cdot p^{(d-1)/d}$$
Therefore, the length of $\gamma$ is at most
$$\begin{array}{lcl} 2p +3ra_1 + N \cdot (d+2)! \cdot p^{(d-1)/d} & \leq & 2p+ 3p \cdot (p^{1/d}+1)^{d-1} + N \cdot (d+2)! \cdot p^{(d-1)/d} \\ \\ & \leq & 2p + d! \cdot p^{(2d-1)/d} + N \cdot (d+2)! \cdot p^{(d-1)/d}  \end{array}$$
Because $\gamma$ starts from $x$, ends at $y$, and passes through each point of $F$, we deduce that 
$$TS(x,F,y) \leq (d+1)! \cdot  p^{(2d-1)/d} + N \cdot (d+2)! \cdot p^{(d-1)/d}.$$
On the other hand, $TC(x,F,y)=2p$, so we conclude thanks to Claim \ref{claim:TSRd} below that there exists some constant $C>0$ such that
$$TS(x,F,y)^d \leq C \cdot (TC(x,F,y)+ \# F)^{2d-1}.$$
Thus, we have proved that $TS(\mathbb{R}^d) \geq \frac{d}{2d-1}$.

\begin{claim}\label{claim:TSRd}
For every $d \geq 1$ and $A,B,p,N \geq 0$, there exists some $C>0$ such that
$$\left(Ap^{(2d-1)/d} +B \cdot N \cdot p^{(d-1)/d} \right)^d \leq C \cdot \left( 2p+ N \right)^{2d-1}.$$
\end{claim}

\noindent
First, notice that
$$\left(Ap^{(2d-1)/d} +B \cdot N \cdot p^{(d-1)/d} \right)^d = \sum\limits_{k=0}^d \left( \begin{matrix} d \\  k \end{matrix} \right) B^kN^k A^{d-k} p^{2d-1-k}.$$
If $C= (AB)^d \cdot d!$, a fortiori
$$\left( \begin{matrix} d \\ k \end{matrix} \right) B^k A^{d-k} \leq C \leq C \cdot 2^{2d-1-k} \left( \begin{matrix} 2d-1 \\ k \end{matrix} \right)$$
for every $0 \leq k \leq d$, hence
$$\begin{array}{lcl} \displaystyle \left(Ap^{(2d-1)/d} +B \cdot N \cdot p^{(d-1)/d} \right)^d & \leq & \displaystyle C \cdot \sum\limits_{k=0}^d \left( \begin{matrix} 2d-1 \\ k \end{matrix} \right) N^k (2p)^{2d-1-k}\\ \\ & \leq & \displaystyle C \cdot \sum\limits_{k=0}^{2d-1}\left( \begin{matrix} 2d-1 \\ k \end{matrix} \right) N^k (2p)^{2d-1-k} \\ \\ & \leq & \displaystyle C \cdot \left(  2p + N \right)^{2d-1} \end{array}$$
This proves our claim. 
\end{proof}

\begin{cor}\label{cor:TSnonhyp}
For every non hyperbolic CAT(0) cube complex $X$, the inequality $TS(X) \leq 2/3$ holds. 
\end{cor}

\begin{proof}
If $X$ is a non hyperbolic CAT(0) cube complex, it follows from \cite[Theorem 3.3]{coningoff} that, for every $n \geq 1$, $X$ contains an isometrically embedded subcomplex isomorphic to $[0,n]^2$. Therefore, $TS(X) \leq TS(\mathbb{R}^2) = 2/3$. 
\end{proof}

\begin{cor}
For any infinite-dimensional CAT(0) cube complex $X$, the equality $TS(X)=1/2$ holds. 
\end{cor}

\begin{proof}
It is well-know that any finite CAT(0) cube complex isometrically embeds into the cube $[0,1]^n$ for some sufficiently large $n \geq 1$; for instance, this is a consequence of Lemma \ref{lem:embedinfiniteprism}. Therefore, if $X$ is an infinite-dimensional CAT(0) cube complex, then, for every $n \geq1$, $X$ contains an isometrically embedded subcomplex isomorphic to $[0,n]^n \subset [0,1]^{n^2}$. Thus, $TS(X) \leq TS(\mathbb{R}^d)$ for every $d \geq 1$, hence $TS(X) \leq 1/2$. On the other hand, we know from Lemma \ref{lem:TS12} that $TS(X) \geq 1/2$, hence $TS(X)=1/2$ finally. 
\end{proof}

\begin{proof}[Proof of Proposition \ref{prop:TS1}.]
If $X$ is not hyperbolic, then $TS(X) \neq 1$ according to Corollary \ref{cor:TSnonhyp}. Conversely, if $X$ is hyperbolic, it follows from Lemma \ref{lem:VleqH} that there exists some constant $C>0$ such that $\# V(Y) \leq C \cdot \# \mathcal{H}(Y)$ for every convex subcomplex $Y \subset X$. On the other hand, if $x,y \in X$ are two vertices and $F \subset X$ a finite collection of vertices, then $TS(x,F,y)$ is at most twice the number of edges of the convex hull $Q$ of $F \cup \{x,y \}$ (indeed, any graph contains a path visiting all its vertices and at most twice each of its edges). Therefore, 
$$TS(x,F,y) \leq 2N \cdot \# V(Q) \leq 2CN \cdot \# \mathcal{H}(F \cup \{x,y \}) \leq 2CN \cdot \left( \# \mathcal{H}(F \cup \{x,y \}) + \# F \right),$$
where $N$ denotes the maximal number of neighbors of the vertices of $X$. Thus, we have proved that $TS(X)=1$. 
\end{proof}

\noindent
By combining Theorem \ref{thm:Wequicompression} with Proposition \ref{prop:TS1}, we obtain the following statement, which is our main application.

\begin{thm}\label{thm:Wcompressionhyp}
Let $H$ be a hyperbolic group acting geometrically on some CAT(0) cube complex. For every finitely generated group $G$ and every $p \geq 1$, 
$$\alpha_p^*(G \wr H) \geq \min \left( \frac{1}{p}, \alpha_p^*(G) \right),$$
with equality if $H$ is non elementary and $p \in [1,2]$.
\end{thm}

\begin{proof}
Let $H$ act geometrically on a CAT(0) cube complex $X$. According to Lemma \ref{lem:modifycubing}, we can suppose without loss of generality that $X$ contains a vertex $x_0$ of trivial stabiliser. Fixing some word metric on $H$ (associated to a finite generating set), it follows from Milnor-\v{S}varc lemma that there exist constants $A,B>0$ such that
$$\frac{1}{A} \cdot d(g,h) - B \leq d(g \cdot x_0, h \cdot x_0) \leq A \cdot d(g,h) +B$$
for every $g,h \in H$. For every distinct $g,h \in H$, notice that
$$d(g \cdot x_0,h \cdot x_0) \leq A \cdot d(g,h) +B \leq (A+B) \cdot d(g,h),$$
and that
$$d(g,h) \leq A \cdot d(g \cdot x_0,h \cdot x_0) +AB \leq A(1+B) \cdot d(g \cdot x_0, h \cdot x_0)$$
so that
$$\frac{1}{A(1+B)} \cdot d(g,h) \leq d(g \cdot x_0,h \cdot x_0).$$
Thus, we have proved that the orbit map $x_0 \mapsto g \cdot x_0$ is Lipschitz and has compression one. Therefore, since $TS(X)=1$ according to Proposition \ref{prop:TS1}, it follows from Theorem \ref{thm:Wequicompression} that
$$\alpha_p^*(G \wr H) \geq \min \left( \frac{1}{p}, \alpha_p^*(G) \right).$$
If $H$ is non elementary, we know that $H$ contains a quasi-isometrically embedded free non abelian subgroup (see for instance \cite{GoulnaraQuasiconvex}), hence $\alpha_p^*(G \wr H) \leq \alpha_p^*(\mathbb{F}_2) = \max \left( \frac{1}{2}, \frac{1}{p} \right)$. We also know that $\alpha_p^*(G \wr H ) \leq \alpha_p^*(G)$, so, if $p \in [1,2]$, our previous inequality turns out to be an equality.
\end{proof}

\begin{remark}
An interesting consequence of the fact that the constant $TS(X)$ of a CAT(0) cube complex $X$ is bounded below by a positive constant is that it is possible to deduce non zero lower bounds for the $\ell^p$-compressions of iterated wreath products. For instance, given a finitely generated group $G$ acting geometrically on some CAT(0) cube complex $X$, set $L_0(G)= G$ and $L_{k+1}(G)= L_k(G) \wr G$ for every $k \geq 0$; and set $R_0(G)= G$ and $R_{k+1}(G)= G \wr R_k$ for every $k \geq 0$. By applying Theorem \ref{thm:Wequicompression}, we find that 
$$\alpha_p^*(L_k(G)) \geq TS(X)^k \cdot \min \left( \frac{1}{p}, \alpha_p(G) \right)$$ 
and that  
$$\alpha_p^*(R_k(G)) \geq \frac{TS(X)}{2^k} \cdot \min \left( \frac{1}{p}, \alpha_p^*(G) \right)$$ 
for every $k,p \geq 1$. However, if $G= \mathbb{Z}$, these lower bounds are not optimal, since it was proved in \cite{NaorPeres} that $\alpha_2^*(L_k(\mathbb{Z}))= \left( 2-2^{1-k} \right)^{-1}$ for every $k \geq 1$; and in \cite{LiWreathProducts} that $\alpha_2^*(R_k(\mathbb{Z})) \geq 1/2k$ for every $k \geq 1$. (Notice however that we recover the lower bound in the latter case for $k=2$.) On the other hand, if $G$ is any cubulable non elementary hyperbolic group, for instance $G= \mathbb{F}_n$ where $n \geq 2$, we deduce from Theorem \ref{thm:Wcompressionhyp} the equality $\alpha_p^*(L_k(G))= \alpha_p^*(G)$ for every $k \geq 1$ and $p \in [1,2]$, so in this situation we know that the previous lower bound is optimal. But we do not know whether our inequality 
$$\alpha^*_p(R_k(\mathbb{F}_n)) \geq \frac{1}{p \cdot 2^k}, \ n \geq 2, \ k,p \geq 1$$
is sharp or not. 
\end{remark}

\section{Application to diagram products}\label{section:appli2}

\subsection{Diagram products}

\noindent
Diagram products were introduced by Guba and Sapir in \cite{MR1725439} as fundamental groups of 2-complexes of groups (see Corollary \ref{cor:2complex}). We use a different point of view in this section, defining diagram products in a similar way as diagram groups, in terms of semigroup diagrams. Our definitions are slightly different from the definitions used in \cite{MR1396957}. Nevertheless, most of the arguments of the first sections of \cite{MR1396957} still hold in our context without any modification, allowing us to refer to them in the sequel.

\medskip \noindent
Let $\mathcal{P}= \langle \Sigma \mid \mathcal{R} \rangle$ be a semigroup presentation. We suppose that, for any relation $u=v$ of $\mathcal{R}$, the relation $v =u$ does not belong to $\mathcal{R}$; in particular, $\mathcal{R}$ contains no relation of the form $u=u$. Let $\mathcal{G}= \{ G_s \mid s \in \Sigma \}$ be a collection of groups, indexed by our alphabet $\Sigma$. For convenience, set $\Sigma (\mathcal{G}) = \bigcup\limits_{s \in \Sigma } \{ s \} \times G_s$.

\medskip \noindent
Given $\mathcal{P}$ and $\mathcal{G}$, a \emph{semigroup diagram} (or a \emph{diagram}, for short) is a finite connected planar graph $\Delta$ whose edges are oriented and labelled by the alphabet $\Sigma(\mathcal{G})$, satisfying the following properties:
\begin{itemize}
	\item $\Delta$ has exactly one vertex-source $\iota$ (which has no incoming edges) and exactly one vertex-sink $\tau$ (which has no outgoing edges);
	\item the boundary of each cell has the form $pq^{-1}$, where $p,q$ are two positive paths respectively labelled by $(a_1,g_1) \cdots (a_n,g_n)$ and $(b_1,h_1) \cdots (b_m,h_m)$, such that the relation $a_1\cdots a_n = b_1 \cdots b_m$ belongs to $\mathcal{R}$;
	\item every vertex belongs to a positive path connecting $\iota$ and $\tau$;
	\item every positive path in $\Delta$ is simple.
\end{itemize}
In particular, $\Delta$ is bounded by two positive paths: the \emph{top path}, denoted by $\mathrm{top}(\Delta)$, and the \emph{bottom path}, denoted by $\mathrm{bot}(\Delta)$. By extension, we also define $\mathrm{top}(\Gamma)$ and $\mathrm{bot}(\Gamma)$ for every subdiagram $\Gamma$ of $\Delta$. In the following, the notations $\mathrm{top}(\cdot)$ and $\mathrm{bot}(\cdot)$ will refer to either the paths or their labels. Also, a \emph{$(u,v)$-cell} (resp. a \emph{$(u,v)$-diagram}) will refer to a cell (resp. a semigroup diagram) whose top path is labelled by $u$ and whose bottom path is labelled by $v$. Finally, we define $\mathrm{top}^-(\cdot)$ and $\mathrm{bot}^-(\cdot)$ as the images of $\mathrm{top}(\cdot)$ and $\mathrm{bot}(\cdot)$, respectively, under the canonical projection $\Sigma(\mathcal{G})^+ \twoheadrightarrow \Sigma^+$ (which just ``forgets'' the second letters in the alphabet $\Sigma(\mathcal{G})$). A semigroup diagram $\Delta$ satisfying $\mathrm{top}^-(\Delta)= \mathrm{bot}^-(\Delta)$ will be said \emph{spherical}. 

\begin{figure}
\begin{center}
\includegraphics[scale=0.7]{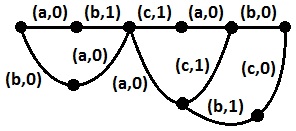}
\end{center}
\caption{Example of a semigroup diagram.}
\label{figure13}
\end{figure}

\medskip \noindent
Figure \ref{figure13} gives an example of diagram $\Delta$ with three cells, constructed from the semigroup presentation
$$\mathcal{P}= \langle a,b,c \mid ab=ba,ac=ca,bc=cb \rangle$$
and the groups $G_a=G_b=G_c= \mathbb{Z}_2$. The edges are supposed oriented from left to right, so that the top path of our diagram is $(a,0)(b,1)(c,1)(a,0)(b,0)$, and its bottom path is $(b,0)(a,0)(a,0)(b,1)(c,0)$. In particular, $\mathrm{top}^-(\Delta)=abcab$ and $\mathrm{bot}^-(\Delta)=baabc$. 

\begin{figure}
\begin{center}
\includegraphics[scale=0.6]{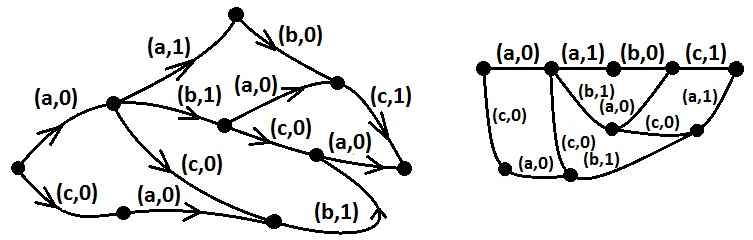}
\end{center}
\caption{Two isotopic semigroup diagrams.}
\label{figure12}
\end{figure}

\medskip \noindent
Since we are only interested in the combinatorics of semigroup diagrams, we will not distinguish isotopic diagrams. For example, the two diagrams of Figure \ref{figure12} will be considered as equal, where we use the semigroup presentation
$$\mathcal{P}= \langle a,b,c \mid ab=ba,ac=ca,bc=cb \rangle$$ 
and $G_a=G_b=G_c= \mathbb{Z}_2$.

\medskip \noindent
Let $\Delta_1, \Delta_2$ be two semigroup diagrams satisfying $\mathrm{bot}^-(\Delta_1)= \mathrm{top}^-(\Delta_2)$, say $\mathrm{bot}(\Delta_1)= (\ell_1,g_1) \cdots (\ell_n,g_n)$ and $\mathrm{top}(\Delta_2)=(\ell_1,h_1) \cdots (\ell_n,h_n)$. The \emph{concatenation} $\Delta_1 \circ \Delta_2$ is defined as the semigroup diagram obtained by gluing $\Delta_1$ and $\Delta_2$ along the paths $\mathrm{bot}(\Delta_1)$ and $\mathrm{top}(\Delta_2)$, which we label by $(\ell_1,g_1h_1) \cdots (\ell_n,g_nh_n)$. See Example \ref{ex:concatenation}.

\medskip \noindent
Suppose that a diagram $\Delta$ contains two cells $\pi_1$ and $\pi_2$ intersecting along $\mathrm{bot}(\pi_1)= \mathrm{top}(\pi_2)$ such that $\mathrm{top}^-(\pi_1)= \mathrm{bot}^-(\pi_2)$ and such that the second coordinates of the letters of the word $\mathrm{bot}(\pi_1)= \mathrm{top}(\pi_2)$ are identities. We say that these two cells form a \emph{dipole}, which we can \emph{reduce} by first removing $\mathrm{bot}(\pi_1)= \mathrm{top}(\pi_2)$ and then identifying the paths $\mathrm{top}(\pi_1)=(\ell_1,g_1) \cdots (\ell_n,g_n)$ and $\mathrm{top}(\pi_2)= (\ell_1,h_1) \cdots (\ell_n,h_n)$ to a single path labelled by $(\ell_1,g_1h_1) \cdots (\ell_n,g_nh_n)$. A diagram is called \emph{reduced} if it does not contain dipoles. The reduced form is unique; see \cite[Theorem 3.17]{MR1396957}. If $\Delta_1, \Delta_2$ are two diagrams for which the concatenation $\Delta_1 \circ \Delta_2$ is well-defined, let us denote by $\Delta_1 \cdot \Delta_2$ the reduced form of $\Delta_1 \circ \Delta_2$. See Example \ref{ex:concatenation}.

\begin{figure}
\begin{center}
\includegraphics[scale=0.55]{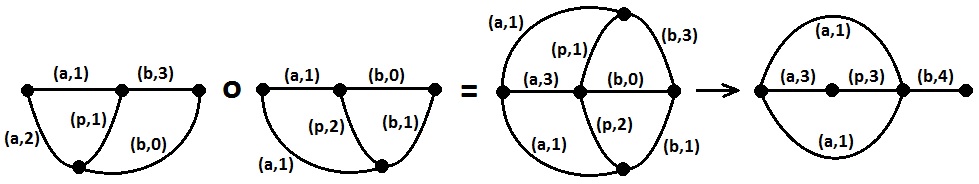}
\end{center}
\caption{A concatenation followed by a reduction.}
\label{figure15}
\end{figure}

\begin{ex}\label{ex:concatenation}
Let $\mathcal{P}= \langle a,b,p \mid a=ap, b=pb \rangle$ and $G_a=G_b=G_p= \mathbb{Z}$. Figure \ref{figure15} gives an example of the concatenation of two semigroup diagrams, followed by the reduction of a dipole. Notice that the produced diagram is reduced. 
\end{ex}

\noindent
If $w=(\ell_1,g_1) \cdots (\ell_n,g_n) \in \Sigma(\mathcal{G})^+$, we denote by $\epsilon(w)$ the diagram consisting in a segment of length $n$ labelled by $w$. Such diagrams are called \emph{linear}. If $u= \ell_1 \cdots \ell_n \in \Sigma^+$, we denote by $\epsilon(u)$ the diagram diagram $\epsilon ( \overline{u})$, where $\overline{u}= (\ell_1,1) \cdots (\ell_n,1)$. Such diagrams are called \emph{trivial}. A diagram is \emph{atomic} if it has a single cell and the letters labelling its edges have the form $(u,1)$ where $u \in \Sigma$. Linear and atomic diagrams will be referred to as \emph{elementary diagrams}. The following lemma is clear. 

\begin{lemma}\label{lem:elementarydiag}
Any semigroup diagram is a concatenation of elementary diagrams.
\end{lemma}

\noindent
In fact, this can be taken as the definition of semigroup diagrams. See the discussion following \cite[Lemma 3.11]{MR1396957}. As a consequence, if we think of the set of (reduced) diagrams as a groupoid following our next proposition, the element diagrams provide a generating set.

\begin{prop}
The set of reduced diagrams $G(\mathcal{P}, \mathcal{G})$, endowed with the product $\cdot$ defined above, is a groupoid.
\end{prop}

\begin{proof}
The product $\circ$ is associative in the sense that, if the concatenations $\Delta_1 \circ \Delta_2$ and $\Delta_2 \circ \Delta_3$ are well-defined, then $(\Delta_1 \circ \Delta_2) \circ \Delta_3= \Delta_1 \circ (\Delta_2 \circ \Delta_3)$. Moreover, the concatenation is compatible with the reduction, in the sense that, if $\Delta_1'$ and $\Delta_2'$ are obtained respectively from $\Delta_1$ and $\Delta_2$ by reducing some dipoles and if the concatenation $\Delta_1' \circ \Delta_2'$ is well-defined, then $\Delta_1 \circ \Delta_2$ is well-defined and $\Delta_1' \circ \Delta_2'$ is obtained from $\Delta_1 \circ \Delta_2$ by reducing some dipoles; see \cite[Lemma 5.1]{MR1396957}. We deduce that the product $\cdot$ is associative. 

\medskip \noindent
It is clear that trivial diagrams are identities. First, we define the inverse of an elementary diagram. Given $u=(\ell_1,g_1) \cdots (\ell_n,g_n) \in \Sigma(\mathcal{G})^+$, set $u'=(\ell_1,g_1^{-1}) \cdots (\ell_n,g_n^{-1})$ and $w=\ell_1 \cdots \ell_n \in \Sigma^+$. Then
$$\epsilon(u) \cdot \epsilon(u') = \epsilon(w) = \epsilon(u') \cdot \epsilon(u),$$
hence $\epsilon(u)^{-1}=\epsilon(u')$. Next, if $A$ is an atomic diagram, let $A^{-1}$ denote the semigroup diagram obtained from taking the mirror image of $A$ with respect to the top path $\mathrm{bot}(A)$, so that the concatenations $A \circ A^{-1}$ and $A^{-1} \circ A$ are well-defined and produce dipoles (see Example \ref{ex:inverse}). We have
$$A \cdot A^{-1}= \epsilon(\mathrm{top}^-(A)) \ \text{and} \ A^{-1} \cdot A = \epsilon( \mathrm{bot}^-(A)).$$
Therefore, elementary diagrams have inverses. Now, let $\Delta$ be an arbitrary semigroup diagram. According to Lemma \ref{lem:elementarydiag}, there exist elementary diagrams $E_1, \ldots, E_n$ such that $\Delta=E_1 \circ \cdots \circ E_n$. Setting $\Delta^{-1}=E_n^{-1} \circ \cdots \circ E_1^{-1}$, we have
$$\Delta \cdot \Delta^{-1}= \epsilon(\mathrm{top}^-(\Delta)) \ \text{and} \ \Delta^{-1} \cdot \Delta = \epsilon( \mathrm{bot}^- (\Delta)).$$
Thus, any semigroup diagram admits left and right inverses. This concludes the proof of our proposition.
\end{proof}

\begin{figure}
\begin{center}
\includegraphics[scale=0.9]{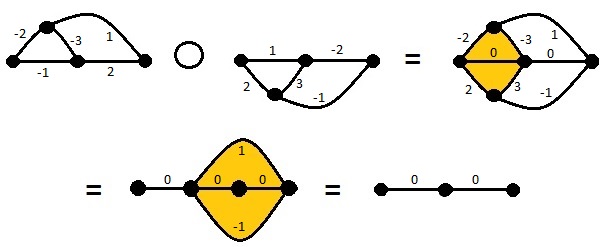}
\end{center}
\label{figure19}
\caption{Reduction of $\Delta^{-1} \circ \Delta$ to $\epsilon(x^2)$.}
\end{figure}

\begin{ex}\label{ex:inverse}
Consider the presentation $\mathcal{P}= \langle x \mid x=x^2 \rangle$, the group $G_x= \mathbb{Z}$ and the base word $x^2$. Figure \ref{figure19} shows how the product $\Delta^{-1} \circ \Delta$, for some diagram $\Delta$, reduces to the trivial diagram $\epsilon(x^2)$ in $D(\mathcal{P}, \mathcal{G},x^2)$. Notice that all the edges of $\Delta$ are necessarily labelled by $x$, so we only indicated the elements of $G_x$. 
\end{ex}

\begin{definition}
Let $w \in \Sigma^+$ be a base word. The \emph{diagram product}\index{Diagram products} $D(\mathcal{P}, \mathcal{G},w)$ is the vertex-group of the groupoid $G(\mathcal{P}, \mathcal{G})$ based at $w$. Equivalently, $D(\mathcal{P},\mathcal{G},w)$ is the set of diagrams $\Delta$ satisfying $\mathrm{top}^-(\Delta)= \mathrm{bot}^-(\Delta)=w$, endowed with the product $\cdot$ defined above.
\end{definition}

\noindent
Notice that if there exists some diagram $\Delta$ such that $\mathrm{top}^-(\Delta)=w$ and $\mathrm{bot}^-(\Delta)=u_1 \cdots u_r av_1 \cdots v_s$ for some letters $u_1, \ldots, u_r,v_1, \ldots, v_s, a \in \Sigma$ then 
$$g \mapsto \Delta \cdot \epsilon((u_1,1) \cdots (u_r,1) (a,g) (v_1,1) \cdots (v_s,1)) \cdot \Delta^{-1}$$
defines an embedding $G_a \hookrightarrow D(\mathcal{P}, \mathcal{G})$ (which is not canonical). Otherwise, if there does not exist such a diagram, then we can remove the letter $a$ from the alphabet $\Sigma$, the relations containing $a$ from $\mathcal{R}$ and the groups $G_a$ from $\mathcal{G}$ without modifying the diagram product. Therefore, up to a natural simplification of $\mathcal{P}$ and $\mathcal{G}$, we can always suppose without loss of generality that the groups of $\mathcal{G}$ embed into $D(\mathcal{P}, \mathcal{G})$. 

\medskip \noindent
Now, we introduce a 2-complex whose fundamental groupoid is isomorphic to $G(\mathcal{P},\mathcal{G})$. In particular, computing its fundamental group will allow us to determine a presentation of $D(\mathcal{P}, \mathcal{G},w)$.

\begin{definition}
The Squier complex $S(\mathcal{P},\mathcal{G})$ is the 2-complex whose set of vertices is $\Sigma^+$, with
\begin{itemize}
	\item an edge $(a, u \to v,b)$ between the words $aub$ and $avb$ whenever $u=v \in \mathcal{R}$; 
	\item a loop $b(g)$ based at $w=\ell_1 \cdots \ell_n \in \Sigma^+$ for every $g \in G(w)=G_{\ell_1} \times \cdots \times G_{\ell_n}$;
	\item a square delimited by the four vertices $aubpc$, $avbpc$, $aubqc$ and $avbqc$ whenever $u=v,p=q \in \mathcal{R}$;
	\item a triangle delimited by the three loops based at $w$ labelled by $g,h,gh \in G(w)$;
	\item a cylinder between the loops labelled by $g \in G(a) \times G(b) \leq G(aub)$ and $h \in G(a) \times G(b) \leq G(avb)$ following the edge $(a,u \to v,b)$, where $u=v \in \mathcal{R}$.
\end{itemize}
\end{definition}

\begin{remark}
If all the groups of $\mathcal{G}$ are trivial, then $S(\mathcal{P},\mathcal{G})$ coincides with the Squier complexes associated to diagram groups in \cite{MR1396957}.
\end{remark}

\begin{thm}\label{thm:fundamentalgroupoid}
$G(\mathcal{P}, \mathcal{G})$ is naturally isomorphic to the fundamental groupoid of $S(\mathcal{P},\mathcal{G})$.
\end{thm}

\begin{proof}
For every edge $(a,u \to v,b)$ of $S(\mathcal{P}, \mathcal{G})$, define $\Phi((a,u \to v,b))$ as the unique atomic diagram satisfying $\mathrm{top}^-(\Delta)= aub$ and $\mathrm{bot}^-(\Delta)=avb$, whose cell corresponds to the relation $u \to v$; and for every loop $b(g)$, where $g=(g_1, \ldots, g_n) \in G(\ell_1 \cdots \ell_n)$, define $\Phi(b(g))$ as the linear diagram $\epsilon(w)$ where $w=(\ell_1,g_1) \cdots (\ell_n,g_n)$. Notice that, for every edge $e$ of $S(\mathcal{P}, \mathcal{G})$, $\mathrm{top}^- (\Phi(e))$ is the initial endpoint of $e$ and $\mathrm{bot}^- (\Phi(e))$ its terminal endpoint, so that it makes sense to extend $\Phi$ from the combinatorial paths of $S(\mathcal{P},\mathcal{G})$ to $G(\mathcal{P}, \mathcal{G})$ by $e_1 \cdots e_n \mapsto \Phi(e_1) \cdots \Phi(e_n)$. Next, we know that the fundamental groupoid $\gamma_1 S(\mathcal{P},\mathcal{G})$ of $S(\mathcal{P}, \mathcal{G})$ is characterized by the relations
\begin{itemize}
	\item $(a,u \to v,bpc)(avb,p \to q,c)=(aub,p \to q,c)(a, u \to v, bqc)$ for every words $a,b,c,p,q \in \Sigma^+$ such that $u=v,p=q \in \mathcal{R}$; 
	\item $b(g_1)b(g_2) = b(g_1g_2)$ for every word $w \in \Sigma^+$ and elements $g_1,g_2 \in G(w)$; 
	\item $(a,u \to v,b) b(g_1) = b(g_2)(a,u \to v,b)$ for every words $a,b,u,v \in \Sigma^+$, where $u=v \in \mathcal{R}$, and elements $g_1 \in G(avb)$, $g_2 \in G(aub)$ representing the same element of $G_a \times G_b$.
\end{itemize}
It is not difficult to verify that the images of these relations by $\Phi$ are relations in $G(\mathcal{P}, \mathcal{G})$. Therefore, $\Phi$ defines a morphism from the fundamental groupoid $\gamma_1 S(\mathcal{P}, \mathcal{G})$ to $G(\mathcal{P}, \mathcal{G})$.

\medskip \noindent
We claim that $\Phi$ turns out to be an isomorphism. Notice that if a diagram $\Delta \in G(\mathcal{P}, \mathcal{G})$ decomposes as a concatenation $E_1 \circ \cdots \circ E_n$ of elementary diagrams, then $\Phi(p)=\Delta$, where $p$ is the path defined as
$$\mathrm{top}^-(E_1), \ \mathrm{bot}^-(E_1)=\mathrm{top}^-(E_2), \ldots, \mathrm{bot}^-(E_{n-1})= \mathrm{top}^-(E_n), \ \mathrm{bot}^-(E_n).$$
Therefore, $\Phi$ is surjective. Finally, let $p$ be a path such that the diagram $\Phi(p)$ contains a dipole. Say that $p$ is a concatenation $e_1\cdots e_r \cdots e_{r+s} \cdots e_n$ of edges such that $\Phi(e_r)$ and $\Phi(e_{r+s})$ define respectively the two cells $\pi_1$ and $\pi_2$ of our dipole. By definition, $e_r= (a,u \to v,b)$ for some $a,b,u,v \in \Sigma^+$ and $u=v \in \mathcal{R}$. The next edge $e_{r+1}$ may be of two types. Either $e_{r+1}$ is a loop $b(g)$ for some $g \in G(avb)$, notice that, since the second coordinates of the letters of $\mathrm{bot}(\pi_1)$ are identifies, necessarily $g \in G(a) \times G(b) \leq G(avb)$; if so, let $e_{r+1}'$ denote the loop $b(g')$ based at $aub$ where $g' \in G(aub)$ represents the same element as $g$ in $G(a) \times G(b)$. Otherwise, the terminal vertex of $e_{r+1}$ is obtained from $avb$ by applying a relation of $\mathcal{R}$ to one of its subword. Since $\pi_2$ is the unique cell of $\Phi(p)$ intersecting an edge of $\mathrm{bot}(\pi_1)$, necessarily this subword must be disjoint from $v$, say $e_{r+1}=(a',x \to y, a''vb)$ for some $a',a'',x,y \in \Sigma^+$ such that $x \to y \in \mathcal{R}$ and $a'xa''=a$ in $\Sigma^+$. If so, let $e_{r+1}'$ denote the edge $(a', x \to y,a''ub)$. Notice that we get a square $(a',x \to y, a'', u \to v,b)$. Similarly, we define $e_{r+i}'$ from $e_{r+i}$ for every $1 \leq i \leq s-1$. Now, replacing the subsegment $e_r \cdots e_{r+s}$ of $p$ with $e_{r+1}' \cdots e_{r+s-1}'$ allows us to shorten $p$ without modifying its homotopy class. Thus, we have proved that the image by $\Phi$ of a path minimizing the length in its homotopy class must be a reduced diagram. As a consequence, any cycle in $S(\mathcal{P}, \mathcal{G})$ with trivial image in $G(\mathcal{P}, \mathcal{G})$ must be homotopically trivial. Necessarily, $\Phi$ is injective.
\end{proof}

\noindent
If $w \in \Sigma^+$, we denote by $S(\mathcal{P}, \mathcal{G},w)$ the connected component of $S(\mathcal{P},\mathcal{G})$ containing $w$. As an immediate consequence of Theorem \ref{thm:fundamentalgroupoid}, we get

\begin{cor}
The diagram product $D(\mathcal{P},\mathcal{G},w)$ is isomorphic to the fundamental group of $S(\mathcal{P},\mathcal{G},w)$.
\end{cor}

\begin{cor}\label{cor:changingword}
If $w_1,w_2 \in \Sigma^+$ are equal modulo $\mathcal{P}$, then the diagram products $D(\mathcal{P}, \mathcal{G},w_1)$ and $D(\mathcal{P}, \mathcal{G}, w_2)$ are isomorphic.
\end{cor}

\begin{proof}
If $w_1$ and $w_2$ are equal modulo $\mathcal{P}$, then they belong to the same connected component of $S(\mathcal{P}, \mathcal{G})$. A path from $w_1$ to $w_2$ in $S(\mathcal{P}, \mathcal{G})$ corresponds naturally to a semigroup diagram $\Delta$ satisfying $\mathrm{top}^-(\Delta)=w_1$ and $\mathrm{bot}^{-}(\Delta)=w_2$. Then
$$A \mapsto \Delta^{-1} \cdot A \cdot \Delta$$
defines an isomorphism $D(\mathcal{P},\mathcal{G},w_1) \to D(\mathcal{P}, \mathcal{G},w_2)$ whic corresponds to changing the basepoint of our fundamental group.
\end{proof}

\noindent
In particular, we recover the original definition of diagram products as 2-complexes of groups due to Guba and Sapir. (In fact, the complex $S(\mathcal{P},\mathcal{G},w)$ is precisely the 2-complex described at the begining of \cite[Section 2]{MR1725439} in the case of the complex of groups described by \cite[Definition 3]{MR1725439}.)

\begin{cor}\label{cor:2complex}
The diagram product $D(\mathcal{P},\mathcal{G},w)$ is isomorphic to the fundamental group of the following 2-complex of groups:
\begin{itemize}
	\item the underlying 2-complex is the 2-skeleton of the Squier complex $S(\mathcal{P},w)$;
	\item to any vertex $u=s_1\cdots s_r \in \Sigma^+$ is associated the group $G_u=G_{s_1} \times \cdots \times G_{s_r}$;
	\item to any edge $e=(a, u \to v,b)$ is associated the group $G_{e}=G_a \times G_b$;
	\item to any square is associated the trivial group;
	\item for every edge $e=(a,u \to v,b)$, the monomorphisms $G_e \to G_{aub}$ and $G_e \to G_{avb}$ are the canonical maps $G_a \times G_b \to G_a \times G_u \times G_b$ and $G_a \times G_b \to G_a \times G_v \times G_b$.
\end{itemize}
\end{cor}

\subsection{Cubical-like geometry}

\noindent
Let $\mathcal{P}= \langle \Sigma \mid \mathcal{R} \rangle$ be a semigroup presentation and $\mathcal{G}$ a collection of groups indexed by our alphabet $\Sigma$. Given a semigroup diagram $\Delta$, define the \emph{length} $\# (\Delta)$ of $\Delta$ as the sum of the number of cells of $\Delta$ with the number of edges of $\Delta$ labelled by a letter of $\Sigma(\mathcal{G})$ with a non trivial second coordinate. For instance, the diagram given by Figure \ref{figure13} has length seven. Say that $\Delta$ is \emph{unitary} if $\#(\Delta)=1$. Notice that unitary diagrams are precisely the atomic diagrams together with the linear diagrams of the form 
$$\epsilon((\ell_1,1) \cdots (\ell_{i-1},1) (\ell_i,g_i) (\ell_{i+1},1) \cdots (\ell_n,1)),$$
where $g_i \neq 1$. Since a linear diagram is cleary a concatenation of unitary linear diagrams, we deduce from Lemma \ref{lem:elementarydiag} that the set of unitary diagrams generate the groupoid $G(\mathcal{P}, \mathcal{G})$. We denote by $X(\mathcal{P},\mathcal{G})$ the Cayley graph of the groupoid $G(\mathcal{P},\mathcal{G})$, with respect to the set of unitary diagrams. Explicitely, $X(\mathcal{P}, \mathcal{G})$ is the graph whose vertices are the reduced semigroup diagrams and whose edges link two diagrams $\Delta_1$ and $\Delta_2$ whenever there exists a unitary diagram $E$ such that $\Delta_2= \Delta_1 \cdot E$. If $w \in \Sigma^+$ is a base word, let $X(\mathcal{P}, \mathcal{G},w)$ denote the connected component of $X(\mathcal{P}, \mathcal{G})$ containing $\epsilon(w)$. 

\medskip \noindent
Notice that two reduced diagrams belong to the same connected component of $X(\mathcal{P}, \mathcal{G})$ if and only if $\mathrm{top}^-(\Delta_1)= \mathrm{top}^-(\Delta_2)$, since the right-multiplication does not modify the top path $\mathrm{top}^-$ of a given diagram, so $X(\mathcal{P}, \mathcal{G},w)$ is the connected component of $X(\mathcal{P}, \mathcal{G})$ corresponding to the diagrams $\Delta$ satisfying $\mathrm{top}^-(\Delta)=w$.

\medskip \noindent
The diagram product $D(\mathcal{P},\mathcal{G},w)$ acts isometrically on $X(\mathcal{P}, \mathcal{G},w)$ by left-multiplication. 

\medskip \noindent
It is worth noticing that, fixing a semigroup diagram $\Delta_0$ satisfying $\mathrm{top}^-(\Delta_0)=w_1$ and $\mathrm{bot}^- (\Delta_0)=w_2$, then the map $\Delta \mapsto \Delta_0^{-1} \cdot \Delta \cdot \Delta_0$ defines an isomorphism $D(\mathcal{P}, \mathcal{G},w_1) \to D(\mathcal{P}, \mathcal{G},w_2)$ (see Corollary \ref{cor:changingword}) and an isometry $X(\mathcal{P}, \mathcal{G}, w_1) \to X(\mathcal{P}, \mathcal{G},w_2)$, which is equivariant when looking at $D(\mathcal{P},\mathcal{G},w_1)$ and $D(\mathcal{P},\mathcal{G}, w_2)$ acting on $X(\mathcal{P},\mathcal{G},w_1)$ and $X(\mathcal{P}, \mathcal{G},w_2)$ respectively. As a consequence, we may always suppose that a fixed base vertex is a trivial diagram up to changing the base word, which does not disturb either the group (up to isomorphism) or the graph (up to isometry) or the action (up to equivariant isometry). This observation will be useful to study the geometry of $X(\mathcal{P}, \mathcal{G},w)$, which turns out to be the quasi-median graph we will be interested in.

\begin{prop}\label{prop:XPquasimedian}
$X(\mathcal{P}, \mathcal{G},w)$ is a quasi-median graph.
\end{prop}

\noindent
The first step toward the proof of this proposition is to understand the geodesics in $X(\mathcal{P}, \mathcal{G},w)$. We begin with the following definition.

\begin{definition}
Let $A_1, \ldots, A_n$ be diagrams such that the concatenation $A_1 \circ \cdots \circ A_k$ is well-defined for every $1 \leq k \leq n$. If 
$$\#(A_1 \circ \cdots \circ A_n) = \# A_1 + \cdots + \# A_n,$$
we say that the concatenation $A_1 \circ \cdots \circ A_n$ is \emph{absolutely reduced}. Equivalently, such a concatenation is absolutely reduced if it is reduced and if two edges of some $\mathrm{bot}(A_i)$ and $\mathrm{top}(A_j)$ are identified in the concatenation then they are never both labelled by letters of $\Sigma(\mathcal{G})$ with non trivial second coordinates. If a diagram $\Delta$ decomposes as an absolutely reduced concatenation $A \circ B$, we say that $A$ is a \emph{prefix} of $\Delta$, written by $A \leq \Delta$, and that $B$ is a \emph{suffix} of $\Delta$.
\end{definition}

\begin{lemma}\label{lem:XPgeodesic}
Let $\Delta_1, \Delta_2 \in X(\mathcal{P}, \mathcal{G})$ be two diagrams. Write $\Delta_1^{-1}\Delta_2$ as an absolutely reduced concatenation $E_1 \circ \cdots \circ E_n$ of unitary diagrams. Then the sequence of vertices $$\Delta_1, \ \Delta_1 \cdot E_1, \ \Delta_1 \cdot E_1 \cdot E_2, \ldots, \Delta_1 \cdot E_1 \cdots E_n= \Delta_2$$ defines a geodesic between $\Delta_1$ and $\Delta_2$ in $X(\mathcal{P}, \mathcal{G})$. Conversely, any geodesic between $\Delta_1$ and $\Delta_2$ is labelled by an absolutely reduced concatenation of unitary diagrams representing the diagram $\Delta_1^{-1}\Delta_2$.
\end{lemma}

\noindent
Formally, the proof of the lemma is the same as the one of Lemma \ref{lem:geodesic}, since a concatenation of minimal length representing a given diagram must be absolutely reduced. This statement was proved for Farley CAT(0) cube complexes in \cite[Lemma 2.3]{arXiv:1505.02053}. As an immediate consequence, we deduce

\begin{cor}\label{cor:XPdist}
$d(\Delta_1,\Delta_2) = \# ( \Delta_1^{-1} \Delta_2)$ for every $\Delta_1, \Delta_2 \in X(\mathcal{P}, \mathcal{G})$. 
\end{cor}

\noindent
Our next preliminary lemma describes how behaves the length of a diagram after right-multiplication by a unitary diagram.

\begin{lemma}\label{lem:XPlength}
Let $\Delta$ be a reduced diagram and $E$ a unitary diagram. If $E$ is atomic
$$\#(\Delta \cdot E) = \left\{ \begin{array}{cl} \#(\Delta)-1 & \text{if $\Delta \circ E$ is not reduced}; \\ \#(\Delta)+1 & \text{if $\Delta \circ E$ is reduced}. \end{array} \right.$$
If $E$ is linear, write $\mathrm{bot}(\Delta)=(\ell_1,h_1) \cdots (\ell_n,h_n)$ and
$$E=\epsilon((\ell_1,1) \cdots (\ell_{i-1},1) (\ell_i,g) (\ell_{i+1},1) \cdots (\ell_n,1)),$$
where $g \in G_{\ell_i} \backslash \{ 1 \}$. Then
$$\#(\Delta \cdot E) = \left\{ \begin{array}{cl} \#(\Delta)+1 & \text{if} \ h_i=1; \\ \#(\Delta)-1 & \text{if} \ g=h_i^{-1}; \\ \#(\Delta) & \text{otherwise}. \end{array} \right.$$ 
\end{lemma}

\begin{proof}
Suppose that $E$ is atomic. Two cases may happen. If $\Delta \circ E$ is not reduced, then $\Delta \cdot E$ is obtained from $\Delta \circ E$ by reducing a dipole: two cells are removed, the removed edges (which are the edges of the intersection between the two cells defining our dipole) are labelled by letters of $\Sigma(\mathcal{G})$ whose second coordinates are trivial, and no edges labelled by letters of $\Sigma(\mathcal{G})$ with non trivial second coordinates are neither created nor removed since the bottom path of the bottom cell defining our dipole does not contain such edges. Therefore, $\#(\Delta \cdot E)=\#(\Delta)-1$. If $\Delta \circ E$ is reduced, then $\Delta \cdot E$ is obtained from $\Delta$ by adding a cell, where the new edges (which are the edges of the bottom path of our new cell) are labelled by letters of $\Sigma(\mathcal{G})$ whose second coordinates are trivial, hence $\#(\Delta \cdot E)=\#(\Delta)+1$.

\medskip \noindent
If $E$ is linear, the statement is clear.
\end{proof}

\noindent
Finally, we introduce a natural family of complete subgraphs in $X(\mathcal{P}, \mathcal{G})$, and we prove that they satisfy properties similar to cliques in quasi-median graphs. In fact, we will see later that they define one of the two families of cliques of $X(\mathcal{P},\mathcal{G})$. 

\begin{definition}
A \emph{coset} of $X(\mathcal{P},\mathcal{G})$ is a complete subgraph generated by the vertices
$$\Delta \cdot \epsilon ( (\ell_1,1) \cdots (\ell_{i-1},1) (\ell_i,g) (\ell_{i+1},1) \cdots (\ell_n,1)), \ g \in G_{\ell_i},$$
the diagram $\Delta$, which satisfies $\mathrm{bot}^-(\Delta)=\ell_1 \cdots \ell_n$, being fixed.
\end{definition}

\begin{lemma}\label{lem:XPcoset}
The intersection between two distinct cosets is either empty or a single vertex.
\end{lemma}

\begin{proof}
Let $C_1,C_2$ be two intersecting cosets. Let $\Delta \in C_1 \cap C_2$. Writing $\mathrm{bot}(\Delta)=(\ell_1,g_1) \cdots (\ell_n,g_n)$, say that
$$C_1= \{ \Delta \cdot \epsilon ( (\ell_1,1) \cdots (\ell_{i-1},1) (\ell_i,g) (\ell_{i+1},1) \cdots (\ell_n,1)), \ g \in G_{\ell_i} \},$$
$$C_2= \{ \Delta \cdot \epsilon ( (\ell_1,1) \cdots (\ell_{j-1},1) (\ell_j,g) (\ell_{i+1},1) \cdots (\ell_n,1)), \ g \in G_{\ell_j} \}.$$
Clearly, if $i=j$ then $C_1=C_2$, and $i \neq j$ then $\Delta$ is the single diagram in $C_1 \cap C_2$.
\end{proof}

\begin{lemma}\label{lem:XPtriangle}
Any triangle of $X(\mathcal{P}, \mathcal{G},w)$ is contained into a coset. 
\end{lemma}

\begin{proof}
Given a triangle, up to translating, we may suppose without loss of generality that $\epsilon(w)$ is one of its vertex; let $E,F$ denote its two other vertices. In particular, $E$, $F$ and $E^{-1}F$ label the three edges of our triangle. 

\medskip \noindent
Notice that, if $E$ and $F$ are two atomic diagrams, then $F \neq E^{-1}$ since there do not exist spherical atomic diagrams (because $\mathcal{R}$ contains no relations of the form $u=u$), so that $d(E,F) = \# (E^{-1}F) \geq 2$. Similarly, if $E$ is atomic and $F$ linear, or if $E$ is linear and $F$ is atomic, then $d(E,F)= \# (E^{-1}F) \geq 2$. Therefore, $E$ and $F$ must be two linear diagrams. Write $w=\ell_1 \cdots \ell_n$ and
$$E= (\ell_1,1) \cdots (\ell_{i-1},1) (\ell_i, g) (\ell_{i+1},1) \cdots (\ell_n,1),$$
$$F= (\ell_1,1) \cdots (\ell_{j-1},1) (\ell_j,h) (\ell_{j+1},1) \cdots (\ell_n,1).$$
Clearly, $d(E,F)=1$ implies that $i=j$. Notice that
$$E^{-1}F= (\ell_1,1) \cdots (\ell_{i-1},1) (\ell_i, g^{-1}h) (\ell_{i+1},1) \cdots (\ell_n,1).$$
We conclude that our triangle is included into a coset.
\end{proof}

\begin{proof}[Proof of Proposition \ref{prop:XPquasimedian}.]
First, we want to verify the triangle condition. So fix three diagrams $\Delta,\Gamma, \Xi \in X(\mathcal{P}, \mathcal{G},w)$ such that $\Delta$ and $\Gamma$ are adjacent and such that the distances $d(\Xi,\Gamma)$ and $d(\Xi,\Delta)$ are equal, say to $n$. Up to translating by $\Xi^{-1}$ (and replacing $w$ with $\mathrm{bot}^-(\Xi)$), we may suppose without loss of generality that $\Xi=\epsilon(w)$. Then, because $\Delta$ and $\Gamma$ are adjacent, there exists some unitary diagram $E$ such that $\Gamma= \Delta \cdot E$. From
$$\# (\Delta \cdot E) = d(\epsilon(w), \Gamma)=d(\epsilon(w),\Delta)= \#\Delta,$$
we deduce thanks to Lemma \ref{lem:XPlength} that $E$ must be a linear diagram. Write 
$$E=\epsilon((\ell_1,1) \cdots (\ell_{i-1},1) (\ell_i,g) (\ell_{i+1},1) \cdots (\ell_n,1),$$
where $\ell_1 \cdots \ell_n= \mathrm{bot}^-(\Delta)$ and $g \in G_{\ell_i} \backslash \{1 \}$. Also, write $\mathrm{bot}(\Delta)= (\ell_1,h_1) \cdots (\ell_n,h_n)$ where $h_j \in G_{\ell_j}$ for every $1 \leq j \leq n$. Because $\#\Delta= \# (\Delta \cdot E)$, necessarily $g \neq h_i^{-1}$ and $h_i \neq 1$. As a consequence, setting
$$\Delta_0= \Delta \cdot \epsilon((\ell_1,1) \cdots (\ell_{i-1},1) (\ell_n,h_i^{-1}) (\ell_{i+1},1) \cdots (\ell_n,1)),$$
as a vertex of $X(\mathcal{P},\mathcal{G},w)$, $\Delta_0$ is adjacent to $\Delta$ and $\Gamma$ and it satisfies 
$$d(\epsilon(w),\Delta_0)=\# \Delta_0= \# \Delta -1=n-1.$$
Thus, $\Delta_0$ is the vertex we are looking for.

\medskip \noindent
Next, we want to verify the quadrangle condition. So let $A,B,C,D \in X(\mathcal{P}, \mathcal{G},w)$ be four diagrams such that $B$ is adjacent to both $A$ and $C$, the distances $d(D,A)$ and $d(D,C)$ are equal, say to $n$, and $d(D,B)=n+1$. Once again, up to translating by $D^{-1}$ and replacing $w$ with $\mathrm{bot}^-(D)$, we will suppose that $D= \epsilon(w)$. Because $B$ is adjacent to $A$ and $C$, there exist two unitary diagrams $E_1,E_2$ such that $B=A \cdot E_1$ and $C=B \cdot E_2$. 

\medskip \noindent
We distinguish four cases. Suppose first that $E_1$ and $E_2$ are atomic diagrams. Thanks to Lemma \ref{lem:XPlength}, we deduce from 
$$\# (A \cdot E_1)= \# B=d(\epsilon(w),B)=d(\epsilon(w),A)+1= \#A+1$$
that $B=A \circ E_1$; similarly, we deduce from
$$\# (C \cdot E_2^{-1}) = \#B=d(\epsilon(w),B)=d(\epsilon(w),C)+1= \#C+1$$
that $B= C \circ E_2^{-1}$. Thus, $C=(A \circ E_1) \cdot E_2$, where the cell of $E_2$ must define a dipole with some cell of $A \circ E_1$. This cannot be the cell of $E_1$, since otherwise we would have $C=A$. Therefore, we can write $A=A_0 \circ E_3$ for some atomic diagram $E_3$, such that $C=(A_0 \circ E_3 \circ E_1) \cdot E_2$ where the cells of $E_2$ and $E_3$ define a dipole. Notice that, necessarily, the top paths of the cells of $E_1$ and $E_3$ intersect $\mathrm{bot}(A_0)$ along two disjoint subpaths, so that we can write $A_0 \circ E_3 \circ E_1$ as $A_0 \circ E_1' \circ E_3'$, where $E_1',E_3'$ are two atomic diagrams whose cells coincide with the celles of $E_1,E_3$ respectively. By reducing the concatenation $(A_0 \circ E_3 \circ E_1) \cdot E_2$, we get $C=A_0 \cdot E_1'$; as a consequence, $A_0$ is adjacent to $C$. Moreover, because $A=A_0 \circ E_3$, we know that $A_0$ and $A$ are adjacent and that
$$d(\epsilon(w),A_0)=\#A_0 = \#A-1=d(\epsilon(w),A)-1=n-1.$$
Consequently, $A_0$ is the vertex we are looking for.

\medskip \noindent
Now, suppose that $E_1$ and $E_2$ are linear. Write $\mathrm{bot}(A)=(\ell_1,g_1) \cdots (\ell_n,g_n)$ and
$$E_1=(\ell_1,1) \cdots (\ell_{i-1},1)(\ell_i,h) (\ell_{i+1},1) \cdots (\ell_n,1),$$
$$E_2=(\ell_1,1) \cdots (\ell_{j-1},1)(\ell_j,k) (\ell_{j+1},1) \cdots (\ell_n,1),$$
where $h \in G_{\ell_i} \backslash \{ 1 \}$ and $k \in G_{\ell_j} \backslash \{ 1 \}$. Thanks to Lemma \ref{lem:XPlength}, we deduce from
$$\#(A \cdot E_1)= d(\epsilon(w),B)=d(\epsilon(w),A)+1= \#A+1$$
that $g_i=1$; similarly, we deduce from
$$\#(B \cdot E_2)= d(\epsilon(w),C)= d(\epsilon(w),B)-1= \#B-1$$
that $k=g_j^{-1}$. In particular, $i \neq j$ so that $E_1 \cdot E_2= E_2 \cdot E_1$. Set $A'=A \cdot E_2$; of course, $A'$ and $A$ are adjacent. Moreover, $A'$ and $C$ are adjacent since
$$C=B \cdot E_2 = A \cdot E_1 \cdot E_2= A \cdot E_2 \cdot E_1 = A' \cdot E_1,$$
and $d(\epsilon(w),A')=\#A'=\#A-1=n-1$, so that $A'$ is the vertex we are looking for.

\medskip \noindent
Next, suppose that $E_1$ is atomic and $E_2$ linear. Write $\mathrm{bot}(B)=(\ell_1,h_1) \cdots (\ell_n,h_n)$ and
$$E_2= \epsilon( (\ell_1,1) \cdots (\ell_{i-1},1) (\ell_i,g) (\ell_{i+1},1) \cdots (\ell_n,1))$$
where $g \in G_{\ell_i} \backslash \{1 \}$. Thanks to Lemma \ref{lem:XPlength}, we deduce from
$$\#(A \cdot E_1)= d(\epsilon(w),B)=d(\epsilon(w),A)+1=\#A+1$$
that $B=A \circ E_1$; similarly, from
$$\#(B \cdot E_2)=d(\epsilon(w),C)=d(\epsilon(w),B)-1=\#B-1$$
we deduce that $g=h_i^{-1}$. As a consequence, $h_i \neq 1$ since $E_2$ not trivial implies that $g \neq 1$. Therefore, because the edges of the bottom path of the cell of $E_1$ are labelled by letters of $\Sigma(\mathcal{G})^+$ whose second coordinates are trivial, necessarily the edge $e_i$ of $\mathrm{bot}(B)$ labelled by $(\ell_i,h_i)$ cannot belong to the top path of the cell of $E_1$. As a consequence, if we write $\mathrm{bot}(A)=(s_1,k_1) \cdots (s_m,k_m)$, then there exists some $1 \leq j \leq m$ such that the edge labelled by $(s_j,k_j)$ is precisely $e_i$. Setting
$$L= \epsilon ( (s_1,1) \cdots (s_{j-1},1) (s_j,k_j^{-1}) (s_{j+1},1) \cdots (s_m,k_m)),$$
we have $E_1 \cdot E_2=L \cdot E_1$. Thus, if $A'=A \cdot L$, then
$$C=A \cdot E_1 \cdot E_2= A \cdot L \cdot E_1= A' \cdot E_1,$$
so that $A'$ and $C$ are adjacent vertices. Moreover, $A$ and $A'$ are clearly adjacent as well, and our choice of $L$ implies, thanks to Lemma \ref{lem:XPlength}, that $d(\epsilon(w),A')=\#A'= \#A-1=n-1$. Thus, $A'$ is the vertex we are looking for.

\medskip \noindent
Finally, if $E_1$ is atomic and $E_2$ linear, the situtation is symmetric to the previous one. This concludes the proof of the quadrangle condition.

\medskip \noindent
Now, notice that, if $\Delta_1,\Delta_2 \in X(\mathcal{P}, \mathcal{G},w)$ are two vertices satisfying $d(\Delta_1,\Delta_2)=2$, then there exists at most two geodesics between $\Delta_1$ and $\Delta_2$. Indeed, according to Lemma \ref{lem:XPgeodesic}, there exists a bijection between the geodesics between $\Delta_1$ and $\Delta_2$ and the reduced concatenation of unitary diagrams representing $\Delta_1^{-1}\Delta_2$. Since $\#(\Delta_1^{-1}\Delta_2)=d(\Delta_1,\Delta_2)=2$, the conclusion follows. As a consequence, $X(\mathcal{P}, \mathcal{G},w)$ cannot contain an induced subgraph isomorphic to $K_{2,3}$.

\medskip \noindent
Finally, let $(A,B,C,D)$ be a square in $X(\mathcal{P}, \mathcal{G},w)$, where $B$ and $D$ are adjacent. According to Lemma \ref{lem:XPtriangle}, the triangles $(A,B,D)$ and $(B,C,D)$ are included into two cosets of $X(\mathcal{P}, \mathcal{G})$. On the other hand, two cosets intersecting along an edge must be equal, so we deduce that the vertices $A$, $B$, $C$ and $D$ belong to the same coset of $X(\mathcal{P}, \mathcal{G})$. In particular, $A$ and $C$ must be adjacent in $X(\mathcal{P},\mathcal{G},w)$. This proves that $X(\mathcal{P},\mathcal{G},w)$ does not contain any induced subgraph isomorphic to $K^-_4$.
\end{proof}

\noindent
If all the groups of $\mathcal{G}$ are trivial, $X(\mathcal{P}, \mathcal{G},w)$ coincides with the 1-skeleton of the cube complex complex introduced by Farley in \cite{MR1978047}. In particular, we recover a result of Farley stating that this cube complex is CAT(0) \cite[Theorem 3.13]{MR1978047}

\begin{cor}
If all the groups of $\mathcal{G}$ are trivial, then $X(\mathcal{P},\mathcal{G},w)$ is a median graph.
\end{cor}

\begin{proof}
If all the groups of $\mathcal{G}$ are trivial, then there do not exist linear unitary diagrams, so that $X(\mathcal{P}, \mathcal{G},w)$ must be triangle-free according to Lemma \ref{lem:XPtriangle}. We conclude that $X(\mathcal{P}, \mathcal{G},w)$ is a median graph thanks to Corollary \ref{cor:whenmedian} combined with Proposition \ref{prop:XPquasimedian}.
\end{proof}

\paragraph{Cliques and prisms of $X(\mathcal{P}, \mathcal{G})$.} Below, our goal is to describe the cliques and the prisms of our quasi-median graph. We begin with its cliques.

\begin{lemma}\label{lem:XPclique}
The cliques of $X(\mathcal{P}, \mathcal{G})$ are the cosets and the edges labelled by atomic diagrams.
\end{lemma}

\begin{proof}
Let $(A,A \cdot E)$ be an edge. If $E$ is atomic, according to Lemma \ref{lem:XPtriangle} this edge is already a clique. Otherwise, if $E$ is linear, our edge belongs to a coset $C$. If this coset is not a clique of $X(\mathcal{P},\mathcal{G})$, there must exist a vertex $\Delta$ adjacent to all the vertices of $C$. Because $C$ has cardinality at least two, this produces a triangle with $\Delta$ as a vertex which shares an edge with $C$. On the other hand, according to Lemma \ref{lem:XPtriangle}, this triangle is included into a coset $C'$, and according to Lemma \ref{lem:XPcoset}, we must have $C=C'$ since $C$ and $C'$ share an edge, hence $\Delta \in C$. This proves that our coset is a clique of $X(\mathcal{P}, \mathcal{G})$.
\end{proof}

\noindent
Now, we focus on prisms of $X(\mathcal{P}, \mathcal{G})$. We first need to introduce a few definitions.

\begin{figure}
\begin{center}
\includegraphics[scale=0.55]{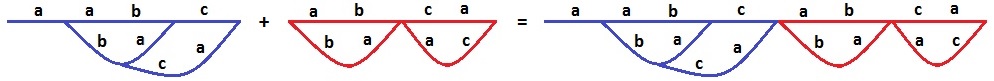}
\end{center}
\caption{A sum of two diagrams.}
\label{figure18}
\end{figure}

\begin{definition}
Given two diagrams $\Delta_1$ and $\Delta_2$, we define the \emph{sum} $\Delta_1+ \Delta_2$ as the diagram obtained by gluing the right-most vertex of $\Delta_1$ with the left-most vertex of $\Delta_2$. See Figure \ref{figure18}. A diagram is \emph{thin} if it is a sum of unitary diagrams.
\end{definition}

\begin{definition}
Let $\Delta$ be a reduced diagram. The \emph{support} of $\Delta$, denoted by $\mathrm{supp}(\Delta)$, is the union of the edges of $\mathrm{top}(\Delta)$ which are either labelled by letters of $\Sigma(\mathcal{G})$ with non trivial second coordinates or included into the top path of some cell of $\Delta$.
\end{definition}

\begin{prop}\label{prop:XPprism}
Let $P$ be a prism of finite cubical dimension of $X(\mathcal{P}, \mathcal{G},w)$, and let $\Delta$ denote the projection of $\epsilon(w)$ onto $P$. There exist edges $e_1,\ldots, e_k$ of $\mathrm{bot}(\Delta)$, respectively labelled by $y_1, \ldots, y_k \in \Sigma$, and a thin $(\mathrm{bot}(\Delta),\ast)$-diagram $A$ such that
$$P = \left\{ \Delta \cdot L_1(g_1) \cdots L_k(g_k) \cdot B \mid g_i \in G_{y_i}, \ B \leq A \right\}$$
and $\mathrm{supp}(A) \cap \{ e_1, \ldots, e_k \}= \emptyset$, where, for every $1 \leq i \leq k$ and every $g \in G_{y_i}$, $L_i(g)$ denotes the diagram $\epsilon(x_i)+ \epsilon((y_i,g)) + \epsilon(z_i)$ if we decompose $\mathrm{bot}^-(\Delta)$ as $x_iy_iz_i$.
\end{prop}

\noindent
Before proving Proposition \ref{prop:XPprism}, we need to understand the squares of $X(\mathcal{P}, \mathcal{G})$.

\begin{lemma}\label{lem:diagproductsquare}
Let $A,B,C,D$ be the vertices of some induced square in $X(\mathcal{P}, \mathcal{G},w)$. Up to permuting $B$ and $D$, there exist words $u,v \in \Sigma(\mathcal{G})^+$, a unitary $(u, \ast)$-diagram $U$ and a unitary $(v,\ast)$-diagram $V$, such that $\mathrm{bot}(A)=uv$, $B=A \circ (\epsilon(u)+V)$, $D= A \circ ( U+ \epsilon(v))$ and $C= A \circ (U+V)$. 
\end{lemma}

\begin{proof}
Because our square is induced, we know from Lemma \ref{lem:XPgeodesic} that $C=A \cdot T$ for some semigroup diagram $T$ of length two. Five cases may happen:
\begin{itemize}
	\item[(i)] $T$ is linear, ie., $T= \epsilon( w_1 \cdot (\ell_1,g_1) \cdot w_2 \cdot (\ell_2,g_2) \cdot w_3)$ for some letters $\ell_1, \ell_2 \in \Sigma$, some non trivial elements $g_1 \in G_{\ell_1}$ and $g_2 \in G_{\ell_2}$, and some words $w_1,w_2,w_3 \in \Sigma(\mathcal{G})^+$ whose letters have trivial second coordinates;
	\item[(ii)] $T$ has a single cell $\pi$, only one of its edge, say $e$, is labelled by a letter of $\Sigma(\mathcal{G})$ with a non trivial second coordinate, and $e$ does not belong to $\mathrm{top}(\pi)$ or $\mathrm{bot}(\pi)$;
	\item[(iii)] $T$ has a single cell $\pi$, only one of its edge, say $e$, is labelled by a letter of $\Sigma(\mathcal{G})$ with a non trivial second coordinate, and $e$ belongs to $\mathrm{top}(\pi)$ or $\mathrm{bot}(\pi)$; 
	\item[(iv)] $T$ has two cells $\pi_1$ and $\pi_2$, its edges are labelled by letters of $\Sigma(\mathcal{G})$ with trivial second coordinates, and the intersection between the boundaries of $\pi_1$ and $\pi_2$ does not contain an edge;
	\item[(v)] $T$ has two cells $\pi_1$ and $\pi_2$, its edges are labelled by letters of $\Sigma(\mathcal{G})$ with trivial second coordinates, and the intersection between the boundaries of $\pi_1$ and $\pi_2$ contains an edge.
\end{itemize}
On the other hand, since there exist (at least) two geodesics between $A$ and $C$, we deduce from Lemma \ref{lem:XPgeodesic} that $T$ can be written as a concatenation of unitary diagrames in (at least) two different ways. This is possible only in the cases $(i)$, $(ii)$, $(iv)$ above. In particular, it is possible to write $T$ as a sum $U+V$, where $\mathrm{top}(T)=uv$ in $\Sigma(\mathcal{G})^+$, and where $U$ is a non trivial unitary $(u,\ast)$-diagram and $V$ a non trivial unitary $(v,\ast)$-diagram. Thus,
$$A, \ A \cdot (\epsilon(u)+V), \ A \cdot (U+V) \ \text{and} \ A, \ A \cdot  (U+ \epsilon(v)), \ A \cdot (U+V)$$
define two geodesics between $A$ and $C$. The conclusion follows from the following observation.

\begin{claim}
In a quasi-median graph, there exist at most two geodesics between two vertices at distance two appart.
\end{claim}

\noindent
Suppose by contradiction that there exist a quasi-median graph $X$ containing two vertices $x,y \in X$ at distance two appart which can be joined by three distinct geodesics. So there exist three distinct vertices $z_1,z_2,z_3$ which adjacent to both $x$ and $y$. By noticing that the these five vertices produce a subgraph isomorphic to $K_{2,3}$, we deduce that either $x$ and $y$ must be adjacent, which is impossible since $d(x,y)=2$, or two vertices among $z_1,z_2,z_3$ must be adjacent, say $z_1$ and $z_2$. Now, the four vertices $x,y,z_1,z_2$ produce a subgraph isomorphic to $K_4^-$, so necessarily $x$ and $y$ must be adjacent, providing a contradiction.
\end{proof}

\begin{proof}[Proof of Proposition \ref{prop:XPprism}.]
Write $P$ as a product of cliques $C_1 \times \cdots \times C_n$ such that $C_1 \cap \cdots \cap C_n = \{ \Delta \}$. According to Lemma \ref{lem:XPclique}, each clique $C_i$ is either a coset or an edge labelled by an atomic diagram. Up to a permutation of the factors, suppose that there exists some $0 \leq k \leq n$ such that $C_1, \ldots, C_k$ are cosets and $C_{k+1}, \ldots, C_n$ are edges labelled by atomic diagrams. So, for every $1 \leq i \leq k$, we can write
$$C_i = \left\{ \Delta \cdot (\epsilon(x_i) + \epsilon((y_i,g)) + \epsilon(z_i)) \mid g \in G_{y_i} \right\},$$
where $x_i,z_i \in \Sigma^+$ are words and $y_i \in \Sigma$ a letter such that $\mathrm{bot}^-(\Delta)=x_iy_iz_i$; and, for every $k+1 \leq i \leq n$, we can write
$$C_i = \left\{ \Delta, \ \Delta \cdot (\epsilon(x_i) + A_i + \epsilon(z_i)) \right\}$$
where $x_i,y_i,z_i \in \Sigma^+$ are words such that $\mathrm{bot}^-(\Delta)=x_iy_iz_i$ and where $A_i$ is a single $(y_i,\ast)$-cell. For every distinct $1 \leq i,j \leq n$, we know that $C_i$ and $C_j$ generate a prism, which contains in particular a square; so it follows from Lemma \ref{lem:diagproductsquare} that $y_i$ and $y_j$ must be disjoint subpaths of $\mathrm{bot}(\Delta)$. Next, for every $1 \leq i \leq k$ and every $g \in G_{y_i}$, set 
$$L_i(g)= \epsilon(x_i)+ \epsilon((y_i,g)) + \epsilon(z_i).$$
Up to permuting the factors $C_{k+1}, \ldots, C_k$, suppose that $\mathrm{bot}^-(\Delta)$ decompose as 
$$a_1y_{k+1} \cdots a_{n-k}y_na_{n-k+1},$$
and set
$$A= \epsilon(a_1)+ A_{k+1}+ \cdots + \epsilon(a_{n-k}) + A_n + \epsilon(a_{n-k+1}).$$
By construction,
$$\{ \Delta \cdot L_1(g_1) \cdots L_k(g_k) \cdot B \mid g_i \in G_{y_i}, \ B \leq A \}$$
is a prism of cubical dimension $n$ which contains the cliques $C_1, \ldots, C_n$. Therefore, this is our prism $P$, concluding the proof.
\end{proof}

\paragraph{Hyperplanes of $X ( \mathcal{P}, \mathcal{G})$.}

\noindent
There are two families of edges in $X(\mathcal{P}, \mathcal{G})$: those which are labelled by atomic diagrams, and those which are labelled by linear diagrams. Interestingly, the combination of Lemmas \ref{lem:XPclique} and \ref{lem:diagproductsquare} shows that the edges of a given hyperplane belong to a single family. Thus, we distinguish two kind of hyperplanes: the \emph{linear hyperplanes} and the \emph{atomic hyperplanes}, depending on the family to which their edges belong. Our first result describe the linear hyperplanes of $X(\mathcal{P},\mathcal{G})$. 

\begin{prop}\label{prop:linearhyp}
Let $J$ be a linear hyperplane of $X(\mathcal{P}, \mathcal{G},w)$ and $\Delta \in N(J)$ a diagram. There exist letters $(a_1,\ell_1), \ldots, (a_r, \ell_r),(b,m), (c_1,n_1), \ldots, (c_s,n_s) \in \Sigma( \mathcal{G})$ such that
$$\mathrm{bot}(\Delta)= (a_1,\ell_1) \cdots (a_r,\ell_r) \cdot (b,m) \cdot (c_1,n_1) \cdots (c_s,n_s),$$
and such that an edge is dual to $J$ if and only if it has the form
$$ [ \Delta \cdot (A+ \epsilon((b,g)) +C), \ \Delta \cdot (A + \epsilon((b,h))+ C)]$$
where $A\in X(\mathcal{P}, \mathcal{G},a_1 \cdots a_r)$ and $C \in X(\mathcal{P}, \mathcal{G},c_1 \cdots c_s)$ are two diagrams, and $g,h \in G_b$ two distinct elements. As a consequence, the fibers of $J$ are the subgraphs generated by
$$\left\{ \Delta \cdot ((A+ \epsilon((b,g)) +C) \left| A\in X(\mathcal{P}, \mathcal{G},a_1 \cdots a_r), C \in X(\mathcal{P}, \mathcal{G},c_1 \cdots c_s) \right. \right\}$$
where $g \in G_b$ is fixed. Moreover, if $\Delta$ is the projection of $\epsilon(w)$ onto $N(J)$ then $\ell_1 = \cdots = \ell_r = n_1 = \cdots =n_s=1$.
\end{prop}

\begin{proof}
Fix some letters $(a_1,\ell_1), \ldots, (a_r, \ell_r),(b,m), (c_1,n_1), \ldots, (c_s,n_s) \in \Sigma( \mathcal{G})$ such that 
$$\mathrm{bot}(\Delta)= (a_1,\ell_1) \cdots (a_r,\ell_r) \cdot (b,m) \cdot (c_1,n_1) \cdots (c_s,n_s),$$
and such that the unique clique dual to $J$ and containing $\Delta$ is the coset
$$\{ \Delta \cdot \epsilon( (a_1,\ell_1) \cdots (a_r,\ell_r) \cdot (b,g) \cdot (c_1,n_1) \cdots (c_s,n_s)) \mid g \in G_b \}.$$
For convenience, set $a=a_1 \cdots a_r$ and $c=c_1 \cdots c_s$, and let $\Delta(A,g,C)$ denote the diagram $\Delta \cdot (A+ \epsilon((b,g))+C)$, where $A$ and $C$ belong to $X(\mathcal{P}, \mathcal{G},a)$ and $X(\mathcal{P}, \mathcal{G},c)$ respectively, and where $g \in G_b$. 

\medskip \noindent
Let $e$ be an edge dual to $J$. According to Claim \ref{claim:chainofsquareshyp} proved below, there exists a sequence of edges 
$$e_0=f, \ e_1, \ldots, \ e_{n-1}, \ e_n=e$$
such that, for every $1 \leq i \leq n$, the edges $e_{i-1}$ and $e_i$ are the opposite edges of some square, where $f= [\Delta( \epsilon(a),g, \epsilon(c)), \Delta(\epsilon(a) , h, \epsilon(c))]$ for some distinct $g,h \in G_b$. Up to inversing $f$, suppose that $h \neq 1$. We want to prove by induction on $n$ that there exist two diagrams $A \in X(\mathcal{P}, \mathcal{G},a)$ and $C \in X(\mathcal{P}, \mathcal{G},c)$  such that 
$$e=[\Delta(A , g, C), \Delta(A,h,C)].$$
If $n=0$, then $e=f$ and there is nothing to prove. Next, suppose that $n \geq 1$. By our induction hypothesis, we know that there exist two diagrams $A \in X(\mathcal{P}, \mathcal{G},a)$ and $C \in X(\mathcal{P}, \mathcal{G},c)$ such that $e_{n-1}=[\Delta(A, g, C), \Delta(A, h,C)]$. Notice that $e_{n-1}$ and $e_n=e$ belong to a common square $Q$. Supposing that our sequence of edges has minimal length (which can be done without loss of generality), $e_{n-1}$ is nearer to $\epsilon(w)$ than $e$, so that $\Delta(A,g,C)$ is the nearest vertex of $Q$ to $\epsilon(w)$. By applying Lemma \ref{lem:diagproductsquare} to $Q$, we deduce that $e$ has the form
$$[\Delta(A' , g, C), \ \Delta(A', h, C)] \ \text{or} \ [\Delta(A,g ,C'), \ \Delta(A, h, C')],$$
where $A'$ (resp. $C'$) is a semigroup diagram obtained from $A$ (resp. $C$) by right-multiplication with a unitary diagram. Conversely, consider an edge
$$e = [ \Delta  (A, g,C), \ \Delta  (A, h, C)]$$
where $A \in X(\mathcal{P}, \mathcal{G},a)$, $C \in X(\mathcal{P}, \mathcal{G},c)$, and where $g,h \in G_b$ are distinct. Write $A$ (resp. $C$) as a concatenation of unitary diagrams $A_1 \circ \cdots \circ A_p$ (resp. $C_1 \circ \cdots \circ C_q$). Now, for every $1 \leq i \leq p$, let $e_i$ denote the edge
$$[\Delta (A_1 \circ \cdots \circ A_i, g, \epsilon(c)), \ \Delta (A_1 \circ \cdots \circ A_i, h, \epsilon(c))],$$
and for every $1 \leq i \leq q$, let $f_i$ denote the edge
$$[\Delta (A,g, C_1 \circ \cdots \circ C_i), \ \Delta (A,h, C_1 \circ \cdots \circ C_i))],$$
and notice that 
$$e_1, \ e_2, \ldots, \ e_{p}, \ f_1, \ f_2, \ldots , \ f_q=e$$
defines a sequence of edges between $e_1 \in J$ and $e$ such that any two consecutive edges are two opposite edges of some square. Therefore, $e$ belongs to $J$. 

\medskip \noindent
This proves the description of the edges of $J$. Notice that this also implies that $\ell_1, \ldots, \ell_r$, $m$, $n_1, \ldots, n_s$ are all trivial is $\Delta$ coincides with the projection of $\epsilon(w)$ onto $N(J)$. Indeed,
$$[\Delta, \ \Delta \cdot \epsilon( (a_1, \ell_1^{-1}) \cdots (a_r,\ell_r^{-1}) \cdot (b,m^{-1}) \cdot (c_1,n_1^{-1}) \cdots (c_s,n_s^{-1}))]$$
is an edge dual to $J$, and the length of the second endpoint, ie., the distance to $\epsilon(w)$, is smaller than the length of $\Delta$ if one of $\ell_1, \ldots, \ell_r, m ,n_1, \ldots, n_s$ is non trivial, which is impossible since $\Delta$ is the unique vertex of $N(J)$ which minimises the distance to $\epsilon(w)$. 

\medskip \noindent
Finally, the description of the fibers of $J$ follows from the following claim, where the clique $\{ \Delta(\epsilon(a),g,\epsilon(c)) \mid g \in G_b\}$ is denoted by $K$. 

\begin{claim}
The restriction to $N(J)$ of the projection onto the clique $K$ is defined by
$$\Delta(A,g,C) \mapsto \Delta(\epsilon(a),g,\epsilon(c))$$
\end{claim}

\noindent
The distance between $\Delta(A,g,C)$ and a vertex $\Delta(\epsilon(a),h,\epsilon(c))$ of $K$, according to Corollary \ref{cor:XPdist}, is equal to
$$\# A + \# C +  \epsilon(g,h), \ \text{where} \ \epsilon(g,h)= \left\{ \begin{array}{cl} 1 & \text{if} \ g \neq h \\ 0 & \text{if} \ g=h \end{array} \right. .$$
Therefore, $\Delta(\epsilon(a),g,\epsilon(c))$ is the vertex of $K$ minimising the distance to $\Delta(A,g,C)$. This concludes the proof of our claim.
\end{proof}

\begin{claim}\label{claim:chainofsquareshyp}
Let $X$ be a quasi-median graph, $C$ a clique and $J$ its dual hyperplane. For every edge $e$ dual to $J$, there exist an edge of $C$ and a chain of adjacent squares from it to $e$. 
\end{claim}

\begin{proof}
Let $x,y$ denote the endpoints of $e$ and $x',y'$ their respective projections onto $C$. Of course, $x'$ belongs to the interval $I(x,y')$ and $y'$ to the interval $I(y,x')$. Moreover, by noticing that $d(x,y')=d(x,x')+1$, we deduce that the concatenation of a geodesic from $y'$ to $y$ with $e$ defines a geodesic from $y'$ to $x$, so that $y$ belongs to the interval $I(x,y')$; similarly, $x \in I(y,x')$. It follows from Lemma \ref{lem:productingflatrectangle} that there exists a flat rectangle $[0,1] \times [0,n] \hookrightarrow X$ such that $[0,1] \times \{ 0 \}=e$, $(0,n)=x'$ and $(1,n)=y'$. This precisely means that there exists a chain of adjacent squares between $e$ and the edge of $C$ linking $x'$ and $y'$. 
\end{proof}

\noindent
As a consequence of Proposition \ref{prop:linearhyp}, we are also able to describe stabilisers of linear hyperplanes.

\begin{cor}\label{cor:XPstabhyp}
Let $J$ be a linear hyperplane. With respect to the notations of the previous proposition, the stabiliser of $J$ is equal to
$$\left\{ \Delta \cdot  \left( A+ \epsilon(b,g) + C \right) \cdot \Delta^{-1} \mid A \in D(\mathcal{P}, \mathcal{G},a), C \in D(\mathcal{P}, \mathcal{G},c), g \in G_b \right\},$$
where $a=a_1 \cdots a_r$ and $c=c_1 \cdots c_s$. In particular, $\mathrm{stab}(J)$ is naturally isomorphic to the product $D(\mathcal{P}, \mathcal{G},a) \times G_b \times D(\mathcal{P}, \mathcal{G},c)$. 
\end{cor}

\begin{proof}
Up to translating by $\Delta^{-1}$, we will suppose for convenience that $\Delta$ is trivial. Let $S \in \mathrm{stab}(J)$. For every $g \in G_b \backslash \{ 1 \}$, there must exist some diagram $A \in X(\mathcal{P},\mathcal{G},a_1 \cdots a_r)$, some diagram $C \in X(\mathcal{P}, \mathcal{G},c_1 \cdots c_s)$, and some distinct elements $h,k \in G_b$ such that
$$S = S \cdot (\epsilon(a)+ \epsilon(b)+ \epsilon(c)) = A+ \epsilon(b,h)+C$$
and
$$S \cdot ( \epsilon(a)+ \epsilon(b,g)+ \epsilon(c)) = A+ \epsilon(b,k) +C.$$
Because $S$ is a spherical diagram, we get the equality
$$abc = \mathrm{top}^-(S) = \mathrm{bot}^-(S)=ubv$$
in $\Sigma^+$, where $u= \mathrm{bot}^-(A)$ and $v= \mathrm{bot}^-(C)$. We claim that $a=u$ and $c=v$, meaning that $A$ and $C$ are spherical. Notice that the equality $a=u$ follows from the equality $abc=ubv$ if we know that $|a|=|u|$. Suppose by contradiction that $|u|< |a|$. From the equality $abc=ubv$, we deduce that $a=up$ for some non empty word $p \in \Sigma^+$. Thus, the equality $abc=ubv$ becomes $pbc=bv$. In particular, $b$ must be the first letter of $p$, say $p=bq$ for some (possibly empty) word $q \in \Sigma^+$. So our equality becomes $qbc=v$. It follows that $\mathrm{top}^-(A)=ubq$, $\mathrm{bot}^-(A)=u$, $\mathrm{top}^-(C)=c$ and $\mathrm{bot}^-(C)=qbc$, so that
$$S \cdot  ( \epsilon(a)+ \epsilon(b,g)+ \epsilon(c)) = A+ \epsilon(b,h)+ (C \cdot (\epsilon(q)+ \epsilon(b,g)+ \epsilon(c))).$$
Thus, our equality
$$S \cdot ( \epsilon(a)+ \epsilon(b,g)+ \epsilon(c)) = A+ \epsilon(b,k) +C.$$
is equivalent to $ \epsilon(b,h)= \epsilon(b,k)$ and $g=1$, which is false. Therefore, $|u| \geq |a|$. The inequality $|u| \leq |a|$ is proved similarly, hence $|u|=|a|$, and finally $a=u$. Next, the equality $abc=ubv$ implies $c=v$, proving our claim. So far, we have proved that $\mathrm{stab}(J)$ is included into
$$\left\{ A+ \epsilon(b,g) + C \mid A \in D(\mathcal{P}, \mathcal{G},a), C \in D(\mathcal{P}, \mathcal{G},c), g \in G_b \right\}.$$
The reverse inclusion is clear.
\end{proof}

\noindent
In the following, a precise description of the atomic diagrams will be not necessary. We only prove the following lemma.

\begin{lemma}\label{lem:DPatomichyp}
An atomic hyperplane $J$ has exactly two fibers. Moreover, $\mathrm{stab}(J)$ stabilises them.
\end{lemma}

\begin{proof}
Let $(\Delta, \Delta \cdot A)$ be an edge dual to $J$, where $A$ is an atomic diagram, say labelled the relation $u \to v$; this means that there exist words $a,b \in \Sigma^+$ such that $\mathrm{top}^-(A)=aub$ and $\mathrm{bot}^-(A)=avb$, where $u=v$ is a relation of $\mathcal{R}$. According to Lemma \ref{lem:XPclique}, this is a clique, so $J$ has only two fibers, say $\partial_-$ and $\partial_+$ containing respectively $\Delta$ and $\Delta \cdot A$. We first notice that

\begin{claim}
If $P \in \partial_-$ and $Q \in \partial_+$ are two adjacent vertices, then $Q=P \cdot B$ for some atomic diagram labelled by the relation $u \to v$.
\end{claim}

\noindent
Because the edge $(P,Q)$ is dual to the same hyperplane as the edge $(\Delta, \Delta \cdot A)$, and that the hyperplane has only two fibers, there must exist a sequence of edges
$$(\Gamma_1, \Xi_1)=(\Delta,\Delta \cdot A), \ (\Gamma_2, \Xi_2), \ldots, \ (\Gamma_{n-1}, \Xi_{n-1}), \ (\Gamma_n, \Xi_n)=(P,Q)$$
such that $(\Gamma_i, \Xi_i)$ and $(\Gamma_{i+1}, \Xi_{i+1})$ are opposite edges of some square for every $1 \leq i \leq n-1$. We argue by induction on $n$. If $n=0$, there is nothing to prove. If $n \geq 1$, then our induction hypothesis implies that $\Xi_{n-1}= \Gamma_{n-1} \cdot C$ for some atomic diagram $C$ labelled by the relation $u \to v$. Now, from description of the squares of $X(\mathcal{P}, \mathcal{G},w)$ given by Lemma \ref{lem:diagproductsquare}, the conclusion follows: there exists an atomic diagram $B$ labelled by $u \to v$ such that $\Xi_n= \Gamma_n \cdot B$. This proves our claim.

\medskip \noindent
Now, suppose by contradiction that there exists some $g \in \mathrm{stab}(J)$ such that $g \cdot \Delta \cdot A \in \partial_-$ and $g \cdot \Delta \in \partial_+$. As a consequence of our claim, $g \cdot \Delta = g \cdot \Delta \cdot A \cdot B$ for some atomic diagram labelled by $u \to v$. It follows that $A \cdot B= \epsilon(z)$ where $z= \mathrm{top}^-(A)$, or equivalently $B=A^{-1}$. Because $A$ is labelled by $u \to v$, it follows that $B$ must be labelled by $v \to u$, a contradiction.
\end{proof}

\subsection{Properties of diagram products}

\noindent
In this section, we want to prove that a diagram product acts topically-transitively on its associated quasi-median graph, and next to apply the criteria stated in Sections \ref{section:topicalactionsI} and \ref{section:topicalactionsII}. 

\begin{prop}
The action of $D(\mathcal{P}, \mathcal{G},w)$ on $X(\mathcal{P}, \mathcal{G},w)$ is topical-transitive. 
\end{prop}

\begin{proof}
Let $C$ be a clique of $X(\mathcal{P}, \mathcal{G},w)$. According to Lemma \ref{lem:DPatomichyp}, if $C$ is an atomic edge, then its stabiliser is trivial, and the stabiliser of the hyperplane $J$ dual to $C$ does not permute the fibers of $J$. 

\medskip \noindent
Next, suppose that $C$ is a linear clique. Fix a vertex $\Delta \in C$. According to Lemma \ref{lem:XPclique}, $C$ is a coset, so that, if we write $\mathrm{bot}^-(\Delta)= \ell_1 \cdots \ell_n$, then
$$C = \left\{ \Delta \cdot \epsilon \left( (\ell_1,1) \cdots ( \ell_{i-1},1) (\ell_i,g) (\ell_{i+1},1) \cdots (\ell_n,1) \right) \mid g \in G_{\ell_i} \right\}$$
for some $1 \leq i \leq n$. It is clear that $\mathrm{stab}(C)$ is equal to
$$ \left\{ \Delta \cdot \epsilon \left( (\ell_1,1) \cdots ( \ell_{i-1},1) (\ell_i,g) (\ell_{i+1},1) \cdots (\ell_n,1) \right) \cdot \Delta^{-1} \mid g \in G_{\ell_i} \right\},$$
so that the action $\mathrm{stab}(C) \curvearrowright C$ is free and transitive on the vertices. Let $S \in \mathrm{stab}(J)$. According to Corollary \ref{cor:XPstabhyp}, 
$$S = \Delta \cdot (A+ \epsilon(\ell_i,g) +C ) \cdot \Delta^{-1}$$
for some $A \in D(\mathcal{P}, \mathcal{G}, \ell_1 \cdots \ell_{i-1})$, $C \in D(\mathcal{P}, \mathcal{G}, \ell_{i+1} \cdots \ell_n)$ and $g \in G_{\ell_i}$. On the other hand, it follows from Proposition \ref{prop:linearhyp} that the set $\mathcal{S}(J)$ of the fibers of $J$ can be naturally identified to $G_{\ell_i}$, so that the action of $S$ on $\mathcal{S}(J)$ corresponds to the action of $g$ on $G_{\ell_i}$ by left-multiplication. Thus, the element
$$\Delta \cdot \epsilon((\ell_1,1) \cdots (\ell_{i-1},1) (\ell_i,g) (\ell_{i+1},1) \cdots (\ell_n,1) \cdot \Delta^{-1}$$
of the stabiliser of $C$ induces the same permutation on $\mathcal{S}(J)$ as $S$.

\medskip \noindent
Thus, we have proved that the action $D(\mathcal{P}, \mathcal{G},w) \curvearrowright X(\mathcal{P}, \mathcal{G},w)$ is topical-transitive. 
\end{proof}

\noindent
Our next preliminary lemma sum up a few easy observations.

\begin{lemma}\label{lem:DPmisc}
The following statements hold.
\begin{itemize}
	\item[(i)] If $\mathcal{P}$ is a finite presentation, then any vertex of $X(\mathcal{P}, \mathcal{G})$ belongs to only finitely many cliques.
	\item[(ii)] Vertex-stabilisers of $X(\mathcal{P}, \mathcal{G})$ are trivial.
	\item[(iii)] If $\mathcal{P}$ is a finite presentation and if the class $[w]_{\mathcal{P}}$ of $w$ modulo $\mathcal{P}$ is finite, then $X(\mathcal{P}, \mathcal{G},w)$ contains finitely many $D(\mathcal{P}, \mathcal{G},w)$-orbits of cliques.
	\item[(iv)] $X(\mathcal{P}, \mathcal{G})$ does not contain an increasing sequence of prisms. Moreover, if $\mathcal{P}$ is a finite presentation and if the class $[w]_{\mathcal{P}}$ is finite, then the cubical dimension of $X(\mathcal{P}, \mathcal{G},w)$ is finite.
	\item[(v)] Suppose that $[w]_{\mathcal{P}}$ is finite. For every prism $P=C_1 \times \cdots \times C_n$, $\mathrm{stab}(P)= \mathrm{stab}(C_1) \times \cdots \times \mathrm{stab}(C_n)$.
\end{itemize}
\end{lemma}

\begin{proof}
We begin by proving Point $(i)$. Let $\Delta \in X(\mathcal{P}, \mathcal{G})$. Clearly, if $\mathrm{bot}^-(\Delta)$ has length $n$, then there exist at most $n$ cosets containing $\Delta$. Moreover, because $\mathcal{P}$ has finitely many relations, necessarily there exist finitely many atomic diagrams $A$ such that the concatenation $\Delta \circ A$ is well-defined. It follows from Lemma \ref{lem:XPclique} that $\Delta$ belongs to finitely many cliques of $X(\mathcal{P}, \mathcal{G})$.

\medskip \noindent
Point $(ii)$ is clear since $X(\mathcal{P}, \mathcal{G})$ is a Cayley graph.

\medskip \noindent
Two vertices $\Delta_1,\Delta_2 \in X(\mathcal{P}, \mathcal{G},w)$ belong to the same $D(\mathcal{P}, \mathcal{G},w)$-orbit if and only if $\mathrm{bot}^-(\Delta_1)= \mathrm{bot}^-(\Delta_2)$. Since we suppose that $[w]_{\mathcal{P}}$ is finite, it follows that $X(\mathcal{P}, \mathcal{G},w)$ contains only finitely many $D(\mathcal{P}, \mathcal{G},w)$-orbits of vertices. On the other hand, we know from Point $(i)$ that a vertex belongs to finitely many cliques, so that Point $(iii)$ follows.

\medskip \noindent
Let $P_1,P_2, \ldots$ be a nondecreasing sequence of prisms, and set $P= \bigcup\limits_{i \geq 1} P_i$. Notice that $P$ is gated since prisms are gated themselves. Let $\Delta$ denote the projection of $\epsilon(w)$ onto $P$, and fix some $j \geq 1$ such that $\Delta \in P_j$. It follows from Proposition \ref{prop:XPprism} that the cubical dimension of any prism containing $\Delta$ is bounded above by the length of $\mathrm{bot}^-(\Delta)$. So $P_1,P_2 , \ldots$ is a nondecreasing sequence of prisms of uniformly bounded cubical dimensions. It follows that the sequence must be eventually constant. This proves the first assertion of Point $(iv)$. The second assertion follows from Point $(i)$ and from the observation that $X(\mathcal{P}, \mathcal{G},w)$ contains finitely $D(\mathcal{P}, \mathcal{G},w)$-orbits of prisms according to Point $(iii)$.

\medskip \noindent
Finally, we want to prove Point $(v)$. Let $P$ be a prism. Up to translating, we can suppose without loss of generality that $\epsilon(w) \in P$. Let
$$P= \{ L_1(g_1) \cdots L_k(g_k) \cdot B \mid g_i \in G_{y_i}, \ B \leq A \}$$
be the description given by Proposition \ref{prop:XPprism}. As a consequence, if $S \in \mathrm{stab}(P)$, then $S = S \cdot \epsilon(w) \in P$ implies that $S= L_1(g_1) \cdots L_k(g_k) \cdot B$ for some elements $g_i \in G_{y_i}$ and some prefix $B \leq A$. A fortiori, $S$ is a thin diagram. On the other hand, $S$ is also a spherical diagram, so we deduce from Claim \ref{claim:thinlinear} below that $S$ must be linear, ie., $B$ must be trivial. Thus, if we write $w=p_1y_1 \cdots p_ky_kp_{k+1}$ where $p_1, \ldots, p_{k+1} \in \Sigma^+$, we have shown that 
$$\mathrm{stab}(P) \subset \{\epsilon(p_1)+ \epsilon((y_1,g_1))+ \cdots + \epsilon(p_k) + \epsilon((y_k,g_k))+ \epsilon(p_{k+1}) \mid g_i \in G_{y_i} \}.$$
The reverse inclusion being clear, this inclusion turns out to be an equality. By noticing that
$$\{\epsilon(p_1y_1 \cdots p_{i-1} y_{i-1}p_i)+ \epsilon((y_i,g))+ \epsilon(p_{i+1}y_{i+1} \cdots p_ky_kp_{k+1}) \mid g \in G_{y_i} \}$$
is the stabiliser of a clique of $P$, we deduce that, if we write our prism $P$ as a product of cliques $C_1 \times \cdots C_n$ such that $C_1, \ldots, C_k$ are the cliques which are cosets and $C_{k+1}, \ldots, C_n$ the cliques labelled by atomic diagrams, then $\mathrm{stab}(P) = \mathrm{stab}(C_1) \times \cdots \times \mathrm{stab}(C_k)$. Since it follows from Lemma \ref{lem:DPatomichyp} and from Point $(ii)$ that $\mathrm{stab}(C_i)$ is trivial for $k+1 \leq i \leq n$, we conclude that $\mathrm{stab}(P)= \mathrm{stab}(C_1) \times \cdots \times \mathrm{stab}(C_n)$.
\end{proof}

\begin{claim}\label{claim:thinlinear}
Suppose $[w]_{\mathcal{P}}$ finite. For any diagrams $A_1, \ldots, A_n \in X(\mathcal{P}, \mathcal{G},w)$, if $A_1 + \cdots + A_n$ is spherical then every $A_i$ must be spherical. 
\end{claim}

\begin{proof}
It is sufficient to prove our claim for $n=2$, the general case following by induction. So let $A,B \in X (\mathcal{P}, \mathcal{G},w)$ be two diagrams such that $A+B$ is spherical. Set $\mathrm{top}^-(A)=a$, $\mathrm{top}^-(B)=b$, $\mathrm{bot}^-(A)=p$ and $\mathrm{bot}^-(B)=q$. So
$$ab= \mathrm{top}^-(A+B)= \mathrm{bot}^-(A+B) = pq$$
holds in $\Sigma^+$. Suppose that $|a| \neq |p|$. Up to considering $(A+B)^{-1}=A^{-1}+B^{-1}$ instead of $A+B$, we can suppose without loss of generality that $|a|< |p|$. So the equality $ab=pq$ implies that $a$ is a proper prefix of $p$, ie., $p=ar$ for some non empty word $r \in \Sigma^+$. On the other, modulo $\mathcal{P}$, we know that $a=\mathrm{top}^-(A)=\mathrm{bot}^-(A)=p$, so the equality $a=ar$ holds modulo $\mathcal{P}$. Therefore, $ar^nb \in [w]_{\mathcal{P}}$ for every $n \geq 0$, so that $[w]_{\mathcal{P}}$ must be infinite since $r$ is non empty. Thus, we have proved that the assumption that $[w]_{\mathcal{P}}$ is finite implies that necessarily $|a|=|p|$. So it follows from the equality $ab=pq$ in $\Sigma^+$ that $a=q$, and next that $b=q$. This precisely means that $A$ and $B$ are spherical.
\end{proof}

\noindent
Now, we are ready to apply the criteria proved in Sections \ref{section:topicalactionsI} and \ref{section:topicalactionsII}. First, we deduce from Proposition \ref{prop:properlydiscontinuous}:

\begin{thm}\label{thm:DPpropercc}
Let $\mathcal{P}= \langle \Sigma \mid \mathcal{R} \rangle$ be a semigroup presentation, $\mathcal{G}$ is collection of groups indexed by $\Sigma$ and $w \in \Sigma^+$ a base word. If the groups of $\mathcal{G}$ act properly on CAT(0) cube complexes, then so does the diagram product $D(\mathcal{P}, \mathcal{G},w)$. 
\end{thm}

\noindent
Next, it follows from Propositions \ref{prop:CAT0metricallyproper}, \ref{aTmenablegroups} and \ref{aBmenablegroups}:

\begin{thm}\label{thm:DPmetricallyproper}
Let $\mathcal{P}= \langle \Sigma \mid \mathcal{R} \rangle$ be a semigroup presentation, $\mathcal{G}$ is collection of groups indexed by $\Sigma$ and $w \in \Sigma^+$ a base word. Suppose that $\mathcal{P}$ is a finite presentation.
\begin{itemize}
	\item If the groups of $\mathcal{G}$ act metrically properly on CAT(0) cube complexes, then so does the diagram product $D(\mathcal{P}, \mathcal{G},w)$;
	\item if the groups of $\mathcal{G}$ are a-T-menable, then so is the diagram product $D(\mathcal{P}, \mathcal{G},w)$;
	\item if the groups of $\mathcal{G}$ are a-$L^p$-menable for some $p \notin 2 \mathbb{Z}$, then so is the diagram product $D(\mathcal{P}, \mathcal{G},w)$.
\end{itemize} 
\end{thm}

\noindent
Finally, we deduce from Proposition \ref{prop:cubulatinggeometrically} and Theorem \ref{thm:producingCAT0groups}:

\begin{thm}\label{thm:DPCAT0}
Let $\mathcal{P}= \langle \Sigma \mid \mathcal{R} \rangle$ be a semigroup presentation, $\mathcal{G}$ is collection of groups indexed by $\Sigma$ and $w \in \Sigma^+$ a base word. Suppose that $\mathcal{P}$ is a finite presentation and that $[w]_{\mathcal{P}}$ is finite. 
\begin{itemize}
	\item If the groups of $\mathcal{G}$ act geometrically on CAT(0) cube complexes, then so does the diagram product $D(\mathcal{P}, \mathcal{G},w)$;
	\item if the groups of $\mathcal{G}$ are CAT(0), then so is the diagram product $D(\mathcal{P}, \mathcal{G},w)$.
\end{itemize} 
\end{thm}

\noindent
We conclude this section by estimating the equivariant $\ell^p$-compressions of some diagram products. 

\begin{thm}\label{thm:DPcompressions}
Let $\mathcal{P}= \langle \Sigma \mid \mathcal{R} \rangle$ be a semigroup presentation, $\mathcal{G}$ is collection of groups indexed by $\Sigma$ and $w \in \Sigma^+$ a base word. Suppose that $\mathcal{P}$ is a finite presentation, that $[w]_{\mathcal{P}}$ is finite, and that the groups of $\mathcal{G}$ are finitely generated. For every $p \geq 1$, the diagram product $D(\mathcal{P}, \mathcal{G},w)$ satifies
$$\alpha_p^*(D(\mathcal{P}, \mathcal{G},w)) \geq \min \left( \frac{1}{p}, \min\limits_{G \in \mathcal{G}} \alpha_p^*(G) \right).$$
\end{thm}

\begin{proof}
It follows from \cite[Proposition 11.2.4]{automaticgroups} that the canonical map $D(\mathcal{P}, \mathcal{G},w) \hookrightarrow X(\mathcal{P},\mathcal{G},w)$ (ie., the orbit map associated to the basepoint $\epsilon(w)$) is a quasi-isometric embedding. Therefore, the conclusion follows from Proposition \ref{prop:equicompression}. 
\end{proof}

\begin{remark}
It is worth noticing that our inequality turns out to be an equality if $p \geq 2$ and if $D(\mathcal{P}, \mathcal{G},w)$ contains a quasi-isometrically embedded non abelian free subgroup; see Remark \ref{rem:compressionGP}. 
\end{remark}

\subsection{When is a diagram product hyperbolic?}

\noindent
In this section, our goal is to determine precisely when a diagram product is hyperbolic, under the assumption that the class of our base word is finite (see the discussion relative to Question \ref{question:twistedGPhyp}). Our criterion is the following:

\begin{thm}\label{thm:DPhyp}
Let $\mathcal{P}= \langle \Sigma \mid \mathcal{R} \rangle$ be a semigroup presentation, $\mathcal{G}$ a collection of finitely generated groups indexed by $\Sigma$, and $w \in \Sigma^+$ a base word whose class $[w]_{\mathcal{P}}$ is finite. The diagram product $D(\mathcal{P}, \mathcal{G},w)$ is hyperbolic if and only if, for every non empty words $u,v \in \Sigma^+$ such that $w$ is equal to $uv$ modulo $\mathcal{P}$, at least one of the two diagram products $D(\mathcal{P}, \mathcal{G}, u)$ and $D(\mathcal{P}, \mathcal{G}, v)$ is finite. 
\end{thm}

\noindent
We want to apply Theorem \ref{thm:qmrelativelyhyp}. First of all, we need to understand when our quasi-median graph $X(\mathcal{P}, \mathcal{G},w)$ is hyperbolic.

\begin{prop}\label{prop:XPhyp}
Suppose that $[w]_{\mathcal{P}}$ is finite. The graph $X(\mathcal{P}, \mathcal{G},w)$ is hyperbolic if and only if, for every non empty words $u,v \in \Sigma^+$ such that $w$ is equal to $uv$ modulo $\mathcal{P}$, at least one of the graphs $D(\mathcal{P}, \mathcal{G}, u)$ and $D(\mathcal{P}, \mathcal{G}, v)$ is bounded. 
\end{prop}

\begin{proof}
Suppose that there exist two words $u,v \in \Sigma^+$ such that $w$ is equal to $uv$ modulo $\mathcal{P}$ and such that both $X(\mathcal{P}, \mathcal{G}, u)$ and $X(\mathcal{P}, \mathcal{G}, v)$ are unbounded. Then the map
$$\left\{ \begin{array}{ccc} X(\mathcal{P},\mathcal{G},u) \times X(\mathcal{P}, \mathcal{G},v) & \to & X(\mathcal{P}, \mathcal{G},v) \\ (\Delta_1, \Delta_2) & \mapsto & \Delta \cdot \left( \Delta_1 + \Delta_2 \right) \cdot \Delta^{-1} \end{array} \right.,$$
where $\Delta$ is a fixed diagram satisfying $\mathrm{top}^-(\Delta)=w$ and $\mathrm{bot}^-(\Delta)=uv$, isometrically embeds a product of two unbounded graphs into $X(\mathcal{P}, \mathcal{G},w)$, so that $X(\mathcal{P}, \mathcal{G},w)$ cannot be hyperbolic.

\medskip \noindent
Conversely, suppose that, for every words $u,v \in \Sigma^+$ such that $w$ is equal to $uv$ modulo $\mathcal{P}$, at most one of the two diagram products $D(\mathcal{P}, \mathcal{G}, u)$ and $D(\mathcal{P}, \mathcal{G}, v)$ is finite. We want to apply Proposition \ref{prop:qmhyp}. So let $R : [0,n] \times [0,m] \hookrightarrow X(\mathcal{P}, \mathcal{G},w)$ be a flat rectangle. Up to translating, we can suppose that $(0,0)= \epsilon(w)$. Let $A$, $B$ and $C$ denote the vertices $(0,m)$, $(n,0)$ and $(n,m)$ respectively. Because $A$ and $B$ belong to some geodesics between $\epsilon(w)$ and $C$, we deduce from Lemma \ref{lem:XPgeodesic} that $A$ and $B$ are prefixes of $C$. Notice that, if $D$ is prefix of both $A$ and $B$, then $D$ belongs to some geodesics between $\epsilon(w)$ and $A$, and between $\epsilon(w)$ and $B$, according to Lemma \ref{lem:XPgeodesic}, so that
$$\begin{array}{lcl} 0 & = & m+n-(m+n) = d(\epsilon(w),A) + d(\epsilon(w),B)- d(A,B) \\ \\ & = & 2d(\epsilon(w),D) + d(D,A)+d(D,B)-d(A,B) \geq 2 \cdot \# D \end{array}$$
which implies that $D$ is trivial. Therefore, $A$ and $B$ do not have a non trivial common prefix. Because $\# C =m+n= \# A+ \# B$, it follows that the prefixes $A$ and $B$ cover completely $C$. Next, we claim that the supports of $A$ and $B$ are disjoint. Indeed, if there exists an edge $e$ which belongs to the supports of $A$ and $B$, four cases may happen:
\begin{itemize}
	\item $e$ is included into the top paths of cells of $A$ and $B$. Since there exists at most one cell below $e$, we deduce that $A$ and $B$ share a cell, which is impossible since we know that $A$ and $B$ have no non trivial common prefix.
	\item $e$ is labelled by a letter of $\Sigma(\mathcal{G})$ with a non trivial second coordinate in both $A$ and $B$. Once again, this is impossible since we know that $A$ and $B$ have no non trivial common prefix.
	\item $e$ is included into the top path of some cell of $B$ and is labelled by a letter of $\Sigma(\mathcal{G})$ with a second coordinate which is non trivial in $A$ and trivial in $B$. As a consequence of Lemma \ref{lem:XPgeodesic}, it is possible to construct $C$ from $B$ by right-concatenating unitary diagrams. Because $e$ is included into the top path of some cell, these concatenations cannot modify the label which has $e$ in $B$, so that $e$ must be labelled in $C$ by a letter of $\Sigma(\mathcal{G})$ which has a trivial second coordinate. On the other hand, we know that this coordinate in not trivial in $A$, so that it cannot be trivial in $C$ since $A$ is a prefix of $C$. 
	\item $e$ is included into the top path of some cell of $A$ and is labelled by a letter of $\Sigma(\mathcal{G})$ with a second coordinate which is non trivial in $B$ and trivial in $A$. The situation is symmetric to the previous one.
\end{itemize}
Thus, we have proved that the supports of $A$ and $B$ are disjoint. We conclude that there exist words $u_1,v_1, \ldots, u_k,v_k \in \Sigma(\mathcal{G})^+$, a $(u_i, \ast)$-diagram $U_i$ and a $(v_i,\ast)$-diagram $V_i$ for every $1 \leq i \leq k$, such that $w=u_1v_1 \cdots u_kv_k$ in $\Sigma(\mathcal{G})^+$, $C=U_1+V_1+ \cdots + U_k+V_k$, $B=U_1+ \epsilon(v_1)+ \cdots + U_k+ \epsilon(v_k)$, and $C=\epsilon(u_1)+ V_1 + \cdots + \epsilon(u_k)+ V_k$. 

\medskip \noindent
Notice that, because $[w]_{\mathcal{P}}$ is finite, $k$ is bounded above by the maximal length $L$ of a word in $[w]_{\mathcal{P}}$. Set also
$$K = \max \left\{ \mathrm{diam}~X(\mathcal{P}, \mathcal{G},u) \mid u \in \mathcal{V}_f \right\}< + \infty,$$
where $\mathcal{V}_f$ denotes the collection of all the subwords $u$ of words of $[w]_{\mathcal{P}}$ such that $X(\mathcal{P}, \mathcal{G},u)$ is bounded. 
Next, by assumption, we know that at most one word $m \in \{ u_1,v_1, \ldots, u_k,v_k \}$ is such that $X(\mathcal{P}, \mathcal{G},m)$ is unbounded. If such a word does not exist, then
$$n = \#B = \# U_1+ \cdots + \# U_k  \leq k \cdot K \leq L \cdot K,$$ 
and similarly $m \leq KL$. Now, if such a word exists, up to switching $B$ and $C$ we can suppose without loss of generality that $\ell=v_i$ for some $1 \leq i \leq k$. So we know that $X(\mathcal{P}, \mathcal{G},u_j)$ has diameter at most $K$ for every $1 \leq j \leq k$, so that $n \leq L \cdot K$.

\medskip \noindent
Thus, we have proved that any flat rectangle of $X(\mathcal{P}, \mathcal{G},w)$ is $KL$-thin, so that $X(\mathcal{P}, \mathcal{G},w)$ must be hyperbolic according to Proposition \ref{prop:qmhyp}. 
\end{proof}

\begin{figure}
\begin{center}
\includegraphics[scale=0.58]{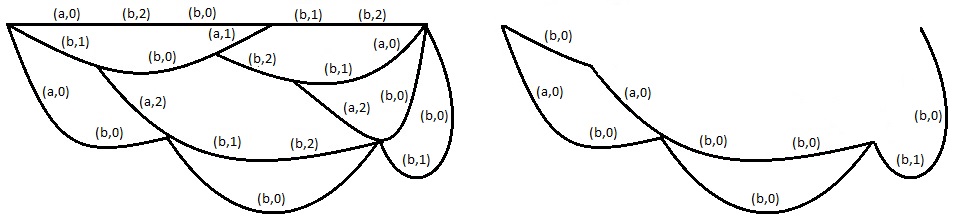}
\end{center}
\label{figure20}
\caption{A diagram and its maximal thin suffix.}
\end{figure}

\noindent
Next, we need to understand further the combinatorics of linear hyperplanes of $X(\mathcal{P}, \mathcal{G})$. This is done by Proposition \ref{prop:hdiaghyp} below, but first we need to introduce preliminary definitions. Our next lemma defines the \emph{maximal thin suffix} of a diagram; see \cite{MR1978047} for the definition in the context of diagram groups, and see Figure \ref{figure20} for an example based on the semigroup presentation
$$\mathcal{P}= \langle a,b \mid ab^2=b^2a, ab=ba, b=b^2 \rangle$$
and the groups $G_a=G_b= \mathbb{Z}$.

\begin{lemma}\label{lem:maxthinsuffix}
Let $\Delta$ be a reduced diagram. There exists a unique thin diagram $T$ such that $\Delta$ decomposes as an absolutely reduced concatenation $\Delta_1 \circ T$, and such that for any other decomposition of $\Delta$ as an absolutely reduced concatenation $\Delta_2 \circ T'$, where $T'$ is a thin diagram, necessarily $\Delta_1 \leq \Delta_2$. 
\end{lemma}

\begin{proof}
In this proof, we work in the interval $I$ between $\epsilon(w)$ and $\Delta$. According to Proposition \ref{prop:intervalmedian}, this is a median graph. Let $\mathcal{H}$ denote the set of all the hyperplanes separating $\epsilon(w)$ and $\Delta$ which contain $\Delta$ in their neighborhoods. Notice that $\mathcal{H}$ is necessarily finite since only finite many hyperplanes separate $\epsilon(w)$ and $\Delta$. Moreover, the hyperplanes of $\mathcal{H}$ are pairwise transverse. Indeed, if there exist two disjoint hyperplanes $J_1,J_2 \in \mathcal{H}$, then cutting $I$ along $J_1$ and $J_2$ produces three connected components, say $A,B,C$ where $B$ is the connected component between $J_1$ and $J_2$, but then $\Delta$ must belong to $B$ because it belongs to the neighborhoods of $J_1$ and $J_2$, so that either $J_1$ or $J_2$ does not separate $\epsilon(w)$ and $\Delta$, depending on whether $\epsilon(w) \in A$ or $\epsilon(w) \in C$, a contradiction. As a consequence, it follows from Proposition \ref{prop:transversehypcube} that there exists a cube $P$ containing $\Delta$ whose dual hyperplanes are precisely the hyperplanes of $\mathcal{H}$. Let $\Delta_1$ denote the projection of $\epsilon(w)$ onto $P$, and set $T= \Delta_1^{-1} \cdot \Delta$. We claim that $T$ is the diagram we are looking for.

\medskip \noindent
First of all, notice that the concatenation $\Delta = \Delta_1 \circ T$ is absolutely reduced as a consequence of Lemma \ref{lem:XPgeodesic}, since there exists a geodesic between $\epsilon(w)$ and $\Delta$ passing through $\Delta_1$.

\medskip \noindent
Next, we claim that $T$ is thin. Notice that, because all the hyperplanes of $\mathcal{H}$ separate $\epsilon(w)$ and $\Delta$, they must also separate $\Delta_1$ and $\Delta$, so that $\Delta_1$ and $\Delta$ are diametrically opposite in $P$. Therefore, our claim follows from

\begin{claim}
Let $P$ be a cube of $X(\mathcal{P},\mathcal{G})$ containing $\epsilon(w)$, and let $A$ denote the vertex of $P$ diametrically opposite to $\epsilon(w)$. Then $A$ is a thin diagram.
\end{claim}

\noindent
Let $A_1, \ldots, A_n$ denote the neighbors of $\epsilon(w)$ in $P$. Notice that they are unitary diagrams since $\# A_i = d(\epsilon(w),A_i)=1$ for every $1 \leq i \leq n$. Moreover, for every $1 \leq i < j \leq n$, the vertices $\epsilon(w)$, $A_i$ and $A_j$ generate a square, so it follows from Lemma \ref{lem:diagproductsquare} that the supports of $A_i$ and $A_j$ are disjoint in $w$. Consequently, we can write $w=x_1a_1 \cdots x_na_nx_{n+1}$ so that (up to permuting $A_1, \ldots, A_n$) the diagram $A_i$ can be written as $\epsilon(x_1a_1 \cdots x_{i-1}a_{i-1}x_i) + \Pi_i + \epsilon(x_{i+1}a_{i+1} \cdots x_na_n)$ where $\Pi_i$ is either a single cell or an edge. Because
$$\{ \epsilon(x_1)+ B_i + \cdots + \epsilon(x_n) + B_n + \epsilon(x_{n+1}) \mid B_i=\Pi_i \ \text{or} \ \epsilon(a_i) \}$$
is a cube generating by the edges $(\epsilon(w), A_1), \ldots, (\epsilon(w),A_n)$, and that a collection of edges generates at most one cube (this is a consequence of the fact that a quasi-median graph does not contain induced subgraphs isomorphic to $K_{2,3}$), it follows that the cube given above is $P$. Therefore, 
$$A= \epsilon(x_1)+ A_1+ \cdots + \epsilon(x_n)+ A_n + \epsilon(x_{n+1}).$$
This concludes the proof of our claim.

\medskip \noindent
Finally, decompose $\Delta$ as an absolutely reduced concatenation $\Delta_2 \circ T'$ where $T'$ is thin. The vertices $\Delta_2$ and $\Delta$ belong to a common cube. Indeed, writting $T'$ as a sum $U_1+ \cdots +U_k$ where $U_i$ is a unitary diagram for every $1 \leq i \leq k$, and setting $\mathrm{top}^-(U_i)=u_i$ for every $1 \leq i \leq k$, then
$$\{ \Delta_2 \circ (V_1+ \cdots + V_k) \mid V_i = U_i \ \text{or} \ \epsilon(u_i) \}$$
defines such a cube. As a consequence, the hyperplanes separating $\Delta_2$ and $\Delta$ define a subcollection $\mathcal{H}' \subset \mathcal{H}$. On the other hand, $P$ contains a vertex such the set of the hyperplanes separating it from $\Delta$ is precisely $\mathcal{H}'$; this vertex must be $\Delta_2$. So we have proved that $\Delta_2 \in P$. A fortiori, there exists a geodesic between $\epsilon(w)$ and $\Delta$ passing through $\Delta_1$ and $\Delta_2$, so that $\Delta_1$ must be a prefix of $\Delta_2$ according to Lemma \ref{lem:XPgeodesic}. 

\medskip \noindent
Thus, we have proved that $T$ is a diagram satisfying the condition we are interested in. Suppose that there exists another such diagram $T'$, so that $\Delta$ can be written as an absolutely reduced concatenation $\Delta_2 \circ T'$. By applying twice the condition satisfied by $T$ and $T'$, we find that $\Delta_1 \leq \Delta_2$ and $\Delta_2 \leq \Delta_1$, hence $\Delta_1= \Delta_2$. It follows that
$$T = \Delta_1^{-1} \cdot \Delta = \Delta_2^{-1} \cdot \Delta = T',$$
which concludes the proof.
\end{proof}

\noindent
Following \cite{arXiv:1505.02053}, we want to encode the linear hyperplanes of $X(\mathcal{P}, \mathcal{G})$ by using specific semigroup diagrams.

\begin{definition}
An \emph{h-diagram} $(\Delta,e)$ is the couple of a semigroup diagram $\Delta$ and an edge $e$ of the bottom path of $\Delta$ such that the maximal thin suffix of $\Delta$ is atomic, such that $e$ belongs to the bottom path of its unique cell, and such that the group corresponding the letter labelling $e$ is non trivial. Two h-diagrams $(\Delta_1,e_1)$ and $(\Delta_2,e_2)$ are \emph{transverse} if there exists a third diagram $\Delta_3$ containing $\Delta_1$ and $\Delta_2$ as prefixes such that the edges $e_1$ and $e_2$ are distinct and both belong to the bottom path of $\Delta_3$.
\end{definition}

\noindent
Given an h-diagram $(\Delta,e)$, let $J(\Delta,e)$ denote the hyperplane dual to the clique
$$\left\{ \Delta \cdot \left( \epsilon(a)+ \epsilon(b,g) + \epsilon(c) \right) \mid g \in G_b \right\},$$
where $a,c \in \Sigma^+$ are words and $b \in \Sigma$ a letter such that $\mathrm{bot}^-(\Delta)=abc$ where $b$ labels the edge $e$.

\begin{prop}\label{prop:hdiaghyp}
The map $(\Delta, e) \mapsto J(\Delta,e)$ defines a bijection from the set of h-diagrams onto the set of linear hyperplanes of $X(\mathcal{P}, \mathcal{G})$. Moreover, two h-diagrams $(\Delta_1,e_1)$ and $(\Delta_2,e_2)$ are transverse if and only if their hyperplanes $J(\Delta_1,e_1)$ and $J(\Delta_2,e_2)$ are transverse.
\end{prop}

\begin{proof}
Let $J$ be a linear hyperplane and let $\Delta$ denote the projection of $\epsilon(w)$ onto $N(J)$. As a consequence of Proposition \ref{prop:linearhyp}, there exist letters $a_1, \ldots, a_r,b,c_1, \ldots, c_s \in \Sigma$ such that $\mathrm{bot}(\Delta)= \epsilon(abc)$ if we write $a=a_1 \cdots a_r$ and $c=c_1 \cdots c_s$, and such that $N(J)$ is the subgraph generated by the vertices
$$\left\{ \Delta \circ \left( A+ \epsilon((b,g))+C \right) \left| A \in X(\mathcal{P},\mathcal{G},a_1 \cdots a_r), C \in X(\mathcal{P}, \mathcal{G},c_1 \cdots c_s), g \in G_b \right. \right\}.$$
In particular, the maximal thin suffix of $\Delta$ must be a sum of atomic diagrams. If $e$ denotes the edge of $\mathrm{bot}(\Delta)$ which is labelled by $b$, suppose by contradiction that there exists a cell in the maximal thin suffix of $\Delta$ which does not contain $e$ in its bottom path. So we can write $\Delta = \Delta_0 \circ (E+ \epsilon((b,1))+ F)$ for some $(\ast, (a_1,1) \cdots (a_r,1))$-diagram $E$ and some $(\ast, (c_1,1) \cdots (c_s,1))$-diagram $F$, such that $\# E + \# F \geq 1$. By noticing that $\Delta \cdot (E^{-1}+ \epsilon((b,1))+ F^{-1}) = \Delta_0$ belongs to $N(J)$ and that $\# \Delta_0 < \# \Delta$, we find a contradiction with the fact that $\Delta$ minimizes the distance to $\epsilon(w)$ in $N(J)$. Thus, the maximal thin suffix of $\Delta$ contains a unique cell, and this cell contains $e$ in its bottom path. This proves that $(\Delta,e)$ is an h-diagram, so that $J=J(\Delta,e)$ since $J$ and $J(\Delta,e)$ are both dual to the clique
$$ C = \left\{ \Delta \circ \left( \epsilon(a)+ \epsilon((b,g))+ \epsilon(c) \right) \mid g \in G_b \right\}.$$
This proves that our map is surjective. Next, for any clique $C'$ dual to $J$ different to $C$, say
$$C'=\left\{ \Delta \circ \left( A+ \epsilon((b,g))+ C \right) \mid g \in G_b \right\}$$
for some $A \in X(\mathcal{P},\mathcal{G},a_1 \cdots a_r)$ and $C \in X(\mathcal{P}, \mathcal{G},c_1 \cdots c_s)$ such that $\#A+ \# C \geq 1$, notice that $\left( \Delta \circ ( A+ \epsilon((b,g))+ C), e \right)$ is not an h-diagram. It follows that our map is injective. So our map turns out to be bijective, but it is worth noticing that we have also proved the following observation:

\begin{fact}\label{fact:hdiagproj}
If $J= J(\Delta,e)$, then $\Delta$ is the projection of $\epsilon(w)$ onto $N(J)$.
\end{fact}

\noindent
Let $(\Delta_1,e_1)$ and $(\Delta_2,e_2)$ be two h-diagrams. Suppose that they are transverse. So there exists a third diagram $\Delta$ containing $\Delta_1$ and $\Delta_2$ as prefixes such that the edges $e_1$ and $e_2$ are distinct and both belong to the bottom path of $\Delta$. We can write $\mathrm{bot}^-(\Delta)= pqrst$ for some words $p,r,t \in \Sigma^+$ and some letters $q,s \in \Sigma$, such that $q$ and $s$ label the edges $e_1$ and $e_2$; up to switching $\Delta_1$ and $\Delta_2$, say that $q$ and $s$ label respectively $e_1$ and $e_2$. We first claim that, fixing some $g \in G_q$ which is different from the inverse of the second coordinate of the letter of $\Sigma(\mathcal{G})^+$ labelling the edge $e_1$ in $\Delta$, the edge between $\Delta$ and $\Delta \cdot (\epsilon(p)+ \epsilon((q,g))+ \epsilon(rst))$ is dual to the hyperplane $J(\Delta_1,e_1)$. Indeed, since $\Delta_1$ is a prefix of $\Delta$ containing the edge $e_1$ in its bottom path, $\Delta$ decomposes as $\Delta_1 \circ (A+ \epsilon(q,h)+B)$ for some diagrams $A,B$ and some element $h \in G_q$. Our claim follows from Proposition \ref{prop:linearhyp}. Similarly, fixing some $k \in G_s$ which is different from the inverse of the second coordinate of the letter of $\Sigma(\mathcal{G})^+$ labelling the edge $e_2$ in $\Delta$, the edge between $\Delta$ and $\Delta \cdot (\epsilon(pqr)+ \epsilon((s,k))+ \epsilon(t))$ must be dual to the hyperplane $J(\Delta_2,e_2)$. On the other hand, the vertices $\Delta$, $\Delta \cdot (\epsilon(p)+ \epsilon((q,g))+ \epsilon(rst))$, $\Delta \cdot (\epsilon(pqr)+ \epsilon((s,k))+ \epsilon(t))$ and $\Delta \cdot ( \epsilon(p) + \epsilon((q,g))+ \epsilon(r) + \epsilon((s,k))+ \epsilon(t))$ defines a square, so that the hyperplanes $J(\Delta_1,e_1)$ and $J(\Delta_2,e_2)$ have to be transverse.

\medskip \noindent
Conversely, suppose that $J(\Delta_1,e_1)$ and $J(\Delta_2,e_2)$ are transverse. According to Proposition \ref{prop:transversehypcube}, there exists a prism $P$ whose dual hyperplanes are precisely $J(\Delta_1,e_1)$ and $J(\Delta_2,e_2)$. Let $\Delta$ denote the projection of $\epsilon(w)$ onto $P$. We know from Fact \ref{fact:hdiagproj} that $\Delta_1$ is the projection of $\epsilon(w)$ onto $N(J(\Delta_1,e_1))$, so that we get a description of $N(J(\Delta_1,e_1))$ from Proposition \ref{prop:linearhyp}. Since $\Delta \in N(J(\Delta_1,e_1))$, it follows that $\Delta_1$ is a prefix of $\Delta$ and $\Delta$ contains the edge $e_1$ in its bottom path. Similarly,  $\Delta_2$ is a prefix of $\Delta$ and $\Delta$ contains the edge $e_2$ in its bottom path. Notice that the clique of $P$ dual to $J(\Delta_1,e_1)$ and containing $\Delta$ corresponds to concatenating linear diagrams along the edge $e_1$ of $\mathrm{bot}(\Delta)$, and similarly, the clique of $P$ dual to $J(\Delta_2,e_2)$ and containing $\Delta$ corresponds to concatenating linear diagrams along the edge $e_2$ of $\mathrm{bot}(\Delta)$, so that $J(\Delta_1,e_1) \neq J(\Delta_2,e_2)$ implies that $e_1 \neq e_2$. This concludes the proof.
\end{proof}

\noindent
Now, we are able to apply Theorem \ref{thm:qmrelativelyhyp} in order to prove the converse of Theorem \ref{thm:DPhyp}. In fact, the conclusion provides the following criterion of relative hyperbolicity:

\begin{prop}\label{prop:DPrelativelyhyp}
Let $\mathcal{P}= \langle \Sigma \mid \mathcal{R} \rangle$ be a semigroup presentation, $\mathcal{G}$ a collection of finitely generated groups indexed by $\Sigma$, and $w \in \Sigma^+$ a base word whose class $[w]_{\mathcal{P}}$ is finite. Suppose that, for every non empty words $u,v \in \Sigma^+$ such that $w$ is equal to $uv$ modulo $\mathcal{P}$, at most one of the two diagram products $D(\mathcal{P}, \mathcal{G}, u)$ and $D(\mathcal{P}, \mathcal{G}, v)$ is finite. Then the diagram product $D(\mathcal{P}, \mathcal{G},w)$ is hyperbolic relatively to a finite collection of subgroups which are commensurable to groups of $\mathcal{G}$.
\end{prop}

\begin{proof}
Notice that the last two conditions of Theorem \ref{thm:qmrelativelyhyp} amount to say that the fibers of any hyperplane dual to an infinite clique must be finite. So we need to prove that the fibers of a linear hyperplane of $X(\mathcal{P}, \mathcal{G},w)$ dual to an infinite coset are finite. In view of Proposition \ref{prop:linearhyp}, this statement is implied by our condition combined with Claim \ref{claim:constantK} proved below. 

\medskip \noindent
Next, suppose that there exist two linear hyperplanes $J_1$ and $J_2$ dual to infinite cosets which are transverse. According to Proposition \ref{prop:hdiaghyp}, there exist two h-diagrams $(\Delta_1,e_1)$ and $(\Delta_2,e_2)$ such that $J_1=J(\Delta_1,e_1)$ and $J_2=J(\Delta_2,e_2)$, and moreover $(\Delta_1,e_1)$ and $(\Delta_2,e_2)$ are transverse, ie., there exists a third diagram $\Delta$ with $\Delta_1$ and $\Delta_2$ as prefixes such that the edges $e_1$ and $e_2$ are included into $\mathrm{bot}(\Delta)$ and distinct. So we can write $\mathrm{bot}^-(\Delta)= pqrst$ for some words $p,r,t \in \Sigma^+$ and some letters $q,s \in \Sigma$, such that $q$ and $s$ label the edges $e_1$ and $e_2$. Notice that, since $J_1$ and $J_2$ are dual to infinite cosets, necessarily $G_q$ and $G_s$ are infinite. On the other hand, the equality $pqrst=w$ holds modulo $\mathcal{P}$. Therefore, our assumptions imply that two hyperplanes dual to infinite cliques cannot be transverse.

\medskip \noindent
The other assumptions of Theorem \ref{thm:qmrelativelyhyp} are satisfied according to Lemma \ref{lem:DPmisc} and Proposition \ref{prop:hdiaghyp}. Therefore, Theorem \ref{thm:qmrelativelyhyp} applies, so that the diagram product $D(\mathcal{P}, \mathcal{G},w)$ is hyperbolic relatively to a finite collection of stabilisers of hyperplanes dual to infinite cliques, which are commensurable to clique-stabilisers, and so which are commensurable to groups of $\mathcal{G}$ as a consequence of Lemma \ref{lem:XPclique}.
\end{proof}

\begin{claim}\label{claim:constantK}
Suppose that $[w]_{\mathcal{P}}$ is finite. For every subword $u$ of every word of $[w]_{\mathcal{P}}$, if $D(\mathcal{P},\mathcal{G},u)$ is finite then so is $X(\mathcal{P},\mathcal{G},u)$. 
\end{claim}

\begin{proof}
Let $u$ be subword of some word of $[w]_{\mathcal{P}}$ such that $D(\mathcal{P}, \mathcal{G},u)$ is finite. Notice that $X(\mathcal{P},\mathcal{G},u)$ is cubically finite, because it contains only finitely many $D(\mathcal{P}, \mathcal{G},u)$-orbits of cliques according to Lemma \ref{lem:DPmisc}. On the other hand, according to Lemma \ref{lem:XPclique}, a clique of $X(\mathcal{P}, \mathcal{G},u)$ is either an edge or a coset. But a coset of $X(\mathcal{P}, \mathcal{G},u)$ must be finite since otherwise $D(\mathcal{P}, \mathcal{G},u)$ would be infinite. Therefore, $X(\mathcal{P}, \mathcal{G},u)$ contains finitely many cliques, and its cliques are all finite, so that $X(\mathcal{P}, \mathcal{G},u)$ is necessarily finite.
\end{proof}

\begin{remark}
In general, a factor of a diagram product (ie., an element of $\mathcal{G}$) does not canonically embed into the diagram product, in the sense that two non conjugated clique-stabilisers may be isomorphic to a common group of $\mathcal{G}$. So the collection of subgroups obtained in Proposition \ref{prop:DPrelativelyhyp} may contain several subgroups isomorphic to the same factor.
\end{remark}

\begin{proof}[Proof of Theorem \ref{thm:DPhyp}.]
If the conditions are all satisfied, we deduce from Proposition \ref{prop:DPrelativelyhyp} that $D(\mathcal{P}, \mathcal{G},w)$ is hyperbolic relatively to hyperbolic groups, so that our diagram product must be hyperbolic (see for instance \cite[Corollary 2.41]{OsinRelativeHyp}). Conversely, suppose that $D(\mathcal{P}, \mathcal{G},w)$ is hyperbolic. It follows from Lemma \ref{lem:QIsubgroup} below that clique-stabilisers are quasiconvex subgroups, so that the groups of $\mathcal{G}$ are necessarily hyperbolic. Next, suppose that there exists a word $u v \in [w]_{\mathcal{P}}$ such that both $D(\mathcal{P}, \mathcal{G},u)$ and $D(\mathcal{P}, \mathcal{G},v)$ are infinite. If $\Delta$ is a diagram satisfying $\mathrm{top}^-(\Delta)=w$ and $\mathrm{bot}^-(\Delta)=u v$, then 
$$\left\{ \Delta \cdot \left( U+V \right) \cdot \Delta^{-1} \mid U \in D(\mathcal{P}, \mathcal{G},u), V \in D(\mathcal{P}, \mathcal{G},v) \right\}$$ 
defines a subgroup of $D(\mathcal{P}, \mathcal{G},w)$ isomorphic to $D(\mathcal{P},\mathcal{G},u) \times D(\mathcal{P}, \mathcal{G},v)$, which is a product of two infinite subgroups. Because a hyperbolic group cannot contain such a subgroup, it follows that $D(\mathcal{P}, \mathcal{G},w)$ is not hyperbolic.
\end{proof}

\begin{lemma}\label{lem:QIsubgroup}
Let $G$ be a finitely generated group acting topically-transitively on a quasi-median graph $X$. Suppose that
\begin{itemize}
	\item vertex-stabilisers are finite and clique-stabilisers are finitely generated;
	\item any vertex of $X$ belongs to finitely many cliques;
	\item $X$ contains finitely many $G$-orbits of cliques.
\end{itemize}
Then clique-stabilisers are quasi-isometrically embedded. 
\end{lemma}

\begin{proof}
Let $\mathcal{C}$ be a collection of cliques such that any $G$-orbit of hyperplanes intersects it along a single clique, and let $\mathcal{C}= \mathcal{C}_1 \sqcup \mathcal{C}_2$ denote the associated decomposition of $\mathcal{C}$. For every $C \in \mathcal{C}_2$, let $\delta_C$ denote the discrete metric $(x,y) \mapsto \left\{ \begin{array}{cl} 1 & \text{if $x \neq y$} \\ 0  & \text{otherwise} \end{array} \right.$ on $C$. If $C \in \mathcal{C}_1$, denote by $x_0(C)$ the projection of a basepoint $x_0 \in X$ onto $C$ and transfer a word metric of $\mathrm{stab}(C)$ (with respect to a finite generating set) to a metric $\delta_C$ on $C$ via the orbit map $g \mapsto g \cdot x_0(C)$. According to Proposition \ref{prop:mtopicalI}, there exists a coherent $G$-invariant system of metrics extending $\{ (C,\delta_C) \mid C \in \mathcal{C} \}$, and, if $\delta$ denotes the corresponding global metric, then $G$ acts geometrically on $(X, \delta)$ according to Corollary \ref{cor:mwhengeometric}. On the other hand, we know that $(C,\delta_C)$ is isometrically embedded into $(X,\delta)$ and that $\mathrm{stab}(C)$ acts geometrically on $(C,\delta_C)$, so it follows that $\mathrm{stab}(C)$ is quasi-isometrically embedded subgroup of $G$. 
\end{proof}

\subsection{Structure of groups acting on quasi-median graphs}

\noindent
In this section, we introduce \emph{rotative actions} on quasi-median graphs, and we state and prove a structure result for groups admitting such actions. We will show in the next section that the actions of diagram products on their associated quasi-median graphs are rotative, which will help us to understand the structure of these products. First of all, we need to introduce some vocabulary. (Recall that rotative stabilisers of hyperplanes were defined in Section \ref{section:GPcurvegraph}.)

\begin{definition}
Let $G$ be a group acting on a quasi-median graph $X$. Let $\mathcal{J}$ denote a $G$-invariant collection of hyperplanes of $X$, and, for every $J \in \mathcal{J}$, let $\mathcal{S}(J)$ denote the set of the sectors delimited by $J$. The action $G \curvearrowright X$ is \emph{$\mathcal{J}$-rotative}\index{Rotative actions} if the action $\mathrm{stab}_{\circlearrowleft}(J) \curvearrowright \mathcal{S}(J)$ is transitive and free for every $J \in \mathcal{J}$. If $\mathcal{J}$ is the set of all the hyperplanes of $X$, a $\mathcal{J}$-rotative action will be referred to as a \emph{fully rotative} action.
\end{definition}
 
\begin{definition}
Let $X$ be a quasi-median graph, $\mathcal{J}$ a collection of hyperplanes and $x_0 \in X$ a base vertex. A subcollection $\mathcal{J}_0 \subset \mathcal{J}$ is \emph{$x_0$-peripheral} if, for every $J \in \mathcal{J}_0$, there does not exist a hyperplane of $\mathcal{J}$ separating $J$ and $x_0$. 
\end{definition}

\noindent
Our structure result is the following:

\begin{thm}\label{thm:rotativeactions}
Let $G$ be a group acting $\mathcal{J}$-rotatively on a quasi-median graph $X$. Fix a basepoint $x_0 \in X$. If $Y \subset X$ denotes the intersection of the sectors containing $x_0$ which are delimited by a hyperplane of $\mathcal{J}$, then
$$G= \mathrm{Rot}(\mathcal{J}) \rtimes \mathrm{stab}(Y), \ \text{where} \ \mathrm{Rot}(\mathcal{J})= \langle  \mathrm{stab}_{\circlearrowleft}(J), \ J \in \mathcal{J} \rangle.$$
Moreover, if $\mathcal{J}_0 \subset \mathcal{J}$ denotes the unique maximal $x_0$-peripheral subcollection of $\mathcal{J}$, then $\mathrm{Rot}(\mathcal{J})$ decomposes as a graph product $\Delta \mathcal{G}$, where $\Delta$ is the graph whose vertices are the hyperplanes of $\mathcal{J}_0$ and whose edges link two hyperplanes which are transverse, and where $\mathcal{G}= \{ \mathrm{stab}_{\circlearrowleft}(J) \mid J \in \mathcal{J}_0 \}$.
\end{thm}

\begin{proof}
Notice that, since $\mathcal{J}$ is $G$-invariant, the subgroup $\mathrm{Rot}(\mathcal{J})$ is normal. 

\medskip \noindent
Let $g \in G$. Fix some $r \in \mathrm{Rot}(\mathcal{J})$ and suppose that $rg \cdot x_0 \notin Y$. Then there exists some $J \in \mathcal{J}$ which separates $rg \cdot x_0$ from $Y$. Since the action $\mathrm{stab}_{\circlearrowleft}(J) \curvearrowright \mathcal{S}(J)$ is transitive, there exists some $s \in \mathrm{stab}_{\circlearrowleft}(J)$ which sends the sector delimited by $J$ which contains $rg \cdot x_0$ to the sector delimited by $J$ which contains $x_0$. Let $a \in N(J)$ denote the projections of $rg \cdot x_0$ onto $N(J)$ and $C$ the clique dual to $J$ containing $a$. Notice that
$$\begin{array}{lcl} d(x_0,rg \cdot x_0) & = & d(rg \cdot x_0,a)+d(a,x_0) = d(srg \cdot x_0,s \cdot a) + d(x_0,s \cdot a)+1 \\ \\ & \geq & d(x_0,srg \cdot x_0)+1 \end{array}$$
Thus, if we choose some $r \in \mathrm{Rot}(\mathcal{J})$ such that
$$d(rg \cdot x_0,g \cdot x_0)= \min \{ d(sg \cdot x_0,g \cdot x_0) \mid s \in \mathrm{Rot}(\mathcal{J}) \},$$ 
we deduce from the previous observation that $rg \cdot x_0 \in Y$. On the other hand, $G$ permutes the connected components of $X$ cutting along the hyperplanes of $\mathcal{J}$, and $Y$ is precisely the connected component which contains $x_0$, so $rg \in \mathrm{stab}(Y)$. Therefore,
$$g \in r^{-1} \cdot \mathrm{stab}(Y) \subset \mathrm{Rot}(\mathcal{J}) \cdot \mathrm{stab}(Y).$$
Thus, we have proved that $G= \mathrm{Rot}(\mathcal{J}) \cdot \mathrm{stab}(Y)$. 

\medskip \noindent
Next, we want to apply Proposition \ref{prop:pingpong} in order to prove that $\mathrm{Rot}(\mathcal{J})$ is isomorphic to the graph product $\Delta \mathcal{G}$. 

\medskip \noindent
The first point to verify is that $\mathrm{Rot}(\mathcal{J})= \langle \mathrm{stab}_{\circlearrowleft}(J), \ J \in \mathcal{J}_0 \rangle$. Let $J_1 \in \mathcal{J}$. We want to prove that there exists some $r \in \langle \mathrm{stab}_{\circlearrowleft}(J), \ J \in \mathcal{J}_0 \rangle$ such that $rJ_1 \in \mathcal{J}_0$, which is sufficient to deduce the previous equality. 

\medskip \noindent
Fix some $r \in \langle \mathrm{stab}_{\circlearrowleft}(J), \ J \in \mathcal{J}_0 \rangle$, and suppose that $rJ_1 \notin \mathcal{J}_0$. Let $y_0$ denote the projection of $x_0$ onto $N(rJ_1)$. If no hyperplane of $\mathcal{J}$ separates $x_0$ and $y_0$, then $\mathcal{J}_0 \cup \{ rJ \}$ defines a new $x_0$-peripheral subcollection of $\mathcal{J}$, contradicting the maximality of $\mathcal{J}_0$. Therefore, there exists some hyperplane $J_2 \in \mathcal{J}$ separating $x_0$ and $y_0$. As a consequence of Lemma \ref{lem:projseparate}, $J_2$ also separates $rJ_1$ and $x_0$. Notice that, if $J_2 \notin \mathcal{J}_0$, then similarly there exists a third hyperplane separating $J_2$ and $x_0$, and so on. Since there exist only finitely many hyperplanes separating $rJ_1$ and $x_0$, we can suppose without loss of generality that $J_2 \in \mathcal{J}_0$. Let $s \in \mathrm{stab}_{\circlearrowleft}(J_2)$ be an element which sends the sector delimited by $J_2$ containing $rJ_1$ to the sector delimited by $J_2$ containing $x_0$. Notice that
$$\begin{array}{lcl} d(x_0,N(rJ_1)) & \geq & d(N(rJ_1),N(J_2)) + d(x_0,N(J_2)) +1 \\ \\ & \geq & d(N(srJ_1),N(J_2)) + d(x_0,N(J_2))+1 \\ \\ & \geq & d(x_0,N(srJ_1)) +1 \end{array}$$
Therefore, if we choose $r$ so that
$$d(x_0,N(rJ_1))= \min \left\{ d(x_0,N(sJ_1)) \mid s \in \langle  \mathrm{stab}_{\circlearrowleft}(J), \ J \in \mathcal{J}_0 \rangle \right\},$$
then $rJ_1 \in \mathcal{J}_0$. This concludes the proof of our first point.

\medskip \noindent
Next, notice that it follows from Lemma \ref{lem:rotativestabcommute} that any element of $\mathrm{stab}_{\circlearrowleft}(J_1)$ commutes with any element of $\mathrm{stab}_{\circlearrowleft}(J_2)$ whenever $J_1,J_2 \in \mathcal{J}_0$ are adjacent in $\Delta$. Now, if $J \in \mathcal{J}_0$, let $X_J$ be the union of all the sectors which do not contain $x_0$. Notice that, for every $J \in \mathcal{J}_0$ and every $g \in \mathrm{stab}_{\circlearrowleft}(J) \backslash \{ 1 \}$, necessarily $g \cdot x_0 \in X_J$ since the action $\mathrm{stab}_{\circlearrowleft} (J) \curvearrowright \mathcal{S}(J)$ is free. Moreover, it follows from Lemma \ref{lem:rotativestab} that, if $J_1,J_2 \in \mathcal{J}_0$ are adjacent, then $g \cdot X_{J_1} \subset X_{J_1}$ for every $g \in \mathrm{stab}_{\circlearrowleft}(J_2) \backslash \{ 1 \}$. Finally, we need to verify that, for every $J_1, J_2 \in \mathcal{J}_0$ which are not transverse and for every $g \in \mathrm{stab}_{\circlearrowleft}(J_1) \backslash \{ 1 \}$, we have $g \cdot X_{J_2} \subset X_{J_1}$. Because $\mathcal{J}_0$ is a $x_0$-peripheral subcollection of $\mathcal{J}$, necessarily $X_{J_2}$ is contained into the sector $S$ delimited by $J_1$ which contains $x_0$. On the other hand, $g \cdot S \subset X_{J_1}$ because the action $\mathrm{stab}_{\circlearrowleft}(J_1) \curvearrowright \mathcal{S}(J_1)$ is free, hence
$$g \cdot X_{J_1} \subset g \cdot S \subset X_{J_2}.$$
Thus, all the hypotheses of Proposition \ref{prop:pingpong} are satisfied, and it follows that $\mathrm{Rot}(\mathcal{J})$ is isomorphic to $\Delta \mathcal{G}$.

\medskip \noindent
Since the hypotheses of Proposition \ref{prop:pingpong} hold, we also know from Fact \ref{fact:pingpong} that, if $g \in \mathrm{Rot}(\mathcal{J})$ is non trivial, then $g \cdot x_0 \in X_{J}$ for some $J \in \mathcal{J}$. On the other hand, $g$ permutes the connected components of $X$ cutting along the hyperplanes of $\mathcal{J}$, $Y$ is precisely the connected component which contains $x_0$, and $X_J$ is a union of connected components, so $g \cdot Y \subset Y$. We deduce from  
$$g \cdot Y \cap Y \subset X_J \cap Y = \emptyset$$
that $g \cdot Y \cap Y \neq \emptyset$, so that $g$ does not stabilise the connected component $Y$, ie., $g \notin \mathrm{stab}(Y)$. Thus, we have proved that $\mathrm{Rot}(\mathcal{J}) \cap \mathrm{stab}(Y)= \{ 1 \}$. This concludes the proof of the decomposition $G = \mathrm{Rot}(\mathcal{J}) \rtimes \mathrm{stab}(Y)$.
\end{proof}

\begin{remark}
It is worth noticing that $\mathcal{J}_0$ is $\mathrm{stab}(Y)$-invariant, so, in some sense, the action $\mathrm{stab}(Y) \curvearrowright \Delta \mathcal{G}$ factors through an action $\mathrm{stab}(Y) \curvearrowright \Delta$. See the discussion related to Question \ref{question:twistedGPhyp}.
\end{remark}

\begin{remark}
If we apply Theorem \ref{thm:rotativeactions} to the action of some graph product $\Gamma \mathcal{G}$ on the quasi-median graph we defined in Section \ref{section:GPgeometry}, we find that the group decomposes as a graph product $\Gamma \mathcal{G}$, so we do not find more structure on graph products. Similarly, if we apply Theorem \ref{thm:rotativeactions} to the action of some wreath product $G \wr H$, where $H$ acts on a CAT(0) cube complex with a vertex of trivial stabiliser, on the graph of wreaths we defined in Section \ref{section:graphofwreaths}, we find that $G \wr H$ decomposes as $\left( \bigoplus\limits_{h \in H} G \right) \rtimes H$. So once again, we obtain the initial decomposition of our group. Theorem \ref{thm:rotativeactions} also applies to right-angled graphs of groups studied in Section \ref{section:appli4}. Given such a graph of groups $\mathfrak{G}$, the semidirect product we find corresponds to the epimorphism $\varphi : \pi_1(\mathfrak{G}) \twoheadrightarrow \pi_1(\Gamma)$, where $\Gamma$ denotes the underlying graph of $\mathfrak{G}$. But in this case, the structure of $\mathrm{ker}(\varphi)$ can be investigated in full generality by using Bass-Serre theory. The situation will be more interesting in the next section, when we will consider diagram products.
\end{remark}

\noindent
As a by-product, we deduce the following characterization of groups splitting as graph products.

\begin{cor}
A group decomposes non trivially as a graph product if and only if it admits a fully rotative and vertex-free action on a quasi-median graph which is not reduced to a clique.
\end{cor}

\begin{proof}
Let $G$ be a group. If $G$ decomposes non trivially as a graph product $\Gamma \mathcal{G}$, then the action $G \curvearrowright \X$ is fully rotative and vertex-free, and $\X$ is not reduced to a clique. Conversely, suppose that $G$ admits a fully rotative and vertex-free action on a quasi-median graph $X$. By applying Theorem \ref{thm:rotativeactions} to the collection $\mathcal{J}$ of all the hyperplanes of $X$, we deduce that $G$ decomposes as $\mathrm{Rot} ( \mathcal{J} )  \rtimes \mathrm{stab} (Y)$, where $Y$ is a vertex of $X$ and where $\mathrm{Rot}(\mathcal{J})$ decomposes as a graph product along the crossing graph of the subcollection $\mathcal{J}_0 \subset \mathcal{J}$ containing the hyperplanes $J$ such that $x_0 \in N(J)$. On the other hand, since the action $G \curvearrowright X$ is vertex-free, vertex-stabilisers are trivial, so that $G = \mathrm{Rot}(\mathcal{J})$. Moreover, since $X$ is not reduced a single clique, we can choose a basepoint $x_0$ which belongs to at least two cliques, so that $\mathcal{J}_0$ contains at least two hyperplanes. Consequently, $G$ decomposes non trivially as a graph product. 
\end{proof}

\subsection{Diagram products as semidirect products}\label{section:DPsemidirect}

\noindent
In this section, we apply Theorem \ref{thm:rotativeactions} to diagram products. Our motivation is twofold: finding a method to compute presentations of diagram products (in addition to Corollary \ref{cor:2complex}), and finding examples of right-angled Artin groups which are diagram groups. We refer to the next section for applications to explicit examples of diagram products.

\medskip \noindent
First of all, we need to introduce some notation. Given an h-diagram $(\Delta,e)$, write $\mathrm{bot}^-(\Delta)=u \ell v$ for some words $u,v \in \Sigma^+$ and some letter $\ell \in \Sigma$ such that $\ell$ labels the edge $e$, and define the subgroup
$$G(\Delta,e) = \left\{ \Delta \cdot \left( \epsilon(u) + g + \epsilon(v) \right) \cdot \Delta^{-1} \mid g \in G_{\ell} \right\}.$$
Next, we define a \emph{$\mathcal{P}$-diagram} as a diagram $\Delta$ whose edges are labelled by letters of $\Sigma(\mathcal{G})$ with trivial second coordinates; alternatively, $\# \Delta$ is equal to the number of cells of $\Delta$. Finally, let $\Gamma(\mathcal{P},\mathcal{G},w)$ denote the graph whose vertices are the h-diagrams $(\Delta,e)$ such that $\Delta$ is a $\mathcal{P}$-diagram satisfying $\mathrm{top}^-(\Delta)=w$; and whose edges link two transverse h-diagrams.

\medskip \noindent
Now, we are ready to state the main result of this section.

\begin{thm}\label{thm:DPsemidirect}
Let $\mathcal{P}= \langle \Sigma \mid \mathcal{R} \rangle$ be a semigroup presentation, $\mathcal{G}$ a collection of groups indexed by $\Sigma$, and $w\in \Sigma^+$ a base word. The diagram product decomposes as 
$$D(\mathcal{P}, \mathcal{G},w)= A(\mathcal{P}, \mathcal{G},w) \rtimes  D(\mathcal{P},w)$$
where $A( \mathcal{P}, \mathcal{G},w)$ denotes the graph product with underlying graph $\Gamma(\mathcal{P},\mathcal{G},w)$ and such that the vertex $(\Delta,e)$ is labelled by the group $G(\Delta,e)$.
\end{thm}

\noindent
The first towards the proof of this theorem is to understand the rotative stabilisers of linear hyperplanes.

\begin{lemma}\label{lem:DProtative}
Let $(\Delta,e)$ be an h-diagram. The rotative stabiliser of the hyperplane $J(\Delta,e)$ is equal to $G(\Delta,e)$. 
\end{lemma}

\begin{proof}
If we write $\mathrm{bot}^-(\Delta)= u \ell v$ for some words $u,v \in \Sigma^+$ and some letter $\ell \in \Sigma$ such that $\ell$ labels the edges $e$, then the hyperplane $J(\Delta,e)$ is dual to the clique
$$\left\{ \Delta \cdot \left( \epsilon(u) + g + \epsilon(v) \right) \mid g \in G_{\ell} \right\}.$$
Consequently, the rotative stabiliser of $J(\Delta,e)$ must be included into the stabiliser of this clique, which is precisely $G(\Delta,e)$. On the other hand, if $C$ is any clique dual to $J(\Delta,e)$, according to Proposition \ref{prop:linearhyp} there must exist diagrams $A$ and $B$ satisfying $\mathrm{top}^-(A)=u$ and $\mathrm{top}^-(B)=v$ such that
$$C = \left\{ \Delta \cdot \left( A+ g+ B \right) \mid g \in G_{\ell} \right\}.$$
Because $G(\Delta,e)$ clearly stabilises this clique, we conclude that $G(\Delta,e)$ is indeed the rotative stabiliser of $J(\Delta,e)$. 
\end{proof}

\begin{proof}[Proof of Theorem \ref{thm:DPsemidirect}.]
Let $\mathcal{J}$ denote the collection of linear hyperplanes of $X(\mathcal{P}, \mathcal{G},w)$, and let $x_0= \epsilon(w)$ be our basepoint. We claim that the action of $D(\mathcal{P}, \mathcal{G},w)$ on $X(\mathcal{P}, \mathcal{G},w)$ is $\mathcal{J}$-rotative. Indeed, if $J \in \mathcal{J}$, then there exists an h-diagram $(\Delta,e)$ such that $J=J(\Delta,e)$ according to Proposition \ref{prop:hdiaghyp}, so that, if we write $\mathrm{bot}^-(\Delta)= u \ell v$ for some words $u,v \in \Sigma^+$ and some letter $\ell \in \Sigma$ such that $\ell$ labels the edges $e$, then it follows from Proposition \ref{prop:linearhyp} that the fibers of $J$ are the subgraphs generated by 
$$\left\{ \Delta \cdot \left( A+ g + B \right) \mid \mathrm{top}^-(A)=u, \mathrm{top}^-(B)=v \right\}$$
where $g \in G_{\ell}$. So $G(\Delta,e)$, which turns out to be the rotative stabiliser of $J(\Delta,e)$ according to Lemma \ref{lem:DProtative}, acts transitive and freely on the set of the fibers of $J(\Delta,e)$, proving our claim.

\medskip \noindent
Thus, Theorem \ref{thm:rotativeactions} applies, so that our diagram products decomposes as
$$D(\mathcal{P}, \mathcal{G},w) = \mathrm{Rot}(\mathcal{J}) \rtimes \mathrm{stab}(Y)$$
where $\mathrm{Rot}(\mathcal{J})= \langle \mathrm{stab}_{\circlearrowleft}(J), \ J \in \mathcal{J} \rangle$ and where $Y$ denotes the intersection of all the sectors delimited by the hyperplanes of $\mathcal{J}$ which contain $x_0$. 

\medskip \noindent
Notice that $Y$ is the subgraph of $X(\mathcal{P}, \mathcal{G},w)$ generated by the $\mathcal{P}$-diagrams, because the edges dual to linear hyperplanes are precisely the edges which are labelled by linear diagrams. As a consequence, $\mathrm{stab}(Y)$ is the subgroup of $D(\mathcal{P}, \mathcal{G},w)$ defined by the spherical $\mathcal{P}$-diagram, which is precisely the diagram group $D(\mathcal{P},w)$. 

\medskip \noindent
Next, we know from Proposition \ref{prop:hdiaghyp} that there is a bijection between $\mathcal{J}$ and the set of the h-diagrams $(\Delta,e)$ satisfying $\mathrm{top}^-(\Delta)=w$. Let $\mathcal{J}_0$ denote the subcollection of $\mathcal{J}$ corresponding to the subset of the h-diagrams $(\Delta,e)$ such that $\Delta$ is a $\mathcal{P}$-diagram. We claim that $\mathcal{J}_0$ is the maximal $x_0$-peripheral subcollection of $\mathcal{J}$.

\medskip \noindent
Let $J=J(\Delta,e)$ be a hyperplane of $\mathcal{J}_0$. According to Fact \ref{fact:hdiagproj}, $\Delta$ is the projection of $\epsilon(w)$ onto $N(J)$. A fortiori, any hyperplane separating $x_0=\epsilon(w)$ and $J$ must separate $x_0$ and $\Delta$. But $\Delta$ is a $\mathcal{P}$-diagram, so that no linear hyperplane can separate $\epsilon(w)$ and $\Delta$ as a consequence of Lemma \ref{lem:XPgeodesic}. This proves that $\mathcal{J}_0$ is a $x_0$-peripheral subcollection of $\mathcal{J}$. Now, let $J=J(\Delta,e)$ be a hyperplane which does not belong to $\mathcal{J}_0$. This means that $\Delta$ is not a $\mathcal{P}$-diagram, so that there exists a linear hyperplane $H$ separating $x_0=\epsilon(w)$ according to Lemma \ref{lem:XPgeodesic}. On the other hand, we know from Fact \ref{fact:hdiagproj} that $\Delta$ is the projection of $x_0=\epsilon(w)$ onto $N(J)$, so we deduce from Lemma \ref{lem:projseparate} that $H$ separates $x_0$ and $J$, so that $\mathcal{J}_0 \cup \{ H \}$ cannot be an $x_0$-peripheral subcollection of $\mathcal{J}$. This proves our claim.

\medskip \noindent
We conclude that the decomposition of $\mathrm{Rot}(\mathcal{J})$ as a graph product provided by Theorem \ref{thm:rotativeactions} is precisely $A(\mathcal{P}, \mathcal{G},w)$ as a consequence of the previous observations, of Proposition \ref{prop:hdiaghyp} and of Lemma \ref{lem:DProtative}. 
\end{proof}

\noindent
If the diagram group $D(\mathcal{P},w)$ is trivial and if the groups of $\mathcal{G}$ are either trivial or infinite cyclic, it follows from Theorem \ref{thm:DPsemidirect} that the diagram product turns out to be a right-angled Artin group. On the other hand, we know from \cite[Theorem 4]{MR1725439} that a diagram product of diagram groups is a diagram group, so that we get a method to produce right-angled Artin groups which are also diagram groups. It is worth noticing that the graph $\Gamma(\mathcal{P}, \mathcal{G},w)$ does not depend heavily on the groups of $\mathcal{G}$. Indeed, we only need to know whether a group of $\mathcal{G}$ is trivial or not. Let us encode this data by a function $\sigma : \Sigma \to \{ 0,1 \}$, such that $\sigma(\ell)=0$ if $G_{\ell}$ is trivial and $\sigma(\ell)=1$ otherwise, and let $\Gamma(\mathcal{P}, \sigma,w)$ denote the corresponding graph. 

\medskip \noindent
The previous discussion implies:

\begin{cor}\label{cor:RAAGDP}
Let $\mathcal{P}= \langle \Sigma \mid \mathcal{R} \rangle$ be a semigroup presentation, $\sigma : \Sigma \to \{ 0,1 \}$ a map and $w \in \Sigma^+$ a base word. If the diagram group $D(\mathcal{P},w)$ is trivial, then the right-angled Artin group $A(\Gamma( \mathcal{P}, \sigma,w))$ turns out to be a diagram group. 
\end{cor}

\noindent
Below are some examples of right-angled Artin groups arising in this way.

\begin{figure}
\begin{center}
\includegraphics[scale=0.7]{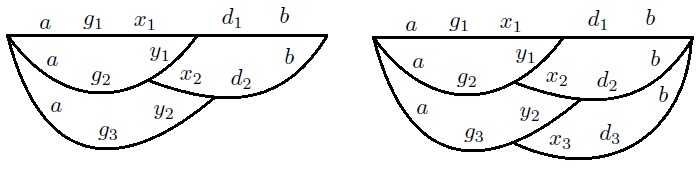}
\end{center}
\caption{The diagrams $G_3$ and $D_3$ respectively.}
\label{figure21}
\end{figure}

\begin{ex}
Let $\Sigma= \{a,b,x_i,y_i,g_i,d_i, i \geq 1 \}$ be our alphabet, and
$$\mathcal{R}= \{ ag_ix_i=ag_{i+1}y_i, \ y_id_ib=x_{i+1}d_{i+1}b, \ i \geq 1 \}$$
our collection of relations. For every $n \geq 1$, let $\mathcal{R}_n$ denote the collection of the first $n$ relations of $\mathcal{R}$, and $\Sigma_n$ the set of letters involved. The semigroup presentation we are interested in is $\mathcal{P}_n = \langle \Sigma_n \mid \mathcal{R}_n \rangle$. Let $\sigma : \Sigma \to \{0,1 \}$ be the map defined by $\sigma(a)=\sigma(b)=0$, $\sigma(x_i)=\sigma(y_i)=0$ for every $i \geq 1$, and $\sigma(g_i)=\sigma(d_i)=1$ for every $i \geq 1$. We claim that the graph $\Lambda_n = \Gamma ( \mathcal{P}_n, \sigma, ag_1x_1d_1b)$ is a segment of length $n-1$. 

\medskip \noindent
Let $G_1$ denote the trivial diagram $\epsilon(ag_1x_1d_1b)$, and, for every $i \geq 2$, let $G_i$ denote the unique $\mathcal{P}$-diagram satisfying $\mathrm{top}^-(G_i)=ag_1x_1d_1b$ and $\mathrm{bot}^-(G_i)=ag_iy_{i-1}d_{i-1}b$. Similarly, for every $i \geq 1$, let $D_i$ denote the unique $\mathcal{P}$-diagram satisfying $\mathrm{top}^-(D_i)=ag_1x_1d_1b$ and $\mathrm{bot}^-(D_i)= ag_ix_id_ib$. For instance, $G_3$ and $D_3$ are illustrated by Figure \ref{figure21}. Notice that the $G_i$'s and the $D_i$'s are the only $\mathcal{P}$-diagrams, and they all have a maximal thin suffix reduced to a single cell. So $\{ (D_i,d_i), (G_i,g_i) \mid i \geq 1 \}$ are the vertices of $\Lambda = \Gamma(\mathcal{P},\sigma, ag_1x_1d_1b)$. Moreover, for every $i \geq 1$, the h-diagram $(G_i,g_i)$ is transverse only to $(D_i,d_i)$ and to $(D_{i-1},d_{i-1})$ (when it exists), and $(D_i,d_i)$ is transverse only to $(G_i,g_i)$ and to $(G_{i+1},g_{i+1})$. Thus, $\Lambda$ is an infinite ray, and $\Lambda_n$ is the initial subsegment of length $n-1$.

\medskip \noindent
Consequently, for every $n \geq 1$, if $P_n$ denotes the segment of length $n$, then the right-angled Artin group $A(P_n)$ turns out to be a diagram group.
\end{ex}

\begin{ex}
To any collection of intervals on the real line is associated a graph whose set of vertices is our collection and whose edges link two intersecting intervals. A graph arising in this way is called an \emph{interval graph}. We proved in \cite{arXiv:1507.01667} that right-angled Artin groups associated to complements of interval graphs are diagram groups. We are able to reprove this result thanks to Corollary \ref{cor:RAAGDP}. So let $\Gamma$ be a finite interval graph. Since it is finite, we may suppose without loss of generality that $\Gamma$ is associated to a collection $\mathcal{C}$ of intervals on $\{ 1, \ldots, n\}$, for some $n \geq 1$. Now, consider the semigroup presentation
$$\mathcal{P}(\mathcal{C}) = \langle x_1, \ldots, x_n, a_I, I \in \mathcal{C} \mid x_I=a_I, \ I \in \mathcal{C} \rangle,$$
where we use the notation $x_I = x_i \cdot x_{i+1} \cdots x_{i+j}$ for every interval $I= \{ i,i+1, \ldots, i+j \}$. Let $\sigma : \Sigma \to \{ 0,1 \}$ be the map defined by $\sigma(x_i)=0$ for every $1 \leq i \leq n$ and $\sigma(a_I)=1$ for every $I \in \mathcal{C}$. We claim that the graph $\Lambda = \Gamma(\mathcal{P}(\mathcal{C}), \sigma,x_1 \cdots x_n)$ is isomorphic to $\bar{\Gamma}$. Indeed, if $A_I$ denotes the atomic diagram satisfying $\mathrm{top}^-(A_I)=x_1 \cdots x_n$ and labelled by the relation $x_I \to a_I$, then $\{ (A_I,a_I) \mid I \in \mathcal{C} \}$ is the set of vertices of our graph $\Lambda$. Moreover, two h-diagrams $(A_I,a_I)$ and $(A_J,a_J)$ are transverse if and only if $I$ and $J$ are disjoint intervals. Consequently, the right-angled Artin group $A\left( \bar{\Gamma} \right)$ turns out to be a diagram group. 
\end{ex}

\subsection{Examples}

\noindent
This section is dedicated to the description of a few explicit families of diagram products.

\begin{ex}
If $\mathcal{P}= \langle a,b,p \mid a=ap, b=pb \rangle$ and $\mathcal{G}=\{ G_a=G_b=\{1 \}, G_p=G\}$ for some group $G$, then the diagram product $D(\mathcal{P}, \mathcal{G}, ab)$ is isomorphic to the wreath product $G \wr \mathbb{Z}$ \cite[Example 10]{MR1725439}. 
\begin{figure}
\begin{center}
\includegraphics[scale=0.7]{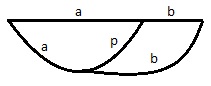}
\end{center}
\caption{The generator $\Delta$ of $D(\mathcal{P},ab)$.}
\label{figure23}
\end{figure}
\begin{figure}
\begin{center}
\includegraphics[scale=0.7]{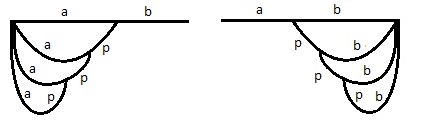}
\end{center}
\caption{From left to right: $A_3$ and $B_3$.}
\label{figure24}
\end{figure}

\medskip \noindent
Indeed, the only h-diagrams we have are $\{(A_n,a), (B_n,b) \mid n \geq 1 \}$, where $A_n$ denotes the unique $\mathcal{P}$-diagram satisfying $\mathrm{top}^-(A_n)=ab$, $\mathrm{bot}^-(ap^nb)$ and $\mathrm{supp}(A_n)=a$, and similarly, $B_n$ the unique $\mathcal{P}$-diagram satisfying $\mathrm{top}^-(B_n)=ab$, $\mathrm{bot}^-(ap^nb)$ and $\mathrm{supp}(B_n)=b$; for instance $A_3$ and $B_3$ are illustrated by Figure \ref{figure24}. Moreover, any two h-diagrams are transverse. On the other hand, by computing the fundamental group of the Squier complex $S(\mathcal{P},ab)$, it follows that the diagram group $D(\mathcal{P},ab)$ is infinite cyclic, generated by the diagram $\Delta$ illustrated by Figure \ref{figure23}, and we notice that $\Delta \cdot A_n=A_{n+1}$ for every $n \geq 1$, $\Delta \cdot B_n= B_{n-1}$ for every $n \geq 2$ and $\Delta \cdot B_1=A_1$. Consequently, by applying Theorem \ref{thm:DPsemidirect}, it follows that
$$D(\mathcal{P}, \mathcal{G},ab)= \left( \bigoplus\limits_{n \in \mathbb{Z}} G \right) \rtimes \mathbb{Z},$$
where $\mathbb{Z}$ acts by shifting the coordinates. This proves our claim.

\medskip \noindent
We do not elaborate more on this example since we dedicated Section \ref{section:appli3} to the study of wreath products. 
\end{ex}

\begin{ex}
If $\mathcal{P}= \langle a,b,p \mid a=ap, b=pb \rangle$ and $\mathcal{G}=\{ G_a=G, G_b=H, G_p=\{1 \} \}$ for some groups $G$ and $H$, then the diagram product $D(\mathcal{P}, \mathcal{G}, ab)$ admits 
$$G \bullet H = \langle G,H,t \mid [g,t^nht^{-n}]=1, g \in G, h \in H,n \geq 0 \rangle$$ 
as a (relative) presentation \cite[Section 8]{MR1396957}. 

\medskip \noindent
Indeed, by computing the fundamental group of the Squier complex $S(\mathcal{P},ab)$ it follows that the diagram group $D(\mathcal{P},ab)$ is infinite cyclic, generated by the diagram $\Delta$ illustrated by Figure \ref{figure23}. Moreover, the only h-diagrams we have are $\{ (\Delta^n,a), (\Delta^m,b) \mid n,m \in \mathbb{Z} \}$, where we identify $(\Delta^n,a)$, when $n \geq 0$, with the h-diagram corresponding to the smallest prefix of $\Delta^n$ containing the edge of $\mathrm{bot}(\Delta^n)$ labelled by $a$, and similarly for $(\Delta^n,b)$, when $n \leq 0$, with respect to the edge of $\mathrm{bot}(\Delta)$ labelled by $b$; for instance, the h-diagrams $(\Delta^3,a)$, $(\Delta^{-3},a)$, $(\Delta^{3},b)$ and $(\Delta^{-3},b)$ are illustrated by Figure \ref{figure22}. Finally, for every $n,m \in \mathbb{Z}$, the h-diagrams $(\Delta^n,a)$ and $(\Delta^m,b)$ are transverse if and only if $n \geq m$. Now, we are ready to apply Theorem \ref{thm:DPsemidirect}. Let $\Gamma(G,H)$ denote the graph product whose underlying graph $\Lambda$ has $\mathbb{Z} \times \{ 0,1 \}$ as vertex-set and $\{ ((n,0),(m,1)) \mid n \geq m \}$ as edge-set, such that the vertices $(n,0)$ are labelled by the group $G$ and the vertices $(n,1)$ by the group $H$. So our diagram product $D(\mathcal{P},\mathcal{G},ab)$ decomposes as the semidirect product $\Gamma(G,H) \rtimes \mathbb{Z}$, where $\mathbb{Z}$ (corresponding to the cyclic subgroup generated by $\Delta$) acts on the graph product $\Gamma(G,H)$ via the translation $(n,i) \mapsto (n+1,i)$ of $\Lambda$. Therefore, given infinitely many copies $G_n,H_m$ of $G,H$ respectively, where $n,m \in \mathbb{Z}$, we find that our diagram product admits
$$\left\langle t, H_n,G_m, n,m \in \mathbb{Z} \left| \begin{array}{l} [g_{(n)},h_{(m)}]=1, n \geq m \\ tg_{(n)}t^{-1}=g_{(n+1)}, th_{(m)}t^{-1} = h_{(m+1)}, n,m \in \mathbb{Z} \end{array}, g \in G,h \in H \right. \right\rangle $$
as a (relative) presentation, where $g_{(n)}$ (resp. $h_{(m)}$) denotes the element $g \in G$ in the copy $G_n$ (resp. the element $h \in H$ in the copy $H_m$). By noticing that $g_{(n)}= t^ng_{(0)}t^{-n}$ for every $n \in \mathbb{Z}$ and every $g \in G$, and $h_{(m)}=t^m h_{(0)}t^{-m}$ for every $m \in \mathbb{Z}$ and every $h \in H$, this presentation simplifies as
$$\left\langle H,G,t \mid \left[ t^ngt^{-n}, t^mht^{-m} \right] = 1, g \in G, h \in H, n \geq m \right\rangle,$$
which gives the presentation mentionned above. 
\begin{figure}
\begin{center}
\includegraphics[scale=0.67]{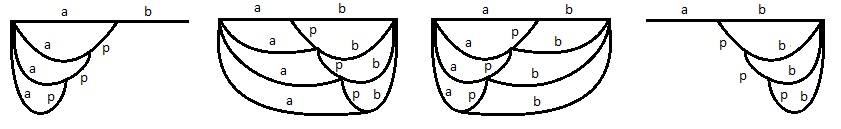}
\end{center}
\caption{From left to right: $(\Delta^3,a)$, $(\Delta^{-3},a)$, $(\Delta^{3},b)$ and $(\Delta^{-3},b)$.}
\label{figure22}
\end{figure}

\medskip \noindent
According to Theorems \ref{thm:DPpropercc} and \ref{thm:DPmetricallyproper}, if $G$ and $H$ act (metrically) properly on CAT(0) cube complexes of dimensions $d_G,d_H \geq 1$ respectively, then the $G \bullet H$ acts (metrically) properly on a CAT(0) cube complex of dimension $d_G+d_H$. Moreover, a $\bullet$-product of two a-T-menable (resp. a-$L^p$-menable, where $p \notin 2 \mathbb{Z}$) groups is a-T-menable (resp. a-$L^p$-menable). On the other hand, $G \bullet H$ is essentially never hyperbolic. More precisely, it is hyperbolic if and only if one of $G,H$ is trivial and the other hyperbolic. Indeed, if $g \in G$ and $h \in H$ are two non trivial elements, then $gtgt^{-1}$ is an infinite-order commuting $t^nht^{-n}$ for every $n \leq 0$, so that $\langle gtgt^{-1} \rangle$ does not have finite index in its centraliser, which implies that $G \bullet H$ cannot be hyperbolic. Conversely, say if $H$ is trivial, then $G \bullet H$ reduces to a free product $G \ast \mathbb{Z}$, which is hyperbolic if so is $G$. (In fact, more generally, we can prove that $G \bullet H$ is relatively hyperbolic if and only if one of $G,H$ is trivial; see for instance \cite[Fact 5.42]{article3} for the torsion-free case.) Nevertheless, the product $G \bullet H$ turns out to be acylindrically hyperbolic \cite[Example 5.41]{article3} whatever are $G$ and $H$. 

\medskip \noindent
About $\ell^p$-compressions of $G \bullet H$, Theorem \ref{thm:DPcompressions} does not apply since $[ab]_{\mathcal{P}}$ is not finite (indeed, it contains $ap^nb$ for every $n \geq 0$). Nevertheless, we expect that the canonical orbit map $G \bullet H \hookrightarrow (X(\mathcal{P}, \mathcal{G},w),\delta)$, where $\delta$ is a global metric constructed by extending word metrics on cliques, is a quasi-isometric embedding, so that, for every $p \geq 1$,
$$\alpha_p^*(G \bullet H) \geq \min \left( \frac{1}{p}, \alpha_p^*(G), \alpha_p^*(H) \right),$$
with equality for $p \in [1,2]$. We already know that this is true for $\mathbb{Z} \bullet \mathbb{Z}$ from \cite{arXiv:1507.01667}, in which we study we study this group. In fact, we also proved a stronger statement: $\mathbb{Z} \bullet \mathbb{Z}$ is isomorphic to the subgroup $\langle a,c,bd \rangle$ in the direct product of two free groups $\langle a,b \mid \ \rangle \times \langle c,d \mid \ \rangle$.
\end{ex}

\begin{ex}\label{ex:DPnonRAAG}
Consider the semigroup presentation  
$$\mathcal{P}= \langle a,b,c,d \mid ab=ac, bd=cd \rangle$$ 
and the collection $\mathcal{G}= \{ G_b=G_c=\{1 \}, G_a=G, G_d=H \}$ for some groups $G,H$. The Squier complex $S(\mathcal{P},abd)$ is a cycle of length, so computing the fundamental group of the graph of groups provided by Corollary \ref{cor:2complex} shows that the diagram product $D(\mathcal{P}, \mathcal{G},abd)$ admits
$$G \square H= \langle G,H, t \mid [g,h]=[g,h^t]=1, g \in G, h \in H \rangle$$
as a (relative) presentation. Notice that $G \square \{ 1 \} \simeq G \ast \mathbb{Z}$, so we will suppose in the sequel that $G$ and $H$ are both non trivial. We introduced the group $\mathbb{Z} \square \mathbb{Z}$ in \cite[Example 5.43]{article3} as an example of a cocompact diagram group which is not a right-angled Artin group. 

\medskip \noindent
According to Theorems \ref{thm:DPpropercc}, \ref{thm:DPmetricallyproper} and \ref{thm:DPCAT0}, we know that
\begin{itemize}
	\item if $G$ and $H$ act properly (resp. metrically properly, geometrically) on some CAT(0) cube complexes of dimensions $d_G,d_H \geq 1$ respectively, then $G \square H$ acts properly (resp. metrically properly, geometrically) on some CAT(0) cube complex of dimension $d_G+d_H$;
	\item if $G$ and $H$ are CAT(0), then so is $G \square H$;
	\item if $G$ and $H$ are a-T-menable (resp. a-$L^p$-menable, where $p \notin 2 \mathbb{Z}$), then $G \square H$ is a-T-menable (resp. a-$L^p$-menable).
\end{itemize}
Moreover, it follows from Theorem \ref{thm:DPcompressions} that
$$\alpha_p^*(G \square H) \geq \min\left( \frac{1}{p}, \alpha_p^*(G), \alpha_p^*(H) \right)$$
for every $p \geq 1$, with equality if $p \in [1,2]$ since $G \square H$ contains a quasi-isometrically embedded non abelian free subgroup. Indeed, $G \square H$ decomposes as an HNN extension $(G \times (H \ast H)) \underset{H}{\ast}$, so that it follows from Britton's lemma that the subgroup $\langle htht^{-1},t \rangle$, where $h$ denotes a non trivial element of $H$, defines a quasi-isometrically embedded non abelian free subgroup of $G \square H$. 

\medskip \noindent
Next, notice $G \square H$ is hyperbolic if and only if $G$ and $H$ are both finite. Indeed, if $G$ is infinite and $h$ denotes a non trivial element of $H$, then $\langle h tht^{-1} \rangle$ defines an infinite cyclic subgroup in the centraliser of $G$ whose intersection with $G$ is reduced to the identity. This implies that $G \square H$ cannot be hyperbolic. On the other hand, if $G$ and $H$ are both finite, then decomposition of $G \square H$ as an HNN extension $(G \times (H \ast H)) \underset{H}{\ast}$ we mentionned above shows that $G \square H$ is virtually free. In fact, if the statement mentionned in the discussion related to Question \ref{question:twistedGPhyp} holds, it follows that no hyperbolic group constructed from a diagram product can be one-ended.

\medskip \noindent
Following the proof of Fact \ref{fact:isometricX}, it is worth noticing that the quasi-median graph associated to $G \square H$ depends only on the cardinalities of $G$ and $H$. In fact, we claim that the arguments of the proof of Theorem \ref{thm:GPlip} apply similarly, so that, if $G_1$ and $G_2$, and $H_1$ and $H_2$, are Lipschitz-equivalent groups, then $G_1 \square H_1$ and $G_2 \square H_2$ are Lipschitz-equivalent. For instance, $\mathbb{F}_n \square \mathbb{Z}$ and $\mathbb{F}_m \square \mathbb{Z}$ are Lipschitz-equivalent for every $n,m \geq 2$. 
\end{ex}

\section{Application to right-angled graphs of groups}\label{section:appli4}

\noindent
We begin this section by fixing the basic definitions and notations related to graph of groups; essentially, we follow \cite{MR1954121}. So far, our graphs were always one-dimensional simplicial complexes, but we need a different definition in order to define graph of groups. To avoid ambiguity, we will refer to the latter as \emph{abstract graphs}. 

\begin{definition}
An \emph{abstract graph} is the data of a set of vertices $V$, a set of arrows $E$, a fixed-point-free involution $e \mapsto \bar{e}$ of $E$, and two maps $s,t : E \to E$ satisfying $t(e)= s \left( \bar{e} \right)$ for every $e \in E$. 
\end{definition}

\noindent
Below, we define graph of groups and their associated fundamental groupoids as introduced in \cite{Higgins}. 

\begin{definition}
A \emph{graph of groups}\index{Graphs of groups} $\mathfrak{G}$ is the data of an abstract graph $(V,E, \bar{\cdot}, s,t)$, a collection of groups indexed by $V \sqcup E$ such that $G_e= G_{\bar{e}}$ for every $e \in E$, and a monomorphism $\iota_e : G_e \hookrightarrow G_{s(e)}$ for every $e \in E$. The \emph{fundamental groupoid} $\mathfrak{F} = \mathfrak{F}(\mathfrak{G})$ of $\mathfrak{G}$ is the groupoid which has vertex set $V$, which is generated by the arrows of $E$ together with $\bigsqcup\limits_{v \in V} G_v$, and satisfying the relations:
\begin{itemize}
	\item for every $v \in V$ and $g,h,k \in G_v$, $gh=k$ if the equality holds in $G_v$;
	\item for every $e \in E$ and $g \in G_e$, $\iota_e(g) \cdot e = e \cdot \iota_{\bar{e}}(g)$. 
\end{itemize}
Notice in particular that, for every $e \in E$, $\bar{e}$ is an inverse of $e$ in $\mathfrak{F}$. Fixing some vertex $v \in V$, the \emph{fundamental group} of $\mathfrak{G}$ (based at $v$) is the vertex-group $\mathfrak{F}_v$ of $\mathfrak{F}$, ie., the loops of $\mathfrak{F}$ based at $1_v$. 
\end{definition}

\noindent
The following normal form, proved in \cite{Higgins}, will be fundamental for our study.

\begin{prop}\label{prop:GFnormalform}
Let $\mathfrak{G}$ be a graph of groups. For every $e \in E$, fix a left-transversal $T_e$ of $\iota_e \left( G_e \right)$ in $G_{s(e)}$ containing $1_{s(e)}$. Any element of $\mathfrak{F}$ can be written uniquely as a word $g_1 \cdot e_1 \cdots g_n \cdot e_n \cdot g_{n+1}$, where 
\begin{itemize}
	\item $(e_1,\ldots, e_n)$ is a direct path of the underlying abstract graph;
	\item $g_i \in T_{e_i}$ for every $1 \leq i \leq n$, and $g_{n+1}$ is an arbitrary element of $G_{t(e_n)}$;
	\item if $e_{i+1} = \bar{e_i}$ for some $1 \leq i \leq n-1$ then $g_{i+1} \neq 1$.
\end{itemize}
\end{prop}

\noindent
Such a word will be referred to as a \emph{normal word}. 

\medskip \noindent
Roughly speaking, we will be interested in graphs of groups gluing graph products. To get something interesting, we need to control the gluings.

\begin{definition}
Given two graph products $\Gamma \mathcal{G}$ and $\Lambda \mathcal{H}$, a morphism $\Phi : \Gamma \mathcal{G} \to \Lambda \mathcal{H}$ is a \emph{graphical embedding} is there exists an embedding $f : \Gamma \to \Lambda$ and an isomorphism $\varphi_v : G_v \to H_{f(v)}$ such that $f(\Gamma)$ is an induced subgraph of $\Lambda$ and $\Phi(g)= \varphi_v(g)$ for every $v \in V(\Gamma)$ and $g \in G_v$. 
\end{definition}

\noindent
Typically, we glue graph products along ``subgraph products'' in a canonical way. We refer to Section \ref{section:RAGGex} for examples.

\begin{definition}
A \emph{right-angled graph of groups}\index{Right-angled graph of groups} is a graph of groups such that each (vertex- and edge-)group has a fixed decomposition as a graph product and and such that each monomorphism of an edge-group into a vertex-group is a graphical embedding (with respect to the structures of graph products we fixed). 
\end{definition}

\noindent
In the following, a \emph{factor} $G$ will refer to a vertex-group of one of these graph products. If $G \subset G_v$ for some $v \in V$, we say that $G$ is labelled by $v$. 

\medskip \noindent
Fix a right-angled graph of groups $\mathfrak{G}$. For every arrow $e \in E$, there exists a natural left-transversal $T_e$ of $\iota_e(G_e)$ in $G_{s(e)}$: the set of reduced words of $G_{s(e)}$ whose tails (see Definition \ref{def:headtail}) do not contain any element of the vertex-groups of $\iota_e(G_e)$. From now on, we fix this choice, and any normal word will refer to this convention.

\medskip \noindent
Let $\mathfrak{S} \subset \mathfrak{F}$ denote the union of the arrows of $E$ together with the vertex-groups (minus the identity) of the graph products $G_v$, $v \in V$. By definition, $\mathfrak{S}$ is a generating set of the fundamental groupoid $\mathfrak{F}$ of $\mathfrak{G}$. Any element $g \in \mathfrak{F}$ can be represented as a word $g_1e_1 \cdots g_ne_ng_{n+1}$ written over $\mathfrak{S}$, where 
\begin{itemize}
	\item for every $1 \leq i \leq n$, $e_i \in E$;
	\item for every $1 \leq i \leq n$, $g_i$ is a reduced word of the graph product $G_{s(e_i)}$, and similarly for $g_{n+1} \in G_{t(e_n)}$;
	\item $g_1e_1 \cdots g_ne_ng_{n+1}$ is a normal word representing $g$ when thought of as a word written over $E \sqcup \bigsqcup\limits_{v \in V} G_v$. 
\end{itemize}
Notice that the length of this word is equal to 
$$n+\sum\limits_{i=1}^{n-1} |g_i|,$$
where $|g_i|$ denotes the length of $g_i$ as an element of the corresponding graph product. In particular, it does not depend on the choice of the reduced elements representing the $g_i$'s. We refer to a word as above as a \emph{reduced word} representing $g$, and denote its length by $|g|$. The following lemma shows that $| \cdot |$ coincides with the word length associated to $\mathfrak{S}$.

\begin{lemma}\label{lem:shortestlength}
For every $g \in \mathfrak{F}$, $|g|$ is the shortest length of a word written on $\mathfrak{S}$ representing $g$.
\end{lemma}

\begin{proof}
Let $g_1e_1 \cdots g_n e_n g_{n+1}$ be a word written over $\mathfrak{S}$. Consider the following operations:
\begin{itemize}
	\item if $g_i$ is trivial for some $2 \leq i \leq n$ and $e_{i-1}= \bar{e_i}$, remove the subword $e_{i-1}g_ie_i$;
	\item if $g_i = a \cdot \iota_{e_i}(b)$ for some $b \in G_{e_i}$ and $1 \leq i \leq n$, replace the subword $g_i \cdot e_i$ with $a \cdot e_i \cdot \iota_{\bar{e_i}}(b)$;
	\item if $g_i$ is not reduced for some $1 \leq i \leq n+1$, replace $g_i$ with a reduced word representing it.
\end{itemize}
Notice that these operations do not increase the length of a word written over $\mathfrak{S}$. On the other hand, it is clear that any word written over $\mathfrak{S}$ can be made reduced by applying finitely many times these operations. Therefore, any word written over $\mathfrak{S}$ representing some $g \in \mathfrak{F}$ has length at least $|g|$. This concludes the proof.
\end{proof}

\noindent
We conclude this section with the following lemma which explains how to reduce the product of a reduced word with a generator. For convenience, if $g$ is a word representing an element of a graph product, we denote by $[g]$ a reduction of $g$.

\begin{lemma}\label{lem:reducingproduct}
Let $g_1e_1 \cdots g_n e_n g_{n+1}$ be a reduced word representing an element $g \in \mathfrak{F}$ and $s \in \mathfrak{S}$ a generator. If $s \in E$, write $g_{n+1}=a \cdot \iota_s(b)$ where $\iota_s(b)$ the suffix of (a reduced word representing) $g_{n+1}$ containing the syllables of its tail which belong to $\iota_s(G_s)$.Then
$$\left\{ \begin{array}{cl} g_1e_1 \cdots  g_n e_n a s \iota_{\bar{s}}(b) & \text{if $a \neq \emptyset$ or $e_n \neq \bar{s}$} \\ \\ g_1e_1 \cdots g_{n-1}e_{n-1} [g_n \iota_{\bar{s}}(b)] & \text{if $a = \emptyset$ and $e_n= \bar{s}$} \end{array} \right.$$
is a reduced word representing $gs$. Otherwise, $g_1e_1 \cdots g_n e_n [g_{n+1}s]$ is a reduced word representing $gs$. 
\end{lemma}

\begin{proof}
Suppose that $s \in E$, $a= \emptyset$, and $e_n= \bar{s}$. It is clear that $g_1e_1 \cdots g_{n-1}e_{n-1}[g_n \iota_{\bar{s}}(b)]$ is a normal word, and it represents $gs$ because 
$$\begin{array}{lcl} g_n \cdot e_n \cdot g_{n+1} \cdot s & = & g_n \cdot e_n \cdot a \cdot \iota_s(b) \cdot s = g_n \cdot e_n \cdot a \cdot s \cdot \iota_{\bar{s}}(b) \\ \\ & = & g_n \cdot \bar{s} \cdot s \cdot \iota_{\bar{s}}(b) = g_n \cdot \iota_{\bar{s}}(b) \end{array}$$
Therefore, $g_1e_1 \cdots g_{n-1}e_{n-1} [g_n \iota_{\bar{s}}(b)]$ is a reduced word representing $gs$. 

\medskip \noindent
Now, suppose that $s \in E$ and either $a \neq \emptyset$ or $e_n \neq \bar{s}$. Notice that, by definition of our decomposition $g_{n+1}= a \cdot \iota_s(b)$, $a$ belongs to $T_s$. Moreover, if $a$ is trivial then $e_n \neq \bar{s}$ by assumption. Therefore, the word $g_1e_1 \cdots  g_n e_n a s \iota_{\bar{s}}(b)$ is normal. And it represents $gs$ because
$$g_{n+1} \cdot s = a \cdot \iota_s(b) \cdot s = a \cdot s \cdot \iota_{\bar{s}}(b).$$
Consequently, $g_1e_1 \cdots  g_n e_n a s \iota_{\bar{s}}(b)$ is a reduced word representing $gs$.

\medskip \noindent
Finally, the case $s \notin E$ follows directly by applying the definition of a reduced word.
\end{proof}

\subsection{Cubical-like geometry}\label{section:RAGGgeometry}

\noindent
Fix a right-angled graph of groups $\mathfrak{G}$, and a vertex $\omega \in V$ of its underlying abstract graph. Define $\mathfrak{X}= \mathfrak{X}(\mathfrak{G},\omega)$ as the connected component of the Cayley graph $\mathfrak{X}(\mathfrak{G})$ of $\mathfrak{F}$, constructed from the generating set $\mathfrak{S}$, which contains the neutral element $1_{\omega}$ based at $\omega$. Equivalently, $\mathfrak{X}$ is the graph whose vertices are the arrow of $\mathfrak{F}$ starting from $\omega$ and whose edges link two elements $g,h \in \mathfrak{F}$ if $g=h \cdot s$ for some $s \in \mathfrak{S}$. It is worth noticing that an edge of $\mathfrak{X}$ is naturally labelled either by an element of $E$, or by an element of $G_v$ for some $v \in V$, and in particular by $v$. Moreover, the fundamental group $\mathfrak{F}_{\omega}$ of $\mathfrak{G}$ (based at $\omega$) acts naturally on $\mathfrak{X}$ by left-multiplication. 

\medskip \noindent
When we work on the Cayley graph of a group, because the action of the group is vertex-transitive, we may always suppose that a fixed base point is the identity element up to a translation. In our situation, where $\mathfrak{X}(\mathfrak{G})$ is the Cayley graph of a groupoid, the situation is slightly different but it is nevertheless possible to make something similar. Indeed, if $x \in \mathfrak{X}(\mathfrak{G},\omega)$ is a base point, then translating by $x^{-1}$ sends $x$ to an identity element, but which does not necessarily belong to $\mathfrak{X}(\mathfrak{G},\omega)$: in general, it will belong to $\mathfrak{X}(\mathfrak{G},\xi)$ where $\xi \in V$ is the terminal point of $x$ when thought of as an element of $\mathfrak{F}$. In the same time, $\mathfrak{F}_{\omega}$ becomes $\mathfrak{F}_{\xi}$, which is naturally isomorphic to $\mathfrak{F}_{\omega}$ via the conjugacy $g \mapsto xgx^{-1}$. Finally, we get a commutative diagram
\begin{displaymath}
\xymatrix{ \mathfrak{F}_{\omega} \ar[rr]^{\text{isomorphism}} \ar[d] & & \mathfrak{F}_{\xi} \ar[d] \\ \mathfrak{X}(\mathfrak{G},\omega) \ar[rr]^{\text{isometry}} & & \mathfrak{X}(\mathfrak{G}, \xi) }
\end{displaymath}
Thus, we may always suppose that a fixed base point is an identity element up to changing the basepoint $\omega$, which does not disturb either the group or the Cayley graph. 

\medskip \noindent
Now, we want to prove that our Cayley graph is quasi-median.

\begin{prop}\label{prop:RAGGqm}
The graph $\mathfrak{X}$ is quasi-median.
\end{prop}

\noindent
First of all, notice that it follows from the definition of $\mathfrak{X}$ and from Lemma \ref{lem:shortestlength} that

\begin{lemma}\label{lem:RAGGdist}
For every vertices $g,h \in \mathfrak{X}$, the equality $d_{\mathfrak{X}}(g,h)=|g^{-1}h|$ holds. 
\end{lemma}

\noindent
It is worth noticing that, by construction of $\mathfrak{X}$, one finds in $\mathfrak{X}$ many copies of the quasi-median graphs associated to the graph products of our graph of groups. In order to apply the results we proved in Section \ref{section:GPgeometry}, we need to show that these copies are nicely embedded into $\mathfrak{X}$. This is the purpose of our next lemma. 

\begin{definition}
A \emph{leaf} of $\mathfrak{X}$ is the subgraph generated by the set of vertices $gG_v$, where $g \in \mathfrak{F}$ is some arrow starting from $\omega$ and ending to some $v \in V$.
\end{definition}

\begin{lemma}\label{lem:GPleaf}
A leaf $\Lambda$ of $\mathfrak{X}$ is a gated subgraph. Moreover, if $\Lambda =g G_v$ for some $g \in \mathfrak{F}$ and $v \in V$, then 
$$\left\{ \begin{array}{ccc} G_v & \to & \Lambda \\ h & \mapsto & gh \end{array} \right.$$
induces an isometry from the quasi-median graph associated to the graph product $G_v$ onto $\Lambda$. A fortiori, $\Lambda$ is a quasi-median graph on its own right.
\end{lemma}

\begin{proof}
Let $\Lambda$ be a leaf and $g \in \mathfrak{X}$ a vertex. Up to translating by $g^{-1}$, suppose that $g=1_{\omega}$. Write $\Lambda=h G_v$, where $h \in \mathfrak{F}$ and $v \in V$, such that $h$ can be written as a reduced word $h_1e_1 \cdots h_ne_nh_{n+1}$ where the tail of $h_{n+1}$ does not contain any syllable of $G_v$. For every $s \in G_v \backslash \{ 1 \}$, it follows from Lemma \ref{lem:reducingproduct} that $h_1e_1 \cdots h_ne_n h_{n+1}s$ is a reduced word representing $hs$, hence $|hs|=|h|+1$. Thus, we have proved that $d(1_{\omega},q) = d(1_{\omega},h)+ d(h,q)$ for every vertex $q \in \Lambda$, so that $h$ is the gate of $1_{\omega}$ in $\Lambda$. Consequently, $\Lambda$ is a gated subgraph.

\medskip \noindent
Now, suppose that $\Lambda = gG_v$ for some $g \in \mathfrak{F}$. Say that $G_v$ is a graph product $\Gamma \mathcal{G}$. The map
$$\left\{ \begin{array}{ccc} G_v & \to & \Lambda \\ h & \mapsto & gh \end{array} \right.$$
induces an isometry $\X \to \Lambda$ according to Lemma \ref{lem:RAGGdist} and Corollary \ref{cor:GPdist}. A fortiori, $\Lambda$ must be a quasi-median graph as a consequence of Proposition \ref{prop:Xquasimedian}.
\end{proof}

\begin{proof}[Proof of Proposition \ref{prop:RAGGqm}.]
Let us begin by verifying the triangle condition. Let $g,h,k \in \mathfrak{X}$ be three vertices such $h$ and $k$ are adjacent, and the distances $d(g,h)$ and $d(g,k)$ are equal, say to $d$. Up to translating by $g^{-1}$, we can suppose without loss of generality that $g=1_{\omega}$. As a consequence, we deduce from Lemma \ref{lem:RAGGdist} that $d(g,h)=|h|$ and $d(g,k)=|k|$. Because $h$ and $k$ are adjacent, there exists some $s \in \mathfrak{S}$ such that $k=h \cdot s$. By noticing that $|h|=|h \cdot s|$, we deduce from Lemma \ref{lem:reducingproduct} that $s \notin E$, ie., $s \in G_v$ for some $v \in V$. Therefore, $h$ and $k$ belong to the leaf $\Lambda = hG_v$. Let $p \in \Lambda$ denote the gate of $1_{\omega}$ in $\Lambda$, which exists according to Lemma \ref{lem:GPleaf}. Notice that
$$d(p,h)=d(1_{\omega},h)- d(1_{\omega},p)= d(1_{\omega},k)- d(1_{\omega},p)= d(p,k).$$
On the other hand, we know from Lemma \ref{lem:GPleaf} that $\Lambda$ is a quasi-median graph, so there exists some $r \in \Lambda$ which is adjacent to both $h$ and $k$, and such that $d(p,r)=d(p,h)-1$. By noticing that
$$d(1_{\omega},r)=d(1_{\omega},p)+ d(p,r)=d(1_{\omega},p)+ d(p,h)-1 = d(1_{\omega},h)-1=d-1,$$
we conclude that $r$ is the vertex we are looking for. Thus, the triangle condition is satisfied.

\medskip \noindent
Next, let us verify the quadrangle condition. Let $g,h,j,k \in \mathfrak{X}$ be four pairwise disctinct vertices such that $j$ is adjacent to both $h$ and $k$, such that the distances $d(g,h)$ and $d(g,k)$ are equal, say to $d$, and such that $d(g,j)=d+1$. Up to translating by $g^{-1}$, we can suppose without loss of generality that $g=1_{\omega}$. 

\medskip \noindent
Suppose first that the edges $(h,j)$ and $(j,k)$ are not labelled by $E$. So there exist $v \in V$ and $r,s \in G_v$ such that $j=h \cdot r$ and $k=j \cdot s$. A fortiori, $h$, $j$ and $k$ belong to the leaf $\Lambda= h G_v$. Let $p$ denote the gate of $1_{\omega}$ in $\Lambda$, which exists according to Lemma \ref{lem:GPleaf}. Notice that
$$d(p,h)= d(1_{\omega},h)- d(1_{\omega},p) = d(1_{\omega},k) - d(1_{\omega},p) = d(p,k),$$
and that
$$d(p,j)=d(1_{\omega},j)- d(1_{\omega},p) = d(1_{\omega},h)-d(1_{\omega},p)-1 = d(h,p)-1.$$
Because $\Lambda$ is quasi-median according to Lemma \ref{lem:GPleaf}, there must exist some $q \in \Lambda$ which is adjacent to both $h$ and $k$, and such that $d(p,q)=d(p,h)-2$. By noticing that
$$d(1_{\omega},q)= d(1_{\omega},p)+d(p,h)-2= d(1_{\omega},h)-2=d-2,$$
we conclude that $q$ is the vertex we are looking for.

\medskip \noindent
Now, suppose that the edges $(h,j)$ and $(j,k)$ are both labelled by $E$. So there exist two edges $f_1,f_2 \in E$ such that $j=h \cdot f_1$ and $k=j \cdot f_2$. Write $h$ as a reduced word $h_1e_1 \cdots h_ne_nh_{n+1}$. Because
$$|h \cdot f_1| = d(1_{\omega}, j) = d+1= d(1_{\omega},h) +1= |h|+1,$$
we deduce from Lemma \ref{lem:reducingproduct} that $h_{n+1}$ can be written as $a \cdot \iota_{f_1}(b)$ and that $h_1 e_1 \cdots h_ne_n a f_1 \iota_{\bar{f_1}}(b)$ is a reduced word representing $h \cdot f_1=j$. Next, because
$$|j \cdot f_2| = d(1_{\omega}, k) = d = d(1_{\omega},j)-1= |j|-1,$$
we deduce from Lemma \ref{lem:reducingproduct} that $f_2= \bar{f_1}$ and that $\iota_{\bar{f_1}}(b)$ is trivial, so that $h_1e_1 \cdots h_ne_na$ is a reduced word representing $j \cdot f_2=k$. On the other hand, since $\iota_{\bar{f_1}}(b)$ is trivial, necessarily $h_{n+1}= a$, hence
$$k= h_1e_1 \cdots h_ne_na = h_1e_1 \cdots h_ne_n h_{n+1}=h.$$
This constracts our assumption that $g,h,j,k$ are four pairwise distinct vertices. 

\medskip \noindent
Finally, suppose that exactly one edge among $(h,j)$ and $(j,k)$ is labelled by $E$, say $(h,j)$. So there exist $e \in E$, $v \in V$ and $s \in G_v$ such that $j=h \cdot e$ and $k = j \cdot s$. Notice that necessarily $v=t(e)$, since otherwise the products would be not defined. Write $h$ as a reduced word $h_1e_1 \cdots h_ne_nh_{n+1}$. Because
$$|h \cdot e| = d(1_{\omega}, j) = d+1= d(1_{\omega},h) +1= |h|+1,$$
we deduce from Lemma \ref{lem:reducingproduct} that $h_{n+1}$ can be written as $a \cdot \iota_{e}(b)$ and that $h_1 e_1 \cdots h_ne_n a e \iota_{\bar{e}}(b)$ is a reduced word representing $h \cdot e=j$. Next, once again according to Lemma \ref{lem:reducingproduct}, we know that $h_1e_1 \cdots h_ne_nae [\iota_{\bar{e}}(b)s]$ is a reduced word representing $j \cdot s=k$. Furthermore,
$$\begin{array}{lcl} 1 & = & d(g,j)-d(g,k)=  |j|-|k| \\ \\ & = & \displaystyle \left(n+1+ \sum\limits_{i=1}^n |h_i| + |a|+ |\iota_{\bar{e}}(b)| \right) - \left( n+1 + \sum\limits_{i=1}^n |h_i| + |a| + |\iota_{\bar{e}}(b) \cdot s| \right) \\ \\ & = & |\iota_{\bar{e}}(b)| - |\iota_{\bar{e}}(b) \cdot s| \end{array}$$
hence $|\iota_{\bar{e}}(b) \cdot s| = |\iota_{\bar{e}}(b)|-1$. Therefore, the product $\iota_{\bar{e}}(b) \cdot s$ is not reduced in the graph product $G_{t(e)}$, so that $s$ must belong to $\iota_{\bar{e}}(G_e)$, say $s= \iota_{\bar{e}}(r)$ where $r \in G_e$. Now, we claim that $p= h_1e_1 \cdots h_ne_n a [\iota_{e}(br)]$ is the vertex we are looking for. By noticing that
$$p \cdot s^{-1} = h_1e_1 \cdots h_ne_n \cdot (a \cdot \iota_{e}(b)) \cdot s \cdot s^{-1} = h_1e_1 \cdots h_ne_nh_{n+1}=h$$
and that
$$p \cdot e = h_1 e_1 \cdots h_ne_na \cdot \iota_e(br) \cdot e = h_1e_1 \cdots h_ne_na \cdot e \cdot \iota_{\bar{e}}(br)=k,$$
we deduce that $p$ is adjacent to both $h$ and $k$. Next, we want to prove that $d(1_{\omega},p)=d-2$. First, notice that
$$|\iota_e(br)| = |\iota_{\bar{e}}(br)| = |\iota_{\bar{e}}(b)s|= |\iota_{\bar{e}}(b)|-1,$$
where the first equality is a consequence of the observation that graphical embeddings preserve the lengths. Next, since we proved that $h_1 e_1 \cdots h_ne_n a e \iota_{\bar{e}}(b)$ is a reduced word representing $j$, we deduce that
$$\begin{array}{lcl} d(1_{\omega},p) & = & \displaystyle |p| \leq n+ \sum\limits_{i=1}^n |h_i| + |a \iota_e(br)| \\ \\ & \leq & \displaystyle n+ \sum\limits_{i=1}^n |h_i| + |a| + |\iota_e(br)| \\ \\ & \leq & \displaystyle n+ \sum\limits_{i=1}^n |h_i| + |a| + |\iota_{\bar{e}}(b)| -1  = d(1_{\omega},j)-2 \end{array}$$
On the other hand,
$$d(1_{\omega},p) \geq d(1_{\omega},j)-d(j,h)-d(h,p) = d(1_{\omega},j)-2.$$
Consequently, $d(1_{\omega},p)=d-2$, which proves that $p$ is indeed the vertex we are looking for. This concludes the proof of the quadrangle condition. 

\medskip \noindent
Finally, if $K$ is a subgraph of $\mathfrak{X}$ isomorphic to $K_{4}^-$, we deduce from Lemma \ref{lem:GPtriangle} below that the edges of $K$ must belong to leaves. But a leaf contains its triangles, since it is gated according to Lemma \ref{lem:GPleaf}, so $K$ must be included into a single leaf. Because a leaf is a quasi-median graph according to Lemma \ref{lem:GPleaf}, we conclude that $K$ cannot be induced in $\mathfrak{X}$. Next, it follows from Lemma \ref{lem:GPdisttwo} below that $\mathfrak{X}$ cannot contain an induced subgraph isomorphic to $K_{3,2}$. Thus, we have proved that $\mathfrak{X}$ is a quasi-median graph. 
\end{proof}

\begin{lemma}\label{lem:GPtriangle}
An edge of $\mathfrak{X}$ which belongs to a triangle cannot be labelled by $E$. 
\end{lemma}

\begin{proof}
Let $g,h,k \in \mathfrak{X}$ be three vertices such that $h$ is adjacent to both $g$ and $k$, and such that the edge $(h,k)$ is labelled by $E$. Up to translating by $g^{-1}$, suppose without loss of generality that $g=1_{\omega}$. It follows from Lemma \ref{lem:reducingproduct} that $|k|= |h| \pm 1$. Therefore,
$$d(g,k) = |k| = |h| \pm 1 = d(g,h) \pm 1 \in \{ 0,2 \}.$$
A fortiori, $g$ and $k$ cannot be adjacent. This concludes the proof.
\end{proof}

\begin{lemma}\label{lem:GPdisttwo}
Between any two vertices of $\mathfrak{X}$ at distance two appart, there exist at most two geodesics.
\end{lemma}

\begin{proof}
It is sufficient to prove that an element of length two in $\mathfrak{F}$ can be written in at most two different ways as a product of two elements of $\mathfrak{S}$. We distinguish four cases:
\begin{itemize}
	\item Let $g \in \mathfrak{F}$ be an element which can be written as a reduced word $e_1e_2$, where $e_1,e_2 \in E$. Because the relations of $\mathfrak{F}$ preserve the generators of $E$, $e_2e_1$ is the only candidate for another reduced word representing $g$. However, this not the case if $e_1 \neq e_2$ according to our normal form. Consequently, $g$ can be written as a product of two elements of $\mathfrak{S}$ in a unique way. 
	\item Let $g \in \mathfrak{F}$ be an element which can be written as a reduced word $he$ where $e \in E$ and $h \in G_v$ for some $v \in V$. Because the relations of $\mathfrak{F}$ preserve the generators of $E$, any other word of length two representing $g$ must have the form $ke$ or $ek'$ where $k \in G_{s(e)}$ and $k' \in G_{t(e)}$. The equality $he=ke$ implies $h=k$, and because $ek'$ is reduced we deduce from our normal form that $he \neq ek'$. Consequently, $g$ can be written as a product of two elements of $\mathfrak{S}$ in a unique way. 
	\item Let $g \in \mathfrak{F}$ be an element which can be written as a reduced word $eh$ where $e \in E$ and $h \in G_v$ for some $v \in V$. Because the relations of $\mathfrak{F}$ preserve the generators of $E$, any other word of length two representing $g$ must have the form $ke$ or $ek'$ where $k \in G_{s(e)}$ and $k' \in G_{t(e)}$. The equality $eh=ek'$ implies $h=k'$. If $eh=ke$, then $ke$ is not a normal word, so $k= \iota_{e}(j)$ for some $j \in G_e$, and we find that $ke=e \cdot \iota_{\bar{e}}(j)$. From the equality $e \cdot \iota_{\bar{e}}(j) = eh$, we deduce that $h= \iota_{\bar{e}}(j)$. Consequently, $g$ can be written as a product of two elements of $\mathfrak{S}$ in at most two ways.
	\item Let $g \in \mathfrak{F}$ be an element which can be written as a reduced word $hk$ where $hk \in G_v$ for some $v \in V$. Because the relations of $\mathfrak{F}$ preserve the generators of $E$, any other reduced word representing $g$ must belong to $G_v$. It follows from the normal form in graph products that $kh$ is the only candidate for another reduced word representing $hk$. Consequently, $g$ can be written as a product of two elements of $\mathfrak{S}$ in at most two ways.
\end{itemize}
This concludes the proof of our lemma.
\end{proof}

\noindent
It is worth noticing that the first point of the previous proof implies the following observation, which will be useful later:

\begin{fact}\label{fact:squareE}
Two consecutive edges of a square of $\mathfrak{X}$ cannot be labelled by $E$. 
\end{fact}

\paragraph{Cliques and prisms.} Now, we want to understand the cliques and prisms of $\mathfrak{X}$. The classification of cliques is given by the following lemma. Recall that a factor is a vertex-group of some graph products in our right-angled graph of groups.

\begin{lemma}\label{lem:RAGGclique}
A clique of $\mathfrak{X}$ is either an edge labelled by $E$ or a complete subgraph $gG$ where $G$ is a factor and $g \in \mathfrak{F}$. 
\end{lemma}

\begin{proof}
Let $Q$ be a clique. If $Q$ contains an edge $e$ labelled by $E$, then it follows from Lemma \ref{lem:GPtriangle} that $Q=e$. Otherwise, $Q$ must belong to a leaf, and the conclusion follows from Lemma \ref{lem:Xclique}, which applies thanks to Lemma \ref{lem:GPleaf}. 
\end{proof}

\noindent
Next, we focus on the prisms of $\mathfrak{X}$. The first step is to understand the squares. In the following lemma, for every edge $e \in E$ we denote by $\varphi_e$ the morphism $\iota_{\bar{e}} \circ \iota_e^{-1} : \iota_e(G_e) \to \iota_{\bar{e}}(G_e)$. 

\begin{lemma}\label{lem:RAGGsquare}
Let $S$ be a square of $\mathfrak{X}$. Fix a vertex $h \in S$ and let $hg,hs \in S$ denote its neighbors where $g,s \in \mathfrak{S}$. Suppose that $g \in G_v$ for some $v \in V$. Let $w\in S$ denote the fourth vertex. Either $s \in E$, so that $g \in \iota_s(G_s)$ and $w=hgs=hs \varphi_s(g)$; or $s \in G_v$, and $g$ and $s$ belong to two adjacent vertex-groups of the graph product $G_v$. 
\end{lemma}

\begin{proof}
Because $w$ is adjacent to both $hg$ and $hs$, there exist $p,q \in \mathfrak{S}$ such that $w=hgp=hsq$. As a consequence, the equality $gp=sq$ must hold. 

\medskip \noindent
Suppose that $s \in E$. Because the relations of $\mathfrak{F}$ preserves $E$, it follows from the equality $gp=sq$ that $p=s \in E$ and that $q \in G_{t(s)}$. Next, because $sq$ is a normal form, we deduce that $gs$ cannot be a normal word, so that $g \in \iota_s(G_s)$ and $gs=s \varphi_s(g)$. From the equality $sq=gs=s \varphi_s(g)$, we deduce that $q= \varphi_s(g)$.

\medskip \noindent
Next, suppose that $s \in G_u$ for some $u \in V$. Because the product $hg$ is well-defined and that $g \in G_v$, we deduce that $u=v$ since the product $hs$ is also well-defined. Notice that $h,hg,hs$ belong to the leaf $\Lambda = h G_v$. Because $\Lambda$ is gated (and, in particular, convex), necessarily $w \in \Lambda$. So $S$ is a square of $\Lambda$. We deduce from Corollary \ref{cor:Xprisms} that $g$ and $s$ belong to two adjacent vertex-groups of the graph product $G_v$. 
\end{proof}

\begin{cor}\label{cor:RAGGsquare}
Let $S$ be a square and $(a,b)$ one of its edge. Suppose that this edge is labelled by $V$. There exists a unique $s \in \mathfrak{S}$ such that the opposite edge of $(a,b)$ in $S$ is $(as,bs)$. 
\end{cor}

\begin{proof}
Let $c$ denote the neighbor of $a$ in $S$ which is distinct from $b$, and let $s \in \mathfrak{S}$ be a generator such that $c=a \cdot s$. According to Lemma \ref{lem:RAGGsquare}, the edge of $S$ which is opposite to $(a,b)$ is $(as,bs)$. If there exists another $r \in \mathfrak{S}$ such that this edge can be written as $(ar,br)$, then the edge $(a,as)=(a,ar)$ is labelled by $s = a^{-1} \cdot (as) = a^{-1} \cdot (ar) = r$, so the uniqueness follows.
\end{proof}

\noindent
Now, we are ready to describe the prisms of $\mathfrak{X}$. Notice that we already understand the prisms which are included into some leaf, as a consequence of Lemma \ref{lem:GPleaf} and Corollary \ref{cor:Xprisms}. The other prisms are characterized by our next lemma.

\begin{lemma}\label{lem:RAGGprism}
For every prism $Q$ of $\mathfrak{X}$ which is not included into a leaf, there exist some $e \in E$ and some prism $P$ which is included into a leaf, such that $Q$ is generated by the set of vertices $\{ g, \ ge \mid g \in P \}$.
\end{lemma}

\begin{proof}
Because $Q$ is not included into a leaf, it must contain an edge $f$ labelled by some $e \in E$. Cutting along the hyperplane $J$ dual to $f$ decomposes $Q$ into two subprisms $Q_1,Q_2$. According to Fact \ref{fact:squareE}, $Q_1$ and $Q_2$ cannot contain edges labelled by $E$, so they are included into some leaves. We claim that the edges of $Q$ which do not belong to neither $Q_1$ nor $Q_2$ are labelled by $e$. Let $a$ be such an edge. There exists a chain of adjacent squares from $f$ to $a$, so we can suppose that $a$ and $f$ are opposite edges of some square $S$, the general case following by induction of the length of the chain. We deduce from Corollary \ref{cor:RAGGsquare} that $a$ is labelled by $e$, proving our claim. As a consequence, if $h \in Q_2$ is a vertex, then $h=g \cdot e$ where $g$ is the unique vertex of $Q_1$ adjacent to $h$. Thus, we have proved that the set of vertices of $Q$ is $\{ g,g \cdot e \mid g \in Q_1 \}$, where $Q_1$ is a prism which is included into a leaf.
\end{proof}

\paragraph{Hyperplanes.} The rest of the section is dedicated to the study of the hyperplanes of $\mathfrak{X}$. It is worth noticing that, as a consequence of Lemma \ref{lem:RAGGclique}, Fact \ref{fact:squareE} and Lemma \ref{lem:RAGGsquare}, all the edges of a given hyperplane are labelled by either some edge of $E$ or by $V$, so that the hyperplanes of $\mathfrak{X}$ are naturally labelled by an edge of $E$ or by $V$; moreover, two hyperplanes labelled by $E$ cannot be transverse. 

\medskip \noindent
For convenience, we will use the following notations. If $e \in E$, denote by $\varphi_e : \iota_e(G_e) \to \iota_{\bar{e}}(G_e)$ the isomorphism $\iota_{\bar{e}} \circ \iota_e^{-1}$. A priori, $\varphi_e$ is not defined on all $G_{s(e)}$, but for every subset $S \subset G_{s(e)}$, we can define $\varphi_e(S)$ as $\varphi_e \left( S \cap \iota_e(G_e) \right)$.  

\medskip \noindent
Roughly speaking, the neighbor of the hyperplane dual to the clique associated to some factor $G$ is generated by the vertices corresponding to elements of $\mathfrak{F}$ which ``commute'' with all the elements of $G$. Because commutation is not well-defined in groupoids, we need to define carefully this idea, which is done by the following definition.

\begin{definition}
Let $G$ be a factor labelled by $v \in V$. An element $h \in \mathfrak{F}$ belongs to the \emph{link} of $G$, denoted by $\mathrm{link}(G)$, if it can be written as normal word $h_1e_1 \cdots h_ne_nh_{n+1}$ such that
\begin{itemize}
	\item for every $1 \leq i \leq n-1$, $\varphi_{e_i} \left( \cdots \left( \varphi_{e_1} (G) \right) \cdots \right) \subset \iota_{e_{i+1}}(G_{e_{i+1}})$;
	\item for every $1 \leq i \leq n$, $h_i$ belongs to the link of the vertex-group $\varphi_{e_{i-1}} \left( \cdots \left( \varphi_{e_{1}}(G) \right) \cdots \right)$ of the graph product $G_{s(e_i)}$;
	\item $h_{n+1}$ belongs to the link of the vertex-group $\varphi_{e_n} \left( \cdots \left( \varphi_{e_1}(G) \right) \cdots \right)$ of the graph product $G_{t(e_n)}$.
\end{itemize}
\end{definition}

\noindent
The ``commutation'' with the elements of $G$ mentionned earlier is made precise by the following lemma. 

\begin{lemma}\label{lem:RAGGcommute}
Let $G$ be a factor and $\ell$ an element of its link. Write $\ell = \ell_1 e_1 \cdots \ell_n e_n \ell_{n+1}$. The restriction $\varphi_{\ell}^G$ of the map $\varphi_{e_n} \circ \cdots \circ \varphi_{e_1}$ to $G$ defines an isomorphism $G \to \varphi_{e_n} \left( \cdots \left( \varphi_{e_1}(G) \right) \cdots \right)$. Moreover, for every $g \in G$, $g \cdot \ell= \ell \cdot \varphi_{\ell}^G(g)$. 
\end{lemma}

\begin{proof}
We argue by induction on $n$. Suppose that $n=0$. Notice that $\varphi_{\ell}^G$ is the identity $G \to G$. Say that $G$ is a vertex-group of the graph product $G_v$ where $v \in V$. We know that $\ell= \ell_1$ belongs to the link of the vertex-group $G$ in the graph product $G_v$, hence $g \cdot \ell = \ell \cdot g$ for every $g \in G$. Now, suppose that $n \geq 1$, and set $m= \ell_1 e_1 \cdots \ell_{n-1} e_{n-1} \ell_n$, so that $\ell = m e_n \ell_{n+1}$. By our induction hypothesis, since $m$ belongs to the link of $G$, the restriction $\varphi_m^G$ of $\varphi_{e_{n-1}} \circ \cdots \circ \varphi_{e_1}$ to $G$ defines an isomorphism $G \to \varphi_{e_{n-1}} \left( \cdots \left( \varphi_{e_1}(G) \right) \cdots \right)$ such that $g \cdot m = m \cdot \varphi_m^G(g)$ for every $g \in G$. On the other hand, we know that the image of $\varphi_{m}^G$ is included into $\iota_{e_n}(G_{e_n})$, which is the domain of the isomorphism $\varphi_{e_n}$, so necessarily $\varphi_{\ell}^G = \varphi_{e_n} \circ \varphi_m^G$ induces an isomorphism $G \to \varphi_{e_n} \left( \cdots \left( \varphi_{e_1}(G) \right) \cdots \right)$. Moreover, for every $g \in G$, 
$$g \cdot \ell = m \cdot \varphi^G_m(g) \cdot e_n \cdot \ell_{n+1} = me_n \cdot \varphi^G_{\ell}(g) \cdot \ell_{n+1} = \ell \cdot \varphi_{\ell}^G(g),$$
where the last equality is justified by the observation that $\ell_{n+1}$ belongs to the link of the vertex-group $\varphi_{e_n}\left( \cdots \left( \varphi_{e_1}(G) \right) \cdots \right)$ of the graph product $G_{t(e_n)}$ and that $\varphi^G_{\ell}(g)$ belongs to $\varphi_{e_n}\left( \cdots \left( \varphi_{e_1}(G) \right) \cdots \right)$. This concludes the proof of our lemma.
\end{proof}

\noindent
Now, we are ready to describe the hyperplanes of $\mathfrak{X}$ which are lablled by factors. In fact, they are the only hyperplanes we need to understand precisely, the others being dual to finite cliques. 

\begin{prop}\label{prop:RAGGhyp}
Let $G$ be a factor and let $J$ denote the hyperplane dual to the clique $G$. An edge $e$ of $\mathfrak{X}$ is dual to $J$ if and only if $e=(h_1 \ell,h_2 \ell)$ for some $h_1,h_2 \in G$ distinct and $\ell \in \mathrm{link}(G)$. As a consequence, $N(J)= G \cdot \mathrm{link}(G)$ and the fibers of $J$ are the $h \cdot \mathrm{link}(G)$ where $h \in G$.
\end{prop}

\begin{proof}
Suppose that $e$ is dual to $J$. According to Claim \ref{claim:chainofsquareshyp}, there exists a chain of adjacent squares from an edge of the clique $G$ and $e$. We argue by induction on the number of squares. If it is zero, there is nothing to prove. Otherwise, its antepenultimate edge is $(h_1 \ell,h_2 \ell)$ where $\ell \in \mathrm{link}(G)$ and $h_1,h_2 \in G$ are distinct. According to Corollary \ref{cor:RAGGsquare}, $e=(h_1 \ell s, h_2 \ell s)$ for some $s \in \mathfrak{S}$. Write $\ell= \ell_1 e_1 \cdots \ell_n e_n \ell_{n+1}$. Notice that $(h_1 \ell,h_2 \ell)= \left( \ell \cdot \varphi_{\ell}^G(h_1), \ell \cdot \varphi_{\ell}^G(h_2) \right)$. We want to prove that $\ell s \in \mathrm{link}(G)$. 

\medskip \noindent
Suppose that $s \in G_v$ for some $v \in V$. Noticing that $\ell_1 e_1 \cdots \ell_n e_n (\ell_{n+1} s)$ is a normal word representing $\ell s$. The only thing to verify in order to prove that $\ell s$ belongs to $\mathrm{link}(G)$ is that $\ell_{n+1}s$ belongs to the link of the vertex-group $\varphi_n \left( \cdots \left( \varphi_{e_1} (G) \right) \cdots \right)$ of the graph product $G_{t(e_n)}$. Notice that, because the product $\ell s$ is well-defined, necessarily $G_{t(e_n)}=G_v$. We already know that our statement is true for $\ell_{n+1}$. On the other hand, we know from Lemma \ref{lem:RAGGsquare} that $s$ and $\varphi_{\ell}^G(h_1)^{-1} \varphi_{\ell}^G(h_2)$ belongs to two adjacent vertex-groups of $G_v$, so that the conclusion follows since $\varphi_{\ell}^G(h_1)^{-1} \varphi_{\ell}^G(h_2)$ belongs to $\varphi_n \left( \cdots \left( \varphi_{e_1} (G) \right) \cdots \right)$. 

\medskip \noindent
Suppose that $s \in E$. Write $\ell_{n+1}$ as $a \cdot \iota_s(b)$, where $\iota_s(b)$ is the suffix of $\ell_{n+1}$ containing the syllables of its tail which belong to $\iota_s(G_s)$. According to Lemma \ref{lem:reducingproduct}, 
$$\ell_1 e_1 \cdots \ell_{n-1} e_{n-1} a s \iota_{\bar{s}}(b)$$
is a normal word representing $\ell s$. In order to prove that $\ell s \in \mathrm{link}(G)$, we have to verify that
\begin{itemize}
	\item[(i)] $a$ belongs to the link of the vertex-group $\varphi_{e_n} \left( \cdots \left( \varphi_{e_1}(G) \right) \cdots \right)$ in in the graph product $G_{t(e_n)}$;
	\item[(ii)] $\varphi_s\left( \varphi{e_n} \left( \cdots \left( \varphi_{e_1}(G) \right) \cdots \right) \right) \subset \iota_s(G_s)$. 
	\item[(iii)] $\iota_{\bar{s}}(b)$ belongs to the link of the vertex-group $\varphi_s\left( \varphi_{e_n} \left( \cdots \left( \varphi_{e_1}(G) \right) \cdots \right) \right)$ in the graph product $G_{t(s)}$;
\end{itemize}
The point $(i)$ follows from the observation that $a$ is a subword of $\ell_{n+1}$ and that we know that $\ell_{n+1}$  belongs to the link of the vertex-group $\varphi_{e_n} \left( \cdots \left( \varphi_{e_1}(G) \right) \cdots \right)$ in in the graph product $G_{t(e_n)}$ since $\ell \in \mathrm{link}(G)$. Next, we deduce from Lemma \ref{lem:RAGGsquare} that $\varphi_{\ell}^G(h_1)^{-1} \varphi_{\ell}^G(h_2) \in \iota_s(G_s)$. Because $\varphi_{\ell}^G(h_1)^{-1} \varphi_{\ell}^G(h_2)$ is a non trivial element of the factor $\varphi_{e_n} \left( \cdots \left( \varphi_{e_1} (G) \right) \cdots \right)$, and that $\iota_s(G_s)$ is a subgroup of the graph product $G_{t(e_n)}$ generated by vertex-groups, we deduce that the point $(ii)$ is satisfied. Finally, notice that, because $\iota_s(b)$ is a subword of $\ell_{n+1}$, it must belong to the link of the vertex-group $\varphi_{e_n} \left( \cdots \left( \varphi_{e_1}(G) \right) \cdots \right)$ in in the graph product $G_{t(e_n)}$. Therefore, $\varphi_s( \iota_s(b)) = \iota_{\bar{s}}(b)$ belongs to the link of the vertex-group $\varphi_s\left( \varphi_{e_n} \left( \cdots \left( \varphi_{e_1}(G) \right) \cdots \right) \right)$ in the graph product $\varphi_s(G_{t(e_n)}) = G_{t(s)}$. 

\medskip \noindent
Conversely, fix $h_1,h_2 \in G$ distinct and $\ell \in \mathrm{link}(G)$. Write $\ell = \ell_1e_1 \cdots \ell_n e_n \ell_{n+1}$. Consider the paths
$$h_1, \ h_1 \ell_1, \ h_1 \ell_1e_1, \ldots, \ h_1 \ell_1 e_1 \cdots \ell_n e_n$$
and
$$h_2, \ h_2 \ell_1, \ h_2 \ell_1e_1, \ldots, \ h_2 \ell_1 e_1 \cdots \ell_n e_n.$$
Notice that the distance between $h_1 \ell_1 e_1 \cdots \ell_i e_i$ and $h_2 \ell_1 e_1 \cdots \ell_i e_i$ is equal to 
$$d\left( \ell_1e_1 \cdots \ell_ie_i \cdot \varphi_{\ell_1e_1 \cdots \ell_ie_i}^G(h_1), \ell_1e_1 \cdots \ell_ie_i \cdot \varphi_{\ell_1e_1 \cdots \ell_ie_i}^G(h_2) \right)=1.$$
Similarly, $h_1 \ell_1 e_1 \cdots \ell_i e_i \ell_{i+1}$ and $h_2 \ell_1 e_1 \cdots \ell_i e_i \ell_i$ are adjacent for every $i$. Therefore, our paths define a chain of adjacent squares from $(h_1,h_2)$ to $(h_1 \ell,h_2 \ell)$. A fortiori, $(h_1 \ell,h_2 \ell)$ is dual to $J$. 

\medskip \noindent
Thus, we have proved the first statement of our proposition. As a consequence, the set of vertices of $N(J)$ is necessarily included into $G \cdot \mathrm{link}(G)$. Conversely, if $g \in G$ and $\ell \in \mathrm{link}(G)$, then, fixing some $h \in G \backslash \{ g \}$, the vertices $g \cdot \ell$ and $h \cdot \ell$ are adjacent since
$$d(g \cdot \ell, h \cdot \ell) = d( \ell \cdot \varphi_{\ell}^G(g), \ell \cdot \varphi_{\ell}^G(h)) = \left| \varphi_{\ell}^G \left( gh^{-1} \right) \right| =1,$$
and the edge $(g \cdot \ell, h \cdot \ell)$ is dual to $J$, so that $g \cdot \ell$ must belong to $N(J)$. Thus, we have proved that $N(J)= G \cdot \mathrm{link}(G)$. 

\medskip \noindent
We claim that $h \cdot \ell$ belongs to the fiber of $J$ passing through $h$. Indeed, for every $g \in G$, 
$$d(h \cdot \ell,g) = d(g^{-1}h \cdot \ell, 1_{\omega}) = d \left( \ell \cdot \varphi_{\ell}^G \left( g^{-1}h \right), 1_{\omega} \right) = \left| \ell \cdot \varphi_{\ell}^G \left( g^{-1}h \right) \right|.$$
Write $\ell = \ell_1 e_1 \cdots \ell_n e_n \ell_{n+1}$. According to Lemma \ref{lem:reducingproduct}, 
$$\ell_1e_1 \cdots \ell_n e_n \cdot \left[ \ell_{n+1} \cdot \varphi_{\ell}^G \left( g^{-1}h \right) \right]$$
is a reduced word representing $\ell \cdot \varphi_{\ell}^G \left( g^{-1}h \right)$. Therefore, $d(h \cdot \ell,g)$ is minimal precisely when the length of $\ell_{n+1} \cdot \varphi_{\ell}^G \left(g^{-1}h\right)$ is minimal. On the other hand, we know that $\ell_{n+1}$ belongs to the link of the vertex-group $\varphi_{\ell}^G(G)$ in the graph product $G_{t(e_n)}$, so that $\ell_{n+1}$ and $ \varphi_{\ell}^G \left(g^{-1}h\right)$ belong to distinct vertex-groups of $\varphi_{\ell}^G(G)$, hence
$$\left| \ell_{n+1} \cdot  \varphi_{\ell}^G \left(g^{-1}h\right) \right| = | \ell_{n+1}| + \left|  \varphi_{\ell}^G \left(g^{-1}h\right) \right|.$$
Consequently, the distance $d(h \cdot \ell,g)$ is minimized precisely when $g=h$.
\end{proof}

\begin{cor}\label{cor:RAGGstabhyp}
Let $G$ be a factor and $J$ the hyperplane dual to the clique $G$. 
$$\mathrm{stab}(J)= \{ g \cdot m \mid g \in G, m \in \mathrm{link}(G), \varphi_m^G(G)=G \}.$$
\end{cor}

\begin{proof}
An element $k$ belongs to $\mathrm{stab}(J)$ if and only if, for every distinct $g_1,g_2 \in G$, there exist distinct $h_1,h_2 \in G$ and $\ell \in \mathrm{link}(G)$ such that $k \cdot (g_1,g_2) = (h_1 \ell , h_2 \ell)$. In particular, taking $g_1=1$, we deduce that there exist $g \in G$ and $m \in \mathrm{link}(G)$ such that $k= gm$. Now we want to prove that $\varphi_m^G(G)=G$. Notice that $m=g^{-1} \cdot gm$ also belongs to $\mathrm{stab}(J)$, since $G$ stabilises the clique $G$ and a fortiori the hyperplane $J$. Therefore, fixing some distinct $g_1,g_2 \in G$, we know that exist distinct $h_1,h_2 \in G$ and $\ell \in \mathrm{link}(G)$ such that $k \cdot (g_1,g_2) = (h_1 \ell , h_2 \ell)$. Write $m= m_1e_1 \cdots m_re_rm_{r+1}$ and $\ell= \ell_1 e_1' \cdots \ell_s e_s' \ell_{s+1}$. From the equality
$$m_1e_1 \cdots m_re_r(m_{r+1}g_i) = \ell_1 e_1' \cdots \ell_s e_s' ( \ell_{s+1} \varphi_{\ell}^G(h_i))$$
between two normal words, we deduce that $r=s$, $m_j=\ell_j$ and $e_j=e_j'$ for every $1 \leq j \leq r=s$, and $m_{r+1}g_i= \ell_{s+1} \varphi_{\ell}^G(h_i)$. Notice that
$$\left\{ \begin{array}{l} m_{r+1}g_1= \ell_{s+1} \cdot \varphi_{\ell}^G(h_1) \\  m_{r+1}g_2= \ell_{s+1} \cdot \varphi_{\ell}^G(h_2) \end{array} \right. \ \text{implies} \ \varphi_{\ell}^G\left( h_2^{-1}h_1 \right) = g_2^{-1}g_1 \in G.$$
Because $\varphi_{\ell}^G(G)$ and $G$ are two vertex-groups of the graph product $G_{\omega}$ which intersect non trivially, we deduce that $\varphi_{\ell}^G(G)=G$. On the other hand, because $e_i=e_i'$ for every $i$, necessarily $\varphi_{\ell}^G= \varphi_m^G$. This proves that $\varphi_m^G(G)=G$. 

\medskip \noindent
Conversely, let $g \in G$ and $m \in \mathrm{link}(G)$  be two elements, such that $\varphi_m^G(G)=G$. Notice that $\varphi_m^G$ defines an automorphism of $G$ so that it is invertible. For every distinct $g_1,g_2 \in G$ and $\ell \in \mathrm{link}(G)$, we have
$$gm \cdot g_i \ell = \left( g \left( \varphi_m^G  \right)^{-1}(g_i) \right) \cdot \left( m \ell \right)$$
where $g \left( \varphi_m^G  \right)^{-1}(g_i) \in G$ and $m \ell \in \mathrm{link}(G)$ according to Lemma \ref{lem:linksubgroupoid} below. Thus, $gm$ sends an edge of the clique $G$ to an edge dual to $J$. A fortiori, $gm \in \mathrm{stab}(J)$. 
\end{proof}

\begin{lemma}\label{lem:linksubgroupoid}
Let $G$ be a factor. For every $m \in \mathrm{link}(G)$ and $\ell \in \mathrm{link} \left( \varphi^G_m(G) \right)$, if the product $m \ell$ is well-defined then it must belong to $\mathrm{link}(G)$. 
\end{lemma}

\begin{proof}
We argue by induction on the length $L$ of a normal word representing $\ell$. If $L=0$, there is nothing to prove. Suppose that $L \geq 1$. Write $\ell$ as a reduced word $\ell_1 f_1 \cdots \ell_n f_n \ell_{n+1}$. Setting $p = \ell_1 f_1 \cdots \ell_{n-1} f_{n-1} \ell_n$, we have $\ell = p f_n \ell_{n+1}$, and, by noticing that $p$ belongs to $\mathrm{link} \left( \varphi_m^G(G) \right)$, we deduce from our induction hypothesis that $mp$ belongs to $\mathrm{link}(G)$. Write $mp$ as a reduced word $a_1 e_1 \cdots a_r e_r a_{r+1}$. Write $a_{r+1} = b \cdot \iota_{f_n}(c)$ so that Lemma \ref{lem:reducingproduct} applies. We distinguish two cases.

\medskip \noindent
Suppose that $b \neq \emptyset$ or $f_n \neq \bar{e_r}$. According to Lemma \ref{lem:reducingproduct}, $a_1e_1 \cdots a_re_rb f_n \iota_{\bar{f_n}}(c)$ is a reduced word representing $mpf_n$, so that, once again according to Lemma \ref{lem:reducingproduct}, 
$$a_1 \cdot e_1 \cdots a_r \cdot e_r \cdot b \cdot f_n \cdot \left[ \iota_{\bar{f_n}}(c) \ell_{n+1} \right]$$
is a reduced word representing $m \ell$. Because $mp \in \mathrm{link}(G)$, necessarily $a_1e_1 \cdots a_r e_r \in \mathrm{link}(G)$. Moreover, because $a_{r+1}= b \cdot \iota_{f_n}(c)$ belongs to the link of the vertex-group $\varphi_{e_r} \left( \cdots \left( \varphi_{e_1}(G) \right) \cdots \right)$ of the graph product $G_{t(e_r)}$, a fortiori $b$ belongs also to this link, and $\iota_{\bar{f_n}}(c)$ belongs to the link of the vertex-group $\varphi_{f_n} \left( \varphi_{e_r} \left( \cdots \left( \varphi_{e_1}(G) \right) \cdots \right) \right)$ of the graph product $G_{t(f_n)}$. Now, notice that, because $\ell \in \mathrm{link}\left( \varphi_m^G(G) \right)$ and because the words $e_1 \cdots e_r$ and $f_1 \cdots f_{n-1}$ contain the same sequence of edges of $E$ up to cancellation in $\mathfrak{F}$, necessarily
$$\varphi_{e_r} \left( \cdots \left( \varphi_{e_1}(G) \right) \cdots \right) = \varphi_{f_{n-1}} \left( \cdots \left( \varphi_{f_1} \left( \varphi^G_m(G) \right) \right) \cdots \right) \subset \iota_{f_n}(G_{f_n}).$$
Finally, we know that $\ell_{n+1}$ belongs to the link of the vertex-group $\varphi_{f_n} \left( \cdots \left( \varphi_{f_1} \left( \varphi_m^G(G) \right) \right) \cdots \right)$ of the graph product $G_{t(f_n)}$. Since the words $mf_1 \cdots f_n$ and $e_1 \cdots e_r f_n$ contain the same sequence of edges of $E$ up to cancellation in $\mathfrak{F}$, it follows that
\begin{equation}\label{eq:vertexgroups}
\varphi_{f_n} \left( \varphi_{e_r} \left( \cdots \left( \varphi_{e_1}(G) \right) \cdots \right) \right)= \varphi_{f_n} \left( \cdots \left( \varphi_{f_1} \left( \varphi_m^G(G) \right) \right) \cdots \right)
\end{equation}
so that $\left[ \iota_{\bar{f_n}}(c) \ell_{n+1} \right]$ belongs to the link of the vertex-group $\varphi_{f_n} \left( \varphi_{e_r} \left( \cdots \left( \varphi_{e_1}(G) \right) \cdots \right) \right)$ of the graph product $G_{t(f_n)}$. Thus, we have proved that $m \ell$ belongs to $\mathrm{link}(G)$.

\medskip \noindent
Now, suppose that the equalities $b = \emptyset$ and $f_n = \bar{e_r}$ both hold. According to Lemma \ref{lem:reducingproduct}, $a_1e_1 \cdots a_{r-1}e_{r-1} \left[ a_r \varphi_{f_n}(a_{r+1}) \right]$ is a reduced word representing $mpf_n$, so that, once again according to Lemma \ref{lem:reducingproduct}, 
$$a_1 \cdot e_1 \cdots a_{r-1} \cdot e_{r-1} \cdot \left[ a_r \varphi_{f_n}(a_{r+1})\ell_{n+1} \right]$$
is a reduced word representing $m \ell$. Notice that, because $mp$ belongs to $\mathrm{link}(G)$, necessarily $a_1e_1 \cdots a_{r-1} e_{r-1}$ belongs $\mathrm{link}(G)$, $a_r$ belongs to the link of the vertex-group $\varphi_{e_{r-1}} \left( \cdots \left( \varphi_{e_1}(G) \right) \cdots \right)$ of the graph product $G_{t(e_{r-1})}$, and $a_{r+1}$ belongs to the link of the vertex-group $\varphi_{e_r}\left( \cdots \left( \varphi_{e_1}(G) \right) \cdots \right)$ of the graph product $G_{t(e_r)}$. Because $f_n= \bar{e_r}$, we deduce that $\varphi_{f_n}(a_{r+1})$ belongs to the link of the vertex-group $\varphi_{e_{r-1}} \left( \cdots \left( \varphi_{e_1}(G) \right) \cdots \right)$ of the graph product $\varphi_{f_n}(G_{t(e_r)}) = G_{s(e_r)} = G_{t(e_{r-1})}$. Finally, because $\ell \in \mathrm{link} \left( \varphi_m^G(G) \right)$, we know that $\ell_{n+1}$ belongs to the vertex-group $\varphi_{f_n} \left( \cdots \left( \varphi_{f_1} \left( \varphi_m^G(G) \right) \right) \cdots \right)$ of the graph product $G_{t(f_n)}=G_{t( \bar{e_r})}=G_{s(e_r)}=G_{t(e_{r-1})}$, which is equal to $\varphi_{f_n} \left( \varphi_{e_r} \left( \cdots \left( \varphi_{e_1}(G) \right) \cdots \right) \right)$ according to the equality \ref{eq:vertexgroups}, but also to $\varphi_{e_{r-1}} \left( \cdots \left( \varphi_{e_1}(G) \right) \cdots \right)$ since $f_n= \bar{e_r}$. Thus, we have proved that the reduction $\left[ a_r \varphi_{f_n}(a_{r+1}) \ell_{n+1} \right]$ belongs the link of the vertex-group $\varphi_{e_{r-1}} \left( \cdots \left( \varphi_{e_1}(G) \right) \cdots \right)$ of the graph product $G_{t(e_{r-1})}$. Consequently, $m \ell \in \mathrm{link}(G)$. 
\end{proof}

\noindent
We conclude this section by proving the following statement about hyperplanes labelled by $E$.

\begin{lemma}\label{lem:RAGGhypE}
Let $J$ be a hyperplane labelled by $E$. Then $J$ has two fibers which are stabilised by $\mathrm{stab}(J)$. 
\end{lemma}

\begin{proof}
Let $(h,h \cdot e)$ be an edge dual to $J$, where $h \in \mathfrak{F}$ and $e \in E$. According to Lemma \ref{lem:RAGGclique}, this edge is a clique, so $J$ has only two fibers, say $\partial_-$ and $\partial_+$ containing respectively $h$ and $h \cdot e$. We first notice that

\begin{claim}
If $a \in \partial_-$ and $b \in \partial_+$ are two adjacent vertices, then $b=a \cdot e$.
\end{claim}

\noindent
Because the edge $(a,b)$ is dual to the same hyperplane as the edge $(h,he)$, and that this hyperplane has only two fibers, there must exist a sequence of edges 
$$(x_1,y_1)=(h,he), \ (x_2,y_2), \ldots, (x_{n-1},y_{n-1}), \ (x_n,y_n)=(a,b)$$
such that $(x_i,y_i)$ and $(x_{i+1},y_{i+1})$ are opposite edges of some square for every $1 \leq i \leq n-1$. We argue by induction on $n$. If $n=0$, there is nothing to prove. If $n \geq 1$, then our induction hypothesis implies that $y_{n-1}=x_{n-1} \cdot e$. Notice that, as a consequence of Fact \ref{fact:squareE}, the edge $(x_{n-1},x_n)$ cannot be labelled by $E$. Therefore, Lemma \ref{lem:RAGGsquare} applies, and we deduce that $y_n=x_n \cdot e$. This concludes the proof of our claim. 

\medskip \noindent
Now, suppose by contradiction that there exists some $g \in \mathrm{stab}(J)$ such that $g h \in \partial_+$ and $ghe \in \partial_-$. According to our previous claim, necessarily $gh=ghee$, hence $ee=1_v$, where $v$ is the ending point of $h$ thought of as an arrow of the groupoid $\mathfrak{F}$. But the normal form provided by Proposition \ref{prop:GFnormalform} asserts that this is not possible. 
\end{proof}

\subsection{When is the action topical-transitive?}\label{section:RAGGtopic}

\noindent
In order to apply our combination results, we need to show that $\mathfrak{F}_{\omega}$ acts topically-transitively on $\mathfrak{X}$. Unfortunately, this is not always the case; see Example \ref{ex:rtimesZ}. Proposition \ref{prop:RAGGtopical} below gives a sufficient condition for our action to be topical-transitive, which is essentially necessary as well. First, given a factor $G$, define
$$\Phi(G) = \{ \varphi_{\ell}^G \mid \ell \in \mathrm{link}(G), \ \varphi_{\ell}^G (G)=G \} \subset \mathrm{Aut}(G).$$

\begin{prop}\label{prop:RAGGtopical}
If $\Phi(G)= \{ \operatorname{Id} \}$ for every factor $G$, then the action $\mathfrak{F}_{\omega} \curvearrowright \mathfrak{X}$ is topical-transitive. 
\end{prop}

\begin{proof}
Let $C$ be a clique of $\mathfrak{X}$ labelled by a factor $G$ and let $J$ denote its dual hyperplane. According to Proposition \ref{prop:RAGGhyp}, there is a bijection between $G$ and the set $\mathcal{S}(J)$ of the fibers of $J$. According to Corollary \ref{cor:RAGGstabhyp}, if $k \in \mathrm{stab}(J)$ then $k=g \cdot m$ for some $g \in G$ and $m \in \mathrm{link}(G)$ such that $\varphi_m^G \in \Phi(G)$. For convenience, set $\varphi = \varphi_m^G$. By using Lemma \ref{lem:RAGGcommute}, we find
$$gm \cdot (h \cdot \mathrm{link}(G))= g \varphi^{-1}(h) \cdot m \cdot \mathrm{link}(G) = g \cdot (h \cdot \mathrm{link}(G)).$$
Thus, the action of $k=gm$ on $\mathcal{S}(J)$ corresponds to the action of $g$ on $G$ by left-multiplication. Therefore, there exists an element of $\mathrm{stab}(C)$, namely $g$, which induces the same permutation of the fibers of $J$ as $k$. Equivalently, $\mathrm{Im}(\rho_C) \subset \rho_C(\mathrm{stab}(C))$. It follows from Lemma \ref{lem:RAGGhypE} that the same inclusion holds if $C$ is a clique labelled by $E$.

\medskip \noindent
Next, if $C$ is a clique labelled by $E$, then $C$ is finite and it follows from Lemma \ref{lem:RAGGhypE} that $\mathrm{stab}(C) = \mathrm{fix}(C)$ (which is trivial). Otherwise, if $C$ is labelled by some factor $G$, then the action $\mathrm{stab}(C) \curvearrowright C$ is isomorphic to the action $G \curvearrowright G$ by left-multiplication. A fortiori, $\mathrm{stab}(C) \curvearrowright C$ is transitive and free on the vertices. 
\end{proof}

\noindent
In order to apply our combination results, we register the following easy observations:

\begin{lemma}\label{lem:RAGGmisc}
The following statements hold.
\begin{itemize}
	\item[(i)] If the underlying abstract graph is locally finite and if the simplicial graphs defining the graph products are finite, then any vertex of $\mathfrak{X}$ belongs to finitely many cliques.
	\item[(ii)] Vertex-stabilisers of $\mathfrak{X}$ are trivial.
	\item[(iii)] If the underlying abstract graph and the simplicial graphs defining the graph products are all finite, then $\mathfrak{X}$ contains finitely many $\mathfrak{F}_{\omega}$-orbits of cliques.
	\item[(iv)] If the underlying abstract graph and the simplicial graphs defining the graph products are all finite, then the cubical dimension of $\mathfrak{X}$ is finite.
	\item[(v)] For every prism $P=C_1 \times \cdots \times C_n$, $\mathrm{stab}(P)= \mathrm{stab}(C_1) \times \cdots \times \mathrm{stab}(C_n)$.
\end{itemize}
\end{lemma}

\begin{proof}
Let $g \in \mathfrak{X}$ be a vertex. As an element of $\mathfrak{F}$, $g$ is an arrow from $\omega$ to some vertex $v \in V$. If the underlying abstract graph of our right-angled graph of groups is locally finite, so that $v$ has finitely many neighbors, then only finitely many cliques labelled by $E$ contain $g$. Moreover, because the cliques labelled by $V$ which contain $g$ are the $gG$'s where $G$ is any vertex-group of the graph product $G_v$, Point $(i)$ follows.

\medskip \noindent
Point $(ii)$ is clear since $\mathfrak{X}$ is a (connected component of a) Cayley graph of $\mathfrak{F}$. 

\medskip \noindent
Suppose that the underlying abstract graph of our right-angled graph of groups is finite. As a consequence, fixing for every $v \in V$ some arrow $g_v \in \mathfrak{F}$ starting from $\omega$ and ending at $v$, the set $R= \{ g_v \mid v \in V \}$ is finite. Now, let $C$ be a clique of $\mathfrak{X}$. According to Lemma \ref{lem:RAGGclique}, if $C$ is labelled by $E$, it is an edge, say $C=(g,ge)$ where $g \in \mathfrak{X}$ and $e \in E$. Say that, as an element of $\mathfrak{F}$, $g$ is an arrow from $\omega$ to some $v \in V$. Then $g_vg^{-1} \cdot C = (g_v,g_ve)$ where $g_vg^{-1} \in \mathfrak{F}_{\omega}$. Thus, we have proved that the edges of the form $(g_u,g_uf)$, where $u \in V$ and $f \in E$, define a finite set of representatives of cliques labelled by $E$. Next, suppose that $C$ is labelled by some factor $G$, ie., $gG$ for some $g \in \mathfrak{F}$. Similarly, if $g$ is an arrow from $\omega$ to some vertex $v \in G$, then $g_vg^{-1} \cdot C = g_vG$ where $g_vg^{-1}$. Therefore, the cliques $g_uG$ define a set of representatives of the cliques labelled by $V$. Point $(iii)$ follows.

\medskip \noindent
According to Lemma \ref{lem:RAGGprism}, the cubical dimension of $\mathfrak{X}$ is at most $1+ \sup\limits_{\Lambda \ \text{leaf}} \dim_{\square} \Lambda$. So Point $(iv)$ follows from Lemma \ref{lem:GPleaf} and Corollary \ref{cor:Xprisms}.

\medskip \noindent
Finally, let $P$ be a prism. If $P$ is included into a leaf, the conclusion of Point $(v)$ follows from Lemma \ref{lem:GPleaf} and Corollary \ref{cor:Xprisms}. Otherwise, according to Lemma \ref{lem:RAGGprism}, $P=P_1 \times P_2$ where $P_2$ is an edge labelled by $E$ and $P_1$ a prism included into a leaf. We deduce from Lemma \ref{lem:RAGGhypE} that $\mathrm{stab}(P)= \mathrm{stab}(P_1) \times \mathrm{stab}(P_2)$, and in fact $\mathrm{stab}(P)= \mathrm{stab}(P_1)$ since $\mathrm{stab}(P_2)$ is trivial. Therefore, the conclusion follows from Lemma \ref{lem:GPleaf} and Corollary \ref{cor:Xprisms}, proving Point $(v)$.
\end{proof}

\noindent
Now, we are ready to apply the criteria proved in Sections \ref{section:topicalactionsI} and \ref{section:topicalactionsII}. First, we deduce from Proposition \ref{prop:properlydiscontinuous}:

\begin{thm}
Let $\mathfrak{G}$ be a right-angled graph of groups such that $\Phi(G)= \{ \operatorname{Id} \}$ for every factor $G$. If the factors act properly on CAT(0) cube complexes, then so does the fundamental group of $\mathfrak{G}$. 
\end{thm}

\noindent
Next, it follows from Propositions \ref{prop:CAT0metricallyproper}, \ref{aTmenablegroups} and \ref{aBmenablegroups}:

\begin{thm}
Let $\mathfrak{G}$ be a right-angled graph of groups such that $\Phi(G)= \{ \operatorname{Id} \}$ for every factor $G$. Suppose that the underlying abstract graph is locally finite and that the simplicial graphs defining the graph products are finite. 
\begin{itemize}
	\item If the factors act metrically properly on CAT(0) cube complexes, then so does the fundamental group of $\mathfrak{G}$;
	\item if the factors are a-T-menable, then so is the fundamental group of $\mathfrak{G}$;
	\item if the factors are a-$L^p$-menable where $p \notin 2 \mathbb{Z}$, then so is the fundamental group of $\mathfrak{G}$.
\end{itemize} 
\end{thm}

\noindent
Finally, we deduce from Proposition \ref{prop:cubulatinggeometrically} and Theorem \ref{thm:producingCAT0groups}:

\begin{thm}\label{thm:RAGGgeom}
Let $\mathfrak{G}$ be a right-angled graph of groups such that $\Phi(G)= \{ \operatorname{Id} \}$ for every factor $G$. Suppose that the underlying abstract graph and the simplicial graphs defining the graph products are all finite. 
\begin{itemize}
	\item If the factors act geometrically on CAT(0) cube complexes, then so does the fundamental group of $\mathfrak{G}$;
	\item if the factors are CAT(0), then so is the fundamental group of $\mathfrak{G}$.
\end{itemize} 
\end{thm}

\noindent
We conclude this section by estimating the equivariant $\ell^p$-compressions of fundamental groups of right-angled graph of groups. 

\begin{thm}
Let $\mathfrak{G}$ be a right-angled graph of groups such that $\Phi(G)= \{ \operatorname{Id} \}$ for every factor $G$. Suppose that the underlying abstract graph and the simplicial graphs defining the graph products are all finite, and that the factors are finitely generated. For every $p \geq 1$, the fundamental group $\mathfrak{F}_{\omega}$ of $\mathfrak{G}$ satifies
$$\alpha_p^*(\mathfrak{F}_{\omega}) \geq \min \left( \frac{1}{p}, \min\limits_{\text{$G$ factor}} \alpha_p^*(G) \right).$$
\end{thm}

\begin{proof}
It follows from \cite[Proposition 11.2.4]{automaticgroups} that the canonical map $\mathfrak{F}_{\omega} \hookrightarrow \mathfrak{X}$ (ie., the orbit map associated to the basepoint $1_{\omega}$) is a quasi-isometric embedding. Therefore, the conclusion follows from Proposition \ref{prop:equicompression}. 
\end{proof}

\begin{remark}
It is worth noticing that our inequality turns out to be an equality if $p \geq 2$ and if $\mathfrak{F}_{\omega}$ contains a quasi-isometrically embedded non abelian free subgroup; see Remark \ref{rem:compressionGP}. By noticing that there is a natural epimorphism $\mathfrak{F}_{\omega} \twoheadrightarrow \pi_1(\Gamma)$, where $\Gamma$ denotes the underlying abstract graph of our right-angled graph of groups, we deduce that there exists such a subgroup if the Euler characteristic of $\Gamma$ is negative.
\end{remark}

\subsection{Hyperbolicity}

\noindent
In this section, we determine precisely when the fundamental groups of our right-angled graphs of groups are hyperbolic. More precisely, we prove:

\begin{thm}\label{thm:RAGGhyp}
Let $\mathfrak{G}$ be a right-angled graph of groups such that $\Phi(G)= \{ \operatorname{Id} \}$ for every factor $G$. Suppose that the factors are finitely generated, and that the underlying abstract graph and the simplicial graphs defining the graph products are all finite. The fundamental group of $\mathfrak{G}$ is hyperbolic if and only if the following conditions are satisfied:
\begin{itemize}
	\item any simplicial graph defining one of our graph products is square-free;
	\item  there do not exist a loop $p$ in the abstract graph, based at some vertex $v \in V$, and two non adjacent vertex-groups $G,H$ in the graph product $G_v$, such that $\varphi_p^G(G)=G$ and $\varphi_p^H(H)=H$;
	\item the link of every factor is finite;
	\item the factors are hyperbolic.
\end{itemize}
\end{thm}

\noindent
In fact, in order to prove the reciprocal, we will show a sufficient condition for being relatively hyperbolic by applying Theorem \ref{thm:qmrelativelyhyp}. More precisely, we prove the following criterion:

\begin{prop}\label{prop:RAGGrelhyp}
Let $\mathfrak{G}$ be a right-angled graph of groups such that $\Phi(G)= \{ \operatorname{Id} \}$ for every factor $G$. Suppose that the underlying abstract graph and the simplicial graphs defining the graph products are all finite, and that the following conditions are satisfied:
\begin{itemize}
	\item any simplicial graph defining one of our graph products is square-free;
	\item  there do not exist a loop $p$ in the abstract graph, based at some vertex $v \in V$, and two non adjacent vertex-groups $G,H$ in the graph product $G_v$, such that $\varphi_p^G(G)=G$ and $\varphi_p^H(H)=H$;
	\item the link of every infinite factor is finite.
\end{itemize}
The fundamental group of $\mathfrak{G}$ is hyperbolic relatively to a finite collection of groups commensurable to factors. 
\end{prop}

\noindent
The first step toward the proof of this proposition is to determine when our quasi-median graph $\mathfrak{X}$ is hyperbolic. This is the purpose of the next proposition. 

\begin{prop}\label{prop:XXhyp}
Suppose that the underlying abstract graph and the simplicial graphs defining the graph products are all finite. Suppose also that there exists an integer $N \geq 1$ such that, for every factor $G$, the automorphisms of $\Phi(G)$ have order at most $N$. The quasi-median graph $\mathfrak{X}$ is not hyperbolic if and only if one the following condition is satisfied:
\begin{itemize}
	\item one of the simplicial graphs defining the graph products contains an induced square;
	\item there exist a loop $p$ in the abstract graph, based at some vertex $v \in V$, and two non adjacent vertex-groups $G,H$ in the graph product $G_v$, such that $\varphi_p^G(G)=G$ and $\varphi_p^H(H)=H$.
\end{itemize}
\end{prop}

\begin{proof}
Suppose that the simplicial graphs defining the graph products are square-free. As a consequence, it follows from Lemma \ref{lem:GPleaf}, Fact \ref{fact:hypiffsquarefree} and Proposition \ref{prop:qmhyp} that there exists a constant $K \geq 1$ such that the flat rectangles of any leaf of $\mathfrak{X}$ must be $K$-thin. Set 
$$L = 1+ \max \left( B, (A^2+1)(K+1) \cdot \# V \right)$$
where $A$ is the maximum number of vertices of the simplicial graph associated to one of the graph products of our right-angle graph of groups, and $B$ the maximal number of vertices of a clique in one of these graphs. If $\mathfrak{X}$ is not hyperbolic, it follows from Proposition \ref{prop:qmhyp} that $\mathfrak{X}$ contains a flat rectangle $R : [0,L] \times [0,L] \hookrightarrow \mathfrak{X}$. Because two transverse hyperplanes cannot be both labelled by an edge of $E$, as a consequence of Fact \ref{fact:squareE}, suppose without loss of generality that the edges of $\{0 \} \times [0,L]$ are not labelled by $E$. Fixing some $0 \leq a \leq b \leq L$, if the edges of $[a,b] \times \{ 0 \}$ are all labelled by $V$, then the subrectangle $[a,b] \times [0,L]$ must be included into some leaf $\Lambda$. By definition of $K$, we deduce that $b-a \leq K$. Thus, $[0,L] \times \{ 0 \}$ contains at most $K$ consecutive edges labelled by $V$, so that it must contain $m \geq (A^2+1) \cdot \# V$ edges labelled by $E$. 

\medskip \noindent
Up to a translating by an element of $\mathfrak{F}$, suppose for convenience that $(0,0)=1_{\omega}$. Write $(0,L)$ as $g_1 \cdots g_n$ where $g_i \in G_{\omega}$ for every $1 \leq i \leq n$. As a consequence of Corollary \ref{cor:RAGGsquare}, if we write $(L,0)$ as a reduced word $h_1e_1 \cdots h_me_mh_{m+1}$, then the vertices of $R$ are
$$g_1 \cdots g_i \cdot h_1e_1 \cdots h_je_j \ \text{and} \ g_1 \cdots g_i \cdot h_1e_1 \cdots h_je_jh_{j+1}, \ 0 \leq i \leq n, 0 \leq j \leq m.$$
Notice that the vertices
$$g_1 \cdots g_i \cdot e_1 \cdots e_j, \ 0 \leq i \leq n, 0 \leq j \leq m$$
define a new flat rectangle. Because $m \geq (A^2+1) \cdot \# V$, the path $e_1,\ldots, e_m$ must contain a loop passing through its basis at least $A^2+1$ times. For convenience, we will supose that $e_1, \ldots, e_m$ is such a loop. 

\medskip \noindent
Because $L > B$, there must be two syllables of $g_1 \cdots g_n$ which belong to two non adjacent vertex-groups $G,H$ of the graph product $G_{\omega}$. By applying Lemma \ref{lem:RAGGsquare}, we deduce that $e_1, \ldots, e_m \in \mathrm{link}(G) \cap \mathrm{link}(H)$. Now, the collection $\left( (\varphi_{e_1\cdots e_i}^G(G), \varphi_{e_1 \cdots e_i}^H(H) ) \right)$ contains at least $A^2+1$ pairs of non adjacent vertex-groups of $G_v$, so there exist $1 \leq r < s \leq m$ such that $\varphi_{e_1 \cdots e_r}^G(G)= \varphi_{e_1 \cdots e_s}^G(G)$ and $\varphi_{e_1 \cdots e_r}^H(H)= \varphi_{e_1 \cdots e_s}^H(H)$, and such that $e_1 \cdots e_r$, $e_{r+1} \cdots e_s$ and $e_{s+1} \cdots e_m$ are loops based at $\omega$. 

\medskip \noindent
Thus, we find that the second condition of the statement of our proposition is satisfied for the loop $e_{r+1}, \ldots, e_s$ and the vertex-groups $\varphi_{e_1 \cdots e_r}^G(G)$, $\varphi_{e_1 \cdots e_r}^H(H)$. 

\medskip \noindent
Conversely, suppose that one the two conditions given in the statement of our proposition is satisfied. If there exists a simplicial graph defining one of the graph products of our right-angled graph of groups, we deduce from Lemma \ref{lem:GPleaf} and Fact \ref{fact:hypiffsquarefree} that $\mathfrak{X}$ is not hyperbolic. Next, suppose that there exist a loop $p$ in the abstract graph underlying our right-angled graph of groups, based at some vertex $v \in V$, and two non adjacent vertex-groups $G,H$ in the graph product $G_v$, such that $\varphi_p^G(G)=G$ and $\varphi_p^H(H)=H$. Fix an arrow $k \in \mathfrak{F}$ starting from $\omega$ and ending at $v$, and two non trivial elements $g \in G$ and $h \in H$. If $n \geq 1$ is any multiple of $N!$ and if $\phi,\psi$ denotes $ \left( \varphi_p^G \right)^{-1}$, $\left( \varphi_p^H \right)^{-1}$ respectively, then 
$$d \left( k(gh)^n, kp^n \right) + d \left( k, kp^n(gh)^n \right) = \left| p^{-n} (gh)^n \right| + \left| p^n(gh)^n \right|= 6n,$$
but
$$d \left( k, k(gh)^n \right) + d \left( kp^n, kp^n(gh)^n \right) = 2 \cdot \left|(gh)^n \right| = 4n$$
and
$$d \left( k,kp^n \right) + d \left( k(gh)^n,kp^n(gh)^n \right) = |p^n| + \left|p^n \left( \phi^n(g) \psi^n(h) \right)^{-n} (gh)^n \right|= 2n$$
since $\phi^n(g)=g$ and $\psi^n(h)=h$. Therefore, if $\mathfrak{X}$ is $\delta$-hyperbolic for some $\delta \geq 0$, necessarily $6n \leq \max(4n,2n) +2 \delta$, hence $\delta \geq n$. Because we can choose $n$ arbitrarily large, we conclude that $\mathfrak{X}$ is not hyperbolic.
\end{proof}

\begin{proof}[Proof of Proposition \ref{prop:RAGGrelhyp}.]
Under our assumptions, it follows from Proposition \ref{prop:XXhyp} that $\mathfrak{X}$ is hyperbolic. Moreover, the three first points of Theorem \ref{thm:qmrelativelyhyp} are satisfied according to Lemma \ref{lem:RAGGmisc}; its two last points are equivalent to the following statement: any fiber of a hyperplane dual to an infinite clique is finite. But $\mathfrak{X}$ satisfies this statement according to Proposition \ref{prop:RAGGhyp} and our assumptions. Therefore, Theorem \ref{thm:qmrelativelyhyp} applies, so that the fundamental group of $\mathfrak{G}$ is hyperbolic relatively to subgroups commensurable to clique-stabilisers of $\mathcal{C}$. The conclusion follows. 
\end{proof}

\begin{proof}[Proof of Theorem \ref{thm:RAGGhyp}.]
If all the conditions of our statement are satisfied, it follows from Proposition \ref{prop:RAGGrelhyp} that the fundamental group $\mathfrak{F}_{\omega}$ of $\mathfrak{G}$ is hyperbolic relatively to hyperbolic groups, so that it must be hyperbolic. Conversely, suppose that $\mathfrak{F}_{\omega}$ is hyperbolic. It follows from the normal form provided by Proposition \ref{prop:GFnormalform} that the factors and the graph products of our right-angled graph of groups are quasi-isometrically embedded, so that they must be hyperbolic. As a consequence, we deduce from Theorem \ref{Meier} that the simplicial graph associated to these graph products must be square-free. 

\medskip \noindent
Assume that there exist a loop $p$ in the abstract graph of our right-angled graph of groups, based at some vertex $v \in V$, and two non adjacent vertex-groups $G,H$ in the graph product $G_v$, such that $\varphi_p^G(G)=G$ and $\varphi_p^H(H)=H$. Fix two non trivial elements $g \in G$ and $h \in H$. Then $p$ and $gh$ define two infinite-order elements of the fundamental group $\mathfrak{F}_v$ based at $v$, and they commute. Therefore, $\mathfrak{F}_v$ contains a subgroup isomorphic to $\mathbb{Z}^2$, and a fortiori it cannot be hyperbolic. The same conclusion holds for $\mathfrak{F}_{\omega}$ since $\mathfrak{F}_{\omega}$ and $\mathfrak{F}_v$ are isomorphic.

\medskip \noindent
Suppose that there exists a factor $G$ whose link is infinite. We want to prove that, under this assumption, the fundamental group of $\mathfrak{G}$ is not hyperbolic. This will conclude the proof of our theorem. We distinguish two cases. 

\medskip \noindent
First, supose that the length of normal words is bounded on $\mathrm{link}(G)$. Because the abstract underlying our right-angled graph of groups is finite, there must exist edges $e_1, \ldots, e_n \in E$ such that infinitely many elements of $\mathrm{link}(G)$ can be written as a reduced word $\ell_1e_1 \cdots \ell_ne_n \ell_{n+1}$. Let $1 \leq i \leq n+1$ be an index such that there are infinitely many choices for $\ell_i$. Up to shortening our initial collection of edges, say that $i=n+1$. So there exist an infinite collection of elements $g_1, g_2, \ldots \in G_{t(e_n)}$ such that $e_1 \cdots e_n g_i \in \mathrm{link}(G)$ for every $i \geq 1$; in particular, the $g_i$'s belong to the link of the vertex-group $\varphi_{e_n}\left( \cdots \left( \varphi_{e_1}(G) \right) \cdots \right)$ of the graph product $G_{t(e_n)}$. Thus, we have found that some graph product of our right-angled graph of groups contains an infinite vertex-groups such that the subgroup generated by its link is infinite. Because a hyperbolic group does not contain two commuting infinite subgroups intersecting trivially, we deduce that $\mathfrak{X}$ is not hyperbolic. 

\medskip \noindent
Next, suppose that the length of normal words is unbounded on $\mathrm{link}(G)$. As a consequence, fixing some integer $m \geq 1$, there exists some $\ell \in \mathrm{link}(G)$ which can be written as a normal word $\ell_1e_1 \cdots \ell_n e_n \ell_{n+1}$ with $n \geq mA \cdot \# V$, where $A$ denotes the maximal number of vertices of a simplicial graph of our right-angled graph of groups. The path $e_1, \ldots, e_n$ must contain a loop, say $e_r, \ldots, e_s$ for some $1 \leq r<s \leq n$, which passes through its basepoint at least $mA$ times. Up to replacing $G$ with $\varphi_{e_1 \cdots e_r}^G(G)$, suppose that $r=1$. Since a graph product of our right-angled graph of groups contains at most $A$ vertex-groups, it follows that $\varphi_{e_i \cdots e_j}^G(G)=G$ for some subloop $e_i, \ldots, e_j$ with the same basepoint as $e_1, \ldots, e_s$ and such that $j-i \geq m$. For convenience, we will suppose that $i=1$; otherwise, replace $G$ with $\varphi_{e_1 \cdots e_i}^G(G)$. Set $p = \ell_1 e_1 \cdots \ell_j e_j \ell_{j+1}$. For every $g \in G$, 
$$g \cdot p = p \cdot \varphi_p^G(g) = p \cdot \varphi_{e_1 \cdots e_j}^G(g)= p \cdot g.$$
Thus, we have proved that, for every $m \geq 1$, there exist a vertex $v \in V$, an infinite factor $H \subset G_v$ and an element $p \in \mathfrak{F}_v$ belongs to the centralizer of $H$ and whose normal form has length at least $m$. Since there exist only finitely many factors, it follows that there exist a vertex $v \in V$ and an infinite factor $H$ such that, for every $m \geq 1$, there exists some $p \in \mathfrak{F}_v$ in the centralizer of $H$ whose normal form has length at least $m$. If $H$ is a torsion subgroup, we already know that $\mathfrak{F}_v$ cannot be hyperbolic since a hyperbolic group does not contain an infinite torsion subgroup. So suppose that $H$ contains an infinite order element $h$. Since the length of normal words is bounded on a coset of $\langle h \rangle$, we deduce from the existence of the previous infinite collection of elements that the centralizer of $h$ cannot be covered by finitely many cosets of $\langle h \rangle$. This means that $\langle h \rangle$ has not finite-index in its centralizer. As a consequence, $\mathfrak{F}_v$ cannot be hyperbolic. 
\end{proof}

\subsection{Examples}\label{section:RAGGex}

\noindent
This section is dedicated to explicit families of right-angled graphs of groups. 

\begin{ex}\label{ex:rtimesZ}
Fix a group $G$ and an automorphism $\varphi \in \mathrm{Aut}(G)$. The semi-direct product $H= G \rtimes_{\varphi} \mathbb{Z}$ is also an HNN extension, so that it turns out to be an example of a right-angled graph of groups. In this case, $\mathfrak{X}$ is the Cayley graph of $G \rtimes_{\varphi} \mathbb{Z}$ with respect to the generating set $G \cup \{ t \}$, where $t$ is a generator of $\mathbb{Z}$. As a consequence, $\mathfrak{X}$ is isometric to the product of a bi-infinite line (corresponding to the powers of $t$) with a clique (corresponding to the cosets of $G$). Let $J$ denote the hyperplane containing the cosets of $G$. The fibers of $J$ are naturally labelled by $G$, and the action of $t^n \in H$ on the set of fibers of $J$ corresponds to the action $g \mapsto \varphi^n(g)$ on $G$. By noticing that the latter action fixes the fiber labelled by the identify element, we deduce that no non trivial element of (the stabiliser of the clique) $G$ may induce the same permutation on the fibers. Consequently, $H$ acts topically-transitively on $\mathfrak{X}$ if and only if $\varphi = \operatorname{Id}$; in this case, $H$ is just the direct product $G \times \mathbb{Z}$.

\medskip \noindent
Thus, our method does not directly apply to a large class of right-angled graphs of groups. But there is a true obstruction. For instance, Gersten noticed in \cite{GerstenAutFn} that the free-by-cyclic group
$$\langle a,b,c,t \mid tat^{-1}=a, \ tbt^{-1}=ba, \ tct^{-1}=ca^2 \rangle$$
cannot be isomorphic to a subgroup of any CAT(0) group. 
\end{ex}

\noindent
A natural class of right-angled graphs of groups satisfying our extra-condition ensuring topical-transitivity of the action is obtained by gluing direct products along factors.

\begin{definition}
A \emph{Cartesian graphs of groups}\index{Cartesian graphs of groups} is a graph of groups such that groups are direct products, edges are factors of vertices, and monomorphisms are canonical embeddings.
\end{definition}

\begin{ex}
Let $A_1,A_2$ be two copies of a group $A$. Consider the HNN extension $(A_1 \times A_2) \underset{A}{\ast}$ associated to the monomorphisms $A \to A_1$ and $A \to A_2$. This is a Cartesian graph of groups, and its fundamental group, which we denote by $A^{\rtimes}$, admits
$$\langle A, \ t \mid [ a,a^t ]=1, \ a \in A \rangle$$
as a (relative) presentation. Notice that, if $A$ is infinite cyclic, we recover the group introduced in \cite{BKS}, which was the first example of fundamental group of a 3-manifold which is not subgroup separable. Following the proof of Fact \ref{fact:isometricX}, it is worth noticing that the quasi-median graph associated to $A^{\rtimes}$ depends only on the cardinality of $A$. In fact, we claim that the arguments used in the proof of Theorem \ref{thm:GPlip} apply similarly, so that, if $A$ and $B$ are two infinite Lipschitz-equivalent groups, then $A^{\rtimes}$ and $B^{\rtimes}$ are Lipschitz-equivalent. As a concrete application, we can say for instance that the groups
$$\mathbb{F}_2^{\rtimes} = \langle a,b,t \mid [a,a^t] = [b,b^t]=1 \rangle$$
and
$$\mathbb{F}_3^{\rtimes} = \langle a,b,c,t \mid [a,a^t]=[b,b^t]=[c,c^t]=1 \rangle$$
are two Lipschitz-equivalent groups acting geometrically on two-dimensional CAT(0) cube complexes, and that 
$$\mathbb{Z}_n^{\rtimes}= \langle x,t \mid x^n=[x,x^t]=1 \rangle$$
is a cubulable hyperbolic group for every $n \geq 2$. (In fact, $A^{\rtimes}$ is (relatively) hyperbolic if and only if $A$ is finite.)

\medskip \noindent
This example can be generalized in the following way. Let $P_n$ be a segment with $n+1$ vertices $v_0, \ldots, v_n$. Given a group $A$, denote by $P_nA$ the graph product associated to $P_n$ where each vertex-group is isomorphic to $A$. Now, consider the graph of groups which a bouquet of $n$ circles, where the vertex is $P_nA$, the edges are $A$, and the monomorphism on the $i$-th edge identifies $A_0$ with $A_i$ for every $1 \leq i \leq n$. This is a right-angled graph of groups such that $\Phi(G) = \{ \operatorname{Id} \}$ for every factor, and its fundamental group admits
$$\langle A,t_1, \ldots, t_n \mid [a,a^{t_i}]=1, \ a \in A, 1 \leq i \leq n \rangle$$
as a (relative) presentation.
\end{ex}

\begin{ex}
Let $A$ be a group and $B_1,B_2$ two copies of a group $B$. Consider the HNN extension $(B_1 \times A \times B_2) \underset{B}{\ast}$ associated to the monomorphisms $B \to B_1$ and $B \to B_2$. The group we obtain, which we denote by $A \Join B$, admits
$$\langle A,B,t \mid [a,b]=[a,b^t]=[b,b^t]=1, \ a \in A, b \in B \rangle $$
as a (relative) presentation. Notice that $A \Join B$ is relatively hyperbolic if and only if $B$ is finite (if so, $A \Join B$ is hyperbolic relatively to $A$). In particular, $A \Join B$ is hyperbolic if and only if $B$ is finite and $A$ hyperbolic. For instance,
$$\mathbb{F}_2 \Join \mathbb{Z}_3 = \langle a,b,c,t \mid [c,a]=[c,b]=[c,a^t]=[c,b^t]=[c,c^t]=c^3=1 \rangle$$
is a hyperbolic group acting geometrically on a two-dimensional CAT(0) cube complex. 
\end{ex}

\noindent
Another possible construction in order to obtain right-angled graph of groups satisfying our extra-condition is to take an HNN extension of a graph product conjugating two isomorphic vertex-groups. Interestingly, many diagram products arise in this way. 

\begin{ex}
Let $A$ be a group and $B_1,B_2$ two copies of a group $B$. Consider the HNN extension $(A \times( B_1 \ast B_2)) \underset{B}{\ast}$ associated to the monomorphisms $B \to B_1$ and $B \to B_2$. We recover the product $A \square B$ from Example \ref{ex:DPnonRAAG}, which admits
$$\langle A,B,t \mid [a,b]=[a,b^t]=1, \ a \in A, b \in B \rangle $$
as a (relative) presentation.

\medskip \noindent
Several generalisations of this construction are possible. For example, let $A$ be a group and, for every $1 \leq i \leq n$, let $B_i^1$ and $B_i^2$ be two copies of a common group $B_i$. Consider the graph product $G= A \times (B_1^1 \ast B_1^2 \ast \cdots \ast B_n^1 \ast B_n^2)$. Now, take the HNN extension $G \underset{B_1 \ast \cdots \ast B_n}{\ast}$ associated to the monomorphism $B_1 \ast \cdots \ast B_n \to B_1^1 \ast \cdots \ast B_n^1$ and $B_1 \ast \cdots \ast B_n \to B_{1}^2 \ast \cdots \ast B_{n}^2$. The group we obtain admits
$$\langle A,B_1, \ldots, B_n ,t \mid [a,b]=[a,b^t]=1, \ a \in A, b \in B_i, 1 \leq i \leq n \rangle$$
as a (relative) presentation. Another possibility would be to make $n$ distinct HNN extensions, to get the (relative) presentation
$$\langle A,B_1, \ldots, B_n, t_1, \ldots, t_n \mid [a,b]=[a,b^{t_i}]=1, \ a \in A, b \in B_i, 1 \leq i \leq n \rangle.$$
\end{ex}

\noindent
Now, suppose that we consider a right-angled graph of groups such that the underlying abstract graph is finite, as well as all the simplicial graphs defining the graph products, and such that all factors are finite. It follows from Lemma \ref{lem:RAGGmisc} that the action of the fundamental group on the associated quasi-median graph is geometric. Thus, although the action is not topical-transitive in general, we are able to deduce some information on the group. In fact, we suspect that the results proved in the previous two sections hold under the weaker hypothesis that $\Phi(G)$ is finite for every factor $G$; see Remark \ref{rem:transitivestrong}.

\begin{prop}\label{prop:RAGGfinitehyp}
Let $\mathfrak{G}$ be a right-angled graph of groups. Suppose that the underlying abstract graph and the simplicial graphs of the graph products are all finite, and that the factors are finite groups. The fundamental group of $\mathfrak{G}$ acts geometrically on some CAT(0) cube complex. Moreover, it is hyperbolic if and only if the following two conditions are satisfied:
\begin{itemize}
	\item the simplicial graphs defining the graph products are square-free;
	\item there do not exist a loop $p$ in the abstract graph, based at some vertex $v \in V$, and two non adjacent vertex-groups $G,H$ in the graph product $G_v$, such that $\varphi_p^G(G)=G$ and $\varphi_p^H(H)=H$.
\end{itemize}
\end{prop}

\noindent
The proposition follows directly from Propositions \ref{prop:quasimedianimplycubical} and \ref{prop:XXhyp}. 

\begin{ex}
Let $A_1,A_2$ be two copies of a finite abelian group $A$. Consider the graph of groups which is a bouquet of two circles such that the vertex is $A_1 \times A_2$, the two edges are $A$, the monomorphisms of one edge are the identities $A \to A_1$ and $A \to A_2$, and the monomorphisms of the other edge are the identity $A \to A_1$ and the inversion $A \to A_2$. For instance, if $A= \mathbb{Z}_3$, we get the group
$$\langle x,r,s \mid x^3=[x,x^r]=1, \ x^s=x^{-1} \rangle.$$
So this is a cubulable hyperbolic group. 

\medskip \noindent
Now, let $B$ be a finite group and $A_1,A_2$ two copies of a finite abelian group $A$. Consider the graph of groups which is a bouquet of two circles such that the vertex is $B \times (A_1 \ast A_2)$, the two edges are $A$, the monomorphisms of one edge are the identities $A \to A_1$ and $A \to A_2$, and the monomorphisms of the other edge are the identity $A \to A_1$ and the inversion $A \to A_2$. For instance, if $A=\mathbb{Z}_3$ and $B=\mathbb{Z}_2$, we get the group
$$\langle x,y,r,s \mid x^3=y^2=[y,x]=[y,x^r]=1, \ x^{rs}=x^{-1} \rangle.$$
This is a cubulable group which is not hyperbolic.
\end{ex}

\section{Open problems}\label{section:open}

\paragraph{Disc diagrams.} Section \ref{section:QMmain} contains essentially all the tools used for CAT(0) cube complexes, except one: disc diagrams. They were used independently in \cite{MR1347406} and in \cite{bigWise} (see also \cite{Hagenthesis, article3}). This is a very convenient tool, and it would be interesting to know whether a similar technology can be developped for quasi-median graphs. For this purpose, the local criterion formulated in Section \ref{section:locallyQM} could be useful. 

\begin{problem}
Study disc diagrams in quasi-median prism complexes.
\end{problem}

\noindent
One motivation is the following. It is known that a cube complex is a locally CAT(0) metric space if and only if it satisfies Gromov's link condition (see \cite{Leary}). Similarly, we showed in Section \ref{section:qmCAT0} that a quasi-median prism complex turns out to be CAT(0), so that a locally quasi-median prism complex must be locally CAT(0); but the converse does not hold. Indeed, two triangles glued along an edge defines a CAT(0) prism complex, but its 1-skeleton is not a quasi-median graph. We expect that, by restricting the gluings in the definition of a prism complex, the converse will hold. More precisely, define a prism complex as a collection of prisms glued together along ``factors'' (thinking of a prism as a product of simplexes). With respect to this definition, 

\begin{question}
Is a prism complex locally CAT(0) if and only if it is locally quasi-median?
\end{question}

\paragraph{Other properties to study.} In Sections \ref{section:topicalactionsI} and \ref{section:topicalactionsII}, we proved combination theorems for several properties, but many properties of interest were not considered. We list below some of them, for which similar arguments should be possible. Typically, it may be expected that a result which holds for graph products and groups acting on CAT(0) cube complexes generalises in the context of quasi-median graphs.

\begin{question}
Let $G$ be a group acting on a quasi-median graph such that the hypotheses of Proposition \ref{prop:cubulatinggeometrically} are satisfied. If clique-stabilisers are semihyperbolic (resp. (bi)automatic), is $G$ semihyperbolic (resp. (bi)automatic)?
\end{question}

\noindent
Extending the normal cube paths from \cite{NibloReeves} to quasi-median graphs, it is possible to choose a canonical chain of cliques between two vertices of a quasi-median graph whose cubical dimension is finite, and next, using Lemma \ref{lem:brokengeod}, this choice leads to a natural way to extend combings from cliques to a combing on the whole quasi-median graph (with respect to a global metric). Notice that combings on graph products were already studied in \cite{HermillerMeier}.

\begin{question}
Let $G$ be a group acting on a quasi-median graph such that the hypothesese of Proposition \ref{prop:cubulatinggeometrically} are satisfied. Does the inequality
$$\mathrm{asdim}(G) \leq \max \left\{ \mathrm{asdim}\left( \mathrm{stab}(P) \right) \mid \text{$P$ prism} \right\}$$
between asymptotic dimensions hold?
\end{question}

\noindent
Notice that the asymptotic dimension of graph products was studied in \cite{AntolinDreesen}, where the above inequality is proved. 

\begin{question}
Let $G$ be a group acting on a quasi-median graph such that the hypotheses of Proposition \ref{prop:cubulatinggeometrically} are satisfied. Does the inequality
$$\delta_G \prec \max \left\{ \delta_{\mathrm{stab}(C)} \mid \text{$C$ clique} \right\}$$
between Dehn functions hold?
\end{question}

\noindent
This inequality for graph products was proved in \cite{GPDehnFunction}.

\begin{question}
Let $G$ be a group acting topically-transitively on a quasi-median graph of finite cubical dimension with finite vertex-stabilisers. If prism-stabilisers satisfy the property RD, does $G$ satisfy property RD?
\end{question}

\noindent
We suspect that, under our assumptions, $G$ satifies the property RD relatively to prism-stabilisers (as defined by Sapir in \cite{SapirRD}), essentially as a consequence of Proposition \ref{prop:quasimedian}. It is worth noticing that this is precisely the meaning of \cite[Theorem 3.24]{SapirRD} in the context of graph products. 

\begin{question}
Let $G$ be a group acting on a quasi-median graph such that the hypotheses of Proposition \ref{prop:cubulatinggeometrically} are satisfied. If clique-stabilisers are coarse median, is $G$ coarse median as well?
\end{question}

\noindent
Coarse median groups were introduced by Bowditch in \cite{Coarsemedian}, including hyperbolic groups, groups acting geometrically on CAT(0) cube complexes, and mapping class groups of surfaces. It is even not known whether being coarse median is stable under graph products. Notice that there is a natural way to extend a collection of coarse medians defined on cliques into a global coarse median thanks to Proposition \ref{prop:quasimedian}. 

\medskip \noindent
The above properties are essentially metric, but similarly to the extensions of metrics and collections of walls, it is possible to define coherent system of topologies and to extend them as a global topology. For instance, we expect that extending contractible topologies leads to a global contractible topology (the idea is that it should be possible to collapse a hyperplane by identifying its fibers if its cliques are contractible). One example of possible application is (see also \cite{GPDehnFunction}):

\begin{question}
Let $\Gamma$ be a simplicial graph and $\mathcal{G}$ a collection of groups indexed by $V(\Gamma)$. Suppose that every group $G \in \mathcal{G}$ is of type $F_{n(G)}$ for some $n(G) \leq + \infty$. Is it true that the graph product $\Gamma \mathcal{G}$ is of type $F_n$ where
$$n= \max \left\{ \sum\limits_{i=1}^k n(G_i) \mid G_1, \ldots, G_k \ \text{pairwise adjacent} \right\}?$$
\end{question}

\paragraph{Some improvements.} We suspect that some of our results may be strengthen. For instance, by weighting the hyperplanes as in \cite{propertyA}, the lower bound in Proposition \ref{prop:compression} might be improved if the cubical dimension of the quasi-median graph is finite.

\begin{question}\label{question:compression}
Under the hypotheses of Proposition \ref{prop:compression}, is it true that
$$\alpha_p(X,\delta) \geq \inf\limits_{C~\text{clique}} \alpha_p(C,\delta_C)$$
if the cubical dimension of $X$ is finite?
\end{question}

\noindent
Theorem \ref{thm:GPlip} proved that Lipschitz equivalent vertex-groups lead to Lipschitz equivalent graph products. A natural problem would be to generalise the proof to a more general context. More precisley: 

\begin{question}
Let $G$ and $H$ be two groups acting topically-transitively on a common quasi-median graphs. If stabilisers of a given clique are quasi-isometric (resp. Lipschitz equivalent, strongly commensurable), when are $G_1$ and $G_2$ quasi-isometric (resp. Lipschitz equivalent, commensurable)?
\end{question}

\noindent
(The question about commensurability is motivated by \cite{GPcommensurability}, whose results are proved in the context of graph products.)

\paragraph{Curve graphs for graph products.} Let us continue the discussion we began in Section \ref{section:GPcurvegraph} about the analogy between graph products and mapping class groups. Although I am conviced that the technics which are mentionned below lead to proofs of the different statements, our discussion remains informal, and the statements need to be showed.

\medskip \noindent
First of all, inspired by Nielsen-Thurston classification of the homeomorphisms of a surface, we can classify the elements of a graph product from their actions on the associated quasi-median graph.

\begin{stat}
Let $\Gamma$ be a simplicial graph and $\mathcal{G}$ a collection of groups indexed by $V(\Gamma)$. If $g \in \Gamma \mathcal{G}$, three cases may happen:
\begin{itemize}
	\item $g$ is \emph{elliptic}, ie., it stabilises a prism of $\X$; if so, the support of $g$ is included into a clique.
	\item $g$ is \emph{reducible}, ie., it is loxodromic and it stabilises a product of two unbounded subcomplexes of $\X$; if so, the support of $g$ is included into a large join but not on a clique.
	\item $g$ is \emph{contracting}, ie., it is contracting isometry of $\X$; if so, the support of $g$ is not included into any large join.
\end{itemize}
\end{stat}

\noindent
This statement follows from the techniques we used in \cite{article3}. Given $\X$, it is natural to define the crossing graph $T(\Gamma, \mathcal{G})$ as the graph whose vertices are the hyperplanes of $\X$ and whose edges link two transverse hyperplanes, and to think of it as the curve graph. However, $T(\Gamma, \mathcal{G})$ is connected only if $\Gamma$ is connected itself. Since this disconnectedness is produced by cut vertices of $\X$, a solution is to add these cut vertices to $T(\Gamma, \mathcal{G})$ and to link them to their adjacent hyperplanes. Alternatively, it is possible to mimic the definition of the contact graph of a CAT(0) cube complex, as defined by Hagen in \cite{MR3217625}. (Notice that these two graphs are quasi-isometric.) But another problem is that contracting elements of $\Gamma \mathcal{G}$ may induce elliptic isometries of $T(\Gamma, \mathcal{G})$ as well as in any of its variations. In fact, if you define a new graph $\Lambda$ as the join of $\Gamma$ with a single vertex and a new collection $\mathcal{H}$ from $\mathcal{G}$ by labelling the vertex of $\Lambda \backslash \Gamma$ with $\mathbb{Z}_2$, then our different graphs turn out to be bounded although the graph product $\Gamma \mathcal{G}$ becomes $\Lambda \mathcal{H}= \Gamma \mathcal{G} \times \mathbb{Z}_2$. Nevertheless, we have:

\begin{stat}
Let $\Gamma$ be a connected simplicial graph and $\mathcal{G}$ a collection of groups indexed by $V(\Gamma)$. Then $T(\Gamma, \mathcal{G})$ is a quasi-tree, unbounded if and only if $\Gamma$ is not clique, on which $\Gamma \mathcal{G}$ acts (non uniformly) acylindrically. Moreover, elliptic and reducible elements of $\Gamma \mathcal{G}$ induce elliptic isometries of $T(\Gamma, \mathcal{G})$, and any element of $\Gamma \mathcal{G}$ inducing a loxodromic isometry of $T(\Gamma, \mathcal{G})$ must be contracting, the converse being true if the groups of $\mathcal{G}$ are all infinite.
\end{stat}

\noindent
We refer to \cite{MR3217625,article3,article4} for the techniques which can be used to prove this statement. As a consequence, although defining a reasonable curve graph in full generality might be more technical and less natural, $T(\Gamma, \mathcal{G})$ is a good candidate when $\Gamma$ is connected and the groups of $\mathcal{G}$ are infinite.

\paragraph{Quasi-isometric classification of graph products.} We expect that the study of the quasi-median graphs associated to graph products will lead to results towards the classification up to quasi-isometry of (some) graph products, in particular when vertex-groups are all finite. A first question is: is the quasi-isometric classification of graph products of finite groups more difficult that the same classification of right-angled Coxeter groups? More precisely:

\begin{question}
Does there exist a graph product of finite groups which is not quasi-isometric to a right-angled Coxeter groups?
\end{question}

\noindent
We suspect a positive answer. Of course, a complete classification is not expected. An interesting family of examples were studied in \cite{Bourdon1997} (whose arguments can be interpreted in a very natural way by using our quasi-median graphs!), where infinitely many pairwise non quasi-isometric hyperbolic graph products of finite cyclic groups are given. A sligthly more general problem would be:

\begin{problem}
Classify cycle of finite groups up to quasi-isometry. 
\end{problem}

\noindent
By a \emph{cycle of groups}, we mean the graph product associated to a cycle. In particular, if $F$ is a group and $n \geq 1$ an integer, we denote by $C_nF$ the graph product associated to a cycle of length $n$ whose vertex-groups are all isomorphic to $F$. 

\begin{question}\label{question:negativeanswerQI}
Let $n \geq 5$. Do there exist two infinite quasi-isometric groups $G,H$ such that the cycles $C_nG$ and $C_nH$ are not quasi-isometric?
\end{question}

\noindent
We expect a positive answer. A case of interest is when $G$ and $H$ are surface groups.

\paragraph{Wreath products.} Theorem \ref{thm:Wequicompression} provides a lower bound on the $\ell^p$-compression of some wreath products which depends on a constant $TS(X)$ associated to a CAT(0) cube complex $X$. It would be interesting to study the possible values of this constant. A (possibly naive) guess is the following:

\begin{question}
Let $X$ be a cocompact CAT(0) cube complex. Let $d$ denote the maximal integer such that $X$ contains an isometrically embedded subcomplex isomorphic to $\mathbb{R}^d$. Is it true that $TS(X)=d/(2d-1)$?
\end{question}

\noindent
Easier questions to test this problem may be the following:

\begin{question}
What is $TS(T_1 \times T_2)$ where $T_1,T_2$ are two trees?
\end{question}

\begin{question}
If $X$ is a finite dimensional CAT(0) cube complex, is is true that $TS(X)>1/2$?
\end{question}

\noindent
Corollary \ref{cor:Wproper} shows that acting properly on a CAT(0) cube complex is stable under wreath products, and according to Proposition \ref{prop:FZnfinitedimension}, the CAT(0) cube complex on which the wreath product acts properly can be chosen finite-dimensional for $F \wr \mathbb{Z}^n$, where $F$ is a finite group, whereas it is typically infinite-dimensional in the general case. 

\begin{question}
When does a wreath product act properly discontinuously on a finite-dimensional CAT(0) cube complex? 
\end{question}

\noindent
One obstruction is provided by the Tits alternative proved in \cite{MR2827012}. A wreath product acting properly on a finite-dimensional CAT(0) cube complex must either contain a non abelian free subgroup or be (locally finite)-by-abelian. For instance, $\mathbb{Z} \wr \mathbb{Z}$ cannot act properly on a finite-dimensional CAT(0) cube complex. Another obstruction is provided by the weak amenability. Indeed, it was noticed in \cite{CSVparticulier} that, for any non trivial finite group $F$, the wreath product $F \wr \mathbb{F}_2$ does not act properly on a finite dimensional CAT(0) cube complex, because a group admitting such an action must must have a Cowling-Haagerup constant equal to one according to \cite{WeakAmenabilityCCMizuta, WeakAmenabilityCCGuentner}, whereas it is proved in \cite[Corollary 2.11]{OzawaPopa} that the Cowling-Haagerup constant of $F \wr \mathbb{F}_2$ is strictly larger than one.

\paragraph{Diagram products.}  In Section \ref{section:appli2}, we were not able to study the stability of being virtually special under diagram products. Following \cite{arXiv:1507.01667}, it may be expected that, at least under the assumption that $[w]_{\mathcal{P}}$ is finite, the action of a diagram group on its quasi-median graph is special. So the main difficulty comes from the last point in the statement of Proposition \ref{prop:cocompactspecial}. Given some semigroup presentation $\mathcal{P}= \langle \Sigma \mid \mathcal{R} \rangle$, some collection of groups $\mathcal{G}$ indexed by $\Sigma$ and some base word $w \in \Sigma^+$, let $u \in \Sigma^+$ be a word equal to $w$ modulo $\mathcal{P}$ which we write as $x_1u_1 \cdots x_nu_nx_{n+1}$ for some $x_1,u_1, \ldots x_n,u_n,x_{n+1} \in \Sigma^+$. The map
$$(U_1, \ldots, U_n) \mapsto \Delta \cdot \left( \epsilon(x_1)+ U_1+ \cdots + \epsilon(x_n)+U_n + \epsilon(x_{n+1}) \right) \cdot \Delta^{-1},$$
where $\Delta$ is a fixed diagram satisfying $\mathrm{top}^-(\Delta)=w$ and $\mathrm{bot}^-(\Delta)=u$, defines a monomorphism $D(\mathcal{P}, \mathcal{G},u_1) \times \cdots \times D(\mathcal{P}, \mathcal{G},u_n) \hookrightarrow D(\mathcal{P},\mathcal{G},w)$. We call a group arising in this way a \emph{canonical subgroup}. 

\begin{question}
Is a canonical subgroup a rectract of the diagram product?
\end{question}

\noindent
In the context of cocompact diagram groups, we showed in \cite{arXiv:1507.01667} that canonical subgroups are virtual retracts. A related question, which may be easier to solve, is:

\begin{question}
If its factors and its underlying diagram group are residually finite, is a diagram product residually finite?
\end{question}

\noindent
In fact, the problem reduces to show that a diagram product of finite groups turns out to be residually finite if its underlying diagram group is itself residually finite. 

\medskip \noindent
In view of Corollary \ref{cor:RAAGDP}, a natural problem is to determine which right-angled Artin groups we are able to obtain in this way. More precisely,

\begin{question}
When is a graph a $\Gamma(\mathcal{P}, \sigma,w)$ for some semigroup presentation $\mathcal{P}= \langle \Sigma \mid \mathcal{R} \rangle$, some map $\sigma : \Sigma \to \{ 0 , 1 \}$ and some base word $w \in \Sigma^+$?
\end{question}

\noindent
Finally, notice that a consequence of Theorem \ref{thm:DPsemidirect}, combined with \cite{OrderGP} and \cite{WiestDiagram}, is that a diagram product of left-orderable groups must be left-orderable. Because diagram groups turn out to be bi-orderable \cite{MR2193191}, a natural question is:

\begin{question}
Is a diagram product of bi-orderable groups bi-orderable?
\end{question}

\noindent
In our opinion, the cocompactness assumption $[w]_{\mathcal{P}}<+ \infty$ in Theorem \ref{thm:DPhyp} can be removed by following the argument used in \cite{arXiv:1505.02053}. To be more precise, the analogue of \cite[Lemma 8]{arXiv:1505.02053} would be that a subgroup $H \leq D(\mathcal{P},\mathcal{G},w)$ either is containted into a prism-stabiliser or splits over a finite subgroup; as a consequence, if $D(\mathcal{P}, \mathcal{G},w)$ is finitely generated, it would decompose as a finite graph of groups whose edges are finite groups and whose vertices are direct products of finitely many groups of $\mathcal{G}$. However, it might be more interesting to study the problem with a more general point of view. Let us introduce the class of \emph{twisted graph products}\index{Twisted graph products}. 

\medskip \noindent
Let $\Gamma$ be a simplicial graph and $\mathcal{G}$ a collection of groups indexed by $V(\Gamma)$. Label the vertices of $\Gamma$ by the isomorphic classes of their associated groups. Let $\mathrm{Isom}^+(\Gamma)$ denote the subgroup of $\mathrm{Isom}(\Gamma)$ preserving this labelling. Fix a group for each of these classes, and, for every $s \in V(\Gamma)$, denote by $G(s)$ the group we fixed in the class of the vertex-group $G_s$, and fix an isomorphism $\varphi_s : G(s) \to G_s$. We deduce a monomorphism
$$\left\{ \begin{array}{ccc} \mathrm{Isom}^+(\Gamma) & \to & \mathrm{Aut}(\Gamma \mathcal{G}) \\ \varphi & \mapsto & \left( g \mapsto \varphi_{\varphi(s)} \circ \varphi_s^{-1}(g), \ g \in G_s \right) \end{array} \right.$$
For every $H \leq \mathrm{Isom}^+(\Gamma)$, we define the \emph{twisted graph product} $\Gamma \mathcal{G} \rtimes H$, where $H$ acts on $\Gamma \mathcal{G}$ via its image in $\mathrm{Aut}(\Gamma \mathcal{G})$ by the previous map. Notice that, up to isomorphism, the group we obtain does not depend on the choices of the isomorphisms $\varphi_s$. We recover the \emph{graph-wreath products} introduced in \cite{GraphWreathProducts} when $\mathcal{G}$ is reduced to a singleton. 

\begin{question}\label{question:twistedGPhyp}
When is a twisted graph products hyperbolic?
\end{question}

\noindent
In this context, Theorem \ref{thm:rotativeactions} shows that any group acting rotatively on a quasi-median graph decomposes as a twisted graph product. Does the converse hold? 

\begin{question}
Does a twisted graph product act naturally on a quasi-median graph?
\end{question}

\paragraph{Right-angled graphs of groups.} In Section \ref{section:appli4}, we were not able to study the stability of being virtually special under right-angled graphs of groups. The main difficulty comes from the last point in the statement of Proposition \ref{prop:cocompactspecial}.

\begin{question}
Consider a right-angled graph of groups satisfying the hypotheses of Theorem \ref{thm:RAGGgeom}. Are hyperplane-stabilisers retracts? 
\end{question}

\noindent
Another problem is that the action on the quasi-median graph is not always special. 

\begin{question}
When is the action of the fundamental group of a right-angled graph of groups on its associated quasi-median graph special?
\end{question}

\noindent
Motivated by Corollary \ref{cor:2complex} and Cartesian graphs of groups defined in Section \ref{section:RAGGex}, define a \emph{Cartesian cubing of groups} as a complex of groups such that its underlying complex is a non positively curved cube complex, its groups are direct products, and its monomorphisms are canonical embeddings of factors into products.

\begin{question}
Does the fundamental group of a Cartesian cubing of groups act topically-transitively on a quasi-median graph?
\end{question}

\addcontentsline{toc}{section}{References}

\bibliographystyle{alpha}
\bibliography{QMgraphs}

\newcommand{\etalchar}[1]{$^{#1}$}
\begin{thebibliography}{BCC{\etalchar{+}}13b}

\bibitem[AD13]{AntolinDreesen}
Y.~Antol\`in and D.~Dreesen.
\newblock The {H}aagerup property is stable under graph products.
\newblock {\em arXiv:1305.6748}, 2013.

\bibitem[Ago13]{MR3104553}
Ian Agol.
\newblock The virtual {H}aken conjecture.
\newblock {\em Doc. Math.}, 18:1045--1087, 2013.
\newblock With an appendix by Agol, Daniel Groves, and Jason Manning.

\bibitem[AGS06]{MR2271228}
G.~N. Arzhantseva, V.~S. Guba, and M.~V. Sapir.
\newblock Metrics on diagram groups and uniform embeddings in a {H}ilbert
  space.
\newblock {\em Comment. Math. Helv.}, 81(4):911--929, 2006.

\bibitem[Akh08]{PerturbationofWP}
Azer Akhmedov.
\newblock Perturbation of wreath products and quasi-isometric rigidity 1.
\newblock {\em International Mathematics Research Notices}, 2008:rnn004, 2008.

\bibitem[Alo96]{GPDehnFunction}
J.M Alonso.
\newblock Dehn functions and finiteness properties of graph products.
\newblock {\em Journal of Pure and Applied Algebra}, 107(1):9 -- 17, 1996.

\bibitem[AM15]{TitsAltGP}
Y.~Antol\`in and A.~Minasyan.
\newblock Tits alternatives for graph products.
\newblock {\em J. Reine Angew. Math.}, 704:55--83, 2015.

\bibitem[ANP09]{WPZwrZ}
T.~Austin, A.~Naor, and Y.~Peres.
\newblock The wreath product of $\mathbb{Z} \wr \mathbb{Z}$ has {H}ilbert
  compression $2/3$.
\newblock {\em Proc. Ameri. Math. Soc.}, 137:85--90, 2009.

\bibitem[Arn13]{labelledpartitions}
S.~Arnt.
\newblock Spaces with labelled partitions and isometric affine actions on
  {B}anach spaces.
\newblock {\em arXiv:1401.0125v2}, 2013.

\bibitem[Arz01]{GoulnaraQuasiconvex}
G.~N. Arzhantseva.
\newblock On quasiconvex subgroups of word hyperbolic groups.
\newblock {\em Geometriae Dedicata}, 87:191--208, 2001.

\bibitem[BC96]{BandeltChepoi}
H.-J. Bandelt and V.~Chepoi.
\newblock A {H}elly theorem in weakly modular space.
\newblock {\em Discrete Math.}, 160(1-3):25--39, 1996.

\bibitem[BC08]{BandeltChepoi2}
H.-J. Bandelt and V.~Chepoi.
\newblock The algebra of metric betweeness {II}: {G}eometry and equational
  characterization of weakly median graphs.
\newblock {\em European Journal of Combinatorics}, 29(3):676--700, 2008.

\bibitem[BCC{\etalchar{+}}13a]{bucolic}
B.~Bresar, J.~Chalopin, V.~Chepoi, T.~Gologranc, and D.~Osajda.
\newblock Bucolic complexes.
\newblock {\em Advances in Mathematics}, 243:127--167, 2013.

\bibitem[BCC{\etalchar{+}}13b]{RetractsChordal}
B.~Bre\v{s}ar, J.~Chalopin, V.~Chepoi, M.~Kov\v{s}e, A.~Labourel, and
  Y.~Vax\`es.
\newblock Retracts of products of chordal graphs.
\newblock {\em Journal of Graph Theory}, 73(2):161--180, 2013.

\bibitem[BH99]{MR1744486}
Martin~R. Bridson and Andr{\'e} Haefliger.
\newblock {\em Metric spaces of non-positive curvature}, volume 319 of {\em
  Grundlehren der Mathematischen Wissenschaften [Fundamental Principles of
  Mathematical Sciences]}.
\newblock Springer-Verlag, Berlin, 1999.

\bibitem[BHSC13]{BHSC}
J.~Behrstock, M.~Hagen, A.~Sisto, and P.-E. Caparace.
\newblock Thickness, relative hyperbolicity, and randomness in {C}oxeter
  groups.
\newblock {\em arXiv:1312.4789}, 2013.

\bibitem[BKS87]{BKS}
R.~G. Burns, A.~Karrass, and D.~Solitar.
\newblock A note on groups with separable finitely generated subgroups.
\newblock {\em Bulletin of the Australian Mathematical Society},
  36(1):153--160, 1987.

\bibitem[BM15]{GPAsymptotic}
G.~Bell and D.~Moran.
\newblock On constructions preserving the asymptotic topology of metric spaces.
\newblock {\em North Carolina Journal of Mathematics and Statistic}, 1:46--57,
  2015.

\bibitem[BMW94]{quasimedian}
H.-J. Bandelt, H.M. Mulder, and E.~Wilkeit.
\newblock Quasi-median graphs and algebras.
\newblock {\em J. Graph Theory}, 18(7):681--703, 1994.

\bibitem[Bou97]{Bourdon1997}
M.~Bourdon.
\newblock Immeubles hyperboliques, dimension conforme et rigidit{\'e} de
  {M}ostow.
\newblock {\em Geometric {\&} Functional Analysis GAFA}, 7(2):245--268, 1997.

\bibitem[Bow12]{relativelyhyperbolic}
Brian Bowditch.
\newblock Relatively hyperbolic groups.
\newblock {\em Internat. J. Algebra and Computation.}, 22(3), 2012.

\bibitem[Bow13]{Coarsemedian}
B.~Bowditch.
\newblock Coarse median spaces and groups.
\newblock {\em Pacific J. Math.}, 261(1):53--93, 2013.

\bibitem[Bow14]{Bowditchcriterion}
Brian Bowditch.
\newblock Uniform hyperbolicity of the curve graphs.
\newblock {\em Pacific J. Math.}, 269:269--280, 2014.

\bibitem[CCHO14]{weaklymoduloar}
J.~Chalopin, V.~Chepoi, H.~Hirai, and D.~Osajda.
\newblock Weakly modular graphs and nonpositive curvature.
\newblock {\em arXiv:1409.3892 (to appear in Memoirs of AMS)}, 2014.

\bibitem[CDE{\etalchar{+}}]{CDEHV}
V.~Chepoi, F.~Dragan, B.~Estellon, M.~Habib, and Y.~Vax\`es.
\newblock Diameters, centers, and approximating trees of $\delta$-hyperbolic
  geodesic spaces and graphs.
\newblock {\em Symposium on Computational Geometry, SoCG'2008}, pages 59--68.

\bibitem[CDH10]{medianviewpoint}
I.~Chatterji, C.~Dru{\c{t}}u, and F.~Haglund.
\newblock Kazhdan and {H}aagerup properties from the median viewpoint.
\newblock {\em Advances in Mathematics}, 225:882--921, 2010.

\bibitem[Che89]{ChepoiTriangles}
V.~Chepoi.
\newblock Classification of graphs by means of metric triangles.
\newblock {\em Metody Diskret. Analiz.}, 49(96):75--93, 1989.

\bibitem[Che00]{mediangraphs}
V.~Chepoi.
\newblock Graphs of some {$\rm CAT(0)$} complexes.
\newblock {\em Adv. in Appl. Math.}, 24(2):125--179, 2000.

\bibitem[Chi12]{OrderGP}
I.~M. Chiswell.
\newblock Ordering graph products of groups.
\newblock {\em International Journal of Algebra and Computation},
  22(04):1250037, 2012.

\bibitem[CMV04]{measuredwallspaces}
P.A. Cherix, F.~Martin, and A.~Valette.
\newblock Spaces with measured walls, the {H}aagerup property and property
  ({T}).
\newblock {\em Ergodic Theory Dynam. Systems}, 24:1895--1908, 2004.

\bibitem[CN05a]{propertyA}
S.~Campbell and G.~Niblo.
\newblock Hilbert space compression and exactness of discrete groups.
\newblock {\em Journal of Functional Analysis}, 22:292--305, 2005.

\bibitem[CN05b]{ChatterjiNibloWallspaces}
I.~Chatterji and G.~Niblo.
\newblock From wall spaces to \rm{CAT}(0) cube complexes.
\newblock {\em International Journal of Algebra and Computation},
  15(5-6):875--885, 2005.

\bibitem[Cor13]{CornulierCommensurated}
Y.~Cornulier.
\newblock Group actions with commensurated subsets, wallings and cubings.
\newblock {\em arxiv:1302.5982}, 2013.

\bibitem[CS11]{MR2827012}
Pierre-Emmanuel Caprace and Michah Sageev.
\newblock Rank rigidity for {CAT}(0) cube complexes.
\newblock {\em Geom. Funct. Anal.}, 21(4):851--891, 2011.

\bibitem[CSV08]{CSVparticulier}
Y.~Cornulier, Y.~Stalder, and A.~Valette.
\newblock Proper actions of lamplighter groups associated with free groups.
\newblock {\em C. {R}. {A}cad. {S}ci. {P}aris, {S}er. {I}}, 346:173--176, 2008.

\bibitem[CSV12]{CSVgeneral}
Y.~Cornulier, Y.~Stalder, and A.~Valette.
\newblock Proper actions of wreath products and generalizations.
\newblock {\em Trans. {A}mer. {M}ath. {S}oc}, 364:3159--3184, 2012.

\bibitem[CW15]{GarsideCurve}
M.~Calvez and B.~Wiest.
\newblock Curve complexes and {G}arside groups.
\newblock {\em arxiv:1503.02482}, 2015.

\bibitem[Dav83]{DavisCovEucl}
M.~Davis.
\newblock Groups generated by reflections and aspherical manifolds not covered
  by {E}uclidean space.
\newblock {\em The Annals of Mathematics, 2nd Ser.}, 117(2):293--324, 1983.

\bibitem[dC06]{CornulierFPWP}
Yves de~Cornulier.
\newblock Finitely presented wreath products and double coset decompositions.
\newblock {\em Geometriae Dedicata}, 122(1):89--108, 2006.

\bibitem[DH89]{Davis1989}
Michael~W. Davis and Jean-Claude Hausmann.
\newblock {\em Aspherical manifolds without smooth or PL structure}, pages
  135--142.
\newblock Springer Berlin Heidelberg, Berlin, Heidelberg, 1989.

\bibitem[DJL12]{DavisJanusLafont}
M.~Davis, T.~Januszkiewicz, and J.-F. Lafont.
\newblock 4-dimensional locally {CAT}(0)-manifolds with no {R}iemannian
  smoothings.
\newblock {\em Duke Math. J.}, 161(1):1--28, 2012.

\bibitem[Djo73]{Djokovic}
D.Z. Djokovic.
\newblock Distance preserving subgraphs of hypercubes.
\newblock {\em Journal of Combinatorial Theory B}, 14:263--267, 1973.

\bibitem[DO01]{DavisOkun}
M.~Davis and B.~Okun.
\newblock Vanishing theorems and conjectures for the $\ell^2$-homology of
  right-angled {C}oxeter groups.
\newblock {\em Geometry and Topology}, 5:7--74, 2001.

\bibitem[Dra08]{AsDimRACG}
A.~Dranishnikov.
\newblock On asymptotic dimension of amalgamated products and right-angled
  {C}oxeter groups.
\newblock {\em Algebraic and Geometric Topology}, 8:1281--1293, 2008.

\bibitem[Dre11]{Dreesen1}
Dennis Dreesen.
\newblock Hilbert space compression for free products and hnn-extensions.
\newblock {\em Journal of Functional Analysis}, 261(12):3585 -- 3611, 2011.

\bibitem[DS04]{RACGproductoftrees}
A.~Dranishnikov and V.~Schroeder.
\newblock Embedding of {C}oxeter groups in aproduct of trees.
\newblock {\em arXiv:0402398}, 2004.

\bibitem[DT14]{DaniThomas}
P.~Dani and A.~Thomas.
\newblock Bowditch's {JSJ} tree and quasi-isometry classification of certain
  {C}oxeter groups.
\newblock {\em arXiv:1402.6224}, 2014.

\bibitem[Dym10]{QIvsLIP}
T.~Dymarz.
\newblock Bilipschitz equivalence is not equivalent to quasi-isometric
  equivalence for finitely generated groups.
\newblock {\em Duke Math. J.}, 154(3):509--526, 2010.

\bibitem[Dyu00]{WreathLip}
Anna Dyubina.
\newblock Instability of the virtual solvability and the property of being
  virtually torsion-free for quasi-isometric groups.
\newblock {\em International Mathematics Research Notices},
  2000(21):1097--1101, 2000.

\bibitem[EFW13]{RigiditySolLamp}
A.~Eskin, D.~Fisher, and K.~Whyte.
\newblock Coarse differentiation of quasi-isometries {II}: {R}igidity for {S}ol
  and lamplighter groups.
\newblock {\em Annals of Mathematics}, 177(3):869--910, 2013.

\bibitem[EPC{\etalchar{+}}92]{automaticgroups}
David B.~A. Epstein, M.~S. Paterson, J.~W. Cannon, D.~F. Holt, S.~V. Levy, and
  W.~P. Thurston.
\newblock {\em Word Processing in Groups}.
\newblock A. K. Peters, Ltd., Natick, MA, USA, 1992.

\bibitem[Far03]{MR1978047}
Daniel~S. Farley.
\newblock Finiteness and {$\rm CAT(0)$} properties of diagram groups.
\newblock {\em Topology}, 42(5):1065--1082, 2003.

\bibitem[Gen15]{arXiv:1507.01667}
A.~Genevois.
\newblock Hyperplanes of {S}quier's cube complexes.
\newblock {\em arXiv:1507.01667}, 2015.

\bibitem[Gen16a]{article4}
A.~Genevois.
\newblock Acylindrical action on the hyperplanes of a {CAT}(0) cube complex.
\newblock {\em arxiv:1610.08759}, 2016.

\bibitem[Gen16b]{coningoff}
A.~Genevois.
\newblock Coning-off \rm{CAT(0)} cube complexes.
\newblock {\em arXiv:1603.06513}, 2016.

\bibitem[Gen16c]{article3}
A.~Genevois.
\newblock Contracting isometries of \rm{CAT}(0) cube complexes and acylindrical
  hyperbolicity of diagram groups.
\newblock {\em arXiv:1610.07791}, 2016.

\bibitem[Gen17a]{arXiv:1505.02053}
A.~Genevois.
\newblock Hyperbolic diagram groups are free.
\newblock {\em Geometriae Dedicata}, 188(1):33--50, Jun 2017.

\bibitem[Gen17b]{aTmW}
A.~Genevois.
\newblock Lampligther groups, median spaces, and a-{T}-menability.
\newblock {\em arXiv:1705.00834}, 2017.

\bibitem[Ger94]{GerstenAutFn}
S.~M. Gersten.
\newblock The automorphism group of a free group is not a \rm{CAT}(0) group.
\newblock {\em Proceedings of the American Mathematical Society},
  121(4):999--1002, 1994.

\bibitem[GH10]{WeakAmenabilityCCGuentner}
Erik Guentner and Nigel Higson.
\newblock Weak amenability of \rm{CAT}(0)-cubical groups.
\newblock {\em Geometriae Dedicata}, 148(1):137--156, 2010.

\bibitem[Gre90]{GreenGP}
E.~Green.
\newblock Graph products of groups.
\newblock {\em PhD Thesis}, 1990.

\bibitem[Gro87]{Gromov1987}
M.~Gromov.
\newblock {\em Hyperbolic Groups}, pages 75--263.
\newblock Springer New York, New York, NY, 1987.

\bibitem[GS97]{MR1396957}
Victor Guba and Mark Sapir.
\newblock Diagram groups.
\newblock {\em Mem. Amer. Math. Soc.}, 130(620):viii+117, 1997.

\bibitem[GS99]{MR1725439}
V.~S. Guba and M.~V. Sapir.
\newblock On subgroups of the {R}. {T}hompson group {$F$} and other diagram
  groups.
\newblock {\em Mat. Sb.}, 190(8):3--60, 1999.

\bibitem[GS06]{MR2193191}
V.~S. Guba and M.~V. Sapir.
\newblock Diagram groups are totally orderable.
\newblock {\em J. Pure Appl. Algebra}, 205(1):48--73, 2006.

\bibitem[Hag08]{MR2413337}
Fr{\'e}d{\'e}ric Haglund.
\newblock Finite index subgroups of graph products.
\newblock {\em Geom. Dedicata}, 135:167--209, 2008.

\bibitem[Hag12]{Hagenthesis}
Mark~F. Hagen.
\newblock Geometry and combinatorics of cube complexes.
\newblock {\em Thesis}, 2012.

\bibitem[Hag14]{MR3217625}
Mark~F. Hagen.
\newblock Weak hyperbolicity of cube complexes and quasi-arboreal groups.
\newblock {\em J. Topol.}, 7(2):385--418, 2014.

\bibitem[Hig76]{Higgins}
P.~J. Higgins.
\newblock The fundamental groupoid of a graph of groups.
\newblock {\em Journal of the London Mathematical Society}, s2-13(1):145--149,
  1976.

\bibitem[HM95]{HermillerMeier}
S.~Hermiller and J.~Meier.
\newblock Algorithms and geometry for graph products of groups.
\newblock {\em J. Algebra}, 171(1):230--257, 1995.

\bibitem[HW99a]{HsuWiseNonRF}
T~Hsu and D.~Wise.
\newblock A non-residually finite square of finite groups.
\newblock {\em Groups St. Andrews 1997 in Bath, I, London Math. Soc., Lecture
  Note Ser., 260}, pages 368--378, 1999.

\bibitem[HW99b]{HsuWise}
T.~Hsu and D.~Wise.
\newblock On linear and residual properties of graph products.
\newblock {\em Michigan Math. J.}, 46:251--259, 1999.

\bibitem[HW08]{MR2377497}
Fr{\'e}d{\'e}ric Haglund and Daniel~T. Wise.
\newblock Special cube complexes.
\newblock {\em Geom. Funct. Anal.}, 17(5):1551--1620, 2008.

\bibitem[HW12]{HruskaWise}
G.~Hruska and D.~Wise.
\newblock Finiteness properties of cubulated groups.
\newblock {\em arXiv:1209.1074}, 2012.

\bibitem[JS01]{GPcommensurability}
T.~Januszkiewicz and J.~Swiatkowski.
\newblock Commensurability of graph products.
\newblock {\em Algebraic and Geometric Topology}, 1:587--603, 2001.

\bibitem[JS03]{HypRACGlargedim}
T.~Januszkiewicz and J.~Swiatkowski.
\newblock Hyperbolic {C}oxeter groups of large dimension.
\newblock {\em Commentarii Mathematici Helvetici}, 78(3):555--583, 2003.

\bibitem[Kim12]{GPsurfacesub}
S.-H. Kim.
\newblock Surface subgroups of graph products of groups.
\newblock {\em International Journal of Algebra and Computation}, 22(8), 2012.

\bibitem[KK13]{embeddingRAAG}
S.-H. Kim and T.~Koberda.
\newblock Embedability between right-angled {A}rtin groups.
\newblock {\em Geom. {T}opol.}, 17(1):493--530, 2013.

\bibitem[KK14]{extensiongraph}
S.-H. Kim and T.~Koberda.
\newblock The geometry of the curve graph of a right-angled {A}rtin group.
\newblock {\em Int. {J}. {A}lgebra {C}omput.}, 24(2):121--169, 2014.

\bibitem[KM16]{GraphWreathProducts}
P.H. Kropholler and A.~Martino.
\newblock Graph-wreath products and finiteness conditions.
\newblock {\em Journal of Pure and Applied Algebra}, 220(1):422 -- 434, 2016.

\bibitem[Kob12]{pingponglemmas}
T.~Koberda.
\newblock Ping-pong lemmas with applications to geoemtry and topology.
\newblock In {\em Geometry, topology and dynamics of character varieties},
  volume~23 of {\em Lect. Notes Ser. Inst. Math. Sci. Natl. Univ. Singap.},
  pages 139--158. Wold Sci. Publ., Hackensack, NJ, 2012.

\bibitem[Lea13]{Leary}
Ian~J. Leary.
\newblock A metric {K}an-{T}hurston theorem.
\newblock {\em Journal of Topology}, 6(1):251--284, 2013.

\bibitem[Lev16]{divCCC}
I.~Levcovitz.
\newblock Divergence of {CAT}(0) cube complexes and {C}oxeter groups.
\newblock {\em arxiv:1611.04378}, 2016.

\bibitem[Li10]{LiWreathProducts}
S.~Li.
\newblock Compression bounds for wreath products.
\newblock {\em Proceedings of the American Mathematical Society},
  138(8):2701--2714, 2010.

\bibitem[Mei94]{MeierGPtopology}
John Meier.
\newblock The topology of graph products of groups.
\newblock {\em Proceedings of the Edinburgh Mathematical Society},
  37(3):539--544, 10 1994.

\bibitem[Mei95]{BNSinvGP}
Holger Meinert.
\newblock The {B}ieri-{N}eumann-{S}trebel invariant for graph products of
  groups.
\newblock {\em Journal of Pure and Applied Algebra}, 103(2):205 -- 210, 1995.

\bibitem[Mei96]{MeierGP}
J.~Meier.
\newblock When is a graph product of hyperbolic groups hyperbolic?
\newblock {\em Geometriae Dedicata}, 61:29--41, 1996.

\bibitem[Miz08]{WeakAmenabilityCCMizuta}
Naokazu Mizuta.
\newblock A {B}ozejko-{P}icardello type inequality for finite-dimensional
  \rm{CAT}(0) cube complexes.
\newblock {\em Journal of Functional Analysis}, 254(3):760--772, 2008.

\bibitem[MO13]{arXiv:1310.6289}
A.~Minasyan and D.~Osin.
\newblock Acylindrical hyperbolicity of groups acting on trees.
\newblock {\em arXiv:1310.6289}, 2013.

\bibitem[Mul80]{Mulder}
H.M. Mulder.
\newblock {\em The interval function of a graph}, volume 132 of {\em Math.
  Centre Tracts}.
\newblock Mathematisch Centrum, Amsterdam, 1980.

\bibitem[Neb71]{NebeskyMedian}
L.~Nebesk\'y.
\newblock Median graphs.
\newblock {\em Commentationes Mathematicae Universitatis Carolinae},
  12(2):317--325, 1971.

\bibitem[Nib02]{SplittingObstruction}
Graham~A. Niblo.
\newblock The singularity obstruction for group splittings.
\newblock {\em Topology and its Applications}, 119:17--31, 2002.

\bibitem[Nic04]{NicaWallspaces}
B.~Nica.
\newblock Cubulating spaces with walls.
\newblock {\em Algebraic and Geometric Topology}, 4:297--309, 2004.

\bibitem[NP08]{NaorPeres}
Assaf Naor and Yuval Peres.
\newblock Embeddings of discrete groups and the speed of random walks.
\newblock {\em International Mathematics Research Notices}, 2008:rnn076, 2008.

\bibitem[NP11]{NaorYuval1}
Assaf Naor and Yuval Peres.
\newblock $l_p$ compression, traveling salesmen, and stable walks.
\newblock {\em Duke Math. J.}, 157(1):53--108, 03 2011.

\bibitem[NR97]{Haagerup}
G.~Niblo and L.~Reeves.
\newblock Groups acting on {CAT}(0) cube complexes.
\newblock {\em Geometry and Topology}, 1:1--7, 1997.

\bibitem[NR98a]{NibloReeves}
G.~Niblo and L.~Reeves.
\newblock The geometry of cube complexes and the complexity of their
  fundamental groups.
\newblock {\em Topology}, 37(3):621--633, 1998.

\bibitem[NR98b]{MR1459140}
Graham~A. Niblo and Martin~A. Roller.
\newblock Groups acting on cubes and {K}azhdan's property ({T}).
\newblock {\em Proc. Amer. Math. Soc.}, 126(3):693--699, 1998.

\bibitem[OP10]{OzawaPopa}
N.~Ozawa and S.~Popa.
\newblock On a class of {II1} factors with at most one {C}artan subalgebra.
\newblock {\em Annals of Mathematics}, 172(1):713--749, 2010.

\bibitem[Osi06]{OsinRelativeHyp}
D.~Osin.
\newblock Relatively hyperbolic groups: intrinsic geometry, algebraic
  properties, and algorithmic problems.
\newblock {\em Mem. Amer. Math. Soc.}, 179(843):vi--100, 2006.

\bibitem[Osi16]{arXiv:1304.1246}
D.~Osin.
\newblock Acylindrically hyperbolic groups.
\newblock {\em Trans. Amer. Math. Soc.}, 368:851--888, 2016.

\bibitem[PS09]{HypBoundRACG}
P.~Przytycki and J.~Swiatkowski.
\newblock Flag-no-square triangulations and {G}romov boundaries in dimension 3.
\newblock {\em Groups Geom. Dyn.}, 3:453--468, 2009.

\bibitem[Rol98]{Roller}
M.~Roller.
\newblock Pocsets, median algebras and group actions: {A}n extended study of
  {D}unwoody's construction and {S}ageev's theorem.
\newblock {\em dissertation}, 1998.

\bibitem[Sag95]{MR1347406}
Michah Sageev.
\newblock Ends of group pairs and non-positively curved cube complexes.
\newblock {\em Proc. London Math. Soc. (3)}, 71(3):585--617, 1995.

\bibitem[Sag14]{SageevCAT(0)}
Michah Sageev.
\newblock \rm{CAT}(0) cube complexes and groups.
\newblock In {\em Geometric Group Theory}, volume~21 of {\em IAS/Park City
  Math. Ser.}, pages 7--54. 2014.

\bibitem[Sap14]{SapirRD}
M.~Sapir.
\newblock The {R}apid {D}ecay property and centroids in groups.
\newblock {\em arxiv:1405.0757}, 2014.

\bibitem[Ser03]{MR1954121}
Jean-Pierre Serre.
\newblock {\em Trees}.
\newblock Springer Monographs in Mathematics. Springer-Verlag, Berlin, 2003.
\newblock Translated from the French original by John Stillwell, Corrected 2nd
  printing of the 1980 English translation.

\bibitem[SW05]{alternative}
M.~Sageev and D.~Wise.
\newblock The {T}its alternative for {CAT}(0) cubical complexes.
\newblock {\em Bulletin of the London Mathematical Society}, 37(5):706--720,
  2005.

\bibitem[Swi16]{RACGboundSier}
J.~Swiatkowski.
\newblock Right-angled {C}oxeter groups with n-dimensional {S}ierpinski
  compacta as boundaries.
\newblock {\em arXiv:1603.02152}, 2016.

\bibitem[Tes11]{TesseraWP}
R.~Tessera.
\newblock Asymptotic isoperimetry on groups and uniform embeddings into
  {B}anach spaces.
\newblock {\em Commentarii Mathematici Helvetici}, 86(3):499--535, 2011.

\bibitem[Why99]{WhyteLip}
K.~Whyte.
\newblock Amenability, bilipschitz equivalence, and the von {N}eumann
  conjecture.
\newblock {\em Duke Math J.}, 99(1):93--112, 1999.

\bibitem[Wie03]{WiestDiagram}
B.~Wiest.
\newblock Diagram groups, braid groups, and orderability.
\newblock {\em Journal of Knot Theory and Its Ramifications}, 12(03):321--332,
  2003.

\bibitem[Wis12]{bigWise}
D.~T. Wise.
\newblock The structure of groups with a quasiconvex hierarchy.
\newblock {\em Preprint}, 2012.

\end{thebibliography}

\addcontentsline{toc}{section}{Index}

\printindex

\end{document}